\documentclass[a4paper,10pt]{smfbook}
\usepackage[french,english]{babel}
\usepackage{smfthm}

\usepackage[T1]{fontenc}
\usepackage[utf8]{inputenc}
\usepackage{amsmath}
\usepackage{array}
\usepackage{amsthm}
\usepackage{amssymb}
\usepackage{stmaryrd}
\input xy
\xyoption{all}
\usepackage{hyperref}
\usepackage{multind}
\makeindex{nota}
\makeindex{termin}

\newcommand{\A}{{\mathcal{A}}}
\newcommand{\B}{{\mathcal{B}}}
\newcommand{\C}{{\mathcal{C}}}
\newcommand{\D}{{\mathcal{D}}}
\newcommand{\F}{{\mathcal{F}}}
\newcommand{\Oo}{{\mathfrak{D}}}
\newcommand{\fct}{{\mathbf{Fct}}}

\newcommand{\Pp}{{\mathcal{P}}}

\newcommand{\pol}{{\mathcal{P}ol}}
\newcommand{\E}{{\mathcal{E}}}

\newcommand{\I}{{\mathcal{I}}}

\newcommand{\M}{{\mathcal{M}}}

\newcommand{\G}{{\mathcal{G}}}
\newcommand{\p}{{\mathfrak{p}}}
\newcommand{\pp}{{\mathfrak{P}}}
\newcommand{\col}{{\rm colim}\,}

\newcommand{\W}{{\mathcal{W}}}

\newcommand{\Md}{\text{-}\mathbf{Mod}}

\newcommand{\Si}{\mathfrak{S}}
\newcommand{\FF}{\mathbb{F}}

\newcommand{\ph}{Hom-polynomial}
\newcommand{\phs}{Hom-polynomiaux}
\newcommand{\nil}{\mathrm{Nil}}
\newcommand{\rj}{\operatorname{J}}
\newcommand{\rf}{\mathfrak{r}}
\newcommand{\AF}{A_{/F}}
\newcommand{\hi}{\underline{\mathrm{Hom}}}
\newcommand{\ti}{\underline{\otimes}}
\newcommand{\Sp}{\mathbf{sp}}
\newcommand{\GL}{\operatorname{GL}}
\newcommand{\SL}{\operatorname{SL}}
\newcommand{\PSL}{\operatorname{PSL}}
\newcommand{\spec}{\operatorname{Spec}}
\newcommand{\cd}{\operatorname{Card}}
\newcommand{\mon}{\mathbf{Mon}}
\newcommand{\ac}{\mathfrak{A}}
\newcommand{\ic}{\mathfrak{I}}
\newcommand{\op}{\mathrm{op}}
\newcommand{\hd}{hyper-diagonalisant}
\newcommand{\hds}{hyper-diagonalisants}
\newcommand{\htr}{hyper-trigonalisant}
\newcommand{\hts}{hyper-trigonalisants}
\newcommand{\cf}{\mathrm{cf}}
\newcommand{\df}{\mathrm{df}}
\newcommand{\Tt}{\mathfrak{T}}
\newcommand{\la}{\operatorname{z}}
\newcommand{\ppol}{\operatorname{p}}
\newcommand{\qpol}{\operatorname{q}}
\newcommand{\pf}{\mathbf{pf}}

\author[A. DJAMENT]{Aur\'elien DJAMENT}
\address{CNRS, laboratoire Analyse, Géométrie et Applications (UMR7539)\\ Institut Galilée\\ Université Sorbonne Paris Nord\\ 99 avenue Jean-Baptiste Clément\\ 93430 VILLETANEUSE\\ FRANCE}
\email{djament@math.cnrs.fr}
\urladdr{https://djament.perso.math.cnrs.fr/}

\author[A. TOUZÉ]{Antoine TOUZÉ}
\address{Université de Lille\\ laboratoire Paul Painlev\'e (UMR8524)\\ Cité scientifique, bât. M2\\ 59655 VILLENEUVE D'ASCQ CEDEX\\ FRANCE}
\email{antoine.touze@univ-lille.fr}
\urladdr{https://pro.univ-lille.fr/antoine-touze/}

\title[Représentations génériques des groupes linéaires infinis]{Sur la structure des représentations génériques des groupes linéaires infinis}

\alttitle{On the structure of generic representations of infinite general linear groups}

\newtheorem{thi}{Th\'eor\`eme}

\newtheorem{pri}[thi]{Proposition}

\newtheorem{cvi}{Convention}[chapter]
\newtheorem{noti}[cvi]{Notation}

\begin{document}
\newenvironment{nota}{\begin{enonce}{Notation}}{\end{enonce}}
\newenvironment{hyp}{\begin{enonce}{Hypoth\`ese}}{\end{enonce}}
\newenvironment{conv}{\begin{enonce}{Convention}}{\end{enonce}}
\newenvironment{prdef}{\begin{enonce}{Proposition et d\'efinition}}{\end{enonce}}

\frontmatter

\begin{abstract}
Nous étudions plusieurs aspects structuraux des catégories de foncteurs d'une petite catégorie additive vers une catégorie de modules, notamment de la catégorie $\F(A,K)$ des foncteurs des modules libres de rang fini sur un anneau commutatif $A$ vers les espaces vectoriels sur un corps commutatif $K$ --- de tels foncteurs sont parfois appelés \textit{représentations génériques} des groupes linéaires sur $A$ à coefficients dans $K$. Nous nous intéressons en particulier aux foncteurs de type fini de $\F(A,K)$ prenant des valeurs de dimensions finies. Nous montrons qu'ils se ramènent moyennant une légère hypothèse supplémentaire (toujours vérifiée si $A$ est un anneau noethérien) à des foncteurs bien mieux compris, à savoir des foncteurs \textit{polynomiaux} (à la Eilenberg-MacLane), ou se factorisant à la source par la réduction modulo un idéal cofini de $A$. Nous en déduisons que de tels foncteurs sont toujours noethériens et que, si l'anneau $A$ est de type fini, ils possèdent des résolutions projectives de type fini.

Nos méthodes reposent principalement sur l'étude de \textit{décompositions en poids} des foncteurs et de leurs effets croisés, sur notre travail antérieur avec Vespa \cite{DTV} et de l'algèbre commutative élémentaire.
\end{abstract}

\begin{altabstract}
We study several structure aspects of functor categories from a small additive category to a module category, in particular the category $\F(A,K)$ of functors from finitely generated free modules over a commutative ring $A$ to vector spaces over a field $K$ --- such functors are sometimes called \textit{generic representations} of linear groups over $A$ with coefficients in $K$. We are especially interested with finitely generated functors of $\F(A,K)$ taking finite dimensional values. We prove that they can, under a mild extra assumption (always satisfied if the ring $A$ is noetherian), be built from much better understood functors, namely \textit{polynomial} functors (in the sense of Eilenberg-MacLane), or factorising at the source through reduction modulo a cofinite ideal of $A$. We deduce that such functors are always noetherian et that, if the ring $A$ is finitely generated, they have finitely generated projective resolutions.

Our methods rely mainly on the study of \textit{weight decompositions} of functors and their cross-effects, our previous work with Vespa \cite{DTV} and elementary commutative algebra. 
\end{altabstract}

\subjclass{18A25, 18E05, 18E10; 13A99, 20C99, 20M14, 20M25.}

\keywords{catégories de foncteurs, foncteurs polynomiaux, décompositions en poids, foncteurs de type fini à valeurs de dimensions finies, noethérianité.}

\altkeywords{functor categories, polynomial functors, weight decompositions, finitely generated functors with finite-dimensional values, noetherianess.}

\thanks{Les auteurs ont bénéficié du soutien partiel de l’Agence Nationale de la Recherche via le Labex CEMPI (ANR-11-LABX-0007-01)}

\maketitle

\tableofcontents

\mainmatter

\chapter*{Introduction}

\section*{Représentations génériques des groupes linéaires}

Soient $A$ un anneau, $k$ un anneau noethérien et $K$ un corps, que nous supposerons tous trois \textbf{commutatifs}. Cet article porte principalement sur la catégorie $\F(A,k)$ des foncteurs depuis la catégorie $\mathbf{P}(A)$ des $A$-modules libres de rang fini vers la catégorie $k\Md$ de tous les $k$-modules. La catégorie $\F(A,k)$ est abélienne et possède de nombreuses propriétés et structures agréables (c'est une catégorie de Grothendieck localement de présentation finie avec assez de projectifs et d'injectifs, munie d'une structure tensorielle donnée par le produit tensoriel au-dessus de $k$ au but) ; on la nomme, depuis la série d'articles de N. Kuhn \cite{Ku1,Ku2,Ku3}, consacrée au cas où $A=k$ est un corps fini, catégorie des {\em représentations génériques des groupes linéaires} sur $A$ à coefficients dans $k$.  De fait, on peut reconstruire un foncteur $F$ de $\F(A,k)$ à partir de la suite des $k$-représentations $F(A^n)$ des groupes $\GL_n(A)$ et de quelques données supplémentaires ; inversement, l'étude de $\F(A,k)$ permet parfois d'obtenir des résultats sur les représentations des groupes $\GL_n(A)$, au moins dans le domaine \guillemotleft~stable~\guillemotright\ (c'est-à-dire pour $n$ \guillemotleft~assez grand~\guillemotright), notamment en matière cohomologique (voir par exemple le travail remarquable de Scorichenko \cite{Sco}).

Ainsi, si $A$ est un corps, il existe une bijection explicite entre la réunion disjointe sur $n\in\mathbb{N}$ de l'ensemble des classes d'isomorphisme de représentations $K$-linéaires irréductibles de $\GL_n(A)$ et l'ensemble des classes d'isomorphisme de foncteurs simples de $\F(A,K)$, donnée par la théorie des recollements. Toutefois, les méthodes employées pour étudier les représentations des groupes linéaires et les représentations génériques diffèrent largement, ce qui donne lieu à des interactions fécondes.

Dans ce mémoire, nous établirons des propriétés de finitude et de structure de vastes classes de foncteurs de $\F(A,K)$ ; nous généraliserons certaines d'entre elles à $\F(A,k)$ dans le dernier chapitre. Nous étudierons en particulier la propriété noethérienne et la propriété de longueur finie, mais aussi la présentation finie et ses généralisations supérieures. Cette étude nous conduira à dégager des résultats de structure beaucoup plus précis, sur lesquels nous reviendrons ultérieurement dans cette introduction.

L'étude fine de \emph{tous} les foncteurs de $\F(A,K)$, sauf parfois lorsque $A$ est un anneau fini (cf. par exemple \cite{Ku1,Ku2,Ku3,GP-cow,Dja-gr,DG}), s'avère en général hors d'atteinte, de même que la compréhension des représentations linéaires d'un groupe infini assez gros (par exemple, contenant un groupe libre non abélien, ce qui est le cas de $\GL_n(A)$ pour $n>1$ si l'anneau $A$ n'est pas localement fini). Une grande partie de ce travail sera consacrée aux foncteurs de type fini prenant des valeurs de dimensions finies, qui possèdent un comportement plus accessible, comme notre récent travail avec C. Vespa \cite{DTV} a commencé à le montrer (au prix d'hypothèses de finitude plus fortes), de même que les représentations \emph{de dimension finie} fournissent une classe importante de représentations linéaires d'un groupe infini plus accessibles que toutes les représentations.

D'autres classes fondamentales de foncteurs joueront un rôle important dans ce mémoire, notamment les foncteurs \emph{polynomiaux}, introduits par S. Eilenberg et S. Mac\-Lane \cite{EML} au début des années 1950, qui ont largement montré leur importance en topologie algébrique et en théorie des représentations. Les foncteurs \emph{antipolynomiaux} introduits dans \cite[§\,4.1]{DTV} comme contrepoint naturel aux foncteurs polynomiaux, et dont la structure a été étudiée en détail par Kuhn \cite{Ku-adv} lorsque l'anneau $A$ est un corps fini, et par Djament et Gaujal \cite{DG} dans le cas général, interviendront également. Enfin, nous définirons de nouvelles classes de foncteurs possédant de bonnes propriétés de stabilité (par sous-quotient, colimite, produit tensoriel...) et pour lesquelles nous démontrerons des théorèmes de structure qui nous permettront d'obtenir nos résultats principaux sur les foncteurs de type fini à valeurs de dimensions finies.

\section*{La propriété localement noethérienne}

La propriété noethérienne pour une algèbre de groupe est notoirement très difficile à étudier en général. Toutefois, on peut voir sans trop de peine que, pour tout entier $n\ge 2$, l'algèbre $k[\GL_n(A)]$ est noethérienne si et seulement si $A$ est \emph{fini} (si l'anneau noethérien $k$ est non nul). Dans le cadre des représentations génériques, ce fait possède l'analogue suivant :

\begin{thi}[Putman-Sam-Snowden]\label{thi-PSS} Supposons l'anneau $k$ non nul. Alors la catégorie $\F(A,k)$ est localement noethérienne si et seulement si l'anneau $A$ est fini.
\end{thi}

La noethérianité locale de $\F(A,k)$ pour $A$ fini constitue un résultat profond établi au milieu des années 2010, par Putman-Sam \cite{PSam} d'une part et Sam-Snowden \cite{SamSn} d'autre part. La réciproque constitue une conséquence facile de la non-noethérianité de $k[\GL_2(A)]$ pour $A$ infini et $k$ non nul (cf. \cite[prop.~11.1]{DTV}).

\paragraph*{Cas des foncteurs polynomiaux} La sous-catégorie des foncteurs additifs de $\F(A,k)$ est équivalente à celle des $A_k$-modules, où $A_k:=A\otimes_\mathbb{Z}k$. Plus généralement, il est facile de voir que, si tout foncteur polynomial de degré $d$ de $\F(A,K)$ est localement noethérien, alors la $k$-algèbre $A_k^{\otimes d}$ (où le produit tensoriel est pris sur $k$) est noethérienne. Il semble toutefois très ardu de donner une condition nécessaire et suffisante simple pour que tout foncteur polynomial de $\F(A,k)$ soit localement noethérien, même si $k$ est un corps (toutefois, dans $\F(K,K)$, \cite[appl.~5.9]{DT-TJM} donne une réponse complète simple à cette question).

Notre travail \cite{DT-TJM} montre néanmoins le résultat suivant, qui fournit une vaste classe d'anneaux infinis $A$ pour lesquels les foncteurs polynomiaux de $\F(A,k)$ sont localement noethériens. On rappelle qu'un anneau (commutatif) est dit \emph{essentiellement de type fini} s'il est isomorphe à un localisé d'un anneau de type fini.

\begin{thi}\label{thi-tfdfln}\cite[th.~5.2]{DT-TJM} Si l'anneau $A$ est essentiellement de type fini, alors tout foncteur polynomial de $\F(A,k)$ est localement noethérien.
\end{thi}

\paragraph*{Cas des foncteurs à valeurs de dimensions finies} L'un des résultats principaux de ce mémoire est le suivant, qui ne nécessite aucune hypothèse (sinon la commutativité) sur l'anneau $A$. 

\begin{thi}[théorème~\ref{th-noethgal}]\label{pri-k1} Tout foncteur de type fini à valeurs de type fini de $\F(A,k)$ est noethérien.
\end{thi}

Le cas particulier de ce théorème où $k$ est un corps répond par l'affirmative à une conjecture de \cite[conj.~11.11.1]{DTV} pour le cas des anneaux \emph{commutatifs} (à la source).

\section*{La propriété localement finie}

(On se limite ici sans perte de généralité au cas de $\F(A,K)$, car la propriété localement finie pour les foncteurs constants de $\F(A,k)$ entraîne que l'anneau $k$ est artinien, auquel cas l'étude de propriétés de finitude dans $\F(A,k)$ se ramène facilement à des questions analogues dans $\F(A,K)$, où $K$ est un corps quotient de $A$.)

Si $G$ est un groupe, il est classique que l'algèbre $K[G]$ est artinienne --- ce qui revient à dire que la catégorie des $K[G]$-modules est localement finie (engendrée par des objets de longueur finie) --- si et seulement si $G$ est \emph{fini} \cite[chap.~10, th.~1.1]{Passman}. Il est de plus trivial que tout $K[G]$-module de dimension finie sur $K$ est fini (i.e. de longueur finie). Pour les représentations génériques, la propriété localement finie est plus difficile à étudier, même si les résultats ici mentionnés demandent un peu moins de préparation que ceux relatifs à la propriété localement noethérienne.

\paragraph*{Cas des foncteurs polynomiaux} La propriété suivante (cf. \cite[lemme~11.10]{DTV}) résulte du fait que, pour tout $d\in\mathbb{N}$, la restriction du foncteur d'évaluation en $A^d$ à la sous-catégorie pleine des foncteurs polynomiaux de degré au plus $d$ de $\F(A,K)$ est fidèle.

\begin{pri}\label{pri-poldff} Tout foncteur polynomial à valeurs de dimensions finies de $\F(A,K)$ est fini.
\end{pri}

En général, la réciproque de cette propriété est fausse, y compris pour les foncteurs additifs, car même un $(A,K)$-bimodule \emph{simple} n'est pas nécessairement de dimension finie sur $K$. On a toutefois le résultat suivant, inclus dans \cite[prop.~6.3]{DTV} :

\begin{pri} Si l'anneau $A$ est de type fini, alors tout foncteur polynomial fini de $\F(A,K)$ est à valeurs de dimensions finies.
\end{pri}

Il n'est pas non plus difficile de déterminer quand \emph{tous} les foncteurs polynomiaux de $\F(A,K)$ sont localement finis :

\begin{pri}[proposition~\ref{pr-tspollf}] Les assertions suivantes sont équivalentes :
\begin{enumerate}
    \item tout foncteur polynomial de type fini de $\F(A,K)$ est fini ;
    \item tout foncteur polynomial de type fini de $\F(A,K)$ est à valeurs de dimensions finies ;
    \item l'algèbre $A\otimes_\mathbb{Z}K$ est de dimension finie sur $K$.
\end{enumerate}
\end{pri}

Il existe également des situations où les foncteurs polynomiaux ne sont pas tous localement noethériens, mais où tous les foncteurs polynomiaux noethériens sont finis. C'est notamment le cas lorsque $A=K$ est la clôture algébrique d'un corps fini, grâce au résultat ci-dessous (qui est inclus dans la proposition~\ref{pr-hdcorlf}), combiné à la proposition~\ref{pri-poldff} et à l'observation que $K\otimes_\mathbb{Z} K$ n'est alors pas un anneau noethérien (de sorte qu'il existe des foncteurs additifs non localement noethériens) :

\begin{pri}\label{pri-polnoeth} Si $A$ est un sous-corps localement fini de $K$, alors tout foncteur polynomial noethérien de $\F(A,K)$ est à valeurs de dimensions finies.
\end{pri}

\paragraph*{Cas des foncteurs antipolynomiaux} On rappelle qu'un foncteur de $\F(A,K)$ est dit \emph{antipolynomial}\index{termin}{antipolynomial \emph{(foncteur)}} s'il se factorise à la source par la réduction modulo un idéal $I$ de $A$ tel que $A/I$ soit de cardinal fini et inversible dans $K$. Dans \cite[cor.~11.8]{DTV}, on observe que l'on peut déduire simplement du théorème~\ref{thi-PSS} que tout foncteur antipolynomial de $\F(A,K)$ est localement fini. Comme réciproquement, si $k$ est un corps fini de même caractéristique que $K$, il est classique que la catégorie $\F(k,K)$ n'est pas localement finie (sa dimension de Krull est infinie), on a :

\begin{pri}\label{prildffac}
Les assertions suivantes sont équivalentes :
\begin{enumerate}
    \item la catégorie $\F(A,K)$ est localement finie ;
    \item l'anneau $A$ est de cardinal fini et inversible dans $K$.
\end{enumerate}
\end{pri}

\paragraph*{Cas des foncteurs à valeurs de dimensions finies} L'un des résultats principaux de ce mémoire est le suivant, qui répond positivement à une autre conjecture de \cite[conj.~11.11.2]{DTV} pour le cas des anneaux commutatifs  :
\begin{thi}[théorème~\ref{th-PAP_df}]\label{thi-fdflf}
Les assertions suivantes sont équivalentes :
\begin{enumerate}
    \item\label{itprif1} tout foncteur à valeurs de dimensions finies de $\F(A,K)$ est localement fini ;
    \item\label{itprif2} l'anneau $A$ n'a pas de quotient fini de même caractéristique que $K$.
\end{enumerate}
\end{thi}

(L'implication \ref{itprif1}$\Rightarrow$\ref{itprif2} découle de la proposition~\ref{prildffac}.)

\section*{Résolutions projectives de type fini}

\begin{thi}[corollaire~\ref{cor-pfi_an_gal}]\label{pri-k2} Si l'anneau $A$ est essentiellement de type fini, alors tout foncteur de type fini à valeurs de type fini de $\F(A,k)$ possède une résolution projective de type fini. 
\end{thi}

Lorsque $k$ est un corps et que le foncteur considéré est \emph{fini}, ce résultat a été établi dans \cite[th.~1]{DT-schw}. Pour des foncteurs \emph{polynomiaux}, il a été démontré, sans hypothèse de finitude sur les valeurs, dans \cite[th.~5.2]{DT-TJM}.

\section*{Structure des foncteurs de type fini à valeurs de dimensions finies}

Notre premier théorème de structure est un résultat de polynomialité qui englobe par exemple le cas où $A$ est une algèbre sur un corps infini. Il vaut dans $\F(A,k)$ pour tout anneau noethérien (non nul) $k$.

\begin{thi}[théorème~\ref{th-asqf}]\label{pri-k3} Supposons l'anneau $k$ non nul.  Les assertions suivantes sont équivalentes :
\begin{enumerate}
    \item tout foncteur de type fini à valeurs de dimensions finies de $\F(A,k)$ est polynomial ;
    \item l'anneau $A$ n'a pas de quotient fini non nul.
\end{enumerate}
\end{thi}

Le résultat suivant est à la base du théorème~\ref{thi-fdflf} et ne nécessite pas non plus de travailler avec des espaces vectoriels au but.

\begin{thi}[corollaire~\ref{cor-pqfe_ann}]\label{pri-k4} Si tous les quotients finis de l'anneau $A$ sont de cardinal inversible dans $k$, alors pour tout foncteur $F$ de type fini et à valeurs de type fini de $\F(A,k)$, il existe un foncteur $B : \mathbf{P}(A)\times\mathbf{P}(A)\to k\Md$, unique à isomorphisme près, tel que :
\begin{enumerate}
\item $F$ est isomorphe à la composée de $B$ et du foncteur diagonale $\mathbf{P}(A)\to\mathbf{P}(A)\times\mathbf{P}(A)$
;
\item $B$ est polynomial par rapport à la première variable ;
\item par rapport à la deuxième variable, $B$ se factorise à travers le foncteur de réduction modulo un idéal cofini de $A$.
\end{enumerate}

De plus, $B$ est de type fini et à valeurs de type fini.
\end{thi}

L'énoncé suivant peut sembler analogue au précédent mais il convient de noter plusieurs différences : il ne se généralise pas (de même que le théorème~\ref{thi-spa} ci-après) au cas où l'on remplace le corps $K$ au but par un anneau noethérien, même dans $\F(\mathbb{Z},\mathbb{Z})$ (cf. exemple~\ref{cex_dec-prim}). De plus, on perd l'unicité du bifoncteur $B$ qui apparaît ; en revanche, il en existe une forme plus précise (cf. théorème~\ref{th-tfdf_gal}) --- également à la base du théorème~\ref{thi-spa} ci-dessous --- qui permet de retrouver une forme d'unicité, et aussi d'obtenir directement des propriétés de finitude de $F$ à partir de propriétés analogues pour $B$, ce qui n'est pas possible avec le seul énoncé ci-dessous.

\begin{thi}[corollaire~\ref{cor-imq}]\label{thi-strtfdf} Si l'anneau $A$ est noethérien, alors pour tout foncteur $F$ de type fini et à valeurs de dimensions finies de $\F(A,K)$, il existe un idéal \textbf{cofini} $I$ de $A$ et un foncteur $B : \mathbf{P}(A)\times\mathbf{P}(A/I)\to K\Md$, de type fini et à valeurs de dimensions finies, tels que :
\begin{enumerate}
\item $F$ est isomorphe à la composée de $B$ et du foncteur canonique $\mathbf{P}(A)\to\mathbf{P}(A)\times\mathbf{P}(A/I)$ 
;
\item $B$ est polynomial par rapport à la première variable.
\end{enumerate}
\end{thi}

Le résultat qui précède (ou plutôt la forme plus précise susmentionnée) permet de démontrer le théorème~\ref{thi-tfdfln} à partir du théorème~\ref{thi-PSS} et de la proposition~\ref{pri-poldff}, après s'être ramené au cas où l'anneau $A$ est noethérien (c'est facile) et où $K$ est un corps (c'est plus difficile).

Le théorème suivant (là encore, le théorème~\ref{th-tfdf_gal} est plus précis et inclut une forme d'unicité) fournit une \guillemotleft~décomposition primaire~\guillemotright\ pour les foncteurs de type fini à valeurs de dimensions finies de $\F(A,K)$.

\begin{thi}[théorème~\ref{th-tfdf_gal}]\label{thi-spa} Soit $F$ un foncteur de type fini et à valeurs de dimensions finies de $\F(A,K)$. Alors il existe des idéaux premiers $\mathfrak{p}_1,\dots,\mathfrak{p}_r$ de $A$, $r\in\mathbb{N}$ et un foncteur $X : \mathbf{P}(A_{\mathfrak{p}_1}/\mathfrak{p}_1^r.A_{\mathfrak{p}_1})\times\dots\times\mathbf{P}(A_{\mathfrak{p}_1}/\mathfrak{p}_n^r.A_{\mathfrak{p}_n})\to K\Md$, de type fini et à valeurs de dimensions finies, tels que $F$ soit isomorphe à la composée de $X$ et du foncteur canonique $\mathbf{P}(A)\to\mathbf{P}(A_{\mathfrak{p}_1}/\mathfrak{p}_1^r.A_{\mathfrak{p}_1})\times\dots\times\mathbf{P}(A_{\mathfrak{p}_1}/\mathfrak{p}_n^r.A_{\mathfrak{p}_n})$.
\end{thi}

Dans certains cas, on peut aller encore plus pour préciser la structure des foncteurs de type fini à valeurs de dimensions finies. C'est notamment le cas pour les endofoncteurs des espaces vectoriels de dimensions finies sur un corps algébriquement clos de caractéristique première, qui se ramènent à des (multi)foncteurs polynomiaux \emph{stricts}, au sens de Friedlander-Suslin \cite{FS} (dont la structure est beaucoup mieux comprise que celle des foncteurs polynomiaux ordinaires) et à de la théorie de Galois grâce au résultat suivant.

\begin{thi}[proposition~\ref{pr-corpspftcarp}]\label{thi-Fpbar} Supposons que $K$ et $L$ sont des corps de même caractéristique première, avec $K$ algébriquement clos et $L$ parfait infini. Si $F$ est un foncteur de type fini et à valeurs de dimensions finies de $\F(L,K)$, alors il existe des morphismes de corps $\varphi_1,\dots,\varphi_n : L\to K$  et un multifoncteur polynomial \textbf{strict} $M$ en $n$ variables sur les $K$-espaces vectoriels tels que $F\simeq M\circ(\varphi_1^*,\dots,\varphi_n^*)$, où $\varphi_i^* : \mathbf
P(L)\to\mathbf{P}(K)$ désigne le foncteur d'extension des scalaires associé à $\varphi_i$.
\end{thi}

\section*{Compléments}

Plusieurs des résultats donnés dans ce mémoire sont plus généraux que les formes présentées dans cette introduction ; notamment, on peut parfois remplacer la catégorie des $A$-modules libres de rang fini par une catégorie additive $A$-linéaire essentiellement petite plus générale, voire arbitraire. On peut également généraliser, pour certains résultats, la catégorie des $K$-espaces vectoriels au but par une catégorie abélienne de Grothendieck $K$-linéaire quelconque.

En revanche, \emph{l'hypothèse de commutativité des anneaux considérés est essentielle} ; l'hypothèse de commutativité du corps $K$ au but (qui permet d'appliquer les techniques usuelles de réduction des endomorphismes des espaces vectoriels de dimension finie) est déjà cruciale dans \cite{DTV}, sur lequel nous nous appuyons à de nombreuses reprises. Ainsi, si $L$ est un corps qui n'est de dimension finie sur aucun sous-corps commutatif, nous ignorons à peu près tout (noethérianité locale, structure des foncteurs simples...) des endofoncteurs des $L$-espaces vectoriels de dimensions finies.

Même dans le cadre de la catégorie de foncteurs $\F(A,K)$ (où $A$ est un anneau commutatif et $K$ un corps commutatif), plusieurs problèmes de finitude sur les foncteurs à valeurs de dimensions finies demeurent ouverts. En particulier, si $A$ est \emph{fini}, les méthodes ici développées ne donnent guère que des résultats triviaux\,\footnote{C'est aussi le cas lorsque le quotient de $A$ par son nilradical est fini et de même caractéristique que $K$ --- cf. remarque~\ref{rq-cas_embetant}.} ; elles servent plutôt à \emph{se ramener} dans certains cas à des foncteurs de $\F(A',K)$ où $A'$ est un anneau fini (et à d'autres types de foncteurs mieux compris comme les foncteurs polynomiaux) --- cf. théorème~\ref{thi-strtfdf}. Lorsque $A$ est fini, deux situations fondamentales très différentes se présentent : celle d'inégale caractéristique, où le cardinal de $A$ est inversible dans $K$, dans laquelle \cite{DG} décrit assez finement la structure des foncteurs de type fini de $\F(A,K)$, et celle où la caractéristique de $A$ est une puissance de celle de $K$, dont la structure reste largement inconnue --- la description hypothétique de sa filtration de Krull donnée dans \cite[conj.~12.2.1 et prop.~12.2.3]{Dja-gr} (lorsque $A=K$ est un corps fini) n'a fait l'objet d'aucune avancée depuis une quinzaine d'années.

Parmi les autres problèmes de finitude que nous ne savons pas résoudre, mentionnons également celui-ci :
la sous-catégorie pleine des foncteurs de $\F(A,K)$ constituée des foncteurs localement à valeurs de dimensions finies (i.e. dont tout sous-foncteur de type fini est à valeurs de dimensions finies) possède-t-elle une famille cogénératrice formée de foncteurs artiniens ? (La question est particulièrement naturelle lorsque $A=K$ est la clôture algébrique d'un corps fini, au vu du théorème~\ref{thi-Fpbar} et du résultat de Kuhn \cite[cor.~5.7]{K97}.)

\section*{Méthodes employées}

Ce mémoire s'inscrit dans le prolongement de \cite{DTV} (auquel nous apportons peu lorsque le corps $K$ est de caractéristique nulle), dont il combine les résultats avec de l'algèbre commutative élémentaire et plusieurs outils nouveaux.

Nous introduisons ainsi, au chapitre~\ref{chap-spec}, les notions de \emph{radical}, de \emph{spectre} et d'\emph{algèbre caractéristique} d'un foncteur de $\F(A,K)$. Mais les outils nouveaux principaux que nous utilisons viennent de la notion de \emph{poids} d'un foncteur, que nous allons maintenant détailler.

\paragraph*{Décomposition en poids des foncteurs}
Ici, un \emph{poids} sera un morphisme \emph{multiplicatif} de $A$ dans $K$. Si $F$ est un foncteur de $\F(A,K)$\,\footnote{Nous nous limitons dans cette introduction au cas fondamental de $\F(A,K)$, mais définissons et étudions les décompositions en poids dans un contexte plus général dans le corps du mémoire.} et $w$ un poids, on dit que $F$ est homogène de poids fort $w$ si $F$ transforme la multiplication par tout scalaire $a\in A$ à la source en la multiplication par $w(a)\in K$ au but.

On dit qu'un foncteur $F$ possède une décomposition en poids forte s'il est somme directe de sous-foncteurs homogènes de certains poids, et qu'il possède une décomposition en poids faible s'il possède une filtration exhaustive dont les sous-quotients sont homogènes de certains poids. Les poids de ces sous-quotients sont appelés simplement poids de $F$.

L'existence d'une décomposition en poids, même faible, est tout sauf automatique. On montre ainsi sans trop de peine que, si tout foncteur de $\F(A,K)$ possède une décomposition en poids faible, alors $A$ est un anneau fini (proposition~\ref{pr-dpftjs}). Réciproquement il est facile de voir que, si l'anneau $A$ est fini et le corps $K$ algébriquement clos, alors tout foncteur de $\F(A,K)$ possède une décomposition en poids faible. Plus important encore est le fait que, si $K$ est algébriquement clos (mais sans hypothèse sur $A$), alors tout foncteur \emph{à valeurs de dimensions finies} de $\F(A,K)$ possède une décomposition en poids faible (corollaire~\ref{cor-tfdp}).

De façon générale, si $A$ est infini, un poids $A\to K$ peut être très \guillemotleft~pathologique~\guillemotright, mais une conséquence importante des résultats de \cite{DTV} est que les poids d'un foncteur à valeurs de dimensions finies de $\F(A,K)$ (qu'on nomme \emph{poids ordinaires}) possèdent une description assez simple, comme produit (essentiellement unique) de morphismes d'anneaux $A\to K$ (quitte à remplacer $K$ par une extension finie de ce corps) et d'un morphisme multiplicatif se factorisant par un quotient fini de $A$. Cela joue un rôle important dans ce travail.

Pour les foncteurs \emph{polynomiaux} de $\F(A,K)$, l'existence d'une décomposition en poids implique des propriétés de structure fortes. Ainsi, on montre qu'un foncteur polynomial noethérien de $\F(A,K)$ possédant une décomposition en poids faible est fini et à valeurs de dimensions finies (théorème~\ref{th-poinoeth}) ; de plus, sous certaines hypothèses supplémentaires (si $A$ est un anneau essentiellement de type fini, ou la clôture algébrique d'un corps fini, par exemple), on peut affaiblir l'hypothèse noethérienne en une hypothèse de type fini. La proposition~\ref{pri-polnoeth} se déduit directement de cela, car l'existence d'une décomposition en poids pour un foncteur polynomial noethérien est automatique dans la situation considérée (cf. corollaire~\ref{cor-poidec}).

\paragraph*{Foncteurs \hds\ et \hts}

Pour des foncteurs quelconques (sans propriété de polynomialité) de $\F(A,K)$, l'existence d'une décomposition en poids (même forte) ne suffit généralement pas à obtenir des résultats de finitude ou de structure significatifs. Par exemple, la structure du plus grand quotient homogène d'un certain poids fort du foncteur projectif de type fini de $\F(A,K)$ composé du foncteur d'oubli de $\mathbf{P}(A)$ vers les ensembles et du foncteur de linéarisation vers les $K$-espaces vectoriels semble en général totalement hors de portée.

Nous montrons en revanche que la condition, pour un foncteur $F$ de $\F(A,K)$, que le bifoncteur $F\circ\oplus$ (où $\oplus : \mathbf{P}(A)\times\mathbf{P}(A)\to\mathbf{P}(A)$ désigne la somme directe), possède une décomposition en poids implique des résultats de structure significatifs. Nous disons que $F$ est \emph{\hd}\ si $F\circ\oplus$ possède une décomposition en poids forte, et que $F$ est \emph{\htr}\ si $F\circ\oplus$ possède une décomposition en poids faible.

La terminologie provient de ce que $F$ est \hd, respectivement \htr, s'il envoie tout endomorphisme diagonalisable d'un $A$-module libre de rang fini sur un endomorphisme diagonalisable, respectivement trigonalisable, de $F(A^n)$, \emph{avec une condition supplémentaire} qui est essentiellement une condition de diagonalisabilité ou de trigonalisabilité \emph{en famille}, laquelle est \emph{automatique si $F$ prend des valeurs de dimensions finies}. On pourra comparer ces considérations à la décomposition à la Steinberg globale de \cite[§\,4]{DTV}, qui est obtenue à partir de la décomposition de Jordan multiplicative de l'image par un foncteur à valeurs de dimensions finies de $\F(A,K)$ des endomorphismes \emph{unipotents} d'un $A$-module libre de rang fini (voir notamment \cite[prop.~4.26]{DTV}).

La propriété \htr e d'un foncteur est liée de près au fait de prendre des valeurs de dimensions finies. On montrera en particulier que
\emph{si l'anneau $A$ est noethérien et le corps $K$ algébriquement clos, alors un foncteur de type fini de $\F(A,K)$ est \htr\ si et seulement s'il est à valeurs de dimensions finies}. Le fait qu'un foncteur à valeurs de dimensions finies de $\F(A,K)$ soit \htr\ si $K$ est algébriquement clos est essentiellement immédiat (et ne nécessite aucune hypothèse sur $A$) ; en revanche, l'énoncé réciproque (pour un foncteur de type fini et $A$ noethérien) est difficile (cf. théorème~\ref{th-htr_ln}) et apparaît comme sous-produit de théorèmes de structure dont se déduisent directement les théorèmes~\ref{thi-strtfdf} et~\ref{thi-spa}.

La structure des foncteurs \hds\ de type fini est beaucoup plus rigide que celle des foncteurs \hts\, (et liée aux foncteurs polynomiaux \emph{stricts}) ; elle est à la base du théorème~\ref{thi-Fpbar}. En effet, lorsque $K$ est un corps algébriquement clos et $L$ un corps parfait de caractéristique première, tout foncteur à valeurs de dimensions finies de $\F(L,K)$ est \hd.

\section*{Organisation du mémoire}

Le premier chapitre est consacré à des rappels algébriques et catégoriques d'usage fréquent dans la suite. Le second reste dans un cadre général, sans résultat original, mais s'approche du c\oe ur du sujet par des rappels sur la version abstraite des décompositions en poids, à savoir la notion de torsion relativement à des idéaux d'un anneau commutatif dans des catégories abéliennes linéaires sur celui-ci.

Le chapitre~\ref{spcf} introduit les décompositions en poids dans les catégories de foncteurs et en donne des propriétés générales élémentaires. Il se place dans le cadre de foncteurs sur une catégorie essentiellement petite munie d'une action d'un monoïde commutatif vers une catégorie abélienne raisonnable linéaire sur un corps commutatif. Le chapitre~\ref{spfca} se concentre sur la situation où la source est additive, linéaire sur un anneau commutatif $A$, pour l'action du monoïde multiplicatif sous-jacent à $A$. Il donne plusieurs propriétés générales des décompositions en poids spécifiques à cette situation, qui est celle à laquelle la suite du mémoire est consacrée, notamment pour des classes fondamentales de foncteurs (comme les foncteurs polynomiaux) dont il rappelle quelques propriétés utiles. Ce chapitre étudie aussi quelques aspects, qui s'apparentent à de l'arithmétique élémentaire, souvent sans intervention de foncteurs, mais qui interviennent à plusieurs reprises dans la suite, des poids dans ce contexte.

C'est au chapitre~\ref{sec-finpolp} qu'apparaissent les premiers résultats nouveaux substantiels. On y montre que l'existence d'une décomposition en poids pour un foncteur polynomial possède d'importantes conséquences en termes de structure et de finitude.

Le chapitre~\ref{chap-spec} introduit des outils fondamentaux, en complément des décompositions en poids, pour l'étude des foncteurs sur une catégorie additive, notamment la notion de \textit{spectre} d'un foncteur, qui s'avérera au c\oe ur de la démonstration du théorème de structure des foncteurs \hds. Le chapitre~\ref{spefcr} définit les notions essentielles de foncteur \hd\ et \htr. Il en dresse les propriétés générales élémentaires et les applique aux foncteurs \phs\ (généralisation des foncteurs polynomiaux introduite dans \cite{DTV}), ce qui permet d'établir un résultat dont le théorème~\ref{thi-fdflf} se déduit facilement. 

Le chapitre~\ref{shd} étudie la structure des foncteurs \hds, permettant d'obtenir comme cas particulier le théorème~\ref{thi-Fpbar}. Le chapitre~\ref{shts} constitue le c\oe ur du présent travail : les résultats de structure des foncteurs \hts\ constituent la base de nos théorèmes les plus importants, déclinés dans les chapitres suivants. On obtient ainsi rapidement, au chapitre~\ref{stfdf}, la noethérianité des foncteurs de type fini et à valeurs de dimensions finies de $\F(A,K)$, déduite d'un théorème de \guillemotleft~décomposition primaire~\guillemotright\ qui entraîne également le théorème~\ref{thi-strtfdf}. Le chapitre~\ref{stftf} s'attache à étendre une partie de ces résultats aux foncteurs de type fini et à valeurs de type fini de $\F(A,k)$. On y montre notamment les théorèmes \ref{pri-k1}, \ref{pri-k2}, \ref{pri-k3} et \ref{pri-k4}.

% TOUTES RÉFÉRENCES À \cite{DT-ext} à supprimer dès que possible !
% ICI, DANS L'INTRO, FAIRE RÉFÉRENCE AUX TROIS NOUVEAUX ARTICLES HOMOLOGIQUES.
L'appendice~\ref{app_act-grp}, qui ne prétend pas à l'originalité, présente des résultats généraux sur une construction de produit semi-direct d'une catégorie abélienne par un groupe fini qui sont utilisés aux chapitres \ref{spfca} et \ref{sec-finpolp} mais nous ont paru trop techniques pour y être exposés. L'appendice~\ref{ch-apExt} aborde, en lien avec \cite{DT-ext}, quelques aspects cohomologiques de décompositions qui apparaissent dans ce mémoire. Enfin, l'appendice~\ref{ap_afmn} complète sur plusieurs points l'étude de $\F(A,K)$ lorsque le quotient de $A$ par son nilradical est fini (situation dans laquelle les décompositions en poids apportent  très peu).

%%%%%%%%%%%
% POINT IMPORTANT À TRAITER
%Penser à revoir les citations à \cite{DT-ext} lorsque son découpage sera effectif.
%%%%%%%%%%%

\chapter{Conventions, notations et rappels généraux}\label{sct-conv}

Ce chapitre préliminaire donne des notations et des conventions qui seront utilisées couramment dans tout ce mémoire. Il rappelle également quelques résultats algébriques et catégoriques pour la plupart classiques et élémentaires (qui seront quelquefois utilisés sans référence par la suite), à l'exception du théorème (\ref{th-PSS}) de finitude récent et profond de Putman-Sam-Snowden, qui constituera l'un des outils essentiels des théorèmes de noethérianité que nous établirons en fin de mémoire.

\paragraph*{Notations ensemblistes générales}

On désigne par $\mathbb{Z}$ l'ensemble des entiers relatifs ; pour $(a,b)\in\mathbb{Z}^2$, on note $\llbracket a, b\rrbracket:=\{i\in\mathbb{Z}\,|\,a\le i\le b\}$. On désigne par $\mathbb{N}$ (resp. $\mathbb{N}^*$) l'ensemble des entiers naturels (resp. strictement positifs).

\paragraph*{Notations catégoriques générales}

Si $\C$ est une catégorie, on note $\C^\op$ la catégorie opposée à $\C$. Si $x$ et $y$ sont des objets de $\C$, on note $\C(x,y)$ ou $\mathrm{Hom}_\C(x,y)$ (voire $\mathrm{Hom}(x,y)$ s'il n'y a pas d'ambiguïté possible) l'ensemble\,\footnote{Comme c'est l'usage le plus courant, les catégories (sans plus de précision) auront pour nous une \emph{classe} d'objets, mais un \emph{ensemble} de morphismes entre deux objets fixés.} des morphismes de $x$ vers $y$.

Si $\C$ est une catégorie essentiellement petite et $\E$ une catégorie, on note $\fct(\C,\E)$\index{nota}{fct@$\fct$|textbf} la catégorie des foncteurs de $\C$ vers $\E$.

\section{Monoïdes}\label{parmo}

\textbf{Tous les monoïdes seront supposés commutatifs, sauf mention expresse du contraire} ; on les notera en général multiplicativement. On désigne par \index{nota}{Mon@$\mon$ \emph{(catégorie des monoïdes)}} $\mon$ la catégorie des monoïdes. Ainsi, pour $M$ et $N$ des objets de $\mon$, l'ensemble des morphismes $\mon(M,N)$ porte une structure naturelle de monoïde.
 
\begin{nota}\label{nota-adj0} Si $M$ est un monoïde, on note $M_0$ le monoïde obtenu en adjoignant à $M$ un élément absorbant externe.

Si $S$ est un semi-groupe, on note $S_+$ le monoïde obtenu en adjoignant à $S$ une unité externe.
\end{nota}

Si $M$ est un monoïde non nécessairement commutatif, on note $\underline{M}$\index{nota}{M@$\underline{M}$ \emph{(catégorie à un objet associée au monoïde $M$)}} la catégorie possédant exactement un objet dont le monoïde d'endomorphismes est $M$. On désigne par $Z(M)$\index{nota}{Z@$Z$ \emph{(centre)}} le centre\index{termin}{centre} de $M$. C'est un monoïde commutatif.

\subsection*{Éléments remarquables d'un monoïde}

Soient $M$ un monoïde et $x$ un élément de $M$. Une et une seule des situations suivantes advient :
\begin{enumerate}
\item la suite $(x^i)_{i\in\mathbb{N}}$ est injective (autrement dit, le sous-monoïde engendré par $x$ est libre) --- on dit alors que $x$ est \emph{sans torsion} ;
\item la suite $(x^i)_{i\in\mathbb{N}}$ est périodique à partir d'un certain rang, ce qui équivaut à dire que le sous-monoïde engendré par $x$ est fini --- on dit alors que $x$ est \emph{périodique}, ou \emph{de torsion}.
\end{enumerate}

Par ailleurs, on dit que $x$ est \emph{régulier}\index{termin}{regulier@régulier} (au sens de von Neumann) s'il existe $y\in M$ tel que $x=x^2 y$. Tout élément inversible de $M$ est régulier.

Comme $M$ est commutatif, ses éléments périodiques (resp. réguliers) en forment un sous-monoïde.

Les propriétés suivantes sont équivalentes :
\begin{enumerate}
\item il existe $i\in\mathbb{N}^*$ tel que $x^{i+1}=x$ ;
\item $x$ est périodique et régulier.
\end{enumerate}

\subsection*{Action d'un monoïde sur une catégorie}

Une \emph{action} d'un monoïde $M$ sur une catégorie $\C$ est la donnée d'une action de $M$ sur chaque ensemble de morphismes $\C(x,y)$ de sorte que les applications de composition soient équivariantes, autrement dit, c'est la donnée d'un enrichissement de $\C$ sur la catégorie monoïdale symétrique des $M$-ensembles (munis du produit cartésien).

 Se donner une action de $M$ sur $\C$ revient également à se donner un morphisme de monoïdes $M\to Z(\C)$, où $Z(\C)$\index{nota}{Z@$Z$ \emph{(centre)}} désigne le \emph{centre}\index{termin}{centre} de $\C$, c'est-à-dire le monoïde (automatiquement commutatif) des transformations naturelles de $\mathrm{Id}_\C$ vers elle-même\,\footnote{Pour que ces transformations naturelles forment un ensemble et pas seulement une classe, on peut supposer que $\C$ est essentiellement petite, ou plus généralement possède un ensemble de générateurs. La notation $Z$ est cohérente avec celle utilisée plus haut au vu de l'identification $Z(\underline{T})\simeq Z(T)$ pour un monoïde non nécessairement commutatif $T$.}. Autrement dit, on se donne pour tout objet $x$ de $\C$ un morphisme $M\to Z(\mathrm{End}_\C(x))$ de sorte que pour tous objets $x$ et $y$ les actions de $M$ sur $\C(x,y)$ données par la composition à la source ou au but coïncident.

\subsection*{Localisation}

Si $M$ est un monoïde et $S$ un sous-monoïde de $M$, on note $M[S^{-1}]$\index{nota}{$\mathfrak{X}[S^{-1}]$ \emph{(localisation d'un monoïde $\mathfrak{X}$, d'un ensemble $\mathfrak{X}$ ou d'une catégorie $\mathfrak{X}$ muni(e) d'une action d'un monoïde par rapport à un sous-monoïde $S$ de celui-ci)}} le monoïde obtenu à partir de $M$ en inversant formellement les éléments de $S$. Si $E$ est un $M$-ensemble à droite, on note $M[S^{-1}]$ l'ensemble quotient $(E\times S)/((x,s)\sim (y,t)\Leftrightarrow\exists u\in S\quad x.(tu)=y.(su))$ : l'action de $M$ sur $E$ induit une action de $M[S^{-1}]$ sur $E[S^{-1}]$ de sorte que l'action de $mt^{-1}$ sur la classe de $(x,s)\in E\times S$ soit la classe de $(x.m,st)$. Ainsi, si l'on note de la même façon un élément de $E$ et son image par l'application canonique (qui est équivariante relativement au morphisme de monoïdes canonique $M\to M[S^{-1}]$, mais n'est pas nécessairement injective) $E\to E[S^{-1}]$ associant à $x$ la classe de $(x,1)$, la classe de $(x,s)$ est $x.s^{-1}$. 

La notation $M[S^{-1}]$, où $M$ est muni de l'action à droite de $M$ par translations, est cohérente avec la précédente.

\smallskip

Supposons maintenant que $M$ agit sur  une catégorie $\C$. On peut définir une catégorie $\C[S^{-1}]$ ayant les mêmes objets que $\C$ et dont les morphismes sont donnés par $\C[S^{-1}](x,y):=\C(x,y)[S^{-1}]$, la composition étant donnée par $(fs^{-1})\circ (gt^{-1}):=(f\circ g)(st)^{-1}$, où $f$ et $g$ sont des flèches composables de $\C$ et $s, t$ des éléments de $S$. La catégorie $\C[S^{-1}]$ est munie d'une action de $M[S^{-1}]$ induite par celle de $M$ sur $\C$, et l'on dispose d'un foncteur canonique $\C\to\C[S^{-1}]$ égal à l'identité sur les objets.

On note que $\underline{M}[S^{-1}]$ s'identifie canoniquement à $\underline{M[S^{-1}]}$.\index{nota}{$\mathfrak{X}[S^{-1}]$ \emph{(localisation d'un monoïde $\mathfrak{X}$, d'un ensemble $\mathfrak{X}$ ou d'une catégorie $\mathfrak{X}$ muni(e) d'une action d'un monoïde par rapport à un sous-monoïde $S$ de celui-ci)}}

Lorsque la catégorie $\C$ possède des colimites filtrantes, le foncteur canonique $\C\to\C[S^{-1}]$ possède un adjoint à gauche qui identifie $\C[S^{-1}]$ à la sous-catégorie pleine de $\C$ des objets sur lesquels les éléments de $S$ agissent par des isomorphismes.

Si l'on suppose $\C$ essentiellement petite, l'action de $M$ sur $\C$ induit, par précomposition, une action de $M$ sur la catégorie de foncteurs $\fct(\C,\E)$ (où est $\E$ une catégorie quelconque).

Le résultat suivant est immédiat mais utile :
\begin{lemm}\label{lm-inversion_fct} Soient $\C$ une catégorie essentiellement petite munie d'une action d'un monoïde $M$, $S$ un sous-monoïde de $M$ et $\E$ une catégorie.

Le foncteur $\fct(\C[S^{-1}],\E)\to\fct(\C,\E)$ de précomposition par le foncteur canonique $\C\to\C[S^{-1}]$ est pleinement fidèle, et son image essentielle est constituée des foncteurs $F$ tels que $F(s)$ soit inversible pour tout $s\in S$. Si $\E$ est cocomplète, alors $\fct(\C[S^{-1}],\E)$ s'identifie à $\fct(\C,\E)[S^{-1}]$.
\end{lemm}

L'action d'un monoïde se propage également au \emph{but} des catégories de foncteurs : toute action d'un monoïde $N$ sur une catégorie $\E$ induit une action sur $\fct(\C,\E)$, si $\C$ est une catégorie essentiellement petite (et, lorsque $\C$ est munie d'une action de $M$, ladite action commute avec l'action de $M$ à la source ; autrement dit, la catégorie $\fct(\C,\E)$ est munie d'une action du monoïde $M\times N$). Si $T$ est un sous-monoïde de $N$, on dispose d'une équivalence canonique $\fct(\C,\E)[T^{-1}]\simeq\fct(\C,\E[T^{-1}])$ lorsque la catégorie $\E$ possède des colimites filtrantes.

\section{Anneaux}

\begin{conv} Sauf mention explicite du contraire, \textbf{tous les anneaux et en particulier tous les corps sont supposés commutatifs}. Les algèbres sont toujours supposées unitaires et associatives, et commutatives sauf mention du contraire.
\end{conv}

Si $k$ est un anneau, non nécessairement commutatif, et $E$ un ensemble, on note $k[E]:=k^{\oplus E}$ le $k$-module à gauche libre sur $E$.

\begin{conv} Dans tout ce mémoire, \textbf{$A$ désigne un  anneau et $K$ un corps}. On note $p$ la caractéristique de $K$. Sauf mention du contraire, les produits tensoriels de base non spécifiée sont pris sur $K$.

On note $A_K$\index{nota}{AK@$A_K$} la $K$-algèbre $A\otimes_\mathbb{Z}K$.

On désigne par $A^\times$ le groupe des unités de $A$.
\end{conv}

On note \index{nota}{Ann@$\mathbf{Ann}$ \emph{(catégorie des anneaux)}} $\mathbf{Ann}$ la catégorie des anneaux.

On désigne comme d'habitude par $\spec(A)$\index{nota}{spec@$\spec$ \emph{(spectre premier d'un anneau)}} l'ensemble des idéaux premiers de $A$.

Un idéal $I$ de $A$ est dit \emph{cofini}\index{termin}{cofini!idéal} si l'anneau quotient $A/I$ est fini.

Le centre\index{termin}{centre} d'un anneau non nécessairement commutatif $R$ est noté $Z(R)$.\index{nota}{Z@$Z$ \emph{(centre)}} 

Une $K$-algèbre $A$ est dite \index{termin}{deploye@déployé!\emph{($K$-algèbre déployée)}} {\em déployée} si, pour tout idéal maximal $\mathfrak{m}$ de $A$, le morphisme de corps $K\to A/\mathfrak{m}$ est un isomorphisme.

On note \index{nota}{nu@$\nu(K)$}\label{pnu} $\nu(K)$ l'ensemble des entiers $n>0$ tels que le polynôme $X^n-1$ de $K[X]$ soit scindé, i.e. tels que l'algèbre $K[X]/(X^n-1)$ soit déployée.

On note $\mathbf{Ab}$ \index{nota}{Ab@$\mathbf{Ab}$} la catégorie des groupes abéliens. Si $k$ est un anneau non nécessairement commutatif, on désigne par \index{nota}{Mod@$k\Md$} $k\Md$ la catégorie des $k$-modules à gauche, et par $\mathbf{P}(k)$ \index{nota}{P@$\mathbf{P}(k)$} la sous-catégorie pleine dont les objets sont les modules libres $k^i$ pour $i\in\mathbb{N}$. On remarquera que dans de nombreuses références (notamment \cite{DTV}, que nous utiliserons abondamment), $\mathbf{P}(k)$ désigne plutôt (un squelette de) la catégorie des modules {\em projectifs} de type fini de $k\Md$ (d'où la notation), mais c'est sans conséquence sur les catégories de foncteurs que nous considérerons, car l'injection des modules libres de type fini dans les modules projectifs de type fini est une équivalence de Morita : elle induit par précomposition une équivalence entre toutes les catégories de foncteurs de but abélien (ou plus généralement où les idempotents se scindent).

\subsection*{Radical}

On note \index{nota}{J@$\rj$ \emph{(radical de Jacobson d'un anneau)}} $\rj(A)$ le \index{termin}{radical!d'un anneau} radical de Jacobson de $A$, c'est-à-dire l'intersection de ses idéaux maximaux, et \index{nota}{nil@$\nil$ \emph{(nilradical d'un anneau)}} $\nil(A)$ son nilradical, c'est-à-dire l'intersection de ses idéaux premiers, ou encore, l'ensemble de ses éléments nilpotents \cite[chap.~II, §\,2, n°\,6, prop.~13]{Bki-AC1}. On dit que $A$ est \index{termin}{reduit@réduit!\emph{(anneau)}} {\em réduit} si $\nil(A)=\{0\}$. 

On rappelle que $A$ est dit \index{termin}{semi-local}{\em semi-local} si l'anneau quotient $A/\rj(A)$ est semi-simple, ce qui équivaut à dire que $A$ n'a qu'un nombre fini d'idéaux maximaux, ou que $A/\rj(A)$ est isomorphe à produit fini de corps \cite[chap.~II, §\,3, n°\,5, prop.~16]{Bki-AC1}.

L'anneau $A$ est dit \index{termin}{semi-parfait} {\em semi-parfait} s'il est semi-local et que les idempotents de $A/\rj(A)$ se relèvent en des idempotents de $A$. Il revient au même de demander que $A$ soit isomorphe à un produit fini d'anneaux locaux.

Contrairement aux autres termes employés dans ce paragraphe pour désigner certains types d'anneaux, le qualificatif d'anneau \emph{quasi-parfait} introduit ci-dessous ne semble pas répertorié, même si la notion est bien connue.

\begin{defi}\label{df-qpft} L'anneau $A$ est dit \index{termin}{quasi-parfait} {\em quasi-parfait} s'il est semi-local et que $\nil(A)=\rj(A)$.
\end{defi}

\begin{prop}\label{pr-qpsp} Un anneau quasi-parfait est semi-parfait.
\end{prop}

\begin{proof} Cela résulte de \cite[§\,9, n°\,4, prop.~7]{Bki}.
\end{proof}

L'énoncé classique et facile suivant, dont nous omettons la démonstration, dresse une liste de relations d'usage courant entre ces différentes propriétés d'un anneau.

\begin{prop}\label{pr-qpsi}\begin{enumerate}
\item Les assertions suivantes sont équivalentes :
\begin{enumerate}
\item l'anneau $A$ est local et quasi-parfait ;
\item l'ensemble $A$ est la réunion disjointe de $\nil(A)$ et de $A^\times$ ;
\item l'anneau $A$ possède exactement un idéal premier.
\end{enumerate}
\item Les assertions suivantes sont équivalentes :
\begin{enumerate}
\item\label{itqpf1} l'anneau $A$ est quasi-parfait ;
\item $A$ est isomorphe à un produit fini d'anneaux locaux quasi-parfaits ;
\item\label{itqpf3} $A$ est semi-local et tout idéal premier de $A$ est maximal.
\end{enumerate}
\item Le spectre d'un anneau quasi-parfait est fini.
\item\label{it-qpfna} Un anneau est artinien si et seulement s'il est quasi-parfait et noethérien.
\end{enumerate}
\end{prop}

L'anneau $A$ est dit
\index{termin}{parfait \emph{(anneau)}} {\em parfait}\,\footnote{Il n'y a pas en pratique de confusion à craindre avec la notion différente de {\em corps parfait}. De toute manière, nous n'utiliserons pratiquement pas la notion d'anneau parfait dans ce mémoire ; elle est mentionnée ici pour la cohérence de cet exposé de rappels sur les anneaux.} s'il est semi-local et que, pour toute suite dénombrable $(x_i)_{i\in\mathbb{N}}$ de $\rj(A)$, il existe $n\in\mathbb{N}$ tel que $\prod_{i\le n}x_i=0$. Un tel anneau est quasi-parfait. La réciproque est toutefois erronée, comme le montre l'exemple suivant.

\begin{exem}\label{ex-qpft} La $\FF_2$-algèbre $\FF_2[x_0,x_1,\dots,x_n,\dots]/(x_0^2,x_1^2,\dots,x_n^2,\dots)$ est locale, et tout élément de son radical est de carré nul. Elle est en particulier quasi-parfaite. Mais la considération de la suite $(x_i)_{i\in\mathbb{N}}$ d'éléments de son radical montre qu'elle n'est pas parfaite.
\end{exem}

L'anneau $A$ est dit {\em semi-primaire}\index{termin}{semi-primaire} s'il est semi-local et que l'idéal $\rj(A)$ est nilpotent. Un tel anneau est parfait. Tout anneau artinien est semi-primaire.

Si $I$ est un idéal \emph{nilpotent} de $A$, alors $A$ est semi-primaire si et seulement si $A/I$ est semi-primaire.

\begin{lemm}\label{lm-sqalgdepl}  Tout produit fini et tout produit tensoriel fini de $K$-algèbres semi-simples (resp. semi-primaires) déployées est semi-simple (resp. semi-primaire) déployé.

Toute sous-algèbre et toute algèbre quotient d'une $K$-algèbre semi-simple (resp. semi-primaire) déployée est semi-simple (resp. semi-primaire) déployée.
\end{lemm}

\begin{proof}%e omettre la dém. ? Faire la remarque que l'assertion relative aux sous-algèbres devient fausse si on supprimer "déployée" ?
L'assertion relative aux produits finis et aux produits tensoriels finis est immédiate.  Celle relative aux sous-algèbres et quotients des algèbres semi-simples déployées est contenue dans \cite[chap.~V, §\,6, prop.~3]{Bki2} (où les algèbres semi-simples déployées sont appelées \emph{algèbres diagonalisables}).

Si $B$ est une sous-algèbre d'une $K$-algèbre semi-primaire déployée $A$, alors l'image du morphisme d'algèbres canonique $B\hookrightarrow A\twoheadrightarrow A/\rj(A)$ est semi-simple déployée en tant que sous-algèbre de l'algèbre semi-simple déployée $A/\rj(A)$. Son noyau est inclus dans $\rj(A)$ ; c'est donc un idéal nilpotent de $B$. Il s'ensuit que cette algèbre est semi-primaire déployée. L'assertion relative aux quotients s'établit de façon analogue.
\end{proof}

En combinant l'assertion du lemme précédente relative aux produits tensoriels finis et celle sur les quotients, on obtient :

\begin{lemm}\label{lm-sqalgdepl2} Toute $K$-algèbre engendrée par un nombre fini de sous-algèbres semi-simples (resp. semi-primaires) déployées est semi-simple (resp. semi-primaire) déployée.
\end{lemm}

\begin{exem}\label{ex-assez_rac1} Si $M$ est un monoïde fini, alors la $K$-algèbre $K[M]$ est de dimension finie, donc semi-primaire. Elle est déployée si et seulement si pour tout $x\in M$, l'ordre $a$ du sous-\emph{semi-groupe} $(x^n)_{n\in\mathbb{N}^*}$ de $M$ engendré par $x$ appartient à $\nu(K)$, et semi-simple si et seulement si $a$ est inversible dans $K$ et $x$ régulier. En effet, les lemmes \ref{lm-sqalgdepl} et \ref{lm-sqalgdepl2} montrent qu'il suffit de traiter le cas où $M$ est engendré par $x$, qui résulte de ce qu'on a alors $K[M]\simeq K[X]/(X^n.(X^a-1))\simeq K[X]/(X^n)\times K[X]/(X^a-1)$ pour un certain $n\in\mathbb{N}$, et $x$ est régulier si et seulement si $n\in\{0,1\}$.
\end{exem}

\begin{lemm}\label{lm-qsp} Supposons que $A$ est un anneau semi-parfait.
\begin{enumerate}
\item\label{itqsp1} Les assertions suivantes sont équivalentes :
\begin{enumerate}
\item tout quotient semi-primaire de $A$ est semi-simple ;
\item pour tout idéal maximal $\mathfrak{m}$ de $A$, $\mathfrak{m}^2=\mathfrak{m}$ ;
\item $\rj(A)^2=\rj(A)$.
\end{enumerate}
\item\label{itqsp2} Les assertions suivantes sont équivalentes :
\begin{enumerate}
\item tout quotient semi-primaire de $A$ est artinien ;
\item pour tout idéal maximal $\mathfrak{m}$ de $A$, $\dim_{A/\mathfrak{m}}\mathfrak{m}/\mathfrak{m}^2<\infty$ ;
\item $\rj(A)/\rj(A)^2$ est un $A/\rj(A)$-module de type fini.
\end{enumerate}
\end{enumerate}
\end{lemm}

\begin{proof} Comme $A$ est un produit fini d'anneaux locaux, il suffit de démontrer le résultat lorsque $A$ est local ; notons $\mathfrak{m}$ son radical.

Un quotient non nul de $A$ est de la forme $A/I$, où $I\subset\mathfrak{m}$ est un idéal, et a pour radical $\mathfrak{m}/I$. Ainsi, $A/I$ est semi-primaire si et seulement s'il existe $n\in\mathbb{N}^*$ tel que $\mathfrak{m}^n\subset I\subset\mathfrak{m}$.

Comme la condition $\mathfrak{m}^2=\mathfrak{m}$ implique $\mathfrak{m}^n=\mathfrak{m}$ pour tout $n\in\mathbb{N}^*$, l'assertion~\ref{itqsp1} s'ensuit.

Comme on dispose d'un morphisme surjectif $A/\mathfrak{m}$-linéaire $S^n_{A/\mathfrak{m}}(\mathfrak{m}/\mathfrak{m}^2)\twoheadrightarrow\mathfrak{m}^n/\mathfrak{m}^{n+1}$ (où $S^n$ désigne la $n$-ième puissance symétrique), l'assertion~\ref{itqsp2} découle également de l'observation précédente.
\end{proof}

\subsection*{Anneaux et monoïdes}

 On note $A_\mu$ \index{nota}{A@$A_\mu$ \emph{(monoïde multiplicatif d'un anneau $A$)}} le monoïde multiplicatif sous-jacent à un anneau $A$. Le foncteur associant à un monoïde $M$ l'anneau (resp. la $K$-algèbre) $\mathbb{Z}[M]$ (resp. $K[M]$) de $M$ est adjoint à gauche au foncteur $\mathbf{Ann}\to\mathbf{Mon}\quad A\mapsto A_\mu$ (resp. à la composée de ce foncteur avec le foncteur d'oubli des $K$-algèbres dans $\mathbf{Ann}$).
 
On désigne par $\bar{K}[M]$\index{nota}{K@$\bar{K}[\cdot]$ \emph{(idéal d'augmentation)}} l'idéal d'augmentation de la $K$-algèbre d'un monoïde $M$.
  
Les exemples ci-dessous illustrent certaines notions introduites au §\,\ref{nota-adj0} dans le monoïde $A_\mu$. 
 
\begin{exem}\label{ex-aa0} Si $A$ est un anneau intègre, le monoïde multiplicatif $A_\mu$ s'identifie à $(A\setminus\{0\})_0$ (cf. notation~\ref{nota-adj0}).
\end{exem}

\begin{exem}\label{exev-mormon} Si $A$ est un anneau local quasi-parfait, on dispose d'un isomorphisme canonique de monoïdes $\mon(A_\mu,K_\mu)\simeq\mathbf{Ab}(A^\times,K^\times)_+$.

En effet, si $\varphi : A_\mu\to K_\mu$ est un morphisme de monoïdes, comme $0\in A$ est idempotent, soit $\varphi$ envoie $0$ sur $1$, auquel cas $\varphi$ est le morphisme trivial (puisque $0$ est absorbant), soit $\varphi$ envoie $0$ sur $0$. Dans ce dernier cas, les éléments non inversibles de $A$, qui sont nilpotents (proposition~\ref{pr-qpsi}), sont nécessairement envoyés sur $0$ par $\varphi$.
\end{exem}

Un élément de $A$ est dit \emph{régulier}\index{termin}{regulier@régulier} (au sens de von Neumann) si c'est un élément régulier du monoïde $A_\mu$. On dit que $A$ est un anneau \index{termin}{absolument plat|textbf} {\em absolument plat} si tous ses éléments sont réguliers. Il revient au même de demander que tout idéal principal de $A$ soit engendré par un idempotent, ou que tout $A$-module de présentation finie soit projectif, ou encore que tout $A$-module soit plat \cite[Théorème~8.6]{Pop}.

La propriété suivante est évidente.

\begin{prop} Tout produit et toute colimite filtrante d'anneaux absolument plats est un anneau absolument plat.
\end{prop}

\begin{coro}\label{cor-vN} Si $G$ est un groupe abélien dont tous les éléments sont d'ordres finis et inversibles dans $K$, alors l'algèbre de groupe $K[G]$ est un anneau absolument plat.
\end{coro}

La notion suivante sera utilisée à deux reprises au chapitre~\ref{shts}.

\begin{nota}\label{notapereg} Le sous-monoïde de $A_\mu$ constitué des éléments périodiques et réguliers de $A$ est noté $A^\mathrm{pereg}$.\index{nota}{Apereg@$A^\mathrm{pereg}$|textbf}
\end{nota}

Nous donnons maintenant deux énoncés classiques reliant les propriétés de finitude des monoïdes et des anneaux. Le premier est établi par exemple dans \cite[Théorème.~7.7]{Gilm}.

\begin{prop}\label{pr-montfn} Soit $M$ un monoïde. La $K$-algèbre $K[M]$ est noethérienne si et seulement si $M$ est un monoïde de type fini.
\end{prop}

Le résultat suivant fut démontré pour la première fois dans \cite{Isb59}.

\begin{theo}[Isbell]\label{th-Isb} Le monoïde $A_\mu$ est de type fini si et seulement si $A$ est un anneau fini.
\end{theo}

Nous utiliserons aussi quelquefois le groupe additif sous-jacent à un anneau $A$, noté $A_\mathrm{add}$. Le $K$-espace vectoriel $K[A]$ hérite ainsi de deux structures de $K$-algèbre : celle de $K[A_\mathrm{add}]$ et celle de $K[A_\mu]$. Le produit donné par cette dernière structure (que nous utiliserons de loin le plus souvent) sera noté par un point (voire l'absence de symbole) : $[a].[b]=[ab]$, tandis que l'autre produit sera noté $\star$\index{nota}{$\star$} : $[a]\star [b]=[a+b]$. Si $U$ et $V$ sont des sous-espaces vectoriels de $K[A]$, nous noterons $U\star V$ le sous-espace vectoriel de $K[A]$ engendré par les éléments $x\star y$ pour $(x,y)\in U\times V$, et nous appellerons $\star$-idéal de $K[A]$ un idéal de $K[A_\mathrm{add}]$. Un sous-ensemble de $K[A]$ qui est un idéal à la fois de $K[A_\mathrm{add}]$ et de $K[A_\mu]$ sera appelé {\em bi-idéal}\index{termin}{biideal@bi-idéal|textbf}\label{pbiid} de $K[A]$. De façon équivalente, si l'on munit l'ensemble $A$ de l'action à gauche canonique du monoïde non commutatif affine $\operatorname{aff}(A):=\{A\xrightarrow{x\mapsto ax+b}A\,|\,(a,b)\in A^2\}\simeq A_\mu\ltimes A_{\mathrm{add}}$, un bi-idéal de $K[A]$ en est un sous-$K[\operatorname{aff}(A)]$-module à gauche.

\begin{rema}\label{rqstev}
\begin{enumerate}
\item Si $V$ est un sous-groupe de $A_\mathrm{add}$, alors le noyau de l'application linéaire $K[A]\twoheadrightarrow K[A/V]$ induite par la projection $A\twoheadrightarrow A/V$ est $K[A]\star\bar{K}[V]$.
\item Si $U$ et $V$ sont des idéaux de $K[A_\mu]$, alors $U\star V$ est aussi un idéal de $K[A_\mu]$. Si de plus $U$ ou $V$ est un bi-idéal de $K[A]$, il s'ensuit que $U\star V$ est un bi-idéal.
\item En particulier, si $I$ est un idéal de $K[A_\mu]$, alors $K[A]\star I$ (c'est-à-dire le $\star$-idéal de $K[A]$ engendré par $I$) est un bi-idéal de $K[A]$, c'est donc le bi-idéal engendré par $I$.
\item L'idéal d'augmentation $\bar{K}[A]$ est un bi-idéal, ses puissances $\bar{K}[A]^{\star d}$ sont donc des bi-idéaux.
\end{enumerate}
\end{rema}

\section{Dérivations}

Si $\varphi : A\to K$ est un morphisme d'anneaux, une \emph{$\varphi$-dérivation}\index{termin}{derivation@dérivation|textbf} de $A$ dans $K$ est un morphisme de groupes abéliens $\mathrm{d} : A_{\mathrm{add}}\to K_{\mathrm{add}}$ tel que
$$\forall (x,y)\in A^2\qquad\mathrm{d}(xy)=\varphi(x)\mathrm{d}(y)+\varphi(y)\mathrm{d}(x).$$
Ces $\varphi$-dérivations forment un sous-espace vectoriel de $\mathbf{Ab}(A,K)$ qu'on notera 
 $\mathrm{Der}_\varphi(A,K)$.\index{nota}{Der@$\mathrm{Der}_\varphi(A,K)$} Nous noterons également $\mathfrak{m}_\varphi$\index{nota}{m@$\mathfrak{m}_\varphi$ \emph{(pour $\varphi\in\mathbf{Ann}(A,K)$)}} l'idéal de la $K$-algèbre $A_K$ noyau du morphisme d'algèbres $A_K\to K$ envoyant $a\otimes\lambda$ sur $\varphi(a)\lambda$. On dispose d'isomorphismes canoniques de $K$-espaces vectoriels
\begin{equation}\label{eq-ident_der}
\mathrm{Der}_\varphi(A,K)\simeq (\mathfrak{m}_\varphi/\mathfrak{m}_\varphi^2)^*\,,
\end{equation}
où l'étoile désigne la dualité, et
\begin{equation}\label{eq-der_Ext}\mathrm{Der}_\varphi(A,K)\simeq\mathrm{Ext}^1_{A_K}(K_\varphi,K_\varphi)\,,
\end{equation}
où $K_\varphi$\index{nota}{K@$K_\varphi$} désigne le $A_K$-module $K$ muni de l'action tautologique de $K$ et de l'action de $A$ via $\varphi$.

Les \emph{conditions de finitude sur les dérivations} suivantes interviendront à quelques reprises dans ce travail :
\begin{equation*}\label{eqcfd}
\text{(CFD)}\qquad\qquad\qquad\forall\varphi\in\mathbf{Ann}(A,K)\quad\dim_K\mathrm{Der}_\varphi(A,K)<\infty\;.
\end{equation*}\index{nota}{CFD@(CFD), (CFD)$^+$ \emph{(conditions de finitude sur les dérivations)}|textbf}
(S'il y a une ambiguïté possible sur $A$ ou $K$, nous dirons que $(A,K)$ vérifie (CFD).)
\begin{equation*}
\text{(CFD)}^+\quad\text{Pour toute extension finie }L\text{ du corps }K, (A,L) \text{ vérifie la condition (CFD)}.
\end{equation*}

On a manifestement les implications suivantes :

\noindent
($A_K$ est un anneau noethérien) $\Rightarrow$ (pour tout $\varphi\in\mathbf{Ann}(A,K)$, $\mathfrak{m}_\varphi$ est un idéal de type fini de $A_K$) $\Rightarrow$ (CFD).

La condition de noethérianité de $A_K$ est en particulier vérifiée dans chacun des cas suivants :
\begin{enumerate}
\item $A$ est un anneau essentiellement de type fini\index{termin}{essentiellement de type fini \emph{(anneau, algèbre)}|textbf} (i.e. un localisé d'un anneau de type fini) ;
\item $A$ est un anneau noethérien et $K$ un corps de type fini.
\end{enumerate}

De plus, la noethérianité de $A_K$ implique celle de $A_L$ pour toute extension finie $L$ du corps $K$, et donc (CFD)$^+$.

\begin{rema}
la condition (CFD)$^+$ est également vérifiée lorsque $p>0$ et que la $\FF_p$-algèbre $A/p$ est un module (ou une algèbre) de type fini sur l'image de son endomorphisme de Frobenius, car une dérivation en caractéristique $p$ est nulle sur les puissances $p$-ièmes. Cette condition est en particulier vérifiée si $A$ est une algèbre de type fini sur un corps parfait de caractéristique $p$. Cela fournit des exemples de situations où (CFD)$^+$ est satisfaite sans que l'anneau $A_K$ soit noethérien (par exemple, si $A=K$ est un corps parfait infini de caractéristique non nulle).
\end{rema}

\begin{rema}\label{rq-CFD_imagefinie}
Supposons que l'image de $\varphi\in\mathbf{Ann}(A,K)$ est \emph{finie}, disons de cardinal $p^r$. Alors la théorie de Galois des corps finis fournit un diagramme commutatif de morphismes d'anneaux
$$\xymatrix{A_K=A\otimes_\mathbb{Z}K\ar[r]^-{\varphi\otimes K}\ar@{->>}[d] & K\otimes_\mathbb{Z}K \\
\FF_{p^r}\otimes_{\FF_p}K\ar[r]^-\simeq\ar@{^{(}->}[ur] & K^r
}$$
dont les composantes $A_K\to K$ de la composée $A_K\to K^r$ sont les $\varphi^{p^i}$ pour $0\le i<r$. Il s'ensuit que, si $I_\varphi$ désigne l'idéal de $A$ noyau de $\varphi$, on a
$$(I_\varphi/I_\varphi^2)\otimes_{\FF_p}K\simeq\bigoplus_{i=0}^{r-1}\mathfrak{m}_{\varphi^{p^i}}/\mathfrak{m}_{\varphi^{p^i}}^2\;.$$
Ainsi, la condition (CFD) implique que $I_\varphi/I_\varphi^2$ est fini, et si tout morphisme d'anneaux $A\to K$ est d'image finie (par exemple, si $K$ est fini), alors la condition (CFD) \emph{équivaut} à la finitude de $I_\varphi/I_\varphi^2$ pour tout $\varphi\in\mathbf{Ann}(A,K)$.
\end{rema}

Les \emph{conditions d'annulation des dérivations} suivantes interviendront également, notamment au chapitre~\ref{sec-finpolp} :
\begin{equation*}\label{eqcad}
\text{(CAD)}\qquad\qquad\qquad\qquad\forall\varphi\in\mathbf{Ann}(A,K)\qquad\mathrm{Der}_\varphi(A,K)=0\;.
\end{equation*}
(S'il y a une ambiguïté possible sur $A$ ou $K$, nous dirons que $(A,K)$ vérifie (CAD).)
\begin{equation*}
\text{(CAD)}^+\quad\text{Pour toute extension finie }L\text{ du corps }K, (A,L) \text{ vérifie la condition (CAD)}.
\end{equation*}
\index{nota}{CAD@(CAD), (CAD)$^+$ \emph{(conditions d'annulation des dérivations)}|textbf}

Ces conditions sont faciles à caractériser lorsque $A=K$ :
\begin{prdef}\label{prdf-sander} Les conditions suivantes sont équivalentes :
\begin{enumerate}
\item\label{itcad1} pour toute extension $L$ du corps $K$, $(K,L)$ vérifie la condition (CAD) ;
\item\label{itcad2} $(K,K)$ vérifie la condition \textnormal{(CAD)}$^+$ ;
\item\label{itcad3} $(K,K)$ vérifie la condition \textnormal{(CAD)} ;
\item\label{itcad4} $K$ est une extension algébrique de $\mathbb{Q}$ ou un corps parfait de caractéristique première.
\end{enumerate}

Lorsqu'elles sont vérifiées, nous dirons que $K$ est \textbf{sans dérivation}.
\end{prdef}

\begin{proof} Les implications \ref{itcad1}$\Rightarrow$\ref{itcad2}$\Rightarrow$\ref{itcad3} sont triviales.

Les implications \ref{itcad3}$\Rightarrow$\ref{itcad4}$\Rightarrow$\ref{itcad1} résultent de \cite[chap.~V, §\,16, th.~1 et son cor.~3]{Bki2} en caractéristique nulle et de \cite[chap.~V, §\,13, prop.~6]{Bki2} en caractéristique première.
\end{proof}

Nous aurons besoin au chapitre~\ref{sec-finpolp} du résultat technique simple suivant.
\begin{lemm}\label{lm-cafd_insep} Supposons que pour toute extension finie \textbf{séparable} $L$ du corps $K$ et tout $\alpha\in\mathbf{Ann}(A,L)$, on a $\dim_L\mathrm{Der}_\alpha(A,L)<\infty$ (resp. $\mathrm{Der}_\alpha(A,L)=0$). Alors la condition \textnormal{(CFD)}$^+$ (resp. \textnormal{(CAD)}$^+$) est vérifiée.
\end{lemm}

\begin{proof}
Il s'agit de montrer que si $L'$ est une extension finie, non nécessairement séparable, de $K$, et $\alpha : A\to L$ un morphisme d'anneaux, on a  $\dim_L\mathrm{Der}_\alpha(A,L')<\infty$ (resp. $\mathrm{Der}_\alpha(A,L')=0$). Soit $L$ la fermeture algébrique séparable \cite[chap.~5, §\,7, n°7]{Bki2} de $K$ dans $L'$. Alors $L'$ est une extension radicielle finie du corps $L$, donc il existe $n\in\mathbb{N}$ tel que la $n$-ième itérée $\varphi$ du morphisme de Frobenius de $L'$ soit à valeurs dans $L$. Ainsi, il existe $\tilde{\alpha}\in\mathbf{Ann}(A,L)$ tel que la composée de $\tilde{\alpha}$ et de l'inclusion $L\subset L'$ coïncide avec $\varphi\circ\alpha$. Ainsi, $\mathrm{Der}_{\varphi\circ\alpha}(A,L')\simeq\mathrm{Der}_{\tilde{\alpha}}(A,L)\underset{L}{\otimes} L'$ est de dimension finie (resp. nul) par hypothèse, puisque $L$ est une extension finie séparable de $K$. Comme la composition par $\varphi$ définit une application $K$-linéaire injective $\mathrm{Der}_\alpha(A,L')\to\mathrm{Der}_{\varphi\circ\alpha}(A,L')$, la conclusion s'ensuit.
\end{proof}

Nous utiliserons au chapitre~\ref{shts} la caractérisation de la condition (CAD) pour des algèbres sur un corps fini donnée par les deux énoncés suivants.
%manifestement classiques : réf. ???

\begin{lemm}\label{lm-deriv_loc} Supposons que $p$ est premier et que $A$ est une $\FF_p$-algèbre locale dont le radical est de carré nul et le corps résiduel $\kappa$ est un sous-corps de $K$. Soit $\alpha$ le morphisme d'anneaux $A\twoheadrightarrow\kappa\hookrightarrow K$. Alors $\mathrm{Der}_\alpha(A,K)$ est nul si et seulement si $A$ est un corps parfait.
\end{lemm}

\begin{proof}
Toute dérivation de $\kappa$ dans $K$ (relativement à l'inclusion $\kappa\hookrightarrow K$) se prolonge, par précomposition par $A\twoheadrightarrow\kappa$, en une $\alpha$-dérivation de $A$ dans $K$, il est donc nécessaire que $\kappa$ soit parfait pour que $\mathrm{Der}_\alpha(A,K)$ soit nul (cf. proposition~\ref{prdf-sander}).

Supposons désormais $\kappa$  parfait : comme $\rj(A)$ est de carré nul, cela implique que l'image de l'endomorphisme de Frobenius de la $\FF_p$-algèbre $A$ induit un isomorphisme de $\kappa$ sur une sous-algèbre de $A$ supplémentaire (comme $\FF_p$-espace vectoriel) de $\rj(A)$, d'où l'on déduit $\mathrm{Der}_\alpha(A,K)\simeq\mathrm{Hom}_{\FF_p}(\rj(A),\FF_p)$ et la conclusion.
\end{proof}

\begin{prop}\label{pr-CAD_fpalg} Supposons que $p$ est premier et que $A$ est une $\FF_p$-algèbre. Alors la condition \textnormal{(CAD)} est vérifiée si et seulement si, pour tout $\mathfrak{p}\in\spec(A)$ tel que le corps $\kappa_\mathfrak{p}:=\operatorname{Frac}(A/\mathfrak{p})$ se plonge dans $K$, l'idéal maximal $\mathfrak{m}_\mathfrak{p}:=\mathfrak{p}.A_\mathfrak{p}$ de l'anneau localisé $A_\mathfrak{p}$ est égal à son carré et le corps $\kappa_\mathfrak{p}$ est parfait.
\end{prop}

\begin{proof}
Soient $\varphi : A\to K$ un morphisme d'anneaux et $\mathfrak{p}$ son noyau : $\mathfrak{p}$ est un idéal premier de $A$ et le corps $\kappa_\mathfrak{p}$ se plonge dans $K$. Notant $\tilde{\varphi} : A_\mathfrak{p}\to K$ l'unique morphisme d'anneaux factorisant $\varphi$ à travers le morphisme canonique $A\to A_\mathfrak{p}$, toute dérivation $\mathrm{d}\in\mathrm{Der}_\varphi(A,K)$ se prolonge de manière unique en une dérivation $\tilde{\mathrm{d}}\in\mathrm{Der}_{\tilde{\varphi}}(A_\mathfrak{p},K)$, donnée par $\tilde{\mathrm{d}}(a.x^{-1})=\alpha(x)^{-1}\mathrm{d}(a)-\alpha(a)\alpha(x)^{-2}\mathrm{d}(x)$ pour $a\in A$ et $x\in A\setminus\mathfrak{p}$, et $\mathrm{d}\mapsto\tilde{\mathrm{d}}$ définit un isomorphisme $\mathrm{Der}_\varphi(A,K)\xrightarrow{\simeq}\mathrm{Der}_{\tilde{\varphi}}(A_\mathfrak{p},K)$. Par ailleurs, comme $\tilde{\varphi}$ est nul sur $\mathfrak{m}_\mathfrak{p}$, tout élément de $\mathrm{Der}_{\tilde{\varphi}}(A_\mathfrak{p},K)$ est nul sur $\mathfrak{m}_\mathfrak{p}^2$, de sorte que la réduction modulo $\mathfrak{m}_\mathfrak{p}^2$ induit un isomorphisme $\mathrm{Der}_{\tilde{\varphi}}(A_\mathfrak{p},K)\xrightarrow{\simeq}\mathrm{Der}_{\bar{\varphi}}(A_\mathfrak{p}/\mathfrak{m}_\mathfrak{p}^2,K)$, où $\bar{\varphi} : A_\mathfrak{p}/\mathfrak{m}_\mathfrak{p}^2\to K$ désigne le morphisme d'anneaux induit par $\tilde{\varphi}$. La conclusion résulte maintenant du lemme~\ref{lm-deriv_loc}.
\end{proof}

Si $M$ est un monoïde et $w : M\to K_\mu$ un morphisme, on appelle \emph{dérivation non additive}\index{termin}{derivation@dérivation!non additive|textbf} de $M$ dans $K$ relative à $w$ toute fonction $\mathrm{d} : M\to K$ telle que
\[\forall(x,y)\in M^2\qquad\mathrm{d}(xy)=w(x).\mathrm{d}(y)+w(y).\mathrm{d}(x).\]

Si $\tilde{w}$ désigne l'image de $w$ par la bijection canonique $\mon(M,K_\mu)\xrightarrow{\simeq}\mathbf{Ann}(\mathbb{Z}[M],K)$, le $K$-espace vectoriel des dérivations non additives $M\to K$ relatives à $w$ s'identifie à $\mathrm{Der}_{\tilde{w}}(\mathbb{Z}[M],K)$. Nous noterons simplement $\mathfrak{m}_w$ pour $\mathfrak{m}_{\tilde{w}}$\label{pagemw}, c'est-à-dire le noyau du morphisme de $K$-algèbres $K[M]\to K$ prolongeant $w$.

\begin{exem}\label{ex-deriv_grp}
 Si $M$ est un \emph{groupe}, alors $w$ est à valeurs dans $K^\times$, et une fonction $\mathrm{d} : M\to K$ est une dérivation non additive relative à $w$ si et seulement si la fonction $\mathrm{d}/w$ est un morphisme de $M$ dans le groupe additif sous-jacent à $K$.
\end{exem}

\begin{rema}\label{rq-derp0}
Une dérivation non additive $\mathrm{d} : M\to K$ relative à $w$ satisfait la relation $\mathrm{d}(x^n)=n w(x)^{n-1}\mathrm{d}(x)$ pour tous $n\in\mathbb{N}^*$ et $x\in M$. En particulier, on a $\mathrm{d}(x^p)=0$ pour tout $x\in M$.
\end{rema}

\section{Catégories linéaires, abéliennes}

\subsection*{Catégories $A$-linéaires}

Une structure \emph{$A$-linéaire} sur une catégorie $\A$ est un enrichissement de $\A$ sur la catégorie monoïdale symétrique $(A\Md,\underset{A}{\otimes},A)$. Autrement dit, les ensembles de morphismes de $\A$ sont munis d'une structure de $A$-module de sorte que les applications de composition soient $A$-bilinéaires. Une catégorie $\mathbb{Z}$-linéaire est encore appelée catégorie \emph{préadditive}. (Une catégorie \emph{additive} est une catégorie préadditive possédant des coproduits finis ; dans ce cas, la structure préadditive est unique.)

Le centre $Z(\A)$\index{nota}{Z@$Z$ \emph{(centre)}}\index{termin}{centre} d'une catégorie préadditive\,\footnote{possédant un ensemble générateur, par exemple, afin d'éviter tout problème ensembliste.} $\A$ est un anneau (commutatif). La donnée d'une structure $A$-linéaire sur une catégorie $\A$ équivaut à la donnée d'une structure préadditive sur $\A$ et d'un morphisme d'anneaux $A\to Z(\A)$.

Si $\A$ est une catégorie additive et $A$-linéaire, on dispose d'un foncteur (unique à isomorphisme près) $\underset{A}{\otimes} : \mathbf{P}(A)\times\A\to\A$ \label{not-ptA} qui est additif et $A$-linéaire par rapport à chaque variable et tel que $A\underset{A}{\otimes}-\simeq\mathrm{Id}_\A$. On dispose de plus d'isomorphismes $\A(V\underset{A}{\otimes}x,y)\simeq\mathrm{Hom}_A(V,\A(x,y))$ et $\A(V^*\underset{A}{\otimes} x,y)\simeq\A(x,V\underset{A}{\otimes}y)$ naturels en les objets $x$, $y$ de $\A$ et $V$ de $\mathbf{P}(A)$, où $V^*:=\mathrm{Hom}_A(V,A)$. 

\label{ppgam} Si $\gamma : A\to B$ est un morphisme d'anneaux et $\A$ une catégorie $A$-linéaire, on note $\A\otimes_A B$, ou $\gamma^*\A$\index{nota}{G@$\gamma^*\A$ \emph{($\A$ catégorie $A$-linéaire, $\gamma : A\to B$ morphisme d'anneaux)}}, la catégorie $B$-linéaire ayant les mêmes objets que $\A$ et dont les morphismes sont donnés par $(\A\otimes_A B)(x,y):=\A(x,y)\otimes_A B$, la composition et les unités étant induites par celles de $\A$. On dispose d'un foncteur canonique $\gamma_* : \A\to\gamma^*\A$ égal à l'identité sur les objets, qui est $\gamma$-linéaire (c'est-à-dire $A$-linéaire, où $\gamma^*\A$ est munie de la structure $A$-linéaire image inverse par $\gamma$ de sa structure $B$-linéaire). Si $\A$ est additive, alors $\gamma^*\A$ est également additive. Pour $\A=\mathbf{P}(A)$, $\gamma_*$ s'identifie au foncteur $\mathbf{P}(A)\to\mathbf{P}(B)$ d'extension des scalaires le long de $\gamma$. Si $\gamma : A\twoheadrightarrow A/I$ est la réduction modulo un idéal $I$ de $A$, on notera $\A/I$ pour $\gamma^*\A$. 

Toute structure $A$-linéaire sur une catégorie fournit une action du monoïde $A_\mu$ sur celle-ci.

Si $S$ est une partie multiplicative de $A$ (c'est-à-dire un sous-monoïde de $A_\mu$), l'anneau de fractions $A[S^{-1}]$ a pour monoïde multiplicatif sous-jacent $A_\mu[S^{-1}]$ ; de même, si $\A$ est une catégorie $A$-linéaire, la catégorie avec action de $A_\mu[S^{-1}]$ sous-jacente à la catégorie $A[S^{-1}]$-linéaire $\A\otimes_A A[S^{-1}]$ est $\A[S^{-1}]$ : toutes les notions de localisation concordent.

\subsection*{Catégories abéliennes et sous-catégories}

Pour les généralités sur les catégories abéliennes, on pourra se référer à \cite{Gabriel,Pop}.

\begin{defi}\label{df-soucat} Soient $\E$ une catégorie abélienne et $\E'$ une sous-catégorie pleine de $\E$.
On dit que $\E'$ est une sous-catégorie :
\begin{itemize}
\item \emph{semi-épaisse}\index{termin}{semiepaisse@semi-épaisse \emph{(sous-catégorie)}} si elle est stable par sous-quotients et par sommes directes finies ;
\item \emph{épaisse}\index{termin}{epaisse@épaisse \emph{(sous-catégorie)}} si $\E'$ est semi-épaisse et stable par extensions ;
\item \emph{prélocalisante}\index{termin}{prelocalisante@prélocalisante \emph{(sous-catégorie)}|textbf} si $\E'$ est semi-épaisse et que le foncteur d'inclusion $\E'\to\E$ admet un adjoint à droite ;
\item \emph{prébilocalisante}\index{termin}{prebilocalisante@prébilocalisante \emph{(sous-catégorie)}|textbf} si $\E'$ est prélocalisante et que l'inclusion $\E'\to\E$ admet un adjoint à gauche ;
\item \emph{localisante}\index{termin}{localisante \emph{(sous-catégorie)}} si $\E'$ est prélocalisante et épaisse ;
\item \emph{bilocalisante}\index{termin}{bilocalisante \emph{(sous-catégorie)}} si $\E'$ est prébilocalisante et épaisse.
\end{itemize}
\end{defi}

Si $\E'$ est une sous-catégorie semi-épaisse de $\E$, l'inclusion $i : \E'\to\E$ admet un adjoint à droite (resp.  à gauche) si et seulement si tout objet $X$ de $\E$ possède un plus sous grand sous-objet (resp. quotient) $\Phi(X)$ appartenant à $\E'$ ; $\Phi$ est alors l'adjoint à droite (resp. à gauche) de $i$, et l'unité $i\Phi\to\mathrm{Id}_\E$ (resp. la coünité $\mathrm{Id}_\E\to i\Phi$) est alors le monorphisme (resp. l'épimorphisme) canonique.

Si $\E'$ est une sous-catégorie semi-épaisse de $\E$, la sous-catégorie pleine de $\E$ des objets possédant une filtration finie dont les sous-quotients appartiennent à $\E'$ est épaisse, c'est la plus petite sous-catégorie épaisse de $\E$ contenant $\E'$, on l'appelle donc sous-catégorie épaisse de $\E$ engendrée par $\E'$.

\subsection*{Catégories de Grothendieck}

On rappelle qu'une catégorie de Grothendieck est une catégorie abélienne cocomplète, possédant un générateur, et dans laquelle les colimites filtrantes sont exactes \cite[p.~63, note~3]{Pop}. Il est bien connu qu'une catégorie de Grothendieck possède toujours des produits \cite[chapitre~3, corollaire~7.10]{Pop} et un cogénérateur injectif \cite[chapitre~3, théorème~10.10 et lemme~7.12]{Pop}. En particulier, on peut y définir des groupes d'extensions. Les sous-objets d'un objet d'une catégorie de Grothendieck forment toujours un ensemble par \cite[chap.~I, prop.~5]{Gabriel}.

Une sous-catégorie semi-épaisse d'une catégorie de Grothendieck est prélocalisante (resp. prébilocalisante) si et seulement si elle est stable par sommes directes (resp. produits) arbitraires. Une sous-catégorie prélocalisante d'une catégorie de Grothendieck est elle-même une catégorie de Grothendieck.

Si $\E'$ est une sous-catégorie semi-épaisse d'une catégorie de Grothendieck $\E$, la sous-catégorie pleine de $\E$ des objets qui sont colimite d'une famille d'objets de $\E'$ est prélocalisante, c'est la sous-catégorie prélocalisante de $\E$ engendrée par $\E'$. Tout objet de cette catégorie est la colimite filtrante de ses sous-objets appartenant à $\E'$.

Différentes propriétés de finitude classiques d'un objet $X$ d'une catégorie de Grothendieck, ou plus généralement d'une catégorie abélienne $\E$ dans laquelle la classe des sous-objets d'un objet donné forme un treillis complet, interviendront de façon importante dans ce travail. Il s'agit notamment des notions d'objet \emph{noethérien}\index{termin}{noetherien@noethérien}, \emph{artinien}\index{termin}{artinien}, \emph{simple}\index{termin}{simple@simple \emph{(objet d'une catégorie abélienne)}}, \emph{de longueur finie} (que nous nommerons simplement \emph{fini}\index{termin}{fini@fini \emph{(objet d'une catégorie abélienne)}} ; cette condition équivaut à \emph{noethérien et artinien}), \emph{de type fini}\index{termin}{type fini@type fini \emph{(objet de type fini d'une catégorie abélienne)}}, ou \emph{de type cofini}\index{termin}{type cofini@type cofini \emph{(objet de type cofini d'une catégorie abélienne)}}\index{termin}{cofini!type}. Toutes ces notions sont classiques et se trouvent par exemple rappelées dans \cite[§\,3.2]{Alb}, à l'exception de la dernière. Un objet de $\E$ est par définition de type cofini si pour toute famille de sous-objets d'intersection nulle, il existe une sous-famille finie d'intersection nulle. Un objet d'une catégorie de Grothendieck est de type cofini si et seulement s'il est extension essentielle d'un sous-objet fini. Comme ces propriétés de finitude d'un objet $X$ de $\E$ ne dépendent que du treillis de ses sous-objets, si $X$ appartient à une sous-catégorie semi-épaisse $\E'$ de $\E$, elles sont équivalentes dans $\E$ et $\E'$.

Un objet d'une catégorie abélienne est dit \emph{localement noethérien} (resp. \emph{localement fini},  \emph{localement de type fini}) s'il est somme d'un ensemble de sous-objets noethériens (resp. finis, de type fini). Une catégorie abélienne est dite \emph{localement noethérienne}\index{termin}{localement noetherienne@localement noethérienne \emph{(catégorie abélienne)}|textbf} (resp. \emph{localement finie}\index{termin}{localement finie@localement finie \emph{(catégorie abélienne)}|textbf}, \emph{localement de type fini}\index{termin}{localement de type fini \emph{(catégorie abélienne)}}) si tous ses objets sont localement noethériens (resp. localement finis, localement de type fini).

Pour $n\in\mathbb{N}\cup\{\infty\}$, on dit qu'un objet $X$ d'une catégorie de Grothendieck $\E$ est de \emph{$n$-présentation finie}\index{termin}{presentation finie@présentation finie \emph{(objet de $n$-présentation finie, $pf_n$)}}\label{pagepfn} --- en abrégé, $pf_n$ --- si le morphisme naturel
$$\underset{J}{\col}\mathrm{Ext}^i_\E(X,-)\circ\Phi\to\mathrm{Ext}^i_\E(X,\underset{J}{\col}\Phi)$$
est bijectif pour $i<n$ et injectif pour $i=n$ (si $n$ est fini), pour tout foncteur $\Phi : J\to\E$, où $J$ est une petite catégorie \emph{filtrante}. Pour davantage de détails sur cette notion, on renvoie à \cite[§\,1]{DT-schw}, dont on suivra les notations --- en particulier, on désignera par $\pf_n(\E)$\index{nota}{pf@$\pf_n$|textbf} la classe des objets de $\E$ vérifiant la propriété $pf_n$.

\subsection*{Catégories abéliennes $k$-linéaires (cf. \cite[début de l'exposé~5]{SGA5})}

Si $k$ est un anneau et $\E$ une catégorie abélienne $k$-linéaire cocomplète, le bifoncteur $\underset{k}{\otimes} : \mathbf{P}(k)\times\E\to\E$ s'étend de façon unique (à isomorphisme près) en un bifoncteur encore noté $\underset{k}{\otimes} : k\Md\times\E\to\E$ qui est cocontinu par rapport à chaque variable. Il est caractérisé par l'isomorphisme d'adjonction
$$\A(V\underset{k}{\otimes}a,b)\simeq\mathrm{Hom}_k(V,\A(a,b))$$
naturel en les objets $a$, $b$ de $\A$ et $V$ de $k\Md$.

Si $k'$ est une $k$-algèbre non nécessairement commutative, et $\E$ une catégorie  abélienne $k$-linéaire, alors on note $\E\underset{k}{\boxtimes} k'$\index{nota}{$\boxtimes$} la catégorie dont les objets sont les couples $(X,f)$ où $X$ est un objet de $\E$ et $f : k'\to\mathrm{End}_\E(X)$ un morphisme de $k$-algèbres (non nécessairement commutatives). Les morphismes $(X,f)\to (Y,g)$ sont les morphismes $u : X\to Y$ de $\E$ tels que le diagramme suivant commute.
$$\xymatrix{k'\ar[r]^-f\ar[d]_g & \mathrm{End}_\E(X)\ar[d]^{u_*}\\
\mathrm{End}_\E(Y)\ar[r]^-{u^*} & \mathrm{Hom}_\C(X,Y)
}$$

La catégorie $\E\underset{k}{\boxtimes} k'$ est abélienne. Si $k'$ est commutative, cette catégorie est $k'$-linéaire. Si $\E$ est de Grothendieck, il en est de même pour $\E\underset{k}{\boxtimes} k'$. Si $\E$ est cocomplète, alors le foncteur d'oubli $\E\underset{k}{\boxtimes} k'\to\E$ possède un adjoint à gauche noté $-\underset{k}{\otimes} k'$ (sa composée avec le foncteur d'oubli s'identifie à l'endofoncteur $k'\underset{k}{\otimes}-$ de $\E$). Si $\E$ est de Grothendieck et que $k'$ est un $k$-module plat, alors $-\underset{k}{\otimes} k'$ est exact.

Si $\E$ est cocomplète, le foncteur $k'\underset{k}{\otimes}-$ induit également un foncteur $k$-linéaire pleinement fidèle  $\E\otimes_k k'\to \E\underset{k}{\boxtimes} k'$, mais celui-ci n'est généralement pas une équivalence (par exemple, pour $\E=\mathbf{Ab}$ et $k=\mathbb{Z}$, son image essentielle est constituée des $k'$-modules libres). Toutefois, si le morphisme d'anneaux $k\to k'$ est un épimorphisme, ce foncteur est une équivalence. C'est en particulier le cas pour un morphisme surjectif ou pour le morphisme canonique $k\to k[S^{-1}]$, où $S$ est une partie multiplicative de $k$.

Un cas particulièrement important de ce type de changement de base est celui où $k'$ est un localisé de l'anneau $k$. Nous utiliserons dans ce cadre le résultat simple suivant :
\begin{lemm}\label{lm-loc-ent}
Soient $S$ une partie multiplicative d'un anneau $k$, $\E$ une catégorie abélienne $k$-linéaire et $X$ un objet de $\E\underset{k}{\boxtimes}k[S^{-1}]$ ; notons $\iota : \E\underset{k}{\boxtimes}k[S^{-1}]\to\E$ le foncteur d'oubli. 

Supposons que $\operatorname{End}(X)$ est une $k$-algèbre (non commutative) \emph{entière}. Alors tout sous-objet de $\iota(X)$ est l'image par $\iota$ d'un sous-objet de $X$. En particulier, si $X$ est de type fini (resp. noethérien, fini, simple...), alors il en est de même pour $\iota(X)$.
\end{lemm}

(On notera que, si $\E$ est cocomplète, alors $\E\underset{k}{\boxtimes}k[S^{-1}]$ s'identifie à $\E[S^{-1}]$.)

\begin{proof}
Un objet $T$ de $\E$ appartient à l'image essentielle de $\iota$ si et seulement l'image de $S$ par le morphisme de $k$-algèbres $k\to Z(\operatorname{End}(T))$ est inversible. Comme par ailleurs $\iota$ est pleinement fidèle, la conclusion découle du lemme élémentaire d'algèbre commutative ci-dessous.
\end{proof}

\begin{lemm}
Soient $S$ une partie mutiplicative d'un anneau $k$ et $R$ une $k[S^{-1}]$-algèbre. Si $R$ est entière comme $k$-algèbre, alors $k$ et $k[S^{-1}]$ ont la même image dans $R$.
\end{lemm}

\begin{proof}
Soit $s\in S$. Comme l'image de $s^{-1}$ dans $R$ annule un polynôme unitaire à coefficients dans $k$, c'est un polynôme (à coefficients dans $k$) en l'image de $s\in k$. 
\end{proof}

\subsection*{Décompositions en somme directe d'objets indécomposables}

Considérons les propriétés suivantes d'un objet $X$ d'une catégorie de Grothendieck $\E$ :

(Dec) {\rm L'objet $X$ est isomorphe à une somme directe d'objets indécomposables de $\E$.}\index{nota}{Dec@(Dec), (DecF), (DecFF)}

(DecF) {\rm L'objet $X$ est isomorphe à une somme directe finie d'objets indécomposables de $\E$.}

(DecFF) {\rm Il existe un entier $d$ tel que toute décomposition de $X$ en somme directe de sous-objets non nuls a au plus $d$ termes.}\index{nota}{Dec@(Dec), (DecF), (DecFF)}

On a manifestement les implications (DecFF)$\Rightarrow$(DecF)$\Rightarrow$(Dec).

\begin{rema} Un objet de type fini de $\E$ qui vérifie la condition (Dec) vérifie nécessairement (DecF).
\end{rema}

\begin{defi}\label{df-goldie} On dit qu'un objet de $\E$ est \emph{de dimension uniforme (ou de Goldie) finie} s'il ne contient pas de somme directe infinie d'objets non nuls.
\end{defi}

\begin{exem} Tous les objets noethériens et tous les objets de type cofini de $\E$ sont de dimension uniforme finie.
\end{exem}

Pour la propriété classique suivante, on pourra consulter \cite[§\,2.2 et 3.2]{Alb}.

\begin{prop} Si $X$ est un objet de $\E$ de dimension uniforme finie, alors il existe $d\in\mathbb{N}$ tel que tout sous-objet de $X$ est somme directe d'au plus $d$ sous-objets non nuls.
\end{prop}

\begin{coro} Tout objet de $\E$ de dimension uniforme finie vérifie la condition $\mathrm{(DecFF)}$.
\end{coro}

\begin{coro}\label{cor-ntcf} Tous les objets noethériens et tous les objets de type cofini de $\E$ vérifient la condition $\mathrm{(DecFF)}$.
\end{coro}

\begin{rema} Un objet de type fini ne vérifie pas nécessairement la condition (Dec). Par exemple, si $A$ est le produit d'une infinité d'anneaux non nuls, le $A$-module de type fini $A$ ne vérifie pas (Dec).
\end{rema}

\begin{prop}\label{pr-decffdf} Si $\E$ est $K$-linéaire, tout objet $X$ de $\E$ tel que $\dim_K\mathrm{End}_\E(X)$ soit finie vérifie la condition $\mathrm{(DecFF)}$.
\end{prop}

\begin{proof} En effet, si $X$ est isomorphe à une somme directe de $d$ objets non nuls, alors $\dim_K\mathrm{End}_\E(X)\ge d$.
\end{proof}

\section{Catégories de foncteurs}

Soient $\C$ et $\E$ des catégories, avec $\C$ essentiellement petite. Si $\E$ est cocomplète (resp. complète, additive, abélienne, de Grothendieck), alors il en est de même pour la catégorie de foncteurs $\fct(\C,\E)$.\index{nota}{fct@$\fct$}

Si $\Phi : \D\to\C$ est un foncteur entre catégories essentiellement petites, on notera $\Phi^* : \fct(\C,\E)\to\fct(\D,\E)$ le foncteur de précomposition par $\Phi$. Si $\Psi : \E\to\E'$ est un foncteur, on notera $\Psi_* : \fct(\C,\E)\to\fct(\C,\E')$ le foncteur de postcomposition par $\Psi$.

Si $k$ est un anneau, on note \index{nota}{FB@$\F$} $\F(\C;k)$ pour $\fct(\C,k\Md)$. On note aussi $\F(A,k)$ pour $\F(\mathbf{P}(A);k)$.

Si $\C$ et $\D$ sont des catégories essentiellement petites, on note
\index{nota}{$\boxtimes$}$$\boxtimes : \F(\C;K)\times\F(\D;K)\to\F(\C\times\D;K)$$
le \emph{produit tensoriel extérieur}\index{termin}{produit tensoriel extérieur} (sur $K$), donné par $(F\boxtimes G)(c,d):=F(c)\otimes G(d)$. Un autre produit tensoriel extérieur apparaîtra quelquefois : si $\E$ est une catégorie de Grothendieck $k$-linéaire et $k\to k'$ un morphisme d'anneaux, on dispose d'une équivalence canonique $\fct(\C,\E)\underset{k}{\boxtimes}k'\simeq\fct(\C,\E\underset{k}{\boxtimes}k')$ de catégories de Grothendieck $k'$-linéaires (qui se spécialise notamment en $\F(\C;k)\underset{k}{\boxtimes}k'\simeq\F(\C;k')$). Nous identifierons souvent ainsi ces deux catégories.

Le produit tensoriel que nous utiliserons le plus souvent sera le produit tensoriel interne à $\F(\C;K)$, calculé au but, et noté simplement $\otimes$ : $(F\otimes G)(c)=F(c)\otimes G(c)$ --- ou encore, $F\otimes G$ est la précomposition de $F\boxtimes G$ par le foncteur diagonale $\C\to\C\times\C$.

Si $c$ est un objet de $\C$, on note \index{nota}{P@$P^c_\C$, $P^c$ \emph{(foncteurs projectifs représentables)}} $P^{c;k}_\C$\label{pproj} l'objet $k[\C(c,-)]$ de $\F(\C;k)$. La mention à la catégorie $\C$ sera souvent omise quand aucune ambiguïté ne peut en résulter, de même que la mention à l'anneau $k$ (surtout lorsque $k=K$). Ce foncteur est projectif de type fini car il représente l'évaluation en $c$, par une variante du lemme de Yoneda. Les $P^{c;k}_\C$ constituent une famille génératrice de la catégorie $\F(\C;k)$ lorsque $c$ parcourt un squelette de $\C$.

Nous utiliserons à plusieurs reprises les notions de {\em support}\index{termin}{support} ou de {\em co-support}\index{termin}{co-support} d'un foncteur, pour lesquelles nous renvoyons à \cite[déf.~1.2]{DTV}. Dans $\F(\C;k)$, un ensemble $E$ d'objets de $\C$ constitue un support d'un foncteur si et seulement s'il est quotient d'une somme directe de foncteurs de la forme $P^c_\C$ avec $c\in E$.

\subsection*{Foncteurs à valeurs de dimension (ou de type) fini(e)}

On note \index{nota}{FB@$\F^\df$} $\F^\df(\C;K)$ la sous-catégorie pleine de $\F(\C;K)$ constituée des foncteurs dont les valeurs sont des espaces vectoriels de dimension finie sur $K$.

L'observation immédiate suivante sera d'usage courant.

\begin{prop}\label{pr-tfdfend} Si $F$ est un objet de type fini de $\F(\C;K)$ et $G$ un objet de $\F^\df(\C;K)$, alors l'espace vectoriel $\mathrm{Hom}(F,G)$ est de dimension finie. En particulier, si $F$ est un objet de type fini de $\F^\df(\C;K)$, la $K$-algèbre (non nécessairement commutative) $\mathrm{End}(F)$ est de dimension finie.
\end{prop}

\begin{proof}
Comme $F$ est de type fini, c'est un quotient d'une somme directe finie de foncteurs de la forme $P^c_\C$, donc $\mathrm{Hom}(F,-)$ est un sous-foncteur d'une somme directe finie de foncteurs d'évaluation. 
\end{proof}

Par la proposition~\ref{pr-decffdf}, on en déduit :

\begin{coro}\label{endofini-dec} Tout foncteur de type fini de $\F^\df(\C;K)$ vérifie la condition $\mathrm{(DecFF)}$.
\end{coro}

La proposition~\ref{pr-tfdfend} admet la généralisation immédiate suivante :
\begin{prop}\label{pr-tftfend} Soient $k$ un anneau noethérien. Si $F$ est un objet de type fini de $\F(\C;k)$ et $G$ un objet de $\F(\C;k)$ à valeurs de type fini, alors le $k$-module $\mathrm{Hom}(F,G)$ est de type fini. En particulier, si $F$ est un objet de type fini et à valeurs de type fini de $\F(\C;k)$, la $k$-algèbre (non nécessairement commutative) $\mathrm{End}(F)$ est finie (i.e. de type fini comme $k$-module).
\end{prop}

Nous aurons besoin en fin de mémoire du résultat suivant :
\begin{prop}\label{pr-loc_fcttftf}
Soient $k$ un anneau noethérien, $M$ un monoïde opérant sur une petite catégorie $\C$ et $F$ un foncteur de type fini et à valeurs de type fini de $\F(\C[M^{-1}];k)$. Alors l'image de $F$ dans $\F(\C;k)$ par la précomposition par le foncteur canonique $\C\to\C[M^{-1}]$ est de type fini.
\end{prop}

\begin{proof}
La catégorie $\F(\C;k)$ est $k[M]$-linéaire, et $\F(\C[M^{-1}];k)$ s'identifie à $\F(\C;k)[M^{-1}]\simeq\F(\C;k)\underset{k[M]}{\boxtimes}k[M,M^{-1}]$. La proposition~\ref{pr-tftfend} montre que l'algèbre $\mathrm{End}(F)$ est finie, donc entière, sur l'anneau noethérien $k$, et a fortiori sur $k[M]$. Par conséquent, le lemme~\ref{lm-loc-ent} donne la conclusion.
\end{proof}

On dit \cite[déf.~3.4]{DTV} que $K$ est un {\em corps de décomposition non additif}\index{termin}{corps de décomposition} de la catégorie $\C$ si, pour tout foncteur \emph{simple} $S$ de $\F^\df(\C;K)$, le morphisme canonique de corps (non commutatif au but) $K\to\mathrm{End}(S)$ est un isomorphisme. Lorsque $\C$ est de plus \emph{additive}, on dit que $K$ est {\em un corps de décomposition} de $\C$ si la condition précédente vaut pour les foncteurs \emph{additifs} simples $S$ de $\F^\df(\C;K)$.

\begin{exem}
\begin{enumerate}
    \item Si $K$ est algébriquement clos, c'est un corps de décomposition non additif de $\C$.
    \item  On peut montrer (mais c'est non trivial, contrairement à l'observation précédente) que tout corps de décomposition $K$ d'une petite catégorie additive qui contient les racines de l'unité (i.e. tel que $\nu(K)=\mathbb{N}^*$) en est également un corps de décomposition non additif \cite[cor.~5.7]{DTV}.
    \item Le corps $K$ est un corps de décomposition de la catégorie $\mathbf{P}(A)$ si et seulement si tout morphisme d'anneaux de $A$ vers un corps qui est une extension \textit{finie} de $K$ est à valeurs dans $K$.
\end{enumerate}
\end{exem}

\subsection*{Propriétés de finitude}

Les résultats ici rappelés constitueront l'un des fondements essentiels des théorèmes de finitude (propriétés localement noethérienne ou localement finie) que nous établirons, dans des situations beaucoup plus générales, à la fin de ce travail.

\begin{theo}[Putman-Sam-Snowden]\label{th-PSS} Si $A$ est un anneau fini (éventuellement non commutatif) et $\E$ une catégorie de Grothendieck localement noethérienne,\index{termin}{localement noetherienne@localement noethérienne \emph{(catégorie abélienne)}} alors la catégorie $\fct(\mathbf{P}(A),\E)$ est localement noethérienne.
\end{theo}

\begin{proof}
Lorsque $\E$ est une catégorie de modules sur un anneau noethérien, le résultat est démontré dans \cite{PSam} ou \cite{SamSn}. Le cas d'une catégorie but localement noethérienne générale s'établit exactement de la même façon, et apparaît explicitement dans \cite{Dja-bki}.
\end{proof}

\begin{prop}\label{pr-PSmodif} Si $A$ est un anneau fini (éventuellement non commutatif) de cardinal inversible dans $K$ et $\E$ une catégorie de Grothendieck $K$-linéaire localement finie, alors la catégorie $\fct(\mathbf{P}(A),\E)$ est localement finie.\index{termin}{localement finie@localement finie \emph{(catégorie abélienne)}}
\end{prop}

\begin{proof} Combiner \cite[Proposition~11.7]{DTV} (qui repose lourdement sur le théorème~\ref{th-PSS}) et \cite[fin du §\,3.2]{DT-schw}.
\end{proof}

\subsection*{Foncteurs additifs}
Si $\A$ et $\B$ sont des catégories préadditives, avec $\A$ essentiellement petite, on note $\mathbf{Add}(\A,\B)$\index{nota}{Add@$\mathbf{Add}$} la sous-catégorie pleine des foncteurs additifs de $\fct(\A,\B)$. Pour $d\in\mathbb{N}$, on note $\mathbf{Add}_d(\A,\B)$\index{nota}{Add@$\mathbf{Add}_d$|textbf} la sous-catégorie pleine de $\fct(\A^d,\B)$ des multifoncteurs qui sont additifs par rapport à chacune des $d$ variables. Ainsi, via l'équivalence de catégories canonique $\fct(\A^d,\B)\simeq\fct(\A,\fct(\A^{d-1},\B))$, on a $\mathbf{Add}_d(\A,\B)\simeq\mathbf{Add}(\A,\mathbf{Add}_{d-1}(\A,\B))$.

Si $R$ est un anneau non nécessairement commutatif et $\E$ une catégorie abélienne cocomplète, l'évaluation en $R$ induit une équivalence de catégories 
\begin{equation}\label{eq-pteAdd}
\mathbf{Add}(\mathbf{P}(R),\E)\xrightarrow{\simeq}\E\underset{\mathbb{Z}}{\boxtimes}R^\op\,,
\end{equation}
d'où pour tout $d\in\mathbb{N}$
\begin{equation}\label{eq-pteAddd}
\mathbf{Add}_d(\mathbf{P}(R),\E)\xrightarrow{\simeq}\E\underset{\mathbb{Z}}{\boxtimes}(R^{\otimes_\mathbb{Z}d})^\op\,,
\end{equation}
et pour tout anneau non nécessairement commutatif $k$ une équivalence
\begin{equation}\label{eq-pteAdda}
\mathbf{Add}(\mathbf{P}(R),k\Md)\simeq (k\underset{\mathbb{Z}}{\otimes}R^\op)\Md
\end{equation}
(\emph{théorème d'Eilenberg-Watts}).

\chapter{Catégories de Grothendieck $A$-linéaires et torsion}\label{storsion}

\begin{cvi} Dans tout le chapitre~\ref{storsion}, la lettre $I$ désigne un idéal de $A$, et $\E$ désigne une catégorie de Grothendieck $A$-linéaire.
\end{cvi}

Les notions et résultats présentés dans ce chapitre sont bien connus et constituent une généralisation immédiate aux catégories de Grothendieck $A$-linéaires de résultats de base sur la $I$-torsion des $A$-modules, qui donne lieu notamment à l'étude de la cohomologie locale (voir à ce sujet l'ouvrage \cite{Sch-Si}). Précisément, toutes les notions de ce chapitre s'obtiennent par produit tensoriel avec $\E$ au-dessus de $A$ de leurs déclinaisons dans $A\Md$, au sens du produit tensoriel de catégories de Grothendieck introduit dans \cite{LRGS}. Nous nous contentons ici d'un point de vue concret et élémentaire.

\section{Objets de $I$-torsion}

\begin{defi} 
\begin{enumerate}
\item On note \index{nota}{U@$U(\E;I)$} $U(\E;I)$ la sous-catégorie pleine des objets $M$ de $\E$ tels que le morphisme d'anneaux canonique $A\to\mathrm{End}_\E(M)$ soit nul sur $I$. Les objets de $U(\E;I)$ sont nommés objets de {\em $I$-torsion forte} de $\E$.

Autrement dit, $U(\E;I)$ s'identifie canoniquement à la catégorie $\E/I$, ou encore à $\E\underset{A}{\boxtimes}A/I$.
\item On note \index{nota}{V@$V^\cf(\E;I)$} $V^\cf(\E;I)$ la sous-catégorie pleine de $\E$ réunion sur $n\in\mathbb{N}$ des sous-catégories $U(\E;I^n)$ (où $I^n$ désigne la puissance $n$-ième de l'idéal $I$).
\item On note \index{nota}{V@$V(\E;I)$} $V(\E;I)$ la sous-catégorie pleine des objets de $\E$ qui sont colimites d'objets appartenant à $V^\cf(\E;I)$. Les objets de $V(\E;I)$ sont nommés objets de {\em $I$-torsion faible} de $\E$, ceux de $V^\cf(\E;I)$ objets de {\em $I$-torsion faible à caractère fini}.
\end{enumerate}
\end{defi}

La proposition suivante est immédiate.

\begin{prop}\label{pr-gal-adj} La sous-catégorie $U(\E;I)$ de $\E$ est prébilocalisante.\index{termin}{prebilocalisante@prébilocalisante \emph{(sous-catégorie)}} On note $X\mapsto X_{A/I}$ (resp. $X\mapsto X^{A/I}$) le foncteur $\E\to U(\E;I)$ adjoint à gauche (resp. à droite) de l'inclusion.

Pour tout objet $X$ de $\E$, l'unité $X\to X_{A/I}$ (resp. la coünité $X^{A/I}\to X$) est un épimorphisme (resp. un monomorphisme) qui identifie $X_{A/I}$ (resp. $X^{A/I}$) au plus grand quotient (resp. sous-objet) de $X$ appartenant à $U(\E;I)$.
\end{prop}

Concrètement, le noyau de l'unité $X\to X_{A/I}$ est le sous-objet $I.X$ de $X$ somme des images des endomorphismes $\lambda_X$ de $X$ (donnant l'action de $\lambda$ sur cet objet) pour $\lambda\in I$.

Les trois énoncés faciles suivants sont laissés en exercice.

\begin{lemm}\label{lm-idproduit} Soient $I, J$ des idéaux de $A$ et $X$ un objet de $\E$. Les assertions suivantes sont équivalentes :
\begin{enumerate}
\item $X$ appartient à $U(\E;IJ)$ ;
\item il existe une suite exacte courte $0\to X'\to X\to X''\to 0$ de $\E$ avec $X'$ dans $U(\E;I)$ et $X''$ dans $U(\E;J)$.
\end{enumerate}
\end{lemm}

\begin{prop}\label{pr-gal-ev2} Soit $X$ un objet de $\E$. Les assertions suivantes sont équivalentes :
\begin{enumerate}
\item $X$ appartient à $U(\E;I^n)$ ;
\item il existe une filtration de longueur $n$
\[0=T_0\subset T_1\subset\dots\subset T_n=X\]
dont les sous-quotients $T_i/T_{i-1}$ (pour $i\in\llbracket 1, n\rrbracket$) appartiennent à $U(\E;I)$.
\end{enumerate}
\end{prop}

\begin{coro}\label{cor-uvcf}
\begin{enumerate}
\item La sous-catégorie $V^\cf(\E;I)$ de $\E$ est la sous-catégorie épaisse engendrée par $U(\E;I)$.
\item La sous-catégorie $V(\E;I)$ de $\E$ est la sous-catégorie prélocalisante\index{termin}{prelocalisante@prélocalisante \emph{(sous-catégorie)}} engendrée par $V^\cf(\E;I)$.
\item\label{ituvcf3} Tout objet de type fini de $V(\E;I)$ appartient à $V^\cf(\E;I)$.
\item Si l'inclusion $I^2\subset I$ d'idéaux de $A$ est une égalité, alors l'inclusion $U(\E;I)\subset V(\E;I)$ est une égalité ; réciproquement, cette égalité pour $\E=A\Md$ implique $I^2=I$.
\end{enumerate} 
\end{coro}

\begin{rema}
La \emph{cohomologie locale} de $\E$ par rapport à $I$ est par définition le foncteur cohomologique obtenu en dérivant à droite l'adjoint à droite de l'inclusion $V(\E;I)\hookrightarrow\E$. Nous n'en aurons pas besoin explicitement dans ce travail.
\end{rema}

La proposition ci-dessous est classique et significative des problèmes de finitude sous-jacents aux situations dans lesquelles nous étudierons la torsion dans ce mémoire, situations où l'hypothèse de type-finitude de l'idéal $I$ ou de noethérianité locale de la catégorie ambiante ne seront presque jamais vérifiées. Elle ne sera toutefois pas utilisée dans l'article, c'est pourquoi nous ne ferons qu'esquisser sa démonstration.

\begin{prop}Si $\E$ est \index{termin}{localement noetherienne@localement noethérienne \emph{(catégorie abélienne)}} localement noethérienne, ou que $I$ est un idéal de type fini de $A$, alors $V(\E;I)$ est une sous-catégorie localisante de $\E$.
\end{prop}

\begin{proof} Le cas où $\E$ est localement noethérienne résulte de la propriété générale bien connue que toute sous-catégorie d'une catégorie localement noethérienne obtenue en saturant par colimites une sous-catégorie épaisse est localisante.

Si $I$ est un idéal de type fini, disons engendré par $n$ éléments, on se ramène d'abord au cas où $n=1$ à l'aide des observations suivantes : si $I$ et $J$ sont deux idéaux de $A$, alors $U(\E;I+J)=U(\E;I)\cap U(\E;J)$ et $V(\E;I+J)=V(\E;I)\cap V(\E;J)$. Par ailleurs, pour tout élément $a$ de $A$, $V(\E;(a))$ est le noyau du foncteur canonique $\E\to\E[a^{-1}]$, qui est \emph{exact} et cocontinu, c'est donc bien une sous-catégorie localisante de $\E$.
\end{proof}

\begin{rema}
Sans hypothèse de finitude sur $I$ ni sur $\E$, la sous-catégorie $V(\E;I)$ n'est généralement pas stable par extensions --- cf. exemple~\ref{ex-Vpasepais} ci-après.
\end{rema}

La propriété suivante, qui constitue une variante du théorème d'Hopkins-Levitzki, nous servira au chapitre~\ref{sec-finpolp} (par l'intermédiaire de son corollaire).

\begin{prop}\label{prqpft} Soit $X$ un objet noethérien de $\E$.

Alors $X$ appartient à $V^\cf(\E;\nil(A))$. Si de plus la catégorie $U(\E;\nil(A))$ est localement finie\index{termin}{localement finie@localement finie \emph{(catégorie abélienne)}}, alors $X$ est fini.
\end{prop}

\begin{proof} Comme $X$ est noethérien, il existe un idéal de type fini $I$ de $A$, inclus dans $\nil(A)$, tel que $\nil(A).X=I.X$. On en déduit l'égalité $\nil(A)^n.X=I^n.X$ par récurrence sur $n\in\mathbb{N}$. Comme $I$ est de type fini et constitué d'éléments nilpotents, il existe $n$ tel que $I^n=0$, d'où $\nil(A)^n.X=0$. Autrement dit, $X$ appartient à $U(\E;\nil(A)^n)\subset V^\cf(\E;\nil(A))$.

La proposition~\ref{pr-gal-ev2} montre que si $U(\E;\nil(A))$ est localement finie, alors $X$ est fini (les sous-quotients de la filtration $(T_i)$ donnée par la proposition sont noethériens et dans $U(\E;\nil(A))$, donc finis).
\end{proof}

\begin{coro}\label{cor-hole} Soient $\E'$ une catégorie de Grothendieck $K$-linéaire et $B$ une $K$-algèbre quasi-parfaite. On suppose également que la $K$-algèbre $B/\rj(B)$ est de dimension finie.

Alors un objet de $\E'\underset{K}{\boxtimes} B$ est noethérien si et seulement si son image par le foncteur d'oubli $\E'\underset{K}{\boxtimes} B\to\E'$ est noethérienne.

Si de plus la catégorie $\E'$ est localement finie, alors tout objet noethérien de $\E'\underset{K}{\boxtimes} B$ est fini.
\end{coro}

%\begin{rema}  Si $I/I^2$ est un $A/I$-module de type fini, alors tout objet de type fini de $V(\E;I)$ a une filtration finie dont les sous-quotients sont de type fini et dans $U(\E;I)$.
%\end{rema}

\section{Propriétés cohomologiques élémentaires}

L'un des principes imprécis mais importants de la théorie de la torsion dans les catégories abéliennes est que \guillemotleft~les objets de torsions assez différentes ne se mélangent pas \guillemotright. La proposition suivante en donne une illustration cohomologique.

\begin{prop}\label{pr-etr} Soient $I$ et $J$ deux idéaux de $A$ tels que $I+J=A$. Si $X$ (resp. $Y$) est un objet de $V(\E;I)$ (resp. $V(\E;J)$), alors $\mathrm{Ext}^*_\E(X,Y)=0$.
\end{prop}

\begin{proof}
Commençons par établir ce résultat lorsque $X$ appartient à $U(\E;I)$. Comme $I+J=A$, il existe $a\in I$ tel que $1-a\in J$. Il s'ensuit que l'image de $a$ dans l'anneau $A/J^n$ est inversible pour tout $n\in\mathbb{N}$, ainsi l'action de $a$ est inversible sur tout objet de $V^\cf(\E;J)$, et donc sur tout objet de $V(\E;J)$. Par ailleurs, les foncteurs $\mathrm{Ext}^i_\E$ se relèvent (au but) aux $A$-modules : ainsi l'action de $a$ sur $\mathrm{Ext}^*_\E(X,Y)$ est nulle car elle l'est sur $X$, qui appartient à $U(\E;I)$, mais elle est aussi inversible, puisqu'elle l'est sur $Y$, qui appartient à $V(\E;J)$. Par conséquent, $\mathrm{Ext}^*_\E(X,Y)=0$.

Montrons maintenant le résultat lorsque $X$ appartient à $V^\cf(\E;I)$. La relation $I+J=A$ entraîne $I^n+J=A$ pour tout $n\in\mathbb{N}$, de sorte que ce cas se déduit directement du précédent.

Montrons maintenant le cas général : $X$ possède une filtration croissante exhaustive par des sous-objets de $V^\cf(\E;I)$, et est donc quotient d'une somme directe d'objets de $V^\cf(\E;I)$. On en déduit par récurrence qu'il existe un complexe homologique de $\E$ dont tous les termes sont des sommes directes d'objets de $V^\cf(\E;I)$ et dont l'homologie est isomorphe à $X$, concentrée en degré $0$. Les suites spectrales d'hypercohomologie correspondantes permettent alors de conclure en se ramenant au cas précédent.
\end{proof}

Un autre aspect de l'étude cohomologique de la torsion consiste à comparer les groupes d'extensions dans $\E$ et $V(\E;I)$ (il n'y a aucun espoir que les $\mathrm{Ext}^1$ coïncident dans $\E$ et dans $U(\E;I)$, sauf si $U(\E;I)=V(\E;I)$, d'après le corollaire~\ref{cor-uvcf}).

Nous aurons besoin pour cela de la version suivante du lemme d'Artin-Rees.

\begin{lemm}\label{artin-rees} Soient $n\in\mathbb{N}$, $Y$ un objet \textbf{noethérien} de $\E$ et $X$ un sous-objet de $Y$. Alors il existe un entier $m\ge n$ tel que $X\cap (I^m.Y)\subset I^n.X$.
\end{lemm}

\begin{proof} Quitte à remplacer $I$ par $I^n$, on peut supposer $n=1$. On peut également supposer que $I$ est un idéal \textbf{de type fini}, car la noethérianité de $X$ et $Y$ entraîne l'existence d'un idéal de type fini $J$ de $A$, inclus dans $I$, tel que $J.X=I.X$ et $J.Y=I.Y$, égalités qui entraînent $J^m.X=I^m.X$ et $J^m.Y=I^m.Y$ par récurrence sur $m\in\mathbb{N}$.

Supposons donc que $I$ est engendré par $r$ éléments de $A$, et montrons le lemme par récurrence sur $r$.

Traitons d'abord le cas $r=1$ : $I$ est engendré par un élément $f$ de $A$. Étant donné $m\in\mathbb{N}$, posons $T_m=(f^m_Y)^{-1}(X)$ : $(T_m)_{m\in\mathbb{N}}$ est une suite croissante de sous-objets de $Y$, elle stationne donc. Par ailleurs, $I^m.T_m=X\cap (I^m.Y)$. Par conséquent, l'égalité $T_{m-1}=T_m$ implique
\[X\cap (I^m.Y)=I.(I^{m-1}.T_{m-1})=I.(X\cap (I^{m-1}.Y))\subset I.X.\]

Supposons maintenant $r$ quelconque ; écrivons $I=J+I'$, où $J$ (resp. $I'$) est un idéal engendré par au plus $r-1$ (resp. $1$) éléments. L'hypothèse de récurrence fournit $m\in\mathbb{N}^*$ tel que $X\cap (J^m.Y)\subset J.X$. Soit $Y'$ l'objet noethérien $Y/(J^m.Y)$ de $\E$ et $X'$ son sous-objet $X/(X\cap (J^m.Y))$. Ce qui précède fournit $l\in\mathbb{N}^*$ tel que $X'\cap (I'^l.Y')\subset I'.X'$, c'est-à-dire $X\cap (I'^l.Y+J^m.Y)\subset I'.X+J^m.Y$, d'où
\[X\cap (I'^l.Y+J^m.Y)\subset X\cap (I'.X+J^m.Y)=I'.X+(X\cap J^m.Y)\subset I'.X+J.X=I.X.\]

Par ailleurs, $I'^l.Y+J^m.Y=(I'^l+J^m).Y$, et $I^{l+m-1}=(I'+J)^{l+m-1}\subset I'^l+J^m$, d'où le lemme.
\end{proof}

\begin{prop}\label{preserv-inj} Supposons que la catégorie $\E$ est \index{termin}{localement noetherienne@localement noethérienne \emph{(catégorie abélienne)}} \textbf{localement noethérienne}. Alors le foncteur d'inclusion $V(\E;I)\to\E$ préserve les objets injectifs ; il induit un isomorphisme entre groupes d'extensions.
\end{prop}

\begin{proof} On note tout d'abord que, d'après le corollaire~\ref{cor-uvcf}, $V(\E;I)$ est une catégorie de Grothendieck, donc en particulier une catégorie abélienne avec assez d'objets injectifs, il suffit donc de montrer que l'inclusion $V(\E;I)\to\E$ (qui est exacte) préserve les objets injectifs.

Soit $L$ un objet injectif de $V(\E;I)$. Pour $n\in\mathbb{N}$, notons $L_n:=L^{A/I^n}$ : alors $(L_n)$ est une suite croissante de sous-objets de $L$, de colimite $L$, et $L_n$ est un objet injectif de la catégorie abélienne $U(\E;I^n)$.

Soient $Y$ un objet noethérien, $X$ un sous-objet de $Y$ et $f\in\E(X,L_n)$. Choisissons un entier $m\ge n$ tel que $X\cap (I^m.Y)\subset I^n.X$ (lemme~\ref{artin-rees}). Alors le morphisme $f$ se factorise à travers la projection $X\twoheadrightarrow X/(X\cap (I^m.Y))$. Comme $L_m$ est injectif dans $U(\E;I^m)$ et que l'inclusion $X\subset Y$ induit un monomorphisme $X/(X\cap (I^m.Y))\hookrightarrow Y/(I^m.Y)\simeq Y_{A/I^m}$ de $U(\E;I^m)$, il s'ensuit qu'on peut former un diagramme commutatif
\[\xymatrix{X\ar[rr]^f\ar@{->>}[rd]\ar[dd]^\subset & & L_n\ar[dd]^\subset \\
& X/(X\cap (I^m.Y)\ar@{-->}[ur]\ar[dd] & \\
Y\ar@{->>}[rd] & & L_m \\
& Y/(I^m.Y)\ar@{-->}[ur] &
}\]

Comme tout morphisme de $X$ vers $L$ se factorise par l'inclusion $L_n\subset L$ pour $n$ assez grand (car $X$ est de type fini), on en déduit que tout morphisme $X\to L$ s'étend à $Y$, dès lors que $Y$ est noethérien. Comme $\E$ est localement noethérienne, cela implique que $L$ est injectif, d'où la proposition.
\end{proof}

L'hypothèse localement noethérienne est forte ; la comparaison cohomologique entre $V(\E;I)$ est difficile en général sans hypothèse de finitude, et ce dès le degré cohomologique $1$ : la sous-catégorie $V(\E;I)$ de $\E$ n'est pas nécessairement épaisse.

\begin{exem}\label{ex-Vpasepais} Considérons l'anneau $A$ quotient de l'algèbre de polynômes $K[x_1,x_2,\dots,x_n,\dots]$ en une infinité dénombrable d'indéterminées sur le corps $K$ par l'idéal engendré par les éléments $x_n x_{i_1}\dots x_{i_n}$, pour tous $n\in\mathbb{N}^*$ et $(i_1,\dots,i_n)\in (\mathbb{N}^*)^n$. Alors $A$ est un anneau local dont le radical $I$, image dans $A$ de l'idéal d'augmentation de $K[x_1,x_2,\dots,x_n,\dots]$, vérifie : $\forall t\in I\quad\exists n\in\mathbb{N}\quad I^n t=0$. Il s'ensuit que le $A$-module $I$ appartient à $V(A\Md;I)$. Par conséquent, le $A$-module $A$ appartient à la sous-catégorie épaisse de $A\Md$ engendrée par $V(A\Md;I)$ (laquelle coïncide donc avec $A\Md$), en vertu de la suite exacte courte de $A$-modules $0\to I\to A\to A/I\to 0$.

Mais le $A$-module $A$ n'appartient pas à $V(A\Md;I)$, car $I$ n'est pas un idéal nilpotent (car l'image de $x_n^n$ dans $A$ n'est pas nulle). 
\end{exem}

\section{Décomposition selon la torsion}\label{sdst}

Au-delà des objets de torsion, ce sont les objets possédant une décomposition en somme directe de sous-objets de torsion relativement à des idéaux appropriés qui s'avéreront particulièrement utiles. L'énoncé ci-dessous donne des propriétés générales simples de telles décompositions.

\begin{prop}\label{poixt} Soit $\M$ un ensemble d'idéaux de $A$ tel que $I+J=A$ dès que $I$ et $J$ sont deux éléments distincts de $\M$ (par exemple, un ensemble d'idéaux \emph{maximaux}). Considérons le foncteur
\[\Phi : \prod_{I\in\M}V(\E;I)\to\E\qquad(X_I)_{I\in\M}\mapsto\bigoplus_{I\in\M}\iota_I(X_I),\]
où $\iota_I : V(\E;I)\to\E$ désigne le foncteur d'inclusion.
\begin{enumerate}
\item Le foncteur $\Phi$ est exact et pleinement fidèle. Il commute aux colimites.
\item Supposons que $\M$ est constitué de tous les idéaux maximaux de $A$ et que $A$ est semi-primaire. Alors $\Phi$ est une équivalence de catégories.
\item Supposons que $\E$ est \index{termin}{localement noetherienne@localement noethérienne \emph{(catégorie abélienne)}} localement noethérienne. Étant donnés des objets $\mathrm{X}=(X_I)_{I\in\M}$ et $\mathrm{Y}=(Y_I)_{I\in\M}$ de $\prod_{I\in\M}V(\E;I)$, le morphisme canonique
\[\prod_{I\in\M}\mathrm{Ext}^*_{V(\E;I)}(X_I,Y_I)\simeq\mathrm{Ext}^*(\mathrm{X},\mathrm{Y})\to\mathrm{Ext}^*_\E(\Phi(\mathrm{X}),\Phi(\mathrm{Y}))\]
est un isomorphisme.
\end{enumerate}
\end{prop}

\begin{proof} Comme chaque foncteur $\iota_I$ est exact, pleinement fidèle, et commute aux colimites, il suffit d'établir la nullité de $\E\Big(\iota_I(X_I),\underset{J\in\M\setminus\{I\}}{\bigoplus}\iota_J(Y_J)\Big)$ pour montrer la première assertion. Comme les colimites filtrantes sont exactes dans $\E$, ce groupe abélien s'injecte dans
\[\E\Big(\iota_I(X_I),\prod_{J\in\M\setminus\{I\}}\iota_J(Y_J)\Big)\simeq\prod_{J\in\M\setminus\{I\}}\E(\iota_I(X_I),\iota_J(Y_J)),\]
qui est nul d'après la proposition~\ref{pr-etr}.

Sous les hypothèses du deuxième point, $A$ est isomorphe à un produit fini d'anneaux locaux de radicaux nilpotents, d'où l'on tire que tout $A$-module $M$ se décompose de façon unique et naturelle en somme directe de sous-modules annulés par une puissance d'un idéal maximal de $A$, ce qui entraîne l'essentielle surjectivité de $\Phi$.

Supposons maintenant que $\E$ est localement noethérienne. La proposition~\ref{preserv-inj} montre que chaque foncteur $\iota_I$ induit un isomorphisme entre groupes d'extensions, il suffit donc pour conclure d'établir la nullité de $\mathrm{Ext}^*_\E\Big(\iota_I(X_I),\underset{J\in\M\setminus\{I\}}{\bigoplus}\iota_J(Y_J)\Big)$. Si $X_I$ est noethérien, ce groupe abélien est isomorphe à $\underset{J\in\M\setminus\{I\}}{\bigoplus}\mathrm{Ext}^*_\E(\iota_I(X_I),\iota_J(Y_J))$ \cite[prop.~1.3]{DT-schw}, donc nul par la proposition~\ref{pr-etr}. Le cas général s'en déduit par un argument formel en écrivant $X_I$ comme la colimite filtrante de ses sous-objets noethériens (cf. fin de la démonstration de la proposition~\ref{pr-etr}).
\end{proof}

\begin{rema} On ne peut se passer d'hypothèse de finitude forte de type noethérianité locale dans la dernière assertion de la proposition~\ref{poixt}. Nous en verrons une illustration dans l'exemple~\ref{exextsdi}.
\end{rema}

\begin{rema} Lorsque $\E=A\Md$, la réciproque de la deuxième assertion est vraie : si $\Phi$ est une équivalence de catégories ($\M$ étant l'ensemble des idéaux maximaux de $A$), alors $A$ est semi-primaire. On vérifie en effet très facilement que l'appartenance du $A$-module $A$ à l'image essentielle de $\Phi$ implique que $A$ est un anneau semi-primaire.
\end{rema}

La proposition suivante donne une situation où l'existence d'une décomposition en sous-objets de torsion relativement aux idéaux maximaux de $A$ est garantie pour tous les objets vérifiant la condition (Dec).

\begin{prop}\label{poivN} Supposons que l'anneau $A$ est absolument plat\index{termin}{absolument plat} et notons $\M$ l'ensemble des idéaux maximaux de $A$. Alors le foncteur
\[\Psi : \prod_{I\in\M}U(\E;I)\to\E\qquad(X_I)_{I\in\M}\mapsto\bigoplus_{I\in\M}\gamma_I(X_I),\]
où $\gamma_I : U(\E;I)\to\E$ désigne le foncteur d'inclusion, est pleinement fidèle et son image essentielle contient tous les objets de $\E$ vérifiant la condition {\rm (Dec)}.
\end{prop}

\begin{proof} Compte-tenu de la proposition~\ref{poixt}, il suffit de montrer que tout objet indécomposable $X$ de $\E$ appartient à l'image essentielle de $\Psi$. L'image du morphisme d'anneaux canonique $A\to Z(\mathrm{End}_\E(X))$ est un anneau absolument plat (car c'est un quotient de $A$), et n'a pas d'idempotent non trivial (car c'est un sous-anneau de $Z(\mathrm{End}_\E(X))$), c'est donc un corps. Son noyau est un idéal maximal $I$ de $A$, et $X$ appartient à $U(\E;I)$.
\end{proof}

\begin{rema} Pour $\E=A\Md$ (où $A$ est un anneau absolument plat), la réciproque est vraie : un $A$-module qui appartient à l'image essentielle de $\Psi$ vérifie la condition {\rm (Dec)}. En particulier, si $A$ n'est pas semi-simple, le $A$-module $A$ n'appartient pas à cette image.
\end{rema}

\chapter{Poids dans les catégories de foncteurs}\label{spcf}

\begin{cvi} Dans tout le chapitre~\ref{spcf}, $\C$ désigne une catégorie essentiellement petite munie d'une action d'un monoïde $M$ et $\E$ une catégorie de Grothendieck $K$-linéaire.
\end{cvi}

Nous nous appliquerons dans ce chapitre à la catégorie de foncteurs $\fct(\C,\E)$, munie de la structure $K[M]$-linéaire pour laquelle les éléments de $M$ opèrent à la source et les éléments de $K$ par multiplication au but, les notions de torsion du chapitre précédent. Le cas particulier le plus important sera celui de la catégorie $\F(\C;K)$. Les \emph{décompositions en poids} (faibles ou fortes) que nous introduirons seront ainsi des cas particuliers des décompositions selon la torsion du §\,\ref{sdst}.

Il s'agit d'une part de donner de nouveaux points de vue sur la théorie esquissée au chapitre précédent, notamment à l'aune de ce que l'algèbre linéaire (réduction simultanée des endomorphismes des $K$-espaces vectoriels) permet, mais aussi d'aller plus loin dans les critères d'existence de telles décompositions. D'autre part, il s'agit de donner plusieurs énoncés de changement de base (changement de corps au but, changement de monoïde agissant à la source...) qui se révéleront importants dans plusieurs de nos considérations ultérieures (au moins pour certains aspects techniques). 

\section{Définitions}

 Les morphismes de monoïdes de $M$ vers le monoïde multiplicatif $K_\mu$ seront aussi appelés \index{termin}{poids} {\em poids}. Un tel morphisme $w$ s'étend de façon unique en un morphisme de $K$-algèbres $\tilde{w} : K[M]\to K$, dont on rappelle (cf. page~\pageref{pagemw}) que $\mathfrak{m}_w$\index{nota}{m@$\mathfrak{m}_w$ \emph{(pour $w\in\mon(M,K_\mu)$)}} désigne le noyau, qui est donc un idéal \emph{maximal} de $K[M]$. 

\begin{defi}\label{defipoi} Soient $w : M\to K_\mu$ un poids et $F$ un foncteur de $\fct(\C,\E)$.
\begin{enumerate}
\item On dit que $F$ est \index{termin}{homogene@homogène!\emph{(foncteur homogène de poids fort ou faible $w$)}}\textbf{homogène de poids fort} $w$ s'il est non nul et appartient à la sous-catégorie pleine $U(\fct(\C,\E);\mathfrak{m}_w)$ de $\fct(\C,\E)$ --- autrement dit, si, pour tout objet $x$ de $\C$ et tout élément $\lambda$ de $M$, l'endomorphisme $\lambda_x$ de $x$ est envoyé par $F$ sur la multiplication par $w(\lambda)$.
La catégorie $U(\fct(\C,\E);\mathfrak{m}_w)$ est notée $\fct(\C,\E)_w$.\index{nota}{F@$\fct(\C,\E)_w$}
\item On dit que $F$ est \textbf{homogène de poids faible} $w$ s'il est non nul et appartient à la sous-catégorie pleine $V(\fct(\C,\E);\mathfrak{m}_w)$ de $\fct(\C,\E)$. Cette sous-catégorie est notée $\fct(\C,\E)_{<w>}$.\index{nota}{F@$\fct(\C,\E)_{<w>}$}
\item On dit que $F$ admet une \textbf{décomposition en poids} forte (resp. faible) s'il est isomorphe à une somme directe, éventuellement infinie, de foncteurs homogènes pour certains poids forts (resp. faibles).
\item On dit que $F$ admet une décomposition en poids forte (resp. faible) finie s'il possède une  décomposition en somme directe finie (resp. une filtration finie) dont chaque terme (resp. chaque sous-quotient) est homogène d'un certain poids fort.
\end{enumerate}
\end{defi}

\begin{rema} La condition de non-nullité est destinée à garantir l'unicité du poids d'un foncteur homogène. Par ailleurs, la première assertion de la proposition~\ref{poixt} montre qu'une décomposition en poids d'un foncteur est nécessairement unique (à isomorphisme près).
\end{rema}

\begin{prdef}\label{poids-evident} Soit $w : M\to K_\mu$ un poids.
\begin{enumerate}
\item La sous-catégorie $\fct(\C,\E)_w$ de $\fct(\C,\E)$ est prébilocalisante.\index{termin}{prebilocalisante@prébilocalisante \emph{(sous-catégorie)}} C'est en particulier une sous-catégorie abélienne de $\fct(\C,\E)$. L'adjoint à droite (resp. à gauche) du foncteur d'inclusion $\fct(\C,\E)_w\to\fct(\C,\E)$ est noté $F\mapsto F^w$ (resp. $F\mapsto F_w$).
\item\label{itprdfr} Plus généralement, pour $r\in\mathbb{N}$, on note  $F\mapsto F_w^{[r]}$ l'adjoint à gauche du foncteur d'inclusion $U(\fct(\C,\E);\mathfrak{m}_w^r)\to\fct(\C,\E)$.
\item La sous-catégorie $\fct(\C,\E)_{<w>}$ de $\fct(\C,\E)$ est prélocalisante.\index{termin}{prelocalisante@prélocalisante \emph{(sous-catégorie)}}
On note $F\mapsto F^{<w>}$ l'adjoint à droite du foncteur d'inclusion $\fct(\C,\E)_{<w>}\to\fct(\C,\E)$.
%Si $M$ est de type fini, ou que $\fct(\C,\E)$ est localement noethérienne, alors $\fct(\C,\E)_{<w>}$ est une sous-catégorie localisante de $\fct(\C,\E)$. En général, si $0\to F\to G\to H\to 0$ est une suite exacte de $\fct(\C,\E)$ avec $F$ et $H$ de poids faible $w$, et que de plus $F$ est de type fini, alors $G$ est de poids faible $w$.
\item Le morphisme canonique
\[\bigoplus_{\pi\in\mon(M,K_\mu)}F^{<\pi>}\to F\]
est un monomorphisme. Le foncteur canonique
$$\prod_{\pi\in\mon(M,K_\mu)}\fct(\C,\E)_{<\pi>}\to\fct(\C,\E)$$
(donné par la somme directe des projections composées avec les inclusions) est pleinement fidèle.
\end{enumerate}
\end{prdef}

\begin{proof} Le caractère prébilocalisant (resp. prélocalisant) de la sous-catégorie $\fct(\C,\E)_w$ (resp. $\fct(\C,\E)_{<w>}$) de $\fct(\C,\E)$ est un cas particulier de la proposition~\ref{pr-gal-adj} (resp. du corollaire~\ref{cor-uvcf}). La dernière assertion provient pour sa part de la proposition~\ref{poixt}.
\end{proof}

\begin{rema}\label{rq-decfini} Une décomposition en poids (faible ou forte) d'un foncteur de type fini est nécessairement finie. Par ailleurs, un foncteur possède une décomposition en poids (faible ou forte) si et seulement s'il est colimite de sous-foncteurs possédant une décomposition en poids (faible ou forte) finie.

En particulier, si la catégorie $\E$ est localement de type fini, alors un foncteur de $\fct(\C,\E)$ possède une décomposition en poids (faible ou forte) si et seulement si tous ses sous-foncteurs de type fini possèdent une décomposition en poids (faible ou forte) finie. On peut ainsi souvent se limiter à étudier les décompositions en poids finies.
\end{rema}

\begin{exem}\label{ex-pgqwP} Soient $c$ un objet de $\C$ et $w\in\mon(M,K_\mu)$. Le foncteur $(P^c_\C)_w$ est donné par $(P^c_\C)_w(t)=K[\C(c,t)]\underset{K[M]}{\otimes}K_w$, où $K_w$ désigne le $K$-espace vectoriel $K$ muni de l'action de $M$ donnée par $w$ ; plus généralement, le plus grand quotient de $P^c_\C$ appartenant à $U(\fct(\C,\E);\mathfrak{m}_w^r)$ (où $r\in\mathbb{N}$) est donné par $t\mapsto K[\C(c,t)]\underset{K[M]}{\otimes}(K[M]/\mathfrak{m}_w^r)$. 
\end{exem}

\begin{rema}
Par la proposition~\ref{pr-gal-ev2}, un foncteur $F$ possède une décomposition en poids faible finie si et seulement s'il existe $r,n\in\mathbb{N}$ et des poids deux à deux distincts $w_1,\dots,w_n : M\to K_\mu$ tels que le morphisme canonique $F\twoheadrightarrow\bigoplus_{i=1}^n F^{[r]}_{w_i}$ soit un isomorphisme.
\end{rema}

\begin{defi}\label{nota-enspoi} Étant donné un foncteur $F$ de $\fct(\C,\E)$, on note \index{nota}{Pi@$\Pi_M$, $\Pi$ \emph{(ensemble des poids d'un foncteur relativement à l'action d'un monoïde $M$ sur sa catégorie source)}}$\Pi_M(F)$, ou simplement $\Pi(F)$, l'ensemble des morphismes de monoïdes $w : M\to K_\mu$ tels que $F^w\ne 0$. Les éléments de $\Pi(F)$ sont appelés poids de $F$.\index{termin}{poids!d'un foncteur}
\end{defi}

\begin{rema} La condition $F^w\ne 0$ équivaut à $F^{<w>}\ne 0$. Si de plus $F$ admet une décomposition en poids faible finie, elle équivaut également à $F_w\ne 0$.
\end{rema}

Si $F$ est un foncteur possédant une décomposition en poids faible, et si $x\in M$ est tel que $w(x)\in K^\times$ pour tout $w\in\Pi_M(F)$, alors $F(x)$ est inversible dans $Z(\mathrm{End}(F))$. Par conséquent, le lemme~\ref{lm-inversion_fct} implique : 
\begin{prop}\label{pr-inversion_fct} Soient $F$ un foncteur de $\fct(\C,\E)$ possédant une décomposition en poids faible (resp. forte) et $S$ un sous-monoïde de $M$ tels que :
$$\forall w\in\Pi_M(F)\quad\forall x\in S\qquad w(x)\in K^\times\;.$$
Alors il existe un foncteur $\tilde{F}$ de $\fct(\C[S^{-1}],\E)$, unique à isomorphisme près, tel que $F$ soit isomorphe à la composée de $\tilde{F}$ et du foncteur canonique $\C\to\C[S^{-1}]$.

De plus, $\tilde{F}$ possède une décomposition en poids faible (resp. forte) relativement à l'action de $M[S^{-1}]$, et l'on a $\Pi_{M[S^{-1}]}(\tilde{F})=\{\tilde{w}\,|\,w\in\Pi_M(F)\}$, où $\tilde{w} : M[S^{-1}]\to K_\mu$ désigne l'unique morphisme de monoïdes tel que $w$ soit la composée de $\tilde{w}$ et du morphisme canonique $M\to M[S^{-1}]$.
\end{prop}

\section{Interprétation en termes d'algèbres d'endomorphismes}

\begin{nota}\label{nota-cZ} Soit $F$ un foncteur de $\fct(\C,\E)$. On note \index{nota}{Z@$\mathfrak{Z}$, $\mathfrak{Z}_M$} $\mathfrak{Z}_M(F)$, ou simplement $\mathfrak{Z}(F)$, la $K$-algèbre image du morphisme canonique $K[M]\to Z(\mathrm{End}(F))$.
\end{nota}

L'énoncé suivant est immédiat.

\begin{lemm}\label{pd1ev} Un foncteur $F$ est homogène d'un certain poids fort si et seulement si $\dim_K\mathfrak{Z}(F)=1$. Son poids est alors donné par la composée du morphisme de monoïdes canonique $M\to K[M]_\mu$, de la projection canonique $K[M]\twoheadrightarrow\mathfrak{Z}(F)$ et de l'unique isomorphisme de $K$-algèbres $\mathfrak{Z}(F)\simeq K$.
\end{lemm}

On en déduit en particulier, en utilisant la proposition~\ref{pr-tfdfend} et le fait que l'anneau (non nécessairement commutatif) d'endomorphismes d'un objet simple est un corps :
\begin{prop}[Lemme de Schur]\label{poisimp} Supposons que le corps $K$ est algébriquement clos. Alors tout foncteur simple à valeurs de dimensions finies de $\F(\C;K)$ est homogène d'un certain poids fort $M\to K_\mu$.
\end{prop}

\begin{lemm}\label{lm-decZalg} Soit $F$ un foncteur de $\fct(\C,\E)$. Supposons que la $K$-algèbre $\mathfrak{Z}(F)$ se décompose en produit fini $\mathfrak{Z}(F)\simeq\prod_{i=1}^n A_i$. Alors il existe une décomposition $F\simeq\bigoplus_{i=1}^n F_i$ dans $\fct(\C,\E)$ avec $\mathfrak{Z}(F_i)\simeq A_i$.
\end{lemm}

\begin{proof}
La décomposition en produit de la sous-algèbre $\mathfrak{Z}(F)$ de $Z(\mathrm{End}(F))$ fournit une famille complète d'idempotents centraux orthogonaux $e_1,\dots,e_n$ de $\mathrm{End}(F)$. La conclusion s'obtient en prenant pour $F_i$ l'image de $e_i$.
\end{proof}

L'utile propriété suivante peut être vue comme un analogue de la deuxième assertion de la proposition~\ref{poixt}.

\begin{prop}\label{pr-poideploye} Soit $F$ un foncteur de $\fct(\C,\E)$. 
\begin{enumerate}
\item Le foncteur $F$ possède une décomposition en poids forte finie si et seulement si la $K$-algèbre $\mathfrak{Z}(F)$ est semi-simple et déployée.
\item Le foncteur $F$ possède une décomposition en poids faible finie si et seulement si la $K$-algèbre $\mathfrak{Z}(F)$ est semi-primaire et déployée. Dans ce cas, il existe un diagramme commutatif d'algèbres
$$\xymatrix{K[M]\ar[rr]^-{(\tilde{w})_{w\in\Pi(F)}}\ar@{->>}[d] && K^{\Pi(F)}\\
\mathfrak{Z}(F)\ar@{->>}[rr] && \mathfrak{Z}(F)/\rj(\mathfrak{Z}(F))\ar[u]^\simeq
}$$
\end{enumerate}
\end{prop}

\begin{proof} Notons $R$ la $K$-algèbre $\mathfrak{Z}(F)$. Les lemmes~\ref{lm-decZalg} et~\ref{pd1ev} et montrent que si $R$ est semi-simple déployée, alors $F$ possède une décomposition en poids forte finie, et que le morphisme $K[M]\xrightarrow{(\tilde{w})_{w\in\Pi(F)}}K^{\Pi(F)}$ induit un isomorphisme d'algèbres $R\xrightarrow{\simeq}K^{\Pi(F)}$. Réciproquement il est clair que si $F$ admet une décomposition en poids forte finie, alors $R$ est semi-simple déployée.

Notons $I$ l'idéal de $K[M]$ noyau du morphisme canonique $K[M]\to Z(\mathrm{End}(F))$ et $J$ la racine de $I$. Si $R\simeq K[M]/I$ est semi-primaire déployée, alors $K[M]/J$ est semi-simple déployée et il existe $n\in\mathbb{N}$ tel que $J^n\subset I$. La première assertion et la proposition~\ref{pr-gal-ev2} montrent alors que $F$ possède une filtration finie dont les sous-quotients admettent une décomposition en poids forte finie, ainsi $F$ admet une décomposition en poids faible finie. De plus, si $G$ est le plus grand sous-foncteur de $F$ ayant une décomposition en poids forte finie, on a $\Pi(G)=\Pi(F)$ et $\mathfrak{Z}(G)=R/\rj(R)$, ce qui permet d'obtenir le diagramme commutatif souhaité.

Réciproquement, si $F$ possède une décomposition en poids faible finie, alors $F$ appartient à la sous-catégorie épaisse engendrée par les foncteurs homogènes d'un certain poids fort (cf. corollaire~\ref{cor-uvcf}), donc il existe un nombre fini de poids $w_1,\dots,w_n : M\to K$ tels que le morphisme canonique $K[M]\to Z(\mathrm{End}(F))$ soit nul sur $\mathfrak{m}_{w_1}\dots\mathfrak{m}_{w_n}$, par le lemme~\ref{lm-idproduit}. Il s'ensuit que $R$ est un quotient de la $K$-algèbre $K[M]/\mathfrak{m}_{w_1}\dots\mathfrak{m}_{w_n}$, qui est semi-primaire déployée, donc $R$ est elle-même une algèbre semi-primaire déployée (cf. lemme~\ref{lm-sqalgdepl}).
\end{proof}

\begin{coro}\label{cor-acfdf} Si la $K$-algèbre $\mathfrak{Z}(F)$ est de dimension finie et que le corps $K$ est algébriquement clos, alors $F$ possède une décomposition en poids faible.
\end{coro}

Le résultat suivant est élémentaire mais fondamental.

\begin{coro}\label{cor-tfdp} Supposons $K$ algébriquement clos. Alors tout foncteur $F$ à valeurs de dimensions finies de $\F(\C;K)$ possède une décomposition en poids faible.
\end{coro}

\begin{proof} Si $F$ est de plus de type fini, alors $\mathfrak{Z}(F)$ est de dimension finie, par la proposition~\ref{pr-tfdfend}, et le corollaire~\ref{cor-acfdf} permet de conclure. Le cas général s'obtient en écrivant $F$ comme colimite de ses sous-foncteurs de type fini (cf. remarque~\ref{rq-decfini}).
\end{proof}

L'énoncé suivant généralise le corollaire~\ref{cor-tfdp}.

\begin{coro}\label{cor-tfdp2} Supposons que $K$ est un corps de décomposition non additif de $\C$. Alors tout foncteur $F$ à valeurs de dimensions finies de $\F(\C;K)$ possède une décomposition en poids faible.
\end{coro}

\begin{proof} Comme dans la démonstration du corollaire~\ref{cor-tfdp}, on peut se contenter de montrer le résultat lorsque $F$ est de plus de type fini. Alors $\mathfrak{Z}(F)$ est de dimension finie, par la proposition~\ref{pr-tfdfend}. En particulier, la $K$-algèbre $\mathfrak{Z}(F)$ est semi-primaire, de sorte qu'il suffit de montrer qu'elle est déployée. Le lemme~\ref{lm-decZalg} permet de se ramener au cas où $\mathfrak{Z}(F)$ est également \emph{locale}.

Notons $I$ et $J$ les idéaux de $K[M]$ noyaux des morphismes canoniques $K[M]\twoheadrightarrow\mathfrak{Z}(F)$ et $K[M]\twoheadrightarrow\mathfrak{Z}(F)\twoheadrightarrow L$ respectivement, où $L$ désigne le corps $\mathfrak{Z}(F)/\rj(\mathfrak{Z}(F))$. Alors $I\subset J$, et il existe $n\in\mathbb{N}$ tel que $J^n\subset I$ puisque le noyau de $\mathfrak{Z}(F)\twoheadrightarrow L$ est nilpotent.

Par définition de $\mathfrak{Z}(F)$, $F$ appartient à $U(\F(\C;K);I)\subset U(\F(\C;K);J^n)$. Par la proposition~\ref{pr-gal-ev2}, il s'ensuit que $F$ possède un quotient non nul $G$ dans $U(\F(\C;K);J)$. Alors $G$ est de type fini et à valeurs de dimensions finies car c'est un quotient de $F$, et $\mathfrak{Z}(G)\simeq L$ car $G$ appartient à $U(\F(\C;K);J)$ (donc $\mathfrak{Z}(G)$ est un quotient de $L$) et est non nul. Ainsi, $G$ possède un quotient simple $S$, et l'on a $\mathfrak{Z}(S)\simeq L$ car $S$ est encore un objet non nul de $U(\F(\C;K);J)$. Mais $S$ est à valeurs de dimensions finies, donc absolument simple puisque $K$ est un corps de décomposition non additif de $\C$. Ainsi $\mathfrak{Z}(S)=K$, d'où $L=K$, ce qui montre que $\mathfrak{Z}(F)$ est bien déployée et achève la démonstration.
\end{proof}

\begin{rema}\label{remacdna} La réciproque du corollaire~\ref{cor-tfdp2} est vraie : si tout foncteur, ou même seulement tout foncteur simple, de $\F^\df(\C;K)$ possède une décomposition en poids faible, alors $K$ est un corps de décomposition non additif de $\C$. Cela résulte du lemme~\ref{pd1ev} et de l'observation qu'un foncteur simple possédant une décomposition en poids faible est nécessairement homogène d'un certain poids fort.
\end{rema}

Le résultat suivant permet d'obtenir des décompositions en poids \emph{à extension (finie) des scalaires (au but) près}.

\begin{coro}\label{cor-excb} Supposons que la $K$-algèbre $\mathfrak{Z}(F)$ est de dimension finie.
\begin{enumerate}
\item Il existe une extension de corps finie $K\to K'$ telle que le foncteur $F\underset{K}{\boxtimes}K'$ de $\fct(\C,\E\underset{K}{\boxtimes}K')$ possède une décomposition en poids faible.
\item Si le corps $K$ est parfait, alors $F$ admet une décomposition en somme directe finie de foncteurs $G_i$ tels qu'il existe une extension de corps finie $K\to L_i$ et des foncteurs $T_i$ de $\fct(\C,\E\underset{K}{\boxtimes} L_i)$ possédant une décomposition en poids faible finie relativement à des poids $M\to L_i$ et que $G_i$ soit isomorphe à la post\-composition de $T_i$ et du foncteur d'oubli $\E\underset{K}{\boxtimes} L_i\to\E$.
\end{enumerate}
\end{coro}

\begin{proof} Comme une $K$-algèbre de dimension finie se décompose en produit fini d'algèbres locales, le lemme~\ref{lm-decZalg} montre qu'on peut supposer sans perte de généralité que l'algèbre $\mathfrak{Z}(F)$ est locale. Son corps résiduel $L$ est une extension finie de $K$. 

Pour la première assertion, l'extension quasi-galoisienne $K'$ de $K$ engendrée par $L$ convient.

L'extension $L$ convient pour la deuxième assertion dès lors qu'elle est séparable --- donc en particulier si $K$ est parfait --- car le morphisme de $K$-algèbres $A\twoheadrightarrow L$ est alors scindé d'après un théorème classique de Wedderburn \cite[§\,13, cor.~1 ou~2 de la prop.~7]{Bki}.
\end{proof}

%%%%%%%%%%%%%%%%%%%%%%%%%%%%%%%%%%%%%%%%%%%%%%%%%%%%%%%%%%%%%%%
%% SI VISÉE ENCYCLOPÉDIQUE, RÉINTÉGRER REMARQUE CI-DESSOUS !
% s'il existe une extension de corps \emph{arbitraire} $K\subset L$ telle que $F\underset{K}{\boxtimes}L$ ait une décomposition en poids faible, alors il existe une sous-extension \emph{finie} $K'$ de $L$ telle que $F\underset{K}{\boxtimes}K'$ ait une décomposition en poids faible, mais il n'est sans doute pas utile de le mentionner. La raison : si $R$ est une $K$-algèbre telle que $R\otimes L$ soit semi-primaire déployée, alors il existe une sous-extension finie $K'$ de $L$ telle que $R\otimes K'$ soit semi-primaire déployée --- ce qui se voit en se ramenant au cas où $R$ est un corps et en traitant à part les extensions radicielles.
%%%%%%%%%%%%%%%%%%%%%%%%%%%%%%%%%%%%%%%%%%%%%%%%%%%%%%%%%%%%%%%

\section[Interprétation en termes de réduction d'endomorphismes]{Interprétation en termes de réduction d'endomorphismes linéaires}\label{sredend}

Lorsque $\C$ est la catégorie à un objet $\underline{M}$ (munie de l'action tautologique de $M$), un foncteur de $\fct(\underline{M},\E)$ consiste en la donnée d'un objet $X$ de $\E$ muni d'une action de $M$, c'est-à-dire d'un morphisme de monoïdes (non nécessairement commutatif au but) $M\to\mathrm{End}_\E(X)$, ou encore d'un endomorphisme de $K$-algèbres (non nécessairement commutative au but) $K[M]\to\mathrm{End}_\E(X)$. En particulier, $\F(\underline{M};K)$ s'identifie à $K[M]\Md$.

Supposons que $M$ est le monoïde (commutatif) libre $\mathbb{N}[E]$ sur un ensemble $E$. Un objet de $\fct(\underline{\mathbb{N}[E]},\E)$ consiste en la donnée d'un objet $X$ de $\E$ et d'une famille $(\varphi_i)_{i\in E}$ d'endomorphismes de $X$ commutant deux-à-deux. Cela nous permet de donner la définition suivante.

\begin{defi}\label{df-diatrig} Une famille commutative $(\varphi_i)_{i\in E}$ d'endomorphismes d'un objet de $\E$ est dite  \emph{collectivement diagonalisable} (resp. \emph{collectivement trigonalisable}) si le foncteur correspondant de $\fct(\underline{\mathbb{N}[E]},\E)$ possède une décomposition en poids forte (resp. faible). Un endomorphisme $\varphi$ de $E$ est dit diagonalisable (resp. trigonalisable) si $\{\varphi\}$ est collectivement diagonalisable (resp. collectivement trigonalisable).
\end{defi}

Lorsque $\E=K\Md$, on retrouve la notion usuelle de diagonalisabilité et de trigonalisabilité, au moins pour les espaces vectoriels de dimension finie. La généralisation aux espaces vectoriels de dimension arbitraire, moins étudiée, présente quelques subtilités ; on pourra se référer à \cite{IMR16} pour ce qui concerne la diagonalisabilité et à \cite{Mes18} pour la trigonalisabilité. L'équivalence de la définition~\ref{df-diatrig}, qui généralise directement la décomposition en sous-espaces caractéristiques, avec celle de \cite[déf.~4]{Mes18}, généralisation directe à la dimension infinie de l'existence d'une base où l'endomorphisme possède une matrice triangulaire, est donnée par \cite[théorème~8]{Mes18}.

\begin{rema}\label{rq-diatricol} Une famille commutative d'endomorphismes diagonalisables (resp. trigonalisables) d'un espace vectoriel de dimension infinie n'est pas nécessairement collectivement diagonalisable (resp. collectivement trigonalisable) --- voir par exemple \cite[ex.~2.13 ou~4.17]{IMR16}. En revanche, une famille commutative collectivement trigonalisable d'endomorphismes diagonalisables d'un objet de $\E$ est  collectivement diagonalisable.
\end{rema}

La proposition ci-dessous, qui constitue un simple jeu  d'écriture, est laissée en exercice.

\begin{prop}\label{pr-tridia-gl} Étant donné un foncteur $F$ de $\fct(\C,\E)$, les assertions suivantes sont équivalentes :
\begin{enumerate}
\item $F$ possède une décomposition en poids forte (resp. faible) ;
\item la famille $\{F(\lambda)\,|\,\lambda\in M\}$ d'endomorphismes de $F$ est collectivement diagonalisable (resp. collectivement trigonalisable) ;
\item pour tout objet $t$ de $\C$, la famille $\{F(\lambda_t)\,|\,\lambda\in M\}$ d'endomorphismes de l'objet $F(t)$ de $\E$ est collectivement diagonalisable (resp. collectivement trigonalisable).
\end{enumerate}
\end{prop}

On ne peut généralement pas ôter l'adverbe \emph{collectivement} de l'énoncé précédent ; on a toutefois le résultat partiel suivant.

\begin{prop}\label{pr-dtf} Étant donné un foncteur $F$ de $\fct(\C,\E)$, les conditions suivantes sont équivalentes :
\begin{enumerate}
\item pour tout objet $t$ de $\C$ et tout élément $\lambda$ de $M$, l'endomorphisme $F(\lambda_t)$ de $F(t)$ est diagonalisable (resp. trigonalisable) ;
\item pour tout $\lambda\in M$, l'endomorphisme $F(\lambda)$ de $F$ est diagonalisable (resp. trigonalisable) ;
\item le foncteur $F$ possède une décomposition en poids forte (resp. faible) relativement à l'action de tout sous-monoïde de type fini de $M$.
\end{enumerate}

Si elles sont vérifiées avec des décompositions fortes et que le foncteur $F$ satisfait à la condition {\rm (Dec)}, alors $F$ possède une décomposition en poids forte.
\end{prop}

\begin{proof}
L'équivalence des trois assertions se déduit directement de \cite[§\,4]{Mes18} dans le cas trigonalisable et de \cite[cor.~4.15]{IMR16} dans le cas diagonalisable. Si elles sont vérifiées avec des décompositions fortes et que $F$ est indécomposable, alors pour tout sous-monoïde de type fini $N$ de $M$, $F$ est homogène d'un certain poids fort $w_N : N\to K_\mu$. Les $w_N$ s'assemblent pour fournir un morphisme de monoïdes $w : M\to K_\mu$ ; il est clair que $F$ est homogène de poids fort $w$. Le cas où $F$ satisfait à la condition {\rm (Dec)} se déduit immédiatement du cas indécomposable.
\end{proof}

%Les conditions avec "diagonalisable" équivalent à $\mathfrak{Z}(F)$ absolument plate déployée et entière sur $K$, ou à réduite, déployée et entière sur $K$, mais peu probable que ce soit bien utile !

\section{Changement de monoïde agissant à la source}

Nous omettrons la démonstration de la proposition suivante, immédiate et d'usage courant.

\begin{prop}\label{pr-cmsource} Soient $N$ un monoïde, $\D$ une catégorie essentiellement petite munie d'une action de $N$, $\varphi : N\to M$ un morphisme de monoïdes, et $\Phi : \D\to\C$ un foncteur $\varphi$-équivariant, c'est-à-dire tel que, pour tout objet $x$ de $\D$, le diagramme
$$\xymatrix{N\ar[r]^-\varphi\ar[d] & M\ar[d]\\
\mathrm{End}_\D(x)\ar[r]^-{\Phi_*} & \mathrm{End}_\C(\Phi(x))
}$$
dont les flèches verticales sont données par l'action de $M$ (resp. $N$) sur $\C$ (resp. $\D$) commute.
\begin{enumerate}
\item Pour tout poids $w : M\to K_\mu$, le foncteur $\Phi^* : \fct(\C,\E)\to\fct(\D,\E)$ envoie $\fct(\C,\E)_w$ (resp. $\fct(\C,\E)_{<w>}$) dans $\fct(\D,\E)_{w\circ\varphi}$  (resp. $\fct(\D,\E)_{<w\circ\varphi>}$).
\item Si un foncteur $F$ de $\fct(\C,\E)$ possède une décomposition en poids forte (resp. faible) relativement à l'action de $M$, alors $\Phi^*F$ possède une décomposition en poids forte (resp. faible) relativement à l'action de $N$. On a de plus $\Pi_N(\Phi^*F)\subset\{w\circ\varphi\,|\,w\in\Pi_M(F)\}$, avec égalité lorsque $\Phi^*$ est fidèle (par exemple, lorsque $\Phi$ est essentiellement surjectif à facteur direct près).
\end{enumerate}
\end{prop}

\begin{lemm}\label{lm-dpfmm} Supposons que le monoïde $M$ est engendré par un nombre fini de sous-monoïdes $N_1,\dots,N_i$. Alors un foncteur de $\fct(\C,\E)$ possède une décomposition en poids forte (resp. faible) finie relativement à l'action de $M$ si et seulement s'il en possède une par rapport à chacun des $N_j$.
\end{lemm}

\begin{proof} La $K$-algèbre $\mathfrak{Z}_M(F)$ est engendrée par les sous-algèbres $\mathfrak{Z}_{N_j}(F)$. La conclusion découle donc de la proposition~\ref{pr-poideploye} et du lemme~\ref{lm-sqalgdepl2}.
\end{proof}

\begin{prop}\label{pr-poiprm} Supposons que le monoïde $M$ est engendré par un nombre fini de sous-monoïdes $N_1,\dots,N_i$. Alors un foncteur de $\fct(\C,\E)$ possède une décomposition en poids forte (resp. faible) relativement à l'action de $M$ si et seulement s'il en possède une par rapport à chacun des $N_j$.
\end{prop}

\begin{proof} Cela découle du lemme~\ref{lm-dpfmm} et de la remarque~\ref{rq-decfini}.
\end{proof}

\begin{exem}\label{ex-adj0poi} Si $M\simeq N_0$ (cf. notation~\ref{nota-adj0}) pour un sous-monoïde $N$, alors $M$ est engendré par $N$ et par le sous-monoïde $\{0,1\}$, par rapport auquel tout foncteur possède une décomposition en poids forte finie, qui correspond au scindement de la partie constante d'un foncteur. Ainsi, un foncteur de $\fct(\C,\E)$ possède une décomposition en poids (faible ou forte) relativement à $M$ si et seulement s'il en possède une relativement à $N$.
\end{exem}

Le résultat suivant est immédiat.

\begin{lemm}\label{lm-postcompop} Soient $\E'$ est une catégorie de Grothendieck $K$-linéaire et $\Phi : \E\to\E'$ un foncteur exact $K$-linéaire.
\begin{enumerate}
\item La post-composition $\Phi_* : \fct(\C,\E)\to\fct(\C,\E')$ par $\Phi$ préserve les foncteurs possédant une décomposition en poids faible (resp. forte) finie.
\item Si $\Phi$ est cocontinu, alors $\Phi_*$ préserve les foncteurs possédant une décomposition en poids faible (resp. forte).
\item Si $\Phi$ est fidèle et que $F$ est un foncteur de $\fct(\C,\E)$ tel que $\Phi_*(F)$ possède une décomposition en poids faible (resp. forte) finie, alors $F$ possède une décomposition en poids faible (resp. forte) finie.
\item Si $\Phi$ est fidèle et cocontinu et que $F$ est un foncteur de $\fct(\C,\E)$ tel que $\Phi_*(F)$ possède une décomposition en poids faible (resp. forte), alors $F$ possède une décomposition en poids faible (resp. forte).
\end{enumerate}
\end{lemm}

Le principe suivant de \guillemotleft~séparation des variables~\guillemotright\ pour les décompositions en poids des bifoncteurs nous sera utile à plusieurs reprises.

\begin{prop}\label{pr-poibif} Soient $\D$ une catégorie essentiellement petite munie d'une action d'un monoïde $N$ et $X : \C\times\D\to\E$ un foncteur. Les assertions suivantes sont équivalentes :
\begin{enumerate}
\item $X$ possède une décomposition en poids forte (resp. faible) relativement à l'action de $M\times N$ ;
\item l'image de $X$ par l'isomorphisme de catégories $\fct(\C\times\D,\E)\simeq\fct(\C,\fct(\D,\E))$ possède une décomposition en poids forte (resp. faible) relativement à l'action de $M$ et est à valeurs dans les foncteurs de $\fct(\D,\E)$ ayant une décomposition en poids forte (resp. faible) relativement à l'action de $N$ ;
\item pour tout objet $c$ de $\C$, le foncteur $X(c,-)$ de $\fct(\D,\E)$ possède une décomposition en poids forte (resp. faible) relativement à l'action de $N$ et pour tout objet $d$ de $\D$, le foncteur $X(-,d)$ de $\fct(\C,\E)$ possède une décomposition en poids forte (resp. faible) relativement à l'action de $M$.
\end{enumerate}

De plus, lorsqu'elles sont satisfaites, l'image de la fonction $\Pi_{M\times N}(X)\hookrightarrow\mon(M\times N,K_\mu)\twoheadrightarrow\mon(M,K_\mu)$ égale $\bigcup_{d\in\mathrm{Ob}\,\D}\Pi_N(X(-,d))$.
\end{prop}

\begin{proof} La proposition~\ref{pr-poiprm} montre que $X$ possède une décomposition en poids forte (resp. faible) relativement à l'action de $M\times N$ si et seulement si l'image de $X$ par l'isomorphisme canonique de catégories $K$-linéaires $\fct(\C\times\D,\E)\simeq\fct(\C,\fct(\D,\E))$ possède une décomposition en poids forte (resp. faible) relativement à l'action de $M$ et l'image de $X$ par l'isomorphisme canonique $\fct(\C\times\D,\E)\simeq\fct(\D,\fct(\C,\E))$ possède une décomposition en poids forte (resp. faible) relativement à l'action de $N$.

Pour conclure à l'équivalence des assertions, il suffit de voir que l'image de $X$ par l'isomorphisme canonique $\fct(\C\times\D,\E)\simeq\fct(\C,\fct(\D,\E))$ possède une décomposition en poids forte (resp. faible) relativement à l'action de $M$ si et seulement si, pour tout objet $d$ de $\D$, le foncteur $X(-,d)$ de $\fct(\C,\E)$ possède une décomposition en poids forte (resp. faible) relativement à l'action de $M$. Cela découle du lemme~\ref{lm-postcompop} appliqué au foncteur exact, cocontinu et $K$-linéaire de $\fct(\D,\E)$ dans $\E^T$ donné par l'évaluation en les différents éléments d'un ensemble $T$ de représentants des classes d'isomorphisme d'objets de $\D$. 

La dernière propriété, relative aux ensembles de poids, est claire lorsque $X$ est homogène d'un certain poids $(w,w')\in\mon(M,K_\mu)\times\mon(N,K_\mu)\simeq\mon(M\times N,K_\mu)$, car $X(-,d)$ est alors soit nul, soit homogène de poids $w$, et est non nul pour au moins un $d$ ; le cas général s'en déduit aussitôt en décomposant $X$ est somme directe de foncteurs homogènes.
\end{proof}

\begin{rema} L'assertion analogue obtenue en remplaçant {\em décomposition en poids} par {\em décomposition en poids finie} dans l'énoncé précédent est erronée. Considérons ainsi le cas où $K=\mathbb{Q}$, $\C=\D=\mathbf{P}(\mathbb{Q})$, $X=\underset{n\in\mathbb{N}}{\bigoplus}\Lambda^n\boxtimes\Lambda^n$ (où $\Lambda^n$ désigne la $n$-ième puissance extérieure sur $\mathbb{Q}$). Les monoïdes considérés sont le groupe multiplicatif $M=N=\mathbb{Q}^\times$. Si $V\simeq\mathbb{Q}^d$ est un objet de $\mathbf{P}(\mathbb{Q})$, alors $X(V,-)\simeq X(-,V)$ admet une décomposition en poids forte finie, puisque $\Lambda^n(V)=0$ pour $n>d$, mais la décomposition en poids de $X$ n'est pas finie.

En revanche, les variantes des deux corollaires suivants obtenues en remplaçant \emph{décomposition en poids} (forte ou faible) par \emph{décomposition en poids finie} sont clairement valides. 
\end{rema}

\begin{coro}\label{cor-ptedp} Sous les mêmes hypothèses, soient $F$ un foncteur de $\F(\C;K)$ et $G$ un foncteur de $\F(\D;K)$. Si $F$ (resp. $G$) possède une décomposition en poids, faible ou forte, relativement à l'action de $M$ (resp. $N$), alors il en est de même pour le foncteur $F\boxtimes G$ de $\F(\C\times\D;K)$ relativement à l'action de $M\times N$. Si de plus $F$ et $G$ sont non nuls, la réciproque est vraie, et $\Pi_{M\times N}(F\boxtimes G)$ est l'image de $\Pi_M(F)\times\Pi_N(G)$ par l'isomorphisme canonique $\mon(M,K_\mu)\times\mon(N,K_\mu)\xrightarrow{\simeq}\mon(M\times N,K_\mu)$.
\end{coro}

\begin{coro}\label{cor-poiptint} Soit $(F_i)_{1\le i\le n}$ une famille finie de foncteurs de $\F(\C;K)$. Supposons que, pour tout $i$, $F_i$ possède une décomposition en poids faible (resp. forte). Alors $F_1\otimes\dots\otimes F_n$ possède une décomposition en poids faible (resp. forte). De plus, $\Pi_M(F_1\otimes\dots\otimes F_n)=\Pi_M(F_1)\dots\Pi_M(F_n)$.
\end{coro}

\section{Critères d'existence de décompositions en poids fortes}

En l'absence de lourdes hypothèses de finitude (sur le foncteur ou sur le monoïde $M$), l'existence d'une décomposition en poids d'un foncteur de $\fct(\C,\E)$ constitue un phénomène plutôt exceptionnel (cf. proposition~\ref{pr-dpftjs} ci-après). Le corollaire~\ref{cor-tfdp} ci-avant et ses variantes comptent parmi les critères les plus importants d'existence de décompositions en poids faibles. Nous donnons maintenant quelques critères généraux garantissant l'existence de décompositions en poids \emph{fortes}.

La proposition suivante indique quand un poids faible est automatiquement fort.

\begin{prop}\label{pr-pforfaib} Soit $w : M\to K_\mu$ un poids. Les assertions suivantes sont équivalentes.
\begin{enumerate}
\item[1.] L'inclusion $\F(\underline{M};K)_w\hookrightarrow\F(\underline{M};K)_{<w>}$ est une égalité.
\item[2.] L'inclusion $\mathfrak{m}_w^2\subset\mathfrak{m}_w$ d'idéaux de $K[M]$ est une égalité.
\item[3.] Toute dérivation non additive\index{termin}{derivation@dérivation!non additive} $M\to K$ relative à $w$ est nulle.
\end{enumerate}

Elles impliquent la propriété suivante :
\begin{enumerate}
\item[4.] l'inclusion $\fct(\C,\E)_w\hookrightarrow\fct(\C,\E)_{<w>}$ est une égalité.
\end{enumerate}
\end{prop}

\begin{proof}
Les implications 1.$\Leftrightarrow$2.$\Rightarrow$4. sont des cas particuliers du corollaire~\ref{cor-uvcf}. L'équivalence entre 2. et 3. résulte de l'isomorphisme de $K$-espaces vectoriels du dual de $\mathfrak{m}_w/\mathfrak{m}_w^2$ vers les dérivations non additives $M\to K$ relatives à $w$ (cf. \eqref{eq-ident_der}, page~\pageref{eq-ident_der}).
\end{proof}

\begin{exem} Si $M$ est un \emph{groupe}, il découle de l'exemple~\ref{ex-deriv_grp} que les conditions précédentes équivalent à \emph{$M$ de torsion} si $p=0$ et à \emph{$M$ $p$-divisible} si $p>0$.
%%%%%%%%%%%%%%%%%%%%%%%%%%%%%%%%%%%%%%%%%%%%%%%%%%%%%%%%%%%%%%%
% UTILISER REMARQUE CI-DESSOUS SI VOLONTÉ ENCYCLOPÉDIQUE !
%On peut aller plus loin et s'affranchir de la condition que $M$ est un groupe en localisant par rapport aux éléments de poids non nul --- voir si cela présente un intérêt de l'écrire. Précisément, si $w(x)\ne 0\Leftrightarrow x\in M^\times$, alors l'absence de $\varphi$-dérivation sur $M$ équivaut à la même chose pour la restriction au groupe des inversibles $M^\times$ et au fait que pour tout élément non inversible $x$, $x$ est le produit de deux éléments non inversibles ou il existe $a\in M^\times$ tel que $ax=x$ et $\varphi(a)\ne 1$.
%%%%%%%%%%%%%%%%%%%%%%%%%%%%%%%%%%%%%%%%%%%%%%%%%%%%%%%%%%%%%%%
\end{exem}

%\begin{rema}\label{rq-recippoiff} La condition 4. de la proposition~\ref{pr-pforfaib} implique les précédentes dans $\F(\C;K)$ (et dans $\fct(\C,\E)$ si $\E$ est non nulle) s'il existe un objet $c$ de $\C$ tel que le morphisme de monoïdes (non commutatif au but) structural $M\to\mathrm{End}_\C(c)$ soit injectif. Cela découle de l'exemple~\ref{ex-pgqwP} et du corollaire~\ref{cor-uvcf}.
%\end{rema}

L'assertion suivante découle de la proposition~\ref{pr-pforfaib} et de la remarque~\ref{rq-derp0}.

\begin{coro}\label{cor-poidsf-carp} Supposons que la caractéristique $p$ du corps $K$ est non nulle et que l'endomorphisme $x\mapsto x^p$ du monoïde $M$ est surjectif. Alors pour tout poids $w : M\to K_\mu$, l'inclusion $\fct(\C,\E)_w\hookrightarrow\fct(\C,\E)_{<w>}$ est une égalité.
\end{coro}

Le résultat suivant montre que, pour un foncteur possédant une décomposition en poids faible, le caractère fort de cette décomposition est une propriété de caractère local par rapport au monoïde $M$.

\begin{prop}\label{pr-dfdfcf} Soit $F$ un foncteur de $\fct(\C,\E)$ possédant une décomposition en poids faible (relativement à $M$). Si $M$ est colimite filtrante de sous-monoïdes $M_i$ tels que $F$ possède une décomposition en poids forte par rapport à l'action de chaque monoïde $M_i$, alors $F$ admet une décomposition en poids forte relativement à $M$.
\end{prop}

\begin{proof}
La $K$-algèbre $\mathfrak{Z}_M(F)$ est semi-primaire déployée, par la proposition~\ref{pr-poideploye}, de sorte qu'il suffit de montrer qu'elle est réduite pour conclure (en utilisant la même proposition). Or $\mathfrak{Z}_M(F)$ est colimite filtrante de ses sous-algèbres $\mathfrak{Z}_{M_i}(F)$, qui sont semi-simples (toujours par la même proposition), donc réduites, d'où la conclusion.
\end{proof}

\begin{rema} On peut aussi démontrer cette proposition avec le point de vue du §\,\ref{sredend}, en combinant la remarque~\ref{rq-diatricol} aux propositions~\ref{pr-tridia-gl} et~\ref{pr-dtf}.
\end{rema}

Nous considérerons pour terminer ce chapitre le cas où $M$ est fini, puis localement fini.

\begin{prop} Supposons que $M$ est {\bf fini} et que le corps $K$ contient assez de racines de l'unité, au sens où tout entier qui est l'ordre du groupe des éléments inversibles du sous-\emph{semi-groupe} de $M$ engendré par un élément de $M$ appartient à $\nu(K)$.
\begin{enumerate}
\item Tout foncteur de $\fct(\C,\E)$ possède une décomposition en poids faible finie. Ainsi, le foncteur canonique
$$\prod_{\pi\in\mon(M,K_\mu)}\fct(\C,\E)_{<\pi>}\to\fct(\C,\E)$$
est une équivalence.
\item Si $M$ est un groupe dont l'ordre est inversible dans $K$, tout foncteur de $\fct(\C,\E)$ possède une décomposition en poids forte finie. Ainsi, le foncteur canonique
$$\prod_{\pi\in\mon(M,K_\mu)}\fct(\C,\E)_\pi\to\fct(\C,\E)$$
est une équivalence.
\end{enumerate}
\end{prop}

\begin{proof} Comme $M$ est fini, la $K$-algèbre $K[M]$ est artinienne, et en particulier semi-primaire. Du fait que $K$ contient assez de racines de l'unité, elle est déployée (cf. exemple~\ref{ex-assez_rac1}). Si de plus $M$ est un groupe d'ordre inversible dans $K$, alors $K[M]$ est semi-simple par le théorème de Maschke. La conclusion provient donc de la proposition~\ref{pr-poideploye}.
\end{proof}

\begin{prop}\label{pr-poidec} Supposons que le corps $K$ contient assez de racines de l'unité et que $M$ est un groupe dont tous les éléments sont d'ordre fini et premier à l'exposant caractéristique de $K$. Alors tout foncteur de $\fct(\C,\E)$ vérifiant la propriété {\rm (Dec)} possède une décomposition en poids forte relativement à l'action de $M$.
\end{prop}

\begin{proof}
Combiner la proposition~\ref{poivN}, le corollaire~\ref{cor-vN} et l'observation que tous les idéaux maximaux de $K[M]$ sont de la forme $\mathfrak{m}_w$ pour un poids $w$ vu que $K$ a assez de racines de l'unité.
\end{proof}

\begin{coro}\label{cor-poidec} Supposons que $k$ est un sous-corps localement fini de $K$ et que $\A$ est une petite catégorie $k$-linéaire. Tout foncteur de $\fct(\A,\E)$ vérifiant la propriété {\rm (Dec)} possède une décomposition en poids forte relativement à l'action de $k_\mu$.
\end{coro}

\begin{proof} Comme $k$ est localement fini et de même caractéristique que $K$, les éléments du groupe multiplicatif $k^\times$ sont d'ordre fini premier à l'exposant caractéristique de $K$. Par ailleurs, comme $K$ contient $k$, $K$ contient assez de racines de l'unité (relativement à $k^\times$). La conclusion découle donc de la proposition~\ref{pr-poidec} et des exemples~\ref{ex-aa0} et~\ref{ex-adj0poi}.
\end{proof}

Les deux énoncés précédents peuvent aussi se déduire de la propriété suivante, qui est très analogue à la proposition~\ref{poivN}.

\begin{prop}\label{pr-col_fil-decpoi} Supposons que $M$ est colimite filtrante de sous-monoïdes $M_i$ et que $F$ est un foncteur de $\fct(\C,\E)$ vérifiant la condition {\rm (Dec)} et possédant une décomposition en poids forte relativement à l'action de chaque monoïde $M_i$. Alors $F$ possède une décomposition en poids forte relativement à l'action de $M$.
\end{prop}

\begin{proof} Il suffit de montrer le résultat lorsque $F$ est indécomposable. Dans ce cas, pour tout $i$, $F$ est homogène d'un certain poids fort $w_i : M_i\to K_\mu$, et les $w_i$ s'assemblent en poids $w : M\to K_\mu$. Alors $F$ est homogène de poids fort $w$, d'où le résultat.
\end{proof}

\chapter{Poids des foncteurs sur une catégorie additive}\label{spfca}

\begin{cvi}
Dans tout le chapitre~\ref{spfca}, $\A$ désigne une catégorie additive $A$-linéaire essentiellement petite. Elle est en particulier munie d'une action du monoïde multiplicatif $A_\mu$ ; toutes les notions de poids considérées dans la suite sont relatives à cette action.
\end{cvi}

\begin{noti} Soit $d\in\mathbb{N}$. On désigne par $\oplus_d : \A^d\to\A$\index{nota}{$\oplus_d$ \emph{ (somme directe $d$-itérée)}} le foncteur de somme directe $d$-itérée, et par $\delta_d : \A\to\A^d$\index{nota}{delta@$\delta$, $\delta_d$ \emph{ (diagonale, diagonale $d$-itérée)}} le foncteur de diagonale $d$-itérée. On note simplement $\delta : \A\to\A\times\A$ pour $\delta_2$.
\end{noti}

Dans toute la suite de ce mémoire, c'est cette situation additive qui nous intéressera ; le présent chapitre donne un certain nombre de conséquences simples des résultats du chapitre~\ref{spcf} dans ce contexte, dans la catégorie $\F(A,K)$ ou, parfois, plus généralement, dans $\fct(\A,\E)$, où $\E$ est une catégorie de Grothendieck $K$-linéaire. Il est aussi l'occasion de rappeler, à la section~\ref{s-rppol}, des définitions et propriétés fondamentales liées à a notion de \emph{foncteur polynomial}, provenant notamment d'Eilenberg-MacLane et Pirashvili, mais aussi de développements plus récents (dus aux auteurs et à Vespa \cite{DTV}) qui introduisent la notion plus générale de \emph{foncteur \ph}, celle de \emph{foncteur antipolynomial} et donnent des \emph{décompositions à la Steinberg} qu'un des objectifs de ce mémoire consiste à préciser.

Une première classe de résultats de ce chapitre consiste à caractériser les situations où \emph{tous} les foncteurs, ou bien tous les foncteurs \emph{polynomiaux}, de $\F(A,K)$, possèdent une décomposition en poids (§\,\ref{par-decPgratos} et~\ref{par-dpfa}). Ils serviront principalement d'échauffement pour la suite, ainsi qu'à illustrer, dans le cas où l'anneau $A$ est infini, à quel point l'existence d'une décomposition en poids est rare sans hypothèse de finitude supplémentaire, et notamment que la condition des valeurs de dimensions finies pour un foncteur de $\F(\A;K)$, qui apparaît au corollaire~\ref{cor-tfdp} et sera cruciale dans de nombreux énoncés de la fin de ce mémoire, n'est nullement surperflue.

La deuxième partie du chapitre est consacrée à l'étude des poids dans notre situation additive, c'est-à-dire à celle du monoïde $\mon(A_\mu,K_\mu)$. Dès lors que le quotient de $A$ par son nilradical est infini, la structure de ce monoïde semble en général hors d'atteinte, et en tout cas difficile à exploiter, par exemple lorsque $A=\mathbb{Z}$ : la compréhension d'une fonction fortement multiplicative à partir de ses valeurs sur les nombres premiers est tout sauf évidente ! C'est pourquoi nous introduisons, à la section~\ref{ppo}, la notion de \emph{poids ordinaire}. Les poids ordinaires forment un sous-monoïde de $\mon(A_\mu,K_\mu)$ dont on peut décrire assez bien la structure. Nous verrons que les poids des foncteurs polynomiaux (ou plus généralement analytiques, ou encore \phs), ainsi que des foncteurs à valeurs de dimensions finies de $\F(\A;K)$, sont toujours ordinaires. Nous verrons qu'il est rare que tous les poids de $\mon(A_\mu,K_\mu)$ soient ordinaires : par exemple, si $A$  est un anneau noethérien, cela survient si et seulement si $A$ est fini (corollaire~\ref{cor-poiord}) ; si $K$ est au plus dénombrable, $\mon(\mathbb{Z}_\mu,K_\mu)$ a la puissance du continu tandis que son sous-monoïde de poids ordinaires est dénombrable.

L'une des propriétés de pure algèbre commutative des poids qui jouera un rôle au chapitre~\ref{shts} est la structure du lieu d'annulation des poids de $\mon(A_\mu,K_\mu)$ : le lieu d'annulation d'un poids \emph{ordinaire} est toujours une réunion \emph{finie} d'idéaux premiers de $A$ (corollaire~\ref{cor-lapord}) --- le cas de $A=\mathbb{Z}$ illustre ici encore que c'est loin de valoir pour un poids quelconque ({\em toutes} les réunions d'idéaux premiers sont des lieux d'annulation de poids).

Enfin, le chapitre se conclut, à la section~\ref{slmar}, par des lemmes d'arithmétique élémentaire sur l'écriture $p$-adique des entiers qui nous serviront pour l'étude des poids (ordinaires) des foncteurs : le lemme~\ref{lm-arit-padic} sert à montrer la proposition~\ref{pr-decfinie} sur la décomposition des poids polynomiaux, propriété de finitude qui consistue l'un des fondements du chapitre~\ref{sec-finpolp}. Quant au lemme~\ref{lmari-sig}, il interviendra au chapitre~\ref{spefcr}.

\section{Décompositions en poids pour le foncteur $P^A$}\label{par-decPgratos}

Les propositions suivantes montrent que l'existence d'une décomposition en poids pour {\em tout} foncteur de $\F(A,K)$ relève de l'exception.

On rappelle que l'ensemble d'entiers $\nu(K)$ est défini page~\pageref{pnu}, et le foncteur projectif $P^A$ page~\pageref{pproj}.

\begin{prop}\label{pr-dpftjs} Les conditions suivantes sont équivalentes.
\begin{enumerate}
\item\label{dpftjs1} Tout foncteur de $\F(A,K)$ possède une décomposition en poids faible finie ;
\item\label{dpftjs2} tout foncteur de $\F(A,K)$ possède une décomposition en poids faible ;
\item\label{dpftjs3} le foncteur $P^A$ de $\F(A,K)$ possède une décomposition en poids faible ;
\item\label{dpftjs4} l'algèbre $K[A_\mu]$ est semi-primaire déployée ;
\item\label{dpftjs4b} l'anneau $A$ est fini et l'algèbre $K[A^\times]$ est déployée ;
\item\label{dpftjs5} $A$ est fini et le groupe multiplicatif $A^\times$ est somme directe de groupes cycliques dont les ordres appartiennent à $\nu(K)$.
\end{enumerate}

Si elles sont vérifiées, alors tout foncteur de $\fct(\A,\E)$, où $\E$ est une catégorie de Grothendieck $K$-linéaire, possède une décomposition en poids faible finie.
\end{prop}

\begin{proof} Les implications \ref{dpftjs1}$\Rightarrow$\ref{dpftjs2}$\Rightarrow$\ref{dpftjs3} sont évidentes.

Comme $P^A$ est un foncteur de type fini, s'il possède une décomposition en poids faible, elle est nécessairement finie. De plus, le morphisme d'algèbres canonique $K[A_\mu]\to\mathrm{End}(P^A)$ est un isomorphisme par le lemme de Yoneda. La proposition~\ref{pr-poideploye} montre donc l'implication \ref{dpftjs3}$\Rightarrow$\ref{dpftjs4}. Comme tout quotient d'une $K$-algèbre semi-primaire déployée est semi-primaire déployé (lemme~\ref{lm-sqalgdepl}), cette même proposition montre que \ref{dpftjs4} entraîne que tout foncteur de $\fct(\A,\E)$ (où $\E$ est une catégorie de Grothendieck $K$-linéaire) possède une décomposition en poids faible finie, et en particulier \ref{dpftjs1}.

Il ne reste donc plus qu'à montrer l'équivalence de \ref{dpftjs4}, \ref{dpftjs4b} et \ref{dpftjs5}.

L'implication \ref{dpftjs5}$\Rightarrow$\ref{dpftjs4b} résulte de l'observation que, pour un groupe abélien fini $G$, l'algèbre $K[G]$ est déployée si et seulement si $G$ est somme directe de groupes cycliques dont les ordres appartiennent à $\nu(K)$ (cf. exemple~\ref{ex-assez_rac1}).

Il suffit d'établir l'implication \ref{dpftjs4b}$\Rightarrow$\ref{dpftjs4} sous l'hypothèse supplémentaire que $A$ est local. En effet, un anneau fini se décompose en produit direct d'anneaux locaux, et la classe des $K$-algèbres semi-primaires déployées est stable par produit tensoriel (lemme~\ref{lm-sqalgdepl}). Il suffit de plus de montrer que $K[A_\mu]$ est déployée, puisque $K[A]$ est de dimension finie, donc semi-primaire, si $A$ est fini. Si $A$ est local et fini, tout élément non inversible de $A$ est nilpotent, de sorte que l'exemple~\ref{ex-assez_rac1} montre que l'algèbre $K[A_\mu]$ est déployée si et seulement si $K[A^\times]$ l'est, d'où \ref{dpftjs4b}$\Rightarrow$\ref{dpftjs4}.

Supposons maintenant l'algèbre $K[A_\mu]$ semi-primaire déployée. Elle est en particulier semi-locale, de sorte que \cite[th.~2.4]{WJ87} implique que $A/\nil(A)$ est fini. Cela implique que $A$ est semi-parfait (proposition~\ref{pr-qpsp}), on peut donc supposer $A$ local. Par ailleurs, l'algèbre de groupe $K[A^\times]$ est un quotient de $K[A_\mu]$, elle est donc semi-primaire, et a fortiori parfaite. On en déduit par un théorème de Woods \cite{Woo} que $A^\times$ est fini. Mais comme $A$ est local, cela implique que $A$ est lui-même fini (si $x\in A$ est non inversible, $1+x$ est inversible). Combiné à la discussion précédente sur le déploiement de $K[A_\mu]$, cela établit l'implication \ref{dpftjs4}$\Rightarrow$\ref{dpftjs5}.
\end{proof}

\begin{rema}
La condition \ref{dpftjs5} de la proposition~\ref{par-decPgratos} peut se reformuler en disant que $A$ est fini et que $K$ est un corps de décomposition non additif de la catégorie $\mathbf{P}(A)$.
\end{rema}

\begin{prop}\label{pr-dpfotjs} Les conditions suivantes sont équivalentes.
\begin{enumerate}
\item Tout foncteur de $\F(A,K)$ possède une décomposition en poids forte finie ;
\item tout foncteur de $\F(A,K)$ possède une décomposition en poids forte ;
\item le foncteur $P^A$ de $\F(A,K)$ possède une décomposition en poids forte ;
\item\label{dpfotjs4} l'algèbre $K[A_\mu]$ est semi-simple déployée ;
\item\label{dpfotjs5} l'anneau $A$ est un produit fini de corps finis dont le cardinal $q$ est tel que $q-1$ appartient à $\nu(K)$ et est inversible dans $K$.
\end{enumerate}

Si elles sont vérifiées, alors tout foncteur de $\fct(\A,\E)$, où $\E$ est une catégorie de Grothendieck $K$-linéaire, possède une décomposition en poids forte finie.
\end{prop}

\begin{proof} La démonstration est entièrement analogue à celle de la proposition~\ref{pr-dpftjs} (et même plus simple pour l'implication \ref{dpfotjs4}$\Rightarrow$\ref{dpfotjs5}, puisqu'une algèbre semi-simple déployée est de dimension finie) ; le seul point qui diffère est le fait suivant : si l'anneau $A$ est fini, alors l'algèbre $K[A_\mu]$ est semi-simple si et seulement si $A$ est un produit fini de corps finis dont le cardinal $q$ est tel que $q-1$ est inversible dans $K$. Il suffit de l'établir lorsque $A$ est local. L'exemple~\ref{ex-assez_rac1}, le caractère cyclique du groupe multiplicatif d'un corps fini et l'observation qu'un élément nilpotent non nul de $A$ n'est pas régulier dans $A_\mu$ permettent de conclure.
\end{proof}

\section{Rappels sur les foncteurs polynomiaux et antipolynomiaux}\label{s-rppol}

\subsection{Foncteurs polynomiaux}

Comme $\A$ possède un objet nul, tout foncteur de $\A$ vers une catégorie abélienne se scinde de façon naturelle et unique en la somme directe d'un foncteur constant (en sa valeur sur l'objet nul $0$ de $\A$) et d'un foncteur \index{termin}{reduit@réduit!\emph{(foncteur)}} réduit (c'est-à-dire nul en $0$). Cette décomposition se généralise à tout multifoncteur en $d$ variables dans $\A$ :  un foncteur $X : \A^d\to\E$, où $\E$ est une catégorie abélienne, se scinde de façon unique et naturelle en somme directe de foncteurs qui par rapport à chaque variable sont soit constants, soit réduits. Le facteur direct réduit par rapport à chaque variable s'appelle partie {\em multiréduite} de $F$. Le $d$-ième \index{termin}{effet croise@effet croisé} {\em effet croisé} d'un foncteur $F$ de $\fct(\A,\E)$ est par définition la partie multiréduite du foncteur $F\circ\oplus_d$ de $\fct(\A^d,\E)$. Cet effet croisé est noté $cr_d(F)$.\index{nota}{cr@$cr_d$ \emph{($d$-ième effet croisé)}}

\begin{rema}\label{rq-dcr}\begin{enumerate}
\item La notion d'effet croisé remonte à Eilenberg-MacLane \cite[chap.~II]{EML}, au moins lorsque la source et le but des foncteurs sont des catégories de modules. Pour davantage de références bibliographiques et de détail sur cette notion classique, on pourra consulter par exemple \cite[§\,2]{DTV}.
\item\label{ppdcr} La décomposition d'un multifoncteur en termes constants ou réduits par rapport à chaque variable peut s'exprimer en termes de poids\,\footnote{Au moins si le but est une catégorie de Grothendieck $K$-linéaire, mais la définition s'étend sans changement au cas d'une catégorie abélienne générale au but.} : elle correspond à l'étude de l'action (variable par variable) du monoïde $\{0,1\}^d$ sur $\A^d$. Dans le cas d'une composition $F\circ\oplus_d$, elle prend la forme
\begin{equation}\label{eq-decCr}
F\circ\oplus_d\simeq\bigoplus_{E\subset\llbracket 1,d\rrbracket}cr_{\cd(E)}(F)\circ\iota_E
\end{equation}
où $\iota_E : \A^d\to\A^{\cd(E)}$ est la composée de la projection $\A^d\simeq\A^{\llbracket 1,d\rrbracket}\to\A^E$ induite par l'inclusion $E\hookrightarrow\llbracket 1,d\rrbracket$ et de l'équivalence $\A^E\simeq\A^{\llbracket 1,\cd(E)\rrbracket}$ obtenue en munissant $E$ de l'ordre induit par celui de $\llbracket 1,d\rrbracket$ (cf. \cite[th.~9.1]{EML}).
\end{enumerate} 
\end{rema}

On obtient ainsi un foncteur continu et cocontinu $cr_d : \fct(\A,\E)\to\fct(\A^d,\E)$ pour tout $d\in\mathbb{N}$. Un foncteur $F$ de $\fct(\A,\E)$ est dit {\em polynomial}\index{termin}{polynomial(e)!foncteur} de degré au plus $d$ si $cr_{d+1}(F)$ est nul. La sous-catégorie pleine des foncteurs polynomiaux de degré au plus $d$ de $\fct(\A,\E)$ est notée \index{nota}{Pol@$\pol_d$ \emph{(catégorie des foncteurs polynomiaux de degré au plus $d$)}} $\pol_d(\A,\E)$ ; elle est bilocalisante. On notera $\pol_d(\A;K)$ pour $\pol_d(\A,K\Md)$ et $\pol_d(A,K)$ pour $\pol_d(\mathbf{P}(A);K)$.

Le foncteur d'inclusion $\pol_d(\A,\E)\to\fct(\A,\E)$ possède un adjoint à droite noté \index{nota}{p@$\ppol_d$ \emph{(plus grand sous-foncteur polynomial de degré au plus $d$)}}$\ppol_d$ et un adjoint à gauche noté \index{nota}{q@$\qpol_d$ \emph{(plus grand quotient polynomial de degré au plus $d$)}} $\qpol_d$. Pour tout foncteur $F$ de $\fct(\A,\E)$, $\ppol_d(F)$ est le plus grand sous-foncteur polynomial de degré au plus $d$ de $F$, et $\qpol_d(F)$ son plus grand quotient polynomial de degré au plus $d$. Un argument formel fondé sur l'adjonction somme/diagonale fournit classiquement des suites exactes naturelles
\begin{equation}\label{eq-qpol}
    \delta_{d+1}^*cr_{d+1}(F)\to F\to\qpol_d(F)\to 0\;;
\end{equation}
\begin{equation}\label{eq-ppol}
    0\to\ppol_d(F)\to F\to\delta_{d+1}^*cr_{d+1}(F)\;.
\end{equation}

L'un des résultats fondamentaux de structure des foncteurs polynomiaux est le suivant, dû à Pirashvili \cite{Pira88}.

\begin{theo}[Pirashvili]\label{th-recPiranv} Soit $\E$ une catégorie de Grothendieck. Le foncteur $cr_d$ induit une équivalence de catégories
$$\pol_d(\A,\E)/\pol_{d-1}(\A,\E)\xrightarrow{\simeq}\mathbf{Add}_d(\A,\E)\rtimes\Si_d\,.$$\index{nota}{Add@$\mathbf{Add}_d$}
\end{theo}

Dans l'énoncé précédent, la notation $\rtimes$ est celle de la définition~\ref{def-psdc} (page~\pageref{def-psdc}) ; le groupe symétrique $\Si_d$ agit sur la catégorie $\mathbf{Add}_d(\A,\E)$ des multifoncteurs multiadditifs $\A^d\to\E$ par permutation des facteurs du produit à la source.

Dans le cas particulier fondamental $\A=\mathbf{P}(A)$, au vu de l'équivalence \eqref{eq-pteAddd} (page~\pageref{eq-pteAddd}) et de la remarque~\ref{rq-wrca} (page~\pageref{rq-wrca}), le théorème précédent prend la forme suivante :
\begin{coro}\label{cor-recPira-gal} Si $\E$ est une catégorie de Grothendieck $K$-linéaire, le foncteur $F\mapsto cr_d(F)(A,\dots,A)$ induit une équivalence
$$\pol_d(\mathbf{P}(A),\E)/\pol_{d-1}(\mathbf{P}(A),\E)\xrightarrow{\simeq}\E\underset{K}{\boxtimes}\Big(A_K\wr\Si_d\Big).$$
\end{coro}
(Pour la signification du symbole $\wr$, voir la notation~\ref{nota-wrt}, page~\pageref{nota-wrt}.)

En particulier, pour $\E=K\Md$, on obtient :
\begin{coro}\label{cor-recPira} Le foncteur $F\mapsto cr_d(F)(A,\dots,A)$ induit une équivalence
$$\pol_d(A,K)/\pol_{d-1}(A,K)\xrightarrow{\simeq}(A_K\wr\Si_d)\Md\,.$$\index{nota}{$\wr$ \emph{(produit en couronne)}}
\end{coro}

Nous aurons besoin du résultat suivant pour établir le corollaire~\ref{cor-carsideg}.
\begin{prop}\label{pr-simpoldfext} Soit $S$ un foncteur polynomial simple de $\F(A,K)$. Il existe une extension de corps $K\subset L$ telle que $S$ soit isomorphe à la composée d'un foncteur polynomial simple de $\F^\df(A,L)$ et du foncteur de restriction des scalaires $L\Md\to K\Md$.
\end{prop}

\begin{proof} Comme $S$ est simple, $L=Z(\mathrm{End}(S))$ est une extension du corps $K$ ; $S$ est évidemment isomorphe à la postcomposition par la restriction des scalaires d'un foncteur $\tilde{S}$, nécessairement simple et polynomial, de $\F(A,L)$. De plus, $\tilde{S}$ a le même corps (non nécessairement commutatif) d'endomorphismes que $S$. La proposition~\ref{pr-endofcentr} (page~\pageref{pr-endofcentr}), \cite[prop.~1.8]{DTV}, le corollaire~\ref{cor-recPira} et le fait que tout objet simple de $A_K\Md$ a un corps d'endomorphismes {\em commutatif} (puisque $A_K$ est un anneau commutatif) montrent que $\mathrm{End}(S)$ est de dimension finie sur $L$. La conclusion découle donc de \cite[prop.~6.3]{DTV}.
\end{proof}

\begin{rema} Ce résultat, qui repose fortement sur la commutativité de $A$ (et de $K$), peut s'exprimer de façon plus intrinsèque par l'\emph{endofinitude centrale} des foncteurs simples polynomiaux de $\F(A,K)$
\end{rema}

\subsection{Foncteurs \phs}

Les effets croisés constituent une version fonctorielle des {\em déviations} d'une fonction entre groupes abéliens --- cf. Eilenberg-MacLane \cite[§\,8]{EML}, qui permet d'introduire la notion de fonction\index{termin}{polynomial(e)!fonction|textbf} polynomiale\,\footnote{Dans \cite{DTV}, cette propriété est nommée EML-polynomiale (pour bien la distinguer des polynômes formels). Dans le présent article, sauf mention du contraire, la polynomialité des fonctions sera toujours cette notion due à Eilenberg et MacLane.} (qui est beaucoup moins rigide que celle de polynôme), et la notion de degré correspondante, entre groupes abéliens. On pourra aussi consulter, par exemple, \cite[chap.~V, §\,1]{Passi} ou \cite[§\,2.2]{DTV} pour les définitions précises et propriétés de base de cette notion classique.

Dans \cite[§\,2.3]{DTV} est introduite la généralisation suivante des foncteurs polynomiaux :
\begin{defi} Un foncteur $F : \A\to\E$, où $\E$ est une catégorie abélienne, est dit \index{termin}{Hom-polynomial}\emph{\ph} si pour tous objets $x$ et $y$ de $\A$, la fonction ensembliste entre groupes abéliens $\A(x,y)\to\E(F(x),F(y))$ qu'induit $F$ est polynomiale.
\end{defi}

Un foncteur polynomial est toujours \ph. Si $\E$ est $\mathbb{Q}$-linéaire, tout foncteur \ph\ de type fini de $\fct(\A,\E)$ est polynomial (cf. proposition~\ref{pr-decpolcar0} ci-après). Un foncteur \ph\ de type fini et de type cofini, ou plus généralement, possédant un support fini et un co-support fini, de $\fct(\A,\E)$ est toujours polynomial \cite[prop.~2.11]{DTV}. En revanche, un foncteur \ph\ de type fini n'est pas nécessairement polynomial \cite[ex.~2.7]{DTV}.

La classe des foncteurs \phs\ de $\fct(\A,\E)$ est stable par sous-quotient et par sommes directes finies. Lorsque le but est une catégorie de modules, elle est stable par produit tensoriel.

L'un des objectifs de ce mémoire consiste à mieux comprendre la structure des foncteurs \phs\ --- voir notamment ci-après le corollaire~\ref{cor-EMLpol-pol} et la section~\ref{phthpol}.

Nous aurons besoin du lien simple suivant entre foncteurs polynomiaux et \phs.

\begin{prop}\label{pr-ppht} Soit $\E$ une catégorie abélienne.
Considérons les propriétés suivantes d'un foncteur $F$ de $\fct(\A,\E)$.
\begin{itemize}
\item[(1)] $F$ est \ph\ ; 
\item[(2a)] pour tout objet $a$ de $\A$, l'épimorphisme $F(a)\twoheadrightarrow\qpol_d(F)(a)$ est un isomorphisme si $d$ est assez grand ;
\item[(2b)] pour tout objet $a$ de $\A$, le monomorphisme $\ppol_d(F)(a)\hookrightarrow F(a)$ est un isomorphisme si $d$ est assez grand.
\end{itemize}
Alors (2a) ou (2b) impliquent chacune (1). Réciproquement, si $F$ est à support (resp. co-support) fini, alors (1) implique (2a) (resp. (2b)).
\end{prop}

\begin{proof} Supposons (2a) vérifié. Soient $a$ et $b$ des objets de $\A$ et $d$ un entier assez grand pour que la projection $F(b)\twoheadrightarrow\qpol_d(F)(b)$ soit un isomorphisme. Considérons le diagramme commutatif
$$\xymatrix{\A(a,b)\ar[r]^-{(i)}\ar[d]_-{(ii)} & \E(F(a),F(b))\ar[d]^\simeq_-{(**)} \\
\E(\qpol_d(F)(a),\qpol_d(F)(b))\ar[r]^-{(*)} & \E(F(a),\qpol_d(F)(b))
}$$
où $(i)$ (resp. $(ii)$) est la fonction ensembliste induite par $F$ (resp. $\qpol_d(F)$) et (*) (resp. (**)) la fonction additive induite par $F(a)\twoheadrightarrow\qpol_d(F)(a)$ (resp. $F(b)\twoheadrightarrow\qpol_d(F)(b)$). Comme $\qpol_d(F)$ appartient à $\pol_d(\A,\E)$, la fonction $(ii)$ est polynomiale de degré au plus $d$ ; il en est de même pour sa composée avec la fonction additive (*). Comme (**) est un isomorphisme de groupes abéliens, il s'ensuit que $(i)$ est également une fonction polynomiale. Ainsi, $F$ est \ph.

L'implication (2b)$\Rightarrow$(1) est duale de la précédente.

Supposons maintenant que $F$ est à support fini et \ph. Soit $a$ un objet de $\A$. Considérons un support fini $\{x_1,\dots,x_n\}$ de $F$ et notons $d$ le maximum des degrés polynomiaux des fonctions $\A(x_i,a)\to\E(F(x_i),F(a))$ induites par $F$. Pour tout objet $t$ de $\A$, la fonction $\A(t,a)\to\E(F(t),F(a))$ induite par $F$ est polynomiale de degré au plus $d$ (ce fait simple se montre exactement comme \cite[lemme~2.12]{DTV}). Prenons $t=a^{\oplus (d+1)}$ : cela implique que le morphisme  $F(a^{\oplus (d+1)})\to F(a)$ somme alternée sur les parties $I$ de $\llbracket 1,d+1\rrbracket$ des $F(\pi_I)$, où $\pi_I : a^{\oplus (d+1)}\to a$ a pour composante $\mathrm{Id}_a$ ou $0_a$ selon que l'indice appartient ou non à $I$, est nul. Comme $\qpol_d(F)(a)$ est le conoyau de ce morphisme (cf. \eqref{eq-qpol}), (2a) est vérifié.

Le fait qu'un foncteur \ph\ à co-support fini vérifie (2b) est dual.
\end{proof}

\subsection{Foncteurs antipolynomiaux}

Suivant \cite[§\,4.1]{DTV}, nous dirons qu'une catégorie additive $\B$ est \emph{$K$-triviale} si pour tous objets $x$ et $y$, le groupe abélien $\B(x,y)$ est fini d'ordre inversible dans $K$. Un foncteur de $\fct(\A,\E)$, où $\E$ est une catégorie abélienne $K$-linéaire, est dit \emph{antipolynomial} s'il se factorise à travers un foncteur additif de $\A$ vers une petite catégorie additive $K$-triviale.

Un foncteur est simultanément \ph\ et antipolynomial si et seulement s'il est constant.

On rappelle également que l'anneau $A$ est dit $K$-trivial si la catégorie $\mathbf{P}(A)$ est $K$-triviale, autrement dit, si $A$ est fini d'ordre inversible dans $K$. Un idéal $I$ de $A$ est dit \emph{$K$-cotrivial}\index{termin}{cotrivial \emph{(idéal)}} si l'anneau $A/I$ est $K$-trivial.

La classe des foncteurs antipolynomiaux de $\fct(\A,\E)$ est stable par sous-quotients et sommes directes finies. Dans $\F(\A;K)$, elle est également stable par produit tensoriel.

\subsection{Décompositions à la Steinberg globales}

La terminologie suivante reprend et étend \cite[déf.~4.3]{DT-ext}.

\begin{defi}\label{def-decSt} Soit $F$ un foncteur de $\fct(\A,\E)$, où $\E$ est une catégorie de Grothendieck $K$-linéaire. On dit que $F$ possède une \emph{décomposition à la Steinberg de type $PAP$ (resp. $HPAP$)} s'il existe un bifoncteur $B : \A\times\A\to\E$ tel que :
\begin{enumerate}
\item $B$ est polynomial (resp. \ph) par rapport à la première variable, i.e. $B(-,x)$ est un foncteur polynomial (resp. \ph) de $\fct(\A,\E)$ pour tout objet $x$ de $\A$ ;
\item $B$ est antipolynomial par rapport à la deuxième variable ;
\item $F$ est isomorphe à la composée de $B$ et du foncteur diagonale $\A\to\A\times\A$.
\end{enumerate}
\end{defi}

\begin{rema}\begin{enumerate}
\item Lorsque $F$ est à support fini, cela implique l'existence d'une petite catégorie additive $K$-triviale $\B$ telle que $B$ se factorise à travers $\A\times\Phi$ pour un foncteur additif $\Phi : \A\to\B$ (cf. \cite[lemme~4.5]{DT-ext}).
\item Une décomposition à la Steinberg d'un des types précédents est unique à isomorphisme près, grâce à \cite[prop.~4.9]{DTV}.
\end{enumerate}
\end{rema}

L'un des résultats principaux de \cite{DTV} est le suivant :

\begin{theo}\cite[th.~4.10]{DTV}\label{th-DTVglob} Tout foncteur de type fini de $\F^\df(\A;K)$ possède une décomposition à la Steinberg de type $HPAP$.
\end{theo}

\begin{coro}\cite[cor.~4.11]{DTV}\label{cor-DTVglob} Tout foncteur de longueur finie de $\F^\df(\A;K)$ possède une décomposition à la Steinberg de type $PAP$.
\end{coro}

(On a en fait besoin seulement de l'hypothèse légèrement plus faible que le foncteur soit de type fini et de type co-fini.)

L'un des objectifs du présent travail est de renforcer les résultats précédents, dans plusieurs directions : soit en remplaçant l'hypothèse de valeurs de dimension finie par une hypothèse a priori plus faible, soit en obtenant une décomposition de type $PAP$ sous des hypothèses plus faibles que celle du corollaire~\ref{cor-DTVglob}, soit en précisant la décomposition de type $HPAP$. Les résultats les plus significatifs seront obtenus dans la catégorie $\F(A,K)$ ; nous introduisons à cette fin la notion suivante.

\begin{defi}\label{def-decFP} On \emph{décomposition de type $FP$} d'un foncteur $F$ de $\F(A,K)$ tout isomorphisme de la forme $F\simeq\alpha^*B$, où $I$ est un idéal cofini de $A$, $B : \mathbf{P}(A)\times\mathbf{P}(A/I)\to K\Md$ un bifoncteur polynomial par rapport à la première variable et $\alpha : \mathbf{P}(A)\to\mathbf{P}(A)\times\mathbf{P}(A/I)$ le foncteur canonique.
\end{defi}

(Les lettres $FP$ abrègent \emph{Fini} et \emph{Polynomial}.)

\section{Décompositions en poids des foncteurs analytiques}\label{par-dpfa}

On rappelle qu'un foncteur est dit \index{termin}{analytique \emph{(foncteur)}}{\em analytique} s'il est colimite (filtrante) de ses sous-foncteurs polynomiaux.

La théorie des foncteurs analytiques (ou \phs) est beaucoup plus simple lorsque $K$ est de caractéristique nulle, en raison du résultat classique suivant :
\begin{prop}\label{pr-decpolcar0} Supposons que $K$ est de caractéristique nulle et que $\E$ est une catégorie de Grothendieck $K$-linéaire. Soit $d\in\mathbb{N}$.
\begin{enumerate}
\item Si $F$ est un foncteur analytique (resp. polynomial de degré au plus $d$) de $\fct(\A,\E)$, alors $F$ possède une décomposition en poids forte relativement à l'action de $\mathbb{Z}_\mu$, et l'ensemble de ses poids est contenu dans celui des $\mathbb{Z}\to\mathbb{Q}\quad\lambda\mapsto\lambda^n$ pour $n\in\mathbb{N}$ (resp. $0\le n\le d$), d'où une décomposition naturelle $F\simeq\bigoplus_{n\in\mathbb{N}}F_{[n]}$ avec $F_{[n]}$ homogène de poids $\lambda\mapsto\lambda^n$. 
\item Tout foncteur \ph\ de $\fct(\A,\E)$ est analytique.
\item\label{itdQrec} Le foncteur $\pol_d(\A,\E)\to\pol_d(\A,\E)/\pol_{d-1}(\A,\E)\times\pol_{d-1}(\A,\E)$ dont les composantes sont le foncteur canonique $\pol_d(\A,\E)\to\pol_d(\A,\E)/\pol_{d-1}(\A,\E)$ et la restriction à $\pol_d(\A,\E)$ de $\qpol_{d-1}$ est une équivalence.

De plus, l'endofoncteur idempotent de $\pol_d(\A,\E)$ composé du foncteur canonique $\pol_d(\A,\E)\to\pol_d(\A,\E)/\pol_{d-1}(\A,\E)$ et de son adjoint à gauche est isomorphe à $F\mapsto F_{[d]}$.
\item Le foncteur canonique $\A\to\A_\mathbb{Q}:=\A\otimes_\mathbb{Z}\mathbb{Q}$ induit une équivalence de catégories $\pol_d(\A_\mathbb{Q},\E)\xrightarrow{\simeq}\pol_d(\A,\E)$
\end{enumerate}
\end{prop}

\begin{proof} On pourra trouver une démonstration des deux premières assertions (celle de la seconde consistant à observer que la décomposition en poids de la première vaut aussi pour les foncteurs \phs) dans \cite[prop.~2.8]{DTV} --- le résultat y est fourni seulement lorsque le but de la catégorie de foncteurs est une catégorie de modules sur une $\mathbb{Q}$-algèbre, mais la démonstration fonctionne pareillement pour toute catégorie abélienne $\mathbb{Q}$-linéaire au but.

La troisième assertion se déduit de la première et du théorème~\ref{th-recPiranv} de Pirashvili, dans lequel il est facile de lire la décomposition en poids du premier point. (Pour montrer l'équivalence $\pol_d(\A,\E)\xrightarrow{\simeq}\pol_d(\A,\E)/\pol_{d-1}(\A,\E)\times\pol_{d-1}(\A,\E)$, on peut aussi raisonner directement à partir du fait que la norme est un isomorphisme entre $\mathbb{Q}[\Si_d]$-modules.)

La dernière assertion se déduit de la première et de la proposition~\ref{pr-inversion_fct}.
\end{proof}

\begin{rema} Lorsque $K$ est de caractéristique $p>0$ et que $\E$ est une catégorie de Grothendieck $K$-linéaire, tout foncteur analytique de $\fct(\A,\E)$ possède une décomposition en poids \emph{faible} relativement à l'action de $\mathbb{Z}_\mu$, mais pas nécessairement une décomposition en poids forte (cf. propositions~\ref{pr-decpfa} et~\ref{pr-decpfoap} ci-après). L'existence d'une décomposition en poids faible peut se voir directement à partir de l'observation élémentaire et bien connue (cf. corollaire~\ref{cor-radPolp} ci-après) que le foncteur canonique $\A\twoheadrightarrow\A/p^d$ induit une équivalence $\pol_d(\A/p^d,\E)\xrightarrow{\simeq}\pol_d(\A,\E)$ lorsque $\E$ est une catégorie abélienne $\mathbb{Z}/p$-linéaire. Ce résultat s'étend aux foncteurs \phs\ (cf. remarque~\ref{rq-dpfg3} ci-après).
\end{rema}

\begin{rema}\label{rq-corNb} On déduit de la proposition~\ref{pr-decpolcar0} que, si $A$ est un corps de nombres, alors tout foncteur polynomial de $\F(A,K)$ possède une décomposition en poids forte si $K$ est un corps de décomposition de $\mathbf{P}(A)$, et que tout foncteur polynomial de type fini de $\F(A,K)$ est à valeurs de dimensions finies.

On en déduit en appliquant la proposition~\ref{pr-col_fil-decpoi} que si $A$ est une extension algébrique de $\mathbb{Q}$ et $K$ un corps de décomposition de $\mathbf{P}(A)$, alors tout foncteur polynomial de $\F(A,K)$ vérifiant la condition {\rm (Dec)} possède une décomposition en poids forte.
\end{rema}

Le résultat suivant améliore le corollaire~\ref{cor-tfdp2} pour les foncteurs \phs.

\begin{prop}\label{pr-vdfdppol} Supposons que $K$ est un corps de décomposition de $\A$. Alors tous les foncteurs analytiques et tous les foncteurs \phs\ de $\F^\df(\A;K)$ possèdent une décomposition en poids faible.
\end{prop}

\begin{proof}
Comme la classe des foncteurs admettant une décomposition en poids faible est stable par colimite, il suffit de montrer qu'un foncteur \ph\ \textit{de type fini} $F$ de $\F^\df(\A;K)$ possède une décomposition en poids faible (un foncteur analytique de type fini étant polynomial).

Commençons par le cas où $F$ est \textit{polynomial} : $F$ est alors fini \cite[lemme~11.10]{DTV}, il suffit donc de traiter le cas où $F$ est simple, qui est contenu dans \cite[th.~5.5]{DTV}.

Dans le cas général où $F$ est seulement de type fini, soit $E$ un support fini de ce foncteur. Il existe $d\in\mathbb{N}$ tel que la projection $F(x)\twoheadrightarrow\qpol_d(F)(x)$ soit un isomorphisme pour tout $x\in E$, par la proposition~\ref{pr-ppht}. On en déduit que le morphisme d'algèbres (non commutatives) $\mathrm{End}(F)\to\mathrm{End}(\qpol_d(F))$ est injectif. Il s'ensuit que $\mathfrak{Z}(F)$ est une sous-algèbre de $\mathfrak{Z}(\qpol_d(F))$. La conclusion découle donc du cas polynomial déjà traité, de la proposition~\ref{pr-poideploye} et du lemme~\ref{lm-sqalgdepl}.
\end{proof}

Dans la suite du §\,\ref{par-dpfa}, on se limite par commodité à la catégorie $\F(A,K)$. On peut facilement généraliser la plupart des résultats à $\fct(\A,\E)$, par exemple en utilisant la notion de spectre d'un foncteur introduite à la définition~\ref{defisp} ci-après (ainsi, la proposition~\ref{pr-decpfa} donnera lieu à la proposition~\ref{pr-ttphtr}).

Le résultat suivant est à comparer à la proposition~\ref{pr-dpftjs}.

\begin{prop}\label{pr-decpfa} Les assertions suivantes sont équivalentes.
\begin{enumerate}
\item\label{itdpfa1} La $K$-algèbre $A_K$ est semi-primaire déployée ;
\item\label{itdpfa2} si $k$ est le sous-corps premier de $K$, la $k$-algèbre $A_k=A\otimes_\mathbb{Z}k$ est semi-primaire et le quotient par son radical est un produit de sous-extensions d'extensions finies galoisiennes de $k$ admettant un plongement dans $K$ ;
\item\label{itdpfa3} tout foncteur analytique de $\F(A,K)$ admet une décomposition en poids faible ;
\item\label{itdpfa4} tout foncteur polynomial de $\F(A,K)$ admet une décomposition en poids faible finie ;
\item\label{itdpfa5} tout foncteur additif de $\F(A,K)$ admet une décomposition en poids faible ;
\item\label{itdpfa6} tout foncteur additif de $\F(A,K)$ admet une décomposition en poids faible finie.
\end{enumerate}
\end{prop}

\begin{proof} De manière générale, si $k\to K$ est une extension de corps séparable et $B$ une $k$-algèbre, il découle de \cite[§\,12, n°4, pr.~3]{Bki} que la $K$-algèbre $B\otimes_k K$ est semi-primaire déployée si et seulement si $B$ est semi-primaire et que $B/\rj(B)$ est un produit de corps $L$ contenant $k$ tels que $L\otimes_k K$ soit isomorphe à $K^n$ pour un entier $n>0$. Cette condition équivaut à dire que $k\to L$ est une extension finie et qu'elle est isomorphe à une sous-extension d'une extension galoisienne finie de $k$ contenue dans $K$. On en déduit l'équivalence entre \ref{itdpfa1} et \ref{itdpfa2}.

L'équivalence de \ref{itdpfa1}, \ref{itdpfa5} et \ref{itdpfa6} se déduit de l'équivalence entre la sous-catégorie pleine des foncteurs additifs de $\F(A,K)$ et $A_K\Md$.

Montrons l'implication \ref{itdpfa6}$\Rightarrow$\ref{itdpfa4}. On suppose donc \ref{itdpfa6} vérifiée et l'on montre par récurrence sur $d\in\mathbb{N}$ que tout foncteur $F$ de $\pol_d(A,K)$ possède une décomposition en poids faible finie. Pour $d=0$, la propriété est immédiate. On suppose donc $d>0$ et l'assertion établie pour les foncteurs de degré au plus $d-1$. Notons $T^d_K$ le foncteur $V\mapsto (V\otimes_\mathbb{Z}K)^{\otimes d}$ de $\F(A,K)$ et posons $G=T^d_K\underset{A_K\wr\Si_d}{\otimes}cr_d(F)(A,\dots,A)$. On dispose d'un morphisme d'adjonction $G\to F$ dont le noyau et le conoyau appartiennent à $\pol_{d-1}(A,K)$ (cf. corollaire~\ref{cor-recPira}). Ceux-ci possèdent donc une décomposition en poids faible finie, par l'hypothèse de récurrence. Par ailleurs, $T^d_K$ possède une décomposition en poids faible finie car c'est le produit tensoriel de $d$ copies d'un foncteur additif, qui admet donc une décomposition en poids faible finie (hypothèse \ref{itdpfa6}), ce qui entraîne la même propriété pour $G$. Il s'ensuit que $F$ possède lui-même une décomposition en poids faible finie.

Comme les implications \ref{itdpfa4}$\Rightarrow$\ref{itdpfa3}$\Rightarrow$\ref{itdpfa5} sont triviales (puisque la classe des foncteurs ayant une décomposition en poids faible est stable par colimite), la proposition est établie.
\end{proof}

\begin{rema}\label{rq-dpfg3} Nous verrons à la proposition~\ref{pr-hpol_htr} que, sous les hypothèses de la proposition~\ref{pr-decpfa}, tous les  foncteurs \phs\ de $\F(A,K)$ possèdent une décomposition en poids faible.
\end{rema}

On rappelle que $p$ désigne la caractéristique du corps $K$.

\begin{prop}\label{pr-decpfoa0} Si $p=0$, les assertions suivantes sont équivalentes.
\begin{enumerate}
\item La $K$-algèbre $A_K$ est semi-simple déployée ;
\item la $\mathbb{Q}$-algèbre $A_\mathbb{Q}$ est un produit fini de sous-extensions d'extensions finies galoisiennes de $\mathbb{Q}$ admettant un plongement dans $K$ ;
\item tout foncteur analytique de $\F(A,K)$ admet une décomposition en poids forte ;
\item tout foncteur polynomial de $\F(A,K)$ admet une décomposition en poids forte finie ;
\item tout foncteur additif de $\F(A,K)$ admet une décomposition en poids forte ;
\item tout foncteur additif de $\F(A,K)$ admet une décomposition en poids forte finie.
\end{enumerate}
\end{prop}

\begin{proof}
Cette proposition s'établit de façon analogue à la proposition~\ref{pr-decpfa} (et parfois plus simplement), en utilisant la proposition~\ref{pr-decpolcar0}.
\end{proof}

Dans le lemme et la proposition qui suivent, si $V$ est un groupe abélien et $d\ge 0$ un entier, on note $\bar{Q}_d(V)$ le quotient de l'idéal d'augmentation de la $K$-algèbre $K[V]$ par sa puissance $(d+1)$-ième.

\begin{lemm}\label{lm-pdeux} Si $p>0$ et que $V$ est un groupe abélien annulé par $p^2$ et tel que la multiplication par $p$ sur $V$ induise $0$ sur $\bar{Q}_p(V)$, alors $V$ est annulé par $p$.
\end{lemm}

\begin{proof} Si $V$ n'est pas annulé par $p$, alors $V$ contient un facteur direct $\mathbb{Z}/p^2$. Il suffit donc de vérifier que la multiplication par $p$ sur $\mathbb{Z}/p^2$ n'induit pas l'endomorphisme nul de $\bar{Q}_p(\mathbb{Z}/p^2)$. En effet, la $K$-algèbre $K[\mathbb{Z}/p^2]$ est isomorphe à $K[x]/(x^{p^2})$, et $\bar{Q}_p(\mathbb{Z}/p^2)$ à son sous-quotient $(x)/(x^{p+1})$ ; à travers cet isomorphisme, la multiplication par $p$ sur $\mathbb{Z}/p^2$ envoie $x$ sur $x^p$, qui n'y est pas nul.
\end{proof}

\begin{prop}\label{pr-decpfoap} Si $p>0$, les assertions suivantes sont équivalentes.
\begin{enumerate}
\item\label{dpfg1} L'anneau $A$ est isomorphe au produit d'une $\mathbb{Z}[1/p]$-algèbre et d'un nombre fini de corps finis admettant un plongement dans $K$ ;
\item\label{dpfg2} tout foncteur analytique de $\F(A,K)$ admet une décomposition en poids forte ;
\item\label{dpfg3} tout foncteur polynomial de $\F(A,K)$ admet une décomposition en poids forte finie ;
\item\label{dpfg4} tout foncteur polynomial de degré au plus $p$ de $\F(A,K)$ admet une décomposition en poids forte ;
\item\label{dpfg5} tout foncteur polynomial de degré au plus $p$ de $\F(A,K)$ admet une décomposition en poids forte finie.
\end{enumerate}
\end{prop}

\begin{proof}
Comme pour les propositions précédentes, on vérifie aisément que tout foncteur additif de $\F(A,K)$ possède une décomposition en poids forte si et seulement si la $\mathbb{Z}/p$-algèbre $A/p$ est isomorphe à un produit fini
de corps finis admettant un plongement dans $K$. 

Supposons maintenant que le foncteur $V\mapsto\bar{Q}_p(V/p^2)$ de $\pol_d(A,K)$ admette une décomposition en poids forte. Comme il est réduit et que $K$ est de caractéristique $p$, cela entraîne que la multiplication par $p$ à la source induit $0$ en l'appliquant. En effet, tout poids polynomial $A\to K$ non constant en $1$ envoie nécessairement $p$ sur $0$. Le lemme~\ref{lm-pdeux} montre alors que $p$ est nul dans $A/(p^2)$, i.e. que l'idéal $(p)$ de $A$ est idempotent, ce qui entraîne que l'anneau $A$ est isomorphe au produit d'une $\mathbb{Z}[1/p]$-algèbre et d'une $\mathbb{Z}/p$-algèbre.

En rassemblant les observations précédentes, on voit que \ref{dpfg4} implique \ref{dpfg1}.

Si l'on suppose maintenant \ref{dpfg1} vérifié, écrivons $A\simeq B\times k$ où $B$ est une $\mathbb{Z}[1/p]$-algèbre et $k$ un produit fini de corps finis admettant un plongement dans $K$. Alors le foncteur canonique $\F(k,K)\to\F(A,K)$ induit par la projection de $A$ sur $k$ à la source induit une équivalence entre les sous-catégories de foncteurs polynomiaux, car une fonction polynomiale d'un $\mathbb{Z}[1/p]$-module vers un $\mathbb{Z}/p$-espace vectoriel est nécessairement constante. Comme tout foncteur de $\F(k,K)$ possède une décomposition en poids forte finie, par la proposition~\ref{pr-dpfotjs}, on en déduit l'implication \ref{dpfg1}$\Rightarrow$\ref{dpfg3}.

Les implications \ref{dpfg3}$\Rightarrow$\ref{dpfg2}$\Rightarrow$\ref{dpfg4} et \ref{dpfg3}$\Rightarrow$\ref{dpfg5}$\Rightarrow$\ref{dpfg4} étant claires, la proposition est démontrée.
\end{proof}

\section{Poids ordinaires}\label{ppo}

Si $w\in\mon(A_\mu,K_\mu)$ est un poids, notons $\mathrm{S}_w$ le foncteur de $\F(A,K)$ prolongement intermédiaire de la représentation $K$-linéaire $K_w$ de dimension $1$ de $A_\mu$ associée à $w$ --- autrement dit, $\mathrm{S}_w$ est l'unique foncteur simple (à isomorphisme près) tel que $\mathrm{S}_w(A)\simeq K_w$ (cf. par exemple \cite[§\,1.3]{DTV}). Explicitement, $\mathrm{S}_w(V)$ est l'image de l'application linéaire
$$K[V]\to K^{V^*}\quad [v]\mapsto\big(l\mapsto w(l(v))\big)$$
où $V^*:=\mathrm{Hom}_A(V,A)$.

Ainsi, $\mathrm{S}_w$ est homogène de poids fort $w$.

\begin{defi}\label{df-poiord} On dit que le poids $w$ est :
\begin{enumerate}
\item \index{termin}{poids!ordinaire}\textbf{ordinaire} si le foncteur $\mathrm{S}_w$ est à valeurs de dimensions finies ;
\item \index{termin}{poids!polynomial}\textbf{polynomial} si le foncteur $\mathrm{S}_w$ est polynomial ;
\item \index{termin}{poids!antipolynomial}\textbf{antipolynomial} si le foncteur $\mathrm{S}_w$ est antipolynomial.
\end{enumerate}
\end{defi}

\begin{prop} Soit $w\in\mon(A_\mu,K_\mu)$.
\begin{enumerate}
\item Le poids $w$ est polynomial si et seulement s'il définit une fonction polynomiale entre les groupes additifs sous-jacents à $A$ et $K$.
\item Le poids $w$ est antipolynomial si et seulement s'il se factorise à travers la réduction de $A$ modulo un idéal $K$-cotrivial.
\end{enumerate}
\end{prop}

\begin{proof} Si $w$ est une fonction polynomiale, \cite[lemme~2.12]{DTV} montre que $\mathrm{S}_w$ est un foncteur polynomial. La réciproque est évidente.

L'assertion relative à la propriété antipolynomiale est immédiate.
\end{proof}

\begin{prop}\label{pr-paphp} Soient $\E$ une catégorie de Grothendieck $K$-linéaire et $F$ un foncteur de $\fct(\A,\E)$.
\begin{enumerate}
\item Si $F$ est \ph\ ou analytique, alors tout poids de $F$ est polynomial.
\item Si $F$ est antipolynomial, alors tout poids de $F$ est antipolynomial.
\end{enumerate}
\end{prop}

\begin{proof} Quitte à remplacer $F$ par son sous-foncteur $F^w$, on peut supposer que $F$ est un foncteur homogène de poids fort $w$. Soit $a$ un objet de $\A$ tel que $F(a)\ne 0$.

Si $F$ est \ph, la fonction $\mathrm{End}_\A(a)\to\mathrm{End}_\E(F(a))$ est polynomiale, donc sa composée avec le morphisme canonique $A\to\mathrm{End}_\A(a)$, qui est additif, l'est également. Cette composée coïncide avec la composée de $w$ et du morphisme canonique $K\to\mathrm{End}_\E(F(a))$, qui est additif et injectif. Par conséquent, $w$ est polynomial.

Tout poids d'un foncteur analytique est poids d'un sous-foncteur polynomial, et donc polynomial par ce qui précède.

Pour montrer l'assertion relative aux foncteurs antipolynomiaux, on peut supposer que $\A$ est $K$-triviale. Le $A$-module $\mathrm{End}_\A(a)$ a alors un ensemble sous-jacent fini d'ordre inversible dans $K$, ce qui implique l'existence d'un idéal $K$-cotrivial $I$ de $A$ qui l'annule. Par conséquent, si $t$ et $u$ sont des éléments de $A$ tels que $t-u\in I$, l'action de $w(t)-w(u)\in K$ sur $F(a)$ est nulle. Comme $K$ est un corps et que $F(a)$ est non nul, il s'ensuit que $w(t)-w(u)$ est nul. Ainsi, $w$ se factorise à travers la projection $A\twoheadrightarrow A/I$, ce qui achève la démonstration.
\end{proof}

\begin{prop}\label{pr-decpord} Notons \index{nota}{Mon@$\mon_{\operatorname{ord}}$, $\mon_{\operatorname{pol}}$, $\mon_{\operatorname{AP}}$} $\mon_{\operatorname{ord}}(A,K)$ (resp. $\mon_{\operatorname{pol}}(A,K)$, $\mon_{\operatorname{AP}}(A,K)$) le sous-ensemble de $\mon(A_\mu,K_\mu)$ constitué des poids ordinaires (resp. polynomiaux, antipolynomiaux).
\begin{enumerate}
\item Les sous-ensembles $\mon_{\operatorname{ord}}(A,K)$, $\mon_{\operatorname{pol}}(A,K)$ et $\mon_{\operatorname{AP}}(A,K)$ sont des sous-monoïdes de $\mon(A_\mu,K_\mu)$.
\item Le morphisme de monoïdes canonique $\mon_{\operatorname{pol}}(A,K)\times\mon_{\operatorname{AP}}(A,K)\to\mon(A_\mu,K_\mu)$ induit un isomorphisme
$$\mon_{\operatorname{pol}}(A,K)\times\mon_{\operatorname{AP}}(A,K)\xrightarrow{\simeq}\mon_{\operatorname{ord}}(A,K)\,.$$
\end{enumerate}
\end{prop}

\begin{proof} La première assertion résulte de ce que $\mathrm{S}_1\simeq K$ et de ce que $\mathrm{S}_{w.w'}$ est un sous-quotient de $\mathrm{S}_w\otimes\mathrm{S}_{w'}$ pour tous poids $w$ et $w'$.

L'inclusion $\mon_{\operatorname{AP}}(A,K)\subset\mon_{\operatorname{ord}}(A,K)$ est évidente ; $\mon_{\operatorname{pol}}(A,K)\subset\mon_{\operatorname{ord}}(A,K)$ découle de \cite[prop.~6.3]{DTV}.

L'injectivité du morphisme de monoïdes induit $\mon_{\operatorname{pol}}(A,K)\times\mon_{\operatorname{AP}}(A,K)\to\mon_{\operatorname{ord}}(A,K)$ se déduit de \cite[prop.~4.9]{DTV} et sa surjectivité du corollaire~\ref{cor-DTVglob}.
\end{proof}

L'importance des poids ordinaires provient notamment du résultat suivant : 

\begin{prop}\label{pr-dfpord} Soit $F$ un foncteur de $\F^\df(\A;K)$. Tout poids $w$ de $F$ est ordinaire.
\end{prop}

\begin{proof} Quitte à remplacer $F$ par un sous-foncteur approprié, on peut supposer que ce foncteur est de type fini et homogène de poids fort $w$. La conclusion résulte donc du théorème~\ref{th-DTVglob} et des propositions~\ref{pr-paphp} et~\ref{pr-decpord}.
\end{proof}

\section{Lieu d'annulation des poids}\label{ssct-la}

\begin{nota} Pour $w\in\mon(A_\mu,K_\mu)$, on note \index{nota}{z@$\la$ \emph{(lieu d'annulation d'un poids)}}$\la(w)$ l'ensemble des éléments $a$ de $A$ tels que $w(a)=0$.
\end{nota}

Le lieu d'annulation $\la(w)$ d'un poids $w$ constitue un outil important, notamment pour les poids ordinaires, que nous étudierons dans les sections suivantes. La considération de lieux d'annulation de poids interviendra au §\,\ref{ss-tshtr} pour démontrer l'un des résultats de structure principaux de ce mémoire.

Les assertions suivantes, immédiates, résument les propriétés formelles qui nous serviront pour manipuler les lieux d'annulation de poids.

\begin{prop}\label{pr-laev} Soient $w$ et $w'$ des éléments de $\mon(A_\mu,K_\mu)$.
\begin{enumerate}
 \item $\la(w.w')=\la(w)\cup\la(w')$ ;
 \item les assertions suivantes sont équivalentes :
 \begin{enumerate}
 \item $\la(w)=\varnothing$ ;
 \item $0\notin\la(w)$ ;
 \item $w=1$.
 \end{enumerate}
 \end{enumerate}
\end{prop}

La notation (LA) introduite ci-dessous abrège \textit{lieu d'annulation}.

\begin{prdef}\label{prdfla} On dit qu'un sous-ensemble $E$ de $A$ vérifie la condition $\mathrm{(LA)}$ s'il satisfait aux trois propriétés suivantes :
\begin{enumerate}
 \item[$\mathrm{(LA1)}$] $1\notin E$ ;
 \item[$\mathrm{(LA2)}$] $\forall x\in E\quad\forall t\in A\qquad x.t\in E$ ;
 \item[$\mathrm{(LA3)}$] $\forall (x,y)\in A^2\quad (x.y\in E)\Rightarrow (x\in E\text{ ou }y\in E)$.
 \end{enumerate}
 
 Alors, pour tout $w\in\mon(A_\mu,K_\mu)$ :
 \begin{enumerate}
 \item $\la(w)$ vérifie $\mathrm{(LA)}$ ;
 \item  si $E$ est une partie de $A$ vérifiant $\mathrm{(LA)}$, la fonction $\xi_E : A\to K$ égale à $0$ sur $E$ et à $1$ ailleurs est un élément de $\mon(A_\mu,K_\mu)$ ;
 \item  tout idéal premier de $A$ vérifie $\mathrm{(LA)}$ ;
 \item  l'ensemble des parties de $A$ vérifiant $\mathrm{(LA)}$ est stable par réunion arbitraire ;
 \item le complémentaire d'une partie de $A$ vérifiant $\mathrm{(LA)}$ est une partie multiplicative de $A$ ;
 \item\label{itlatil} si $S=A\setminus\la(w)$, ou plus généralement, si $S$ est une partie multiplicative de $A$ incluse dans $A\setminus\la(w)$, alors il existe un unique $\tilde{w}\in\mon(A[S^{-1}]_\mu,K_\mu)$ tel que $w$ soit la composée de $\tilde{w}$ et du morphisme d'anneaux canonique $A\to A[S^{-1}]$.
 \end{enumerate}
\end{prdef}
 
On déduit par ailleurs de la proposition~\ref{pr-inversion_fct} l'assertion suivante :

\begin{prop}\label{pr-invad} Soit $F$ un foncteur de $\fct(\A,\E)$ possédant une décomposition en poids forte (resp. faible). On pose $$\tilde{A}_F:=A\Big[\Big(\bigcap_{w\in\Pi(F)}(A\setminus\la(w))\Big)^{-1}\Big].$$\index{nota}{A@$\tilde{A}_F$|textbf}

 Il existe un foncteur $\tilde{F}$ de $\fct(\A\otimes_A\tilde{A}_F,\E)$, unique à isomorphisme près, tel que $F$ soit isomorphe à la composée de $\tilde{F}$ et du foncteur canonique $A\to\A\otimes_A\tilde{A}_F$. De plus, $\tilde{F}$\index{nota}{F@$\tilde{F}$} possède une décomposition en poids forte (resp. faible) et $\Pi_{\tilde{A}_F}(\tilde{F})=\{\tilde{w}\,|\,w\in\Pi_A(F)\}$, avec la notation de la proposition~\ref{prdfla}.\ref{itlatil}.
\end{prop}
 
\section{Poids antipolynomiaux}

Le sous-monoïde $\mon_{\operatorname{AP}}(A,K)$ de $\mon(A_\mu,K_\mu)$ est la colimite filtrante des sous-monoïdes $\mon((A/I)_\mu,K_\mu)$, où $I$ parcourt l'ensemble des idéaux $K$-cotriviaux de $A$.

Commençons donc par étudier le monoïde $\mon(A_\mu,K_\mu)$ lorsque l'anneau $A$ est $K$-trivial. Un tel anneau est en particulier semi-parfait : c'est un produit fini d'anneaux locaux, disons $A\simeq A_1\times\dots\times A_n$. Cela fournit un isomorphisme de monoïdes $\mon(A_\mu,K_\mu)\simeq\mon((A_1)_\mu,K_\mu)\times\dots\times\mon((A_n)_\mu,K_\mu)$. Autrement dit, on se ramène au cas où $A$ est local. Comme les $A_i$ sont quasi-parfaits, on dispose d'isomorphismes canoniques $\mon((A_i)_\mu,K_\mu)\simeq\mathbf{Ab}(A_i^\times,K^\times)_+$ (exemple~\ref{exev-mormon}).

Si l'on revient au cas où $A$ est un anneau quelconque, il résulte de la discussion précédente que :

\begin{prop}\label{pr-zerantipol} Pour $w\in\mon_{\operatorname{AP}}(A_\mu,K_\mu)$, $\la(w)$ est une réunion finie d'idéaux maximaux $K$-cotriviaux de $A$. Si $w$ est non trivial, la réunion est non vide.
\end{prop}

La proposition~\ref{pr-zerantipol} va nous aider à caractériser les anneaux $A$ dont tous les poids sont antipolynomiaux. Commençons par un énoncé facile :

\begin{prop}\label{prpoiantipol_fac} Si $A$ possède un idéal $K$-cotrivial de carré nul, alors tous les poids de $\mon(A_\mu,K_\mu)$ sont antipolynomiaux.
\end{prop}

\begin{proof} Un tel anneau $A$ est semi-primaire ; comme ci-dessus, on peut se ramener au cas où $A$ est local, et l'on a alors $\mon(A_\mu,K_\mu)\simeq\mathbf{Ab}(A^\times,K^\times)_+$. Il s'agit donc de montrer que, si $\varphi : A^\times\to K^\times$ est un morphisme de groupes, il existe un idéal $K$-cotrivial $I$ de $A$ tel que $\varphi$ se factorise à travers le morphisme canonique $A^\times\to (A/I)^\times$.

Soit $\mathfrak{r}$ un idéal $K$-cotrivial de $A$ tel que $\mathfrak{r}^2=0$. Soient $N:=\{x\in\mathfrak{r}\,|\,\varphi(1+x)=1\}$ et $a_1,\dots,a_n$ des relevés dans $A$ des éléments du groupe multiplicatif $(A/\mathfrak{r})^\times$ ; posons
$$I:=\bigcap_{i=1}^n a_i^{-1}N.$$

Comme $\mathfrak{r}^2=0$, on dispose d'une suite exacte de groupes abéliens
$$0\to\mathfrak{r}_\mathrm{add}\xrightarrow{x\mapsto 1+x} A^\times\to (A/\mathfrak{r})^\times\to 0\;.$$
Étant donné que $A/\mathfrak{r}$ est fini, on en déduit l'existence de $N\in\mathbb{N}^*$ tel que $a^N=1$ pour tout $a\in A^\times$. Ainsi, l'image de $\varphi$ est incluse dans le sous-groupe $\mu_N(K)$ des racines $N$-ièmes de l'unité de $K$. Il s'ensuit que $\mathrm{Ker}\,\varphi$ est un sous-groupe d'indice fini de $A^\times$, puis que $N$ est un sous-groupe (additif) d'indice fini de $\mathfrak{r}$, et donc de $A$.

Vérifions que $I$ est un idéal de $A$ : pour $x\in I$ et $t\in A$, soit $t\in\mathfrak{r}$ et $xt\in\mathfrak{r}^2=0$, soit il existe $i\in\llbracket 1,n\rrbracket$ tel que $t-a_i\in\mathfrak{r}$, auquel cas $xt=xa_i\in N$. Ainsi $A.I\subset N$, d'où l'on déduit $A.I=I$, ce qui achève la démonstration.
\end{proof}

La réciproque est vraie si $K$ est assez gros :

\begin{theo}\label{th-ttantipol} Supposons que tout poids de $\mon(A_\mu,K_\mu)$ est antipolynomial.
\begin{enumerate}
\item\label{it-thPfac} L'anneau $A$ est quasi-parfait et son quotient $A/\rj(A)$ est $K$-trivial.
\item Si de plus $K$ contient assez de racines de l'unité au sens où toutes les puissances de la caractéristique de $A$ appartiennent à $\nu(K)$, alors :
\begin{enumerate}
    \item\label{it-thPsp} $A$ est semi-primaire ;
    \item\label{it-thPdur} $A$ possède un idéal $K$-cotrivial de carré nul.
\end{enumerate}
\end{enumerate}
\end{theo}

Ce résultat, qui figure avant tout ici à titre culturel (nous ne nous en servirons pas dans la suite), sera établi dans l'appendice~\ref{ap_afmn} (§\,\ref{par-thIc2}). En effet, bien qu'élémentaire, sa démonstration requiert un certain nombre d'étapes, et est sans rapport avec les catégories de foncteurs.

\begin{rema} On ne peut pas se passer de l'hypothèse que $K$ contient assez de racines de l'unité pour conclure que $A$ est semi-primaire. Ainsi, si $k$ est un corps fini de caractéristique $l\ne p$, et $A$ un anneau local quasi-parfait de corps résiduel $k$, et que le polynôme $X^l-1\in K[X]$ n'a pas d'autre racine que $1$ (ce qui est possible sans que $A$ soit semi-primaire --- il suffit d'adapter légèrement l'exemple~\ref{ex-qpft}), alors l'injection canonique $\mon(k_\mu,K_\mu)\to\mon(A_\mu,K_\mu)$ est bijective, de sorte que tout poids de $\mon(A_\mu,K_\mu)$ est antipolynomial.
\end{rema}

\begin{coro}\label{cor-poiantipol} Les assertions suivantes sont équivalentes :
\begin{enumerate}
\item $A$ est noethérien et tous les poids de $\mon(A_\mu,K_\mu)$ sont antipolynomiaux ;
\item $A$ est $K$-trivial.
\end{enumerate}
\end{coro}

\section{Poids polynomiaux}

\subsection{Décomposition à la Steinberg des poids polynomiaux}

Tout morphisme d'anneaux $A\to K$ définit un poids polynomial (de degré $1$ si $A$ est non nul) ; par conséquent, tout produit (fini) de morphismes d'anneaux $A\to K$ est un poids polynomial.

\begin{defi}\label{df-poidepl} Un poids polynomial de $\mon(A_\mu,K_\mu)$ est dit \index{termin}{deploye@déployé!\emph{(poids polynomial)}}\index{termin}{poids!polynomial!déployé}\textbf{déployé} si c'est un produit fini de morphismes d'anneaux de $A$ dans $K$.

On dit qu'une extension de corps $K\to L$ \textbf{déploie} un poids polynomial $w\in\mon(A_\mu,K_\mu)$ si l'image de $w$ par le monomorphisme de monoïdes $\mon(A_\mu,K_\mu)\to\mon(A_\mu,L_\mu)$ est un poids polynomial déployé.
\end{defi}

Comme conséquence d'un théorème de décomposition à la Steinberg des foncteurs simples polynomiaux de $\F(A,K)$, \cite{DTV} établit le résultat suivant :

\begin{prop}\label{pr-dppol}\cite[prop.~10.1]{DTV} Tout poids polynomial $A\to K$ se déploie sur une extension finie de $K$.

 De plus, la décomposition d'un poids en produit de morphismes d'anneaux est unique si $K$ est de caractéristique nulle, et unique à l'action du morphisme de Frobenius de $K$ près sinon.
\end{prop}

En particulier, si $K$ est algébriquement clos, tout poids polynomial $A\to K$ est déployé. Plus généralement, si $K$ est un corps de décomposition de la catégorie $\mathbf{P}(A)$ (i.e. si $K$ est un corps de décomposition de la $K$-algèbre $A_K$ --- cf. \cite[rem.~3.5]{DTV}), alors tout poids polynomial $A\to K$ est déployé. Par exemple, si $A$ est un corps fini se plongeant dans $K$, alors tout poids polynomial $A\to K$ est déployé.

En revanche, on ne peut généralement pas se passer d'une extension finie de corps dans la décomposition précédente si $K$ n'est pas un corps de décomposition de $\mathbf{P}(A)$ (condition qui revient à dire que tout morphisme d'anneaux de $A$ dans une extension finie de $K$ est à valeurs dans $K$) en vertu du résultat de théorie des corps élémentaire suivant, dont nous nous servirons au chapitre~\ref{sec-finpolp} :

\begin{lemm}\label{lm-extpgg} Soient $L$ une extension finie du corps $K$ et $\varphi : A\to L$ un morphisme d'anneaux. Il existe une extension finie $L'$ du corps $L$ et une famille finie $\alpha_1,\dots,\alpha_n$ de morphismes d'anneaux $A\to L'$ telle que $\alpha_1$ coïncide avec la composée de $\varphi$ et de $L\hookrightarrow L'$ et que le produit des $\alpha_i$ soit à valeurs dans $K$.

De plus, si $L$ est une extension séparable de $K$, on peut supposer qu'il en est de même pour $L'$.
\end{lemm}

\begin{proof}
Soit $L'$ l'extension quasi-galoisienne de $K$ engendrée par $L$ (dans une clôture algébrique de $L$) : c'est une extension finie de $L$ \cite[chap.~5, §\,9, cor.~1 de la prop.~5]{Bki2}, et une extension séparable de $K$ si c'est le cas de $L$ \cite[chap.~5, page~55, discussion précédant la prop.~1]{Bki2}. Cette extension $L'$  convient en prenant pour les $\alpha_i$ des composés de $\varphi$ et des $K$-automorphismes du corps $L$ (éventuellement répétés, si $L$ n'est pas une extension séparable de $K$) par \cite[chap.~5, prop.~4 du §\,8 et prop.~3 du §\,9]{Bki2}.
\end{proof}

La propriété d'unicité de la proposition~\ref{pr-dppol} entraîne aussitôt :

\begin{coro}\label{cor-depl_pft} Supposons le corps $K$ parfait. Si $\pi : A\to K$ est un poids polynomial déployé, pour toute décomposition de $\pi$ en produit de morphismes d'anneaux $\alpha_i : A\to L$, où $L$ est une extension du corps $K$, les $\alpha_i$ sont à valeurs dans $K$.
\end{coro}

\begin{coro}\label{cor-lappol} Pour tout $\pi\in\mon_{\operatorname{pol}}(A,K)$, $\la(\pi)$ est une réunion finie --- non vide si $\pi\ne 1$ --- d'idéaux premiers $\p$ de $A$ tels que $A/\p$ se plonge dans une extension finie du corps $K$. Si de plus $\pi$ est de torsion, alors ces idéaux premiers sont cofinis.
\end{coro}

\begin{coro}\label{cor-lapord} Pour tout $\pi\in\mon_{\operatorname{ord}}(A,K)$, $\la(\pi)$ est une réunion finie d'idéaux premiers de $A$.
\end{coro}

\begin{proof} Combiner les propositions~\ref{pr-decpord},~\ref{pr-zerantipol} et le corollaire~\ref{cor-lappol}.
\end{proof}

Le corollaire suivant s'applique en particulier lorsque le corps $K$ est algébriquement clos et que $F$ est un foncteur de type fini de $\F^\df(\A;K)$, par le corollaire~\ref{cor-tfdp} et la proposition~\ref{pr-dfpord}.

\begin{coro}\label{cor-asl} Soient $\E$ une catégorie de Grothendieck $K$-linéaire et $F$ un foncteur de $\fct(\A,\E)$ possédant une décomposition en poids faible. Supposons que $\Pi(F)$ est fini et constitué de poids ordinaires. Alors l'anneau $\tilde{A}_F$\index{nota}{A@$\tilde{A}_F$} est semi-local.
\end{coro}

\begin{coro}\label{cor-ppol0} Si $K$ est algébriquement clos de caractéristique nulle, le monoïde (commutatif) $\mon_{\operatorname{pol}}(A,K)$ est libre sur $\mathbf{Ann}(A,K)$.
\end{coro}

Dans ce qui suit, pour un entier $p>0$, on note $\mathbb{N}[1/p]$ le sous-monoïde des éléments positifs (pour l'ordre usuel) du groupe additif $\mathbb{Z}[1/p]$.

\begin{coro}\label{cor-ppolp} Supposons que $K$ est algébriquement clos de caractéristique $p>0$. Notons $\sim$ la relation d'équivalence sur $\mathbf{Ann}(A,K)$ engendrée par l'action de l'endomorphisme de Frobenius de $K$ et $\varphi\mapsto [\varphi]$ la fonction canonique $\mathbf{Ann}(A,K)\to\mathbf{Ann}(A,K)/\sim$.

Le monoïde $\mon_{\operatorname{pol}}(A,K)$ est isomorphe à la somme directe sur $[\varphi]\in\mathbf{Ann}(A,K)/\sim$ du sous-monoïde engendré par les éléments de $[\varphi]$. Celui-ci est isomorphe à :
\begin{enumerate}
\item $(\mathbb{Z}/(q-1))_+$ si l'image de $\varphi$ est de cardinal fini $q$ ;
\item $\mathbb{N}[1/p]$ si $\varphi$ est d'image infinie.
\end{enumerate}
\end{coro}

Dans la définition qui suit, la relation d'équivalence $\sim$ sur $\mathbf{Ann}(A,K)$ est celle du corollaire précédent si $p>0$, et l'égalité si $p=0$.

\begin{defi}\label{def-fima} Une famille $(\alpha_i)$ de morphismes d'anneaux $A\to K$ est dite \textbf{indépendante}\index{termin}{independante@indépendante \emph{(famille de morphismes d'anneaux)}|textbf} si les deux conditions suivantes sont réalisées :
\begin{enumerate}
\item chaque $\alpha_i$ est d'image infinie ;
\item les images des $\alpha_i$ dans $\mathbf{Ann}(A,K)/\sim$ sont deux à deux distinctes.
\end{enumerate}
\end{defi}

(En caractéristique $p=0$, cela se réduit à demander que les $\alpha_i$ soient deux à deux distincts.)

Cette terminologie est motivée par l'indépendance linéaire des caractères, qui implique que la famille $(\alpha_i)_{i\in E}$ est indépendante au sens précédent si et seulement si la famille $(\alpha_i^j)_{(i,j)\in E\times\mathbb{N}^*}$ de ses puissances est linéairement indépendante dans l'espace vectoriel $K^A$.

\subsection{Critères de polynomialité des poids}

Si $p$ est un nombre premier, on appelle \emph{$p$-anneau fini} tout anneau fini $A$ vérifiant les conditions équivalentes suivantes :
\begin{enumerate}
\item le cardinal de $A$ est une puissance de $p$ ;
\item la caractéristique de $A$ est une puissance de $p$ ;
\item l'anneau semi-simple $A/\rj(A)$ est de caractéristique $1$ ou $p$.
\end{enumerate}

\begin{prop}\label{pr-anfppol} Supposons que $p$ est non nul et que $A$ est un $p$-anneau fini. Alors tout poids $A\to K$ est polynomial.
\end{prop}

\begin{proof} En effet, toute fonction d'un $p$-groupe fini vers un corps de caractéristique $p$ est polynomiale \cite[chap.~VI, th.~1.2]{Passi}.
\end{proof}

\begin{prop}\label{pr-poirednil} Supposons que $p$ est non nul et que la caractéristique de $A$ est une puissance de $p$. Alors le morphisme canonique $\mon((A/\nil(A))_\mu,K_\mu)\to\mon(A_\mu,K_\mu)$ est un isomorphisme.
\end{prop}

\begin{proof} Soient $x, y\in A$ tels que $x-y\in\nil(A)$ : il existe $n\in\mathbb{N}$ tel que $(x-y)^{p^n}=0$. Si la caractéristique de $A$ est $p^m$, on en déduit $x^{p^{n+m-1}}=y^{p^{n+m-1}}$, d'où $\pi(x)^{p^{n+m-1}}=\pi(y)^{p^{n+m-1}}$ pour tout $\pi\in\mon(A_\mu,K_\mu)$. Comme $K$ est de caractéristique $p$, il s'ensuit que $\pi(x)=\pi(y)$.
\end{proof}

\begin{lemm}\label{lm-car0pol} Soit $V$ un groupe abélien non nul. Si la fonction $V\to K$ envoyant $0$ sur $0$ et $V\setminus\{0\}$ sur $1$ est polynomiale, alors $p\ne 0$ et $V$ est un $p$-groupe fini.
\end{lemm}

\begin{proof} Si $V$ est infini, on peut construire par récurrence une suite $(v_i)_{i\in\mathbb{N}}$ d'éléments de $V$ telle que $\sum_{i\in E}v_i\ne 0$ pour toute partie finie non vide $E$ de $\mathbb{N}$. On en déduit que $V$ est nécessairement fini. L'hypothèse entraîne maintenant que \emph{toute} fonction $V\to K$ est polynomiale, puisque les fonctions polynomiales forment un sous-espace vectoriel de $K^E$ invariant par translation à la source. La conclusion résulte donc de \cite[chap.~VI, th.~1.2]{Passi}.
\end{proof}

L'énoncé et la démonstration du résultat suivant présentent de nombreuses similitudes avec son pendant antipolynomial, le théorème~\ref{th-ttantipol}.

\begin{prop}\label{pr-poitjspol} Les assertions suivantes sont équivalentes :
\begin{enumerate}
\item tout poids de $\mon(A_\mu,K_\mu)$ est polynomial ;
\item $A/\nil(A)$ est un produit fini de corps finis de même caractéristique que $K$.
\end{enumerate}
\end{prop}

\begin{proof} Si $A/\nil(A)$ est un produit fini de corps finis de même caractéristique que $K$, $p\in\nil(A)$, donc les propositions~\ref{pr-poirednil} et~\ref{pr-anfppol} montrent que tout poids de $\mon(A_\mu,K_\mu)$ est polynomial.

Réciproquement, le lemme~\ref{lm-car0pol} montre que pour tout $\p\in\spec(A)$, $A/\p$ est un $p$-anneau fini (considérer le poids $\xi_\p$ introduit à la proposition-définition~\ref{prdfla}), et donc un corps fini de caractéristique $p$. En particulier, tout idéal premier de $A$ est maximal, et il suffit pour conclure de vérifier que $A$ ne possède qu'un nombre fini d'idéaux maximaux. Supposons par l'absurde qu'il existe une suite infinie $(\mathfrak{m}_i)_{i\in\mathbb{N}}$ d'idéaux maximaux deux à deux distincts de $A$. La polynomialité du poids $\xi_{\cup_i\mathfrak{m}_i}$ et le corollaire~\ref{cor-lappol} fournissent un nombre fini d'idéaux premiers $\p_1,\dots,\p_n\in\spec(A)$ tels que $\underset{i\in\mathbb{N}}{\bigcup}\mathfrak{m}_i=\bigcup_{j=1}^n\p_j$. Comme les $\mathfrak{m}_i$ sont maximaux, on a forcément $\mathfrak{m}_i\in\{\p_1,\dots,\p_n\}$ pour chaque $i$, par \cite[ch.~II, §\,1, prop.~2]{Bki-AC1}, contradiction qui achève la démonstration.
\end{proof}

Pour un prolongement de l'énoncé précédent, on renvoie au théorème~\ref{th-ts_simples_pol}.

\begin{prop}\label{pr-ttpoiord} Supposons que tout poids de $\mon(A_\mu,K_\mu)$ est ordinaire. Alors $A$ est quasi-parfait, $A/\rj(A)$ est fini et $A$ est le produit d'un anneau $A'$ tel que tout poids de $\mon(A'_\mu,K_\mu)$ est polynomial et d'un anneau $A''$ tel que tout poids de $\mon(A''_\mu,K_\mu)$ est antipolynomial.
\end{prop}

\begin{proof} Soit $\p\in\spec(A)$. Le poids $\xi_\p$ est le produit d'un poids antipolynomial $\alpha$ et d'un poids polynomial $\pi$ (proposition~\ref{pr-decpord}). Si $\alpha\ne 1$, alors $\p=\la(\xi_\p)$ contient un idéal maximal $K$-cotrivial de $A$ (proposition~\ref{pr-zerantipol}), donc $\p$ est maximal et $A/\p$ fini. Sinon, $\xi_\p$ est polynomial, et on voit comme dans la démonstration de la proposition précédente que $\p$ est maximal et que $A/\p$ est fini.

Ainsi, tous les idéaux premiers de $A$ sont maximaux et cofinis. Par ailleurs, en utilisant le corollaire~\ref{cor-lapord} et en raisonnant comme dans la démonstration précédente, on voit que $A$ n'a qu'un nombre fini d'idéaux maximaux. Il s'ensuit que $A$ est quasi-parfait et $A/\rj(A)$ fini.

En particulier, $A$ est isomorphe à un produit fini d'anneaux locaux, de sorte que, pour démontrer la dernière propriété, il suffit de vérifier que, si $A$ est local, alors soit tous les poids de $\mon(A_\mu,K_\mu)$ sont polynomiaux, soit ils sont tous antipolynomiaux. Si la caractéristique de $A$ est une puissance de $p$, alors $A$ n'a pas d'idéaux cotriviaux autres que $A$, donc il n'y a pas de poids antipolynomial non trivial dans $\mon(A_\mu,K_\mu)$, et $\mon_{\operatorname{pol}}(A,K)=\mon_{\operatorname{ord}}(A,K)=\mon(A_\mu,K_\mu)$ par la proposition~\ref{pr-decpord}. Si la caractéristique de $A$ est étrangère à l'exposant caractéristique de $K$, alors il n'y a pas de morphisme d'anneaux de $A$ dans une extension de $K$, donc pas de poids polynomial non trivial dans $\mon(A_\mu,K_\mu)$ par la proposition~\ref{pr-dppol}, de sorte qu'on a alors $\mon_{\operatorname{AP}}(A,K)=\mon_{\operatorname{ord}}(A,K)=\mon(A_\mu,K_\mu)$, cqfd.
\end{proof}

\begin{coro}\label{cor-poiord} Les assertions suivantes sont équivalentes :
\begin{enumerate}
\item $A$ est noethérien et tout poids de $\mon(A_\mu,K_\mu)$ est ordinaire ;
\item $A$ est fini.
\end{enumerate}
\end{coro}

Ce corollaire implique en particulier que, si $A$ est un anneau noethérien infini, alors il existe dans $\F(A,K)$ au moins un foncteur simple n'appartenant pas à $\F^\df(A,K)$.

\section{Lemmes arithmétiques}\label{slmar}

\begin{nota}
Soit $l>1$ un entier. On note \index{nota}{s@$s_l$ \emph{(somme des chiffres dans l'écriture en base $l$)}}
$s_l(n)$ la somme des chiffres d'un élément $n$ de $\mathbb{N}[1/l]$ dans son écriture $l$-adique.
\end{nota}

Le résultat suivant, classique et facile, est laissé en exercice.

\begin{lemm}\label{lm-som_ebp} Pour toute famille presque nulle $(a_i)_{i\in\mathbb{N}}$ d'entiers naturels et tout entier $l>1$, on a
$$\sum_{i\in\mathbb{N}}a_i\ge s_l\Big(\sum_{i\in\mathbb{N}}a_i l^i\Big)\,.$$
\end{lemm}

\begin{prop}\label{pr-degdp} Soit $(\alpha_1,\dots,\alpha_n)$ une famille finie indépendante de morphismes d'anneaux $A\to K$. Alors pour tous entiers naturels $i_1,\dots,i_n$, le degré polynomial de la fonction $\prod_{r=1}^n\alpha_r^{i_r} : A\to K$ est $\sum_{r=1}^n i_r$ si $p=0$ et $\sum_{r=1}^n s_p(i_r)$ sinon.
\end{prop}

\begin{proof}
Cela découle du résultat d'unicité de la proposition~\ref{pr-dppol} et, pour $p$ premier, du lemme~ \ref{lm-som_ebp}.
\end{proof}

Si $p>1$ est un entier, on définit une fonction $C_p : \mathbb{N}\times\mathbb{N}\to\mathbb{N}$ par récurrence comme suit : $C_p(0,n)=1$ si $n=0$ et $0$ si $n>0$, et, pour $d>0$,
$$C_p(d,n)=n.\sum_{0\le u\le\frac{d-n}{p-1}}C_p(d-1,n-1+u(p-1)).$$
On note que $C_p(d,n)=0$ pour $n>d$. Une estimation grossière montre qu'on a toujours $C_p(d,n)\le (d!)^2$.

\begin{lemm}\label{lm-arit-padic} Soient $d\ge 0$ et $p\ge 2$ des entiers et $r\in\mathbb{N}[1/p]$. L'ensemble des suites presque nulles d'entiers naturels $(a_i)_{i\in\mathbb{Z}}$ telles que $\sum_{i\in\mathbb{Z}}a_i=d$ et $\sum_{i\in\mathbb{Z}}a_i.p^i=r$ est fini, de cardinal majoré par $C_p(d,s_p(r))$, donc par $(d!)^2$.
\end{lemm}

\begin{proof} Notons $E(d,r)$ l'ensemble en question et, pour $j\in\mathbb{Z}$, $E_j(d,r):=\{(a_i)\in E(d,r)\,|\,a_j>0\}$.

On montre par récurrence sur $d$ que $\operatorname{Card}E(d,r)\le C_p(d,s_p(r))$ pour tout $r\in\mathbb{N}$. La propriété est évidente pour $d=0$, on suppose donc $d>0$ et l'inégalité $\operatorname{Card}E(d-1,r)\le C_p(d-1,s_p(r))$ vérifiée.

Soient $j\in\mathbb{Z}$ et $(a_i)_{i\in\mathbb{Z}}$ un élément de $E_j(d,r)$. On a donc $p^j\le r$ ; notons $t$ le plus petit entier supérieur à $j$ tel que $d_t\ne 0$, où $r=\sum_{i\in\mathbb{Z}}d_i.p^i$ (avec $0\le d_i\le p-1$) est l'écriture $p$-adique de $r$. On a alors $s_p(r-p^j)=s_p(r)-1+u(p-1)$ où $u:=t-j$. En considérant la suite $(a'_i)$ donnée par $a'_j=a_{j}-1$ et $a'_i=a_i$ pour $i\ne j$, on en tire $\operatorname{Card}E_j(d,r)\le C_p(d-1,s_p(r)-1+u(p-1))$ grâce à l'hypothèse de récurrence. Si cet entier est non nul, on a $s_p(r)-1+u(p-1)\le d-1$, d'où $0\le u\le\frac{d-s_p(r)}{p-1}$. Ainsi
$$\operatorname{Card}E(d,r)\le\sum_{j\in\mathbb{Z}}\operatorname{Card}E_j(d,r)\le\underset{d_t\ne 0}{\sum_{t\in\mathbb{Z}}}\sum_{0\le u\le\frac{d-s_p(r)}{p-1}}C_p(d-1,s_p(r)-1+u(p-1)).$$
Comme l'ensemble $\{t\in\mathbb{Z}\,|\,d_t\ne 0\}$ est de cardinal majoré par $s_p(r)$, cela termine la démonstration.
\end{proof}

\begin{rema} La majoration est très éloignée de l'optimalité, mais seul nous importe de disposer d'une borne ne dépendant que de $d$.
\end{rema}

Le lemme~\ref{lm-arit-padic} nous permet d'établir le résultat suivant, qui jouera un rôle important dans la section~\ref{sec-finpolp}.

\begin{prop}\label{pr-decfinie} Soient $d\in\mathbb{N}$ et $\pi\in\mon(A_\mu,K_\mu)$. L'ensemble
$$\{(\alpha_1,\dots,\alpha_d)\in\mathbf{Ann}(A,K)^d\,|\,\pi=\prod_{i=1}^d\alpha_i\}$$
est fini.
\end{prop}

\begin{proof} Quitte à agrandir le corps $K$, on ne perd pas en généralité à le supposer algébriquement clos. Le cas où sa caractéristique est nulle découle du corollaire~\ref{cor-ppol0}. Sinon, cela provient du corollaire~\ref{cor-ppolp} et du lemme~\ref{lm-arit-padic}.
\end{proof}

Nous aurons également usage, pour $p\ge 2$, de la fonction \index{nota}{sig@$\sigma_p$}
$\sigma_p : \mathbb{N}\to\mathbb{N}$ \label{psig} définie par $\sigma_p(n)=\max\{s_p(i)\,|\,i\in\llbracket 0, n\rrbracket\}$. On vérifie aussitôt que $\sigma_p(n)$ est le plus grand $N\in\mathbb{N}$ tel que $\sigma_p^!(N)\le n$, où $\sigma_p^!(N):=(r+1)p^a-1$ si $a$ (resp. $r$) est le quotient (resp. le reste) de la division euclidienne de $N$ par $p-1$. Par exemple, $\sigma_2^!(N)=2^N-1$ et $\sigma_2(n)=[\log_2(n+1)]$ (où les crochets désignent la partie entière) ; en général $p^{[\frac{N}{p-1}]}-1\le\sigma_p^!(N)\le p^{\lceil\frac{N}{p-1}\rceil}-1$ (où $\lceil x\rceil:=-[-x]$) et $(p-1)[\log_p(n+1)]\le\sigma_p(n)\le (p-1)[\log_p(n+1)]+p-2$. Nous n'aurons pas besoin de ces relations, mais elles permettent de donner une borne explicite au degré polynomial dans le théorème~\ref{th-chipol} ci-après.

Le lemme suivant interviendra dans la démonstration de la proposition~\ref{pr-degpoipar}.

\begin{lemm}\label{lmari-sig} Soient $n$, $r_0,\dots,r_n$ et $a_0,\dots,a_n$ des entiers naturels. On a
$$\sigma_p\Big(\sum_{i=0}^n a_i\Big)\le\max\Big\{s_p\Big(\sum_{j=0}^n t_j p^{r_j}\Big)\,|\,0\le t_j\le a_j\Big\}.$$
\end{lemm}

\begin{proof} Quitte à réordonner les $a_i$ et les $r_i$, on peut supposer que la suite finie $(r_i)$ est croissante. Quitte à retrancher $r_0$ à tous les $r_i$, on peut également supposer $r_0=0$.

On montre le résultat par récurrence sur $n$. Il est trivial pour $n=0$, on suppose donc $n>0$ et l'assertion établie pour $n-1$.

Supposons d'abord $a_0<p^{r_1}$. Quitte à remplacer la suite $(a_i)$ par une suite $(a'_i)$ avec $0\le a'_i\le a_i$ et $s_p\Big(\sum_{i=0}^n a'_i\Big)=\sigma_p\Big(\sum_{i=0}^n a_i\Big)$, on voit qu'il suffit de montrer qu'il existe des entiers $0\le t_j\le a_j$ (pour $0\le j\le n$) tels que $s_p\Big(\sum_{i=0}^n a_i\Big)\le s_p\Big(\sum_{j=0}^n t_j p^{r_j}\Big)$. L'hypothèse de récurrence fournit des entiers naturels $t_j\le a_j$ pour $1\le j\le n$ tels que $s_p\Big(\sum_{i=1}^n a_i\Big)\le s_p\Big(\sum_{j=1}^n t_j p^{r_j}\Big)$.
Par ailleurs, comme $a_0<p^{r_1}$, on a $s_p(a_0+u)=s_p(a_0)+s_p(u)$ pour tout entier naturel $u$ divisible par $p^{r_1}$. On a donc
$$s_p\Big(\sum_{i=0}^n a_i\Big)\le s_p(a_0)+s_p\Big(\sum_{i=1}^n a_i\Big)\le s_p(a_0)+s_p\Big(\sum_{j=1}^n t_j p^{r_j}\Big)=s_p\Big(\sum_{j=0}^n t_j p^{r_j}\Big)$$
comme souhaité.

Supposons maintenant $a_0\ge p^{r_1}$. L'hypothèse de récurrence montre qu'il existe des entiers $0\le t'\le a_0+a_1$ et $0\le t_j\le a_j$ pour $2\le j\le n$ tels que
$$\sigma_p\Big(\sum_{i=0}^n a_i\Big)\le s_p\Big(t'+\sum_{j=2}^n t_j p^{r_j}\Big).$$
Si $t'\le a_0$, le choix $t_0=t'$ et $t_1=0$ convient. Supposons désormais $t'>a_0$.
Comme $a_0\ge p^{r_1}$, il existe un entier $0\le t_0\le a_0$ congru à $t'$ modulo $p^{r_1}$ ; considérons le plus grand d'entre eux : on a $t_0\le a_0<t_0+p^{r_1}$, et $t'=t_0+t_1 p^{r_1}$ pour un $t_1\in\mathbb{N}$. On a donc
$$t_1=p^{-r_1}(t'-t_0)\le p^{-r_1}(a_1+a_0-t_0)<p^{-r_1}(a_1+p^{r_1})\le a_1+1,$$
d'où $t_1\le a_1$. Cela termine la démonstration.
\end{proof}

\chapter[Foncteurs polynomiaux à poids]{Propriétés de finitude des foncteurs polynomiaux à poids}\label{sec-finpolp}

\begin{cvi}
Dans tout le chapitre~\ref{sec-finpolp}, $\pi\in\mon(A_\mu,K_\mu)$ désigne un poids, $d$ un entier naturel, $r$ un entier strictement positif, $\A$ une catégorie additive $A$-linéaire essentiellement petite et $\E$ une catégorie de Grothendieck $K$-linéaire.
\end{cvi}

L'objectif principal de ce chapitre consiste à démontrer que tout foncteur polynomial noethérien de $\F(A,K)$ possédant une décomposition en poids est à valeurs de dimension finie. De plus, l'hypothèse noethérienne peut être affaiblie en une hypothèse de type-finitude sous une hypothèse de finitude raisonnable sur $A$ (par exemple, si $A$ est un anneau essentiellement de type fini), et la conclusion s'étend aussitôt aux foncteurs \phs\ possédant une décomposition en poids.

Ces résultats, établis dans la section~\ref{parfpp}, reposent lourdement sur des considérations de pure algèbre commutative exposées dans la section~\ref{swd}. Celles-ci ne sont toutefois pas dénués de liens avec la théorie élémentaire des poids des foncteurs de $\F(A,K)$ du chapitre~\ref{spfca}, car elles utilisent de façon fondamentale la proposition~\ref{pr-dppol}, qui est déduite de résultats fonctoriels de \cite{DTV}.

La section~\ref{sliencr} présente le lien, direct, entre les constructions de la section~\ref{swd} et les foncteurs polynomiaux à poids de $\F(A,K)$ (ou plus généralement $\fct(\A,\E)$) ; elles seront également utilisées plus tard dans ce mémoire, pour l'étude des foncteurs polynomiaux \hds\ et \hts.

Nous terminons ce chapitre en donnant une condition suffisante, à la section~\ref{sfdp}, pour qu'un foncteur polynomial de type fini de $\F(A,K)$ possède une décomposition en poids faible si et seulement s'il est fini.

\section{L'algèbre $W_d(A,\pi)$}\label{swd}

\begin{nota} On désigne par \index{nota}{W@$W_d$, $W_d(A,\pi)$, $W_d^{[r]}(A,\pi)$} $W_d(A,\pi)$ le quotient de la $K$-algèbre $A_K^{\otimes d}\simeq K\otimes_\mathbb{Z}A^{\otimes_\mathbb{Z}d}$ par l'idéal $\I_d(A,\pi)$ engendré par les éléments $1\otimes a^{\otimes d}-\pi(a)\otimes 1^{\otimes d}$ pour $a\in A$. Plus généralement, on note $W_d^{[r]}(A,\pi)$ le quotient de $A_K^{\otimes d}$ par la puissance $r$-ième de l'idéal $\I_d(A,\pi)$.
\end{nota}

Cette construction est fonctorielle en $A$ au sens suivant : tout morphisme d'anneaux $\varphi : B\to A$ induit un morphisme de $K$-algèbres $W_d^{[r]}(B,\pi\circ\varphi)\to W_d^{[r]}(A,\pi)$. Il existe une fonctorialité analogue en le corps $K$.

On note que $\I_d(A,\pi)$ est stable par l'action canonique du groupe symétrique $\Si_d$ sur l'algèbre $A_K^{\otimes d}$, de sorte que $W_d(A,\pi)$ possède une action naturelle de $\mathfrak{S}_d$.

\begin{rema}\label{rq-Wprod} Pour $\alpha, \beta\in\mon(A_\mu,K_\mu)$ et $i, j\in\mathbb{N}$, l'isomorphisme canonique $A_K^{\otimes (i+j)}\xrightarrow{\simeq}A_K^{\otimes i}\otimes A_K^{\otimes j}$ envoie $\I_{i+j}(A,\alpha\beta)$ dans $\I_i(A,\alpha)\otimes A_K^{\otimes j}+A_K^{\otimes i}\otimes\I_j(A,\beta)$, d'où un morphisme surjectif de $K$-algèbres
$W_{i+j}(A,\alpha\beta)\twoheadrightarrow W_i(A,\alpha)\otimes W_j(A,\beta)$ et plus généralement $W_{i+j}^{[r+s-1]}(A,\alpha\beta)\twoheadrightarrow W_i^{[r]}(A,\alpha)\otimes W_j^{[s]}(A,\beta)$ pour tous entiers $r, s>0$.

On en déduit par récurrence sur $n\in\mathbb{N}$, pour tous poids $\rho_1,\dots,\rho_n\in\mon(A_\mu,K_\mu)$ et tous $i_1,\dots,i_n\in\mathbb{N}$, un morphisme surjectif de $K$-algèbres
\begin{equation}\label{eq-Wprod}
W_{i_1+\dots+i_n}^{[n(r-1)+1]}(A,\rho_1\dots\rho_n)\twoheadrightarrow W_{i_1}^{[r]}(A,\rho_1)\otimes\dots\otimes W_{i_n}^{[r]}(A,\rho_n).
\end{equation}
\end{rema}

\begin{prop}\label{prw-fond}
\begin{enumerate}
\item L'anneau $W_d^{[r]}(A,\pi)$ est quasi-parfait.
\item\label{itprwdf} L'algèbre $W_d^{[r]}(A,\pi)/\rj(W_d^{[r]}(A,\pi))$ est de dimension finie sur $K$.
\item\label{pr-wdepl} La $K$-algèbre $W_d(A,\pi)$ est déployée dans chacun des deux cas suivants :
\begin{enumerate}
    \item $K$ est un corps de décomposition de $\mathbf{P}(A)$ ;
    \item  $\pi$ est un poids polynomial déployé et $K$ est un corps parfait.
\end{enumerate}
\end{enumerate}
\end{prop}

\begin{proof} Un idéal premier $\p$ de $W_d^{[r]}(A,\pi)$ est le noyau d'un morphisme de $K$-algèbres de $W_d^{[r]}(A,\pi)$ vers une extension $L$ du corps $K$. Un tel morphisme est donné par des morphismes d'anneaux $\alpha_i : A\to L$ pour $1\le i\le d$ tels que le morphisme $K\otimes_\mathbb{Z}A^{\otimes_\mathbb{Z}d}\to L\quad \lambda\otimes a_1\otimes\dots\otimes a_d\mapsto\lambda\alpha_1(a_1)\dots\alpha_d(a_d)$ soit nul sur $\I_d(A,\pi)$, c'est-à-dire tels que $\pi=\prod_{i=1}^d\alpha_i$. La proposition~\ref{pr-dppol} montre alors que les $\alpha_i$ prennent leurs valeurs dans une sous-extension finie de l'extension $L$ du corps $K$. Cela entraîne que $W_d^{[r]}(A,\pi)/\p$ est de dimension finie sur $K$, et est donc un corps. Ainsi, tout idéal premier de $W_d^{[r]}(A,\pi)$ est maximal. La proposition~\ref{pr-decfinie} montre que $W_d^{[r]}(A,\pi)$ n'a qu'un nombre fini d'idéaux maximaux, ce qui achève de prouver les deux premières assertions.

Si $K$ est un corps de décomposition de $\mathbf{P}(A)$, alors tout morphisme d'anneaux de $A$ vers une extension finie du corps $K$ est à valeurs dans $K$, donc l'algèbre $W_d^{[r]}(A,\pi)$ est déployée. La même conclusion vaut lorsque $\pi$ est un poids polynomial déployé et que le corps $K$ est parfait grâce au corollaire~\ref{cor-depl_pft}.
\end{proof}

\begin{rema}\label{rqw-corpremier} La démonstration précédente montre aussi que la $K$-algèbre $W_d(A,\pi)$ est non nulle si et seulement si $\pi$ possède une décomposition en produit de $d$ morphismes d'anneaux de $A$ vers une extension de $K$. Dans ce cas, le poids $\pi$ se factorise de façon unique en un poids $\tilde{\pi} : A_k=k\otimes_\mathbb{Z} A\to K$, où $k$ est le sous-corps premier de $K$, à travers le morphisme canonique $A\to A_k$. Il est clair que $W_d(A,\pi)$ s'identifie à $W_d(A_k,\tilde{\pi})$. Autrement dit, on peut sans perte de généralité supposer que $A$ est une $k$-algèbre.
\end{rema}

\begin{coro}\label{cor-wdfgal} Si $W_d^{[r]}(A,\pi)$ est un anneau noethérien, alors $\dim_K W_d^{[r]}(A,\pi)$ est finie. 
\end{coro}

\begin{proof}
Si $W_d^{[r]}(A,\pi)$ est un anneau noethérien, la première assertion de la proposition~\ref{prw-fond} et la proposition~\ref{pr-qpsi}.\,\ref{it-qpfna} montrent que la $K$-algèbre $W_d^{[r]}(A,\pi)$ est artinienne. Comme son quotient par son radical est de dimension finie (grâce à la proposition~\ref{prw-fond} également), cela donne la conclusion.
\end{proof}

\begin{coro}\label{cor-algfinie} Si la $K$-algèbre $A_K$ est essentiellement de type fini\index{termin}{essentiellement de type fini \emph{(anneau, algèbre)}|textbf} (c'est-à-dire isomorphe à une localisation d'une $K$-algèbre de type fini), alors $W_d^{[r]}(A,\pi)$ est de dimension finie sur $K$.
\end{coro}

\begin{proof} La $K$-algèbre
$A_K^{\otimes d}$ est essentiellement de type fini comme $A_K$, elle est donc noethérienne. L'algèbre $W_d^{[r]}(A,\pi)$, qui en est un quotient, est également noethérienne. On conclut par le corollaire~\ref{cor-wdfgal}.
\end{proof}

\begin{coro}\label{cor-noethp} Le foncteur canonique $\E\underset{K}{\boxtimes}{W_d^{[r]}(A,\pi)}\to\E$ préserve les objets noethériens.
\end{coro}

\begin{proof} Cela résulte de la proposition~\ref{prw-fond} et du corollaire~\ref{cor-hole}.
\end{proof}

Dans l'énoncé suivant, si $p$ est un nombre premier, nous disons qu'une $\FF_p$-algèbre est \emph{$p$-parfaite}\,\footnote{Cette notion généralise celle de corps parfait mais n'a rien à voir avec celle d'anneau parfait !}\index{termin}{p-parfait@$p$-parfait(e) \emph{(anneau, $\FF_p$-algèbre)}} si son endomorphisme de Frobenius est surjectif.

\begin{prop}\label{pr-wss} Supposons que $p$ est non nul et que la $\FF_p$-algèbre $A/p$ est $p$-parfaite (par exemple, que $A$ est un corps parfait de caractéristique $p$).

Alors l'algèbre $W_d^{[r]}(A,\pi)$ est semi-simple et de dimension finie sur $K$.
\end{prop}

\begin{proof} On peut supposer sans perte de généralité que $A$ est une $\FF_p$-algèbre (remarque~\ref{rqw-corpremier}), et que le corps $K$ est parfait, quitte à le plonger dans un corps parfait et étendre les scalaires.

Il suffit de montrer que $W_d^{[r]}(A,\pi)$ est semi-simple, puisqu'on sait, grâce à la proposition~\ref{prw-fond}, que le quotient de cette algèbre par son radical est de dimension finie sur $K$. De plus, comme $W_d^{[r]}(A,\pi)$ est quasi-parfait, on peut se contenter de montrer que cet anneau est réduit.

Montrons d'abord la conclusion sous l'hypothèse supplémentaire que l'endomorphisme de Frobenius de $A=A/p$ est \emph{bijectif}. C'est alors aussi le cas des endomorphismes de Frobenius de $A_K^{\otimes d}$ et de $W_d^{[r]}(A,\pi)$, par fonctorialité de la construction $W_d^{[r]}$. En particulier, ces algèbres sont réduites.

Si l'on suppose maintenant seulement que l'endomorphisme de Frobenius de $A$ est surjectif, il existe un morphisme surjectif de $\FF_p$-algèbres $B\twoheadrightarrow A$ où le Frobenius de $B$ est bijectif (prendre par exemple pour $B$ la colimite de copies de la $\FF_p$-algèbre libre sur $A$ le long du Frobenius). Comme le foncteur $W_d^{[r]}$ préserve les morphismes surjectifs et qu'un quotient d'une algèbre semi-simple est semi-simple, la conclusion s'ensuit.
\end{proof}

Afin de généraliser la proposition précédente, nous utiliserons le lemme facile suivant :
\begin{lemm}\label{lm-etfW} Soit $\iota : B\to A$ un morphisme d'anneaux. Si $A$ est une $B$-algèbre de type fini (resp. essentiellement de type fini), alors $W_d^{[r]}(A,\pi)$ est une $W_d^{[r]}(B,\pi\circ\iota)$-algèbre de type fini (resp. essentiellement de type fini).
\end{lemm}

\begin{proof}
$A^{\otimes_\mathbb{Z} d}$ est une $B^{\otimes_\mathbb{Z} d}$-algèbre de type fini (resp. essentiellement de type fini). Le produit tensoriel par $K$ et le passage au quotient qui définissent $W_d^{[r]}$ préservent également la propriété (essentiellement) de type fini pour les algèbres, d'où le lemme.
\end{proof}

Le résultat suivant généralise à la fois le corollaire~\ref{cor-algfinie} et la proposition~\ref{pr-wss}.

\begin{prop}\label{pr-Wdf} Supposons que $p$ est non nul et que $A/p$ est une algèbre essentiellement de type fini sur une $\FF_p$-algèbre $p$-parfaite $B$.

Alors $W_d^{[r]}(A,\pi)$ est de dimension finie sur $K$.
\end{prop}

\begin{proof}
Par la proposition~\ref{pr-wss}, $W_d^{[r]}(B,\pi)$ est de dimension finie sur $K$. Il s'ensuit grâce au lemme~\ref{lm-etfW} que $W_d^{[r]}(A,\pi)$ est une $K$-algèbre essentiellement de type fini, donc en particulier noethérienne. Le corollaire~\ref{cor-wdfgal} permet de conclure.
\end{proof}

\begin{prop}\label{pr-derWd} Considérons les assertions suivantes.
\begin{enumerate}
\item\label{itWdd1} Pour toute extension finie $L$ du corps $K$ et toute décomposition de $A\xrightarrow{\pi}K\hookrightarrow L$ en produit de $d$ morphismes d'anneaux $\alpha_i : A\to L$ ($i\in\llbracket 1,d\rrbracket$), on a $\dim_K\mathrm{Der}_{\alpha_i}(A,L)<\infty$ (resp. $\mathrm{Der}_{\alpha_i}(A,L)=0$) pour tout $i\in\llbracket 1,d\rrbracket$. 
\item\label{itWdd2} Pour tout $r\in\mathbb{N}$ et tout $\mathfrak{p}\in\mathrm{Spec}(W_d^{[r]}(A,\pi))$, on a $\dim_K\mathfrak{p}/\mathfrak{p}^2<\infty$ (resp. $\mathfrak{p}^2=\mathfrak{p}$).
\item\label{itWdd3} Pour tout $r\in\mathbb{N}$, tout quotient semi-primaire de la $K$-algèbre $W_d^{[r]}(A,\pi)$ est de dimension finie sur $K$ (resp. semi-simple de dimension finie).
\end{enumerate}
Alors \ref{itWdd1}$\Rightarrow$\ref{itWdd2}$\Rightarrow$\ref{itWdd3}, et \ref{itWdd3}$\Rightarrow$\ref{itWdd1} si toutes les extensions de corps $K\subset L$ qui apparaissent dans \ref{itWdd1} sont séparables --- en particulier, si le corps $K$ est parfait.
\end{prop}

\begin{proof}
Pour $\mathfrak{p}\in\mathrm{Spec}(W_d^{[r]}(A,\pi))$, notons $\tilde{\mathfrak{p}}$ l'image réciproque de $\mathfrak{p}$ par la projection $A_K^{\otimes d}\twoheadrightarrow W_d^{[r]}(A,\pi)$ : il existe une extension finie $L$ du corps $K$ et une décomposition de $A\xrightarrow{\pi}K\hookrightarrow L$ en produit de $d$ morphismes d'anneaux $\alpha_i : A\to L$ ($i\in\llbracket 1,d\rrbracket$) telles que $\tilde{\mathfrak{p}}$ soit le noyau du morphisme de $K$-algèbres $A_K^{\otimes d}\to L$ donné par $a_1\otimes\dots\otimes a_d\mapsto\tilde{\alpha}_1(a_1)\dots\tilde{\alpha}_d(a_d)$, où $\tilde{\alpha}_i : A_K\to L$ est le morphisme de $K$-algèbres prolongeant $\alpha_i\in\mathbf{Ann}(A,L)$. Il s'ensuit que
$$\tilde{\mathfrak{p}}/\tilde{\mathfrak{p}}^2\simeq\bigoplus_{i=1}^d\mathfrak{m}_{\alpha_i}/\mathfrak{m}_{\alpha_i}^2$$
(où $\mathfrak{m}_{\alpha_i}=\mathrm{Ker}\,\tilde{\alpha}_i$).
Comme $\mathfrak{p}/\mathfrak{p}^2$ est un quotient de $\tilde{\mathfrak{p}}/\tilde{\mathfrak{p}}^2$, on en déduit l'implication \ref{itWdd1}$\Rightarrow$\ref{itWdd2} en utilisant \eqref{eq-ident_der} (page~\pageref{eq-ident_der}).

L'implication \ref{itWdd2}$\Rightarrow$\ref{itWdd3} découle du lemme~\ref{lm-qsp} ainsi que des propositions~\ref{pr-qpsp} et~\ref{prw-fond}.

Enfin, supposons \ref{itWdd3} vérifiée et soient $L$, $\alpha_1,\dots,\alpha_d$ comme dans \ref{itWdd1}. Si $L$ est une extension séparable de $K$, alors tout quotient semi-primaire de $W_d^{[r]}(A,\pi)\otimes L$ est de dimension finie (resp. semi-simple de dimension finie), parce qu'on a alors $\rj(W_d^{[r]}(A,\pi)\otimes L)\simeq\rj(W_d^{[r]}(A,\pi))\otimes L$ \cite[§\,12, cor. de la prop.~10]{Bki}, ce qui permet de transporter les critères du lemme~\ref{lm-qsp} par extension des scalaires de $K$ à $L$. On peut donc, quitte à étendre les scalaires, supposer $K=L$. Le morphisme surjectif~\eqref{eq-Wprod} de la remarque~\ref{rq-Wprod} montre alors que tout quotient semi-primaire de $W_1^{[2]}(A,\alpha_i)$ (pour tout $i\in\llbracket 1,d\rrbracket$) est de dimension finie (resp. semi-simple de dimension finie). Comme $\I_1(A,\alpha)$ est l'idéal $\mathfrak{m}_\alpha$ de $A_K$ pour tout $\alpha\in\mathbf{Ann}(A,K)$, la $K$-algèbre $W_1^{[2]}(A,\alpha)$ est semi-primaire, et elle est de dimension finie (resp. semi-simple de dimension finie) si et seulement si $\dim_K\mathrm{Der}_\alpha(A,K)<\infty$ (resp. $\mathrm{Der}_\alpha(A,K)=0$), d'où \ref{itWdd1}.
\end{proof}

\begin{coro}\label{cor-qspfW} Les assertions suivantes sont équivalentes :
\begin{enumerate}
\item\label{itspfW1} pour tous $d\in\mathbb{N}$, $r\in\mathbb{N}^*$ et $\pi\in\mon(A_\mu,K_\mu)$, tout quotient semi-primaire de la $K$-algèbre $W_d^{[r]}(A,\pi)$ est de dimension finie (resp. semi-simple de dimension finie) ;
\item\label{itspfW2} la condition \textnormal{(CFD)}$^+$ (resp. \textnormal{(CAD)}$^+$)\index{nota}{CAD@(CAD), (CAD)$^+$ \emph{(conditions d'annulation des dérivations)}}\index{nota}{CFD@(CFD), (CFD)$^+$ \emph{(conditions de finitude sur les dérivations)}} introduite page \pageref{eqcfd} (resp. \pageref{eqcad}) est vérifiée.
\end{enumerate}
\end{coro}

\begin{proof}
La proposition~\ref{pr-derWd} montre l'implication \ref{itspfW2}$\Rightarrow$\ref{itspfW1}.

Supposant maintenant la condition \ref{itspfW1} vérifiée, grâce à la proposition~\ref{pr-derWd} et au lemme~\ref{lm-extpgg}, on voit que $\mathrm{Der}_\alpha(A,L)=0$ si $\alpha$ est un morphisme d'anneaux de $A$ vers une extension finie \emph{séparable} $L$ de $K$. Le lemme~\ref{lm-cafd_insep} permet de conclure.
\end{proof}

\begin{prop}\label{pr-pradw}\begin{enumerate}
\item Si $p$ est nul ou supérieur ou à égal à $d$, alors $t^{d!}=0$ pour tout $t\in\rj(W_d(A,\pi))$.
\item Si $p$ est nul ou strictement supérieur à $d!$, alors $\rj(W_d(A,\pi))^{d!}=0$, et l'anneau $W_d^{[r]}(A,\pi)$ est semi-primaire.
\end{enumerate}
\end{prop}

\begin{proof} Tout élément de $A_K^{\otimes d}$ annule un polynôme unitaire de degré $d!$ à coefficients dans la sous-algèbre des invariants $(A_K^{\otimes d})^{\Si_d}$. Si $p$ est nul ou supérieur à $d$, cette sous-algèbre est engendrée par les éléments $1\otimes a^{\otimes d}$ où $a$ parcourt $A$ (et les tenseurs sont pris sur $\mathbb{Z}$). Comme l'image d'un tel élément dans $W_d(A,\pi)$ appartient à $K$, on en déduit que tout élément $t$ de $W_d(A,\pi)$ annule un polynôme unitaire de degré $d!$ à coefficients dans $K$. Si $t$ appartient à $\rj(W_d(A,\pi))$, il s'ensuit que $t^{d!}=0$.

La nullité de  $\rj(W_d(A,\pi))^{d!}$ pour $p=0$ ou $p>d!$ découle de la première assertion et de l'observation générale que si un idéal $I$ d'un anneau est tel que $x^n=0$ pour tout $x\in I$ et que $n!$ est inversible dans cet anneau, alors $I^n=0$. Comme $W_d^{[r]}(A,\pi)$ est manifestement semi-primaire dès lors que $W_d(A,\pi)$ l'est, cela termine la démonstration.
\end{proof}

\begin{rema}
Pour $d>2$, l'exposant $d!$ qui apparaît dans l'énoncé précédent n'est peut-être pas optimal. Mais on ne peut pas se passer de l'hypothèse que $p$ est nul ou assez grand pour conclure que $W_d(A,\pi)$ est semi-primaire, même si $A$ est un corps (cf. exemple~\ref{ex-Wcorimp} ci-après).
\end{rema}

En combinant la proposition~\ref{pr-pradw} et le corollaire~\ref{cor-qspfW}, il vient :

\begin{coro}\label{cor-wssder} Supposons $p=0$. Alors la condition \textnormal{(CFD)}$^+$ (resp. \textnormal{(CAD)}$^+$)\index{nota}{CAD@(CAD), (CAD)$^+$ \emph{(conditions d'annulation des dérivations)}}\index{nota}{CFD@(CFD), (CFD)$^+$ \emph{(conditions de finitude sur les dérivations)}} est vérifiée si et seulement si la $K$-algèbre $W_d^{[r]}(A,\pi)$ est de dimension finie (resp. semi-simple) pour tout poids $\pi\in\mon(A_\mu,K_\mu)$ et tous entiers $d\ge 0$ et $r>0$.
\end{coro}

\begin{rema}
Si $A$ est un corps s'injectant dans $K$, $(A,K)$ vérifie la condition (CAD)$^+$ si et seulement si $A$ est sans dérivation (proposition-définition~\ref{prdf-sander}). La proposition~\ref{pr-wss} montre donc qu'on peut se passer de toute hypothèse sur la caractéristique $p$ de $K$ dans l'énoncé précédent \emph{si $A$ est un corps}. On prendra toutefois garde que ce n'est pas le cas si $A$ est un anneau quelconque, comme nous le verrons dans l'exemple~\ref{ex-Wbiz}.
\end{rema}

\begin{prop}\label{pr-Wpcarp} Supposons que $K$ et $A$ sont de même caractéristique $p$ première ; soit $\alpha\in\mathbf{Ann}(A,K)$. Alors il existe un morphisme surjectif de $K$-algèbres $W_p(A,\alpha)\twoheadrightarrow A\,_{\mathrm{Fr}}{\underset{A}{\otimes}}{}_\alpha K$, qui est un isomorphisme si $p=2$, où $\mathrm{Fr}$ désigne l'endomorphisme de Frobenius de la $\FF_p$-algèbre $A$.
\end{prop}

\begin{proof}
Le produit $p$-itéré tensorisé par $K$ est un morphisme surjectif d'anneaux $K\otimes A^{\otimes p}\twoheadrightarrow K\otimes A$ (où les produits tensoriels sont pris sur $\mathbb{Z}$) ; il induit un morphisme surjectif de $W_p(A,\alpha)$ sur le quotient de $K\otimes A$ par l'idéal engendré par les éléments $\alpha(x)\otimes 1-1\otimes x^p$, quotient qui n'est autre que $K\,_\alpha{\underset{A}{\otimes}}{}_{\mathrm{Fr}} A$.

Si $p=2$, comme $\alpha$ est additif, $1\otimes a\otimes 1+1\otimes 1\otimes a$ appartient à $\I_2(A,\alpha)$ pour tout $a\in A$, de sorte que $W_2(A,\alpha)$ s'identifie au quotient de $K\otimes_\mathbb{Z} A$ par l'idéal engendré par les éléments $1\otimes a^2+\alpha(a)\otimes 1$, d'où le résultat. 
\end{proof}

\begin{exem}\label{ex-Wcorimp} Supposons que $A$ est un corps de caractéristique $p>0$ de degré infini sur l'image $A^p$ de son endomorphisme de Frobenius. On vérifie facilement que l'anneau $A\,_{\mathrm{Fr}}{\underset{A}{\otimes}}{}_{\mathrm{Fr}}A$ n'est pas semi-primaire (en utilisant une \emph{$p$-base} de $A$ sur  $A^p$ --- cf. \cite[chap.~V, §\,13, déf.~1 et th.~2]{Bki2}), il en est donc de même pour $W_p(A,\alpha)$, si $\alpha : A\to K$ est un morphisme de corps avec $K$ parfait.
\end{exem}

\begin{coro}\label{cor-Wpgros} Supposons que $p$ est premier et que $A$ est une $\FF_p$-algèbre augmentée dont l'endomorphisme de Frobenius est nul sur l'idéal d'augmentation. Notons $\alpha\in\mathbf{Ann}(A,K)$ la composée de l'augmentation $A\to\FF_p$ et du morphisme de corps $\FF_p\to K$. Alors il existe un morphisme surjectif de $K$-algèbres $W_p(A,\alpha)\twoheadrightarrow A_K$, qui est un isomorphisme si $p=2$.
\end{coro}

\begin{exem}\label{ex-Wbiz} Soit $P$ l'ensemble des parties vides ou infinies de $\mathbb{N}$. Munissons le $\FF_p$-espace vectoriel $A:=\FF_p[P]$ d'une structure de $\FF_p$-algèbre à l'aide de la multiplication donnée par $[I].[J]=[I\sqcup J]$ si  $I$ et $J$ sont des éléments disjoints de $P$, et $[I].[J]=0$ sinon. L'application linéaire $A\to\FF_p$ envoyant $[I]$ sur $0$ pour $I\ne\varnothing$ (et $[\varnothing]=1_A$ sur $1$) est une augmentation sur $A$ ; $A$ est un anneau local dont le radical est l'idéal d'augmentation, noté $\mathfrak{m}$. Comme $\mathfrak{m}$ est engendré linéairement par des éléments de carré nul, l'endomorphisme de Frobenius de $A$ est nul sur $\mathfrak{m}$. On remarque que $\mathfrak{m}^2=\mathfrak{m}$ (car toute partie infinie de $\mathbb{N}$ est la réunion disjointe de deux parties infinies), d'où l'on déduit que $(A,K)$ vérifie la condition (CAD)$^+$. Ainsi, tout quotient semi-primaire de $W_p(A,\alpha)$ (où $\alpha$ est défini comme dans le corollaire~\ref{cor-Wpgros}) est isomorphe à $K$, par le corollaire~\ref{cor-qspfW}. Mais l'anneau quasi-parfait $W_p(A,\alpha)$ n'est pas semi-primaire, ni même parfait. En effet, par le corollaire~\ref{cor-Wpgros}, il suffit de vérifier la même propriété pour l'anneau $A$, ce qui provient de l'observation suivante : si $I_l$ désigne, pour $l\in\mathbb{N}$ premier, l'ensemble des puissances de $l$ autres que $1$, alors $([I_l])$ est une suite de $\mathfrak{m}$ dont aucune sous-suite finie n'est de produit nul.
\end{exem}

\begin{rema}\label{rq-I2I} Poursuivons l'analyse de l'exemple précédent, en supposant pour simplifier $p=2$ et $K=\FF_2$ (cette dernière condition n'étant pas restrictive, quitte à étendre les scalaires). Tout d'abord, de manière générale, si $A$ est une $\FF_2$-algèbre munie d'une augmentation $\alpha : A\to\FF_2$ dont le noyau $\mathfrak{m}$ est constitué d'éléments de carré nul, alors l'idéal $\I_2(A,\alpha)$ de $A\otimes A$ est engendré par les éléments $x\otimes 1+1\otimes x$ pour $x\in\mathfrak{m}$. En effet, un tel élément appartient à $\I_2(A,\alpha)$ puisqu'il égale la somme de $(x+1)^{\otimes 2}+\alpha(x+1)$ et $x^{\otimes 2}+\alpha(x)$, et la relation
$$t\otimes t+\alpha(t)=(t\otimes 1).(\bar{t}\otimes 1+1\otimes\bar{t})+t^2\otimes 1+\alpha(t)=(t\otimes 1).(\bar{t}\otimes 1+1\otimes\bar{t})$$
pour $t\in A$, où $\bar{t}:=t+\alpha(t)$ appartient à $\mathfrak{m}$, montre que l'idéal engendré par les éléments de ce type contient $\I_2(A,\alpha)$ (on a $t^2=\alpha(t)$ car $\bar{t}^2=0$ par hypothèse).

Vérifions que, dans la situation de l'exemple~\ref{ex-Wbiz} (avec $p=2$), l'idéal $\I_2(A,\alpha)$ de $A\otimes A$ n'est pas égal à son carré. Si $E$ désigne l'ensemble des parties infinies de $\mathbb{N}$, les $x_I:=[I]$, pour $I\in E$, constituent une base du $\FF_2$-espace vectoriel $\mathfrak{m}$. Dans l'algèbre $R:=A\otimes A$, notons $a_I:=x_I\otimes 1$ et $b_I:=1\otimes x_I$ pour $I\in E$ : $1$, les $a_I$, les $b_J$ et les $a_I b_J$, pour $(I,J)\in E^2$, constituent une base de $R$ sur $\FF_2$. Considérons l'application $\FF_2$-linéaire $\psi : R\to A$ envoyant $1$ sur $1$, $a_I$ sur $x_I$, $b_J$ sur $0$, et $a_I b_J$ sur $x_{I\sqcup J}$ si $I\cap J=\varnothing$ et $\min I<\min J$ (pour l'ordre usuel de $\mathbb{N})$ et sur $0$ sinon. Alors $\psi$ n'est pas nulle sur $\I_2(A,\alpha)$, car $a_I+b_I\in\I_2(A,\alpha)$ est envoyé sur $x_I\ne 0$. En revanche, $\psi$ est nul sur l'idéal $\I_2(A,\alpha)^2$ de $R$ : cet idéal est engendré par les éléments $(a_I+b_I)(a_J+b_J)$ (pour $(I,J)\in E^2$) ; une famille génératrice de $\I_2(A,\alpha)^2$ comme $\FF_2$-espace vectoriel est donc constituée des éléments $(a_I+b_I)(a_J+b_J)$, $a_T(a_I+b_I)(a_J+b_J)$, $b_U(a_I+b_I)(a_J+b_J)$ et $a_T b_U(a_I+b_I)(a_J+b_J)$ pour $(I,J,T,U)\in E^4$. Montrons par exemple que $\psi$ annule les éléments de la forme $a_T b_U(a_I+b_I)(a_J+b_J)$ (les autres cas sont similaires) : si $I$, $J$, $T$ et $U$ ne sont pas deux à deux disjoints, alors $\psi(a_T b_U(a_I+b_I)(a_J+b_J))=0$ car $\psi$ envoie sur $0$ tous les termes obtenus en développant le produit. Si $I$, $J$, $T$ et $U$ sont deux à deux disjoints, alors
$$a_T b_U(a_I+b_I)(a_J+b_J)=a_{I\sqcup J\sqcup T}b_U+a_{I\sqcup T}b_{J\sqcup U}+a_{J\sqcup T}b_{I\sqcup U}+a_T b_{I\sqcup J\sqcup U}$$
fait apparaître un nombre pair de termes dans lesquels $I$, $J$, $T$ ou $U$ figure en indice de $a$, de sorte qu'on a bien $\psi(a_T b_U(a_I+b_I)(a_J+b_J))=0$.

On a donc bien $\I_d(A,\alpha)^2\ne\I_2(A,\alpha)$.
\end{rema}

\section{Lien avec les effets croisés}\label{sliencr}

Comme la catégorie $\A$ (resp. $\E$) est $A$-linéaire (resp. $K$-linéaire), la catégorie de multifoncteurs multiadditifs $\mathbf{Add}_d(\A,\E)$ est $A^{\otimes_\mathbb{Z}d}\otimes_\mathbb{Z}K\simeq A_K^{\otimes d}$-linéaire, de même que la catégorie $\mathbf{Add}_d(\A,\E)\rtimes\Si_d$.

\begin{nota}
On note $\pol_d(\A,\E)_\pi^{[r]}$ la sous-catégorie pleine de $\fct(\A,\E)$ constituée des foncteurs polynomiaux de degré au plus $d$ de $U(\fct(\A,\E);\mathfrak{m}_\pi^r)$.

Lorsque $\A=\mathbf{P}(A)$ et $\E=K\Md$, on note simplement $\pol_d(A,K)_\pi^{[r]}$ cette catégorie.
\end{nota}

Le résultat suivant, qui justifie les considérations de la section~\ref{swd},  est immédiat.
\begin{prop}\label{pr-Wcrimm} 
\begin{enumerate}
\item Le foncteur $\mathbf{Add}_d(\A,\E)\to\pol_d(\A,\E)$ de précomposition par la diagonale $d$-itérée envoie $U(\mathbf{Add}_d(\A,\E);\I_d(A,\pi)^r)=\mathbf{Add}_d(\A,\E)\underset{A_K^{\otimes d}}{\boxtimes}W_d^{[r]}(A,\pi)$ dans $\pol_d(\A,\E)_\pi^{[r]}$.
\item Le foncteur $\pol_d(\A,\E)\to\mathbf{Add}_d(\A,\E)$ qu'induit $cr_d$ envoie $\pol_d(\A,\E)_\pi^{[r]}$ dans $U(\mathbf{Add}_d(\A,\E);\I_d(A,\pi)^r)$.
\end{enumerate}
\end{prop}

Avant de donner, à la section~\ref{parfpp}, les conséquences principales des résultats de la section~\ref{swd} sur les foncteurs polynomiaux à poids déduites de l'observation précédente, nous abordons la question générale de savoir quand toute décomposition en poids faible pour un foncteur polynomial de $\fct(\A,\E)$ est automatiquement forte, obtenant à la proposition~\ref{pr-tppfor} ci-après un analogue partiel de la proposition~\ref{pr-pforfaib} pour les foncteurs \emph{polynomiaux}.

\begin{coro}\label{cor-wId} Considérons les assertions suivantes :
\begin{enumerate}
\item\label{eqitw1} pour toute petite catégorie additive $A$-linéaire $\A$ et toute catégorie de Grothendieck $K$-linéaire $\E$, tout foncteur polynomial homogène de poids faible $\pi$ de $\fct(\A,\E)$ est homogène de poids fort $\pi$ ;
\item\label{eqitw2} tout foncteur polynomial de $\F(A,K)$ homogène de poids faible $\pi$ est de poids fort $\pi$ ;
\item\label{eqitw3} pour tout $d\in\mathbb{N}$, l'idéal $\I_d(A,\pi)$ de $A_K^{\otimes d}$ est égal à son carré.
\end{enumerate}
Alors \ref{eqitw1}$\Rightarrow$\ref{eqitw2}$\Rightarrow$\ref{eqitw3}, et \ref{eqitw3}$\Rightarrow$\ref{eqitw1} si $p=0$.
\end{coro}

\begin{proof} L'implication \ref{eqitw1}$\Rightarrow$\ref{eqitw2} est triviale.

Utilisant l'équivalence du corollaire~\ref{cor-recPira} (dû à Pirashvili), la sous-catégorie $\pol_d(A,K)_\pi^{[r]}/\pol_{d-1}(A,K)_\pi^{[r]}$ de $\pol_d(A,K)/\pol_{d-1}(A,K)$ s'identifie à $W_d^{[r]}(A,\pi)\rtimes\Si_d\Md$, par la proposition~\ref{pr-Wcrimm}. Ainsi, la condition~\ref{eqitw2} implique que le foncteur $W_d(A,\pi)\rtimes\Si_d\Md\to W_d^{[2]}(A,\pi)\rtimes\Si_d\Md$ de restriction le long de la projection $W_d^{[2]}(A,\pi)\rtimes\Si_d\twoheadrightarrow W_d(A,\pi)\rtimes\Si_d=(A_K^{\otimes d}/\I_d(A,\pi)^2\twoheadrightarrow A_K^{\otimes d}/\I_d(A,\pi))\rtimes\Si_d$ est une équivalence, c'est-à-dire que l'inclusion $\I_d(A,\pi)^2\subset\I_d(A,\pi)$ d'idéaux de $A_K^{\otimes d}$ est une égalité. Ainsi \ref{eqitw2}$\Rightarrow$\ref{eqitw3}.

Si la condition \ref{eqitw3} est vérifiée et que $p$ est nul, alors la proposition~\ref{pr-Wcrimm} et le scindement par le degré des foncteurs polynomiaux de but $\mathbb{Q}$-linéaire (proposition~\ref{pr-decpolcar0}.\ref{itdQrec}) montrent que l'inclusion $\pol_d(\A,\E)_\pi^{[r]}\subset\pol_d(\A,\E)_\pi^{[1]}$ est une égalité, d'où \ref{eqitw1}.
\end{proof}

\begin{rema} Nous ignorons si l'implication \ref{eqitw3}$\Rightarrow$\ref{eqitw1} persiste en caractéristique $p$ première.
\end{rema}

On observe par ailleurs (à partir de l'équivalence~\eqref{eq-pteAdda}, page~\pageref{eq-pteAdda}) que la condition $\mathrm{Der}_\varphi(A,K)=0$, pour $\varphi\in\mathbf{Ann}(A,K)$, équivaut au fait qu'un foncteur additif de $\F(A,K)$ homogène de poids faible $\varphi$ est nécessairement homogène de poids fort $\varphi$, et implique la même propriété pour les foncteurs additifs de $\fct(\A,\E)$.

Ainsi, la condition (CAD) équivaut au fait qu'un poids d'un foncteur additif de $\F(A,K)$ est nécessairement fort, mais elle ne suffit pas à garantir la même propriété pour les foncteurs polynomiaux de degré supérieur (c'est toutefois le cas si $K$ est assez gros, par exemple algébriquement clos, et de caractéristique nulle, comme il résultera de la proposition~\ref{pr-tppfor}).

\begin{prop}\label{pr-tppfor} Si toute décomposition en poids d'un foncteur polynomial de $\F(A,K)$ est forte, alors la condition \textnormal{(CAD)}$^+$ est vérifiée.\index{nota}{CAD@(CAD), (CAD)$^+$ \emph{(conditions d'annulation des dérivations)}}

Réciproquement, si \textnormal{(CAD)}$^+$ est vérifiée et que la caractéristique $p$ de $K$ est nulle, alors toute décomposition en poids d'un foncteur polynomial de $\fct(\A,\E)$ est forte.
\end{prop}

\begin{proof} Supposons que toute décomposition en poids d'un foncteur polynomial de $\F(A,K)$ est forte. Soient $L$ une extension finie du corps $K$ et $\varphi\in\mathbf{Ann}(A,L)$. Choisissons une extension finie $L'$ du corps $L$ et $\alpha_1,\dots,\alpha_n\in\mathbf{Ann}(A,L')$ comme dans le lemme~\ref{lm-extpgg}, et notons $\pi$ le produit des $\alpha_i$, qu'on peut voir comme élément de $\mon(A_\mu,K_\mu)$. 

Supposons $\mathrm{Der}_\varphi(A,L)$ non nul : il existe alors un foncteur additif homogène de poids faible non fort $\varphi$ dans $\F(A,L)$, et donc dans $\F(A,L')$ ; soit $X_1$ un tel foncteur. Par ailleurs, notons $X_i$, pour $i\in\llbracket 2,n\rrbracket$ le foncteur de $\F(A,L')$ d'extension des scalaires le long de $\alpha_i$ ; c'est un foncteur additif homogène de poids $\alpha_i$, à valeurs toujours non nulles, sauf en $0$. L'image du foncteur $\bigotimes_{i=1}^n X_i$ dans $\F(A,K)$ par la postcomposition par la restriction des scalaires de $L'$ à $K$ est polynomiale, homogène de poids faible $\pi$. Comme $X_1$ n'est pas homogène de poids \emph{fort} $\alpha_1$ et que les $X_i$, pour $i\ge 2$, sont à valeurs non nulles hors de l'objet nul, il s'ensuit que notre foncteur ne peut être homogène de poids \emph{fort} $\pi$. Cette contradiction montre que (CAD)$^+$ est nécessairement vérifiée.

Supposons maintenant (CAD)$^+$ vérifiée et $p$ nul. Alors le corollaire~\ref{cor-wssder} entraîne que l'idéal $\I_d(A,\pi)$ de $A_K^{\otimes d}$ est égal à son carré. Le corollaire~\ref{cor-wId} permet de conclure.
\end{proof}

\begin{rema} Le sens réciproque de l'énoncé précédent peut tomber en défaut si l'on omet l'hypothèse de caractéristique nulle, comme l'illustre le cas de $A=\mathbb{Z}/p^2$ (qui vérifie (CAD)$^+$ alors qu'il existe des foncteurs de $\F(\mathbb{Z}/p^2,K)$ admettant une décomposition en poids faible mais non forte, par les propositions \ref{pr-decpfa} et \ref{pr-decpfoap}).

Nous verrons toutefois ultérieurement (corollaire~\ref{cor-polpoiforgratp}) que si $p$ est non nul et que $A$ est une \emph{$\FF_p$-algèbre} vérifiant une hypothèse de finitude raisonnable (par exemple, $A$ de type fini comme anneau) et la condition (CAD)$^+$, alors toute décomposition en poids d'un foncteur polynomial de $\F(A,K)$ est nécessairement forte.

Néanmoins, même si $A$ est une $\FF_p$-algèbre, sans hypothèse de finitude supplémentaire, la condition (CAD)$^+$ ne suffit généralement pas à garantir que toute décomposition en poids d'un foncteur polynomial de $\F(A,K)$ soit forte, comme le montre l'exemple~\ref{ex-Wbiz}, combiné à la remarque~\ref{rq-I2I} et au corollaire~\ref{cor-wId}.
\end{rema}

\section{Applications aux foncteurs polynomiaux à poids}\label{parfpp}

Commençons par donner une conséquence générale simple des résultats rappelés dans la section~\ref{s-rppol} :
\begin{lemm}\label{lm-crdtfn}
\begin{enumerate}
\item Les foncteurs $cr_d : \fct(\A,\E)\to\fct(\A^d,\E)$ et $\delta_d^* : \fct(\A^d,\E)\to\fct(\A,\E)$ préservent les objets de type fini.
\item Si $F$ est un objet noethérien de $\pol_d(\A,\E)$, alors $cr_d(F)$ est un objet noethérien de $\mathbf{Add}_d(\A,\E)$.
\end{enumerate}
\end{lemm}

\begin{proof}
Le premier point résulte formellement de ce que les foncteurs $\oplus_d^*$, dont $cr_d$ est facteur direct, et $\delta_d^*$ sont adjoints \emph{des deux côtés} (en particulier, ils sont exacts et commutent aux colimites) et de ce que la type-finitude d'un objet $X$ d'une catégorie de Grothendieck équivaut à la commutation de $\operatorname{Hom}(X,-)$ aux colimites filtrantes de monomorphismes.

Si $F$ est un objet noethérien de $\pol_d(\A,\E)$, alors \cite[chapitre~5, lemme~8.3]{Pop} et le théorème~\ref{th-recPiranv} de Pirashvili montrent que $cr_d(F)$ est un objet noethérien de $\mathbf{Add}_d(\A,\E)\rtimes\Si_d$.  Comme $\Si_d$ est un groupe \emph{fini}, il s'ensuit, par la proposition~\ref{pr-onpt}, que $cr_d(F)$ est noethérien dans $\mathbf{Add}_d(\A,\E)$.
\end{proof}

Dans la suite de ce chapitre, nous nous plaçons sur la catégorie source $\mathbf{P}(A)$, sur laquelle les foncteurs (multi)additifs se décrivent bien.

\begin{prop}\label{pr-vppoi} Soit $F$ un foncteur de $\pol_d(\mathbf{P}(A),\E)^{[r]}_\pi$.
\begin{enumerate}
\item Supposons que $F$ est de type fini et que la $K$-algèbre $W_d^{[r]}(A,\pi)$ est de dimension finie. Alors le multifoncteur $cr_d(F)$ est à valeurs dans les objets de type fini de $\E$.
\item Supposons que $F$ est noethérien. Alors le multifoncteur $cr_d(F)$ est à valeurs dans les objets noethériens de $\E$.
\end{enumerate}
\end{prop}

\begin{proof} L'équivalence de catégories~\eqref{eq-pteAddd} (page~\pageref{eq-pteAddd}), la proposition~\ref{pr-Wcrimm} et le lemme~\ref{lm-crdtfn} montrent que $cr_d(F)$ est un objet de type fini dans le premier cas, noethérien dans le second, de $\mathbf{Add}_d(\mathbf{P}(A),\E)\underset{A_K^{\otimes d}}{\boxtimes}W_d^{[r]}(A,\pi)\simeq\E\underset{K}{\boxtimes}W_d^{[r]}(A,\pi)$. Si $W_d^{[r]}(A,\pi)$ est de dimension finie sur $K$, il définit donc un objet de type fini de $\E$ ; autrement dit, le multifoncteur $cr_d(F)$ est à valeurs dans les objets de type fini de $\E$. Si $F$ est noethérien, on conclut grâce au corollaire~\ref{cor-noethp}.
\end{proof}

\begin{rema} De la même façon, si la $K$-algèbre $A_K$ est absolument plate et déployée (ou plus généralement s'il existe une extension finie de $L$ de $K$ telle que $A_L$ soit déployée), alors pour tout foncteur polynomial noethérien $F$ de $\fct(\mathbf{P}(A),\E)$, $cr_d(F)$ est à valeurs noethériennes. Il suffit pour le voir de raisonner comme précédemment, en utilisant la proposition~\ref{poivN} (et le corollaire~\ref{cor-ntcf}) au lieu du corollaire~\ref{cor-noethp}.

Ainsi, la deuxième assertion du théorème~\ref{th-poinoeth_abst} ci-dessous vaut pour tous les foncteurs polynomiaux noethériens (sans condition de décomposition en poids) sous l'hypothèse précédente sur $A$ et $K$.
\end{rema}

L'un des résultats principaux de ce chapitre est le suivant.

\begin{theo}\label{th-poinoeth_abst} Soit $F$ un foncteur de $\pol_d(\mathbf{P}(A),\E)_\pi^{[r]}$.
\begin{enumerate}
\item Si $F$ est de type fini et que pour tout entier $n\le d$, l'algèbre $W_n^{[r]}(A,\pi)$ est de dimension finie sur $K$, alors $F$ est à valeurs dans les objets de type fini de $\E$.
\item Si $F$ est noethérien, alors $F$ prend ses valeurs dans les objets noethériens de $\E$.
\end{enumerate}
\end{theo}

\begin{proof} On montre le théorème par récurrence sur $d$. Le cas $d=0$ étant trivial, on suppose $d>0$ et les assertions établies pour les foncteurs de degré au plus $d-1$.

Le foncteur $\qpol_{d-1}(F)$ est de type fini dans le premier cas, noethérien dans le second, car c'est un quotient de $F$. Il appartient par ailleurs à $\pol_{d-1}(\mathbf{P}(A),\E)_\pi^{[r]}$. Ainsi, l'hypothèse de récurrence montre qu'il est à valeurs dans les objets de type fini de $\E$ dans le premier cas, noethériens dans le second.

Il en est de même pour $\delta^*_d cr_d(F)$ par la proposition~\ref{pr-vppoi}. Par conséquent, la suite exacte $\delta_d^*cr_d(F)\to F\to\qpol_{d-1}(F)\to 0$ (cf. \eqref{eq-qpol}, page~\pageref{eq-qpol}) achève la démonstration.
\end{proof}

L'énoncé précédent s'étend aussitôt aux foncteurs \phs\,:
\begin{coro}\label{cor-phpoinoeth} Soit $F$ un foncteur \ph\ de type fini de $\fct(\mathbf{P}(A),\E)$, homogène de poids faible $\pi$.
 Supposons que pour tous entiers $n, r>0$, l'algèbre $W_n^{[r]}(A,\pi)$ est de dimension finie sur $K$. Alors $F$ est à valeurs dans les objets de type fini de $\E$.
\end{coro}

\begin{coro}\label{cor-noethvn} Tout foncteur \ph\ noethérien de $\fct(\mathbf{P}(A),\E)$ possédant une décomposition en poids faible est à valeurs noethériennes.
\end{coro}

\begin{proof}[Démonstration des corollaires~\ref{cor-phpoinoeth} et \ref{cor-noethvn}] Comme toute décomposition en poids faible d'un foncteur de type fini est finie, il suffit de démontrer que si $F$ est un foncteur \ph\ de type fini de $U(\fct(\mathbf{P}(A),\E);\mathfrak{m}_\pi^r)$, et que les $W_n^{[r]}(A,\pi)$ sont de dimension finie (resp. que $F$ est noethérien), alors $F$ est à valeurs de type fini (resp. noethériennes).

Pour tout $d\in\mathbb{N}$, le foncteur $\qpol_d(F)$ est polynomial, de type fini (resp. noethérien), et appartient à $\pol_d(\mathbf{P}(A),\E)_\pi^{[r]}$, il est donc à valeurs de type fini (resp. noethériennes) par le théorème~\ref{th-poinoeth_abst}. La conclusion résulte alors de la proposition~\ref{pr-ppht}.
\end{proof}

Comme l'évaluation sur $A^d$ définit un foncteur exact et \emph{fidèle} $\pol_d(\mathbf{P}(A),\E)\to\E$, elle détecte les objets noethériens et les objets finis (cf. \cite[lm~11.10]{DTV} et \cite[lm~6.20]{DT-schw}). Par conséquent, le théorème~\ref{th-poinoeth_abst} entraîne :

\begin{coro}\label{cor-poinoeth_abst} Soit $F$ un foncteur de $\pol_d(\mathbf{P}(A),\E)_\pi^{[r]}$.
\begin{enumerate}
\item Supposons le foncteur $F$ noethérien et la catégorie $\E$ localement finie. Alors $F$ est fini (dans $\fct(\mathbf{P}(A),\E)$).
\item Supposons que $F$ est de type fini, que la catégorie $\E$ est localement noethérienne (resp. localement finie) et que pour tout entier $n\le d$, l'algèbre $W_n^{[r]}(A,\pi)$ est de dimension finie sur $K$. Alors $F$ est noethérien (resp. fini).
\end{enumerate}
\end{coro}

Dans le cas particulier fondamental de la catégorie $\F(A,K)$, le théorème~\ref{th-poinoeth_abst} et le corollaire~\ref{cor-poinoeth_abst} fournissent :

\begin{theo}\label{th-poinoeth} Soit $F$ un foncteur de $\pol_d(A,K)_\pi^{[r]}$. 

Supposons que l'une des conditions suivantes est satisfaite :
\begin{enumerate}
\item $F$ est noethérien ;
\item pour tout entier $n\le d$, l'algèbre $W_n^{[r]}(A,\pi)$ est de dimension finie sur $K$, et $F$ est de type fini.
\end{enumerate}

Alors $F$ est fini et à valeurs de dimension finie.
\end{theo}

\begin{rema} Ce théorème dit en particulier qu'un foncteur polynomial \emph{simple} à poids de $\F(A,K)$ est à valeurs de dimension finie, ce qui permet de lui appliquer les résultats de \cite[§\,5]{DTV}.

Ce résultat peut également se déduire des propositions~\ref{pr-simpoldfext} et~\ref{pr-dppol}. Inversement, le cas particulier du théorème~\ref{th-poinoeth} (dont la démonstration n'utilise pas la proposition~\ref{pr-simpoldfext}) où $F$ est simple donne une nouvelle démonstration de la proposition~\ref{pr-simpoldfext}, car si $S$ est un foncteur polynomial simple de $\F(A,K)$, l'action tautologique sur $S$ du surcorps $L:=Z(\mathrm{End}(S))$ de $K$ fournit un foncteur $\tilde{S}$ de $\F(A,L)$ polynomial, simple, et qui possède un poids par le lemme~\ref{pd1ev}.  
\end{rema}

En caractéristique $p=0$, la structure des foncteurs polynomiaux à poids se déduit plus directement de celle des algèbres $W_n^{[r]}(A,\pi)$, car le scindement de la filtration polynomiale (proposition~\ref{pr-decpolcar0}.\ref{itdQrec}), combiné au corollaire~\ref{cor-recPira} de Pirashvili et à la proposition~\ref{pr-Wcrimm}, fournit une équivalence de catégories
\begin{equation}\label{eq-scW}
\pol_d(A,K)_\pi^{[r]}\simeq\prod_{n=0}^d(W_n^{[r]}(A,\pi)\rtimes\Si_n)\Md\quad(\text{si }p=0).
\end{equation}

\begin{coro}\label{cor-fppc0ss} En caractéristique $p=0$, les assertions suivantes sont équivalentes :
\begin{enumerate}
\item la condition  \textnormal{(CAD)}$^+$\index{nota}{CAD@(CAD), (CAD)$^+$ \emph{(conditions d'annulation des dérivations)}} est vérifiée ;
\item tout foncteur analytique de $\F(A,K)$ possédant une décomposition en poids faible est semi-simple.
\end{enumerate}
\end{coro}

\begin{proof}
Cela résulte du corollaire~\ref{cor-wssder}, de la proposition~\ref{pr-psspsd} et de l'équivalence~\eqref{eq-scW}.
\end{proof}

\begin{coro}\label{cor-fppc0lf} En caractéristique $p=0$, les assertions suivantes sont équivalentes :
\begin{enumerate}
\item\label{itc01} la condition \textnormal{(CFD)}$^+$ \index{nota}{CFD@(CFD), (CFD)$^+$ \emph{(conditions de finitude sur les dérivations)}} est vérifiée ;
\item\label{itc02} tout foncteur analytique de $\F(A,K)$ possédant une décomposition en poids faible est localement fini ;
\item\label{itc03} tout foncteur analytique de $\F(A,K)$ possédant une décomposition en poids faible est localement noethérien ;
\item\label{itc04} tout foncteur polynomial de type fini de $\F(A,K)$ possédant une décomposition en poids faible est à valeurs de dimension finie.
\end{enumerate}
\end{coro}

\begin{proof}
Le corollaire~\ref{cor-wssder} et l'équivalence~\eqref{eq-scW} montrent que les conditions \ref{itc01} et \ref{itc04} sont équivalentes. L'implication \ref{itc04}$\Rightarrow$\ref{itc02} découle d'un argument d'évaluation (cf. corollaire~\ref{cor-poinoeth_abst}) ; l'implication \ref{itc02}$\Rightarrow$\ref{itc03} est évidente. L'implication \ref{itc03}$\Rightarrow$\ref{itc04} résulte du théorème~\ref{th-poinoeth}. 
\end{proof}

Nous spécifions maintenant des conséquences directes du théorème~\ref{th-poinoeth_abst} et de résultats antérieurs de ce chapitre en caractéristique $p$ première.

\begin{coro}\label{cor-poids-ln_gal} Supposons que $p$ est premier et que l'anneau $A/p$ est une algèbre essentiellement de type fini sur une $\FF_p$-algèbre $p$-parfaite\index{termin}{essentiellement de type fini \emph{(anneau, algèbre)}}.
\begin{enumerate}
\item Tout foncteur \ph\ de type fini de $\fct(\mathbf{P}(A),\E)$ possédant une décomposition en poids faible est à valeurs de type fini (dans $\E$).
\item Si la catégorie $\E$ est localement noethérienne (resp. localement finie), alors tout foncteur polynomial de $\fct(\mathbf{P}(A),\E)$ possédant une décomposition en poids faible est localement noethérien (resp. localement fini).
\end{enumerate}
\end{coro}

\begin{proof} Cela résulte de la proposition~\ref{pr-Wdf} et des corollaires~\ref{cor-phpoinoeth} et ~\ref{cor-poinoeth_abst}.
\end{proof}

En particulier, pour $\E=K\Md$, on a les résultats suivants.

\begin{coro}\label{cor-phdf} Supposons que $p$ est premier et que l'anneau $A/p$ est une algèbre essentiellement de type fini sur une $\FF_p$-algèbre $p$-parfaite. Alors tout foncteur \ph\ de type fini de $\F(A,K)$ possédant une décomposition en poids faible est à valeurs de dimensions finies.
\end{coro}

\begin{coro}\label{cor-poids-ln} Supposons que $p$ est premier et que l'anneau $A/p$ est une algèbre essentiellement de type fini sur une $\FF_p$-algèbre $p$-parfaite\index{termin}{essentiellement de type fini \emph{(anneau, algèbre)}}. Alors tout foncteur polynomial de $\F(A,K)$ possédant une décomposition en poids faible est localement fini et localement à valeurs de dimensions finies.
\end{coro}

\begin{rema} On dispose de la réciproque partielle suivante au corollaire~\ref{cor-poids-ln}.

Supposons la $\FF_p$-algèbre $A/p$ de type fini (ici, essentiellement de type fini ne suffit pas, et l'on ne peut pas se contenter de supposer que $A/p$ est de type fini sur un corps parfait arbitraire). Alors tout foncteur polynomial simple de $\F(A,K)$ est à valeurs de dimension finie \cite[proposition~7.10]{DTV} (ladite proposition suppose que l'anneau $A$ est de type fini, mais la même démonstration fonctionne dès que la $\FF_p$-algèbre $A/p$ est de type fini), donc possède un poids, au moins si $K$ est algébriquement clos (proposition~\ref{poisimp}), ou plus généralement si $K$ est un corps de décomposition de la catégorie additive $\mathbf{P}(A)$ (proposition~\ref{pr-vdfdppol}). Tout foncteur polynomial localement fini de $\F(A,K)$ possède alors une décomposition en poids faible.
\end{rema}

\begin{exem} Si $k$ est un corps de caractéristique $p>0$ qui est une extension de type fini d'un sous-corps parfait, alors tout foncteur polynomial de type fini de $\F(k,K)$ possédant une décomposition en poids faible est fini et à valeurs de dimension finie.
\end{exem}

\section[Finitude des dimensions prises par un foncteur polynomial]{La condition de finitude des dimensions des valeurs d'un foncteur polynomial de $\F(A,K)$}\label{sfdp}
 
Nous commençons cette section par un résultat qui n'utilise pas les résultats antérieurs du présent chapitre, mais qui possède un intérêt intrinsèque et met en perspective la proposition~\ref{pr-poipoldfgra} ci-après.

\begin{prop}\label{pr-tspollf} Les assertions suivantes sont équivalentes : 
\begin{enumerate}
    \item\label{it-pris1} $\dim_K A_K<\infty$ ;
    \item\label{it-pris2} tout foncteur additif de type fini de $\F(A,K)$ est à valeurs de dimensions finies ;
    \item\label{it-pris3} tout foncteur polynomial de type fini de $\F(A,K)$ est à valeurs de dimensions finies ;
     \item\label{it-pris4} tout foncteur polynomial de type fini de $\F(A,K)$ est fini ;
    \item\label{it-pris5} tout foncteur polynomial de type fini de $\F(A,K)$ est artinien ;
    \item\label{it-pris6} tout foncteur polynomial de type fini de degré au plus $2$ de $\F(A,K)$ est artinien ;
    \item\label{it-pris7} les $K$-algèbres $A_K$ et $A_K^{\otimes 2}$ sont artiniennes.
\end{enumerate}
\end{prop}

\begin{proof} L'équivalence entre \ref{it-pris1} et \ref{it-pris2} est immédiate (cf. \eqref{eq-pteAdda}, page~\pageref{eq-pteAdda}).

La proposition~\ref{pr-decpolcar0} et le corollaire~\ref{cor-recPira} montrent l'implication  \ref{it-pris1}$\Rightarrow$\ref{it-pris3} lorsque $K$ est de caractéristique nulle. En caractéristique $p>0$, la condition \ref{it-pris1} signifie que l'anneau $A/p$ est fini, ce qui entraîne que $A/p^n$ est fini pour tout $n\in\mathbb{N}$. Comme tout foncteur polynomial de $\F(A,K)$ se factorise par la réduction modulo $p^n$ pour $n\in\mathbb{N}$ assez grand (ce fait classique et élémentaire est redémontré au corollaire~\ref{cor-radPolp} ci-après, sans utiliser la présente proposition), cela montre \ref{it-pris1}$\Rightarrow$\ref{it-pris3} dans le cas de la caractéristique première.

L'implication \ref{it-pris3}$\Rightarrow$\ref{it-pris4} résulte de \cite[lemme~11.10]{DTV} ; les implications \ref{it-pris4}$\Rightarrow$\ref{it-pris5}$\Rightarrow$\ref{it-pris6} sont triviales. 

Si la condition \ref{it-pris6} est vérifiée, alors le corollaire~\ref{cor-recPira} et \cite[chap.~5, cor.~8.4]{Pop} montrent que les algèbres $A_K$ et $A_K\wr\Si_2$ (cette dernière étant non commutative) sont artiniennes à gauche. Comme $\Si_2$ est un groupe fini, on en déduit que $A_K^{\otimes 2}$ est également artinienne par \cite[thm~2]{Grz85} (ou la proposition~\ref{pr-psspsd}), d'où l'implication \ref{it-pris6}$\Rightarrow$\ref{it-pris7}.

Enfin, \ref{it-pris7} entraîne \ref{it-pris1} grâce au lemme~\ref{lm-artDF} ci-dessous, d'où la proposition.
\end{proof}

\begin{lemm}\label{lm-artDF} Soit $R$ une $K$-algèbre. Si les algèbres $R$ et $R\otimes_K R$ sont artiniennes, alors $\dim_K R$ est finie.
\end{lemm}

\begin{proof} Il est classique que, si $K\subset L$ est une extension de corps, alors l'anneau $L\otimes_K L$ n'est artinien que si cette extension est \emph{finie} (\cite{Sharp} montre que l'extension est nécessairement algébrique, et \cite[th.~11]{Vamos} qu'elle est nécessairement de type fini). Par conséquent, si $\mathfrak{m}$ est un idéal maximal d'une $K$-algèbre $R$ telle que $R\otimes_K R$ soit artinienne, alors $R/\mathfrak{m}$ est une extension finie de $K$. Si de plus l'algèbre $R$ est artinienne, cela entraîne que $\dim_K (R/\rj(R))<\infty$. Comme $\rj(R)$ est nilpotent et que les $R/\rj(R)$-modules $\rj(R)^n/\rj(R)^{n+1}$ sont de type fini pour tout $n\in\mathbb{N}$, cela démontre le lemme.
\end{proof}

Nous utilisons maintenant les résultats de la section~\ref{parfpp} pour établir le lien important suivant entre l'existence d'une décomposition en poids et des propriétés de finitude fortes pour un foncteur polynomial de $\F(A,K)$, sous une hypothèse de finitude sur $A$.

\begin{prop}\label{pr-poipoldfgra} Supposons que $K$ est un corps de décomposition de $\mathbf{P}(A)$ et que la $K$-algèbre $A_K$ est de type fini. Soit $F$ un foncteur polynomial de $\F(A,K)$. Les assertions suivantes sont équivalentes :
\begin{enumerate}
\item\label{itpripoi1} $F$ est fini ;
\item\label{itpripoi2} $F$ est à valeurs de dimensions finies ;
\item\label{itpripoi3} $F$ est de type fini et admet une décomposition en poids faible.
\end{enumerate}
\end{prop}

\begin{proof} L'implication \ref{itpripoi1}$\Rightarrow$\ref{itpripoi2} découle de \cite[prop.~6.3]{DTV} (cet énoncé suppose que $A$ est un anneau de type fini, mais la même démonstration s'applique au cas plus général où $A_K$ est une $K$-algèbre de type fini).

Si $F$ est un foncteur polynomial à valeurs de dimensions finies, alors il est fini, donc de type fini, par \cite[lemme~11.10]{DTV} ; comme $K$ est un corps de décomposition de $\mathbf{P}(A)$, $F$ possède une décomposition en poids faible par la proposition~\ref{pr-vdfdppol}, d'où l'implication \ref{itpripoi2}$\Rightarrow$\ref{itpripoi3}.

Enfin, \ref{itpripoi3} entraîne \ref{itpripoi1} grâce au théorème~\ref{th-poinoeth_abst} et au corollaire~\ref{cor-algfinie}.
\end{proof}

\chapter{Spectre d'un foncteur}\label{chap-spec}

\begin{cvi}
Dans tout ce chapitre, $\A$ désigne une catégorie additive $A$-linéaire essentiellement petite et $\E$ une catégorie de Grothendieck $K$-linéaire.
\end{cvi}

Nous introduisons dans ce chapitre trois invariants d'un foncteur de $\fct(\A,\E)$ : son \textit{radical}, qui est un idéal de $A$, son \textit{algèbre caractéristique}, qui est un quotient de l'algèbre $K[A_\mu]$, et son \textit{spectre}, qui est un foncteur de $\F(A,K)$, quotient du foncteur projectif $P^A$. Ces invariants sont de plus en plus précis : l'algèbre caractéristique détermine le radical et le spectre détermine l'algèbre caractéristique. Le spectre, qui constitue de loin la notion la plus importante de ce chapitre, et l'un des outils les plus importants de ce mémoire avec les décompositions en poids, détecte certaines propriétés déjà rencontrées, comme la polynomialité, l'existence d'une décomposition en poids finie ou l'existence de certaines décompositions à la Steinberg. Si le spectre ne détecte pas la noethérianité, il nous servira de façon cruciale pour établir, ultérieurement, la noethérianité locale de $\F^\df(A,K)$.

\section{Radical d'un foncteur}

\begin{prdef} Soit $F$ un foncteur de $\fct(\A,\E)$.
Notons\index{nota}{r@$\rf$ \emph{(radical d'un foncteur)}}
$$\rf(F):=\{\lambda\in A\,|\,\forall (x,y)\in\mathrm{Ob}\,\A^2\quad\forall (f,g)\in\A(x,y)^2\quad F(f+\lambda.g)=F(f)\}.$$
Alors $\rf(F)$ est un idéal de $A$.
De plus, $F$ se factorise de manière à travers le foncteur canonique $\A\twoheadrightarrow\A/\rf(F)$, et tout idéal $I$ de $A$ tel que $F$ se factorise à travers $\A\twoheadrightarrow\A/I$ est inclus dans $\rf(F)$.

L'idéal $\rf(F)$ est appelé \index{termin}{radical!d'un foncteur}\textbf{radical} de $F$ relativement à $A$ ; l'anneau quotient $A/\rf(F)$ sera noté \index{nota}{A@$\AF$} $\AF$.
\end{prdef}

La notation $\AF$ est analogue à celle de l'anneau $A_M$ des homothéties d'un $A$-module $M$ utilisée dans \cite{Bki}. La terminologie de {\em radical} est inspirée de celle du radical de Jacobson, qui est constitué des éléments d'un anneau qui agissent trivialement sur tout module simple. Il n'y aura pas de confusion possible avec la notion générale de radical d'un objet d'une catégorie abélienne (l'intersection de ses sous-objets stricts maximaux), que nous n'utiliserons jamais dans ce mémoire.

\begin{rema}[Changement de catégories source ou but]\label{rq-radcbev}
\begin{enumerate}
    \item Dans cette définition, le fait que la catégorie but $\E$ soit $K$-linéaire, ou même abélienne, ne joue aucun rôle. Nous n'aurons toutefois besoin de cette notion que dans ce contexte.
    \item En particulier, si $\Phi : \E\to\E'$ est un foncteur arbitraire, alors $\rf(\Phi\circ F)\supset\rf(F)$ pour tout foncteur $F$ de $\fct(\A,\E)$, avec égalité si $\Phi$ est fidèle.
    \item\label{itrqradev} Si $\B$ est une catégorie additive $A$-linéaire essentiellement petite, $\Psi : \B\to\A$ un foncteur $A$-linéaire et $F$ un foncteur de $\fct(\A,\E)$, alors $\rf(F\circ\Psi)\supset\rf(F)$, avec égalité si $\Psi$ est plein et essentiellement surjectif.
\end{enumerate}
\end{rema}

\begin{exem}\label{ex-radPAI}
Tout idéal $I$ de $A$ est le radical d'un foncteur de $\F(A,K)$. Plus précisément, si $\pi_I : \mathbf{P}(A)\to\mathbf{P}(A/I)$ désigne la réduction modulo $I$, on vérifie aussitôt que $\rf(\pi_I^*P^{A/I}_{\mathbf{P}(A/I)})=I$ et que, de plus, $\pi_I^*P^{A/I}_{\mathbf{P}(A/I)}$ est le plus grand quotient du foncteur $P^A$ de $\F(A,K)$ dont le radical égale $I$.
\end{exem}

Comme le foncteur canonique $\A\twoheadrightarrow\A/I$ est plein et essentiellement surjectif pour tout idéal $I$ de $A$, la précomposition par celui-ci, qui, comme tout foncteur de précomposition, est bicontinue, et commute au produit tensoriel pour $\E=K\Md$, est pleinement fidèle, et son image essentielle est stable par sous-quotient. Ainsi :

\begin{prop}\label{pr-prebiloc} Pour tout idéal $I$ de $A$, la sous-catégorie pleine des foncteurs $F$ de $\fct(\A,\E)$ tels que $\rf(F)\supset I$ est prébilocalisante.\index{termin}{prebilocalisante@prébilocalisante \emph{(sous-catégorie)}} Pour $\E=K\Md$, elle est également stable par produit tensoriel.
\end{prop}

\begin{rema}
Le théorème d'excision \cite[th.~5.1]{DT-ext} montre que $\F(A/I,K)$ est une sous-catégorie épaisse de $\F(A,K)$ si et seulement si $I/I^2\otimes_\mathbb{Z}K=0$, et que cette condition implique plus généralement que $\fct(\A/I,\E)$ est épaisse dans $\fct(\A,\E)$.
\end{rema}

\begin{nota}
Si $a$ est un objet de $\A$, on note $\tau_a$\index{nota}{t@$\tau_a$ \emph{(foncteur de décalage par $a$)}} le \emph{foncteur de décalage} par $a$ dans $\fct(\A,\E)$, c'est-à-dire son endofoncteur de précomposition par $a\oplus -$.
\end{nota}

Dans le deuxième point de l'énoncé suivant, si $I$ est un idéal de $A$, on voit le sous-espace vectoriel $\bar{K}[I]$ (idéal d'augmentation de $K[I]$, vu comme $K$-algèbre du groupe \emph{additif} sous-jacent à $I$) de $K[A]$ comme un idéal de l'algèbre $K[A_\mu]$, qui agit sur la catégorie abélienne $\fct(\A,\E)$.

\begin{prop}\label{pr-radevtau} Soient $I$ un idéal de $A$ et $F$ un foncteur de $\fct(\A,\E)$. Les assertions suivantes sont équivalentes :
\begin{enumerate}
\item\label{itrad1} $\rf(A)\supset I$ ;
\item\label{itrad2} pour tout objet $a$ de $\A$, $\tau_a(F)\in U(\fct(\A,\E);\bar{K}[I])$ ;
\item\label{itrad3} pour tout objet $x$ de $\A$ et tout $\lambda\in I$, le morphisme 
$F(x\oplus x)\xrightarrow{F(1\quad\lambda)-F(1\quad 0)} F(x)$
de $\E$ est nul ;
\item\label{itrad4} pour tout objet $x$ de $\A$ et tout $\lambda\in I$, le morphisme 
$F(x)\xrightarrow{F\left(\begin{array}{c} 1 \\
                                                \lambda 
                                               \end{array}\right)-F\left(\begin{array}{c} 1 \\
                                                0 
                                               \end{array}\right)} F(x\oplus x)$
de $\E$ est nul.
\end{enumerate}
\end{prop}

\begin{proof}
Comme $\bar{K}[I]$ est engendré comme $K$-espace vectoriel par les éléments $[\lambda]-[0]$ pour $\lambda\in I$, la condition \ref{itrad2} signifie que pour tous objets $a$ et $x$ de $\A$ et tout $\lambda\in I$, l'endomorphisme $F(1_a\oplus\lambda.1_x)-F(1_a\oplus 0_x)$ de $F(a\oplus x)$ est nul. La définition de $\rf(F)$ montre donc l'implication \ref{itrad1}$\Rightarrow$\ref{itrad2}.

Comme le morphisme $F(x\oplus x)\xrightarrow{F(1\quad\lambda)-F(1\quad 0)} F(x)$ égale la composée de l'endomorphisme $F(1_x\oplus\lambda.1_x)-F(1_x\oplus 0_x)$ de $F(x\oplus x)$ et du morphisme $F(x\oplus x)\to F(x)$ induit par la somme $x\oplus s\to x$, \ref{itrad2} entraîne \ref{itrad3} ; de même, \ref{itrad2} entraîne \ref{itrad4}.

Soient $a$, $b$ des objets de $\A$ et $f$, $g$ des éléments de $\A(a,b)$. Posons $x:=a\oplus b$. Le morphisme $F(f+\lambda g)-F(f) : F(a)\to F(b)$ de $\E$ égale la composée du morphisme $F(a)\xrightarrow{F\left(\begin{array}{c} 1 \\
                                                f\\
                                                0\\
                                                g
                                               \end{array}\right)}F(a\oplus b\oplus a\oplus b)=F(x\oplus x)$,
de $F(x\oplus x)\xrightarrow{F(1\quad\lambda)-F(1\quad 0)} F(x)=F(a\oplus b)$ et de $F(a\oplus b\twoheadrightarrow b)$. Il s'ensuit que \ref{itrad3} entraîne \ref{itrad1}. La démonstration de l'implication \ref{itrad4}$\Rightarrow$\ref{itrad1} est entièrement analogue, d'où la proposition.
\end{proof}

\begin{coro}
Soient $F$ et $G$ des foncteurs non nuls de $\F(\A;K)$. On a
$$\rf(F\otimes G)=\rf(F)\cap\rf(G)\,.$$
\end{coro}

\begin{proof} Soit $a$ un objet de $\A$ tel que $G(a)$ soit non nul. Alors le foncteur $\tau_a(G)$ contient un facteur direct $K$ (foncteur constant), donc $\tau_a(F\otimes G)=\tau_a(F)\otimes\tau_a(G)$ contient un facteur direct $\tau_a(F)$. La proposition~\ref{pr-radevtau} montre donc que $\rf(F\otimes G)\subset\rf(F)$. De même $\rf(F\otimes G)\subset\rf(G)$, d'où $\rf(F\otimes G)\subset\rf(F)\cap\rf(G)$.

L'inclusion inverse, immédiate, figure déjà dans la proposition~\ref{pr-prebiloc}.
\end{proof}

\begin{coro}\label{cor-radExt}
Soit $0\to G\to F\to H\to 0$ une suite exacte courte de $\fct(\A,\E)$. On a une inclusion
$$\rf(G).\rf(H)\subset\rf(F)\subset\rf(G)\cap\rf(H)$$
d'idéaux de $A$.
\end{coro}

\begin{proof}
La proposition~\ref{pr-radevtau} et le lemme~\ref{lm-idproduit} montrent que, pour $(\alpha,\beta)\in\rf(G)\times\rf(H)$, l'élément $[\alpha\beta]-[0]=([\alpha]-[0])([\beta]-[0])$ de $K[A_\mu]$ appartient à $\rf(F)$, d'où $\alpha\beta\in\rf(F)$, toujours par la proposition~\ref{pr-radevtau}. Comme $\rf(F)$ est un idéal de $A$, il s'ensuit que $\rf(G).\rf(H)\subset\rf(F)$. L'autre inclusion, immédiate, figure déjà dans la proposition~\ref{pr-prebiloc}.
\end{proof}

Le résultat suivant est classique mais nous n'en connaissons pas de référence explicite dans la littérature (il est mentionné, sans démonstration, dans \cite[lemme~5.7.2]{DT-TJM}).

\begin{coro}\label{cor-radPolp}
Si $F$ est un foncteur de $\pol_d(\A,\E)$ et que $K$ est de caractéristique $p>0$, alors $\rf(F)$ contient l'idéal $(p^d)$ de $A$.
\end{coro}

\begin{proof}
Si $X$ est un foncteur de $\mathbf{Add}_d(\A,\E)$, la multiplication par $p$ sur n'importe laquelle des $d$ variables à la source induit, par multiaddtivité, la multiplication par $p$, c'est à dire $0$, au but, puisque $\E$ est $\mathbb{Z}/p$-linéaire. On en déduit que le radical du foncteur $\delta_d^*(X)$ de $\fct(\A,\E)$ contient l'idéal $(p)$.

On montre maintenant l'assertion par récurrence sur $d$ : on peut supposer $d>0$ et l'assertion établie pour les foncteurs de $\pol_{d-1}(\A,\E)$, puisqu'elle est triviale en degré $0$. Il existe une suite exacte courte $0\to G\to F\to\qpol_{d-1}(F)\to 0$ où $G$ est un quotient de $\delta_d^*(cr_d(F))$. L'hypothèse de récurrence montre l'inclusion $(p^{d-1})\subset\rf(\qpol_{d-1}(F))$ ; comme $G$ est un quotient de $\delta_d^*(cr_d(F))$ et que $cr_d(F)$ appartient à $\mathbf{Add}_d(\A,\E)$ (puisque $F$ est polynomial de degré au plus $d$), ce qui précède (et la proposition~\ref{pr-prebiloc}) montre l'inclusion $(p)\subset\rf(G)$. Le corollaire~\ref{cor-radExt} permet de conclure.
\end{proof}

\begin{rema}\begin{enumerate}
    \item Le corollaire~\ref{cor-radPolp} est le pendant en caractéristique première de la dernière assertion de la proposition~\ref{pr-decpolcar0}, qui implique pour sa part que, en caractéristique $p=0$, si $F$ est un foncteur analytique de $\fct(\A,\E)$, alors $\rf(F)$ contient l'idéal des éléments de $\mathbb{Z}$-torsion de $A$.
    \item Le corollaire~\ref{cor-radPolp}, qu'on mentionne ici comme simple illustration de la notion de radical d'un foncteur (on ne s'en servira pas par la suite), n'est pas optimal : en examinant les quotients de la $K$-algèbre du groupe $\mathbb{Z}$ par les puissances de son idéal d'augmentation, on peut montrer qu'on a en fait $(p^r)\subset\rf(F)$ pour $F$ dans $\pol_d(\A,\E)$ dès lors que $p^r>d$.
\end{enumerate}
\end{rema}

Nous terminons cette section par d'autres considérations élémentaires sur le radical des foncteurs, liées aux foncteurs antipolynomiaux. 

\begin{prop}\begin{enumerate}
    \item Un foncteur $F$ de $\F(A,K)$ est antipolynomial si et seulement si l'anneau $A_{/F}$ est $K$-trivial.
    \item Si le $A$-module $\A(a,b)$ est de type fini pour tous objets $a$ et $b$ de $\A$, alors tout foncteur $F$ de $\fct(\A,\E)$ tel que l'anneau $A_{/F}$ soit $K$-trivial est antipolynomial.
\end{enumerate}
\end{prop}

\begin{proof}
La première assertion résulte de la définition. La seconde découle de ce que, si $I$ est un idéal $K$-cotrivial de $A$ et $M$ un $A$-module de type fini, alors $M/I.M$ est un $A/I$-module de type fini, donc un groupe abélien fini dont le produit tensoriel sur $\mathbb{Z}$ avec $K$ est trivial.
\end{proof}

Les deux énoncés suivants sont immédiats.

\begin{prop}\label{pr-radLin} Soient $M$ un objet non nul de $\E$ et $X : \A\to\mathbf{Ab}$ un foncteur additif. L'anneau (resp. l'idéal) $A_{/M[X]}$ (resp. $\rf(M[X])$) est l'image (resp. le noyau) du morphisme d'anneaux $A\to Z(\mathrm{End}(X))$ donné par la structure $A$-linéaire de $\A$.
\end{prop}

\begin{coro}\label{cor-radLin}
Soient $a$ un objet de $\A$ et $M$ un objet non nul de $\E$. L'anneau (resp. l'idéal) $A_{/M[\A(a,-)]}$ (resp. $\rf(M[\A(a,-)])$) est l'image (resp. le noyau) du morphisme d'anneaux $A\to Z(\mathrm{End}_\A(a))$ donné par la structure $A$-linéaire de $\A$.
\end{coro}

\begin{coro}
Si $F$ est un foncteur antipolynomial à support fini de $\fct(\A,\E)$, alors l'anneau $A_{/F}$ est $K$-trivial.
\end{coro}

\begin{proof}
Quitte à remplacer $\A$ par $\A/\I$ pour un idéal $\I$ approprié de la catégorie additive $\A$ (cf. remarque~\ref{rq-radcbev}.\ref{itrqradev}), on peut supposer que celle-ci est $K$-triviale. Comme $F$ est à support fini, il existe un ensemble fini $E$ d'objets de $\A$ et une famille $(M_a)_{a\in E}$ d'objets de $\A$ telle que $F$ soit un quotient de $\bigoplus_{a\in E} M_a[\A(a,-)]$. La proposition~\ref{pr-prebiloc} montre alors l'inclusion
$$\rf(F)\supset\bigcap_{a\in E}\rf(M_a[\A(a,-)])$$
tandis que le corollaire~\ref{cor-radLin} montre que $\rf(M_a[\A(a,-)])$ est un idéal $K$-cotrivial de $A$ (le quotient de $A$ par cet idéal s'injecte dans $Z(\mathrm{End}_\A(a))$, qui est un anneau $K$-trivial puisque la catégorie $\A$ est $K$-triviale). Le corollaire résulte donc de ce qu'une intersection finie d'idéaux $K$-cotriviaux de $A$ est $K$-cotriviale.
\end{proof}

\begin{rema} Si $F$ est un foncteur antipolynomial de $\F(\A;K)$, l'anneau $A_{/F}$ n'est généralement pas $K$-trivial, même si tous les $A$-modules $\A(a,b)$ sont de type fini. Il suffit pour le voir de considérer, pour $A=\mathbb{Z}$ et $\A$ la sous-catégorie pleine de $\mathbf{Ab}$ constituée des groupes abéliens finis d'ordre inversible dans $K$, la linéarisation $F$ du foncteur d'oubli $X : \A\to\mathbf{Ab}$ : le morphisme d'anneaux $\mathbb{Z}\to Z(\mathrm{End}(X))$ est manifestement injectif, ainsi
$A_{/F}=\mathbb{Z}$ par la proposition~\ref{pr-radLin}.
\end{rema}

\section{Algèbre et fonction caractéristiques d'un foncteur}

\begin{defi} Soit $F$ un foncteur de $\fct(\A,\E)$. On note \index{nota}{chi@$\chi$, $\chi_F$, $\bar{\chi}_F$ \emph{(fonction caractéristique du foncteur $F$)}}$\chi_F : A_\mu\to Z(\mathrm{End}(F\circ\oplus))_\mu$ le morphisme de monoïdes donné par $a\mapsto F(1\oplus a)$. On l'appelle \index{termin}{fonction caracteristique@fonction caractéristique}\textbf{fonction caractéristique} de $F$ (relativement à $A$).

L'image du morphisme de $K$-algèbres $K[A_\mu]\to Z(\mathrm{End}(F\circ\oplus))$ qui prolonge $\chi_F$ est appelée \index{termin}{algebre caracteristique@algèbre caractéristique}\textbf{algèbre caractéristique} de $F$ (relativement à $A$) et notée \index{nota}{A@$\ac$, $\ac_A$ \emph{(algèbre caractéristique d'un foncteur)}} $\ac_A(F)$, ou simplement $\ac(F)$. Son noyau est noté \index{nota}{I@$\ic$, $\ic_A$ \emph{(idéal caractéristique d'un foncteur)}} $\ic_A(F)$ (ou $\ic(F)$) et appelé \index{termin}{ideal caracteristique@idéal caractéristique}\textbf{idéal caractéristique} de $F$.

Le morphisme $\chi_F$ étant à valeurs dans $\ac_A(F)$, on notera $\bar{\chi}_F : A_\mu\to\ac_A(F)_\mu$ le morphisme de monoïdes induit par $\chi_F$.
\end{defi}

Si elles sont loin de caractériser le type d'isomorphisme d'un foncteur $F$ de $\F(A,K)$ (ou $\fct(\A,\E)$), la fonction et l'algèbre caractéristiques de $F$ contiennent beaucoup d'informations sur ce foncteur, au moins si $A$ est assez gros, comme nous le verrons par la suite (cf. par exemple théorème~\ref{th-chipol}). La terminologie employée est inspirée de celle de \emph{caractère} en théorie classique des représentations.

\begin{rema}[Changement de catégories source ou but]\label{rq-radcbevic} Soit $\rho : B\to A$ un morphisme d'anneaux. Si $\B$ est une catégorie additive essentiellement petite $B$-linéaire, $\E'$ une catégorie de Grothendieck $K$-linéaire et $\Phi : \E\to\E'$ (resp. $\Psi : \B\to\A$) un foncteur $K$-linéaire (resp. $B$-linéaire, où $B$ agit sur $\A$ via $\rho$), alors pour tout foncteur $F$ de $\fct(\A,\E)$, on a $\ic_B(\Phi\circ F\circ\Psi)\supset K[\rho]^{-1}(\ic_A(F))$, avec égalité si $\Phi$ est pleinement fidèle et $\Psi$ plein et essentiellement surjectif. (Cf. remarque~\ref{rq-radcbev}.)

Au lieu de supposer que $\Psi$ est plein et essentiellement surjectif, on peut plus généralement supposer que la précomposition par $\Psi$ est pleinement fidèle, ce qui est aussi le cas lorsque $A=B[S^{-1}]$ pour une partie multiplicative $S$ de $A$, que $\A=\B[S^{-1}]$ et que $\rho : B\to B[S^{-1}]$ (resp. $\Psi : \B\to\B[S^{-1}]$) est le morphisme (resp. foncteur) canonique.
\end{rema}

Le résultat suivant, analogue à la proposition~\ref{pr-prebiloc}, est immédiat.

\begin{prop}\label{pr-prebilocIdCar} Pour tout sous-ensemble $E$ de $K[A]$, la sous-catégorie pleine des foncteurs $F$ de $\fct(\A,\E)$ tels que $\ic(F)\supset E$ est prébilocalisante.\index{termin}{prebilocalisante@prébilocalisante \emph{(sous-catégorie)}}
\end{prop}

\begin{rema} En revanche, la sous-catégorie pleine des foncteurs $F$ de $\F(\A;K)$ des foncteurs dont l'idéal caractéristique contient un sous-ensemble donné de $A$ n'est généralement pas stable par produit tensoriel, voir l'exemple~\ref{ex-icpt} ci-dessous.
\end{rema}

La démonstration de la proposition suivante est facile et très semblable à celle de la proposition~\ref{pr-radevtau}, nous la laissons donc en exercice.

\begin{prop}\label{pr-rfevtau} Soient $\xi=\sum_{i=1}^n\lambda_i[a_i]$ un élément de $K[A]$ (où $a_i\in A$ et $\lambda_i\in K$) et $F$ un foncteur de $\fct(\A,\E)$. Les assertions suivantes sont équivalentes :
\begin{enumerate}
\item\label{itrfev1} $\xi\in\ic(F)$ ;
\item\label{itrfeveq} pour toutes flèches parallèles $f$ et $g$ de $\A$, $\sum_{i=1}^n\lambda_i F(f+a_i.g)=0$ ;
\item\label{itrfev2} pour tout objet $x$ de $\A$, l'image  de $\xi$ par la projection $K[A_\mu]\twoheadrightarrow\mathfrak{Z}(\tau_x(F))$ est nulle ;
\item\label{itrfev3} pour tout objet $x$ de $\A$, le morphisme 
$F(x\oplus x)\xrightarrow{\sum_{i=1}^n\lambda_i F(1\quad a_i)} F(x)$
de $\E$ est nul ;
\item\label{itrfev4} pour tout objet $x$ de $\A$ et tout $\lambda\in I$, le morphisme 
$F(x)\xrightarrow{\sum_{i=1}^n\lambda_i F\left(\begin{array}{c} 1 \\
                                                a_i
                                               \end{array}\right)} F(x\oplus x)$
de $\E$ est nul.
\end{enumerate}
\end{prop}

\begin{exem}\label{ex-icpt} En caractéristique $p=0$, on vérifie à partir de l'équivalence \ref{itrfev1}$\Leftrightarrow$\ref{itrfev2} ci-dessus que $\ac_K(T^d)\simeq K^{d+1}$ ($K$-algèbre produit), où $T^d$ désigne le foncteur $d$-ième puissance tensorielle dans $\F(K,K)$. La fonction caractéristique de $T^d$ s'identifie via cet isomorphisme à $K\to K^{d+1}\quad\lambda\mapsto (\lambda^i)_{0\le i\le d}$.

En particulier, si $d>0$, $\ac_K(T^d\otimes T^d)$ n'est pas un quotient de $\ac_K(T^d)$, de sorte que $\ic_K(T^d\otimes T^d)$ ne peut pas contenir $\ic_K(T^d)$.
\end{exem}

Avant de donner un autre exemple simple de calcul d'idéal caractéristique (proposition~\ref{pr-algcarP}), nous donnons trois résultats préliminaires.

\begin{lemm}\label{lm-tauP} Soit $w\in\mon(A_\mu,K_\mu)$ un poids tel que
$$\forall a\in A\qquad (w(a)\ne 0)\Rightarrow (a\in A^\times)\,.$$
Alors le foncteur $P^A$ est facteur direct de $\tau_A(P^A_w)$ dans $\F(A,K)$. 
\end{lemm}

\begin{proof}
Considérons le morphisme $\varphi : P^A\to\tau_A(P^A_w)$ composé de la flèche $P^A\to\tau_A(P^A)$ donnée sur $V$ par $P^A(V)\to\tau_A(P^A)(V)=P^A(A\oplus V)\quad [v]\mapsto [(1,v)]$ pour tout $v\in V$ et de la projection $\tau_A(P^A)\twoheadrightarrow\tau_A(P^A_w)$. Par ailleurs, le morphisme $\tau_A(P^A)\to P^A$ donné sur $V$ par $[(a,v)]\mapsto w(a)[a^{-1}.v]$ si $a\in A^\times$ et $[(a,v)]\mapsto 0$ si $a\notin A^\times$ induit un morphisme $\psi : \tau_A(P^A_w)\to P^A$ grâce à l'hypothèse faite sur $w$. La composée $\psi\circ\varphi$ égale l'identité de $P^A$, d'où le lemme.
\end{proof}

\begin{lemm}\label{lm-ictau} Soient $F$ un foncteur de $\fct(\A,\E)$ et $a$ un objet de $\A$. Les idéaux $\ic(F)$ et $\ic(\tau_a(F))$ de $K[A_\mu]$ sont égaux.
\end{lemm}

\begin{proof}
Comme $F$ est facteur direct de $\tau_a(F)$, on a $\ic(\tau_a(F))\subset\ic(F)$ (proposition~\ref{pr-prebilocIdCar}). L'inclusion réciproque résulte de l'équivalence \ref{itrfev1}$\Leftrightarrow$\ref{itrfev2} de la proposition~\ref{pr-rfevtau}.
\end{proof}

\begin{prop}\label{pr-algcarP} Pour tout $w\in\mon(A_\mu,K_\mu)$, on dispose d'un morphisme d'algèbres $\ac(P^A_w)\simeq K[(A/I)_\mu]$ où $I$ est l'idéal des éléments de $A$ annulant un élément de $A\setminus\la(w)$.
\end{prop}

\begin{proof} L'idéal $I$ est le noyau du morphisme d'anneaux canonique $A\to A[(A\setminus\la(w))^{-1}]$. Par conséquent, la proposition~\ref{pr-invad} montre que, quitte à remplacer $A$ par $A[(A\setminus\la(w))^{-1}]$ (cf. remarque~\ref{rq-radcbevic}), on peut supposer $\la(w)=A\setminus A^\times$. La conclusion résulte alors des lemmes~\ref{lm-tauP} et~\ref{lm-ictau} et de la nullité, immédiate, de $\ic(P^A)$.
\end{proof}

La propriété suivante nous servira notamment au §\,\ref{pp0}. On rappelle que la notion de bi-idéal est introduite page~\pageref{pbiid}.

\begin{prop}\label{pr-idc} L'idéal caractéristique de tout foncteur de $\fct(\A,\E)$ est un bi-idéal\index{termin}{biideal@bi-idéal} de $K[A]$.
\end{prop}

\begin{proof} Soit $F$ un foncteur de $\fct(\A,\E)$.
Par définition, $\ic(F)$ est un idéal de $K[A_\mu]$. Si $\sum_{i=1}^n\lambda_i[a_i]$ est un élément de $\ic(F)$ (où $a_i\in A$ et $\lambda_i\in K$) et $t$ un élément de $A$, en substituant $f+t.g$ à $f$ dans l'égalité du point~\ref{itrfeveq} de la proposition~\ref{pr-rfevtau}, on voit que $\sum_{i=1}^n\lambda_i[a_i+t]$ appartient également à $\ic(F)$, qui est donc aussi un idéal de $K[A_\mathrm{add}]$.
\end{proof}

Réciproquement, tout bi-idéal de $K[A]$ est l'idéal caractéristique d'un foncteur de $\F(A,K)$. Plus précisément, on a l'analogue suivant de l'exemple~\ref{ex-radPAI} :

\begin{prop}\label{pr-biid} Soit $\I$ un bi-idéal de $K[A]$. Définissons un quotient $P^A_\I$ de $P^A$ par la formule $P^A_\I(V)=K[V]/(K[V]\star (\I.K[V]))$. (Ici, on voit $K[V]$ comme un $K[A_\mu]$-module d'une part, d'où la notation $\I.K[V]$, et comme $K$-algèbre du groupe additif sous-jacent à $V$ d'autre part, avec la notation $\star$\index{nota}{$\star$} pour la multiplication.)

Alors $\I$ est l'idéal caractéristique de $P^A_\I$, et $P^A_\I$ est le plus grand quotient de $P^A$ ayant cette propriété.
\end{prop}

\begin{proof} De manière générale, il résulte de l'équivalence \ref{itrfev1}$\Leftrightarrow$\ref{itrfev3} de la proposition~\ref{pr-rfevtau} que tout foncteur $F$ de $\fct(\A,\E)$ possède un plus grand quotient $F_\I$ tel que $\I\subset\ic(F_\I)$, donné par 
$$F_\I(x)=\mathrm{Coker}\,\left(\bigoplus_{\xi=\sum_{i=1}^n\lambda_i[a_i]\in\I}F(x\oplus x)\xrightarrow{\sum_{i=1}^n\lambda_i F(1\quad a_i)}F(x)\right)\,,$$
en particulier, $P^A_\I$ est le plus grand quotient de $P^A$ tel que $\I\subset\ic(P^A_\I)$. Il ne reste donc plus qu'à montrer l'inclusion $\ic(P^A_\I)\subset\I$.

Comme $\I$ est un bi-idéal de $K[A]$, on a $K[A]\star(\I.K[A])=\I$ et donc $P^A_\I(A)=K[A]/\I$. Du fait que $P^A_\I$ est un quotient de $P^A$, le morphisme $\mathrm{End}(P^A_\I)\to P^A_\I(A)=K[A]/\I$ induit par l'évaluation en $A$ est injectif. La composée du morphisme canonique $K[A]\to\mathrm{End}(P^A_\I)$ et de l'injection précédente égale la projection $K[A]\twoheadrightarrow K[A]/\I$. Comme son noyau $\I$ contient $\ic(P^A_\I)$, la proposition est démontrée.
\end{proof}

\begin{rema}La structure des foncteurs $P^A_\I$ est difficile à comprendre en général. Par exemple, si $\I=\bar{K}[A]^{\star n}$ pour un $n\in\mathbb{N}$, $P^A_\I$ est polynomial si $A$ n'a pas de quotient fini de caractéristique $p$ par le théorème~\ref{th-chipol} ci-après, mais nous n'en connaissons pas de démonstration directe, et $P^A_{\bar{K}[A]^{\star n}}$ n'est pas polynomial si $n$ est assez grand et que $A$ possède un quotient fini de caractéristique $p$.
\end{rema}

\subsection*{Propriétés de finitude de l'algèbre caractéristique}

Les foncteurs dont l'algèbre caractéristique possède certaines propriétés de finitude ou de régularité joueront un rôle considérable dans ce mémoire (cf. proposition-définition~\ref{pr-htdcf} ci-après). Une observation fondamentale est la généralisation presque immédiate suivante de la proposition~\ref{pr-tfdfend}, qui s'applique notamment aux foncteurs de type fini et à valeurs de dimension finie.
 
\begin{prop}\label{pr-dfchi} Soit $F$ un foncteur de $\F(\A;K)$. Supposons qu'il existe un ensemble fini $E$ d'objets de $\A$ qui constitue un support ou un co-support de $F$ et tel que $\dim_K F(x\oplus x)<\infty$ pour tout $x\in E$. Alors $\mathrm{End}(F\circ\oplus)$, donc $\ac(F)$, est de dimension finie sur $K$.
\end{prop}

\begin{proof} Si $E$ est un support de $F$ (le cas d'un co-support s'établit de façon analogue), alors la flèche canonique $\F(\A;K)(F,G)\to\bigoplus_{x\in E}\mathrm{Hom}_K(F(x),G(x))$ est injective pour tout foncteur $G$. La conclusion découle donc de l'isomorphisme d'adjonction somme/diagonale $\mathrm{End}(F\circ\oplus)\simeq\F(\A;K)(F,F\circ\oplus\circ\delta)$.
\end{proof}

\begin{coro}\label{cor-dfchi} Si $F$ est un foncteur de type fini de $\F^\df(\A;K)$, alors $\dim_K\ac(F)$ est finie.
\end{coro}

Nous verrons un peu plus tard (cf. corollaire~\ref{cor-spchi2}) que, réciproquement, tout foncteur de type fini de $\F(A,K)$ (ou plus généralement de $\F(\A;K)$, sous une hypothèse de finitude sur les morphismes de $\A$) dont l'algèbre caractéristique est de dimension finie est à valeurs de dimension finie.

La proposition suivante nous servira au chapitre~\ref{stfdf}, par l'intermédiaire de son cas particulier fondamental qu'est le corollaire~\ref{cor-acdf}.

\begin{prop}\label{pr-acdf} Supposons que $\A=\B\otimes_\mathbb{Z}A$ pour une petite catégorie additive $\B$, et que $F$ est un foncteur de type fini de $\fct(\A,\E)$ tel que $\ac_A(F)$ soit une $K$-algèbre de type fini. Alors il existe un sous-anneau de type fini $A'$ de $A$ tel que l'image de $F$ dans $\fct(\B\otimes_\mathbb{Z} A',\E)$ par la précomposition par le foncteur canonique $\B\otimes_\mathbb{Z}A'\to\B\otimes_\mathbb{Z}A=\A$ soit de type fini.
\end{prop}

\begin{proof} Soit $E$ un sous-ensemble fini de $A$ tel que la $K$-algèbre $\ac_A(F)$ soit engendrée par les $\chi_F(a)$ pour $a\in E$. Notons $A'$ le sous-anneau de $A$ engendré par $E$. Pour tout $x\in A$ et toute paire $(f,g)$ de flèches parallèles de $\A$, $F(f+x.g)$ est ainsi déterminé par les $F(f+a.g)$ avec $a\in A'$ : il existe $a_1,\dots,a_n\in A'$ et $\lambda_1,\dots,\lambda_n\in K$ tels que $F(f+x.g)=\sum_{i=1}^n\lambda_i F(f+a_i.g)$. Par conséquent, le foncteur de restriction à la source $\fct(\A,\E)\to\fct(\B\otimes_\mathbb{Z}A',\E)$ induit un isomorphisme entre l'ensemble ordonné des sous-objets de $F$ et celui de son image, d'où la proposition.
\end{proof}

\begin{coro}\label{cor-acdf} Soit $F$ un foncteur de type fini de $\F(A,K)$ tel que $\ac(F)$ soit une $K$-algèbre de type fini. Il existe un sous-anneau de type fini $A'$ de $A$ tel que l'image de $F$ dans $\F(A',K)$ (par la précomposition par le changement de base $\mathbf{P}(A')\to\mathbf{P}(A)$) soit de type fini.
\end{coro}

\subsection*{Idéal caractéristique et radical}

La proposition élémentaire suivante nous servira à de multiples reprises.

\begin{prop}\label{pr-acrad} Soient $F$ un foncteur de $\fct(\A,\E)$ et $a$, $b$ des éléments de $A$. On a $\chi_F(a)=\chi_F(b)$ si et seulement si $a-b\in\rf(F)$. 

Par conséquent, $\chi_F$ induit une fonction multiplicative \emph{injective} $A_{/F}\hookrightarrow\ac_A(F)$.
\end{prop}

\begin{proof} En utilisant la définition de $\chi_F$ et $\ic(F)$ puis la proposition~\ref{pr-idc}, on obtient les équivalences
$$(\chi_F(a)=\chi_F(b))\Leftrightarrow ([a]-[b]\in\ic(F))\Leftrightarrow ([a-b]-[0]\in\ic(F))\,.$$

La proposition~\ref{pr-rfevtau} (équivalence \ref{itrfev1}$\Leftrightarrow$\ref{itrfeveq}) montre que cette dernière condition équivaut à $a-b\in\rf(F)$.
\end{proof}

\begin{coro}\label{cor-tailleA} Soit $F$ un foncteur de $\fct(\A,\E)$ tel que $\ac_A(F)$ est de dimension finie sur $K$. Alors $\AF$ est fini si $K$ est fini, et $\operatorname{Card}\AF\le\operatorname{Card}K$ sinon.
\end{coro}

\begin{coro}\label{cor-Kfiniln} Si le corps $K$ est fini, la catégorie $\F^\df(A,K)$ est \index{termin}{localement noetherienne@localement noethérienne \emph{(catégorie abélienne)}} localement noethérienne.
\end{coro}

\begin{proof} Il suffit de montrer que tout foncteur $F$ de $\F(A,K)$ de type fini et à valeurs de dimensions finies est noethérien. La proposition~\ref{pr-dfchi} et le corollaire~\ref{cor-tailleA} montrent que $\AF$ est un anneau fini ; autrement dit, on peut supposer que $A$ est {\em fini}. La conclusion est alors donnée par le théorème~\ref{th-PSS} de Putman-Sam-Snowden.
\end{proof}

\begin{rema} Nous verrons ultérieurement (théorème~\ref{th-noeth_cor_gal}) qu'on peut s'affranchir de l'hypothèse de finitude de $K$, mais la démonstration est beaucoup plus difficile.
\end{rema}

\section{Les foncteurs adjoints $\hi$ et $\ti$}

Les notions rappelées dans cette section sont dues au second auteur --- voir par exemple \cite[§\,4.1 et 4.2]{T-AIF}.

Pour $F$ dans $\fct(\A,\E)$ et $V$ dans $\mathbf{P}(A)$, on note $F_V:=F\circ (V\underset{A}{\otimes}-)$ (voir page~\pageref{not-ptA} pour la définition du produit tensoriel $\underset{A}{\otimes} : \mathbf{P}(A)\times\A\to\A$) et $F^V:=F_{V^*}$, où $V=\mathrm{Hom}_A(V,A)$.

\subsection{Définitions}

On définit un foncteur \index{nota}{Hom@$\hi$|(} $\hi : \fct(\A,\E)^\op\times\fct(\A,\E)\to\F(A,K)$ par
$$\hi(F,G)(V)=\fct(\A,\E)(F,G_V)(\simeq\fct(\A,\E)(F^V,G))$$
et un foncteur \index{nota}{$\ti$|(} $\ti : \F(A,K)\times\fct(\A,\E)\to\fct(\A,\E)$ par
$$(T\ti X)(a):=\mathrm{Coend}(\mathbf{P}(A)^\op\times\mathbf{P}(A)\xrightarrow{(U,V)\mapsto (U^*,V\underset{A}{\otimes}a)}\mathbf{P}(A)\times\A\xrightarrow{T\boxtimes X}\E)$$
(pour les rappels utiles sur $\mathrm{Coend}$ --- et $\mathrm{End}$ --- voir par exemple \cite[chap.~IX, §\,5-6]{ML-cat}), où $\boxtimes$ désigne le produit tensoriel extérieur sur $K$ : $(T\boxtimes X)(a,b):=T(a)\otimes X(b)$ (où $\otimes : K\Md\times\E\to\E$ est le bifoncteur cocontinu par rapport à chaque variable tel que $K\otimes -=\mathrm{Id}_\E$).

On dispose d'un isomorphisme d'adjonction
$$\fct(\A,\E)(T\ti X,Y)\simeq\F(A,K)(T,\hi(X,Y))$$
naturel en les foncteurs $T$ de $\F(A,K)$, $X$ et $Y$ de $\fct(\A,\E)$.

Le bifoncteur $\hi$ (resp. $\ti$) commute aux limites (resp. colimites) par rapport à chaque variable. De plus, on dispose d'isomorphimes naturels $P^V\ti F\simeq F_V$, en particulier $P^A\ti -\simeq\mathrm{Id}_{\fct(\A,\E)}$, et $\hi(P^a_\A,F)\simeq F\circ (-\underset{A}{\otimes} a)$ (en particulier, dans $\F(A,K)$, $\hi(P^A,-)\simeq\mathrm{Id}_{\F(A,K)}$).

Sur $\F(A,K)$, $\ti$ définit une structure monoïdale symétrique d'unité $P^A$.\index{nota}{$\ti$|)} De façon plus générale, la dissymétrie de la définition de $\ti : \F(A,K)\times\fct(\A,\E)\to\fct(\A,\E)$ n'est qu'apparente, comme le montre le résultat suivant.

\begin{prop}\label{pr-EKti}
Soient $T$ un foncteur de $\F(A,K)$ et $F$ un foncteur de $\fct(\A,\E)$. Notons $\tilde{T}$ l'extension de Kan à gauche de $T$ le long de l'inclusion de $\mathbf{P}(A)\to A\Md$.

Pour tout objet $a$ de $\A$, on dispose dans $\E$ d'un isomorphisme naturel
$$(T\ti F)(a)\simeq\mathrm{Coend}\big(\A^\op\times\A\xrightarrow{\A(-,a)\times\A}A\Md\times\A\xrightarrow{\tilde{T}\boxtimes F}\E\big)\,.$$
\end{prop}

\begin{proof}Cela découle de ce que le bifoncteur $\ti$ est caractérisé, à isomorphisme canonique près, par le fait qu'il commute aux colimites par rapport à chaque variable et que l'on dispose dans $\fct(\A,\E)$ d'un isomorphisme $P^V\ti M[\A(a,-)]\simeq M[V\underset{A}{\otimes}\A(a,-)]$ naturel en les objets $V$ de $\mathbf{P}(A)$, $a$ de $\A$ et $M$ de $\E$.
\end{proof}

On dispose également d'un foncteur, encore noté $\hi$ (il coïncide avec le précédent pour $\A=\mathbf{P}(A)$ et $\E=K\Md$), $\F(A,K)^\op\times\fct(\A,\E)\to\fct(\A,\E)$ donné par
$$\hi(F,X)(a)=\mathrm{End}(\mathbf{P}(A)^\op\times\mathbf{P}(A)\xrightarrow{F\times X\circ (-\underset{A}{\otimes}a)}(K\Md)^\op\times\E\xrightarrow{\operatorname{Hom}_K}\E).$$
(Pour $\E=K\Md$, on a donc  $\hi(F,X)(a)=\F(A,K)(F,X\circ (-\underset{A}{\otimes}a))$.)
Ce foncteur $\hi$ possède des propriétés formelles analogues au précédent.\index{nota}{Hom@$\hi$|)}

\subsection{Foncteurs $\ti$ et $\hi$ et finitude des valeurs}

\begin{prop}\label{pr-vdffac}
\begin{enumerate}
\item Soient $F$ et $G$ des foncteurs de $\F(\A;K)$. Supposons que $F$ est de type fini et $G$ à valeurs de dimension finie. Alors $\hi(F,G)$ appartient à $\F^\df(A,K)$.
\item Soient $T$ un foncteur de type fini de $\F(A,K)$ et $F$ un foncteur de $\F^\df(A,K)$. Alors $T\ti F$ appartient à $\F^\df(A,K)$.
\end{enumerate}
\end{prop}

\begin{proof} Si $F$ est un foncteur de type fini de $\F(\A;K)$, c'est un quotient d'une somme directe finie de foncteurs de la forme $P^a$ avec $a\in\mathrm{Ob}\,\A$. Il s'ensuit que $\hi(F,G)$ est un sous-foncteur d'une somme directe finie de foncteurs de la forme $\hi(P^a,G)\simeq G\circ (-\underset{A}{\otimes} a)$. En particulier, si $G$ est à valeurs de dimensions finies, il en est de même pour $\hi(F,G)$.

La deuxième assertion s'établit de façon analogue.
\end{proof}

Comme l'extension de Kan à gauche le long de l'inclusion $\mathbf{P}(A)\to A\Md$ ne préserve pas les foncteurs à valeurs de dimensions finies, on voit facilement, à partir de la proposition~\ref{pr-EKti}, que, si $T$ est un foncteur de $\F^\df(A,K)$ et $F$ un foncteur de type fini de $\F(\A;K)$, $T\ti F$ n'appartient généralement pas à $\F^\df(\A;K)$. Pour y remédier, nous introduisons la condition de finitude suivante sur les morphismes de $\A$.

(MTF)\quad Pour tous objets $a$ et $b$ de $\A$, le $A$-module $\A(a,b)$ est de type fini.\index{nota}{MTF@(MTF) \emph{(condition de finitude sur les morphismes)}|textbf}

(La notation abrège {\em Morphismes} et {\em Type Fini}.)

\begin{lemm}\label{lm-vdfEK} Soit $T$ un foncteur de $\F(A,K)$, et $\tilde{T}$ son extension de Kan à gauche à $A\Md$.
\begin{enumerate}
    \item Le foncteur $\tilde{T}$ préserve les épimorphismes.
    \item Si $T$ appartient à $\F^\df(A,K)$, alors $\tilde{T}$ prend des valeurs de dimensions finies sur les $A$-modules de type fini.
\end{enumerate}
\end{lemm}

\begin{proof}
Le foncteur $T$ est un quotient d'une somme directe de foncteurs de la forme $P^{A^n}$ pour $n\in\mathbb{N}$, dont l'extension de Kan à gauche à $A\Md$ est $V\mapsto K[V^n]$. Un tel foncteur préserve les épimorphismes. Il s'ensuit que $\tilde{T}$ est quotient d'un foncteur qui préserve les épimorphismes, et donc préserve lui-même les épimorphismes.

La deuxième assertion résulte de la première et de ce que la restriction de $\tilde{T}$ à $\mathbf{P}(A)$ est $T$, puisque l'inclusion $\mathbf{P}(A)\to A\Md$ est pleinement fidèle.
\end{proof}

\begin{prop}\label{pr-valtdf} Supposons que la condition \textnormal{(MTF)} est satisfaite. Soit $T$ un foncteur de $\F^\df(A,K)$ et $F$ un foncteur de type fini de $\F(\A;K)$. Alors $T\ti F$ appartient à $\F^\df(\A;K)$. 
\end{prop}

\begin{proof}
Ce résultat s'établit à partir de la proposition~\ref{pr-EKti} et de la deuxième assertion du lemme~\ref{lm-vdfEK} comme la proposition~\ref{pr-vdffac}.
\end{proof}

\subsection{Foncteur $\ti$ et bifoncteurs}

Dans l'énoncé suivant (comme à plusieurs autres endroits de ce mémoire), qui constitue une variation autour de l'adjonction somme/diagonale, on identifie $\F(A^2,K)$ à $\F(\mathbf{P}(A)^2;K)$.

\begin{prop}\label{pr-asdti} Soient $X$ un foncteur de $\F(A^2,K)$ et $F$ un foncteur de $\fct(\A,\E)$. On dispose d'un isomorphisme naturel
$\delta^*(X)\ti F\simeq\delta^*(X\ti(\oplus^*F))$ dans $\fct(\A,\E)$.
\end{prop}

(Le foncteur $X\ti(\oplus^*F)$ de $\fct(\A^2,\E)$ est défini en utilisant la structure $A^2$-linéaire de $\A^2$.)

\begin{proof}
Soient $U$ et $V$ des objets de $\mathbf{P}(A)$. On dispose d'isomorphismes naturels $\delta^*(P^{(U,V)})\simeq P^{U\oplus V}$ dans $\F(A,K)$ donc $\delta^*(P^{(U,V)})\ti F\simeq F\circ((U\oplus V)\underset{A}{\otimes} -)$ dans $\fct(\A,\E)$. On a par ailleurs un isomorphisme naturel
$$P^{(U,V)}\ti (\oplus^*F)\simeq (\oplus^*F)\circ\big((U,V)\underset{A\times A}{\otimes}-\big)$$
dans $\fct(\A\times\A,\E)$, i.e. $\big(P^{(U,V)}\ti (\oplus^*F)\big)(x,y)\simeq F\Big((U\underset{A}{\otimes}x)\oplus (V\underset{A}{\otimes}y)\Big)$, d'où en particulier $\big(P^{(U,V)}\ti (\oplus^*F)\big)(x,x)\simeq F\Big((U\oplus V)\underset{A}{\otimes}x\Big)$. Ainsi, on a dans $\fct(\A,\E)$ un isomorphisme
$$\delta^*(P^{(U,V)}\ti(\oplus^*F))\simeq F\circ((U\oplus V)\underset{A}{\otimes} -)\simeq\delta^*(P^{(U,V)})\ti F$$
naturel en $F$, $U$ et $V$. Comme les $P^{(U,V)}$ engendrent la catégorie $\F(A^2,K)$ et que les foncteurs $X\mapsto\delta^*(X)\ti F$ et $X\mapsto\delta^*(X\ti(\oplus^*F))$ commutent aux colimites, la conclusion s'ensuit.
\end{proof}

Nous aurons également besoin de la propriété élémentaire suivante de compatibilité entre le produit tensoriel interne $\ti$ et le produit tensoriel extérieur $\boxtimes$.
\begin{prop}\label{pr-pte_ti} Soient $F,G,F',G'$ des foncteurs de $\F(A,K)$. On dispose dans $\F(A^2,K)$ d'un isomorphisme naturel
$(F\boxtimes G)\ti(F'\boxtimes G')\simeq (F\ti F')\boxtimes (G\ti G')$ dans $\fct(\A,\E)$.
\end{prop}

(Dans le terme de gauche, c'est encore la structure $A^2$-linéaire de $\mathbf{P}(A)^2$ qui est utilisée.)

\begin{proof}
C'est une conséquence immédiate de l'isomorphisme naturel $(X\boxtimes Y)\underset{\C\times\D}{\otimes}(X'\boxtimes Y')\simeq (X\underset{\C}{\otimes} X')\otimes (Y\underset{\D}{\otimes} Y')$, où $\C$ et $\D$ sont des petites catégories et $X$, $X'$, $Y$, $Y'$ des foncteurs de $\F(\C^\op;K)$, $\F(\C;K)$, $\F(\D^\op;K)$ et $\F(\D;K)$ respectivement. 
\end{proof}

\section{Spectre d'un foncteur}

\subsection{Définition ; premières propriétés}

Nous introduisons maintenant l'une des notions nouvelles les plus importantes du présent mémoire, le \emph{spectre} d'un foncteur. La plupart de nos théorèmes principaux s'exprimeront en termes de spectre, ou se déduiront facilement de résultats concernant le spectre. La terminologie employée provient notamment de ce que, si $F$ est un foncteur de type fini de $\F^\df(\A;K)$ et que le corps $K$ est algébriquement clos, alors les valeurs propres de l'image par $F$ des endomorphismes diagonaux des $A$-modules $A^n$ peuvent se déduire du spectre de $F$, comme nous le verrons au chapitre suivant (cf. propositions \ref{pr-df_htr} et \ref{pr-poidspar_ac}).

\begin{defi}\label{defisp} Soit $F$ un foncteur de $\fct(\A,\E)$. On note \index{nota}{sp@$\Sp$, $\Sp_A$ \emph{(spectre d'un foncteur)}}$\Sp_A(F)$, ou simplement $\Sp(F)$, le quotient du foncteur $P^A$ de $\F(A,K)$ image du morphisme $P^A\xrightarrow{\eta_F}\hi(F,F)$ correspondant à l'élément $\mathrm{Id}_F\in\fct(\A,\E)(F,F)$ via les isomorphismes canoniques $\F(A,K)(P^A,\hi(F,F))\simeq\hi(F,F)(A)\simeq\fct(\A,\E)(F,F)$. On dit que $\Sp(F)$ est le \index{termin}{spectre \emph{(d'un foncteur)}}\textbf{spectre} du foncteur $F$ (relativement à $A$).
\end{defi}

On note que $\eta_F$ est adjoint à l'isomorphisme canonique $P^A\ti F\xrightarrow{\simeq} F$.

La notation $\Sp(F)$ désignera tantôt le foncteur correspondant de $\F(A,K)$, tantôt un élément de l'ensemble ordonné des quotients du foncteur $P^A$ ; nous ne préciserons explicitement ce qu'il en est que lorsque ce sera nécessaire pour éviter toute confusion.

\begin{rema} L'espace vectoriel $\Sp(F)(0)$ est nul si $F$ est réduit, et de dimension $1$ sinon.

L'espace vectoriel $\Sp(F)(A)$ est nul si et seulement si $F$ est nul ; il est de dimension $1$ si et seulement si $F$ est homogène d'un certain poids fort.
\end{rema}

On dispose de la caractérisation concrète suivante du spectre d'un foncteur :

\begin{prop}\label{pr-sp_concr}
Soient $F$ un foncteur de $\fct(\A,\E)$ et $n\in\mathbb{N}$.
\begin{enumerate}
    \item Comme quotient de $K[A^n]=P^A(A^n)$, $\Sp(F)(A^n)$ est l'image du morphisme canonique
    \begin{equation}\label{eqXi}
\Xi^{(n)}_F : K[A^n]\to\mathrm{End}_{\fct(\A^n,\E)}(F\circ\oplus_n)\qquad [a_1,\dots,a_n]\mapsto F(a_1\oplus\dots\oplus a_n).
\end{equation}
    \item Soient $r\in\mathbb{N}$, $(\lambda_i)_{1\le i\le r}\in K^r$ et $(a^i_j)_{1\le i\le r,1\le j\le n}\in A^{rn}$. Alors l'élément $\sum_{i=1}^r\lambda_i[a_1^i,\dots,a_n^i]$ de $K[A^n]=P^A(A^n)$ appartient au noyau de la projection $K[A^n]=P^A(A^n)\twoheadrightarrow\Sp(F)(A^n)$ si et seulement si, pour tout $n$-uplet $(f_1,\dots,f_n)$ de flèches parallèles de $\A$, on a $\sum_{i=1}^r\lambda_i F(\sum_{j=1}^n a_j^i f_j)=0$.
\end{enumerate}
\end{prop}

\begin{proof} La première assertion provient de l'isomorphisme
$$\hi(F,F)(A^n)\simeq\mathrm{Hom}(F,F\circ (A^n\underset{A}{\otimes}-))\simeq\mathrm{Hom}(F,F\circ \oplus_n\circ\delta_n)\simeq\mathrm{Hom}(F\circ\oplus_n,F\circ\oplus_n)$$
déduit de l'adjonction la somme itérée $\oplus_n$ et la diagonale itérée $\delta_n$, par lequel l'élément $\eta_F([a_1,\dots,a_n])$ est envoyé sur $F(a_1\oplus\dots\oplus a_n)$.

La deuxième assertion découle de la première.
\end{proof}

L'énoncé précédent, qui est à rapprocher des propositions~\ref{pr-radevtau} et~\ref{pr-rfevtau}, dit essentiellement que le spectre d'un foncteur code ses relations sur les morphismes qui proviennent de la structure $A$-linéaire à la source.

\begin{rema} Le spectre \emph{ne} définit \emph{pas} de foncteur de $\fct(\A,\E)$ dans $\F(A,K)$ (ni dans l'ensemble ordonné des quotients de $P^A$), comme on le déduit aisément de la deuxième assertion de la proposition~\ref{pr-sp_concr}. Certaines des propriétés que nous présenterons dans la section~\ref{parspcft} s'apparenteront toutefois à une forme de fonctorialité du spectre.
\end{rema}

\begin{prop}\label{pr-spdf} Soit $F$ un foncteur de $\fct(\A,\E)$. La projection $P^A\twoheadrightarrow\Sp(F)$ induit des isomorphismes
$$\hi(\Sp(F),F)\xrightarrow{\simeq}\hi(P^A,F)\xrightarrow{\simeq}F\quad\text{et}\quad F\simeq P^A\ti F\xrightarrow{\simeq}\Sp(F)\ti F.$$
De plus, $\Sp(F)$ est le plus petit quotient de $P^A$ possédant cette propriété : si $f : P^A\twoheadrightarrow T$ induit un isomorphisme $\hi(T,F)\to\hi(P^A,F)\xrightarrow{\simeq}F$ ou un isomorphisme $P^A\ti F\xrightarrow{\simeq}T\ti F$, alors $P^A\twoheadrightarrow\Sp(F)$ se factorise à travers $f$.
\end{prop}

\begin{proof}
Comme le morphisme $\eta_F$ est adjoint à l'isomorphisme canonique $P^A\ti F\xrightarrow{\simeq}F$, ce dernier se factorise par la projection $P^A\ti F\twoheadrightarrow\Sp(F)\ti F$, laquelle est donc un isomorphisme.

Si $P^A\twoheadrightarrow T$ induit un isomorphisme $F\simeq P^A\ti F\xrightarrow{\simeq}T\ti F$, alors son inverse induit par adjonction un morphisme $T\to\hi(F,F)$, et la composée $P^A\twoheadrightarrow T\to\hi(F,F)$ coïncide avec $\eta_F$, d'où la factorisation cherchée.

Les assertions relatives à $\hi$ s'établissent de façon analogue.
\end{proof}

\begin{prop}\label{pr-dfsp} Soit $F$ un foncteur de type fini de $\F(\A;K)$.
\begin{enumerate}
    \item Si $F$ appartient à $\F^\df(\A;K)$, alors $\Sp(F)$ appartient à $\F^\df(A,K)$.
    \item Supposons que la catégorie $\A$ vérifie la condition \textnormal{(MTF)}\index{nota}{MTF@(MTF) \emph{(condition de finitude sur les morphismes)}}. Si $\Sp(F)$ appartient à $\F^\df(A,K)$, alors $F$ appartient à $\F^\df(\A;K)$.
\end{enumerate}
\end{prop}

\begin{proof}
Si $F$ est un foncteur de type fini de $\F^\df(\A;K)$, la proposition~\ref{pr-vdffac} montre que $\hi(F,F)$ est à valeurs de dimensions finies, c'est donc a fortiori le cas de son sous-foncteur $\Sp(F)$.

Sous l'hypothèse (MTF), si $\Sp(F)$ appartient à $\F^\df(A,K)$, alors $\Sp(F)\ti F$ est à valeurs de dimensions finies, puisque $F$ est de type fini, grâce à la proposition~\ref{pr-valtdf}. La proposition~\ref{pr-spdf} permet alors de conclure.
\end{proof}

\subsection{Lien avec l'algèbre caractéristique}\label{par-spac}

Soit $F$ un foncteur de $\fct(\A,\E)$. Notons $\Theta : K[A]\to K[A^2]$ l'application $K$-linéaire telle que $\Theta([a])=[1,a]$. Il résulte des propositions~\ref{pr-rfevtau} (équivalence \ref{itrfev1}$\Leftrightarrow$\ref{itrfeveq}) et~\ref{pr-sp_concr} (deuxième point) qu'un élément $\xi$ de $K[A]$ appartient à $\ic(F)$ si et seulement si $\Theta(\xi)$ appartient au noyau de la projection $K[A^2]=P^A(A^2)\twoheadrightarrow\Sp(F)(A^2)$ ; en particulier, \emph{le spectre d'un foncteur (et même sa seule évaluation sur $A^2$) détermine son algèbre caractéristique}.

Nous allons maintenant voir que, réciproquement, les algèbres caractéristiques des différents foncteurs $F\circ\oplus_n$ de $\fct(\A^n,\E)$, pour $n\in\mathbb{N}$, déterminent au moins l'effet sur les objets de $\Sp(F)$. Le point de départ de ces considérations est le morphisme $\Xi^{(n)}_F : K[A^n]\to\mathrm{End}_{\fct(\A^n,\E)}(F\circ\oplus_n)$ (cf. \eqref{eqXi}) de la proposition~\ref{pr-sp_concr}.

Ce morphisme, a priori $K$-linéaire, définit en fait un morphisme d'algèbres $K[A_\mu^n]\to Z(\mathrm{End}_{\fct(\A^n,\E)}(F\circ\oplus_n))$. En particulier, $\Sp(F)(A^n)$ est naturellement une sous-algèbre de $Z(\mathrm{End}_{\fct(\A^n,\E)}(F\circ\oplus_n))$ (et un quotient de $K[A_\mu^n]$) --- $\Sp(F)(A^n)$ s'identifie à l'algèbre $\mathfrak{Z}_{A^n_\mu}(F\circ\oplus_n)$, selon la notation~\ref{nota-cZ}. On prendra toutefois garde que $\Sp(F)$ ne prend pas \textbf{naturellement} ses valeurs dans les $K$-algèbres : la structure d'algèbre ainsi obtenue sur $\Sp(F)(V)$, où $V$ est un $A$-module libre de rang fini, dépend du choix d'une base de $V$.

Pour $1\le i\le n$, notons $\Upsilon_i^n : \A^n\to\A^2$ le foncteur $(a_1,\dots,a_n)\mapsto (\bigoplus_{j\ne i} a_j,a_i)$. On dispose d'un isomorphisme évident $\oplus\circ\Upsilon^i_n\simeq\oplus_n$, ainsi la précomposition par $\Upsilon^i_n$ induit un morphisme d'algèbres (non nécessairement commutatives)
$$\alpha^n_i : \mathrm{End}_{\fct((\A^2,\E)}(F\circ\oplus)\to\mathrm{End}_{\fct(\A^n,\E)}(F\circ\oplus_n)$$
(qui pour $n\ge 2$ est injectif, car $\Upsilon^i_n$ est alors un épimorphisme scindé). On a ainsi
\begin{equation}\label{eq-chin}
\Xi^{(n)}_F([a_1,\dots,a_n])=\prod_{i=1}^n\alpha^n_i(\chi_F(a_i)).
\end{equation}

Par conséquent :
\begin{prop}\label{pr-valsp} Pour tout $n\in\mathbb{N}$ et tout foncteur $F$ de $\fct(\A,\E)$, on a $\Sp(F)(A^n)\simeq\mathfrak{Z}_{A^n_\mu}(F\circ\oplus_n)$ ; cette $K$-algèbre est un quotient de $\ac(F)^{\otimes n}$.
\end{prop}

De plus, $\ac(F)$ s'identifie manifestement à une sous-algèbre de $\Sp(F)(A^2)$ (et lui est isomorphe si $F$ est homogène d'un certain poids fort). Il en découle :

\begin{coro}\label{cor-spchi} Soit $F$ un foncteur de $\fct(\A,\E)$. Les assertions suivantes sont équivalentes :
\begin{enumerate}
\item\label{itspac1} l'algèbre caractéristique $\ac(F)$ est de dimension finie sur $K$ ;
\item\label{itspac2}  $\Sp(F)$ appartient à $\F^\df(A,K)$ ;
\item\label{itspac3} $\Sp(F)(A^2)$ est un $K$-espace vectoriel de dimension finie.
\end{enumerate}
\end{coro}

\begin{proof}
L'implication \ref{itspac1}$\Rightarrow$\ref{itspac2} résulte de la proposition~\ref{pr-valsp}. L'implication \ref{itspac2}$\Rightarrow$\ref{itspac3} est triviale, et \ref{itspac3}$\Rightarrow$\ref{itspac1} provient également de la proposition~\ref{pr-valsp} car $\ac(F)$ est par définition une sous-algèbre de $\mathfrak{Z}_{A^2_\mu}(F\circ\oplus)$.
\end{proof}

En combinant le corollaire précédent avec la proposition~\ref{pr-dfsp}, il vient :

\begin{coro}\label{cor-spchi2} Supposons que $\A$ vérifie la condition \textnormal{(MTF)}\index{nota}{MTF@(MTF) \emph{(condition de finitude sur les morphismes)}}. Soit $F$ un foncteur de type fini de $\F(\A;K)$. Alors les assertions suivantes sont équivalentes :
\begin{enumerate}
\item $\dim_K\ac(F)<\infty$ ;
\item  $\Sp(F)$ appartient à $\F^\df(A,K)$ ;
\item $F$ appartient à $\F^\df(\A;K)$.
\end{enumerate}
\end{coro}

Dans l'énoncé suivant, qui découle directement des considérations du début de ce §\ref{par-spac} (ou même simplement de la première assertion de la proposition~\ref{pr-sp_concr}), on identifie les catégories $\F(\mathbf{P}(A)^n,K)$ et $\F(A^n,K)$ à l'aide de l'équivalence induite par l'inclusion $\mathbf{P}(A^n)\to\mathbf{P}(A)^n$.

\begin{prop}\label{pr-spcr} Soient $F$ un foncteur de $\fct(\A,\E)$ et $n\in\mathbb{N}$. On dispose dans $\F(A^n,K)$ d'isomorphismes naturels $\Sp_A(F)\circ\oplus_n\simeq\Sp_{A^n}(F\circ\oplus_n)$ et $cr_n(\Sp_A(F))\simeq\Sp_{A^n}(cr_n(F))$.
\end{prop}

\begin{coro}\label{cor-sppol} Soient $F$ un foncteur de $\fct(\A,\E)$ et $d\in\mathbb{N}$. Alors $F$ appartient à $\pol_d(\A,\E)$ si et seulement si $\Sp(F)$ appartient à $\pol_d(A,K)$.
\end{coro}

Nous verrons ultérieurement (corollaire~\ref{cor-phsp}) que la propriété Hom-polynomiale est également détectée par le spectre, sous des hypothèses raisonnables.

\section{Les catégories $\fct_T(\A,\E)$}\label{parspcft}

Il est généralement très difficile de dire quoi que ce soit du treillis des sous-foncteurs d'un foncteur donné de $\fct(\A,\E)$, ou même seulement de $\F(A,K)$. Toutefois, pour un foncteur $F$ dont le spectre $T$ possède de bonnes propriétés, la sous-catégorie $\F_T(A,K)$ que nous allons définir ci-dessous se ramène parfois à des catégories de foncteurs beaucoup plus simples ; comme c'est une sous-catégorie prélocalisante de $\F(A,K)$ qui contient $F$, le treillis des sous-foncteurs de $F$ est le même dans $\F(A,K)$ et dans $\F_T(A,K)$. En particulier, si $\F_T(A,K)$ est localement noethérienne et $F$ de type fini, alors $F$ est noethérien. C'est ainsi que nous montrerons, au chapitre~\ref{stfdf}, que la catégorie $\F^\df(A,K)$ est localement noethérienne, à partir de l'étude du spectre des foncteurs de type fini de $\F^\df(A,K)$.

\subsection{Définition, généralités}\label{pfTg}

\begin{prop}\label{pr-spqp} Soit $X$ un quotient de $P^A$ dans $\F(A,K)$. Alors $X$ est égal à son spectre.
\end{prop}

\begin{proof} De façon générale, si $\alpha : F\to G$ est un morphisme de $\F(A,K)$, le diagramme
$$\xymatrix{P^A\ar[rr]^-{\eta_F}\ar[d]_-{\eta_G} & & \hi(F,F)\ar[d]^{\hi(F,\alpha)} \\
\hi(G,G)\ar[rr]^-{\hi(\alpha,G)} & & \hi(F,G)
}$$
commute. Si $q : P^A\twoheadrightarrow X$ désigne la projection, on a donc un diagramme commutatif
$$\xymatrix{P^A\ar[rr]^-{\eta_{P^A}}_-\simeq\ar[d]_-{\eta_X} & & \hi(P^A,P^A)\ar[d]_{\hi(P^A,q)}\ar[r]^-\simeq & P^A\ar@{->>}[d]^q \\
\hi(X,X)\ar@{^(->}[rr]^-{\hi(q,X)} & & \hi(P^A,X) \ar[r]^-\simeq & X
}$$
qui donne la conclusion.
\end{proof}

\begin{coro}\label{cor-acsp} Pour tout foncteur $F$ de $\fct(\A,\E)$, on a $\rf(\Sp(F))=\rf(F)$, $\ac(\Sp(F))=\ac(F)$ et $\Sp(\Sp(F))=\Sp(F)$.
\end{coro}

(L'égalité des algèbres caractéristiques s'entend comme égalité de quotients de l'algèbre $K[A_\mu]$, elle est équivalente à l'égalité des idéaux caractéristiques ; l'égalité des spectres s'entend comme égalité de quotients du foncteur $P^A$ de $\F(A,K)$.)

\begin{proof}
Comme $\Sp(F)$ est un quotient de $P^A$, l'égalité $\Sp(\Sp(F))=\Sp(F)$ découle de la proposition~\ref{pr-spqp}. Comme l'algèbre caractéristique d'un foncteur est déterminée par son spectre (cf. début du §\,\ref{par-spac}), l'égalité $\ac(\Sp(F))=\ac(F)$ s'en déduit. Enfin, le radical d'un foncteur est déterminé par son idéal caractéristique (grâce aux propositions~\ref{pr-radevtau} et~\ref{pr-rfevtau}), de sorte que cette dernière égalité entraîne $\rf(\Sp(F))=\rf(F)$.
\end{proof}

\begin{defi} Soit $T$ un quotient de $P^A$. On note $\fct_T(\A,\E)$\index{nota}{fct@$\fct_T$} la sous-catégorie pleine de $\fct(\A,\E)$ constituée des foncteurs $F$ tels que $\Sp(F)\le T$. Lorsque $\A=\mathbf{P}(A)$ et $\E=K\Md$, on notera simplement $\F_T(A,K)$\index{nota}{FB@$\F_T$} cette sous-catégorie.
\end{defi}

D'après la proposition~\ref{pr-spdf}, un foncteur $F$ de $\fct(\A,\E)$ appartient à $\fct_T(\A,\E)$ si et seulement si l'épimorphisme canonique $F\twoheadrightarrow T\ti F$ est un isomorphisme, ou encore si et seulement si le monomorphisme canonique $\hi(T,F)\hookrightarrow F$ est un isomorphisme. On en déduit en particulier les deux énoncés suivants.

\begin{prop}\label{pr-adj_FT} Soit $T$ un quotient de $P^A$. Alors les endofoncteurs $\hi(T,-)$ et $T\ti -$ de $\fct(\A,\E)$ sont à valeurs dans $\fct_T(\A,\E)$. Les foncteurs $\fct(\A,\E)\to\fct_T(\A,\E)$ qu'ils induisent sont adjoints respectivement à droite et à gauche au foncteur d'inclusion $\fct_T(\A,\E)\hookrightarrow\fct(\A,\E)$.
\end{prop}

(Les morphismes canoniques précédents ne sont autres que l'unité et la coünité des adjonctions.)

\begin{prop}\label{pr-catfindT} La sous-catégorie $\fct_T(\A,\E)$ de $\fct(\A,\E)$ est prébilocalisante et stable par l'opération $F\mapsto F_V$ pour tout $V\in\mathbf{P}(A)$.

De plus, pour tout foncteur $X$ de $\F(A,K)$, $\fct_T(\A,\E)$ est stable par $\hi(X,-)$ et $X\ti -$.
\end{prop}

On dispose de la réciproque suivante :

\begin{prop}\label{pr-catfindT2} Soit $\G$ une sous-catégorie prébilocalisante\index{termin}{prebilocalisante@prébilocalisante \emph{(sous-catégorie)}} de $\F(A,K)$ stable par $F\mapsto F_V$ pour tout $V\in\mathbf{P}(A)$. Alors il existe un unique quotient $T$ de $P^A$ tel que $\G=\F_T(A,K)$.
\end{prop}

\begin{proof} Soit $q : \F(A,K)\to\G$ l'adjoint à gauche de l'inclusion : $q(X)$ est le plus grand quotient de $X$ appartenant à $\G$, pour tout $X\in\mathrm{Ob}\,\F(A,K)$. Si $\G=\F_T(A,K)$, on a nécessairement $T=q(P^A)$, d'où l'unicité de $T$. Vérifions que $\G=\F_{q(P^A)}(A,K)$.

Si $F$ appartient à $\F_{q(P^A)}(A,K)$, alors $F\simeq q(P^A)\ti F$ (proposition~\ref{pr-spdf}). Comme $F$ est un quotient d'une somme directe de foncteurs de la forme $P^V$, il s'ensuit que $F$ est quotient d'une somme directe, donc sous-quotient d'un produit direct, de foncteurs de la forme $q(P^A)\ti P^V\simeq q(P^A)_{V^*}$, et appartient donc à $\G$.

Réciproquement, si $F$ est un foncteur de $\G$, alors $\hi(F,F)$ est un sous-foncteur d'un produit de foncteurs du type $\hi(P^V,F)\simeq F_V$, donc appartient à $\G$. Par conséquent, $\Sp(F)$, qui s'injecte dans $\hi(F,F)$, appartient à $\G$. Comme $q(P^A)$ est le plus grand quotient de $P^A$ appartenant à $\G$, on a $\Sp(F)\le q(P^A)$ i.e. $F\in\F_{q(P^A)}(A,K)$, d'où la proposition.
\end{proof}

Une structure monoïdale symétrique sur un foncteur $F$ de $\F(A,K)$ est la donnée de morphismes $F(U)\otimes F(V)\to F(U\otimes_A V)$ naturels en $U$ et $V$ et d'un élément de $F(A)$ vérifiant les conditions d'unité, associativité et symétrie usuelles. Par adjonction, il revient au même de se donner des morphismes $F\ti F\to F$ et $P^A\to F$ de $\F(A,K)$ faisant de $F$ un monoïde commutatif dans la catégorie monoïdale symétrique $(\F(A,K),\ti,P^A)$.

Si $X$ est un foncteur de $\F(A,K)$ muni d'une structure monoïdale symétrique, son extension de Kan à gauche $\tilde{X} : A\Md\to K\Md$ le long de l'inclusion $\mathbf{P}(A)\to A\Md$ hérite d'une structure monoïdale symétrique. On peut définir une catégorie $K$-linéaire $X\A$ ayant les mêmes objets que $\A$, dont les morphismes sont donnés par $(X\A)(a,b)=\tilde{X}(\A(a,b))$, où la composition est donnée par
$$\tilde{X}(\A(a,b))\otimes\tilde{X}(\A(b,c))\to\tilde{X}(\A(a,b)\otimes_A\A(b,c))\to\tilde{X}(\A(a,c)),$$
la première flèche venant de la structure monoïdale de $\tilde{X}$ et la seconde de la composition de $\A$. On définit de même l'unité.

Si $X$ est un quotient de $P^A$, les propositions~\ref{pr-spdf} et~\ref{pr-spqp} fournissent un isomorphisme $X\ti X\xrightarrow{\simeq}X$ qui définit une structure monoïdale symétrique sur $X$. La projection $P^A\twoheadrightarrow X$ en constitue l'unité ; elle est en particulier compatible aux structures monoïdales.

\begin{prop}\label{pr-fctX_d2} Soit $X$ un quotient de $P^A$. La précomposition par le foncteur $K$-linéaire $\Phi : K[\A]=P^A\A\twoheadrightarrow X\A$ induit par la transformation naturelle monoïdale $P^A\twoheadrightarrow X$ identifie la sous-catégorie $\fct_X(\A,\E)$ de $\fct(\A,\E)$ (identifiée à la catégorie des foncteurs $K$-linéaires de $K[\A]$ vers $\E$) à la catégorie des foncteurs $K$-linéaires de $X\A$ vers $\E$.
\end{prop}

\begin{proof}
Comme $\Phi$ est plein et essentiellement surjectif, $\Phi^*$ identifie la catégorie des foncteurs $K$-linéaires de $X\A$ vers $\E$ à la sous-catégorie pleine des foncteurs $F$ de $\fct(\A,\E)$ tels que, pour tous objets $x, y$ de $\A$ et tout élément $\xi$ de $K[\A(x,y)]$ dont l'image dans $\tilde{X}(\A(x,y))$ est nulle, $F(\xi)=0$. En considérant la suite exacte canonique $\bigoplus_{n\in\mathbb{N},\alpha\in\mathrm{Ker}(P^A(A^n)\twoheadrightarrow X(A^n))} P^{A^n}\to P^A\to X\to 0$ de $\F(A,K)$, on voit qu'il revient au même de demander que, pour tout élément $\alpha=\sum_{i=1}^r\lambda_i[a_1^i,\dots,a_n^i]$ de $\mathrm{Ker}\,(K[A^n]=P^A(A^n)\twoheadrightarrow X(A^n))$ (avec $\lambda_i\in K$ et $a_j^i\in A$), le morphisme $K[\A(x,y)^n]\xrightarrow{\alpha^*}K[\A(x,y)]\xrightarrow{F_*}\E(F(x),F(y))$ soit nul. Si $f_1,\dots,f_n$ sont des éléments de $\A(x,y)$, on a $\alpha^*([f_1,\dots,f_n])=\sum_{i=1}^r\lambda_i F(\sum_{j=1}^n a_j^i f_j)$, de sorte que la deuxième assertion de la proposition~\ref{pr-sp_concr} fournit donc la conclusion.
\end{proof}

La description des catégories $\fct_X(\A,\E)$ donnée par la proposition précédente nous servira, au chapitre~\ref{shd}, à d'établir un lien entre spectre et structures polynomiales {\em strictes} généralisées.

\subsection{Comparaison de catégories de la forme $\fct_T(\A,\E)$}

La propriété suivante constitue une conséquence directe de la proposition~\ref{pr-spdf} et de l'associativité du produit tensoriel $\ti$.

\begin{prop}\label{pr-intersect_FT}
Soient $T$ et $U$ des quotients de $P^A$. Alors $T\ti U$ est la borne inférieure de $\{T,U\}$ dans l'ensemble ordonné des quotients de $P^A$. La sous-catégorie $\fct_{T\ti U}(\A,\E)$ de $\fct(\A,\E)$ est l'intersection des sous-catégories $\fct_T(\A,\E)$ et $\fct_U(\A,\E)$.
\end{prop}

\begin{prop}\label{pr-sp_sdir}
Soient $T$ et $U$ des quotients de $P^A$. Les assertions suivantes sont équivalentes :
\begin{enumerate}
    \item le morphisme $P^A\to T\oplus U$ dont les composantes sont les projections $P^A\twoheadrightarrow T$ et $P^A\twoheadrightarrow U$ est un épimorphisme ;
    \item le foncteur $T\ti U$ de $\F(A,K)$ est nul.
\end{enumerate}

Lorsque ces conditions sont vérifiées, le foncteur $\fct_T(\A,\E)\times\fct_U(\A,\E)\to\fct(\A,\E)$ composée de l'inclusion $\fct_T(\A,\E)\times\fct_U(\A,\E)\to\fct(\A,\E)\times\fct(\A,\E)$ et de la somme directe sur $\fct(\A,\E)$ induit une équivalence de catégories  $\fct_T(\A,\E)\times\fct_U(\A,\E)\xrightarrow{\simeq}\fct_{T\oplus U}(\A,\E)$.
\end{prop}

\begin{proof}
Si le morphisme $P^A\to T\oplus U$ est un épimorphisme, comme $P^A\twoheadrightarrow T$ se factorise à travers celui-là, la proposition~\ref{pr-spdf} et le corollaire~\ref{cor-acsp} montrent que l'épimorphisme canonique $T\simeq P^A\ti T\twoheadrightarrow (T\oplus U)\ti T\simeq (T\ti T)\oplus (T\ti U)$ est un isomorphisme. Comme la première composante de ce morphisme est l'isomorphisme canonique $T\simeq T\ti T$ (donné également par la proposition~\ref{pr-spdf} et le corollaire~\ref{cor-acsp}), il s'ensuit que $T\ti U$ est nul.

Supposons réciproquement $T\ti U$ nul et notons $X$ le conoyau du morphisme $P^A\to T\oplus U$. Comme $P^A\to T$ est un épimorphisme, $X$ est un quotient de $U$, de sorte que $X\ti U\simeq X$. De même, $X\ti T\simeq X$. Il s'ensuit que $X\simeq X\ti (T\ti U)$ est nul. De plus, si $F$ (resp. $G$) est un foncteur de $\fct_T(\A,\E)$ (resp. $\fct_U(\A,\E)$), alors 
$$\hi(F,G)\simeq\hi(T\ti F,\hi(U,G))\simeq\hi(U\ti T\ti F,G)=0$$
et en particulier $\fct(\A,\E)(F,G)\simeq\hi(F,G)(A)$ est nul ; de même $\fct(\A,\E)(G,F)=0$. On en déduit que la somme directe induit un foncteur pleinement fidèle $\fct_T(\A,\E)\times\fct_U(\A,\E)\to\fct_{T\oplus U}(\A,\E)$. Montrons son essentielle surjectivité : si $F$ est un foncteur de $\fct_{T\oplus U}(\A,\E)$, on a $F\simeq (T\oplus U)\ti F\simeq (T\ti F)\oplus (U\ti F)$, or $T\ti F$ appartient à $\fct_T(\A,\E)$ et $U\ti F$ à $\fct_U(\A,\E)$, d'où la conclusion.
\end{proof}

\begin{rema} On vérifie aussitôt que, de manière générale, si $T$ et $U$ sont des quotients de $P^A$, l'image $S$ du morphisme canonique $P^A\to T\oplus U$ est leur borne supérieure dans l'ensemble ordonné des quotients de $P^A$. Il ne semble toutefois pas exister de description simple de $\fct_S(\A,\E)$ à partir de $\fct_T(\A,\E)$ et $\fct_U(\A,\E)$ hors de la situation de la proposition~\ref{pr-sp_sdir}.
\end{rema}

Dans l'énoncé suivant, on identifie $\F(A^2,K)$ à $\F(\mathbf{P}(A)^2;K)$.

\begin{prop}\label{pr-spec_diag}
Soit $X$ un quotient de $P^{A\times A}\simeq (P^A)^{\boxtimes 2}$ dans $\F(A^2,K)$.  Les assertions suivantes sont équivalentes :
\begin{enumerate}
    \item le morphisme $\varphi : P^A\to (P^A)^{\otimes 2}\simeq (P^A)^{\boxtimes 2}\circ\delta\twoheadrightarrow X\circ\delta$ (où $P^A\to (P^A)^{\otimes 2}$ est le coproduit) de $\F(A,K)$ est un épimorphisme ;
    \item la restriction du foncteur $\delta^* : \F(A^2,K)\to\F(A,K)$ à la sous-catégorie $\F_X(A^2,K)$ est pleinement fidèle et son image essentielle est stable par quotient.
\end{enumerate}
Lorsqu'elles sont vérifiées, le foncteur $\delta^* : \fct(\A^2,\E)\to\fct(\A,\E)$ induit une équivalence de catégories $\fct_X(\A^2,\E)\xrightarrow{\simeq}\fct_{\delta^*X}(\A,\E)$.
\end{prop}

\begin{proof} En composant l'adjonction somme/diagonale avec celle de la proposition~\ref{pr-adj_FT}, on voit que le foncteur $\Phi :\fct_X(\A^2,\E)\hookrightarrow\fct(\A^2,\E)\xrightarrow{\delta^*}\fct(\A,\E)$ est adjoint à droite au foncteur $\Psi : F\mapsto X\ti\oplus^*F$ et que l'unité $\eta$ de cette adjonction s'identifie au morphisme
$F\simeq P^A\ti F\to\delta^*(X\ti\oplus^*F)\simeq\delta^*X\ti F$ (grâce à la proposition~\ref{pr-asdti}) obtenu en appliquant $-\ti F$ au morphisme $\varphi : P^A\to\delta^*X$.

Si le foncteur $\Phi$ est pleinement fidèle et que son image essentielle est stable par quotient, alors $\eta$ est un épimorphisme ($\eta_F$ est la projection de $F$ sur son plus grand quotient appartenant à cette image essentielle). En particulier, la deuxième condition implique la première.

Réciproquement, supposons que $\varphi$ est un épimorphisme, ce qui entraîne que $\eta$ est un épimorphisme. Si $\varepsilon$ désigne la coünité de l'adjonction, pour tout objet $G$ de $\fct_X(\A^2,\E)$, la composée $\Phi G\xrightarrow{\eta_{\Psi(G)}}\Phi\Psi\Phi G\xrightarrow{\Phi(\varepsilon_G)}\Phi G$ est l'identité. Comme $\eta_{\Psi(G)}$ est un épimorphisme, il s'ensuit que c'est un isomorphisme, de même que $\Phi(\varepsilon_G)$. Puisque $\Phi$ est exact et fidèle, cela montre que $\varepsilon$ est un isomorphisme. Cela entraîne pour des raisons formelles que $\Phi$ est pleinement fidèle, et que son image essentielle est constituée des objets $F$ de $\fct(\A,\E)$ tels que $\eta_F$ soit un isomorphisme, c'est-à-dire $\fct_{\delta^*X}(\A,\E)$, qui est une sous-catégorie de $\fct(\A,\E)$  stable par quotient. Cela termine la démonstration.
\end{proof}

\begin{prop}\label{pr-sp_pte} Soient $T$ et $U$ des quotients de $P^A$. L'image essentielle de la sous-catégorie $\fct_{T\boxtimes U}(\A\times\A,\E)$ de $\fct(\A\times\A,\E)$ par l'équivalence canonique $\fct(\A\times\A,\E)\xrightarrow{\simeq}\fct(\A,\fct(\A,\E))$ égale $\fct_T(\A,\fct_U(\A,\E))$.
\end{prop}

\begin{proof}
Il découle par exemple de la proposition~\ref{pr-sp_concr} que l'image de $\fct_{T\boxtimes P^A}(\A\times\A,\E)$ (resp. $\fct_{P^A\boxtimes U}(\A\times\A,\E)$) par $\fct(\A\times\A,\E)\xrightarrow{\simeq}\fct(\A,\fct(\A,\E))$ égale $\fct_T(\A,\fct(\A,\E))$ (resp. $\fct(\A,\fct_U(\A,\E))$). Comme $(T\boxtimes P^A)\ti(P^A\boxtimes U)\simeq T\ti U$ par la proposition~\ref{pr-pte_ti}, la conclusion découle de la proposition~\ref{pr-intersect_FT}.
\end{proof}

Le résultat suivant découle des propositions~\ref{pr-spec_diag} et~\ref{pr-sp_pte}.

\begin{coro}\label{cor-sp_pt}
Soient $T$ et $U$ des quotients de $P^A$.
 Les assertions suivantes sont équivalentes :
\begin{enumerate}
    \item le morphisme $P^A\to (P^A)^{\otimes 2}\twoheadrightarrow T\otimes U$ (où $P^A\to (P^A)^{\otimes 2}$ est le coproduit) de $\F(A,K)$ est un épimorphisme ;
    \item la restriction du foncteur $\delta^* : \F(A^2,K)\to\F(A,K)$ à la sous-catégorie $\fct_T(\mathbf{P}(A),\F_U(A,K))$ de $\fct(\mathbf{P}(A),\F(A,K))\simeq\F(A\times A,K)$ est pleinement fidèle et son image essentielle est stable par quotient.
\end{enumerate}
Lorsqu'elles sont vérifiées, le foncteur $\fct(\A,\fct(\A,\E))\simeq\fct(\A^2,\E)\xrightarrow{\delta^*}\fct(\A,\E)$ induit une équivalence de catégories $\fct_T(\A,\fct_U(\A,\E))\xrightarrow{\simeq}\fct_{T\otimes U}(\A,\E)$.
\end{coro}

\subsection{Quelques exemples fondamentaux}

\begin{exem} Le corollaire~\ref{cor-sppol} peut se reformuler par l'égalité $\pol_d(\A,\E)=\fct_{\qpol_d(P^A)}(\A,\E)$ pour tout $d\in\mathbb{N}$.
\end{exem}

\begin{exem}\label{ex-radSpS} Si $I$ est un idéal de $A$, alors $\fct_{\pi_I^*P^{A/I}}(\A,\E)\simeq\fct(\A/I,\E)$ (où $\pi_I : \mathbf{P}(A)\to\mathbf{P}(A/I)$ désigne la réduction modulo $I$) est la sous-catégorie de $\fct(\A,\E)$ des foncteurs $F$ tels que $\rf(F)\supset I$. 
\end{exem}

On peut combiner les deux exemples précédents pour obtenir :
\begin{prop} Soient $I$ un idéal de $A$ tel que $(A/I)\underset{\mathbb{Z}}{\otimes}K=0$ et $d\in\mathbb{N}$. Alors le morphisme $P^A\to P^A\otimes P^A\twoheadrightarrow\pi_I^*(P^{A/I})\otimes\qpol_d(P^A)$ composé du coproduit et du produit tensoriel des projections est un épimorphisme. De plus, $\fct_{\pi_I^*(P^{A/I})\otimes\qpol_d(P^A)}(\A,\E)$ s'identifie à $\fct(\A/I,\pol_d(\A,\E))$ (catégorie plongée dans $\fct(\A\times\A,\E)$ par précomposition par la diagonale). 
\end{prop}

\begin{proof}
Le fait que $P^A\to P^A\otimes P^A\twoheadrightarrow\pi_I^*(P^{A/I})\otimes\qpol_d(P^A)$ soit un épimorphisme découle, par dualité, de \cite[lemme~4.20]{DTV}. Le reste s'en déduit alors par le corollaire~\ref{cor-sp_pt} et les exemples précédents.
\end{proof}

\begin{coro} Soit $F$ un foncteur à support fini de $\F(A,K)$. Les assertions suivantes sont équivalentes :
\begin{enumerate}
\item $F$ possède une décomposition à la Steinberg de type $PAP$ ;
\item $\Sp(F)$ possède une décomposition à la Steinberg de type $PAP$ ;
\item il existe un idéal $K$-cotrivial $I$ de $A$ et un entier naturel $d$ tels que $\Sp(F)\le\pi_I^*(P^{A/I})\otimes\qpol_d(P^A)$.
\end{enumerate}
\end{coro}

\begin{exem}\label{exidbiid} Si $\I$ est un bi-idéal de $K[A]$, alors $\fct_{P^A_\I}(\A,\E)$ est la sous-catégorie des foncteurs $F$ de $\fct(\A,\E)$ tels que $\ic(F)\supset\I$ (cf. proposition~\ref{pr-biid}).
\end{exem}

\begin{exem} Soit $w\in\mon(A_\mu,K_\mu)$ un poids. On a  $\fct_{P^A_w}(\A,\E)=\fct(\A,\E)_w$.  Plus généralement, pour $n\in\mathbb{N}$, la sous-catégorie $U(\fct(\A,\E),\mathfrak{m}_w^n)$ de $\fct(\A,\E)$ égale $\fct_{P^A_{w,n}}(\A,\E)$, où $P^A_{w,n}$ désigne le plus grand quotient de $P^A$ appartenant à $U(\F(A,K),\mathfrak{m}_w^n)$.
\end{exem}

\begin{exem}\label{ex-sdirpq} Si $w$ et $w'$ sont deux poids distincts de $\mon(A_\mu,K_\mu)$, alors $P^A_w\ti P^A_{w'}=0$ (s'il était non nul, ce foncteur serait homogène de poids fort à la fois $w$ et $w'$, d'après l'exemple précédent). Il découle donc d'une application itérée de la proposition~\ref{pr-sp_sdir} et de l'exemple précédent que si $w_1,\dots,w_r$ est une famille finie de poids deux à deux distincts de $\mon(A_\mu,K_\mu)$, alors $\fct_{P^A_{w_1}\oplus\dots\oplus P^A_{w_r}}(\A,\E)\simeq\fct(\A,\E)_{w_1}\times\dots\times\fct(\A,\E)_{w_r}$ s'identifie à la sous-catégorie de $\fct(\A,\E)$ des foncteurs $F$ admettant une décomposition en poids forte et tels que $\Pi(F)\subset\{w_1,\dots,w_r\}$.

De même, plus généralement, étant donné $n\in\mathbb{N}$, on a $\fct_{P^A_{w_1,n}\oplus\dots\oplus P^A_{w_r,n}}(\A,\E)\simeq U(\fct(\A,\E),\mathfrak{m}_{w_1}^n)\times\dots\times U(\fct(\A,\E),\mathfrak{m}_{w_r}^n)$.
\end{exem}

L'exemple précédent implique directement le résultat suivant.

\begin{prop}\label{pr-sppoi} Soient $F$ un foncteur de $\fct(\A,\E)$ et $w\in\mon(A_\mu,K_\mu)$ un poids.
\begin{enumerate}
\item Le foncteur $F$ de $\fct(\A,\E)$ est homogène de poids fort $w$ si et seulement si le foncteur $\Sp(F)$ de $\F(A,K)$ est homogène de poids fort $w$.
\item Le foncteur $F$ possède une décomposition en poids forte (resp. faible) finie si et seulement si $\Sp(F)$ en possède une. Le cas échéant, on a $\Pi(F)=\Pi(\Sp(F))$.
\end{enumerate}
\end{prop}

\begin{rema} La deuxième assertion tombe en défaut si l'on omet l'hypothèse de finitude des décompositions en poids.

Considérons par exemple le cas $A=K=\mathbb{Q}$ et la catégorie $\F(\mathbb{Q},\mathbb{Q})$. Le foncteur $F:=\bigoplus_{d\in\mathbb{N}}S^d$ (où $S^d$ désigne la $d$-ième puissance symétrique) possède une décomposition en poids forte, et son ensemble de poids est $\{x\mapsto x^d\,|\,d\in\mathbb{N}\}\subsetneq\mon(\mathbb{Q}_\mu,\mathbb{Q}_\mu)$. On vérifie aisément que $\Sp(F)=P^\mathbb{Q}$, foncteur qui ne possède pas de décomposition en poids (même faible). De plus, $\Pi(P^\mathbb{Q})$ est réduit au poids trivial (mais tout élément de $\mon(\mathbb{Q}_\mu,\mathbb{Q}_\mu)$ est le poids d'un \emph{quotient} de $\Pi(P^\mathbb{Q})$).
\end{rema}

\chapter{Poids et effets croisés}\label{spefcr}

\begin{cvi}
Dans ce chapitre encore, $\A$ désigne une catégorie additive $A$-linéaire essentiellement petite et $\E$ une catégorie de Grothendieck $K$-linéaire.
\end{cvi}

Hors des foncteurs polynomiaux (cf. chapitre~\ref{sec-finpolp}), l'existence d'une décomposition en poids pour un foncteur $F$ de $\fct(\A,\E)$ ne permet pas de dire beaucoup sur sa structure. Le lemme~\ref{lm-tauP} en donne une illustration : remplacer le foncteur $P^A$ de $\F(A,K)$ par son plus grand quotient homogène d'un certain poids fort $w$ nul sur les éléments non inversibles de $A$ ne diminue pas réellement sa complexité. En revanche, l'existence d'une décomposition en poids pour le bifoncteur $F\circ\oplus$ (relativement à l'action de $A_\mu\times A_\mu$) implique des propriétés de structure et de finitude remarquables, comme nous le verrons ultérieurement. Dans le présent chapitre, nous donnons d'abord les définitions et premières propriétés relatives aux foncteurs $F$ dits \emph{\hts}\ ou \emph{\hds}, c'est-à-dire tels que $F\circ\oplus$ possède une décomposition en poids faible ou forte respectivement. Cela nous permettra, en utilisant également les notions développées au chapitre~\ref{chap-spec}, de démontrer à la section~\ref{safphs} que tout foncteur \ph\ de type fini de $\fct(\A,\E)$ est polynomial lorsque l'anneau $A$ n'a pas de quotient fini de même caractéristique que $K$ (corollaire~\ref{cor-EMLpol-pol}).

\section[Foncteurs \hds\ et \hts]{Définition des foncteurs \hds\ et \hts}

Les propositions~\ref{def-htr} et~\ref{defo} ci-après découlent d'une application itérée de la proposition~\ref{pr-poibif}, de la proposition~\ref{pr-tridia-gl} et de la décomposition de la précomposition d'un foncteur par $\oplus_n$ à l'aide des effets croisés \cite[th.~9.1]{EML}.

\begin{prdef}\label{def-htr} Soit $F$ un foncteur de $\fct(\A,\E)$. Les assertions suivantes sont équivalentes :
\begin{enumerate}
\item le foncteur $F\circ\oplus$ de $\fct(\A^2,\E)$ possède une décomposition en poids faible (relativement à l'action de $A_\mu^2$) ;
\item\label{itt2} pour tout $n\in\mathbb{N}$, le foncteur $F\circ\oplus_n$ de $\fct(\A^n,\E)$ possède une décomposition en poids faible (relativement à l'action de $A_\mu^n$) ;
\item pour tous objets $a_1,\dots,a_n$ de $\A$, l'ensemble $\{F(x_1\oplus\dots\oplus x_n)\,|\,(x_1,\dots,x_n)\in A^n\}$ d'endomorphismes de $F(a_1\oplus\dots\oplus a_n)$ est collectivement trigonalisable ;
\item\label{itt3} pour tout $n\in\mathbb{N}$, le foncteur $cr_n(F)$ de $\fct(\A^n,\E)$ possède une décomposition en poids faible (relativement à l'action de $A_\mu^n$) ;
\item pour tous objets $a_1,\dots,a_n$ de $\A$, l'ensemble $\{cr_n(F)(x_1,\dots,x_n)\,|\,(x_1,\dots,x_n)\in (A\setminus\{0\})^n\}$ d'endomorphismes de $cr_n(F)(a_1,\dots,a_n)$ est collectivement trigonalisable ;
\item pour tout objet $a$ de $\A$, le foncteur $\tau_a(F)$ de $\fct(\A,\E)$ possède une décomposition en poids faible.
\end{enumerate}
Si ces conditions sont vérifiées, on dira que $F$ est un foncteur \index{termin}{hyper-trigonalisant}\textbf{\htr} (relativement à $A$).

La sous-catégorie pleine des foncteurs \hts\ de $\fct(\A,\E)$ (resp. $\F(A,K)$) sera notée \index{nota}{T@$\Tt$, $\Tt_A$ \emph{(catégorie des foncteurs \hts)}} $\Tt_A(\A,\E)$ (resp. $\Tt(A,K)$).
\end{prdef}

\begin{prdef}\label{defo} Soit $F$ un foncteur de $\fct(\A,\E)$. Les assertions suivantes sont équivalentes :
\begin{enumerate}
\item le foncteur $F\circ\oplus$ de $\fct(\A^2,\E)$ possède une décomposition en poids forte (relativement à l'action de $A_\mu^2$) ;
\item pour tout $n\in\mathbb{N}$, le foncteur $F\circ\oplus_n$ de $\fct(\A^n,\E)$ possède une décomposition en poids forte (relativement à l'action de $A_\mu^n$) ;
\item pour tous objets $a_1,\dots,a_n$ de $\A$, l'ensemble $\{F(x_1\oplus\dots\oplus x_n)\,|\,(x_1,\dots,x_n)\in A^n\}$ d'endomorphismes de $F(a_1\oplus\dots\oplus a_n)$ est collectivement diagonalisable ;
\item pour tout $n\in\mathbb{N}$, le foncteur $cr_n(F)$ de $\fct(\A^n,\E)$ possède une décomposition en poids forte (relativement à l'action de $A_\mu^n$) ;
\item pour tous objets $a_1,\dots,a_n$ de $\A$, l'ensemble $\{cr_n(F)(x_1,\dots,x_n)\,|\,(x_1,\dots,x_n)\in (A\setminus\{0\})^n\}$ d'endomorphismes de $cr_n(F)(a_1,\dots,a_n)$ est collectivement diagonalisable ;
\item pour tout objet $a$ de $\A$, le foncteur $\tau_a(F)$ de $\fct(\A,\E)$ possède une décomposition en poids forte.
\end{enumerate}
Si ces conditions sont vérifiées, on dira que $F$ est un foncteur \index{termin}{hyper-diagonalisant}\textbf{\hd} (relativement à $A$).

La sous-catégorie pleine des foncteurs \hds\ de $\fct(\A,\E)$ (resp. $\F(A,K)$) sera notée $\Oo_A(\A,\E)$ (resp. $\Oo(A,K)$).\index{nota}{D@$\Oo$, $\Oo_A$ \emph{(catégorie des foncteurs \hds)}}
\end{prdef}

Les caractérisations en termes de familles collectivement diagonalisables ou trigonalisables d'endomorphismes prennent une forme particulièrement frappante, qui justifie la terminologie employée, dans le cas des catégories $\F(k,K)$, où $k$ est un corps :

\begin{coro}\label{cor-justif_termi} Si $k$ est un corps, un foncteur de $\F(k,K)$ est \hd\ (resp. \htr) si et seulement si pour tout $n\in\mathbb{N}^*$, l'image par $F$ de l'ensemble des matrices diagonales de $\GL_n(k)$ est un ensemble collectivement diagonalisable (resp. collectivement trigonalisable) d'endomorphismes de $F(k^n)$.
\end{coro}

Nous introduisons également la définition suivante, en général strictement moins forte que les notions de foncteur \hd\ ou \htr\ (cf. propositions \ref{diago-fpbar} et \ref{pr-hdcorlf}) ; elle jouera un rôle beaucoup moins important car elle donne lieu à bien moins de bonnes propriétés que les précédentes.

\begin{defi}\label{df-diatri} On dit qu'un foncteur $F$ de $\fct(\A,\E)$ est \textbf{diagonalisant} (resp. \textbf{trigonalisant}) si, pour tous objets $a_1,\dots,a_n$ de $\A$ et tous éléments $x_1,\dots,x_n$ de $A$, l'endomorphisme $F(x_1\oplus\dots\oplus x_n)$ de $F(a_1\oplus\dots\oplus a_n)$ est diagonalisable (resp. trigonalisable).
\end{defi}

\begin{rema}
La définition de foncteur diagonalisant apparaît dès les années $1960$ dans le travail d'Epstein \cite[1.4]{Ep69} (amélioré dans Epstein-Kneser \cite{EK69}), sous le nom de foncteur \emph{diagonal}. Toutefois, Epstein ne considère que des foncteurs entre espaces vectoriels de dimensions finies (à la source comme au but, mais sur des corps éventuellement différents), et n'obtient des résultats significatifs qu'en caractéristique nulle.

Epstein \cite[3.3-3.4]{Ep69} utilise une décomposition en poids. Les méthodes d'Epstein-Kneser \cite{EK69} sont pour leur part plus proches de celles employées ultérieurement dans \cite[§\,4.4]{DTV}.
\end{rema}

\begin{rema} Les définitions et propriétés précédentes relatives aux décompositions en poids des effets croisés, ainsi que la plupart des propriétés formelles des quelques pages qui suivent, valent plus généralement pour une catégorie source non nécessairement additive, mais monoïdale symétrique, dont l'unité (monoïdale) est objet nul, et munie d'une action d'un monoïde avec un élément absorbant (agissant par le morphisme nul sur la catégorie source) compatible avec sa structure monoïdale.
\end{rema}

Le critère le plus important pour garantir qu'un foncteur est \htr\ est le suivant, qui découle du corollaire~\ref{cor-tfdp2} et de la proposition~\ref{pr-vdfdppol}.

\begin{prop}\label{pr-df_htr} 
\begin{enumerate}
    \item Si $K$ est un corps de décomposition non additif de $\A$ (par exemple, si $K$ est algébriquement clos), alors tout foncteur de $\F^\df(\A;K)$ est \htr.
    \item Si $K$ est un corps de décomposition de $\A$, alors tout foncteur analytique ou \ph\ de $\F^\df(\A;K)$ est \htr.
\end{enumerate}
\end{prop}

Lorsque le corps $K$ n'est pas assez gros, nous utiliserons parfois la variante suivante :

\begin{prop}\label{pr-tfdf_htr-ext} Soit $F$ un foncteur de type fini de $\F^\df(\A;K)$. Il existe une extension finie de corps $K\subset L$ telle que le foncteur $F\otimes_K L$ de $\F(\A;L)$ soit \htr. 
\end{prop}

\begin{proof}
Le foncteur $F\circ\oplus$ de $\F(\A\times\A;K)$ est également de type fini et à valeurs de dimensions finies. Le résultat découle donc du corollaire~\ref{cor-excb} et de la proposition~\ref{pr-tfdfend}.
\end{proof}

La propriété suivante découle quant à elle du corollaire~\ref{cor-poidsf-carp}.

\begin{prop}\label{prto} Supposons $K$ de caractéristique $p>0$ et la fonction $A\to A\quad x\mapsto x^p$ surjective. Alors $\Tt_A(\A,\E)=\Oo_A(\A,\E)$.
\end{prop}

En combinant les propositions~\ref{pr-df_htr} et~\ref{prto}, on obtient :
\begin{coro}\label{cor-df_HD} Supposons que $K$ est un corps de décomposition non additif de $\A$ (par exemple, que $K$ est algébriquement clos) de caractéristique $p>0$ et que la fonction $A\to A\quad x\mapsto x^p$ est surjective. Alors tout foncteur de $\F^\df(\A;K)$ est \hd.
\end{coro}

Nous travaillerons le plus souvent sur des foncteurs \hts\ ou \hds\ vérifiant la propriété de finitude introduite dans l'énoncé suivant.

\begin{prdef}\label{pr-htdcf} Soit $F$ un foncteur de $\fct(\A,\E)$. Les assertions suivantes sont équivalentes :
\begin{enumerate}
\item le foncteur $F\circ\oplus$ de $\fct(\A^2,\E)$ possède une décomposition en poids forte (resp. faible) {\em finie} relativement à l'action de $A_\mu^2$ ;
\item pour tout $n\in\mathbb{N}$, le foncteur $F\circ\oplus_n$ de $\fct(\A^n,\E)$ possède une décomposition en poids forte (resp. faible) {\em finie}  relativement à l'action de $A_\mu^n$ ;
\item pour tout $n\in\mathbb{N}$, le foncteur $cr_n(F)$ de $\fct(\A^n,\E)$ possède une décomposition en poids forte (resp. faible) {\em finie}  relativement à l'action de $A_\mu^n$ ;
\item\label{tcf-it4} l'algèbre caractéristique $\ac_A(F)$ est semi-simple (resp. semi-primaire) déployée ;
\item\label{itsansf} le spectre $\Sp(F)$ est un foncteur \hd\ (resp. \htr) de $\F(A,K)$.
\end{enumerate}
Lorsqu'elles sont vérifiées, on dira que $F$ est \index{termin}{caractere fini@ caractère fini \emph{(foncteur \hd\ ou \htr\ à caractère fini)}} {\em \hd\ à caractère fini} (resp. {\em \htr\ à caractère fini}).

On note $\Oo_A^\cf(\A,\E)$\index{nota}{Dcf@$\Oo^\cf$, $\Oo^\cf_A$ \emph{(catégorie des foncteurs \hts\ à caractère fini)}} (resp. \index{nota}{Tcf@$\Tt^\cf$, $\Tt^\cf_A$ \emph{(catégorie des foncteurs \hts\ à caractère fini)}} $\Tt^\cf_A(\A,\E)$) la sous-catégorie pleine de $\fct(\A,\E)$ des foncteurs \hds\ (resp. \hts) à caractère fini. On pose $\Oo^\cf(A,K)=\Oo_A^\cf(\mathbf{P}(A),K\Md)$ et $\Tt^\cf(A,K)=\Tt_A^\cf(\mathbf{P}(A),K\Md)$. 

Tout foncteur \hd\ (resp. \htr) de type fini de $\fct(\A,\E)$ est \hd\ à caractère fini (resp. \htr) à caractère fini.

Tout foncteur \hd\ (resp. \htr) de $\fct(\A,\E)$ est colimite filtrante de ses sous-foncteurs \hds\ (resp. \hts) à caractère fini.
\end{prdef}

\begin{proof} L'équivalence entre les quatre premières conditions se déduit de la proposition~\ref{pr-poideploye}, de la discussion qui précède la proposition~\ref{pr-valsp}, de la décomposition de $F\circ\oplus_n$ à l'aide des effets croisés et du fait que la classe des $K$-algèbres semi-simples (resp. semi-primaires) déployées est stable par sous-quotient et produit tensoriel.

La précomposition par $\oplus$ préserve les foncteurs de type fini, grâce à l'adjonction somme/diagonale. Si un foncteur de type fini $F$ est \htr, la décomposition en poids de $F\circ\oplus$ est donc nécessairement finie : $F$ est à caractère fini. Comme $\Sp(F)$ est toujours de type fini, et a même algèbre caractéristique que $F$ (corollaire~\ref{cor-acsp}), cela implique l'équivalence de \ref{itsansf} avec les conditions précédentes.

Supposons maintenant que $F$ est un foncteur de $\Tt_A(\A,\E)$ (resp. $\Oo_A(\A,\E)$). Le bifoncteur $F\circ\oplus$ est colimite filtrante de sous-foncteurs $X_i$ possédant une décomposition en poids faible (resp. forte) {\bf finie} relativement à l'action de $A^2_\mu$. Si l'on note $G_i$ l'image du morphisme $X_i\circ\delta\to F$ de $\fct(\A,\E)$ déduit de l'inclusion $X_i\hookrightarrow F\circ\oplus$ par adjonction somme/diagonale, on vérifie aisément que $F$ est la colimite filtrante des sous-foncteurs $X_i$ et que ceux-ci appartiennent à $\Tt^\cf_A(\A,\E)$ (resp. $\Oo^\cf_A(\A,\E)$).
\end{proof}

Le résultat suivant donne une liste de propriétés élémentaires d'usage courant des classes de foncteurs précédemment introduites. Il se déduit aussitôt des généralités sur les décompositions en poids exposées au chapitre~\ref{spcf} (par exemple, de la proposition~\ref{pr-poibif} pour le point~5).

\begin{prop}\label{pr-elt-ord}
\begin{enumerate}
\item[1.] Les sous-catégories $\Oo_A(\A,\E)$ et $\Tt_A(\A,\E)$ de $\fct(\A,\E)$ sont prélocalisantes.\index{termin}{prelocalisante@prélocalisante \emph{(sous-catégorie)}}
\item[2.] La sous-catégorie $\Oo^\cf_A(\A,\E)$ de $\fct(\A,\E)$ est semi-épaisse ; la sous-catégorie $\Tt^\cf_A(\A,\E)$ est épaisse.
\item[3.] Les sous-catégories $\Oo_A(\A;K)$, $\Oo_A^\cf(\A;K)$, $\Tt_A(\A;K)$ et $\Tt^\cf_A(\A;K)$ de $\F(\A;K)$ sont stables par produit tensoriel.
\end{enumerate}
 Soient $B$ un anneau et $\B$ une catégorie additive $B$-linéaire essentiellement petite.
\begin{enumerate}
\item[4.] Si $\varphi : B\to A$ est un morphisme d'anneaux et $\Phi : \B\to\A$ un foncteur additif $\varphi$-linéaire, le foncteur de précomposition $\Phi^* : \fct(\A,\E)\to\fct(\B,\E)$ envoie $\Oo_A(\A,\E)$ (resp. $\Oo^\cf_A(\A,\E)$, $\Tt_A(\A,\E)$, $\Tt^\cf_A(\A,\E)$) dans $\Oo_B(\B,\E)$ (resp. $\Oo^\cf_B(\B,\E)$, $\Tt_B(\B,\E)$, $\Tt^\cf_B(\B,\E)$).
\item[5.] Un bifoncteur $F$ de $\fct(\A\times\B,\E)$ est \hd\ (resp. \htr) si et seulement s'il est par rapport à chaque variable. Plus précisément, $F$ appartient à $\Oo_{A\times B}(\A\times\B,\E)$ (resp. $\Tt_{A\times B}(\A\times\B,\E)$) si et seulement si, pour tout $a\in\mathrm{Ob}\,\A$ et tout $b\in\mathrm{Ob}\,\B$, le foncteur $F(a,-)$ appartient à $\Oo_B(\B,\E)$ (resp. $\Tt_B(\B,\E)$) et le foncteur $F(-,b)$ appartient à $\Oo_A(\A,\E)$ (resp. $\Tt_A(\A,\E)$).
\end{enumerate}
\end{prop}

\begin{rema} Si $F$ est un foncteur de $\Oo(K,K)$ et $X$ un foncteur de $\Oo_A(\A;K)$, alors $F\circ X$ appartient à $\Oo_A(\A;K)$. (Lorsque $F$ prend des valeurs de dimension infinie, il faut remplacer $F$ par son extension de Kan à gauche $\tilde{F}$ à $K\Md$ pour définir la composition $F\circ X$.)
\end{rema}

\begin{rema} Des propriétés analogues aux précédentes valent pour les classes de foncteurs diagonalisants et trigonalisants.
\end{rema}

\begin{rema} La sous-catégorie $\Tt_A(\A,\E)$ de $\fct(\A,\E)$ n'est pas nécessairement épaisse (cf. exemple~\ref{exextsdi}).
\end{rema}

\section{Poids partiels et multipoids}\label{sppm}

\begin{defi}Soit $F$ un foncteur de $\Tt_A(\A,\E)$.
Pour $d\in\mathbb{N}$, on note \index{nota}{W@$\widetilde{\W}_d$, $\W_d$ \emph{(ensemble des $d$-multipoids d'un foncteur)}} $$\W_d(F):=\Pi_{A^d_\mu}(cr_d(F))\subset\mon(A_\mu^d,K^\mu)\simeq\mon(A_\mu,K_\mu)^d$$ et $\widetilde{\W}_d(F):=\Pi_{A^d_\mu}(F\circ\oplus_d)$. Les éléments de $\W_d(F)$ sont appelés \index{termin}{multipoids} $d$-\textbf{multipoids} (relativement à $A$) de $F$.

On note\index{nota}{P@$\pp$ \emph{(ensemble des poids partiels d'un foncteur)}}
$$\pp(F):=\W_1(F)\cup\{w\in\mon(A_\mu,K_\mu)\,|\,\exists w'\in\mon(A_\mu,K_\mu)\quad (w,w')\in\W_2(F)\}$$
(ou quelquefois $\pp_A(F)$ si une ambiguïté sur $A$ est possible) et
$$\widetilde{\pp}(F):=\{w\in\mon(A_\mu,K_\mu)\,|\,\exists w'\in\mon(A_\mu,K_\mu)\quad (w,w')\in\widetilde{\W}_2(F)\}.$$
Les éléments de $\pp(F)$ sont appelés \index{termin}{poids!partiel}\textbf{poids partiels} de $F$ (relativement à $A$).
\end{defi}

Les deux propositions suivantes seront très souvent utilisées dans la suite de ce mémoire, parfois sans même être citées.

\begin{prop}\label{pr-poidspar} Soient $F$ un foncteur de $\Tt_A(\A,\E)$, $d, n, m\in\mathbb{N}$, et $f : \llbracket 1,n\rrbracket\to\llbracket 1,m\rrbracket$ une fonction surjective.
\begin{enumerate}
\item Les sous-ensembles $\W_d(F)$ et $\widetilde{\W}_d(F)$ de $\mon(A_\mu,K_\mu)^d$ sont stables par l'action de $\mathfrak{S}_d$ par permutation des facteurs du produit cartésien.
\item\label{pttild} On a $\W_d(F)=\widetilde{\W}_d(F)\cap (\mon(A_\mu,K_\mu)\setminus\{1_{A\to K}\})^d$, et
$\widetilde{\W}_d(F)$ est le saturé par l'action de $\Si_d$ de $\bigsqcup_{i\le d}\W_i(F)\times\{1_{A\to K}\}^{d-i}$ dans $\mon(A_\mu,K_\mu)^d$.
\item Pour tout $\mathbf{w}=(w_1,\dots,w_n)\in\widetilde{\W}_n(F)$, le $m$-uplet $f_*(\mathbf{w})=(v_j)_{1\le j\le m}$ défini par $v_j=\underset{f(i)=j}{\prod}w_i$ appartient à $\widetilde{\W}_m(F)$.
\item Pour tout multipoids $\mathbf{w}\in\W_n(F)$, on a $f_*(\mathbf{w})\in\W_m(F)$. En particulier, le produit des coordonnées d'un multipoids non vide de $F$ appartient à $\W_1(F)$.
\end{enumerate}
\end{prop}

\begin{proof} La première propriété découle de l'invariance (à isomorphisme près) de $F\circ\oplus_d$ et $cr_d(F)$ par l'action de $\Si_d$ sur $\A^d$. 

Comme le multifoncteur $cr_d(F)$ est réduit par rapport à chaque variable, on a $\W_d(F)\subset (\mon(A_\mu,K_\mu)\setminus\{1_{A\to K}\})^d$. En utilisant la décomposition de $F\circ\oplus_d$ en effets croisés \eqref{eq-decCr} (page~\pageref{eq-decCr}), on en déduit l'assertion~\ref{pttild}.

Pour établir l'assertion suivante, on considère les foncteurs
$$\Phi : \A^m\to\A^n\quad (a_j)_{1\le j\le m}\mapsto(a_{f(i)})_{1\le i\le n}\qquad\text{et}$$
$$\Psi : \A^m\to\A^m\quad (a_j)_{1\le j\le m}\mapsto (a_j^{\oplus\cd f^{-1}(\{j\})})_{1\le j\le m}\;.$$
Alors $\Phi$ est $f^*$-linéaire, où $f^* : A^m\to A^n$ désigne le morphisme d'anneaux $(a_j)_{1\le j\le m}\mapsto(a_{f(i)})_{1\le i\le n}$, et essentiellement surjectif à facteur direct près parce que $f$ est surjectif, de sorte que $\Pi(F\circ\oplus_n\circ\Phi)=\widetilde{\W}_n(F)\circ f^*=f_*(\widetilde{\W}_n(F))$, tandis que $\Psi$ est $A$-linéaire, et $\Psi$  de sorte que $\Pi(F\circ\oplus_m\circ\Psi)\subset\widetilde{\W}_m(F)$, grâce à la proposition~\ref{pr-cmsource}. Comme le diagramme
$$\xymatrix{\A^m\ar[r]^-\Phi\ar[d]_-\Psi & \A^n\ar[d]^-{\oplus_n} \\
\A^m\ar[r]^-{\oplus_m} & \A
}$$
commute à isomorphisme près, on en déduit la troisième assertion.

Enfin, comme un poids est distinct de $1_{A\to K}$ si et seulement s'il envoie $0$ sur $0$, tout produit fini non vide de poids distincts de $1_{A\to K}$ est distinct de $1_{A\to K}$. Par conséquent, la dernière assertion se déduit des deux précédentes.
\end{proof}

L'énoncé suivant donne les caractérisations générales simples des poids partiels d'un foncteur \htr.

\begin{prop}\label{pr-poidspar2} Soient $F$ un foncteur de $\Tt_A(\A,\E)$ et $d\in\mathbb{N}$.
\begin{enumerate}
\item\label{itppst0} Les coordonnées d'un multipoids de $F$ sont des poids partiels de $F$. Plus généralement, si $(w_1,\dots,w_d)$ est un multipoids de $F$, pour toute partie non vide $I$ de $\llbracket 1,d\rrbracket$, $w_I:=\prod_{i\in I}w_i$ est un poids partiel de $F$.
\item\label{itppst1} On a $\pp(F)=\widetilde{\pp}(F)\cap (\mon(A_\mu,K_\mu)\setminus\{1_{A\to K}\})$. Si $F$ est non nul, on a $\widetilde{\pp}(F)=\pp(F)\cup\{1_{A\to K}\}$.
\item\label{itppst2} On a $\widetilde{\pp}(F)=\bigcup_{a\in {\rm Ob}\,\A}\Pi(\tau_a(F))\subset\mon(A_\mu,K_\mu)$.
\end{enumerate}
\end{prop}

\begin{proof}
Soient $(w_1,\dots,w_d)$ un multipoids de $F$, $I$ une partie non vide de $\llbracket 1,d\rrbracket$ et $J$ la partie complémentaire. La dernière assertion de la proposition~\ref{pr-poidspar} montre que $w_I$ appartient à $\W_1(F)\subset\pp(F)$ si $J$ est vide, et que $(w_I,w_J)\in\W_2(F)$, donc $w_I\in\pp(F)$ si $J$ est non vide. Cela établit la première assertion.

La deuxième propriété se déduit directement de la proposition~\ref{pr-poidspar}, tandis que l'assertion~\ref{itppst2} est un cas particulier de la proposition~\ref{pr-poibif}.
\end{proof}

La proposition suivante établit un lien simple mais important entre poids partiels et algèbre caractéristique. 

\begin{prop}\label{pr-poidspar_ac} Soient $F$ un foncteur de $\Tt^\cf_A(\A,\E)$ et $d\in\mathbb{N}$.
\begin{enumerate}
\item\label{itppac1} Un poids $w\in\mon(A_\mu,K_\mu)$ appartient à $\widetilde{\pp}(F)$ si et seulement s'il existe un morphisme de $K$-algèbres $\varphi : \ac_A(F)\to K$ tel que $w=\varphi\circ\bar{\chi}_F$.
\item\label{itppac2} On a $\pp(F)=\pp(\Sp(F))$ et $\W_d(F)=\W_d(\Sp(F))$.
\item\label{itppac3} Le morphisme canonique de $K$-algèbres $K[A_\mu]\to K^{\widetilde{\pp}(F)}$ induit un isomorphisme $\ac_A(F)/\rj(\ac_A(F))\xrightarrow{\simeq}K^{\widetilde{\pp}(F)}$.
\item\label{itppac4} Si de plus $F$ est \hd, alors ce morphisme induit un isomorphisme $\ac_A(F)\xrightarrow{\simeq}K^{\widetilde{\pp}(F)}$.
\end{enumerate}
\end{prop}

\begin{proof} Si l'on munit la catégorie $\A\times\A$ de l'action du monoïde $A_\mu$ \emph{sur le deuxième facteur}, on a $\ac_A(F)=\mathfrak{Z}_{A_\mu}(F\circ\oplus)$ et $\widetilde{\pp}(F)=\Pi_{A_\mu}(F\circ\oplus)$. Par conséquent, les assertions \ref{itppac1}, \ref{itppac3} et \ref{itppac4} se déduisent de la proposition~\ref{pr-poideploye}. L'assertion~\ref{itppac2} se déduit quant à elle des propositions \ref{pr-spcr} et \ref{pr-sppoi}.
\end{proof}

La conséquence suivante de la proposition~\ref{pr-poidspar_ac} interviendra au chapitre~\ref{shts}.

\begin{coro}\label{cor-kerpp} Soient $F$ un foncteur de $\Tt^\cf_A(\A,\E)$ et $I$ un idéal de $A$. On suppose que $I\subset\la(\pi)$ pour tout $\pi\in\pp(F)$. Alors il existe $n\in\mathbb{N}$ tel que $I^n\subset\rf(F)$.
\end{coro}

\begin{proof} Soit $n$ l'indice de nilpotence du radical de l'algèbre caractéristique $\ac_A(F)$. La proposition~\ref{pr-poidspar_ac}.\ref{itppac3} montre que l'image de $[x]-[0]$ dans $\ac_A(F)$ appartient à $\rj(\ac_A(F))$ pour tout $x\in I$. Ainsi, pour $x_1,\dots,x_n\in I$, l'image de $\prod_{i=1}^n ([x_i]-[0])=[x_1\dots x_n]-[0]$ dans $\ac_A(F)$ est nulle. La proposition~\ref{pr-acrad} permet d'en déduire $x_1\dots x_n\in\rf(F)$, d'où le résultat.
\end{proof}

Comme tout poids non trivial de $\mon(A_\mu,K_\mu)$ est nul sur les éléments nilpotents de $A$, on en déduit en particulier :

\begin{coro}\label{cor-ht_nil} Si $F$ est un foncteur de $\Tt^\cf_A(\A,\E)$, alors l'idéal $\nil(\AF)$ de $\AF$ est nilpotent.
\end{coro}

\begin{rema}\label{rq-nilcar} Une autre conséquence du corollaire~\ref{cor-kerpp} et de la nullité d'un poids non trivial de $\mon(A_\mu,K_\mu)$ sur $\nil(A)$ est la suivante.

Supposons que $I$ est un nilidéal de $A$ égal à son carré (par exemple, $A=k[\mathbb{Z}[1/p]/\mathbb{Z}]$ où $k$ est un anneau de caractéristique $p$, et $I$ est l'idéal d'augmentation). Alors $I\subset\rf(F)$ pour foncteur $F$ de $\Tt_A(\A,\E)$ (le cas des foncteurs de  $\Tt^\cf_A(\A,\E)$ s'étend à tous les foncteurs \hts\ grâce à la dernière assertion de la proposition~\ref{pr-htdcf}). On en déduit aussi que $I\subset\rf(F)$ pour foncteur $F$ de $\F^\df(\A;K)$ (le cas où $K$ est algébriquement clos se déduit de ce qui précède et de la proposition~\ref{pr-df_htr} ; le cas général s'y ramène par extension des scalaires au but, car cette opération ne change pas le radical d'un foncteur).
\end{rema}

La conséquence directe suivante de l'assertion~\ref{itppac1} de la proposition~\ref{pr-poidspar2} va nous servir à établir le corollaire~\ref{cor-carsideg}.

\begin{coro}\label{cor-chi_pp} Soit $F$ un foncteur de $\Tt_A^\cf(\A,\E)$ et $d\in\mathbb{N}$. Si la fonction caractéristique $\chi_F$ est polynomiale de degré au plus $d$, alors tous les poids partiels de $F$ sont des fonctions polynomiales de degré au plus $d$.
\end{coro}

La proposition suivante nous sera utile dans la section~\ref{ssct-dcps}.

\begin{prop}\label{pr-ppdep} Tout poids partiel d'un foncteur \ph\ de $\Tt_A(\A,\E)$ est polynomial déployé.
\end{prop}

\begin{proof} Comme la classe des foncteurs \phs\ et celle des foncteurs \hts\ sont stables par sous-quotient et par translation, et que tout foncteur simple \ph\ est polynomial, la proposition~\ref{pr-poidspar2}.\ref{itppst2} montre qu'il suffit de vérifier que l'existence d'un foncteur polynomial $F$ de degré $d$, \htr\ et homogène de poids faible $w$ entraîne que le poids polynomial $w$ est déployé (remplacer le foncteur de départ par un facteur de composition du plus grand quotient homogène de poids fort $w$ d'un de ses décalages). Le multifoncteur $cr_d(F)$ est non nul, additif par rapport à chaque variable, et possède une décomposition en poids faible. Soit $(\alpha_1,\dots,\alpha_d)$ un poids de $cr_d(F)$ : les poids $\alpha_i : A\to K$ sont des morphismes d'anneaux (par multiadditivité de $cr_d(F)$), et l'on a  $w=\alpha_1\dots\alpha_d$ par la proposition~\ref{pr-poidspar}, comme souhaité.
\end{proof}

La propriété suivante, qui complète la quatrième assertion de la proposition~\ref{pr-elt-ord}, découle de la proposition~\ref{pr-cmsource}.

\begin{lemm}\label{lm-asqf1} Soient $B$ un anneau, $\B$ une catégorie additive $B$-linéaire essentiellement petite, $\varphi : B\to A$ un morphisme d'anneaux et $\Phi : \B\to\A$ un foncteur $\varphi$-linéaire. Si $F$ est un foncteur de $\Tt_A(\A,\E)$, alors $\pp_B(\Phi^*(F))\subset\pp_A(F)\circ\varphi$, avec égalité si $\Phi$ est essentiellement surjectif à facteur direct près.
\end{lemm}

\begin{lemm}\label{lm-asqf2} Soient $n\in\mathbb{N}$ et $F_1,\dots,F_n$ des foncteurs de $\Tt_A(\A,\E)$. Alors $\widetilde{\pp}(F_1\otimes\dots\otimes F_n)=\widetilde{\pp}(F_1)\dots\widetilde{\pp}(F_n)$ et
$$\pp(F_1\otimes\dots\otimes F_n)=\bigcup_{\varnothing\ne I\subset\llbracket 1,n\rrbracket}\prod_{i\in I}\pp(F_i)\subset\mon(A_\mu,K_\mu)\,.$$
\end{lemm}

\begin{proof} La proposition~\ref{pr-poidspar2}.\ref{itppst2} et
corollaire~\ref{cor-poiptint} montrent la première égalité. La deuxième s'en déduit en utilisant la proposition~\ref{pr-poidspar2}.\ref{itppst1} et fait qu'un produit non vide d'éléments de $\mon(A_\mu,K_\mu)$ distincts de $1$ est distinct de $1$.
\end{proof}

\begin{lemm}\label{lm-asqf3} Soient $d\in\mathbb{N}^*$ et $F$ un foncteur non nul de $\F(K,K)$ qui est un quotient ou un sous-foncteur de la $d$-ième puissance tensorielle $T^d$. Alors $F$ est \hd\ et
$$\pp(F)=\{\lambda\mapsto\lambda^i\,|\,i\in\llbracket 1,d\rrbracket\}\,.$$
\end{lemm}

\begin{proof} Les décalages de $T^d$ sont des sommes directes des $T^i$ avec $i\le d$, qui apparaissent tous dans $\tau_K(T^d)$. Comme $T^i$ est homogène de poids fort $\lambda\mapsto\lambda^i$, $T^d$, donc $F$, est \hd\ et $\pp(F)\subset\pp(T^d)=\{\lambda\mapsto\lambda^i\,|\,i\in\llbracket 1,d\rrbracket\}$ grâce à la proposition~\ref{pr-poidspar2}.\ref{itppst2}.

Le foncteur $F$ est de degré polynomial exactement $d$ (conséquence directe et classique de l'adjonction somme/diagonale itérées), donc $cr_d(F)$ est un quotient ou un sous-foncteur non nul de $cr_d(T^d)$. Or $cr_d(T^d)(V_1,\dots,V_d)\simeq\bigoplus_{\sigma\in\Si_d}V_{\sigma(1)}\otimes\dots\otimes V_{\sigma(d)}$, d'où $\W_d(T^d)=\{(\mathrm{Id}_K,\dots,\mathrm{Id}_K)\}\subset\mon(K_\mu,K_\mu)^d$. Sous-ensemble non vide de $\W_d(T^d)$, $\W_d(F)$ est donc également égal à $\{(\mathrm{Id}_K,\dots,\mathrm{Id}_K)\}$. La proposition~\ref{pr-poidspar2}.\ref{itppst0} implique donc l'inclusion $\pp(F)\supset\{\lambda\mapsto\lambda^i\,|\,i\in\llbracket 1,d\rrbracket\}$, d'où le résultat.
\end{proof}

Nous sommes désormais prêts pour établir deux énoncés qui sont à la base de la démonstration de l'un des principaux de ce chapitre, le théorème~\ref{th-chipol}. La notation $\sigma_p$ qui y apparaît est celle introduite page~\pageref{psig}.

\begin{prop}\label{pr-degpoipar} Supposons que l'anneau $A$ n'a pas de quotient fini de même caractéristique que $K$, et que ce corps est algébriquement clos. Soit $S$ un foncteur simple de $\F^\df(A,K)$, polynomial de degré $d>0$. Alors le degré maximal d'un poids partiel de $S$ est $d$ si la caractéristique $p$ de $K$ est nulle, et vaut au moins $\sigma_p(d)$ si $p>0$.
\end{prop}

(On sait déjà, par la proposition~\ref{pr-df_htr}, que $S$ est \htr.)

\begin{proof} Comme $K$ est algébriquement clos, \cite[th.~5.5 et prop.~6.6]{DTV} montre que $S$ possède une décomposition de la forme $S\simeq\bigotimes_{i=1}^n\alpha_i^*E_i$ où les $\alpha_i : A\to K$ sont des morphismes d'anneaux deux à deux distincts et les $E_i$ sont des foncteurs \emph{élémentaires} non constants de $\F(K,K)$, c'est-à-dire des quotients simples d'une puissance tensorielle $T^{d_i}$ avec $d_i>0$. Ici, la notation $\alpha^*$ désigne la précomposition par l'extension des scalaires $\mathbf{P}(A)\to\mathbf{P}(K)$ le long de $\alpha$, foncteur qui est essentiellement surjectif.

On déduit alors des lemmes \ref{lm-asqf1}, \ref{lm-asqf2} et \ref{lm-asqf3} que
$$\pp(S)=\Big\{\prod_{i=1}^n\alpha_i^{r_i}\,|\,(r_1,\dots,r_n)\in\prod_{i=1}^n\llbracket 0,d_i\rrbracket\setminus\{(0,\dots,0)\}\Big\}\,.$$

Si $p=0$, $\prod_{i=1}^n\alpha_i^{r_i}$ est de degré polynomial $\sum_{i=1}^n r_i$, tandis que $S$ est de degré $d=\sum_{i=1}^n d_i$, ce qui permet de conclure. On suppose désormais $p>0$ et l'on regroupe les $\alpha_i$ qui coïncident à une puissance du morphisme de Frobenius près et l'on utilise que $A$ est sans quotient fini de caractéristique $p$ (les $\alpha_i$ sont donc d'images infinies) pour obtenir une famille \emph{indépendante}\,\footnote{Rappelons que cette notion d'indépendance est introduite dans la définition~\ref{def-fima}.} de morphismes d'anneaux $(\beta_1,\dots,\beta_m)$ telle que $\{\alpha_i\,|\,i\in\llbracket 1,n\rrbracket\}=\{\beta_j^{p^{c(t,j)}}\,|\,j\in\llbracket 1,m\rrbracket, t\in\llbracket 1,u(j)\rrbracket\}$ où, pour chaque $j$, $(c(t,j))_{1\le t\le u(j)}$ est une famille finie non vide d'entiers naturels deux à deux distincts. On a alors
$$\pp(S)=\Big\{\prod_{j=1}^m\beta_j^{\sum_{t=1}^{u(j)}r(t,j)p^{c(t,j)}}\,|\,r(t,j)\in\llbracket 0,d(t,j)\rrbracket\Big\}\setminus\{1_{A\to K}\}$$
pour des entiers $d(t,j)>0$ de somme $d$.

Par la proposition~\ref{pr-degdp}, $\prod_{j=1}^m\beta_j^{\sum_{t=1}^{u(j)}r(t,j)p^{c(t,j)}}$ est de degré polynomial
$$\sum_{j=1}^m s_p\Big(\sum_{t=1}^{u(j)}r(t,j)p^{c(t,j)}\Big).$$
Ainsi, le degré maximal d'un poids partiel de $S$ est
$$\max\Big\{\sum_{j=1}^m s_p\Big(\sum_{t=1}^{u(j)}r(t,j)p^{c(t,j)}\Big)\,|\,0\le r(t,j)\le d(t,j)\Big\}=$$
$$\sum_{j=1}^m\max\Big\{s_p\Big(\sum_{t=1}^{u(j)}r(t,j)p^{c(t,j)}\Big)\,|\,0\le r(t,j)\le d(t,j)\Big\}$$
qui vaut au moins, par le lemme~\ref{lmari-sig},
$$\sum_{j=1}^m\sigma_p\Big(\sum_{t=1}^{u(j)}d(t,j)\Big)\ge\sigma_p\Big(\sum_{j=1}^m\sum_{t=1}^{u(j)}d(t,j)\Big)=\sigma_p(d)\,.$$
\end{proof}

\begin{coro}\label{cor-carsideg} Supposons que l'anneau $A$ n'a pas de quotient fini de même caractéristique que $K$. Soit $S$ un foncteur simple de $\F(A,K)$, polynomial de degré $d>0$. Alors le degré polynomial de la fonction $\chi_S$ est au moins $d$ si $p=0$, et au moins $\sigma_p(d)$ si $p>0$.
\end{coro}

\begin{proof} Quitte à agrandir le corps $K$, on peut supposer que $S$ est à valeurs de dimensions finies, par la proposition~\ref{pr-simpoldfext}. Le foncteur $S\otimes\bar{K}$ (où $\bar{K}$ désigne une clôture algébrique de $K$) de $\F^\df(A,\bar{K})$ possède un quotient simple polynomial $T$ de degré $d$. Par la proposition~\ref{pr-degpoipar}, $T$ possède un poids partiel de degré au moins $d$ si $p=0$, et au moins $\sigma_p(d)$ si $p>0$. Il s'ensuit, par le corollaire~\ref{cor-chi_pp}, que le degré polynomial de la fonction caractéristique $\chi_T$ vaut au moins $d$ ou $\sigma_p(d)$ selon que $p$ est nul ou non. Grâce au lemme~\ref{lm-chisq} ci-dessous, la même propriété vaut pour $\chi_{S\otimes\bar{K}}$ puisque $T$ est un quotient de $S\otimes\bar{K}$, et donc pour $\chi_S$, comme souhaité.
\end{proof}

\begin{lemm}\label{lm-chisq} Pour tout $n\in\mathbb{N}$, la sous-catégorie pleine des foncteurs $F$ de $\fct(\A,\E)$ dont la fonction caractéristique est polynomiale de degré au plus $n$ est prébilocalisante.
\end{lemm}

\begin{proof} La fonction $\chi_F$ est polynomiale de degré au plus $n$ si et seulement si l'idéal caractéristique $\ic_A(F)$ contient $\bar{K}[A]^{\star (n+1)}$ (puissance $(n+1)$-ième de l'idéal d'augmentation de $K[A]$). Le résultat se déduit donc de la proposition~\ref{pr-prebilocIdCar}.
\end{proof}

\begin{rema} En utilisant les résultats de \cite{DG}, on peut aussi décrire les poids partiels des foncteurs simples \textit{antipolynomiaux} $S$ de $\F(A,K)$, qui sont \hts\ (grâce à la proposition~\ref{pr-df_htr}) si $K$ est un corps de décomposition non additif de $\A$ --- ce que nous supposerons dans la suite de cette remarque --- car à valeurs de dimensions finies. Par \cite[th.~10.2]{DG}, on a $S\simeq\mathrm{Q}(T,M)$ où $T$ est un $A$-module fini et $M$ un $K[\mathrm{Aut}(T)]$-module simple, uniquement déterminés à isomorphisme près (pour les notations $\mathrm{Q}^{T,M}$ et $\mathrm{Q}(T,M)$, on renvoie à \cite[§\,5]{DG}). Comme le radical de $\mathrm{Q}(T,M)$ contient l'annulateur $I$ de $T$ dans $A$, quitte à remplacer $A$ par $A/I$, on peut supposer que $T$ est un $A$-module fidèle.

 Il existe une suite exacte $0\to\mathrm{Q}(T,M)\to\mathrm{Q}^{T,M}\to F\to 0$ avec $F$ dans $\F_{d-1}(A,K)$ (sous-catégorie de $\F(A,K)$ introduite dans \cite[déf.~8.1]{DG}), grâce aux résultats de \cite[§\,6]{DG}, et un objet $x$ de $\mathbf{P}(A)$ tel que $K[T\underset{A}{\otimes}-]$ soit facteur direct de $\tau_x(\mathrm{Q}^{T,M})$ par \cite[prop.~5.7]{DG} (et les résultats de dualité \cite[cor.~3.6 et prop.~4.3]{DG}). Dans la suite exacte $0\to\tau_x(\mathrm{Q}(T,M))\to\tau_x(\mathrm{Q}^{T,M})\to\tau_x(F)\to 0$, $\tau_x(F)$ appartient à $\F_{d-1}(A,K)$ par \cite[prop.~8.3]{DG}. Comme tout morphisme de $\mathrm{Q}^T\subset K[T\underset{A}{\otimes}-]$ vers un foncteur de $\F_{d-1}(A,K)$ est nul par \cite[th.~9.9]{DG}, il s'ensuit que $\tau_x(\mathrm{Q}^{T,M})$ contient $\mathrm{Q}^T$. En utilisant de nouveau \cite[prop.~5.7]{DG}, on déduit de la proposition~\ref{pr-poidspar2}.\ref{itppst2} que tout élément de $\mon(A_\mu,K_\mu)$ est poids partiel de $S$.
\end{rema}

\section{Application aux foncteurs \phs}\label{safphs}

Cette section applique les résultats de la fin de la précédente ainsi que les notions de spectre et de fonction caractéristique aux foncteurs \phs\ pour établir un énoncé important, le corollaire~\ref{cor-EMLpol-pol}. Celui-ci garantit que certains foncteurs \phs\ sont polynomiaux en caractéristique $p>0$ --- à comparer aux énoncés faciles de \cite[prop.~2.8 et 2.11]{DTV}.

Commençons, pour fixer les idées, par un énoncé simple sur les foncteurs polynomiaux, laissé en exercice.

\begin{prop} Si $F$ est un foncteur polynomial de degré $d$ de $\fct(\A,\E)$, la fonction caractéristique $\chi_F : A\to K$ est polynomiale de degré au plus $d$.
\end{prop}

\begin{prop}\label{pr-chi_ph} Soit $F$ un foncteur \ph\ de $\fct(\A,\E)$. Si $F$ est à support fini, ou à co-support fini, alors $\chi_F$ est une fonction polynomiale.
\end{prop}

\begin{proof} Supposons que $E$ est un ensemble fini d'objets de $\A$ constituant un support, ou un co-support, de $F$. Alors le morphisme de $K$-algèbres (non commutatives) canonique
$$\mathrm{End}(F\circ\oplus)\to\prod_{x\in E}\mathrm{End}_K(F(x\oplus x))$$
est injectif.

On en déduit par composition avec $\chi_F$ une fonction multiplicative
\begin{equation}\label{eqchipol}
A\to\mathrm{End}(F\circ\oplus)\to\prod_{x\in E}\mathrm{End}_K(F(x\oplus x))
\end{equation}
dont chaque composante $A\to\mathrm{End}_K(F(x\oplus x))$ égale
$$A\xrightarrow{a\mapsto 1\oplus a}\mathrm{End}_\A(x\oplus x)\to\mathrm{End}_K(F(x\oplus x))\,,$$
qui est polynomiale (la première flèche est affine, et la seconde est polynomiale car $F$ est \ph). Comme $E$ est fini, \eqref{eqchipol} est polynomiale, ce qui entraîne que $\chi_F$ l'est également, puisque \eqref{eqchipol} est la composée de $\chi_F$ suivie d'une fonction additive injective.
\end{proof}

\begin{rema} On ne peut s'affranchir de toute hypothèse de finitude sur $F$ dans l'énoncé précédent. Ainsi, dans $\F(K,K)$, le foncteur $\bigoplus_{d\in\mathbb{N}}\Lambda^d$ est \ph, mais on vérifie facilement que sa fonction caractéristique n'est pas polynomiale lorsque le corps $K$ est infini (cela se déduit aussi directement du théorème~\ref{th-chipol} ci-après).
\end{rema}

\begin{lemm}\label{lm-phEK}
Soit $X$ un foncteur \ph\ de $\F(A,K)$. L'extension de Kan à gauche $\tilde{X}$ de $X$ le long de l'inclusion de $\mathbf{P}(A)$ dans la catégorie des $A$-modules de type fini est \ph e.
\end{lemm}

\begin{proof} Soient $M$ et $N$ des $A$-modules de type fini. Choissisons des surjections $A$-linéaires $u : U\twoheadrightarrow M$ et $v : V\twoheadrightarrow N$ avec $U$ et $V$ dans $\mathbf{P}(A)$ et considérons le diagramme commutatif
$$\xymatrix{\mathrm{Hom}_A(M,N)\ar[r]^-{u^*}\ar[d]^{\mathrm{(1)}} & \mathrm{Hom}_A(U,N)\ar[d]^{\mathrm{(2)}} & \mathrm{Hom}_A(U,V)\ar[l]_-{v_*}\ar[d]^{\mathrm{(3)}} \\
\mathrm{Hom}_K(\tilde{F}(M),\tilde{F}(N))\ar[r]^-{\tilde{F}(u)^*} & \mathrm{Hom}_K(F(U),\tilde{F}(N)) & \mathrm{Hom}_K(F(U),F(V))\ar[l]_-{\tilde{F}(v)^*}
}$$
dont les flèches verticales sont induites par $\tilde{F}$ (qui coïncide avec $F$ sur $\mathbf{P}(A)$). La fonction (3) est polynomiale car $F$ est \ph ; comme $\tilde{F}(v)^*$ est additive, et $v_*$ additive surjective (car $v$ est surjective et $U$ projectif), il s'ensuit que la fonction (2) est également polynomiale. Comme $u^*$ est additive, et $\tilde{F}(u)^*$ additive injective ($\tilde{F}(u)$ est surjective par le lemme~\ref{lm-vdfEK}), il s'ensuit que la fonction (1) est polynomiale, donc que le foncteur $\tilde{F}$ est \ph. 
\end{proof}

\begin{prop}\label{pr-hisfhp}\begin{enumerate}
\item Soient $F$ et $G$ des foncteurs de $\fct(\A,\E)$. On suppose que $F$ est à support fini et $G$ \ph. Alors le foncteur $\hi(F,G)$ de $\F(A,K)$ est \ph.
\item Soient $X$ un foncteur à support fini de $\F(A,K)$ et $F$ un foncteur \ph\ de $\fct(\A,\E)$. Alors $\hi(X,F)$ et $X\ti F$ sont des foncteurs \phs\ de $\fct(\A,\E)$.
\item Supposons que la condition \textnormal{(MTF)}\index{nota}{MTF@(MTF) \emph{(condition de finitude sur les morphismes)}} est satisfaite. Soient $X$ un foncteur \ph\ de $\F(A,K)$ et $F$ un foncteur à support fini de $\fct(\A,\E)$. Alors $X\ti F$ est un foncteur \ph\ de $\fct(\A,\E)$.
\end{enumerate}
\end{prop}

\begin{proof} Montrons la première assertion (très analogue à la proposition~\ref{pr-vdffac}). Si $F$ est à support fini, c'est un quotient d'une somme directe finie de foncteurs de la forme $M[\A(a,-)](:=M^{\oplus\A(a,-)}])$, pour des objets $M$ de $\E$ et $a$ de $\A$. Or $\hi(M[\A(a,-)],G)\simeq\E(M,-)\circ G\circ(-\underset{A}{\otimes}a)$. Comme les foncteurs $\E(M,-)$ et $-\underset{A}{\otimes}a$ sont additifs, il s'ensuit que, si $G$ est \ph, alors $\hi(M[\A(a,-)],G)$ et $\hi(F,G)$ sont \phs.

La deuxième assertion s'établit de façon entièrement similaire. La dernière se montre également de la même manière, en utilisant la proposition~\ref{pr-EKti} et le lemme~\ref{lm-phEK}.
\end{proof}

\begin{coro}\label{cor-phsp} Soit $F$ un foncteur à support fini de $\fct(\A,\E)$. 
\begin{enumerate}
\item Si $F$ est \ph, alors le foncteur $\Sp(F)$ de $\F(A,K)$ est \ph.
\item Supposons que la condition \textnormal{(MTF)}\index{nota}{MTF@(MTF) \emph{(condition de finitude sur les morphismes)}} est satisfaite. Si $\Sp(F)$ est \ph, alors $F$ est \ph.
\end{enumerate}
\end{coro}

\begin{proof} Si $F$ est à support fini et \ph, la première assertion de la proposition~\ref{pr-hisfhp} montre que $\hi(F,F)$, et donc son sous-foncteur $\Sp(F)$, sont \phs.

Sous l'hypothèse \textnormal{(MTF)}, si $F$ est à support fini et $\Sp(F)$ \ph, la troisième assertion de la même proposition montre que $\Sp(F)\ti F$, qui est isomorphe à $F$ par la proposition~\ref{pr-spdf}, est \ph.
\end{proof}

\begin{prop}\label{pr-chiph} Soit $F$ un foncteur de $\fct(\A,\E)$. Sa fonction caractéristique $\chi_F$ est polynomiale si et seulement si son spectre $\Sp(F)$ est un foncteur \ph\ de $\F(A,K)$.
\end{prop}

\begin{proof} De manière générale, comme le foncteur $P^A$ de $\F(A,K)$ est engendré par $[\mathrm{Id}_A]\in P^A(A)$, un quotient $X$ de $P^A$ est \ph\ si et seulement si, pour tout objet $V$ de $\mathbf{P}(A)$, la fonction
$$V\xrightarrow{v\mapsto [v]}K[V]\simeq P^A(V)\twoheadrightarrow X(V)$$
est polynomiale (cette observation simple se démontre exactement comme \cite[lemme~2.12]{DTV}). Il s'ensuit que $\Sp(F)$ est \ph\ si et seulement si, pour tout $n\in\mathbb{N}$, la fonction
$$A^n\xrightarrow{v\mapsto [v]}K[V]\simeq P^A(A^n)\xrightarrow{\eta_F(A^n)}\hi(F,F)(A^n)\simeq\mathrm{End}_{\fct(\A^n,\E)}(F\circ\oplus_n)$$
(où $\eta_F$ est la transformation naturelle de la définition~\ref{defisp}) est polynomiale. Or cette fonction coïncide avec $(a_1,\dots,a_n)\mapsto\Xi^{(n)}_F([a_1,\dots,a_n])$ (cf. \eqref{eqXi}, page \pageref{eqXi}), avec les notations du §\,\ref{par-spac}, de sorte que la formule \eqref{eq-chin} (page \pageref{eq-chin}) et le fait qu'un produit de fonctions polynomiales est polynomial permettent de conclure.
\end{proof}

Le dernier ingrédient nécessaire à la démonstration du résultat principal de cette section est le lemme général suivant sur les foncteurs polynomiaux (qui est indépendant des autres considérations de ce chapitre).

\begin{lemm}\label{lm-degfc} Soient $i\in\mathbb{N}$ et $F$ un foncteur de $\F(\A;K)$ dont tous les sous-quotients simples appartiennent à $\pol_i(\A;K)$. Alors $F$ appartient à $\pol_i(\A;K)$.
\end{lemm}

\begin{proof} Comme $\pol_i(\A;K)$ est stable par colimites, on peut supposer que $F$ est de type fini. Le plus grand quotient $\qpol_i(F)$ de $F$ appartenant à $\pol_i(\A;K)$ est le conoyau d'un morphisme naturel $F\circ\oplus_{i+1}\circ\delta_{i+1}\to F$ (cf. \eqref{eq-qpol}, page~\pageref{eq-qpol}). Or $F\circ\oplus_{i+1}\circ\delta_{i+1}$ est, comme $F$, de type fini, donc le noyau $G$ de la projection $F\twoheadrightarrow\qpol_i(F)$ également. Supposons $G$ non nul : ce foncteur possède alors un quotient simple $S$, qui appartient par hypothèse à $\pol_i(\A;K)$. En formant la somme amalgamée $T$ des morphismes $G\twoheadrightarrow S$ et $G\hookrightarrow\qpol_i(F)$, on obtient une suite exacte courte $0\to S\to T\to\qpol_i(F)\to 0$, où $T$ est un quotient de $F$. Comme $\pol_i(\A;K)$ est une sous-catégorie épaisse de $\F(\A;K)$, $T$ lui appartient. Cela contredit le fait que $\qpol_i(F)$ est le plus grand quotient de $F$ dans $\pol_i(\A;K)$, de sorte que l'hypothèse que $G$ est non nul est absurde. Ainsi $F\simeq\qpol_i(F)$ appartient bien à $\pol_i(\A;K)$.
\end{proof}

Si $p$ est premier, $\sigma_p$ et $\sigma_p^!$ désignent les fonctions $\mathbb{N}\to\mathbb{N}$ introduites page~\pageref{psig}, et pour $p=0$ on pose $\sigma_0=\sigma_0^!=\mathrm{Id}_\mathbb{N}$. En particulier, on a $\sigma_p(n)\le\sigma_p(n+1)\le\sigma_p(n)+1$ et $\sigma_p^!(d)=\max\{n\in\mathbb{N}\,|\,\sigma_p(n)\le d\}$ pour tous $n, d\in\mathbb{N}$.

\begin{theo}\label{th-chipol} Supposons que l'anneau $A$ n'a pas de quotient fini de caractéristique $p$. Si $F$ est un foncteur de $\F(A,K)$ dont la fonction caractéristique est polynomiale de degré $d$, alors $F$ est polynomial de degré au plus $\sigma_p^!(d)$.
\end{theo}

\begin{proof} D'après la proposition~\ref{pr-chiph} et le corollaire~\ref{cor-phsp}, $F$ est \ph. Il s'ensuit que tous les facteurs de composition de $F$ sont polynomiaux, car un foncteur \ph\ simple est polynomial, et un sous-quotient d'un foncteur \ph\ est \ph. Le corollaire~\ref{cor-carsideg} et le lemme~\ref{lm-chisq} montrent que le degré des facteurs de composition de $F$ ne peut pas excéder $\sigma_p^!(d)$. Le lemme~\ref{lm-degfc} fournit alors la conclusion.
\end{proof}

\begin{coro}\label{cor-EMLpol-pol} Supposons que l'anneau $A$ n'a pas de quotient fini de caractéristique $p$. Alors tout foncteur \ph\ à support fini $F$ de $\fct(\A,\E)$ est polynomial.
\end{coro}

\begin{proof} Le foncteur $\Sp(F)$ est \ph\ (corollaire~\ref{cor-phsp}) et à support fini, donc sa fonction caractéristique est polynomiale (proposition~\ref{pr-chi_ph}). On conclut alors par le théorème précédent et le corollaire~\ref{cor-sppol}.
\end{proof}

\begin{rema} Il résultera du corollaire~\ref{cor-imq} ci-après que, si $A$ est noethérien (ou vérifie plus généralement la condition (IMQ) introduite page~\pageref{pageIMQ}) et que $F$ est un foncteur \ph\ de type fini de $\F(\A;K)$ \textit{prenant des valeurs de dimensions finies}, alors il existe un idéal cofini $I$ de $A$ tel que $F$ se factorise comme composée du foncteur canonique $\A\to\A/I\times\A$ et d'un foncteur de $\F(\A/I\times\A;K)$ qui est polynomial par rapport à la deuxième variable. Nous ignorons si cette généralisation partielle du corollaire~\ref{cor-EMLpol-pol} subsiste lorsqu'on supprime l'hypothèse que $F$ prend des valeurs de dimensions finies.
\end{rema}

\chapter{Foncteurs \hds}\label{shd}

\begin{cvi}
Dans ce chapitre, $\A$ désigne toujours une catégorie additive $A$-linéaire essentiellement petite et $\E$ une catégorie de Grothendieck $K$-linéaire.
\end{cvi}

Nous donnerons, au théorème~\ref{th-struct_ord}, une caractérisation complète des foncteurs \hds\ à caractère fini à l'aide de quotients finis de $A$ et de (multi)foncteurs \emph{strictement} polynomiaux sur $K$. Ceux-ci, introduits dans les années $1990$ par Friedlander et Suslin \cite{FS} (nous en rappellerons la définition dans la section~\ref{ssct-dcps}), sont beaucoup mieux compris que les foncteurs polynomiaux et possèdent une structure beaucoup plus rigide (qui ne dépend presque que de la caractéristique de $K$). Nous donnerons des applications de notre théorème principal dans la section~\ref{sapkb}, notamment à la structure de $\F^\df(\bar{\FF}_p,\bar{\FF}_p)$, frappantes dans la mesure où celle de $\F(\bar{\FF}_p,\bar{\FF}_p)$ semble hors d'atteinte, à partir de l'observation que tout foncteur à valeurs de dimensions finies de cette catégorie est \hd.

\section{Propriétés élémentaires}\label{ssct-hdelem}

Nous donnons dans cette section des conditions suffisantes, et quelquefois nécessaires, pour qu'un foncteur soit \hd\,; nous traiterons aussi parfois de foncteurs diagonalisants (à titre essentiellement culturel, afin d'illustrer les similitudes et différences entre foncteurs \hds\ et diagonalisants).

Notre premier énoncé repose sur les résultats du chapitre~\ref{sec-finpolp}.

\begin{prop}\label{pr-hdpolpoi}
 Supposons que $K$ est un corps de décomposition de $\mathbf{P}(A)$ et que l'une des conditions suivantes est réalisée :
\begin{enumerate}
\item $p=0$ et la condition \textnormal{(CAD)}\index{nota}{CAD@(CAD), (CAD)$^+$ \emph{(conditions d'annulation des dérivations)}} est vérifiée ;
\item $p>0$ et la fonction $A\to A\quad x\mapsto x^p$ est surjective.
\end{enumerate}

Alors tout foncteur \ph\ de $\F(A,K)$ possédant une décomposition en poids faible est \hd.
\end{prop}

\begin{proof} On note tout d'abord que, du fait que $K$ est un corps de décomposition de $\mathbf{P}(A)$, la condition (CAD) implique (CAD)$^+$.

Comme $\Oo(A,K)$ est stable par colimites, il suffit de montrer que tout foncteur \ph\ de type fini $F$ de $\F(A,K)$ est \hd. Le corollaire~\ref{cor-phpoinoeth} et la proposition~\ref{pr-wss} (pour $p>0$) ou le corollaire~\ref{cor-wssder} (pour $p=0$) montrent que $F$ est à valeurs de dimensions finies. Comme $F$ est \ph\ et que $K$ est un corps de décomposition de $\mathbf{P}(A)$, la proposition~\ref{pr-df_htr} permet d'en déduire que $F$ est \htr. On conclut que $F$ est \hd\ en utilisant la proposition~\ref{pr-tppfor} dans le cas $p=0$ ou la proposition~\ref{prto} dans le cas $p>0$.
\end{proof}

\begin{prop}\label{pr-red} Si $F$ est un foncteur diagonalisant de $\fct(\A,\E)$, alors l'algèbre caractéristique $\ac_A(F)$ et l'anneau $\AF$ sont réduits.
\end{prop}

\begin{proof} Soit $\xi=\sum_{i=1}^n\lambda_i[a_i]$ un élément de $K[A]$, où les $a_i$ (resp. $\lambda_i$) sont des éléments de $A$ (resp. $K$). Puisque $F$ est diagonalisant, les endomorphismes $F(1\oplus a_i)$ de $F\circ\oplus$ sont diagonalisables. Comme ils commutent deux à deux, la sous-algèbre de $\mathrm{End}(F\circ\oplus)$ qu'ils engendrent est constituée d'éléments diagonalisables (cela se montre comme dans le cas classique des espaces vectoriels de dimension finie ; cf. par exemple \cite[thm~4.14]{IMR16}), elle est donc réduite. Il s'ensuit que l'image de $\xi$ dans $\ac(F)$ ne peut être nilpotente que si elle est nulle, ainsi $\ac(F)$  est réduit.

La fonction caractéristique $\chi_F$ induit un morphisme de monoïdes multiplicatifs {\em injectif} (proposition~\ref{pr-acrad}) $\varphi : \AF\to\ac(F)$. Si $a\in \AF$ est tel que $a^n=0$, alors $(\varphi(a)-\varphi(0))^n=\varphi(a^n)-\varphi(0)=0$, d'où $\varphi(a)-\varphi(0)=0$ et $a=0$, comme souhaité.
\end{proof}

Nous aurons besoin des deux énoncés généraux suivants pour établir la proposition~\ref{diago-corfini}.

\begin{lemm}\label{lm-matdiagELTR} Soient $n\in\mathbb{N}$ et $x\in A$. Supposons que $x$ et $x-1$ sont inversibles.
Alors toute représentation de $\GL_n(A)$ sur laquelle les matrices diagonales dont les valeurs propres sont des puissances de $x$ opèrent trivialement a une restriction au sous-groupe engendré par les matrices élémentaires triviale.
\end{lemm}

\begin{proof} Cela découle de l'identité suivante dans $\GL_2(A)$ :
$$\left(\begin{array}{cc} 1 & a\\
                                                0 & 1
                                               \end{array}\right)=\left(\begin{array}{cc} 1 & 0\\
                                                0 & x^{-1}
                                               \end{array}\right)\left(\begin{array}{cc} 1 & a\\
                                                0 & x-1
                                               \end{array}\right)\left(\begin{array}{cc} 1 & 0\\
                                                0 & x
                                               \end{array}\right)\left(\begin{array}{cc} 1 & a\\
                                                0 & x-1
                                               \end{array}\right)^{-1}\,.$$
\end{proof}

\begin{lemm}\label{lm-actELT} Soit $F$ un foncteur de $\F(A,K)$ tel que, pour tout $n\in\mathbb{N}$, le sous-groupe $\operatorname{E}_n(A)$ de $\GL_n(A)$ engendré par les matrices élémentaires opère trivialement sur $F(A^n)$. Alors $F$ est constant.
\end{lemm}

\begin{proof}
Soient $U=A^i$ et $V=A^j$ des objets et $f : U\to V$ un morphisme de $\mathbf{P}(A)$. Alors l'automorphisme $\varphi$ de $U\oplus V$ donné matriciellement par $\left(\begin{array}{cc} 1_U & 0\\
                                               f & 1_V
                                               \end{array}\right)$ appartient à $\operatorname{E}_{i+j}(A)$, ainsi $F(\varphi)$ est l'identité. Comme $f$ égale la composée
$$U\hookrightarrow U\oplus V\xrightarrow{\varphi}U\oplus V\twoheadrightarrow V\,,$$
il s'ensuit que $F(f)=F(0_{U\to V})$, ce qui entraîne que $F$ est constant.
\end{proof}

\begin{prop}\label{diago-corfini} Supposons que $A$ est un corps fini de cardinal $q$.
\begin{enumerate}
\item Si $q-1$ est inversible dans $K$ et appartient à $\nu(K)$ (par exemple, si $K$ contient $k$), alors $\Oo^\cf(A,K)=\F(A,K)$.
\item Inversement, si $q-1$ est nul dans $K$, ou n'appartient pas à $\nu(K)$, alors tout foncteur diagonalisant de $\F(A,K)$ est constant.
\end{enumerate}
\end{prop}

\begin{proof} Le premier point résulte de la proposition~\ref{pr-dpfotjs}.

Sous les hypothèses du deuxième point, il existe un élément $x$ de $A^\times$ différent de $1$ et appartenant au noyau de tout morphisme de groupes $A^\times\to K^\times$. Si $F$ est diagonalisant, cela entraîne que l'image par $F$ de tout automorphisme diagonal de $A^n$, pour tout $n\in\mathbb{N}$, dont les valeurs propres sont des puissances de $x$ est l'identité. Les lemmes \ref{lm-matdiagELTR} et \ref{lm-actELT} permettent alors de conclure.
\end{proof}

\begin{prop}\label{pr-antipol_hd} Soit $F$ un foncteur de $\F(A,K)$ tel que $\AF$ soit fini. Alors $F$ est \hd\ si et seulement si $\AF$ est isomorphe à un produit fini de corps finis $\FF_q$ tels que $q-1$ appartienne à $\nu(K)$ et soit inversible dans $K$.
\end{prop}

\begin{proof}
Quitte à remplacer $A$ par $\AF$, on peut supposer que l'anneau $A$ est fini. En le décomposant en produit d'anneaux locaux et en utilisant la proposition~\ref{pr-elt-ord}.5, on voit qu'on peut supposer également que $A$ est local. Comme un anneau local fini réduit est un corps, la proposition~\ref{pr-red} montre que, si $F$ est \hd, alors $\AF$ est un corps fini $\FF_q$. La proposition~\ref{diago-corfini} permet alors de conclure.
\end{proof}

\begin{rema} La proposition~\ref{pr-antipol_hd} s'applique en particulier aux foncteurs antipolynomiaux de $\F(A,K)$. Elle permet de voir qu'il existe des foncteurs simples antipolynomiaux (et donc en particulier à valeurs de dimensions finies) qui ne sont pas \hds, même si $K$ est algébriquement clos. Cela contraste avec la situation des foncteurs simples polynomiaux de $\F^\df(A,K)$ (cf. proposition~\ref{pr-simplesPolHD} ci-après).
\end{rema}

\begin{prop} Soit $F$ un foncteur diagonalisant de $\fct(\A,\E)$.
Supposons de plus que l'une des conditions suivantes est réalisée :
\begin{enumerate}
\item $F$ est \htr\ ;
\item le bifoncteur $F\circ\oplus$ vérifie la condition {\rm (Dec)} ;
\item le bifoncteur $F\circ\oplus$ est noethérien ;
\item $\E=K\Md$ et $F$ appartient à $\F^\df(\A;K)$.
\end{enumerate}
Alors le foncteur $F$ est \hd.
\end{prop}

\begin{proof}  Le premier cas est clair (cf. remarque~\ref{rq-diatricol}). Le deuxième découle de la proposition~\ref{pr-dtf}. Le troisième s'y ramène par le corollaire~\ref{cor-ntcf}. Le dernier s'y ramène également lorsque $F$ est de type fini par le corollaire~\ref{endofini-dec}, et le cas général s'en déduit par stabilité de $\Oo_A(\A,K\Md)$ par colimites.
\end{proof}

L'énoncé simple suivant s'applique notamment à $\F(\bar{\FF}_p,\bar{\FF}_p)$, qui contient par contraste de nombreux foncteurs qui \emph{ne} sont \emph{pas} \hds\ (cf. la proposition~\ref{pr-dpfotjs} et la fin de ce chapitre).

\begin{prop}\label{diago-fpbar} Supposons que $A$ est un sous-corps localement fini de $K$. Alors tout foncteur de $\fct(\A,\E)$ est diagonalisant.
\end{prop}

\begin{proof}
Le corps $A$ est colimite \emph{filtrante} de ses sous-corps finis $k$, or l'algèbre $K[k_\mu]$ est semi-simple déployée (cf. l'implication \ref{dpfotjs5}$\Rightarrow$\ref{dpfotjs4} de la proposition~\ref{pr-dpfotjs}), ce qui permet de conclure.
\end{proof}

\section{Structures polynomiales strictes}\label{ssct-dcps}

Pour $d\in\mathbb{N}$, on désigne par $\Gamma^d_K$, ou simplement $\Gamma^d$, l'endofoncteur des $K$-espaces vectoriels de $d$-ième puissance divisée, autrement dit, $\Gamma^d(V)$ est le sous-espace de $V^{\otimes d}$ des éléments invariants sous l'action du groupe symétrique $\Si_d$. Le foncteur $\Gamma^d$ possède une structure monoïdale symétrique canonique (relativement au produit tensoriel sur $K$). Par conséquent, si $\B$ est une petite catégorie $K$-linéaire, on peut former une catégorie $K$-linéaire $\Gamma^d\B$ ayant les mêmes objets que $\B$ (cf. §\,\ref{pfTg}). Les foncteurs strictement polynomiaux (relativement à $K$) homogènes de degré $d$ de $\B$ vers $\E$ sont par définition les foncteurs $K$-linéaires de $\Gamma^d\B$ dans $\E$, dont la catégorie est notée \index{nota}{P@$\Pp_d$ \emph{(catégorie de foncteurs strictement polynomiaux)}} $\Pp_d(\B,\E)$. On note $\Pp(\B,\E)$ la catégorie $\prod_{d\in\mathbb{N}}\Pp_d(\B,\E)$ ; ses objets sont appelés \textit{foncteurs strictement polynomiaux} (relativement à $K$) de $\B$ dans $\E$. Si $B$ et $C$ sont des $K$-algèbres, on note $\Pp(B,C)$ pour $\Pp(\mathbf{P}(B),C\Md)$.

On note $\Pp^\cf(\B,\E)$ la sous-catégorie pleine des foncteurs de  $\Pp(\B,\E)$ dont seul un nombre fini de composantes homogènes est non nul, ses objets sont appelés foncteurs strictement polynomiaux à caractère fini de $\B$ dans $\E$. (Attention, dans certaines références, ce sont ces foncteurs qui sont nommés foncteurs strictement polynomiaux.)

Le morphisme naturel $K[V]\to\Gamma^d(V)\quad [v]\mapsto v^{\otimes d}$ induit un foncteur $K$-linéaire $K[\B]\to\Gamma^d\B$ égal à l'identité sur les objets ; par précomposition, on en déduit un foncteur canonique $\Pp_d(\B,\E)\to\fct(\B,\E)$ qui est fidèle, continu et cocontinu. Ces foncteurs s'assemblent de façon unique en un foncteur cocontinu $\Pp(\B,\E)\to\fct(\B,\E)$ qui est également exact et fidèle.

L'image de $\Pp_d(\B,\E)$ est incluse dans les sous-catégories $\pol_d(\B,\E)$ et $\fct(\B,\E)_{w_d}$ de $\fct(\B,\E)$, où $w_d : K\to K$ est le poids $x\mapsto x^d$, ce qui explique la terminologie.

\begin{rema}Les foncteurs strictement polynomiaux ont été introduits par Friedlander et Suslin \cite[§\,2]{FS} dans le cas où $\B=\mathbf{P}(K)$ et $\E=K\Md$ ; on pourra consulter aussi \cite[§\,2.1]{TouzeRingel} pour des généralités sur ces foncteurs. Nous suivons ici le point de vue de \cite[§\,2]{DT-TJM}.

Dans \cite{T-torsion}, le second auteur utilise des considérations de poids relativement à l'action de $\mathbb{Z}_\mu$ pour contrôler la torsion (sur $\mathbb{Z}$) des valeurs prises par les foncteurs strictement polynomiaux entre $R$-modules, où $R$ est un \emph{anneau}, avec des applications en topologie algébrique et en algèbre homologique.
\end{rema}

Les notations $\gamma^*$ et $\gamma_*$ qui apparaissent dans ce qui suit ont été introduites page~\pageref{ppgam}.

\begin{defi}\label{def-decStPolStr} On appelle {\em structure polynomiale stricte faible} (relativement à $K$) sur un foncteur $F$ de $\fct(\A,\E)$ tout isomorphisme entre $F$ et la composition de l'image dans $\fct(\gamma_1^*\A\times\dots\times\gamma_n^*\A,\E)$ d'un objet $X$ de $\Pp^\cf(\gamma_1^*\A\times\dots\times\gamma_n^*\A,\E)$ et du foncteur canonique 
$$\A\to\gamma_1^*\A\times\dots\times\gamma_n^*\A$$
dont les composantes sont les $(\gamma_i)_* :\A\to\gamma_i^*\A$,
où $n\in\mathbb{N}$ et $\gamma_1,\dots,\gamma_n : A\to K$ sont des morphismes d'anneaux.

Si $X$ est un produit tensoriel extérieur de foncteurs de $\Pp^\cf(\gamma_1^*\A,\E),\dots,\Pp^\cf(\gamma_n^*\A,\E)$, on parle de {\em structure polynomiale stricte faible tensorielle}.

Si $(\gamma_1,\dots,\gamma_n)$ est une famille \textbf{indépendante}\index{termin}{independante@indépendante \emph{(famille de morphismes d'anneaux)}} de morphismes d'anneaux, on parle de {\em structure polynomiale stricte} (cf. définition~\ref{def-fima}).
\end{defi}

L'intérêt spécifique de cette dernière situation est illustré par les deux énoncés suivants.

Le lemme ci-dessous est la forme duale d'un résultat de Borel-Tits \cite[prop.~2.1]{Borel-Tits} (voir aussi \cite[corollaire~9.13]{DTV}).

\begin{lemm}\label{lm-BT} Si $(\gamma_1,\dots,\gamma_n)$ est une famille indépendante\index{termin}{independante@indépendante \emph{(famille de morphismes d'anneaux)}} de morphismes d'anneaux $A\to K$, alors pour tout $m\in\mathbb{N}$, le morphisme canonique
$$P^A\to\bigoplus_{1\le j_1,\dots,j_n\le m}\Big(\bigotimes_{i=1}^n\gamma_i^*\Gamma^{j_i}\Big)$$
de $\F(A,K)$ est un épimorphisme.
\end{lemm}

(Les composantes de ce morphisme sont obtenues à partir des morphismes  $K[V]\to\gamma^*\Gamma^d(V)=\Gamma^d(V\,_A\otimes_\gamma K)\quad [v]\mapsto (v\otimes 1)^{\otimes d}$ en utilisant le coproduit $n$-itéré sur $P^A$.)

On en déduit aussitôt, en utilisant les propositions \ref{pr-fctX_d2} et \ref{pr-sp_sdir} et le corollaire~\ref{cor-sp_pt} :

\begin{prop}\label{pr-strPBT} Si $(\gamma_1,\dots,\gamma_n)$ est une famille indépendante\index{termin}{independante@indépendante \emph{(famille de morphismes d'anneaux)}} de morphismes d'anneaux $A\to K$, alors le foncteur canonique $\Pp(\gamma_1^*\A\times\dots\times\gamma_n^*\A,\E)\to\fct(\A,\E)$ est pleinement fidèle et son image essentielle est stable par sous-quotient.

De plus, l'image essentielle de sa restriction à $\Pp^\cf(\gamma_1^*\A\times\dots\times\gamma_n^*\A,\E)$ est constituée des foncteurs $F$ de $\fct(\A,\E)$ tels qu'il existe $m\in\mathbb{N}$ tel que
$$\Sp(F)\le\bigoplus_{1\le j_1,\dots,j_n\le m}\Big(\bigotimes_{i=1}^n\gamma_i^*\Gamma^{j_i}\Big)\ .$$
\end{prop}

\begin{rema} Pour $\A=\mathbf{P}(A)$ et $\E=K\Md$, on vérifie sans peine la réciproque : si le foncteur canonique $\Pp(\gamma_1^*\mathbf{P}(A)\times\dots\times\gamma_n^*\mathbf{P}(A),K\Md)\to\F(A,K)$ est pleinement fidèle, alors $(\gamma_1,\dots,\gamma_n)$ est une famille indépendante de morphismes d'anneaux.
\end{rema}

\begin{prop}\label{pr-polstrHD} Tout foncteur de $\fct(\A,\E)$ possédant une structure polynomiale stricte faible est \hd.
\end{prop}

\begin{proof}
L'image essentielle du foncteur canonique $\Pp(K,K)\to\F(K,K)$ est constituée de foncteurs ayant une décomposition en poids forte. Comme elle est stable par translation, il s'ensuit qu'elle est incluse dans $\Oo(K,K)$. La conclusion provient alors de la proposition~\ref{pr-elt-ord}.
\end{proof}

Le résultat suivant, qui repose lourdement sur \cite{DTV}, constitue l'une des clefs de la démonstration du théorème de structure des foncteurs \hds\ (théorème~\ref{th-struct_ord} ci-après).

\begin{prop}\label{pr-simplesPolHD} Soit $S$ un foncteur simple de $\F(A,K)$, polynomial de degré $d$. Les conditions suivantes sont équivalentes :
\begin{enumerate}
\item\label{dsts1} $S$ possède une structure polynomiale stricte faible tensorielle ;
\item\label{dsts1b} $S$ possède une structure polynomiale stricte faible ;
\item\label{dsts2} les facteurs de composition de $cr_d(S)$ dans $\F(A^d,K)$ sont absolument simples ;
\item\label{dsts3} $S$ est \hd\,;
\item\label{dsts4} $S$ est \htr\,;
\item\label{dsts5} le multifoncteur $cr_d(S)$ possède une décomposition en poids faible.
\end{enumerate}

Si l'on suppose de plus $K$ parfait, il suffit que $S$ possède un poids (global) polynomial déployé pour que les conditions précédentes soient satisfaites.
\end{prop}

\begin{proof} L'équivalence entre \ref{dsts1} et \ref{dsts2} est donnée par \cite[th.~5.3]{DTV}. 

La proposition~\ref{pr-polstrHD} fournit l'implication \ref{dsts1b}$\Rightarrow$\ref{dsts3}.

Les implications \ref{dsts1}$\Rightarrow$\ref{dsts1b} et \ref{dsts3}$\Rightarrow$\ref{dsts4}$\Rightarrow$\ref{dsts5} sont triviales.

L'implication \ref{dsts5}$\Rightarrow$\ref{dsts2} provient de ce qu'un foncteur simple de $\F(A^d,K)$ additif par rapport à chacune des $d$ variables, qui s'identifie à un $A_K^{\otimes d}$-module simple, possède un poids si et seulement s'il est absolument simple.

Supposons maintenant que $K$ est parfait et que $S$ possède un poids $\pi$ déployé. Par le théorème~\ref{th-poinoeth}, $S$ appartient à $\F^\df(A,K)$. Il s'ensuit (proposition~\ref{pr-tfdf_htr-ext}) qu'il existe une extension finie de corps $K\subset L$ telle que le foncteur $S_L:=S\otimes L$ de $\F(A,L)$ soit \htr. Les propositions~\ref{pr-Wcrimm} et~\ref{prw-fond}.\ref{pr-wdepl} montrent alors que la condition \ref{dsts5} est vérifiée, d'où la conclusion.
\end{proof}

\begin{rema} On ne peut s'affranchir de l'hypothèse de perfection de $K$ dans la dernière partie de la proposition. Supposons par exemple que $K$ est de degré fini $n>1$ sur l'image de son endomorphisme de Frobenius $\varphi$. Notons $\Phi\in\F(K,K)$ le foncteur de restriction le long de $\varphi$. On vérifie sans peine que $\Phi^{\otimes p}$ est homogène de poids faible $\mathrm{Id}_K$. Ce foncteur est polynomial de degré $p$, à valeurs de dimension finie, et donc fini. Il possède en particulier au moins un facteur de composition $S$ polynomial de degré $p$. Ce foncteur est homogène de poids fort $\mathrm{Id}_K$, qui est polynomial déployé. Mais $cr_d(S)$ ne peut admettre de décomposition en poids, car il n'existe pas d'endomorphismes $\alpha_1,\dots,\alpha_p$ du corps $K$ dont le produit soit l'identité, puisque $K$ est imparfait (utiliser encore l'unicité au Frobenius près de la proposition~\ref{pr-dppol}).
\end{rema}

Le résultat suivant constitue une variante de la proposition~\ref{pr-hdpolpoi} qui se démontre de la même façon, en utilisant la dernière assertion de la proposition~\ref{pr-simplesPolHD} (et le fait que tout foncteur polynomial de $\F^\df(A,K)$ est fini) au lieu de la proposition~\ref{pr-vdfdppol}.

\begin{prop}\label{pr-hdpolpoi2}
 Supposons que le corps $K$ est parfait et que l'une des conditions suivantes est réalisée :
\begin{enumerate}
\item $p=0$ et la condition \textnormal{(CAD)}\index{nota}{CAD@(CAD), (CAD)$^+$ \emph{(conditions d'annulation des dérivations)}} est vérifiée ;
\item $p>0$ et la fonction $A\to A\quad x\mapsto x^p$ est surjective.
\end{enumerate}

Alors tout foncteur polynomial de $\F(A,K)$ homogène d'un certain poids polynomial déployé fort est \hd.
\end{prop}

\section{Théorème de structure}

L'observation simple suivante, qui se rapproche de la proposition~\ref{pr-sp_pte}, est laissée en exercice.

\begin{prdef} Soient $\B$ une petite catégorie $K$-linéaire, $\C$ une petite catégorie et $F$ un foncteur de $\fct(\C\times\B,\E)$. Les assertions suivantes sont équivalentes :
\begin{enumerate}
\item l'image de $F$ par l'équivalence canonique $\fct(\C\times\B,\E)\simeq\fct(\B,\fct(\C,\E))$ appartient à l'image essentielle de $\Pp(\B,\fct(\C,\E))$ (resp. $\Pp^\cf(\B,\fct(\C,\E))$) ;
\item l'image de $F$ par l'équivalence canonique $\fct(\C\times\B,\E)\simeq\fct(\C,\fct(\B,\E))$ est à valeurs dans l'image essentielle du foncteur canonique $\Pp(\B,\E)\to\fct(\B,\E)$ (resp. $\prod_{d\le n}\Pp_d(\B,\E)\to\fct(\B,\E)$ pour un certain $n\in\mathbb{N}$) ;
\item vu comme foncteur $K$-linéaire $K[\C\times\B]\to\E$, $F$ est isomorphe à la somme directe sur $d\in\mathbb{N}$ (resp. sur $d\le n$ pour un certain $n\in\mathbb{N}$) de foncteurs se factorisant à travers le foncteur canonique $K[\C\times\B]\simeq K[\C]\otimes K[\B]\to K[\C]\otimes\Gamma^d\B$.
\end{enumerate}
Lorsqu'elles sont satisfaites, on dira que $F$ est \emph{strictement polynomial} (resp. \emph{strictement polynomial à caractère fini}) \emph{par rapport à $\B$}.
\end{prdef}

La définition suivante renforce la notion de décomposition de type $FP$ introduite à la définition~\ref{def-decFP}.

\begin{defi}\label{def-decStGlobStr} On appelle \emph{décomposition de type $FP$ stricte} d'un foncteur $F$ de $\fct(\A,\E)$ tout isomorphisme entre $F$ et la composition du foncteur canonique
$$\A\to\C\times\B\quad\text{où}\quad\B=\prod_{j=1}^m\gamma_j^*\A\quad\text{et}\quad\C=\prod_{i=1}^n\A/\mathfrak{m}_i$$
et d'un foncteur de $\fct(\C\times\B,\E)$ strictement polynomial à caractère fini par rapport à $\B$, où $n, m$ sont des entiers naturels, les $\mathfrak{m}_i$ ($1\le i\le n$) sont des idéaux maximaux deux à deux distincts de $\A$ tels que chaque corps $A/\mathfrak{m}_i$ soit de cardinal fini $q_i$, où $q_i-1$ appartient à $K^\times$ et à $\nu(K)$, et les $\gamma_j : A\to K$ ($1\le j\le m$) forment une famille \textbf{indépendante}\index{termin}{independante@indépendante \emph{(famille de morphismes d'anneaux)}} de morphismes d'anneaux.
\end{defi}

\begin{rema} Si l'on impose de plus que le multifoncteur qui apparaît dans la définition précédente ne soit constant par rapport à aucune variable, alors les idéaux $\mathfrak{m}_i$ sont entièrement déterminés par $F$ (à l'ordre près), et les morphismes d'anneaux $\gamma_j$ sont entièrement déterminés par $F$ à l'ordre \emph{et à l'action du morphisme de Frobenius}, si $p>0$, près.
\end{rema}

Le morphisme canonique qui apparaît dans l'énoncé suivant s'obtient en utilisant le coproduit du foncteur $P^A_{\mathbf{P}(A)}$, le morphisme $P^A_{\mathbf{P}(A)}\twoheadrightarrow\pi_I^*P^{A/I}_{\mathbf{P}(A/I)}$ de réduction modulo $I$ et le morphisme canonique du lemme~\ref{lm-BT}.

\begin{lemm} \label{lm-modsp} Soient $n, m\in\mathbb{N}$, $\gamma_1,\dots,\gamma_n : A\to K$ une famille indépendante\index{termin}{independante@indépendante \emph{(famille de morphismes d'anneaux)}} de morphismes d'anneaux et $I$ un idéal cofini de $A$. Alors le morphisme canonique
$$P^A_{\mathbf{P}(A)}\to\pi_I^*P^{A/I}_{\mathbf{P}(A/I)}\otimes\left(\bigoplus_{1\le j_1,\dots,j_n\le m}\Big(\bigotimes_{i=1}^n\gamma_i^*\Gamma^{j_i}\Big)\right)=:\mathrm{Q}_{A,I,(\gamma_1,\dots,\gamma_n),m}$$
de $\F(A,K)$ est un épimorphisme, où $\pi_I$ désigne la réduction modulo $I$.
\end{lemm}

\begin{proof} Le résultat se déduit par dualité de \cite[lemme~4.20 et corollaire~9.13]{DTV}.
\end{proof}

On en déduit, par le corollaire~\ref{cor-sp_pt} :

\begin{prop}\label{pr-strglobsp} Dans cet énoncé, les lettres $\mathfrak{m}_1,\dots,\mathfrak{m}_n$, $\gamma_1,\dots,\gamma_m$, $\B$ et $\C$ ont la même signification qu'à la définition~\ref{def-decStGlobStr}.

La restriction aux foncteurs strictement polynomiaux à caractère fini par rapport à $\B$ du foncteur canonique $\fct(\C\times\B,\E)\to\fct(\A,\E)$ est pleinement fidèle et son image essentielle est stable par sous-quotient.

De plus, son image essentielle est constituée des foncteurs $F$ de $\fct(\A,\E)$ tels qu'il existe $r\in\mathbb{N}$ tel que
$$\Sp(F)\le\mathrm{Q}_{A,\mathfrak{m}_1\cap\dots\cap\mathfrak{m}_n,(\gamma_1,\dots,\gamma_m),r}\ .$$
\end{prop}

La définition générale suivante nous servira à montrer, à la proposition~\ref{pr-hypdia-sp}, que le spectre d'un foncteur \hd\ est déterminé par ses ensembles de multipoids.

\begin{defi}\label{def-sature} On dit qu'un foncteur \htr\ de $\F(A,K)$ est {\em saturé} si, pour tout $d\in\mathbb{N}^*$, on a $$\W_d(F)=\{(w_1,\dots,w_d)\in\pp(F)^d\,|\,\forall \varnothing\ne E\subset\llbracket 1,d\rrbracket\qquad\prod_{i\in E}w_i\in\pp(F)\}.$$
\end{defi}

D'après la proposition~\ref{pr-poidspar2}, on a toujours
$$\W_d(F)\subset\{(w_1,\dots,w_d)\in\pp(F)^d\,|\,\forall \varnothing\ne E\subset\llbracket 1,d\rrbracket\qquad\prod_{i\in E}w_i\in\pp(F)\}\,,$$
de sorte que le foncteur \htr\ $F$ est saturé si et seulement si ses ensembles de multipoids sont les plus grands possibles, l'ensemble de poids partiels $\pp(F)$ étant fixé, ce qui justifie la terminologie.  La proposition~\ref{pr-poidspar}.\ref{pttild} montre également que la condition de saturation équivaut à 
$$\forall d\in\mathbb{N}\quad\widetilde{\W}_d(F)=\{(w_1,\dots,w_d)\in\widetilde{\pp}(F)^d\,|\,\forall E\subset\llbracket 1,d\rrbracket\quad\prod_{i\in E}w_i\in\widetilde{\pp}(F)\}.$$

Le résultat suivant est simple mais crucial.

\begin{lemm}\label{lm-sat} Soient $n, m\in\mathbb{N}$, $\gamma_1,\dots,\gamma_n : A\to K$ une famille indépendante\index{termin}{independante@indépendante \emph{(famille de morphismes d'anneaux)}} de morphismes d'anneaux et $I$ un idéal cofini de $A$. On suppose également que le groupe multiplicatif $(A/I)^\times$ est somme directe de groupes cycliques dont les ordres appartiennent à $\nu(K)$. Alors le foncteur $\mathrm{Q}_{A,I,(\gamma_1,\dots,\gamma_n),m}$ du lemme~\ref{lm-modsp} est \htr\ et saturé. Si l'anneau $A/I$ est un produit de corps de cardinal $q$ tel que $q-1$ appartienne à $\nu(K)$ et à $K^\times$, il est même \hd. De plus,
$$\widetilde{\pp}(\mathrm{Q}_{A,I,(\gamma_1,\dots,\gamma_n),m})=\mon((A/I)_\mu,K_\mu).\left\{\prod_{i=1}^n \gamma_i^{e_i}\,|\,0\le e_i\le m\right\}$$
(où l'on désigne encore, par abus, par $\mon((A/I)_\mu,K_\mu)$ le sous-monoïde de $\mon(A_\mu,K_\mu)$ des morphismes se factorisant par la réduction modulo $I$).
\end{lemm}

\begin{proof}
Le foncteur $P^{A/I}_{\mathbf{P}(A/I)}$ est \htr\ grâce à la proposition~\ref{pr-dpftjs}, et \hd\ sous l'hypothèse supplémentaire sur $A/I$ par la proposition~\ref{pr-antipol_hd}. Les propositions~\ref{pr-elt-ord} et~\ref{pr-polstrHD} permettent d'en déduire que $\mathrm{Q}_{A,I,(\gamma_1,\dots,\gamma_n),m}$ est \htr, et \hd\ sous l'hypothèse supplémentaire. La conclusion résulte donc des observations générales suivantes :
\begin{enumerate}
\item $\widetilde{\pp}(F\otimes G)=\widetilde{\pp}(F).\widetilde{\pp}(G)$ pour tous foncteurs $F$ et $G$ de $\F(A,K)$ (lemme~\ref{lm-asqf2}), d'où $\widetilde{\W}_d(F\otimes G)=\widetilde{\W}_d(F).\widetilde{\W}_d(G)$ pour tout $d\in\mathbb{N}$ ;
\item $\widetilde{\pp}(P^k_{\mathbf{P}(k)})=\Pi(P^k_{\mathbf{P}(k)})=\mon(k_\mu,K_\mu)$ pour tout anneau fini $k$, et $P^k_{\mathbf{P}(k)}$ est saturé (à cause de la propriété exponentielle $P^k_{\mathbf{P}(k)}\circ\oplus_d\simeq (P^k_{\mathbf{P}(k)})^{\boxtimes d}$) ;
\item $\widetilde{\pp}(\gamma^*\Gamma^i)=\{\gamma^e\,|\,0\le e\le i\}$ pour tout $e\in\mathbb{N}$ et tout morphisme d'anneaux $\gamma : A\to K$ (grâce aux lemmes \ref{lm-asqf3} et \ref{lm-asqf1}), et $\gamma^*\Gamma^i$ est saturé (à cause de la propriété exponentielle graduée : $\gamma^*\Gamma^i\circ\oplus_d\simeq\underset{j_1+\dots+j_d=i}{\bigoplus}\gamma^*\Gamma^{j_1}\boxtimes\dots\boxtimes\gamma^*\Gamma^{j_d}$) ;
\item le morphisme d'anneaux canonique $\mon((A/I)_\mu,K_\mu)\times\mathbb{N}^n\to\mon(A_\mu,K_\mu)$ correspondant à l'inclusion $\mon((A/I)_\mu,K_\mu)\hookrightarrow\mon(A_\mu,K_\mu)$ et aux $n$ éléments $\gamma_1,\dots,\gamma_n$ est injectif grâce aux propositions~\ref{pr-decpord} et~\ref{pr-dppol}, ce qui permet d'obtenir le caractère saturé de $\mathrm{Q}_{A,I,(\gamma_1,\dots,\gamma_n),m}$.
\end{enumerate}
\end{proof}

\begin{prop}\label{pr-hypdia-sp}  Soient $F$ et $G$ des foncteurs  de $\Oo_A^\cf(\A,\E)$.
\begin{enumerate}
\item\label{it-hypdia-sp} Pour tout $d\in\mathbb{N}$, la dimension du $K$-espace vectoriel $\Sp(F)(A^d)$ égale le cardinal de $\widetilde{\W}_d(F)$. En particulier, $\Sp(F)$ appartient à $\F^\df(A,K)$.
\item Les assertions suivantes sont équivalentes :
\begin{enumerate}
\item\label{it2a} $\Sp(F)=\Sp(G)$ (resp. $\Sp(F)\le\Sp(G)$) dans l'ensemble ordonné des quotients de $P^A$ ;
\item\label{it2b} pour tout $d\in\mathbb{N}$, $\widetilde{\W}_d(F)=\widetilde{\W}_d(G)$ (resp. $\widetilde{\W}_d(F)\subset\widetilde{\W}_d(G)$) ;
\item\label{it2c} pour tout $d\in\mathbb{N}$, $\W_d(F)=\W_d(G)$ (resp. $\W_d(F)\subset\W_d(G)$).
\end{enumerate}
\end{enumerate}
\end{prop}

\begin{proof} Pour $w\in\widetilde{\W}_d(F)$, notons $e_w$ l'élément de $\hi(F,F)(A^d)$ composé du morphisme $F\simeq F\circ (A\otimes_A -)\to F\circ (A^d\otimes_A -)$ induit par la diagonale itérée $A\to A^d$ et de l'idempotent de $F\circ (A^d\otimes_A -)\simeq F\circ\oplus_d\circ\delta_d$ associé au poids $w$. Les $e_w$ forment une famille libre dans le $K$-espace vectoriel $\hi(F,F)(A^d)$ en vertu de l'isomorphisme d'adjonction
\[\fct(\A,\E)(F,F\circ (A^d\otimes_A -))\simeq \fct(\A^d,\E)(F\circ\oplus_d,F\circ\oplus_d)\]
(induit par précomposition par l'adjonction entre somme et diagonale itérées).

Comme $F$ appartient à $\Oo_A^\cf(\A,\E)$, les $e_w$ forment une base de l'espace vectoriel $\Sp(F)(A^d)$. En effet, si $a=(a_1,\dots,a_d)$ est un élément de $A^d$, l'image de $[a_1,\dots,a_d]\in K[A^d]$ par la projection canonique $P^A\twoheadrightarrow\Sp(F)$ (évaluée sur $A^d$) est
\[\sum_{w\in\widetilde{\W}_d(F)}\prod_{i=1}^d w_i(a_i)\cdot e_w.\]
Ainsi, $\Sp(F)(A^d)$ est inclus dans l'espace vectoriel engendré par les $e_w$ ; il lui est exactement égal car les poids de $F\circ\oplus_d$, vus comme éléments de $\mon(A_\mu^d,K_\mu)\subset K^{A^d}$, forment une famille libre de ce $K$-espace vectoriel (indépendance linéaire des morphismes de monoïdes --- cf. par exemple \cite[chap.~V, §\,6.1, corollaire]{Bki2}). On en déduit l'égalité $\dim_K\Sp(F)(A^d)=\cd\widetilde{\W}_d(F)$, d'où la première assertion. Cela montre également que $\Sp(F)$ est entièrement déterminé, comme quotient de $P^A$, par les ensembles $\widetilde{\W}_d(F)$, d'où l'on tire l'équivalence entre \ref{it2a} et \ref{it2b}. L'équivalence de \ref{it2b} et \ref{it2c} découle quant à elle de la proposition~\ref{pr-poidspar}.\ref{pttild}.
\end{proof}

\begin{lemm}\label{lm-diagsp} Soit $B : \mathbf{P}(A)\times\mathbf{P}(A)\to K\Md$ un bifoncteur de type fini, anti-polynomial par rapport à la première variable et \ph\ par rapport à la deuxième. Si la composée de $B$ avec la diagonale $\delta$ de $\mathbf{P}(A)$ appartient à $\Oo(A,K)$, alors $B$ appartient à $\Oo_{A\times A}(\mathbf{P}(A)\times\mathbf{P}(A),K\Md)$.
\end{lemm}

\begin{proof} Comme la précomposition par $\delta$ préserve les foncteurs de type fini et qu'un foncteur de type fini est \hd\ si et seulement si son algèbre caractéristique est semi-simple déployée (proposition~\ref{pr-htdcf}), il suffit de voir que $B$ et $B\circ\delta$ ont même algèbre caractéristique. La proposition~\ref{pr-chi_ph} montre que la fonction caractéristique de $B$ est polynomiale, disons de degré $d$, par rapport à la deuxième variable.

Il s'ensuit que, pour un idéal $K$-cotrivial $I$ approprié de $A$, le morphisme d'algèbres $K[A_\mu\times A_\mu]\to Z(\mathrm{End}(B\circ\oplus))$ qui prolonge $\chi_B$ se factorise par la projection canonique $K[A\times A]\twoheadrightarrow K[A/I]\otimes Q_d(A)$, où $Q_d:=\qpol_d P^A$. Or \cite[lemme~4.20]{DTV} montre (par dualité) que la composée du morphisme $K[A]\to K[A\times A]$ induit par la diagonale avec ce morphisme $K[A\times A]\twoheadrightarrow K[A/I]\otimes Q_d(A)$ est surjective, d'où le lemme.
\end{proof}

Nous pouvons maintenant démontrer le résultat principal de ce chapitre :

\begin{theo}\label{th-struct_ord} Un foncteur de $\fct(\A,\E)$ appartient à $\Oo^\cf_A(\A,\E)$ si et seulement s'il possède une décomposition de type $FP$ stricte.
\end{theo}

\begin{proof} Si $F$ est un foncteur possédant une décomposition de type $FP$ stricte, alors la proposition~\ref{pr-strglobsp} et le lemme~\ref{lm-sat} montrent que $\Sp(F)$ appartient à $\Oo^\cf(A,K)$. La proposition~\ref{pr-htdcf} implique ensuite que $F$ appartient à $\Oo^\cf_A(\A,\E)$.

Réciproquement, soit $F$ un foncteur de $\Oo^\cf(\A,\E)$. Le foncteur $\Sp(F)$ de $\F(A,K)$ est de type fini, et à valeurs de dimensions finies par la proposition~\ref{pr-hypdia-sp}. Par \cite[th.~4.10]{DTV} (cité au théorème~\ref{th-DTVglob} du présent mémoire), il s'ensuit qu'il est isomorphe à la composée de la diagonale $\mathbf{P}(A)\to\mathbf{P}(A)\times\mathbf{P}(A)$ et d'un bifoncteur $B$ (nécessairement de type fini) antipolynomial par rapport à la première variable et \ph\ par rapport à la seconde. Le lemme~\ref{lm-diagsp} montre que $B$ appartient à $\Oo_{A\times A}(\mathbf{P}(A)\times\mathbf{P}(A),K\Md)$.

La proposition~\ref{pr-antipol_hd} montre donc que, par rapport à la première variable, $B$ se factorise par la réduction modulo un idéal $I$ de $A$ tel que $A/I$ soit isomorphe à un produit fini de corps finis $\FF_q$ de caractéristiques différentes de $p$, avec $q-1$ dans $\nu(K)$ et dans $K^\times$.

Chaque poids partiel de $F$ s'exprime comme produit d'un morphisme de $\mon((A/I)_\mu,K_\mu)$ et d'un poids polynomial déployé (par la proposition~\ref{pr-ppdep}). Le lemme~\ref{lm-sat} permet d'en déduire l'inclusion $\pp(F)\subset\pp(\mathrm{Q}_{A,J,(\gamma_1,\dots,\gamma_n),m})$ pour des entiers naturels convenables $n$ et $m$, une famille indépendante convenable $(\gamma_i)$ de $\mathbf{Ann}(A,K)$ et un idéal à quotient fini $J$ intersection de $I$ et d'un idéal $J'$ tel que $A/J'$ soit un produit fini de corps finis de caractéristique $p$ se plongeant dans $K$ ($J'$ s'obtient comme intersection des noyaux des morphismes d'anneaux d'image finie dans une décomposition de la partie polynomiale des poids partiels de $F$ en produit de morphismes d'anneaux $A\to K$). Comme $\mathrm{Q}_{A,J,(\gamma_1,\dots,\gamma_n),m}$ est saturé (lemme~\ref{lm-sat}), cela entraîne $\W_d(F)\subset\W_d(\mathrm{Q}_{A,J,(\gamma_1,\dots,\gamma_n),m})$ pour $d$, d'où $\Sp(F)\le\Sp(\mathrm{Q}_{A,J,(\gamma_1,\dots,\gamma_n),m})$ par la proposition~\ref{pr-hypdia-sp}. La proposition~\ref{pr-strglobsp} permet de conclure.
\end{proof}

\begin{rema}
 Alors que la notion de foncteur \hd\ est liée à l'effet du foncteur sur les endomorphismes diagonalisables, la démonstration précédente fait un usage crucial de la décomposition à la Steinberg globale de \cite{DTV}, qui repose sur l'examen de l'effet d'un foncteur sur les endomorphismes unipotents.
\end{rema}

Les résultats suivants constituent des corollaires directs du théorème~\ref{th-struct_ord} (en utilisant la dernière assertion de la proposition~\ref{pr-htdcf} pour passer de $\Oo^\cf$ à $\Oo$).

\begin{coro}\label{cor-ttHDpol} Supposons qu'il n'existe pas de morphisme d'anneaux surjectif de $A$ sur un corps fini de cardinal $q$ tel que $q-1$ appartienne à $\nu(K)$ et $K^\times$. Alors tout foncteur de $\Oo_A^\cf(\A,\E)$ possède une structure polynomiale stricte.

Réciproquement, si tout foncteur de $\Oo^\cf(A,K)$ est polynomial, alors il n'existe pas de morphisme d'anneaux surjectif de $A$ sur un corps fini de cardinal $q$ tel que $q-1$ appartienne à $\nu(K)$ et $K^\times$.
\end{coro}

Le résultat suivant précise la proposition~\ref{diago-corfini}.

\begin{coro}\label{cor-HDcst} Si les deux conditions suivantes sont réalisées, alors tout foncteur de $\Oo_A(\A,\E)$ est constant.
\begin{enumerate}
\item Il n'existe pas de morphisme d'anneaux surjectif de $A$ sur un corps fini de cardinal $q$ tel que $q-1$ appartienne à $\nu(K)$ et $K^\times$.
\item Il n'existe pas de morphisme d'anneaux de $A$ dans $K$.
\end{enumerate}
Réciproquement, si tout foncteur de $\Oo^\cf(A,K)$ est constant, alors les conditions précédentes sont remplies.
\end{coro}

\section{Applications}\label{sapkb}

\subsection{Propriétés de finitude des foncteurs \hds}

\begin{theo}\label{th-finitude-diago}
\begin{enumerate}
\item Si $K$ est de caractéristique nulle, alors la catégorie $\Oo(A,K)$ est semi-simple.
\item La catégorie $\Oo(A,K)$ est \index{termin}{localement noetherienne@localement noethérienne \emph{(catégorie abélienne)}} localement noethérienne. Elle est localement finie\index{termin}{localement finie@localement finie \emph{(catégorie abélienne)}} si et seulement s'il n'existe pas de morphisme d'anneaux $A\to K$ d'image finie.
\end{enumerate}
\end{theo}

\begin{proof} Pour montrer la première assertion, il suffit de voir que si $K$ est un corps \textit{algébriquement clos} de caractéristique $0$, alors la catégorie $\Oo(A,K)$ est semi-simple. En effet, en général, si $\bar{K}$ est une clôture algébrique de $K$ et $F$ un foncteur de $\Oo(A,K)$, alors $F\otimes\bar{K}$ appartient à $\Oo(A,\bar{K})$, et la semi-simplicité de $F\otimes\bar{K}$ dans $\F(A,\bar{K})$ entraîne celle de $F$ dans $\F(A,K)$.

Si $K$ est de caractéristique $0$ et que $k$ est un corps fini, Kuhn \cite{Ku-adv} a montré que $\F(k,K)$ est semi-simple, donc équivalente à un produit (infini) de copies de la catégorie des $K$-espaces vectoriels si $K$ est algébriquement clos (puisque les corps d'endomorphismes des simples de $\F(k,K)$ sont de dimension finie sur $K$, $k$ étant fini). On en déduit par récurrence sur $n\in\mathbb{N}$ que, si $k_1,\dots,k_n$ sont des corps finis, alors $\F(\mathbf{P}(k_1)\times\dots\mathbf{P}(k_n);K)$ est équivalente à un produit de copies de la catégorie des $K$-espaces vectoriels.

Par ailleurs, il est classique (voir par exemple \cite[th.~3.2]{FS} et \cite[(2.6e)]{Green}) que la catégorie $\Pp(K,K)$ est semi-simple.

En utilisant la proposition \ref{pr-strglobsp}, le corollaire \ref{cor-sp_pt} et le théorème \ref{th-struct_ord}, on déduit de ce qui précède que tout foncteur de $\Oo^\cf(A,K)$ est semi-simple, donc également tout foncteur de $\Oo(A,K)$.

Montrons maintenant la deuxième assertion ($K$ est donc de caractéristique quelconque). Le fait qu'un foncteur de type fini de $\Oo(A,K)$ soit noethérien (resp. fini s'il n'y a pas de morphisme d'anneaux $A\to K$ d'image finie) résulte du théorème \ref{th-struct_ord}, de la proposition \ref{pr-strglobsp}, du corollaire \ref{cor-sp_pt} et du théorème~\ref{th-PSS} (resp. de la proposition~\ref{pr-PSmodif}). Si $A$ possède un quotient fini $\FF_q$ se plongeant dans $K$, le foncteur canonique $\F(\FF_q,K)\to\F(A,K)$ est pleinement fidèle et a une image essentielle prébilocalisante ; comme $\F(\FF_q,K)=\Oo(\FF_q,K)$ (cf. proposition~\ref{diago-corfini}) n'est pas localement finie, $\F(A,K)$ ne peut pas non plus l'être.
\end{proof}

On peut également établir que les foncteurs de type fini de $\Oo(A,K)$ vérifient la propriété $pf_\infty$ (dans $\F(A,K)$), en utilisant le théorème~\ref{th-struct_ord} et \cite{DT-schw}. Nous montrerons en fait ultérieurement un résulat plus général (grâce aux propositions \ref{pr-hypdia-sp}.\ref{it-hypdia-sp} et \ref{pr-dfsp}), le théorème~\ref{th-pfi_corps}.

\subsection{Structure de $\F^\df(k,K)$ ($k$ corps parfait infini de caractéristique $p$)}

Lorsque $K$ est de caractéristique nulle, si $k$ est un corps infini, alors tout foncteur de type fini de $\F^\df(k,K)$ est polynomial (et donc constant si $k$ est de caractéristique première), par \cite[cor.~4.19]{DTV}. Si l'on suppose de plus que $k$ est une extension algébrique de $\mathbb{Q}$, alors $(k,K)$ vérifie la condition (CAD)$^+$, de sorte que le corollaire~\ref{cor-fppc0ss} montre que tout foncteur de $\F^\df(k,K)$ est semi-simple. De plus, dans ce cadre, tout foncteur analytique de $\F(k,K)$ admettant une décomposition en poids est \hd\ si $K$ est un corps de décomposition de $\mathbf{P}(k)$, grâce à la proposition~\ref{pr-hdpolpoi}.

Nous donnons ci-après des résultats analogues lorsque la caractéristique $p$ de $K$ est première, qui sont alors plus difficiles à établir directement.

\begin{prop}\label{pr-corpspftcarp} Soit $k$ un corps parfait infini de caractéristique $p>0$. Supposons que $K$ est un corps de décomposition de $\mathbf{P}(k)$. Si $F$ est un foncteur de type fini de $\F(k,K)$, alors les conditions suivantes sont équivalentes :
\begin{enumerate}
\item\label{itkK1} $F$ possède une structure polynomiale stricte ;
\item\label{itkK2} $F$ est \hd\,;
\item\label{itkK3} $F$ est à valeurs de dimensions finies ;
\item\label{itkK4} $F$ est polynomial et possède une décomposition en poids faible.
\end{enumerate}
\end{prop}

\begin{proof}
L'équivalence entre \ref{itkK1} et \ref{itkK2} est un cas particulier du théorème~\ref{th-struct_ord}, et l'implication \ref{itkK1}$\Rightarrow$\ref{itkK3} est immédiate. Si $F$ est à valeurs de dimensions finies, alors $F$ est \ph\ par le théorème~\ref{th-DTVglob}, donc polynomial par le corollaire~\ref{cor-EMLpol-pol}. Comme $K$ est un corps de décomposition de $\mathbf{P}(k)$, il s'ensuit que $F$ possède une décomposition en poids faible par la proposition~\ref{pr-vdfdppol}. Ainsi, on a bien \ref{itkK3}$\Rightarrow$\ref{itkK4}. Enfin, la proposition~\ref{pr-hdpolpoi} fournit l'implication \ref{itkK4}$\Rightarrow$\ref{itkK2}.
\end{proof}

\begin{rema}
Il découlera des résultats du chapitre~\ref{shts} que, sous les hypothèses de la proposition~\ref{pr-corpspftcarp}, les conditions de l'énoncé sont également équivalentes au caractère \htr\ de $F$, mais nous ne savons pas le montrer avec les seuls résultats du présent chapitre. En revanche, si $k$ est localement fini, cette équivalence est claire, comme nous allons le voir dans la proposition suivante.
\end{rema}

L'énoncé suivant illustre à quel point il est plus fort pour un foncteur d'être \hd\ que diagonalisant --- en effet, sous ses hypothèses, \emph{tout} foncteur de $\F(k,K)$ est diagonalisant par la proposition~\ref{diago-fpbar}.

\begin{prop}\label{pr-hdcorlf} Supposons que $k$ est un sous-corps infini et localement fini de $K$.
Étant donné un foncteur $F$ de type fini de $\F(k,K)$, les assertions suivantes sont équivalentes :
\begin{enumerate}
\item\label{itkKlf1} $F$ possède une structure polynomiale stricte ;
\item\label{itkKlf2} $F$ est \hd\,;
\item\label{itkKlf3} $F$ est \htr\,;
\item\label{itkKlf4} $F$ est à valeurs de dimensions finies ;
\item\label{itkKlf5} $F\circ\oplus$ est fini ;
\item\label{itkKlf6} $F\circ\oplus$ est noethérien ;
\item\label{itkKlf7} $F\circ\oplus$ vérifie la condition $\mathrm{(Dec)}$ ;
\item\label{itkKlf8} $F$ est polynomial et possède une décomposition en poids faible ;
\item\label{itkKlf9} $F$ est polynomial et fini ;
\item\label{itkKlf10} $F$ est polynomial et noethérien ;
\item\label{itkKlf11} $F$ est polynomial et vérifie la condition $\mathrm{(Dec)}$.
\end{enumerate}
\end{prop}

\begin{proof} Comme $k$ est une extension algébrique de $\FF_p$, c'est un corps parfait, et le fait que $k$ soit un sous-corps de $K$ entraîne que $K$ est un corps de décomposition de $\mathbf{P}(k)$. Ainsi, l'équivalence entre \ref{itkKlf1}, \ref{itkKlf2}, \ref{itkKlf4} et \ref{itkKlf8} est un cas particulier de la proposition~\ref{pr-corpspftcarp}. L'équivalence entre \ref{itkKlf2} et \ref{itkKlf3} découle pour sa part de la proposition~\ref{prto}.

L'implication \ref{itkKlf5}$\Rightarrow$\ref{itkKlf6} est triviale, et \ref{itkKlf6}$\Rightarrow$\ref{itkKlf7} découle du corollaire~\ref{cor-ntcf}, tandis que \ref{itkKlf7}$\Rightarrow$\ref{itkKlf2} se déduit du corollaire~\ref{cor-poidec}, qui montre aussi que \ref{itkKlf11} entraîne \ref{itkKlf8} .

Les implications \ref{itkKlf1}$\Rightarrow$\ref{itkKlf9}$\Rightarrow$\ref{itkKlf10}$\Rightarrow$\ref{itkKlf11} sont claires.

Pour conclure, il suffit de noter que \ref{itkKlf1} entraîne \ref{itkKlf5}, parce que $F\circ\oplus$ possède une structure polynomiale stricte et est de type fini si c'est le cas de $F$.
\end{proof}

\begin{rema}Sous les hypothèses de la proposition précédente, il existe dans $\F(k,K)$ des foncteurs additifs non \hds\ (par exemple $V\mapsto V\otimes_\mathbb{Z}K$), et des foncteurs simples à poids non \hds\ (il existe dans $\mon(k_\mu,K_\mu)$ des poids $w$ non polynomiaux, et le foncteur $P^k_w$ de $\F(k,K)$ possède au moins un quotient simple, homogène de poids $w$, qui ne peut être \hd).
\end{rema}

\begin{rema}
Lorsque $k$ est une extension algébrique de $\mathbb{Q}$ et $K$ un surcorps de la clôture galoisienne de $k$ sur $\mathbb{Q}$, plusieurs des équivalences de la proposition~\ref{pr-hdcorlf} subsistent (cf. le début de ce paragraphe). Nous ne savons toutefois pas si sa condition \ref{itkKlf7}, ou même \ref{itkKlf5}, entraîne \ref{itkKlf4}.
\end{rema}

\subsection{Fonctions polynomiales homogènes}

Nous terminons cette section par une application à un énoncé où aucune catégorie de foncteurs n'apparaît. Comme dans tout le reste du texte, la notion de fonction polynomiale (entre groupes abéliens) s'entend au sens d'Eilenberg-MacLane \cite[§\,8]{EML}.

\begin{defi}\label{def-Fonction_polhom} Soient $\pi : A\to K$ un poids, $V$ un $A$-module et $E$ un $K$-espace vectoriel. On dit qu'une fonction $f : V\to E$ est \index{termin}{homogene@homogène!\emph{(fonction polynomiale homogène)}|(} \emph{homogène de poids $\pi$} si $\forall (\lambda,v)\in A\times V\quad f(\lambda.v)=\pi(\lambda).f(v)$.

Lorsque $A=K$, on dit que $f$ est $d$-homogène, pour $d\in\mathbb{N}$, si $f$ est homogène de poids $\lambda\mapsto\lambda^d$.
\end{defi}

\begin{prop}\label{fonctions-pol-homog} Supposons que le corps $K$ est parfait et que $\pi : A\to K$ est un poids polynomial déployé. Supposons  également que l'une des conditions suivantes est réalisée :
\begin{enumerate}
\item $p=0$ et la condition \textnormal{(CAD)}\index{nota}{CAD@(CAD), (CAD)$^+$ \emph{(conditions d'annulation des dérivations)}} est vérifiée ;
\item $p>0$ et la fonction $A\to A\quad x\mapsto x^p$ est surjective.
\end{enumerate}

 Soient $n\in\mathbb{N}$ et $f : A^n\to K$ une \index{termin}{polynomial(e)!fonction|(} fonction polynomiale homogène de poids $\pi$. Alors il existe des morphismes d'anneaux $\alpha_1,\dots,\alpha_r : A\to K$ et un polynôme $P\in K[X_{i,j}]_{1\le i\le r,1\le j\le n}$ tels que $\alpha_1\dots\alpha_r=\pi$ et
\[\forall (v_1,\dots,v_n)\in A^n\qquad f(v_1,\dots,v_n)=P(\alpha_i(v_j)).\] 
\end{prop}

\begin{proof} Soit $F\in\mathrm{Ob}\,\F(A,K)$ le sous-foncteur de $J:=K^{\mathbf{P}(A)(-,A)}$ engendré par $f\in K^{A^n}\simeq J(A^n)$. Comme $f$ est une fonction polynomiale, \cite[lemme~2.12]{DTV} montre que $F$ est un foncteur polynomial. Du fait que la fonction $f$ est homogène de poids $\pi$, le foncteur $F$ est homogène de poids fort $\pi$. La proposition~\ref{pr-hdpolpoi2} montre donc que $F$ est \hd. Comme $F$ est de type fini, il possède une décomposition de type $FP$ stricte, par le théorème~\ref{th-struct_ord}.
La conclusion s'ensuit à partir des observations suivantes :
\begin{enumerate}
\item sauf dans le cas trivial où $f$ est nulle, l'espace vectoriel $F(A)$ est de dimension $1$, engendré par $\pi\in K^A$ ;
\item l'isomorphisme $\phi : K\xrightarrow{\simeq}F(A)$ qu'on en déduit, et le morphisme $\psi : K\to F(A^n)$ donné par $f$ fournissent une fonction
$$A^n\simeq\mathbf{P}(A)(A^n,A)\xrightarrow{F_*}\mathrm{Hom}_K(F(A^n),F(A))\xrightarrow{\mathrm{Hom}_K(\psi,\phi^{-1})}\mathrm{Hom}_K(K,K)\simeq K$$
égale à $f$ ;
\item si $k$ est un corps fini se plongeant dans $K$, toute fonction $k^n\to K$ est induite par un polynôme en $n$ variables sur $K$ et un plongement de $k$ dans $K$.
\end{enumerate}
\end{proof}

\begin{coro}\label{cor-aprith-car0} Supposons que $K$ est une extension algébrique de $\mathbb{Q}$. Soient $n, d\in\mathbb{N}$ et $f : K^n\to K$ une fonction polynomiale $d$-homogène. Alors il existe un polynôme $P\in K[X_1,\dots,X_n]$, homogène de degré $d$, tel que
$$\forall (v_1,\dots,v_n)\in K^n\qquad f(v_1,\dots,v_n)=P(v_1,\dots,v_n).$$
\end{coro}

\begin{rema}\label{rq-car0der} En caractéristique $p=0$, l'hypothèse d'algébricité de  $K$ sur $\mathbb{Q}$ équivaut à la nullité de toutes les dérivations absolues de $K$ ; elle est également \emph{nécessaire} pour que la conclusion sur les fonctions polynomiales homogènes vaille. En effet, si $\mathrm{d} : K\to K$ est une dérivation, la fonction  $K^2\to K\quad (x,y)\mapsto x.\mathrm{d}(y)-y.\mathrm{d}(x)$ est $2$-homogène et polynomiale de degré au plus $2$, et l'on vérifie sans peine qu'elle ne satisfait la conclusion du corollaire~\ref{cor-aprith-car0} que si $\mathrm{d}$ est nulle.

Plus généralement, les fonctions polynomiales homogènes \guillemotleft~exotiques~\guillemotright\ de degré quelconque semblent liées aux dérivations supérieures.
\end{rema}

\begin{coro}\label{cor-aprith-carp} Supposons que $K$ est parfait de caractéristique $p>0$. Soient $n, d\in\mathbb{N}$ et $f : K^n\to K$ une fonction polynomiale $d$-homogène.  \index{termin}{homogene@homogène!\emph{(fonction polynomiale homogène)}|)} Alors il existe $r\in\mathbb{N}$ et un polynôme $P\in K[X_1,\dots,X_n]$, homogène de degré $dp^r$, de sorte que
$$\forall (v_1,\dots,v_n)\in K^n\qquad f(v_1,\dots,v_n)=P(v_1^{p^{-r}},\dots,v_n^{p^{-r}}).$$\index{termin}{polynomial(e)!fonction|)}
\end{coro}

\begin{rema} En caractéristique $p>0$, l'hypothèse de perfection du corps $K$ équivaut  à la nullité de toutes les dérivations absolues de $K$, et l'on peut formuler les mêmes observations qu'à la remarque~\ref{rq-car0der}.
\end{rema}

\chapter{Foncteurs \hts}\label{shts}

\begin{cvi}
Dans ce chapitre, $\A$ désigne toujours une catégorie additive $A$-linéaire essentiellement petite et $\E$ une catégorie de Grothendieck $K$-linéaire.
\end{cvi}

Notre objectif principal consiste à caractériser les foncteurs de $\Tt_A^\cf(\A,\E)$ d'une manière analogue à ce qui a été fait pour ceux de $\Oo_A^\cf(\A,\E)$ au chapitre précédent. Pour y parvenir, nous commencerons par traiter deux cas particuliers fondamentaux :
\begin{itemize}
\item celui des foncteurs  de $\Tt_A^\cf(\A,\E)$ dont tous les poids partiels sont sans torsion (section~\ref{par-fhtmst}) : un argument combinatoire simple permet de voir qu'ils sont toujours \emph{polynomiaux} ;
\item celui des foncteurs $F$ de $\Tt_A^\cf(\A,\E)$ dont tous les poids partiels sont de torsion (section~\ref{par-fhtmf}), pour lesquels on montre que $\AF$ est fini si $K$ est de caractéristique nulle ou $A$ noethérien (voir le corollaire~\ref{cor-carhtper} pour une caractérisation dans le cas général).
\end{itemize}

On montre enfin à la section~\ref{ss-tshtr} comment le cas général peut se ramener aux deux cas particuliers précédents, en utilisant un argument de lieu d'annulation des poids partiels et un principe formel de \guillemotleft~séparation~\guillemotright\ de certains poids partiels établi à la section~\ref{ssdp}.

Les trois premières sections contiennent des préliminaires sur la propriété hyper-trigonalisante pour des foncteurs $F$ possédant une propriété supplémentaire forte telle que la polynomialité, l'Hom-polynomialité ou la finitude de l'anneau $\AF$.

\section{Propriétés élémentaires}\label{s1hts}

La proposition~\ref{pr-dpftjs} permet de déterminer quand tous les foncteurs de $\F(A,K)$ sont \hts\ ; plus généralement :

\begin{prop}\label{pr-tthtr} Si les conditions de la proposition~\ref{pr-dpftjs} sont satisfaites, alors tout foncteur de $\fct(\A,\E)$ appartient à $\Tt_A^\cf(\A,\E)$.
\end{prop}

On rappelle que $p$ désigne la caractéristique du corps $K$.

\begin{lemm}\label{lm-finjppp} Supposons $p>0$. Soit $F$ un foncteur de $\Tt_A(\A,\E)$. Considérons la fonction multiplicative canonique $\alpha : \AF\to K^{\pp(F)}$ dont les composantes $\AF\to K$ sont induites par les éléments de $\pp(F)$.
\begin{enumerate}
\item Soit $(x,y)\in\AF^2$. On a $\alpha(x)=\alpha(y)$ si et seulement s'il existe $r\in\mathbb{N}$ tel que $x^{p^r}=y^{p^r}$.
\item  Si de plus $F$ appartient à $\Tt_A^\cf(\A,\E)$, alors il existe $r\in\mathbb{N}$ tel que, pour tout $(x,y)\in\AF^2$, $\alpha(x)=\alpha(y)$ si et seulement si $x^{p^r}=y^{p^r}$.
\end{enumerate}
\end{lemm}

\begin{proof} Comme $F$ est la colimite filtrante de ses sous-foncteurs $G$ appartenant à $\Tt_A^\cf(\A,\E)$ et que $\pp(F)$ (resp. $\rf(F)$) est la réunion (resp. l'intersection) filtrante croissante (resp. décroissante) des $\pp(G)$ (resp. $\rf(G)$), il suffit de montrer la deuxième assertion.

On suppose donc $F$ dans $\Tt_A^\cf(\A,\E)$ ; ainsi, l'algèbre $\ac_A(F)$ est semi-primaire. Quitte à remplacer $A$ par $\AF$ et $\A$ par $\A/\rf(F)$, on peut supposer également $\rf(F)=0$. On choisit pour $r$ un entier naturel tel que $\rj(\ac(F))^{p^r}=0$.

La composée de la fonction $\alpha$ et de l'injection canonique $K^{\pp(F)}\hookrightarrow K^{\widetilde{\pp}(F)}$ (dont l'image est constituée des fonctions valant $1$ sur le morphisme trivial $A_\mu\to K_\mu$) égale la composée $A\xrightarrow{\bar{\chi}_F}\ac(F)\twoheadrightarrow\ac(F)/\rj(\ac(F))\xrightarrow{\simeq}K^{\widetilde{\pp}(F)}$, où le dernier isomorphisme est donné par la proposition~\ref{pr-poidspar_ac}.\ref{itppac3}. Comme la fonction $\bar{\chi}_F$ est injective (proposition~\ref{pr-acrad}), il s'ensuit qu'on a $\alpha(x)=\alpha(y)$ si et seulement si $([x]-[y])^{p^r}=0$ dans $K[A_\mu]$, relation équivalente à $[x^{p^r}]=[y^{p^r}]$ puisque $K$ est de caractéristique $p$, et donc à $x^{p^r}=y^{p^r}$.
\end{proof}

L'énoncé suivant, qui utilise la notation~\ref{notapereg} (page~\pageref{notapereg}), constitue une variante du précédent en caractéristique nulle.

\begin{lemm}\label{lm-finjpp} Supposons $p=0$. Soit $F$ un foncteur de $\Tt_A(\A,\E)$. Alors la fonction multiplicative canonique
$$\alpha : \AF^{\mathrm{pereg}}\hookrightarrow \AF\xrightarrow{(\pi)_{\pi\in\pp(F)}} K^{\pp(F)}$$
est injective.
\end{lemm}

\begin{proof} Raisonnant comme dans la démonstration du lemme~\ref{lm-finjppp}, on voit qu'il suffit de montrer l'énoncé lorsque $F$ appartient à $\Tt_A^\cf(\A,\E)$ et que $\rf(F)$ est nul. On voit également de même que la relation $\alpha(x)=\alpha(y)$ entraîne que l'image de $[x]-[y]$ dans $\ac(F)$ est nilpotente.

Si $x$ est un élément de $A^{\mathrm{pereg}}$\index{nota}{Apereg@$A^\mathrm{pereg}$}, alors la sous-algèbre $K[x]$ de $\ac(F)$ engendrée par $\chi_F(x)$ est finie séparable, car $x$ annule un polynôme $X^n-X$ (avec $n\ge 2$) de $K[X]$, qui est séparable puisque $K$ est de caractéristique $0$. Si $x$ et $y$ sont des éléments de $A^{\mathrm{pereg}}$, la sous-algèbre $K[x,y]$ de $\ac(F)$ qu'ils engendrent est donc également séparable, et en particulier réduite. Ainsi, si $\alpha(x)=\alpha(y)$, l'image de $[x]-[y]$ dans $\ac(F)$ est nilpotente, donc nulle, d'où $\chi_F(x)=\chi_F(y)$ et $x=y$ par la proposition~\ref{pr-acrad}, comme souhaité.
\end{proof}

\begin{lemm}\label{lm-htrqAf} Soient $F$ un foncteur de $\Tt_A(\A,\E)$ et $n\in\mathbb{N}^*$. Si le groupe $\AF^\times$ possède un élément $x$ d'ordre $n$, alors $n\in\nu(K)$.
\end{lemm}

\begin{proof} Si $p>0$, on peut supposer $n$ étranger à $p$, vu que $l\in\nu(K)\Leftrightarrow pl\in\nu(K)$ pour tout $l\in\mathbb{N}^*$.

Supposons $n\notin\nu(K)$. Alors il existe un élément non nul de $\mathbb{Z}/n$ appartenant au noyau de tout morphisme de groupes $\mathbb{Z}/n\to K^\times$. Il existe donc une puissance $y$ de $x$, différente de $1$, telle que $\pi(y)=1$ pour tout $\pi\in\pp(F)$. Les lemmes~\ref{lm-finjpp} (si $p=0$) ou \ref{lm-finjppp} (si $p>0$) montrent que c'est impossible, d'où la conclusion.
\end{proof}

Le résultat suivant fait écho à la proposition~\ref{pr-antipol_hd} et s'applique notamment aux foncteurs antipolynomiaux de $\F(A,K)$.

\begin{prop}\label{pr-htr_AFfini} Soit $F$ un foncteur de $\fct(\A,\E)$ tel que l'anneau $\AF$ soit fini. Les propriétés suivantes sont équivalentes : 
\begin{enumerate}
\item\label{ihtr1} $F$ appartient à $\Tt_A^\cf(\A,\E)$ ;
\item\label{ihtr2} $F$ est \htr\;;
\item\label{ihtr3} la $K$-algèbre $K[\AF^\times]$ est déployée ;
\item\label{ihtr4} la $K$-algèbre $K[(\AF)_\mu]$ est déployée ;
\item\label{ihtr5} le groupe $\AF^\times$ est somme directe de groupes cycliques dont les ordres appartiennent à $\nu(K)$.
\end{enumerate}
\end{prop}

\begin{proof} L'implication \ref{ihtr1}$\Rightarrow$\ref{ihtr2} est triviale. Le lemme~\ref{lm-htrqAf} montre que \ref{ihtr2} entraîne \ref{ihtr5}. Les implications \ref{ihtr5}$\Rightarrow$\ref{ihtr4}$\Rightarrow$\ref{ihtr3}$\Rightarrow$\ref{ihtr1} découlent de la proposition~\ref{pr-dpftjs}.
\end{proof}

\section{Hyper-trigonalisabilité et polynomialité}\label{phtpol}

La proposition~\ref{pr-decpfa} indique quand tout foncteur polynomial de $\F(A,K)$ est \htr\ ; plus généralement :

\begin{prop}\label{pr-ttphtr} Si les conditions de la proposition~\ref{pr-decpfa} sont satisfaites, alors tout foncteur polynomial (resp. analytique) de $\fct(\A,\E)$ appartient à $\Tt_A^\cf(\A,\E)$ (resp. $\Tt_A(\A,\E)$).
\end{prop}

\begin{proof} Il suffit de démontrer l'assertion relative aux foncteurs polynomiaux, puisque $\Tt_A(\A,\E)$ est stable par colimites. Elle résulte de la proposition~\ref{pr-decpfa} appliquée au spectre du foncteur considéré (cf. \ref{pr-htdcf}.\ref{itsansf}) précomposé avec le foncteur de somme directe --- en effet, si l'algèbre $A_K$ est semi-primaire déployée, il en est de même pour $(A\times A)_K\simeq A_K\times A_K$. 
\end{proof}

\begin{lemm}\label{lm-Tdhypertrig} Soient $d\in\mathbb{N}$ et $F$ un foncteur de $\pol_d(\A,\E)$. Les assertions suivantes sont équivalentes :
\begin{enumerate}
\item $F$ appartient à $\Tt_A^\cf(\A,\E)$ ;
\item $cr_d(F)$ possède une décomposition en poids faible finie (relativement à l'action de $A^d$) et $\ppol_{d-1}(F)$ appartient à $\Tt_A^\cf(\A,\E)$ ;
\item $cr_d(F)$ possède une décomposition en poids faible finie (relativement à l'action de $A^d$) et $\qpol_{d-1}(F)$ appartient à $\Tt_A^\cf(\A,\E)$.
\end{enumerate}

Si ces conditions sont vérifiées et que de plus $\ppol_{d-1}(F)$ ou $\qpol_{d-1}(F)$ possède une filtration finie dont le gradué associé est \hd, alors $F$ possède une filtration finie dont le gradué associé est \hd.
\end{lemm}

\begin{proof} Cela découle des suites exactes \eqref{eq-qpol} et \eqref{eq-ppol} (page~\pageref{eq-qpol}) et de ce que si $X$ est un multifoncteur de $\fct(\A^d,\E)$ additif par rapport à chacune des $d$ variables et possédant une décomposition en poids faible finie (relativement à l'action de $A^d_\mu$), alors $X$ possède une filtration finie dont le gradué associé appartient à $\Oo^\cf_{A^d}(\A^d,\E)$, donc $\delta_d^*X$ possède une filtration finie dont le gradué associé appartient à $\Oo^\cf_A(\A,\E)$.
\end{proof}

\begin{prop}\label{pr-polhd_ht} Tout foncteur polynomial de $\Tt_A^\cf(\A,\E)$ possède une filtration finie dont le gradué associé appartient à $\Oo_A^\cf(\A,\E)$.
\end{prop}

\begin{proof} La propriété se déduit du lemme~\ref{lm-Tdhypertrig} par récurrence sur le degré polynomial.
\end{proof}

\begin{rema}\label{rq-hthdext} Ce phénomène est spécifique aux foncteurs polynomiaux.

Supposons par exemple que $q$ est une puissance d'un nombre premier, différente de $2$, que $p$ est un diviseur premier de $q-1$ et $K$ un corps de caractéristique $p$ contenant assez de racines de l'unité. Alors  $\F(\mathbb{F}_q,K)$ ne contient aucun foncteur \hd\ non constant (proposition~\ref{diago-corfini}), tandis que tous ses foncteurs sont \hts\ (proposition~\ref{pr-tthtr}). En particulier, cette catégorie contient des foncteurs simples \hts\ mais pas \hds.

Il peut même exister des foncteurs \emph{\phs}, noethériens, à valeurs de dimensions finies et \hts\ mais n'ayant aucune filtration dont le gradué associé est \hd. Pour l'illustrer, supposons que $A$ est un $p$-anneau fini non semi-simple, et que $K$ contient assez de racines de l'unité. Alors tout foncteur de $\F(A,K)$ est \htr\ (proposition~\ref{pr-tthtr}), et un foncteur $F$ de $\F(A,K)$ est \hd\ si et seulement si $A_{/F}$ est un quotient de $A/\rj(A)$ (proposition~\ref{pr-antipol_hd}). On en déduit facilement que le plus grand sous-foncteur \hd\ $G$ de $P^A$ est $V\mapsto K[I.V]$, où $I$ est l'idéal annulateur de $\rj(A)$ dans $A$, et que $P^A/G$ est un foncteur non nul de $\F(A,K)$ ne possédant aucun sous-foncteur \hd\ non nul. Ainsi $P^A$ est un foncteur \ph\ (comme tout foncteur de $\F(A,K)$), noethérien (car $\F(A,K)$ est localement noethérienne), à valeurs de dimensions finies et \htr\ ne possédant aucune filtration dont le gradué associé est \hd.
\end{rema}

L'énoncé suivant utilise la construction $W_d$ introduite au début du chapitre~\ref{sec-finpolp}.

\begin{prop}\label{pr-polpoihtr} Soient $d\in\mathbb{N}$ et $\pi\in\mon(A_\mu,K_\mu)$. Les assertions suivantes sont équivalentes.
\begin{enumerate}
\item\label{ithtrp1} Tout foncteur de $\pol_d(A,K)$ homogène de poids faible $\pi$ est \htr.
\item\label{ithtrp2} Pour tout entier $n\le d$, la $K$-algèbre $W_n(A,\pi)$ est semi-primaire déployée.
\end{enumerate}
\end{prop}

\begin{proof} Si $T^n_\pi$ désigne le plus grand quotient homogène de poids fort $\pi$ de $V\mapsto (V\otimes_\mathbb{Z}K)^{\otimes d}$, on a
$$cr_n(T^n_\pi)(V_1,\dots,V_n)\simeq\bigoplus_{\sigma\in \Si_n}(V_{\sigma(1)}\otimes_\mathbb{Z}\dots\otimes_\mathbb{Z} V_{\sigma(n)}\otimes_\mathbb{Z}K)\underset{A_K^{\otimes n}}{\otimes} W_n(A,\pi)\,,$$
d'où $\mathfrak{Z}_{A^n_\mu}(cr_n(T^n_\pi))\simeq W_n(A,\pi)$ et l'implication \ref{ithtrp1}$\Rightarrow$\ref{ithtrp2}.

Pour montrer la réciproque, on peut se contenter de montrer que, sous l'hypothèse \ref{ithtrp2}, tout foncteur polynomial $F$ de degré $n\le d$, de type fini, et homogène de poids fort $\pi$ est \htr. La proposition~\ref{pr-Wcrimm} et le caractère semi-primaire déployé de $W_n(A,\pi)$ montrent que $cr_n(F)$ possède une décomposition en poids faible, nécessairement finie puisque $F$ est de type fini. Par récurrence sur $n$, on peut supposer que tout foncteur polynomial de degré $<n$ de $\F(A,K)$ est \htr, ce qui implique que le foncteur polynomial de type fini $\qpol_{n-1}(F)$ appartient à $\Tt^\cf(A,K)$. On conclut grâce au lemme~\ref{lm-Tdhypertrig}.
\end{proof}

\begin{coro}\label{cor-polpoihtr} Supposons que $K$ est de caractéristique nulle et que $\pi : A_\mu\to K_\mu$ est un poids polynomial déployé. Alors tout foncteur polynomial homogène de poids faible $\pi$ de $\F(A,K)$ est \htr.
\end{coro}

\begin{proof} Cela découle des propositions~\ref{pr-polpoihtr}, \ref{pr-pradw} et~\ref{prw-fond}.\ref{pr-wdepl}.
\end{proof}

\begin{rema} En revanche, sans hypothèse polynomiale, il est exceptionnel que tout foncteur homogène de poids $\pi$ soit \htr, même en caractéristique nulle, comme le montrent les propositions~\ref{pr-algcarP} et~\ref{pr-dpftjs}.
\end{rema}

En caractéristique $p>0$, on a l'analogue suivant, qui nécessite davantage d'hypothèses puisque les algèbres $W_n(A,\pi)$ ne sont alors plus automatiquement semi-primaires (cf. exemple~\ref{ex-Wcorimp}) :

\begin{prop}\label{pr-polpoihtrp} Supposons $p>0$ et l'anneau $A$ de type fini ou, plus généralement, que $A/p$ est une algèbre essentiellement de type fini sur une $\FF_p$-algèbre $p$-parfaite. Supposons également que $K$ est un corps de décomposition de $\mathbf{P}(A)$.

Alors tout foncteur polynomial de $\fct(\A,\E)$ possédant une décomposition en poids faible est \htr.
\end{prop}

\begin{proof}
Il suffit de montrer que tout foncteur $F$ de $\fct(\A,\E)$ possédant une décomposition en poids faible \textit{finie} est \htr\ (cf. remarque~\ref{rq-decfini}). Le foncteur $\Sp(F)$ de $\F(A,K)$ possède alors une décomposition en poids faible finie (proposition~\ref{pr-sppoi}), et est polynomial (corollaire~\ref{cor-sppol}). Le corollaire~\ref{cor-phdf} implique que $\Sp(F)$ appartient à $\F^\df(A,K)$. La proposition~\ref{pr-df_htr} montre alors que $\Sp(F)$ est \htr. Il s'ensuit que $F$ lui-même est \htr, d'où la conclusion.
\end{proof}

\begin{rema}\label{rq-complPAScordec}
\begin{enumerate}
    \item Si l'on supprime l'hypothèse que $K$ est un corps de décomposition de $\mathbf{P}(A)$, la conclusion de la proposition~\ref{pr-polpoihtrp} subsiste sous la forme affaiblie suivante : si $F$ est un foncteur polynomial de $\fct(\A,\E)$ possédant une décomposition en poids faible \emph{finie}, alors il existe une extension finie $L$ du corps $K$ tel que le foncteur $F\underset{K}{\boxtimes} L$ de $\fct(\A,\E\underset{K}{\boxtimes} L)$ soit \htr. La démonstration est la même sinon qu'il faut remplacer, à la fin, l'usage de la proposition~\ref{pr-df_htr} par celui du corollaire~\ref{cor-excb} (appliqué au bifoncteur $\Sp(F)\circ\oplus$).
    \item Dans la proposition~\ref{pr-polpoihtrp} ou la variante mentionnée ci-dessus, on peut supposer seulement que $F$ est \ph, dès lors que la condition (MTF) est satisfaite (par exemple, dans $\F(A,K)$). La même démonstration s'applique, en remplaçant le corollaire~\ref{cor-sppol} par le corollaire~\ref{cor-phsp}.
\end{enumerate}
\end{rema}

La proposition suivante, analogue pour les foncteurs polynomiaux \hts\ de résultats de la section~\ref{parfpp}, nous servira notamment pour établir le théorème~\ref{th-htr_ln}.

\begin{prop}\label{pr-polhtfini} Les assertions suivantes sont équivalentes :
\begin{enumerate}
\item la condition $\mathrm{(CFD)}$\index{nota}{CFD@(CFD), (CFD)$^+$ \emph{(conditions de finitude sur les dérivations)}} (introduite page~\pageref{eqcfd}) est vérifiée ;
\item\label{itG2} tout foncteur additif de type fini de $\F(A,K)$ possédant une décomposition en poids faible est noethérien ;
\item tout foncteur additif de type fini de $\F(A,K)$ possédant une décomposition en poids faible est fini ;
\item\label{itG4} tout foncteur additif de type fini de $\F(A,K)$ possédant une décomposition en poids faible est à valeurs de dimensions finies ;
\item\label{itG5} pour toute catégorie de Grothendieck $K$-linéaire $\E$, tout foncteur additif de type fini de $\fct(\mathbf{P}(A),\E)$ possédant une décomposition en poids faible est à valeurs de type fini ;
\item\label{itG6} pour toute catégorie de Grothendieck $K$-linéaire $\E$, tout foncteur polynomial \htr\ de type fini de $\fct(\mathbf{P}(A),\E)$ est à valeurs de type fini ;
\item\label{itG7} pour toute catégorie de Grothendieck $K$-linéaire localement noethérienne (resp. localement finie) $\E$, tout foncteur polynomial \htr\ de type fini de $\fct(\mathbf{P}(A),\E)$ est noethérien (resp. fini).
\end{enumerate}
\end{prop}

\begin{proof}
Utilisant l'équivalence de catégories canonique $\mathbf{Add}(A,K)\simeq A_K\Md$ (cf. \eqref{eq-pteAdda}, page~\pageref{eq-pteAdda}) et la proposition~\ref{pr-poideploye}, on voit qu'un foncteur additif de type fini de $\F(A,K)$ possédant une décomposition en poids faible correspond à une somme directe finie de $A_K$-modules pour lesquels l'action de $A_K$ est nulle sur $\mathfrak{m}_\varphi^r$\index{nota}{m@$\mathfrak{m}_\varphi$ \emph{(pour $\varphi\in\mathbf{Ann}(A,K)$)}} pour un $\varphi\in\mathbf{Ann}(A,K)$ et un $r\in\mathbb{N}$, on obtient alors l'équivalence entre les quatre premières conditions grâce au corollaire~\ref{cor-hole} et à l'isomorphisme~\eqref{eq-ident_der} (page~\pageref{eq-ident_der}).

Pour la même raison, la première  assertion implique la cinquième, un foncteur de type fini de $\mathbf{Add}(\mathbf{P}(A),\E)\simeq\E\underset{K}{\boxtimes}A_K$ (cf. \eqref{eq-pteAdd}, page~\pageref{eq-pteAdd}) possédant une décomposition en poids faible étant une somme directe finie d'objets provenant de $\E\underset{K}{\boxtimes}A_K/\mathfrak{m}_\varphi^r$ pour un $\varphi\in\mathbf{Ann}(A,K)$ et un $r\in\mathbb{N}$.

Montrons l'implication \ref{itG5}$\Rightarrow$\ref{itG6}. Supposant \ref{itG5} vérifiée, on montre d'abord que tout foncteur \htr\ de type fini de $\mathbf{Add}_d(\mathbf{P}(A),\E)$ est à valeurs de type fini. Cela découle par récurrence sur $d$ de \ref{itG5} et des équivalences de catégories canoniques  $\mathbf{Add}_d(\mathbf{P}(A),\E)\simeq\mathbf{Add}(\mathbf{P}(A),\mathbf{Add}_{d-1}(\mathbf{P}(A),\E))$. On montre maintenant par récurrence sur $d$ qu'un foncteur \htr\ de type fini $F$ de $\pol_d(\mathbf{P}(A),\E)$ est à valeurs de type fini. On peut supposer $d>1$. Le foncteur $cr_d(F)$ de $\mathbf{Add}_d(\mathbf{P}(A),\E)$ est comme $F$ de type fini et \htr, il est donc à valeurs de type fini d'après ce qu'on vient de voir. Le foncteur $\qpol_{d-1}(F)$ est également à valeurs de type fini d'après l'hypothèse de récurrence. La conclusion résulte donc de la suite exacte naturelle $\delta_d^*cr_d(F)\to F\to\qpol_{d-1}(F)\to 0$.

L'implication \ref{itG6}$\Rightarrow$\ref{itG7} découle de \cite[lm~11.10]{DTV} et \cite[lm~6.20]{DT-schw} ; l'implication \ref{itG7}$\Rightarrow$\ref{itG2} s'obtient immédiatement en prenant $\E=K\Md$, d'où la proposition.
\end{proof}

\section{Hyper-trigonalisabilité et Hom-polynomialité}\label{phthpol}

La proposition~\ref{pr-decpfa} donne des conditions suffisantes --- et également nécessaires dans $\F(A,K)$ --- pour que tout foncteur \ph\ soit \htr\ :

\begin{prop}\label{pr-hpol_htr} Supposons les conditions de la proposition~\ref{pr-decpfa} satisfaites. Alors tout foncteur \ph\ $F$ de $\fct(\A,\E)$ est \htr.
\end{prop}

\begin{proof} Comme $\Tt_A(\A,\E)$ est stable par colimite, on peut supposer que $F$ est à support fini. Par la proposition~\ref{pr-chi_ph}, sa fonction caractéristique $\chi_F$ est polynomiale. Autrement dit, il existe $d\in\mathbb{N}$ tel que l'application linéaire $K[A]\to Z(\mathrm{End}(F\circ\oplus))$ prolongeant $\chi_F$ se factorise par la projection $K[A]\twoheadrightarrow Q_d(A)$, où $Q_d:=\qpol_d(P^A)\in\mathrm{Ob}\,\F(A,K)$.

La condition \ref{itdpfa4} de la proposition~\ref{pr-decpfa} montre que $Q_d$ possède une décomposition en poids faible finie, donc $\mathfrak{Z}(Q_d)\simeq Q_d(A)$ est une algèbre semi-primaire déployée. Ce qui précède montre donc que l'algèbre caractéristique $\ac_A(F)$ est semi-primaire déployée, ainsi $F$ appartient à $\Tt_A^\cf(\A,\E)$ (caractérisation \ref{tcf-it4} de la proposition et définition~\ref{pr-htdcf}), d'où le résultat.
\end{proof}

\begin{rema} Nous ignorons si, sous les hypothèses de la proposition~\ref{pr-decpfa}, tout foncteur \ph\ $F$ de $\fct(\A,\E)$ appartient à $\Tt_A^\cf(\A,\E)$.
%probablement NON, mais contre-ex. ?
\end{rema}

Nous allons maintenant donner des conditions suffisantes pour que tout foncteur de $\Tt^\cf(A,K)$ soit \ph. Nous nous limitons ici au cas d'une caractéristique $p$ première ; le corollaire~\ref{cor-hpap} ci-après montrera que, si $K$ est de caractéristique nulle, alors tout foncteur de $\Tt^\cf(A,K)$ est \ph\ (ou polynomial) si et seulement si l'anneau $A$ est sans quotient fini, au moins si $K$ contient assez de racines de l'unité.

Notre raisonnement reposera essentiellement sur des propriétés de l'algèbre $K[A_\mu]$, fournies par les deux lemmes suivants.

\begin{lemm}\label{lm-radnl} Si $p>0$ et que la caractéristique de $A$ est une puissance de $p$, alors la projection $K[A]\twoheadrightarrow K[A/\nil(A)]$ induit un isomorphisme d'algèbres $K[A_\mu]/\nil(K[A_\mu])\xrightarrow{\simeq} K[(A/\nil(A))_\mu]$.
\end{lemm}

\begin{proof} Pour $x, y\in A$, on a $x-y\in\nil(A)$ si et seulement si $x^{p^n}=y^{p^n}$ pour un $n\in\mathbb{N}$.
%Fait général à ériger en lemme séparé ? En tout début d'article ??? Ou remplacer par une réf. ?
Précisément, si $A$ est de caractéristique $p^r$, et si $x$ et $t$ sont des éléments de $A$ tels que $t^i=0$, alors la formule du binôme montre que $x^{p^n}=(x+t)^{p^n}$ dès que $p^n\ge i$ et $n\ge r+\frac{i-1}{p-1}$ ; réciproquement, si $x^{p^n}=y^{p^n}$, comme $A/\nil(A)$ est de caractéristique $p$ (ou $1$), on a $(x-y)^{p^n}\in\nil(A)$ donc $x-y\in\nil(A)$. La conclusion résulte maintenant de ce que l'élévation à la puissance $p^n$ est pour tout $n\in\mathbb{N}$ un endomorphisme de l'anneau $K[A_\mu]$ qui se restreint en un endomorphisme injectif de $K$.
\end{proof}

Dans le lemme suivant, on fait un usage répété de la notion de bi-idéal et de la notation $\star$ introduites page~\pageref{pbiid}, ainsi que de la remarque~\ref{rqstev}.

\begin{lemm}\label{lm-piaRad} Supposons $A$ semi-primaire et $A/\rj(A)$ fini de caractéristique $p$. Alors pour tout $n\in\mathbb{N}$, il existe $d\in\mathbb{N}$ tel que le bi-idéal\index{termin}{biideal@bi-idéal} de $K[A]$ engendré par $\rj(K[A_\mu])^n$ contienne $\bar{K}[A]^{\star d}$.
\end{lemm}

\begin{proof} Posons $J:=\rj(A)$, $N:=\min\{i\in\mathbb{N}\,|\,J^{i+1}=0\}$ et $I:=J^N\subset J$. 

Comme $A/J$ est fini de caractéristique $p$, l'idéal d'augmentation de $K[A/J]$ est $\star$-nilpotent \cite[chap.~VI, th.~1.2]{Passi}. Du fait que le noyau de l'épimorphisme canonique $K[A]\twoheadrightarrow K[A/J]$ est $K[A]\star\bar{K}[J]$, il existe donc $r\in\mathbb{N}$ tel que
\begin{equation}\label{eqstar}
\bar{K}[A]^{\star r}\subset K[A]\star\bar{K}[J]\;.
\end{equation}

Le lemme~\ref{lm-radnl} montre par ailleurs que $\rj(K[A_\mu])=K[A]\star\bar{K}[J]$.

Raisonnant par récurrence, on peut supposer la conclusion vraie
\begin{itemize}
\item pour $n-1$ : on dispose donc de $a\in\mathbb{N}$ tel que
\begin{equation}\label{recst1}
\bar{K}[A]^{\star a}\subset K[A]\star (K[A]\star\bar{K}[J])^{n-1}\; ;
\end{equation}
\item pour l'anneau $A/I$ : on dispose donc de $b\in\mathbb{N}$ tel que $\bar{K}[A/I]^{\star b}\subset K[A/I]\star (K[A/I]\star\bar{K}[J/I])^n$ ; comme le noyau du morphisme canonique $\bar{K}[A]^{\star b}\twoheadrightarrow\bar{K}[A/I]^{\star b}$ est $\bar{K}[A]^{\star b-1}\star\bar{K}[I]$, on en déduit
\begin{equation}\label{recst2}
\bar{K}[A]^{\star b}\subset (\bar{K}[A]^{\star b-1}\star\bar{K}[I])+(K[A]\star (K[A]\star\bar{K}[J])^n).
\end{equation}
\end{itemize}

Comme $I.J=0$, on a $(1+x).(1+t)=1+x+t$ pour $x\in I$ et $t\in J$. En utilisant \eqref{recst1}, on en déduit que
\begin{equation*}
K[A]\star (K[A]\star\bar{K}[J])^n\supset K[A]\star\big((K[A]\star\bar{K}[I]).(K[A]\star\bar{K}[J])^{n-1}\big)\supset
\end{equation*}
\begin{equation*}K[A]\star\big((K[A]\star\bar{K}[I]).(K[A]\star\bar{K}[J]^{\star a})\big)\supset K[A]\star\bar{K}[I]\star\bar{K}[J]^{\star a}\, ,
\end{equation*}
d'où, grâce à \eqref{eqstar},
\begin{equation}
K[A]\star (K[A]\star\bar{K}[J])^n\supset \bar{K}[A]^{\star ar}\star\bar{K}[I].
\end{equation}
Comme \eqref{recst2} entraîne
\begin{equation}
\bar{K}[A]^{\star c}\subset (\bar{K}[A]^{\star c-1}\star\bar{K}[I])+(K[A]\star (K[A]\star\bar{K}[J])^n)
\end{equation}
pour tout entier $c\ge b$, pour $d=\max\{b,ar+1\}$, on a 
$$\bar{K}[A]^{\star d}\subset K[A]\star (K[A]\star\bar{K}[J])^n\,,$$
ce qui achève la démonstration.
\end{proof}

\begin{prop}\label{pr-htr_impl_ph} Si $A$ est semi-primaire et que $A/\rj(A)$ est fini de caractéristique $p$, alors tout foncteur $F$ de $\Tt^\cf(A,K)$ est \ph.
\end{prop}

\begin{proof}
L'algèbre $\ac_A(F)$ est semi-primaire déployée, donc l'idéal caractéristique $\ic_A(F)$ contient une puissance du radical de $K[A_\mu]$. Comme $\ic_A(F)$ est un bi-idéal de $K[A]$, le lemme~\ref{lm-piaRad} fournit $d\in\mathbb{N}$ tel que $\ic_A(F)\supset\bar{K}[A]^{\star d}$. Autrement dit, $\ic_A(F)$ est un quotient de $K[A]/\bar{K}[A]^{\star d}$, ce qui signifie que la fonction caractéristique $\chi_F$ est polynomiale de degré au plus $d$. On conclut par la proposition~\ref{pr-chiph} et le corollaire~\ref{cor-phsp}.
\end{proof}

\section{Foncteurs de $\Tt^\cf(\A,\E)$ aux poids partiels sans torsion}\label{par-fhtmst}

\begin{lemm}\label{lm-ppst} Soit $F$ un foncteur de $\Tt_A^\cf(\A,\E)$. Il existe un entier $N$ tel que, pour tout multipoids $(w_1,\dots,w_d)\in\W_d(F)$, le nombre d'indices $i\in\llbracket 1,d\rrbracket$ tels que le poids partiel $w_i$ soit sans torsion dans le monoïde $\mon(A_\mu,K_\mu)$ soit majoré par $N$.
\end{lemm}

\begin{proof} Comme $F$ est \htr\ \emph{à caractère fini}, l'ensemble $\pp(F)$ de ses poids partiels est fini ; soient $n$ son cardinal et $N:=n^2$.

Si $t$ est un poids partiel sans torsion de $F$ et $(w_1,\dots,w_d)$ un multipoids de $F$, notons $r$ le nombre d'indices $i$ tels que $w_i=t$. La proposition~\ref{pr-poidspar}.\ref{itppst0} montre que $t^j$ appartient à $\pp(F)$ pour tout $j\in\llbracket 1,i\rrbracket$. Comme $t$ est sans torsion, ces éléments de $\mon(A_\mu,K_\mu)$ sont deux à deux distincts, d'où $i\le n$. Par conséquent, le nombre de composantes d'un multipoids de $F$ qui sont des poids partiels sans torsion est majoré par $n$ fois le nombre total de poids partiels sans torsion, donc par $N$, d'où le lemme.
\end{proof}

\begin{prop}\label{pr-raispol}  Soit $F$ un foncteur de $\Tt_A^\cf(\A,\E)$. Supposons que chaque élément de  $\pp(F)$ est sans torsion dans le monoïde $\mon(A_\mu,K_\mu)$.

Alors $F$ est polynomial.
\end{prop}

\begin{proof} Le lemme précédent fournit $N\in\mathbb{N}$ tel que tout multipoids de $F$ a au plus $N$ éléments. Autrement dit, le foncteur $cr_{N+1}(F)$, qui admet une décomposition en poids faible, ne possède aucun poids : il est nul, donc $F$ est polynomial de degré au plus $N$.
\end{proof}

\begin{rema} Les mêmes énoncés valent sans changement dans le cas plus général où $\A$, plutôt qu'une catégorie additive $A$-linéaire, est une petite catégorie monoïdale symétrique dont l'unité est objet nul munie d'une action d'un monoïde commutatif.
\end{rema}

\begin{rema} On peut a posteriori améliorer la borne sur le degré d'un foncteur $F$ vérifiant les hypothèses de la proposition~\ref{pr-raispol} déduite des démonstrations précédentes, à savoir $n^2$, où $n=\cd(\pp(F))$ : on a en fait $\deg(F)\le n$. Pour le voir, supposons que $S$ est un facteur de composition de $F$. Alors $S$ est, comme $F$, polynomial et \htr, de sorte que la proposition~\ref{pr-simplesPolHD} montre qu'il possède une structure polynomiale stricte (la décomposition n'est pas faible car l'hypothèse sur $F$ exclut tout morphisme d'anneaux $A\to K$ d'image finie comme poids partiel de $F$). On déduit alors directement de la proposition~\ref{pr-dppol} que le degré de $S$ est au plus égal au nombre de ses poids partiels, lui-même majoré par $n$. Le lemme~\ref{lm-degfc} fournit alors l'inégalité $\deg(F)\le n$ souhaitée.
\end{rema}

\section{Foncteurs de $\Tt^\cf(\A,\E)$ aux poids partiels de torsion}\label{par-fhtmf}

L'objectif de cette section consiste à établir le théorème suivant.

\begin{theo}\label{th-pptor_gal} Soit $F$ un foncteur de $\Tt_A^\cf(\A,\E)$ dont les poids partiels sont périodiques. Alors $\AF$ est isomorphe au produit d'un anneau fini et d'un anneau semi-primaire dont le quotient par le radical est un produit fini de sous-corps finis de $K$.
\end{theo}

La démonstration nécessitera un certain nombre d'étapes ; on traitera indépendamment le cas où la caractéristique de $K$ est première (au §\,\ref{parccp}) et celui où $K$ est de caractéristique nulle (au §\,\ref{pp0}). Curiseusement, c'est ce dernier cas le plus délicat, à l'opposé de la plupart des considérations sur les foncteurs de source additive et de but abélien (y compris celles du présent mémoire).

Avant d'aborder la démonstration du théorème~\ref{th-pptor_gal}, nous en tirons quelques conséquences directes.

\begin{coro}\label{cor-carppoiord} Sous les hypothèses du théorème~\ref{th-pptor_gal}, tous les poids partiels de $F$ sont ordinaires.
\end{coro}

\begin{proof}
Cela découle des propositions \ref{pr-poitjspol} et \ref{pr-decpord} ainsi que du théorème~\ref{th-pptor_gal}.
\end{proof}

\begin{coro}\label{cor-carhtper} Étant donné un foncteur $F$ de $\F(A,K)$, les conditions suivantes sont équivalentes :
\begin{enumerate}
\item\label{icth1} $F$ appartient à $\Tt^\cf(A,K)$ et ses poids partiels sont périodiques ;
\item\label{icth2} $F$ possède une décomposition à la Steinberg de type $HPAP$ et $\AF$ est isomorphe au produit d'un anneau fini dont le groupe multiplicatif est somme directe de groupes cycliques d'ordres dans $\nu(K)$ et d'un anneau semi-primaire dont le quotient par le radical est un produit fini de sous-corps finis de $K$.
\end{enumerate}
\end{coro}

\begin{proof} L'implication \ref{icth1}$\Rightarrow$\ref{icth2} résulte du théorème~\ref{th-pptor_gal}, du lemme~\ref{lm-htrqAf} et de la proposition~\ref{pr-htr_impl_ph}. La réciproque se déduit des propositions~\ref{pr-htr_AFfini} et~\ref{pr-hpol_htr}.
\end{proof}

La conclusion du théorème~\ref{th-pptor_gal} nous conduit à introduire la condition suivante (qui porte sur les Idéaux Maximaux à quotient Fini, d'où la notation) :
\smallskip

\noindent
(IMF)\;\;\; Si $\mathfrak{m}$ est un idéal maximal de $A$ tel que $A/\mathfrak{m}$ soit un sous-corps fini de $K$, alors $\mathfrak{m}/\mathfrak{m}^2$ est fini.\index{nota}{IMF@(IMF) \emph{(condition de finitude sur les idéaux maximaux à quotient fini)}|textbf}\label{pageimf}
\smallskip

Cette condition est en particulier vérifiée lorsque $p=0$ ou que $A$ est un anneau noethérien.

\begin{rema}\label{rq-CFD_IMF}
La remarque~\ref{rq-CFD_imagefinie} montre l'implication (CFD)$\Rightarrow$(IMF).\index{nota}{CFD@(CFD), (CFD)$^+$ \emph{(conditions de finitude sur les dérivations)}}
\end{rema}

Le corollaire~\ref{cor-carhtper} entraîne aussitôt le résultat suivant.

\begin{coro}\label{cor-IMF} Si la condition \textnormal{(IMF)} est vérifiée, alors pour tout foncteur $F$ de $\Tt_A^\cf(\A,\E)$ dont les poids partiels sont périodiques, l'anneau $\AF$ est fini.

Réciproquement, si tout foncteur $F$ de $\Tt^\cf(A,K)$ dont les poids partiels sont périodiques, l'anneau $\AF$ est fini, alors la condition \textnormal{(IMF)} est vérifiée.
\end{coro}

\subsection{Cas de la caractéristique première}\label{parccp}

\begin{hyp} Dans tout le §\,\ref{parccp}, on suppose que $p$ est un nombre premier.
\end{hyp}

Les résultats de ce paragraphe reposent sur l'énoncé d'algèbre commutative suivant.

\begin{lemm}\label{ann-puissances-finies} Soit $r\in\mathbb{N}$, posons $N:=p^r$. Supposons que le sous-ensemble $S:=\{x^N\,|\,x\in A\}$ de $A$ est fini. Alors $A$ est isomorphe au produit direct d'un anneau fini et d'un anneau $B$ possédant les propriétés suivantes : 
\begin{enumerate}
\item la caractéristique de $B$ est une puissance de $p$ ;
\item il existe un entier $i$ tel que, pour tout élément $x$ du radical $\mathfrak{m}$ de $B$, $x^i$ soit nul ;
\item l'anneau quotient $B/\mathfrak{m}$ est fini.
\end{enumerate}

En particulier, $A$ est quasi-parfait, et si $\mathfrak{m}$ est un idéal de type fini --- par exemple, si $A$ est noethérien --- alors $A$ est fini.
\end{lemm}

\begin{proof} Comme l'ensemble $\{x^N\,|\,x\in\mathbb{Z}\}$ est infini, $A$ ne peut pas être de caractéristique nulle. Par décomposition primaire, on peut supposer que la caractéristique de $A$ est une puissance d'un nombre premier $l$ ; on distinguera plus tard selon que $l$ égale ou non $p$.

Notons $R$ le sous-anneau de $A$ engendré par $S$. Comme $S$ est stable par multiplication et fini, et que $A$ est de caractéristique positive, $R$ est fini. Par conséquent, $R$ est un produit fini d'anneaux {\em locaux} finis, et cette décomposition induit une décomposition en produit de $A$, ce qui permet de supposer que $R$ est un anneau local. Comme il est artinien, il existe un entier $j>0$ tel que tout élément non inversible $t$ de $R$ vérifie $t^j=0$. Il s'ensuit que, si $x$ est un élément non inversible de $A$, on a $(x^N)^j=0$. Ainsi, $A$ est local et la puissance $i$-ème de tout élément de $\rj(A)$ est nulle, pour $i:=Nj$.

Le corps $A/\rj(A)$ est fini, car l'ensemble de ses puissances $N$-ièmes est fini.

Il ne reste plus qu'à démontrer que, si $p$ est inversible dans $A$, alors cet anneau est lui-même fini. L'endomorphisme $x\mapsto x^N$ du groupe multiplicatif $A^\times$ est d'image finie, et son noyau est également fini, car il s'injecte dans $(A/\rj(A))^\times$ --- en effet, si $t$ est un élément nilpotent non nul de $A$, on ne peut avoir $(1+t)^N=1$ (relation qui implique $N.t\in (t^2)$, donc $t\in (t^2)$ puisque $N\in A^\times$). Ainsi $A^\times$ est fini ; comme $1+x$ est inversible pour tout élément non inversible $x$ de $A$, il s'ensuit que $A$ est fini, ce qui achève la démonstration.
\end{proof}

\begin{rema} Sous les hypothèses du lemme précédent, $A$ n'est pas nécessairement semi-primaire, comme le montre l'exemple de
$$A:=k[x_1,\dots,x_n,\dots]/(x_1^p,\dots,x_n^p,\dots)\,,$$
où $k$ est un corps fini de caractéristique $p$.
\end{rema}

\begin{prop}\label{cr-poids-faibles-torsion} Supposons que $F$ est un foncteur de $\Tt_A^\cf(\A,\E)$ dont les poids partiels sont périodiques. Alors il existe $r\in\mathbb{N}$ tel que, en posant $N:=p^r$, le sous-ensemble $\{x^N\,|\,x\in \AF\}$ de l'anneau $\AF$ soit fini.
\end{prop}

\begin{proof} Le lemme~\ref{lm-finjppp} fournit $r\in\mathbb{N}$ tel que le morphisme multiplicatif canonique $\alpha : \AF\to K^{\pp(F)}$ vérifie : $\forall (x,y)\in\AF^2\qquad\alpha(x)=\alpha(y)\Rightarrow x^{p^r}=y^{p^r}$.

Il suffit donc de montrer que l'image de $\alpha$ est finie. Cela découle de la finitude de $\pp(F)$ et du fait que les poids partiels de $F$ sont périodiques, donc d'images finies dans $K$, d'où la proposition.
\end{proof}

\begin{proof}[Démonstration du théorème~\ref{th-pptor_gal} (pour $p>0$)] Combinant la proposition~\ref{cr-poids-faibles-torsion} au lemme~\ref{ann-puissances-finies}, on voit que $\AF$ est isomorphe au produit d'un anneau fini et d'un anneau quasi-parfait $A'$ dont le quotient par le radical est fini de caractéristique $p$ ou $1$. Le corollaire~\ref{cor-ht_nil} montre alors que $A'$ est semi-primaire. Le fait que les corps qui sont quotients de $A'$ s'injectent dans $K$ se déduit du lemme~\ref{lm-finjppp}. En effet, si $x$ est un élément de $A'\setminus\rj(A')$ dont la classe dans $A'/\rj(A')$ est d'ordre $n$ dans le groupe $(A'/\rj(A'))^{\times}$, alors $x$ est d'ordre $n.p^s$ pour un $s\in\mathbb{N}$ dans le groupe $(A')^{\times}$. Cela termine la démonstration.
\end{proof}

\begin{rema}\label{rq-poiptpAFfini} Sous les hypothèses du théorème~\ref{th-pptor_gal}, $\AF$ n'est pas nécessairement fini, même si l'on suppose de plus que $F$ est fini, polynomial et à valeurs de dimensions finies.
\end{rema}

\begin{exem}\label{ex-poiptpAFfini} Supposons que $K$ est infini et que $A$ est le quotient d'un anneau de polynômes sur $\FF_p$ en une infinité dénombrable d'indéterminées $(x_i)$ par le carré de son idéal d'augmentation. Notons $\varphi : A\to K$ le morphisme d'anneaux composé de la projection $A\twoheadrightarrow\FF_p$ et de l'inclusion $\FF_p\hookrightarrow K$. Ainsi, $\varphi$ définit un $A_K$-module simple $K_\varphi$, de dimension $1$ sur $K$. Pour toute famille $(\xi_i)$ d'éléments de $K$, l'application $\FF_p$-linéaire $\mathrm{d} : A\to K$ envoyant $1$ sur $0$ et $x_i$ sur $\xi_i$ est une $\varphi$-dérivation\index{termin}{derivation@dérivation} $A\to K$. Ainsi, $\mathrm{d}$ définit un $A_K$-module $M_\mathrm{d}$ de dimension $2$ sur $K$, extension de $K_\varphi$ par lui-même (cf. \eqref{eq-der_Ext}, page~\pageref{eq-der_Ext}). Mais si la famille $(\xi_i)$ est libre sur $\FF_p$ --- ce qui est possible puisque $K$ est infini --- alors le morphisme d'anneaux (non nécessairement commutatif au but) structural $A\to\mathrm{End}_{A_K}(M_\mathrm{d})$ est injectif.

Par conséquent, le foncteur additif $F:=M_\mathrm{d}\underset{A}{\otimes}-$ de $\F(A,K)$ vérifie les propriétés suivantes :
\begin{enumerate}
\item $\rf(F)=0$, donc $\AF=A$ est infini ;
\item $F$ est polynomial, fini, à valeurs de dimensions finies ;
\item $F$ appartient à $\Tt^\cf(A,K)$ avec pour unique poids partiel $\varphi$, qui est de torsion car d'image finie.
\end{enumerate}
\end{exem}

\subsection{Cas de la caractéristique nulle}\label{pp0}

\begin{hyp} Dans tout le §\,\ref{pp0}, on suppose que $K$ est de caractéristique nulle.
\end{hyp}

\begin{lemm}\label{lm-icZ} Pour tout foncteur $F$ de $\fct(\A,\E)$, soit l'idéal caractéristique de $F$ relativement à $\mathbb{Z}$ est nul, soit $\ac_\mathbb{Z}(F)$ est de dimension finie sur $K$.
\end{lemm}

\begin{proof} En effet, $\ic_\mathbb{Z}(F)$ est un bi-idéal de $K[\mathbb{Z}]$ (proposition~\ref{pr-idc}), donc en particulier un idéal de l'algèbre $K[\mathbb{Z}_\mathrm{add}]\simeq K[t,t^{-1}]$, dont tous les quotients stricts sont de dimension finie.
\end{proof}

\begin{lemm}\label{lm-chiessai} Soient $F$ un foncteur de $\fct(\A,\E)$, $(x_i)_{i\in I}$ une famille finie d'éléments de $A$ et, pour tout $i\in I$, $(a_{i,j})_{j\in J(i)}$ une famille d'éléments de $A$ telle que $\chi_F(x_i)\in\operatorname{Vect}_K\big\{\chi_F(a_{i,j})\,|\,j\in J(i)\big\}$. Alors pour tout $t\in A$, on a
$$\chi_F\biggl(t+\prod_{i\in I}x_i\biggl)\in\operatorname{Vect}_K\left\{\chi_F\biggl(t+\prod_{i\in I}a_{i,\varphi(i)}\biggl)\,\Biggl|\,\varphi\in\prod_{i\in I}J(i)\right\}\;.$$
\end{lemm}

\begin{proof}
Pour $t=0$, cela découle de la multiplicativité de $\chi_F$. Comme le noyau $\ic_A(F)$ du morphisme $K$-linéaire $K[A]\to\ac_A(F)$ est un $\star$-idéal de $K[A]$ (proposition~\ref{pr-idc}), on en déduit le cas général.
\end{proof}

\begin{lemm}\label{lm-algent} Supposons que $A$ est un anneau de type fini et que $F$ est un foncteur de $\fct(\A,\E)$ tel que l'algèbre $\ac_A(F)$ soit entière sur $K$. Alors  $\ac_A(F)$ est de dimension finie sur $K$.
\end{lemm}

\begin{proof} Le lemme~\ref{lm-icZ} montre déjà que $\ac_\mathbb{Z}(F)$ est de dimension finie : il existe une partie finie $E$ de $\mathbb{N}$ telle que les $\chi_F(i)$ pour $i\in E$ engendrent linéairement $\ac_\mathbb{Z}(F)$.

Considérons par ailleurs des éléments $a_1,\dots,a_n$ de $A$ engendrant cet anneau. Par hypothèse, pour tout $i$, il existe un polynôme unitaire, disons de degré $d_i$, de $K[X]$ qui annule $\chi_F(a_i)\in\ac_A(F)$. En utilisant de façon itérée le lemme~\ref{lm-chiessai}, on en déduit que les éléments
$$\chi_F\left(\sum_{0\le s_i<d_i} r(s_1,\dots,s_n)\prod_{i=1}^n a_i^{s_i}\right)$$
de $\ac_A(F)$ l'engendrent comme $K$-espace vectoriel lorsque $r$ parcourt l'ensemble fini des fonctions de $\prod_{i=1}^n\llbracket 0,d_i-1\rrbracket$ dans $E$.
\end{proof}

\begin{lemm}\label{lm-cardanf} Si $A$ est un anneau fini, alors $\cd(A)\le\cd(A^{\mathrm{pereg}})^2$.\index{nota}{Apereg@$A^\mathrm{pereg}$}
\end{lemm}

\begin{proof} Si l'anneau fini $A$ est local, alors $\cd(A^{\mathrm{pereg}})\ge 2$ car $0$ et $1$ sont des éléments distincts de $A^{\mathrm{pereg}}$, et $\cd(A^{\mathrm{pereg}})\ge\cd(A)/2$ car, si $x$ est un élément non inversible de $A$, alors $1-x\in A^\times\subset A^{\mathrm{pereg}}$. Il s'ensuit que $\cd(A)\le\cd(A^{\mathrm{pereg}})^2$.

Le cas général s'en déduit en décomposant $A$ en produit d'anneaux locaux.
\end{proof}

\begin{proof}[Démonstration du théorème~\ref{th-pptor_gal} (pour $p=0$)] Quitte à remplacer $A$ par $\AF$ et $\A$ par $\A/\rf(F)$, on peut supposer sans perte de généralité que $\rf(F)$ est nul.

Comme les poids partiels de $F$ sont en nombre fini, disons $n$, et périodiques, il existe $i\in\mathbb{N}^*$ tel que $\pi^{i+1}=\pi$ pour tout $\pi\in\pp(F)$. Cela signifie que chaque $\pi$ est à valeurs dans $\{0\}\cup\mu_i(K)\subset K$, où $\mu_i(K)$ désigne le sous-groupe de $K^\times$ des racines $i$-èmes de l'unité. Cela implique par ailleurs (cf. proposition~\ref{pr-poidspar_ac}.\ref{itppac3}) que, pour tout $x\in A$, $(\chi_F(x)^{i+1}-\chi_F(x))^j=0$ dans $\ac_A(F)$, où $j$ désigne l'indice de nilpotence du radical de cette algèbre semi-primaire. En particulier, l'algèbre $\ac_A(F)$, qui est engendrée par les $\chi_F(x)$, est entière sur $K$.

Supposons dans un premier temps l'anneau $A$ de type fini. Le lemme~\ref{lm-algent} permet d'en déduire la finitude de $\dim_K\ac_A(F)$. Il s'ensuit que $\Sp(F)$ est un foncteur de type fini et à valeurs de dimensions finies de $\F(A,K)$ (corollaire~\ref{cor-spchi}). Par conséquent, $\Sp(F)$ possède une décomposition à la Steinberg de type $PAP$ \cite[cor.~4.19]{DTV}. Comme les poids partiels de $F$ sont de torsion, cela entraîne que $\Sp(F)$ est antipolynomial (utiliser la proposition~\ref{pr-decpord} et la propriété d'unicité de la proposition~\ref{pr-dppol}). Ainsi, l'anneau $A=\AF$ est fini. En utilisant le fait que les poids partiels de $F$ sont à valeurs dans $\{0\}\cup\mu_i(K)\subset K$ et le lemme~\ref{lm-finjpp}, on obtient $\cd(A^{\mathrm{pereg}})\le (i+1)^n$, puis, par le lemme~\ref{lm-cardanf}, $\cd(A)\le (i+1)^{2n}$.

Traitons maintenant le cas général, où $A$ n'est pas nécessairement de type fini. Si $B$ est un sous-anneau de type fini de $A$, considérant $\A$ comme une catégorie $B$-linéaire, $F$ appartient à $\Tt_B^\cf(\A,\E)$, son ensemble de poids partiels n'a pas plus d'éléments que $\pp_A(F)$, et la relation $\pi^{i+1}=\pi$ vaut aussi pour $\pi$ dans $\pp_B(F)$. Le raisonnement précédent procure donc l'inégalité $\cd(B)\le (i+1)^{2n}$. Comme elle est vraie pour tout sous-anneau de type fini $B$ de $A$, cet anneau est lui-même fini de cardinal au plus $(i+1)^{2n}$.
\end{proof}

\section[Décomposition selon les poids partiels]{Décomposition des foncteurs \hts\ selon les poids partiels}\label{ssdp}

\begin{nota}
Dans tout le §\,\ref{ssdp}, $n$ désigne un entier naturel et $\mathfrak{M}_1,\dots,\mathfrak{M}_n$ des sous-monoïdes de $\mon(A_\mu,K_\mu)$.
\end{nota}

\begin{nota} Si $\mathfrak{M}$ est un sous-monoïde de $\mon(A_\mu,K_\mu)$, on désigne par $\Tt_{A,[\mathfrak{M}]}(\A,\E)$ (resp. $\Tt^\cf_{A,[\mathfrak{M}]}(\A,\E)$) la sous-catégorie pleine de $\Tt_A(\A,\E)$ (resp. $\Tt^\cf_A(\A,\E)$) constituée des foncteurs $F$ tels que $\pp(F)\subset\mathfrak{M}$.

On note simplement $\Tt_{[\mathfrak{M}_1\times\dots\times\mathfrak{M}_n]}(\A^n,\E)$ (resp. $\Tt^\cf_{[\mathfrak{M}_1\times\dots\times\mathfrak{M}_n]}(\A^n,\E)$) pour $\Tt_{A^n,[\mathfrak{M}_1\times\dots\times\mathfrak{M}_n]}(\A^n,\E)$ (resp. $\Tt^\cf_{A^n,[\mathfrak{M}_1\times\dots\times\mathfrak{M}_n]}(\A^n,\E)$), où $\mathfrak{M}_1\times\dots\times\mathfrak{M}_n$ est vu comme sous-monoïde de $\mon(A_\mu^n,K_\mu)\simeq\mon(A_\mu,K_\mu)^n$.
\end{nota}

Grâce à la proposition~\ref{pr-poibif}, on dispose d'équivalences canoniques
\begin{equation}\label{eq-decevttp}
\Tt_{[\mathfrak{M}_1\times\mathfrak{M}_2]}(\A^2,\E)\simeq\Tt_{[\mathfrak{M}_1]}(\A,\Tt_{[\mathfrak{M}_2]}(\A,\E))\simeq\Tt_{[\mathfrak{M}_2]}(\A,\Tt_{[\mathfrak{M}_1]}(\A,\E))
\end{equation}
qui envoient $\Tt^\cf_{[\mathfrak{M}_1\times\mathfrak{M}_2]}(\A^2,\E)$ dans $\Tt^\cf_{[\mathfrak{M}_1]}(\A,\Tt^\cf_{[\mathfrak{M}_2]}(\A,\E))$ et $\Tt^\cf_{[\mathfrak{M}_2]}(\A,\Tt^\cf_{[\mathfrak{M}_1]}(\A,\E))$.

\begin{defi}\label{def-sep} On dit que $\mathfrak{M}_1,\dots,\mathfrak{M}_n$ vérifient la condition de séparation des poids partiels $\mathrm{(Sep)}$ si les deux propriétés suivantes sont satisfaites :
\begin{enumerate}
\item le morphisme canonique $\mathfrak{M}_1\times\dots\times\mathfrak{M}_n\to\mon(A_\mu,K_\mu)$ est injectif ;
\item pour tout foncteur $F$ non nul de $\Tt^\cf_{[\mathfrak{M}_1\dots\mathfrak{M}_n]}(\A,\E)$ (où $\mathfrak{M}_1\dots\mathfrak{M}_n$ désigne le sous-monoïde de $\mon(A_\mu,K_\mu)$ engendré par les $\mathfrak{M}_i$, qui est isomorphe à leur produit vu l'hypothèse précédente), on a
$$\widetilde{\W}_n(F)\cap (\mathfrak{M}_1\times\dots\times\mathfrak{M}_n)\ne\varnothing\quad\text{ dans }\quad\mon(A_\mu,K_\mu)^n.$$
\end{enumerate}
\end{defi}

On a alors $\widetilde{\W}_n(F)\cap (\mathfrak{M}_1\times\dots\times\mathfrak{M}_n)\ne\varnothing$ pour tout foncteur $F$ non nul de $\Tt_{[\mathfrak{M}_1\dots\mathfrak{M}_n]}(\A,\E)$.

\begin{rema}\label{rq-sepoiev}
\begin{enumerate}
    \item Si $\mathfrak{M}'_1\subset\mathfrak{M}_1$ et $\mathfrak{M}'_1\subset\mathfrak{M}_1$ sont des sous-monoïdes de $\mon(A_\mu,K_\mu)$ et que $\mathfrak{M}_1$ et $\mathfrak{M}_2$ vérifient la condition (Sep), alors il en est de même pour $\mathfrak{M}'_1$ et $\mathfrak{M}'_2$.
    \item Si les $n$ sous-monoïdes $\mathfrak{M}_1,\dots,\mathfrak{M}_n$ vérifient la condition (Sep), alors, pour tout $i\in\llbracket 1,n-1\rrbracket$, les $2$ sous-monoïdes $\mathfrak{M}_1\dots\mathfrak{M}_i$ et $\mathfrak{M}_{i+1}\dots\mathfrak{M}_n$ la vérifient également.
\end{enumerate}
\end{rema}

\begin{theo}\label{th-decPPart} Supposons que le morphisme canonique $\mathfrak{M}_1\times\dots\times\mathfrak{M}_n\to\mon(A_\mu,K_\mu)$ est injectif.
\begin{enumerate}
\item\label{ithp1} La restriction à $\Tt_{[\mathfrak{M}_1\times\dots\times\mathfrak{M}_n]}(\A^n,\E)$ du foncteur $\fct(\A^n,\E)\to\fct(\A,\E)$ de précomposition par la diagonale itérée $\delta_n : \A\to\A^n$ est pleinement fidèle.
\item\label{ithp2} Le foncteur $\Tt_{[\mathfrak{M}_1\times\dots\times\mathfrak{M}_n]}(\A^n,\E)\to\Tt_{[\mathfrak{M}_1\dots\mathfrak{M}_n]}(\A,\E)$ qu'induit $\delta_n^*$ est une équivalence si et seulement si la condition $\mathrm{(Sep)}$ est vérifiée.
\item\label{ithp3} Si $X$ est un foncteur de $\Tt_{[\mathfrak{M}_1\times\dots\times\mathfrak{M}_n]}(\A^n,\E)$ et $Y$ un foncteur de $\Tt^\cf_{[\mathfrak{M}_1\times\dots\times\mathfrak{M}_n]}(\A^n,\E)$, alors le morphisme
\[\mathrm{Ext}^*_{\fct(\A^n,\E)}(X,Y)\to\mathrm{Ext}^*_{\fct(\A,\E)}(\delta_n^*X,\delta_n^*Y)\]
induit par $\delta_n^*$ est un isomorphisme.
\end{enumerate}
\end{theo}

\begin{proof} Établissons d'abord \ref{ithp1} et \ref{ithp3}. Soient $X$ et $Y$ des foncteurs de $\Tt_{[\mathfrak{M}_1\times\dots\times\mathfrak{M}_n]}(\A^n,\E)$. L'adjonction somme/diagonale fournit un isomorphisme naturel $\mathrm{Ext}^*_{\fct(\A,\E)}(\delta_n^*X,\delta_n^*Y)\simeq\mathrm{Ext}^*_{\fct(\A^n,\E)}(X,\oplus_n^*\delta_n^*Y)$, il s'agit donc de montrer que le morphisme canonique $Y\to\oplus_n^*\delta_n^*Y$ induit un isomorphisme $\mathrm{Ext}^*_{\fct(\A^n,\E)}(X,Y)\xrightarrow{\simeq}\mathrm{Ext}^*_{\fct(\A^n,\E)}(X,\oplus_n^*\delta_n^*Y)$ lorsque le degré cohomologique est nul ou que $Y$ appartient à $\Tt^\cf_{[\mathfrak{M}_1\times\dots\times\mathfrak{M}_n]}(\A^n,\E)$. Comme $X$ et $\oplus_n^*\delta_n^*Y$ possèdent des décompositions en poids faibles, il suffit pour cela, d'après la proposition~\ref{pr-etr}, de vérifier que le morphisme canonique $Y\to\oplus_n^*\delta_n^*Y$ s'identifie à l'injection dans $\oplus_n^*\delta_n^*Y$ de la somme directe de ses parties homogènes de poids faibles appartenant à $\mathfrak{M}_1\times\dots\times\mathfrak{M}_n$.

Par hypothèse, $Y\circ\delta_n^{\A^n}$ (où l'exposant précise la catégorie additive sur laquelle où considère la diagonale $n$-itérée) possède une décomposition en poids faible relativement à l'action de $(A_\mu^n)^n\simeq A_\mu^{n^2}$ dont les composantes $\pi_{i,j}\in\mon(A_\mu,K_\mu)$ (pour $1\le i,j\le n$), où $\pi_{i,j}$ correspond à l'action de $A_\mu$ sur le $j$-ième terme dans la $i$-ème variable du $n$-multifoncteur $Y$, vérifient $\pi_{i,j}\in\mathfrak{M}_i$.

Considérant le diagramme commutatif (à isomorphisme près)
\[\xymatrix{\A^n\ar[r]^-{\oplus_n^\A}\ar[d]_-{\delta_n^{\A^n}} & \A\ar[d]^-{\delta_n^\A} \\
(\A^n)^n\ar[r]^-{\oplus_n^{\A^n}} & \A^n
}\,,\]
on en déduit que les poids de $\oplus_n^*\delta_n^*Y$ sont de la forme $(\prod_{i=1}^n\pi_{i,j})_{1\le j\le n}$ avec $\pi_{i,j}\in\mathfrak{M}_i$. Vu l'hypothèse faite sur les $\mathfrak{M}_j$, on a $\prod_{i=1}^n\pi_{i,j}\in\mathfrak{M}_j$ si et seulement si $\pi_{i,j}=1$ pour $i\ne j$. Comme la partie de poids $1$ d'un foncteur défini sur $\A$ est exactement son terme constant, cela termine la démonstration de \ref{ithp1} et \ref{ithp3}.

Montrons maintenant \ref{ithp2}. Tout d'abord, par \ref{ithp1}, $\delta_n^*$ induit un foncteur pleinement fidèle $\Phi : \Tt_{[\mathfrak{M}_1\times\dots\times\mathfrak{M}_n]}(\A^n,\E)\to\Tt_{[\mathfrak{M}_1\dots\mathfrak{M}_n]}(\A,\E)$. Utilisant la proposition~\ref{poids-evident} et l'adjonction somme/diagonale, on voit que $\Phi$ possède un adjoint à droite $\Psi$ donné par
\[F\mapsto\bigoplus_{w\in\mathfrak{M}_1\times\dots\times\mathfrak{M}_n}(\oplus_n^*F)^{<w>}.\]
Comme $\Phi$ est pleinement fidèle, ce foncteur est une équivalence si et seulement si $\Psi$ est fidèle. Or $\Psi$ est exact, car la restriction à la sous-catégorie pleine des foncteurs admettant une décomposition en poids faible d'un foncteur de la forme $X\mapsto X^{<w>}$ est exacte (puisque facteur direct de l'identité). Ainsi, $\Psi$ est fidèle si et seulement si la relation $\Psi(F)=0$ entraîne $F=0$. Comme $\Psi(F)=0$ équivaut, par définition, à $\widetilde{\W}_n(F)\cap (\mathfrak{M}_1\times\dots\times\mathfrak{M}_n)=\varnothing$, on en déduit \ref{ithp2}, ce qui achève la démonstration.
\end{proof}

\begin{rema} L'énoncé et une partie de la démonstration du théorème précédent sont très analogues à ceux de \cite[th.~4.4]{DT-ext}. Par ailleurs, si l'on se restreint aux foncteurs \hds, le théorème~\ref{th-decPPart} recoupe largement les propositions~\ref{pr-strPBT} et~\ref{pr-strglobsp}.
\end{rema}

L'énoncé suivant constitue, compte-tenu du théorème~\ref{th-decPPart}, une variation du corollaire~\ref{cor-sp_pt}.

\begin{prop}\label{pr-sepepi} Soient $E$ un foncteur exponentiel de $\F(\A;K)$ (c'est-à-dire un foncteur {\em monoïdal symétrique fort} de $(\A,\oplus,0)$ vers $(K\Md,\otimes,K)$ --- cf. \cite{T-expo}) et $\rho_i : E\twoheadrightarrow X_i$ des épimorphismes, pour $1\le i\le n$. On suppose que $X_i$ appartient à $\Tt_{[\mathfrak{M}_i]}(\A;K)$ pour chaque $i$ et que la condition $\mathrm{(Sep)}$ est vérifiée. Alors le morphisme
$$f : E\to E^{\otimes n}\xrightarrow{\rho_1\otimes\dots\otimes\rho_n}X_1\otimes\dots\otimes X_n$$
(où la première flèche est le coproduit itéré déduit de la structure exponentielle de $E$) est un épimorphisme.
\end{prop}

\begin{proof}
Notons $Y:=\mathrm{Coker}\,f$ : le théorème~\ref{th-decPPart} montre l'existence d'un foncteur $T$ de $\F(\A^n;K)$ et d'un épimorphisme $\varphi : X_1\boxtimes\dots\boxtimes X_n\twoheadrightarrow T$ tels que $Y\simeq\delta_n^*T$ et que la projection $X_1\otimes\dots\otimes X_n\twoheadrightarrow Y$ s'identifie à $\delta_n^*(\varphi)$. 

Comme l'adjoint de $f$ (pour l'adjonction entre $\oplus_n^*$ et $\delta_n^*$) est $\rho_1\boxtimes\dots\boxtimes\rho_n : E^{\boxtimes n}\twoheadrightarrow X_1\boxtimes\dots\boxtimes X_n$, l'adjoint de la composée $E\xrightarrow{f}X_1\otimes\dots\otimes X_n\twoheadrightarrow Y\simeq\delta_n^*T$, qui est nulle, est l'épimorphisme $$E^{\boxtimes n}\xrightarrow{\rho_1\boxtimes\dots\boxtimes\rho_n}X_1\boxtimes\dots\boxtimes X_n\xrightarrow{\varphi}T,$$
d'où la nullité de $T$ et donc de $Y$.
\end{proof}

\section{Théorèmes de structure}\label{ss-tshtr}

\subsection{Premier théorème de décomposition}

\begin{nota} Si $w=(w_1,\dots,w_n)\in\mon(A_\mu,K_\mu)^n$ est une suite finie de poids, on note $\ell(w)$ le nombre d'indices $i\in\llbracket 1,n\rrbracket$ tels que $w_i$ soit un élément sans torsion du monoïde $\mon(A_\mu,K_\mu)$.
\end{nota}

\begin{prop}\label{pr-pphtrmax} Soit $F$ un foncteur non nul de $\Tt^\cf_A(\A,\E)$.
\begin{enumerate}
\item\label{imev} Le sous-ensemble $\{\ell(w)\,|\,w\in\bigsqcup_{n\in\mathbb{N}}\W_n(F)\}$ de $\mathbb{N}$ possède un plus grand élément $s$.
\item\label{imp} Si $w=(\alpha_1,\dots,\alpha_r,\beta_1,\dots,\beta_s)$ est un multipoids de $F$ tel que les $\beta_i$ soient sans torsion, alors les poids $\alpha_i$ sont ordinaires.
\item\label{imp2} Il existe un multipoids $w=(\alpha_1,\dots,\alpha_r,\beta_1,\dots,\beta_s)$ de $F$ tel  que les $\beta_i$ soient sans torsion et que si $w'=(\alpha'_1,\dots,\alpha'_t,\beta'_1,\dots,\beta'_s)$ est un autre multipoids de $F$ avec les $\beta'_i$ sans torsion et $\la(\beta'_i)\subset\la(\beta_i)$ pour tout $i\in\llbracket 1,s\rrbracket$, alors $\la(\beta'_i)=\la(\beta_i)$ pour tout $i\in\llbracket 1,s\rrbracket$.
\item\label{impf} Si $w$ est un multipoids comme dans \ref{imp2}, alors les $\beta_i$ sont des morphismes d'anneaux et le multifoncteur $cr_{r+s}(F)^{<w>}$ est additif par rapport à chacune des $s$ dernières variables.
\end{enumerate}
\end{prop}

\begin{proof}
La première assertion résulte du lemme~\ref{lm-ppst}.

Comme les ensembles de multipoids sont invariants par permutation des coordonnées, on peut se contenter, afin d'alléger l'écriture, de montrer que $\alpha_1$ est un poids ordinaire, dans la situation \ref{imp}. Soient $a_2,\dots,a_r,b_1,\dots,b_s$ des objets de $\A$ tels que le foncteur $G:=cr_{r+s}(F)^{<w>}(-,a_2,\dots,a_r,b_1,\dots,b_s)$ de $\fct(\A,\E)$ soit non nul ; il est donc homogène de poids faible $\alpha_1$ et \htr. De plus, si $(\gamma_1,\dots,\gamma_n)$ est un multipoids de $G$, alors $(\gamma_1,\dots,\gamma_n,\alpha_2,\dots,\alpha_r,\beta_1,\dots,\beta_s)$ est un multipoids de $F$. La maximalité de $s$ montre que les $\gamma_i$ sont nécessairement tous périodiques. Par conséquent, le corollaire~\ref{cor-carppoiord} montre que $\alpha_1$ est un poids de torsion ordinaire.

L'assertion~\ref{imp2} résulte de \ref{imev} et de la finitude de $\pp(F)$. Pour montrer l'assertion~\ref{impf}, il suffit de montrer que $cr_{r+s}(F)^{<w>}$ est additif par rapport à la $(r+1)$-ième variable (si $s>0$), ce qui entraîne automatiquement que $\beta_1$ est un morphisme d'anneaux. Supposons le contraire : il existe alors des objets $a_1,\dots,a_r,b_2,\dots,b_s$ de $\A$ tels que le foncteur $H:=cr_{r+s}(F)^{<w>}(a_1,\dots,a_r,-,b_2,\dots,b_s)$ de $\fct(\A,\E)$ soit non additif. Comme $H$ est \htr, cela entraîne que $\W_2(H)$ est non vide. Soit $(\gamma,\gamma')$ un de ses éléments : comme $H$ est homogène de poids faible $\beta_1$, on a $\beta_1=\gamma\gamma'$. Par ailleurs, $(\alpha_1,\dots,\alpha_r,\gamma,\gamma',\beta_2,\dots,\beta_s)$ est un multipoids de $F$. La maximalité de $s$ montre que $\gamma$ et $\gamma'$ ne peuvent pas être tous deux sans torsion. Ils ne peuvent pas non plus être tous deux périodiques, car sinon leur produit $\beta_1$ le serait aussi. On peut donc supposer par exemple que $\gamma$ est périodique et $\gamma'$ sans torsion. La relation $\beta_1=\gamma\gamma'$ entraîne $\la(\gamma')\subset\la(\beta_1)$, de sorte que l'hypothèse de minimalité sur les $\la(\beta_i)$ fournit $\la(\gamma')=\la(\beta_1)$, puis $\la(\gamma)\subset\la(\gamma')=\la(\beta_1)$. Par ailleurs, $\gamma$ est un poids non trivial de torsion ordinaire par \ref{imp}, donc $\la(\gamma)$ contient un idéal maximal cofini $\mathfrak{m}$ de $A$, grâce aux propositions \ref{pr-zerantipol}, \ref{pr-decpord} et au corollaire~\ref{cor-lappol}. Comme $\pp(H)$ est fini, il s'ensuit qu'il existe une intersection finie $I$ d'idéaux maximaux cofinis de $A$ telle que $I\subset\la(\pi)$ pour tout $\pi\in\pp(H)$. Le corollaire~\ref{cor-kerpp} montre alors que les poids partiels de $H$ se factorisent par la projection $A\twoheadrightarrow A/I^n$ pour un certain $n\in\mathbb{N}$. L'anneau $A/I^n$ est semi-primaire et le quotient par son radical est fini, ce qui implique que tout poids de $\mon((A/I^n)_\mu,K_\mu)$ est de torsion.
% à ériger en lemme, au chapitre \ref{spfca} ?
% sinon, cette démonstration est un peu touffue...
(En effet, on peut se ramener au cas local ; si $A$ est un anneau local dont le radical $J$ est nilpotent et le corps résiduel $\kappa$ fini, $\mon(A_\mu,K_\mu)\simeq\mathbf{Ab}(A^\times,K^\times)_+$, et $\mathbf{Ab}(A^\times,K^\times)\simeq\mathbf{Ab}(\kappa^\times,K^\times)\oplus\mathbf{Ab}((1+J)_\mu,K^\times)$, or $(1+J)_\mu$ est un groupe annulé par $l^r$, où $l$ est la caractéristique de $\kappa$ et $r$ l'indice de nilpotence de $J$.) Comme $\gamma$ est sans torsion mais se factorise par la réduction modulo $I^n$, c'est absurde, d'où la proposition.
\end{proof}

\begin{nota}\label{not-tordpoli}
On note $\mathfrak{M}_{\mathrm{tord}}$ le sous-monoïde $\mon(A_\mu,K_\mu)$ constitué des poids ordinaires de torsion, et $\mathfrak{M}_{\mathrm{poli}}$ le sous-monoïde engendré par les morphismes d'anneaux $A\to K$ d'image infinie.
\end{nota}

L'indice tord (resp. poli) abrège {\em de Torsion ORDinaire} (resp. {\em POLynomial d'image Infinie}).

\begin{coro}\label{cor-ppord} Tous les poids partiels d'un foncteur \htr\ sont ordinaires. De plus, de tels poids partiels s'expriment comme produit d'un élément de $\mathfrak{M}_{\mathrm{tord}}$ et d'un élément de $\mathfrak{M}_{\mathrm{poli}}$, et ce de façon unique.
\end{coro}

\begin{proof} Tout poids partiel $\pi$ d'un foncteur \htr\ de $\fct(\A,\E)$ est le poids d'un foncteur $F$ de $\Tt_A^\cf(\A,\E)$ homogène de poids faible $\pi$. Appliquant les assertions~\ref{impf} et~\ref{imp} de la proposition~\ref{pr-pphtrmax}, on obtient un multipoids $(\alpha_1,\dots,\alpha_r,\beta_1,\dots,\beta_s)$ de $F$ tel quel les $\alpha_i$ soient des poids de torsion ordinaires et les $\beta_j$ des morphismes d'anneaux. Il s'ensuit que $\pi=\alpha_1\dots\alpha_r\beta_1\dots\beta_s$ est un poids ordinaire possédant une décomposition du type requis. L'unicité résulte des propositions~\ref{pr-decpord} et~\ref{pr-dppol}.
\end{proof}

\begin{coro}\label{cor-septorpol} Les sous-monoïdes $\mathfrak{M}_{\mathrm{tord}}$ et $\mathfrak{M}_{\mathrm{poli}}$ de $\mon(A_\mu,K_\mu)$ vérifient la condition de séparation des poids partiels $\mathrm{(Sep)}$ (cf. définition~\ref{def-sep}).
\end{coro}

\begin{proof}
Cela découle de l'unicité du corollaire précédent et des assertions~\ref{imp2} et~\ref{impf} de la proposition~\ref{pr-pphtrmax}.
\end{proof}

En combinant le théorème~\ref{th-decPPart}, les corollaires~\ref{cor-ppord} et~\ref{cor-septorpol}, et les équivalences \eqref{eq-decevttp} (page~\pageref{eq-decevttp}), on obtient :

\begin{theo}\label{th-decHTR} La précomposition par la diagonale $\A\to\A\times\A$ induit une équivalence de catégories
$$\Tt_{[\mathfrak{M}_{\mathrm{poli}}\times\mathfrak{M}_{\mathrm{tord}}]}(\A\times\A,\E)\simeq\Tt(\A,\E)\;,$$
d'où des équivalences
$$\Tt(\A,\E)\simeq\Tt_{[\mathfrak{M}_{\mathrm{poli}}]}(\A,\Tt_{[\mathfrak{M}_{\mathrm{tord}}]}(\A,\E))\simeq\Tt_{[\mathfrak{M}_{\mathrm{tord}}]}(\A,\Tt_{[\mathfrak{M}_{\mathrm{poli}}]}(\A,\E))\;.$$

Ces équivalences se restreignent en une équivalence
$$\Tt^\cf(\A,\E)\simeq\Tt^\cf_{[\mathfrak{M}_{\mathrm{poli}}\times\mathfrak{M}_{\mathrm{tord}}]}(\A\times\A,\E)$$
et envoient ces catégories respectivement dans
$\Tt^\cf_{[\mathfrak{M}_{\mathrm{poli}}]}(\A,\Tt^\cf_{[\mathfrak{M}_{\mathrm{tord}}]}(\A,\E))$ et \linebreak
$\Tt^\cf_{[\mathfrak{M}_{\mathrm{tord}}]}(\A,\Tt^\cf_{[\mathfrak{M}_{\mathrm{poli}}]}(\A,\E))$.
\end{theo}

\begin{coro}\label{cor-hpap} Tout foncteur de $\Tt^\cf(A,K)$ possède une décomposition à la Steinberg de type $HPAP$ (cf. définition~\ref{def-decSt}).

Si la condition \textnormal{(IMF)}\index{nota}{IMF@(IMF) \emph{(condition de finitude sur les idéaux maximaux à quotient fini)}} est vérifiée, tout foncteur de $\Tt^\cf(A,K)$ possède une décomposition de type $FP$ (cf. définition~\ref{def-decFP}).
\end{coro}

\begin{proof}
Ce résultat découle du théorème~\ref{th-decHTR}, de la proposition~\ref{pr-raispol} et des corollaires~\ref{cor-carhtper} et~\ref{cor-IMF}.
\end{proof}

\subsection{Propriétés de finitude des foncteurs \hts\ de $\F(A,K)$}

\begin{coro}\label{cor-htn_df} Tout foncteur \htr\ noethérien $F$ de $\F(A,K)$ appartient à $\F^\df(A,K)$.
\end{coro}

\begin{proof} Comme $F$ appartient à $\Tt^\cf(A,K)$, puisque $F$ est de type fini, le corollaire précédent montre qu'il existe un quotient $K$-trivial $A'$ de $A$ et un bifoncteur $B$ de $\F(\mathbf{P}(A)\times\mathbf{P}(A');K)$ tels que $F$ soit isomorphe à la composée de $B$ et du foncteur canonique $\phi : \mathbf{P}(A)\to\mathbf{P}(A)\times\mathbf{P}(A')$, avec $B$ \ph\ par rapport à la première variable et possédant une décomposition en poids. Comme $\phi^*B$ est noethérien, \cite[prop.~4.9]{DTV} montre que $B$ est noethérien. Le corollaire~\ref{cor-noethvn} montre alors que l'image de $B$ par l'équivalence de catégories canonique $\F(\mathbf{P}(A)\times\mathbf{P}(A');K)\simeq\fct(\mathbf{P}(A),\F(A',K))$ est à valeurs noethériennes. Or tout foncteur de type fini de $\F(A',K)$ est à valeurs de dimensions finies, puisque $A'$ est fini. Il s'ensuit que $F$ lui-même est à valeurs de dimensions finies.
\end{proof}

Le théorème~\ref{th-decHTR} permet de déterminer facilement quand tous les foncteurs \hts\ de $\F(A,K)$ sont localement noethériens :
\begin{theo}\label{th-htr_ln} Les propriétés suivantes sont équivalentes :
\begin{enumerate}
\item[(a)] la catégorie $\Tt(A,K)$ est localement noethérienne ;
\item[(b)] tout foncteur additif de $\Tt(A,K)$ est localement noethérien ;
\item[(c)] la condition \textnormal{(CFD)}\index{nota}{CFD@(CFD), (CFD)$^+$ \emph{(conditions de finitude sur les dérivations)}} est vérifiée ;
\item[(d)] tout foncteur de type fini de $\Tt(A,K)$ est à valeurs de dimensions finies.
\end{enumerate}
\end{theo}

\begin{proof} La proposition~\ref{pr-polhtfini} montre que les conditions (b) et (c) sont équivalentes et que (d), qui implique la condition \ref{itG4} de la proposition~\ref{pr-polhtfini}, entraîne (c).

Utilisant l'équivalence $\Tt(A,K)\simeq\Tt_{[\mathfrak{M}_\mathrm{poli}]}(\mathbf{P}(A),\Tt_{[\mathfrak{M}_\mathrm{tord}]}(A,K))$ donnée par le théorème~\ref{th-decHTR}, la proposition~\ref{pr-raispol}, qui montre que tout foncteur de type fini de cette catégorie est polynomial (dans $\fct(\mathbf{P}(A),\Tt_{[\mathfrak{M}_\mathrm{tord}]}(A,K))$) et la proposition~\ref{pr-polhtfini}, on voit que, sous la condition (CFD), tout foncteur de type fini de $\Tt(A,K)$ est à valeurs de type fini dans $\Tt_{[\mathfrak{M}_\mathrm{tord}]}(A,K)$ (via l'équivalence précédente). Comme la condition (CFD) entraîne (IMF) (remarque~\ref{rq-CFD_IMF}), le corollaire~\ref{cor-IMF} montre que, pour tout foncteur $F$ de type fini de $\Tt_{[\mathfrak{M}_\mathrm{tord}]}(A,K)$, l'anneau $\AF$ est fini. Cela implique que $F$ est à valeurs de dimensions finies, d'où l'implication (c)$\Rightarrow$(d), mais aussi que $F$ est noethérien par le théorème~\ref{th-PSS}. Comme un foncteur polynomial à valeurs noethériennes de source $\mathbf{P}(A)$ est noethérien \cite[lm~6.20]{DT-schw}, cela montre aussi l'implication (c)$\Rightarrow$(a).

Comme (a) entraîne trivialement (b), ceci achève la démonstration.
\end{proof}

\begin{rema} Le corollaire~\ref{cor-htn_df} donne une autre démonstration, plus directe, de l'implication (a)$\Rightarrow$(d) du théorème précédent.
\end{rema}

\subsection{Décomposition primaire des foncteurs \hts}

On peut préciser les résultats précédents en introduisant le raffinement suivant de la notation~\ref{not-tordpoli} :

\begin{nota}\label{nota-decprimpoi}
\begin{enumerate}
    \item Si $I$ est un idéal maximal cofini de $A$, on note $\mathfrak{M}_{\mathrm{tord},I}$ le sous-monoïde des poids $w$ de $\mathfrak{M}_{\mathrm{tord}}$ se factorisant à travers la projection $A\twoheadrightarrow A/I^n$ pour un certain $n\in\mathbb{N}$.
    \item Si $\alpha : A\to K$ est un morphisme d'anneaux d'image infinie, on note $\mathfrak{M}_{\mathrm{poli},\alpha}$ le sous-monoïde de $\mathfrak{M}_{\mathrm{poli}}$ engendré par $\alpha$ lorsque $p=0$, et le sous-monoïde constitué des $w\in\mathbf{Mon}(A_\mu,K_\mu)$ tels qu'il existe $(r,n)\in\mathbb{N}^2$ tel que $w^{p^r}=\alpha^n$ si $p>0$.
\end{enumerate}
\end{nota}

L'énoncé suivant découle des propositions \ref{pr-decpord}, \ref{pr-dppol} et de la décomposition d'un anneau fini en produit d'anneaux locaux. 

\begin{prop}\label{pr-indep_imF}
\begin{enumerate}
    \item Si $I_1,\dots,I_n$ sont des idéaux maximaux cofinis deux à deux distincts de $A$, alors les monoïdes $\mathfrak{M}_{\mathrm{tord},I_1},\dots,\mathfrak{M}_{\mathrm{tord},I_n}$ vérifient la condition de séparation des poids partiels $\mathrm{(Sep)}$.
    \item Un foncteur $F$ de $\fct(\A,\E)$ appartient à $\Tt^\cf_{[\mathfrak{M}_{\mathrm{tord}}]}(\A,\E)$ si et seulement s'il existe des idéaux maximaux cofinis deux à deux distincts $I_1,\dots,I_n$ de $A$ tels que $F$ appartienne à $\Tt^\cf_{[\mathfrak{M}_{\mathrm{tord},I_1}\dots\mathfrak{M}_{\mathrm{tord},I_n}]}(\A,\E)$.
\end{enumerate}
\end{prop}

Le résultat suivant se déduit quant à lui de la proposition~\ref{pr-dppol} et des assertions~\ref{imp2} et~\ref{impf} de la proposition~\ref{pr-pphtrmax}.

\begin{prop}\label{pr-indep_alpha}
\begin{enumerate}
    \item Si $(\alpha_1,\dots,\alpha_n)$ est une famille indépendante de morphismes d'anneaux $A\to K$, alors les monoïdes $\mathfrak{M}_{\mathrm{poli},\alpha_1},\dots,\mathfrak{M}_{\mathrm{poli},\alpha_n}$ vérifient la condition de séparation des poids partiels $\mathrm{(Sep)}$.
    \item Un foncteur $F$ de $\fct(\A,\E)$ appartient à $\Tt^\cf_{[\mathfrak{M}_{\mathrm{poli}}]}(\A,\E)$ si et seulement s'il existe une famille indépendante $(\alpha_1,\dots,\alpha_n)$ de $\mathbf{Ann}(A,K)$ telle que $F$ appartienne à $\Tt^\cf_{[\mathfrak{M}_{\mathrm{poli},\alpha_1}\dots\mathfrak{M}_{\mathrm{poli},\alpha_n}]}(\A,\E)$.
\end{enumerate}
\end{prop}

Les deux propositions suivantes précisent la structure des foncteurs \hts\ dont les poids partiels appartiennent tous à l'une des \guillemotleft~briques élémentaires~\guillemotright\, introduites dans la notation~\ref{nota-decprimpoi}.

\begin{prop}\label{pr-tordI} Soient $I$ un idéal maximal cofini de $A$ et $F$ un foncteur de $\Tt^\cf_{[\mathfrak{M}_{\mathrm{tord},I}]}(\A,\E)$.
\begin{enumerate}
\item Il existe $r\in\mathbb{N}$ tel que $I^r\subset\rf(F)$.
\item Si $F$ est non constant, alors le cardinal $q$ du corps fini $A/I$ vérifie $q-1\in\nu(K)$.
\end{enumerate}
\end{prop}

\begin{proof}
La première assertion provient du corollaire~\ref{cor-kerpp}, et la deuxième se déduit de la proposition~\ref{pr-htr_AFfini}.
\end{proof}

\begin{prop}\label{pr-poli_alpha} Soient $\alpha : A\to K$ un morphisme d'anneaux d'image infinie, $\mathfrak{p}$ son noyau et $F$ un foncteur de $\Tt^\cf_{[\mathfrak{M}_{\mathrm{poli},\alpha}]}(\A,\E)$.
\begin{enumerate}
\item Il existe $r\in\mathbb{N}$ tel que $\mathfrak{p}^r\subset\rf(F)$.
\item Il existe un foncteur de $G$ de $\Tt^\cf_{[\mathfrak{M}_{\mathrm{poli},\alpha}]}(\A\otimes_A (A_\mathfrak{p}/\mathfrak{p}^r A_\mathfrak{p}),\E)$, unique à isomorphisme près, tel que $F$ soit la composée de $G$ et du foncteur additif canonique $\A\to\A\otimes_A (A_\mathfrak{p}/\mathfrak{p}^r A_\mathfrak{p})$.
\end{enumerate}
\end{prop}

\begin{proof} La première assertion découle du corollaire~\ref{cor-kerpp}. La deuxième s'en déduit en utilisant la proposition~\ref{pr-invad}.
\end{proof}

\begin{rema} Dans les conditions de la proposition~\ref{pr-poli_alpha}, si $F$ est non constant, l'anneau $\tilde{A}_F$\index{nota}{A@$\tilde{A}_F$} (avec la notation de la proposition~\ref{pr-invad}) est local, semi-primaire, et son corps résiduel $\mathrm{Frac}(A/\mathfrak{p})$ est un sous-corps infini de $K$. 

De façon analogue, dans les conditions de la proposition~\ref{pr-tordI}, si $F$ est non constant, l'anneau $\AF$ est local, semi-primaire ; son corps résiduel $A/I$ est fini, et s'il est de caractéristique $p$, il se plonge dans $K$.
\end{rema}

On peut résumer une grande partie des résultats précédents les plus importants en l'énoncé suivant :

\begin{theo}\label{th-struct_htr}
\begin{enumerate}
\item\label{itpfstr} Si $(\alpha_1,\dots,\alpha_n)$ est une famille finie indépendante de $\mathbf{Ann}(A,K)$ et que $I_1,\dots,I_m$ sont des idéaux maximaux cofinis deux à deux distincts de $A$, alors le foncteur
$$\Tt_{[\mathfrak{M}_{\mathrm{poli},\alpha_1}\times...\times\mathfrak{M}_{\mathrm{poli},\alpha_n}\times\mathfrak{M}_{\mathrm{tord},I_1}\times\dots\times\mathfrak{M}_{\mathrm{tord},I_m}]}(\A^{n+m},\E)\to\fct(\A,\E)$$
restriction de $\delta_{n+m}^*$ est pleinement fidèle et a pour image essentielle\linebreak $\Tt_{[\mathfrak{M}_{\mathrm{poli},\alpha_1}\dots\mathfrak{M}_{\mathrm{poli},\alpha_n}.\mathfrak{M}_{\mathrm{tord},I_1}\dots\mathfrak{M}_{\mathrm{tord},I_m}]}(\A,\E)$.
\item\label{itexist} Si $F$ est un foncteur de $\Tt^\cf(\A,E)$, il existe un unique $(n,m)\in\mathbb{N}^2$, une famille $(\alpha_1,\dots,\alpha_n)$ indépendante de $\mathbf{Ann}(A,K)$, unique à l'ordre, et à l'action du morphisme de Frobenius si $p>0$, près, une famille $(I_1,\dots,I_m)$ d'idéaux maximaux cofinis deux à deux distincts de $A$, unique à l'ordre près, et un foncteur $X$ de $\Tt^\cf_{[\mathfrak{M}_{\mathrm{poli},\alpha_1}\times\dots\times\mathfrak{M}_{\mathrm{poli},\alpha_n}\times\mathfrak{M}_{\mathrm{tord},I_1}\times\dots\times\mathfrak{M}_{\mathrm{tord},I_m}]}(\A^{n+m},\E)$ qui n'est constant par rapport à aucune des $n+m$ variables, unique à isomorphisme près, tels que $F\simeq\delta_{n+m}^*(X)$.
\item Ce multifoncteur $X$ vérifie de plus les propriétés suivantes :
\begin{enumerate}
\item\label{it_sup1} $X$ est polynomial par rapport aux $n$ premières variables ;
\item\label{it_sup2} pour $j\in\llbracket 1,m\rrbracket$, si $A/I_j$ est de caractéristique $p$, et que la condition \textnormal{(MTF)}\index{nota}{MTF@(MTF) \emph{(condition de finitude sur les morphismes)}} est vérifiée, alors $X$ est \ph\ par rapport à la $(n+j)$-ième variable ;
\item\label{it_sup3} pour $j\in\llbracket 1,m\rrbracket$, si $A/I_j$ est de caractéristique différente de $p$, et que la condition \textnormal{(MTF)} est vérifiée, alors $X$ est antipolynomial par rapport à la $(n+j)$-ième variable.
\end{enumerate}
\end{enumerate}
\end{theo}

\begin{proof} L'assertion~\ref{itpfstr} et l'énoncé d'existence de l'assertion~\ref{itexist} découlent des théorèmes \ref{th-decHTR} et \ref{th-decPPart} ainsi que des propositions \ref{pr-indep_imF} et \ref{pr-indep_alpha} (et de la remarque~\ref{rq-sepoiev}). \textit{La famille indépendante $(\alpha_1,\dots,\alpha_n)$ de $\mathbf{Ann}(A,K)$ et la famille $(I_1,\dots,I_m)$ d'idéaux maximaux cofinis deux à deux distincts de $A$ étant fixées}, l'unicité à isomorphisme près de $X$ dans $\Tt^\cf_{[\mathfrak{M}_{\mathrm{poli},\alpha_1}\times\dots\times\mathfrak{M}_{\mathrm{poli},\alpha_n}\times\mathfrak{M}_{\mathrm{tord},I_1}\times\dots\times\mathfrak{M}_{\mathrm{tord},I_m}]}(\A^{n+m},\E)$ tel que $F\simeq\delta_{n+m}^*(X)$ résulte de l'assertion~\ref{itpfstr}. La propriété d'unicité plus forte énoncée en~\ref{itexist} s'en déduit en regroupant les familles d'idéaux maximaux cofinis et les familles indépendantes de morphismes d'anneaux apparaissant dans deux décompositions distinctes et en en supprimant les doublons (y compris au Frobenius près pour les morphismes d'anneaux, dans le cas où $p>0$) afin d'obtenir une famille d'idéaux maximaux cofinis englobant celles des deux décompositions et une famille indépendante de morphismes d'anneaux englobant (à l'action du Frobenius près si $p>0$) celles des deux décompositions.

L'assertion~\ref{it_sup1} résulte de la proposition~\ref{pr-raispol}, l'assertion~\ref{it_sup2} des propositions \ref{pr-htr_impl_ph} (qui s'étend de $\Tt^\cf(A,K)$ à $\Tt^\cf(\A,\E)$ sous l'hypothèse (MTF) grâce au corollaire~\ref{cor-phsp}) et \ref{pr-tordI}, et l'assertion~\ref{it_sup3} du théorème~\ref{th-pptor_gal}.
\end{proof}

\begin{rema}
Les idéaux premiers $\mathrm{Ker}\,\alpha_i$ et les idéaux maximaux cofinis $I_j$ apparaissant dans l'énoncé précédent peuvent être caractérisés comme étant exactement les idéaux premiers qui sont le lieu d'annulation d'un poids partiel de $F$. C'est pourquoi on peut penser à l'isomorphisme du théorème~\ref{th-struct_htr} comme à une \textit{décomposition primaire} pour les foncteurs \hts\ à caractère fini.
\end{rema}

\subsection{Spectre des foncteurs \hts}

Afin d'exprimer en termes de spectre les résultats précédents, nous introduisons la notation suivante (on rappelle que, si $F$ est un foncteur de $\F(A,K)$, $w\in\mon(A_\mu,K_\mu)$ un poids et $r$ un entier naturel, $F^{[r]}_w$ désigne le plus grand quotient de $F$ appartenant à $U(\F(A,K);\mathfrak{m}_w^r)$ --- cf. définition~\ref{poids-evident}.\ref{itprdfr}).

\begin{nota}\label{not-Tspc}Pour $d,r,n\in\mathbb{N}$, $w_1,\dots,w_n : A\to K$ des poids et $I$ un idéal de $A$, on définit dans $\F(A,K)$
$$\mathrm{T}_{A,I,d,r,(w_1,\dots,w_n)}:=\pi_I^*P^{A/I}_{\mathbf{P}(A/I)}\otimes\Big(\bigoplus_{i=1}^n\qpol_d(P^A)_{w_i}^{[r]}\Big)\,,$$
où $\pi_I : \mathbf{P}(A)\to\mathbf{P}(A/I)$ désigne la réduction modulo $I$.
\end{nota}

\begin{prop}\label{pr-spht} Soient $d,r,n\in\mathbb{N}$, $w_1,\dots,w_n$ des éléments deux à deux distincts de $\mathfrak{M}_{\mathrm{poli}}$ et $I$ un idéal cofini de $A$.

Alors le morphisme $P^A\to\mathrm{T}_{A,I,d,r,(w_1,\dots,w_n)}$
de $\F(A,K)$ composé de la diagonale $P^A\to P^A\otimes P^A$ et du morphisme canonique $P^A\to\pi_I^*P^{A/I}_{\mathbf{P}(A/I)}$ tensorisé par le morphisme $P^A\to\bigoplus_{i=1}^n\qpol_d(P^A)_{w_i}^{[r]}$ dont les composantes sont les projections est un épimorphisme.
\end{prop}

\begin{proof} Quitte à agrandir $K$, on peut supposer ce corps algébriquement clos.

Il existe un sous-anneau $A'$ \textit{de type fini} de $A$ tel que :
\begin{enumerate}
\item les restrictions à $A'$ des $w_i$ sont deux à deux distinctes ; 
\item chaque $w_i$ est produit de morphismes d'anneaux $A\to K$ dont la restriction à $A'$ a une image de cardinal strictement supérieur à celui de $A/I$.
\end{enumerate}

En vertu de  la proposition~\ref{pr-polpoihtrp}, tout foncteur polynomial de $\F(A,K)$ possédant une décomposition en poids est \htr \emph{ relativement à l'action de $A'$}.

Soit $\mathfrak{M}_1$ le sous-monoïde de $\mon(A'_\mu,K_\mu)$ constitué des morphismes se factorisant par la réduction modulo $I\cap A'$, et soit $\mathfrak{M}_2$ le sous-monoïde constitué engendré par des morphismes d'anneaux $A'\to K$ d'image de cardinal strictement supérieur à celui de $A/I$. Le corollaire~\ref{cor-septorpol}, la proposition~\ref{pr-indep_imF} et la remarque~\ref{rq-sepoiev} montrent que $\mathfrak{M}_1$ et $\mathfrak{M}_2$ vérifient la condition (Sep). En effet, si $I'_1,\dots,I'_r$ désignent les idéaux maximaux (nécessairement cofinis) de $A'$ contenant $I\cap A'$, on a $\mathfrak{M}_1\subset\mathfrak{M}_{\mathrm{tord},I'_1}\dots\mathfrak{M}_{\mathrm{tord},I'_r}$ et $\cd(A'/I'_i)\le\cd(A/I)$, tandis que $\mathfrak{M}_2$ est inclus dans le sous-monoïde de $\mon(A'_\mu,K_\mu)$ engendré par $\mathfrak{M}_{\mathrm{poli}}$ et par les $\mathfrak{M}_{\mathrm{tord},\mathfrak{m}}$, où $\mathfrak{m}$ est un idéal maximal cofini de $A'$ tel que $\cd(A'/\mathfrak{m})>\cd(A/I)$ et que $A'/\mathfrak{m}$ s'injecte dans $K$.

Comme le morphisme canonique $P^A\to\bigoplus_{i=1}^n\qpol_d(P^A)_{w_i}^{[r]}$ est un épimorphisme (cf. proposition~\ref{pr-sp_sdir} et exemple~\ref{ex-sdirpq}), la conclusion se déduit alors de la proposition~\ref{pr-sepepi}.
\end{proof}

Le corollaire suivant donne dans certains cas importants --- par exemple lorsque $K$ est de caractéristique nulle, ou que $K$ est un corps parfait et $A$ un anneau de type fini --- une caractérisation complète en termes de spectre des foncteurs \hts\ à caractère fini.

\begin{coro}\label{cor-hpapfor} Soit $F$ un foncteur de $\fct(\A,\E)$.
\begin{enumerate}
\item Si la condition \textnormal{(IMF)}\index{nota}{IMF@(IMF) \emph{(condition de finitude sur les idéaux maximaux à quotient fini)}} est vérifiée et que $F$ appartient à $\Tt^\cf_A(\A,\E)$, alors il existe un idéal cofini $I$ de $A$ tel que $(A/I)^\times$ soit somme directe de groupes cycliques d'ordres appartenant à $\nu(K)$, des entiers naturels $r$, $d$, $n$ et des éléments $w_1,\dots,w_n$ de $\mathfrak{M}_{\mathrm{poli}}$ deux à deux distincts tels que
$$\Sp(F)\le\mathrm{T}_{A,I,d,r,(w_1,\dots,w_n)}\,.$$
\item Supposons que l'une des deux conditions suivantes est vérifiée :
\begin{enumerate}
    \item\label{itccar1} $K$ est de caractéristique nulle ;
    \item\label{itccar2} $K$ est un corps parfait de caractéristique $p>0$ et $A/p$ est une algèbre essentiellement de type fini sur une $\FF_p$-algèbre $p$-parfaite (par exemple, $A$ un anneau essentiellement de type fini).
\end{enumerate}
Alors tout foncteur $F$ vérifiant la condition précédente appartient à $\Tt^\cf_A(\A,\E)$. 
\end{enumerate}
\end{coro}

(On notera que chacune des deux conditions \ref{itccar1} et \ref{itccar2} implique (IMF).)

\begin{proof}
Le corollaire~\ref{cor-IMF} et la proposition~\ref{pr-htr_AFfini} montrent que, si (IMF) est vérifiée, alors pour tout foncteur $F$ de $\Tt^\cf_A(\A,\E)$ tel que $\pp(F)\subset\mathfrak{M}_{\mathrm{tord}}$, il existe un idéal cofini $I$, tel que $(A/I)^\times$ soit somme directe de groupes cycliques d'ordres appartenant à $\nu(K)$ et que $\Sp(F)\le\pi_I^*P^{A/I}_{\mathbf{P}(A/I)}$.

Par ailleurs, si $F$ est un foncteur de $\Tt^\cf_A(\A,\E)$ tel que $\pp(F)\subset\mathfrak{M}_{\mathrm{poli}}$, alors $F$ possède en particulier une décomposition en poids faible finie, avec des poids $w_1,\dots,w_n$ dans $\mathfrak{M}_{\mathrm{poli}}$, donc il existe $r\in\mathbb{N}$ tel que $\Sp(F)\le\bigoplus_{i=1}^n\qpol_d(P^A)_{w_i}^{[r]}$. La première assertion résulte de ces deux cas particuliers, du théorème~\ref{th-decHTR} et du corollaire~\ref{cor-sp_pt}.

Sous l'hypothèse \ref{itccar1}, le foncteur $\bigoplus_{i=1}^n\qpol_d(P^A)_{w_i}^{[r]}$ appartient à $\Tt^\cf(A,K)$ par le corollaire~\ref{cor-polpoihtr}. C'est également le cas sous l'hypothèse \ref{itccar2}, car le corollaire~\ref{cor-poids-ln} montre que ce foncteur est alors fini et à valeurs de dimensions finies ; comme ses poids sont déployés et que $K$ est parfait, la proposition~\ref{pr-simplesPolHD} montre que ses facteurs de composition sont \hds, ce qui implique que notre foncteur est \htr.

Par ailleurs, si $I$ est un idéal cofini de $A$ tel que $(A/I)^\times$ soit somme directe de groupes cycliques d'ordres dans $\nu(K)$, alors $\pi_I^*P^{A/I}_{\mathbf{P}(A/I)}$ appartient à $\Tt^\cf(A,K)$ par la proposition~\ref{pr-htr_AFfini}. Il s'ensuit que, si l'une des hypothèses \ref{itccar1} ou \ref{itccar2} est remplie, le foncteur $\mathrm{T}_{A,I,d,r,(w_1,\dots,w_n)}$ est \htr. En conséquence, tout foncteur $F$ de $\fct(\A,\E)$ tel que $\Sp(F)\le\mathrm{T}_{A,I,d,r,(w_1,\dots,w_n)}$ appartient à $\Tt^\cf_A(\A,\E)$, ce qui achève la démonstration.
\end{proof}

\subsection{Retour sur les foncteurs polynomiaux \hts}

\begin{prop}\label{pr-htshd_polp} Supposons que $p$ est premier, que $A$ est de caractéristique $p$ et que la condition \textnormal{(CAD)}\index{nota}{CAD@(CAD), (CAD)$^+$ \emph{(conditions d'annulation des dérivations)}} est vérifiée. Alors tout foncteur \ph\ \htr\ de $\fct(\A,\E)$ est \hd.
\end{prop}

\begin{proof}
Il suffit de montrer que tout foncteur \ph\ $F$ de $\Tt_A^\cf(\A,\E)$ est \hd. En utilisant le théorème~\ref{th-struct_htr} et la proposition~\ref{pr-elt-ord}, on voit qu'il suffit de le faire dans chacun des deux cas suivants :
\begin{itemize}
\item $F$ appartient à $\Tt^\cf_{[\mathfrak{M}_{\mathrm{poli},\alpha}]}(\A,\E)$, où $\alpha : A\to K$ est un morphisme d'anneaux d'image infinie ;
\item $F$ appartient à $\Tt^\cf_{[\mathfrak{M}_{\mathrm{tord},I}]}(\A,\E)$, où $I$ est un idéal maximal cofini de $A$ tel que $\mathrm{car}(A/I)=p$.
\end{itemize}

Dans chacun des deux cas, il existe $\mathfrak{p}\in\spec(A)$ tel que le corps résiduel $\kappa:=\mathrm{Frac}(A/\mathfrak{p})$ se plonge dans $K$ et $r\in\mathbb{N}$ tel que $\tilde{A}_F$ soit un quotient de $A_\mathfrak{p}/\mathfrak{p}^r A_\mathfrak{p}$, par les propositions \ref{pr-tordI} et \ref{pr-poli_alpha}. Comme la condition (CAD) est vérifiée, la proposition~\ref{pr-CAD_fpalg} montre que $\tilde{A}_F$ est un corps parfait (ou l'anneau nul). Le corollaire~\ref{cor-poidsf-carp} permet alors de conclure.
\end{proof}

\begin{rema}
Dans $\F(A,K)$, la condition (CAD) est aussi évidemment \textit{nécessaire} pour que tout foncteur \ph\ \htr\ soit \hd\ --- en effet, elle signifie exactement que tout foncteur \textit{additif} de $\F(A,K)$ possédant une décomposition en poids faible (i.e. \htr) admet une décomposition en poids forte (i.e. est \hd).
\end{rema}

\begin{rema}
On peut aussi déduire directement le cas particulier fondamental suivant de la proposition~\ref{pr-htshd_polp} d'un résultat cohomologique général \cite[th.~7.5]{DT-str_pol} : si $p$ est non nul, que $A$ est de caractéristique $p$ et que (CAD) est satisfaite, alors tout foncteur polynomial \htr\ de $\F^\df(A,K)$ est \hd.
\end{rema}

\begin{rema}
On déduit directement du théorème~\ref{th-struct_htr} que tout foncteur \htr\ de $\F(A,K)$ est \hd\ si et seulement si tout foncteur \ph\ \htr\ de $\F(A,K)$ est \hd\ et que tout idéal maximal $K$-cotrivial de $A$ est égal à son carré.
\end{rema}

\begin{coro}\label{cor-polpoihdgrat} Supposons que $K$ est un corps de décomposition de la catégorie $\mathbf{P}(A)$ et que la condition \textnormal{(CAD)}\index{nota}{CAD@(CAD), (CAD)$^+$ \emph{(conditions d'annulation des dérivations)}} est vérifiée. Supposons de plus que l'une des conditions suivantes est satisfaite :
\begin{enumerate}
    \item $p=0$ ;
    \item $p>0$ et $A$ une algèbre essentiellement de type fini sur une $\FF_p$-algèbre $p$-parfaite (par exemple, $A$ est un anneau de type fini de caractéristique $p$).
\end{enumerate}
Alors tout foncteur polynomial de $\fct(\A,\E)$ possédant une décomposition en poids faible est \hd.
\end{coro}

\begin{proof} Le corollaire~\ref{cor-polpoihtr} (si $p=0$) ou la proposition~\ref{pr-polpoihtrp} (si $p>0$) montrent qu'un tel foncteur est \htr. On conclut qu'il est \hd\ en utilisant la proposition~\ref{pr-tppfor} (si $p=0$) ou la proposition~\ref{pr-htshd_polp}.
\end{proof}

Pour établir le corollaire~\ref{cor-polpoiforgratp} ci-après, nous aurons besoin de la variante suivante du corollaire précédent :

\begin{coro}\label{cor-polpoihdgratbis} Supposons que la condition \textnormal{(CAD)}$^+$\index{nota}{CAD@(CAD), (CAD)$^+$ \emph{(conditions d'annulation des dérivations)}} est vérifiée, ainsi que l'une des conditions suivantes :
\begin{enumerate}
    \item $p=0$ ;
    \item $p>0$ et $A$ une algèbre essentiellement de type fini sur une $\FF_p$-algèbre $p$-parfaite (par exemple, $A$ est un anneau de type fini de caractéristique $p$).
\end{enumerate}

Soit $F$ un foncteur polynomial de $\fct(\A,\E)$ possédant une décomposition en poids faible finie. Alors il existe une extension finie $L$ du corps $K$ telle que le foncteur $F\underset{K}{\boxtimes} L$ de $\fct(\A,\E\underset{K}{\boxtimes} L)$ soit \hd.

Si la condition \textnormal{(MTF)}\index{nota}{MTF@(MTF) \emph{(condition de finitude sur les morphismes)}} est également satisfaite, la même conclusion vaut pour les foncteurs \phs\ de $\fct(\A,\E)$ possédant une décomposition en poids faible finie.
\end{coro}

\begin{proof} Le corollaire~\ref{cor-polpoihtr} et la proposition~\ref{pr-dppol} (si $p=0$) ou la proposition~\ref{pr-polpoihtrp} et la remarque~\ref{rq-complPAScordec} (si $p>0$) montrent l'existence d'une extension finie $L$ de $K$ telle que $F\underset{K}{\boxtimes} L$ soit \htr. Les propositions \ref{pr-tppfor} (si $p=0$) ou \ref{pr-htshd_polp} permettent de conclure.
\end{proof}

\begin{rema}
Si le corps $K$ est parfait et que les poids de $F$ sont déployés, alors $F$ lui-même est \hd, sans nécessité d'une extension de corps au but (cf. proposition~\ref{pr-simplesPolHD}).
\end{rema}

L'énoncé suivant, qui découle directement du corollaire~\ref{cor-polpoihdgratbis}, complète le sens réciproque dans la proposition~\ref{pr-tppfor} :

\begin{coro}\label{cor-polpoiforgratp} Supposons $p>0$ et la condition \textnormal{(CAD)}$^+$\index{nota}{CAD@(CAD), (CAD)$^+$ \emph{(conditions d'annulation des dérivations)}} vérifiée. Supposons également que $A$ une algèbre essentiellement de type fini sur une $\FF_p$-algèbre $p$-parfaite.

Alors toute décomposition en poids d'un foncteur polynomial de $\fct(\A,\E)$ est forte. C'est également vrai pour les foncteurs \phs\ si la condition \textnormal{(MTF)} est satisfaite.
\end{coro}

\chapter{Foncteurs de type fini à valeurs de dimensions finies}\label{stfdf}

\begin{cvi} Dans ce chapitre comme dans les précédents, $\A$ désigne une catégorie additive $A$-linéaire essentiellement petite.
\end{cvi}

Tout foncteur de type fini de $\F^\df(\A;K)$ étant \htr\ quitte à agrandir le corps $K$ (proposition~\ref{pr-tfdf_htr-ext}), il n'est pas difficile de déduire des considérations des chapitres précédents, notamment du chapitre~\ref{shts}, des résultats significatifs pour un tel foncteur. Nous démontrerons ainsi deux conjectures de \cite[conj.~11.11]{DTV}, toutefois sous l'hypothèse supplémentaire, générale dans ce travail, que l'anneau à la source soit commutatif. Ces conjectures traitent des propriétés localement finie et localement noethérienne de $\F^\df(A,K)$ ; nous verrons à la section~\ref{sAsqf} que la première est vérifiée si et seulement si $A$ n'a pas de quotient fini de caractéristique $p$, et, à la section~\ref{sAgal}, que la seconde toujours satisfaite. Nous utiliserons pour cela un théorème de décomposition primaire directement déduit du théorème~\ref{th-struct_htr} et présenté à la section~\ref{sAtf}, qui l'applique également à l'existence de résolutions projectives de type fini pour les foncteurs de type fini de $\F^\df(A,K)$ (utilisant et généralisant notablement les résultats principaux de \cite{DT-schw}), sous une hypothèse de finitude sur $A$.

\section{Cas d'un anneau $A$ sans quotient fini}\label{sAsqf}

On rappelle que l'anneau $A$ est dit \emph{sans quotient fini} s'il n'admet pas d'idéal propre cofini.

\begin{theo}\label{th-asqfc}
Supposons que l'anneau $A$ est sans quotient fini. Alors tout foncteur $F$ de type fini de $\F^\df(\A;K)$ est polynomial.
\end{theo}

Nous donnerons, exceptionnellement, deux démonstrations de ce résultat. Indépendantes, elles établissent en fait chacune un résultat plus général que le théorème précédent.

\begin{proof}[Première démonstration]
Comme $A$ est sans quotient fini, la catégorie $\A$ ne possède aucun idéal $K$-cotrivial strict. Il s'ensuit, par le théorème~\ref{th-DTVglob}, que $F$ est \ph. La conclusion résulte alors du corollaire~\ref{cor-EMLpol-pol}.
\end{proof}

\begin{proof}[Deuxième démonstration]
Il existe une extension de corps $K\subset L$ telle que le foncteur $F_L:=F\otimes_K L$ de $\F(\A;L)$ soit \htr\, (par exemple, une clôture algébrique de $K$, par la proposition~\ref{pr-df_htr}). Comme $F$ est polynomial si et seulement si $F_L$ l'est, on peut supposer que $F$ lui-même est \htr. La proposition~\ref{pr-dfpord} montre que tous les poids partiels de $F$ sont ordinaires. Le fait que $A$ est sans quotient fini, combiné à la deuxième assertion de la proposition~\ref{pr-decpord} et à la proposition~\ref{pr-dppol}, montre que tous les poids partiels de $F$ sont sans torsion. La proposition~\ref{pr-raispol} donne alors la conclusion.
\end{proof}

Nous abordons maintenant une situation plus générale que celle où l'anneau $A$ est sans quotient fini, celle où $A$ n'a pas de quotient fini de caractéristique $p$. Le théorème suivant répond par l'affirmative à la conjecture \cite[conj.~11.11\,(2)]{DTV}, dans le cas d'un anneau \emph{commutatif} à la source.

\begin{theo}\label{th-PAP_df} Les assertions suivantes sont équivalentes :
\begin{enumerate}
\item l'anneau $A$ ne possède pas de quotient fini de caractéristique $p$ ;
\item tout foncteur de type fini de $\F^\df(A,K)$ possède une décomposition à la Steinberg de type $PAP$ ;
\item la catégorie $\F^\df(A,K)$ est localement finie.
\end{enumerate}

Si ces conditions sont vérifiées et que $\mathrm{(CFD)}$\index{nota}{CFD@(CFD), (CFD)$^+$ \emph{(conditions de finitude sur les dérivations)}} l'est également, alors la catégorie $\Tt(A,K)$ est localement finie. 
\end{theo}

\begin{proof}
Le théorème~\ref{th-DTVglob} et le corollaire~\ref{cor-EMLpol-pol} montrent que la première condition entraîne la seconde. La proposition~\ref{pr-PSmodif} et \cite[prop.~4.9 et lemme~11.10]{DTV} montrent que la seconde implique la troisième. Il est classique que la troisième assertion entraîne la première, car, si $k$ est un corps fini de caractéristique $p$, la catégorie $\F(k,K)$ n'est pas localement finie (le foncteur projectif de type fini $P^k$ n'est pas fini, car sa filtration copolynomiale est infinie).

Si la condition $\mathrm{(CFD)}$ est satisfaite, le théorème~\ref{th-htr_ln} montre que tout foncteur de type fini de $\Tt(A,K)$ appartient à $\F^\df(A,K)$, et donc est fini si $A$ n'a pas de quotient fini de caractéristique $p$, d'après ce qu'on vient de voir. Cela termine la démonstration.
\end{proof}

Le résultat précédent s'étend à certains anneaux non commutatifs :

\begin{theo}\label{th-PAP_df-noncom} Supposons que l'anneau $A$ n'a pas de quotient fini de caractéristique $p$ et que $k$ est une $A$-algèbre non nécessairement commutative. Alors tout foncteur de type fini $F$ de $\F^\df(k,K)$ admet une décomposition à la Steinberg de type $PAP$ et est noethérien.
\end{theo}

\begin{proof}
Le théorème~\ref{th-DTVglob} et le corollaire~\ref{cor-EMLpol-pol} montrent que $F$ possèdent une décomposition à la Steinberg de type $PAP$. Comme la proposition~\ref{pr-PSmodif} et \cite[prop.~4.9 et lemme~11.10]{DTV} valent également pour un anneau non commutatif à la source, on en déduit que $F$ est noethérien comme dans la démonstration précédente.
\end{proof}

\section{Décomposition primaire}\label{sAtf}

\subsection{Préliminaires}

Commençons par introduire une variante de la condition (IMF) (un peu plus forte que cette dernière) introduite page~\pageref{pageimf} :

\label{pageIMQ}\textnormal{(IMQ) Si $\mathfrak{m}$ est un idéal maximal de $A$ tel que $A/\mathfrak{m}$ soit un corps fini de caractéristique $p$, alors $\mathfrak{m}/\mathfrak{m}^2$ est fini.}

Cette condition est en particulier vérifiée lorsque $K$ est de caractéristique $p=0$ ou lorsque $A$ est noethérien.

Il nous sera par ailleurs utile l'introduire les notations suivantes.

\begin{nota} On note $\pol^\df_d(\A;K)$ la sous-catégorie pleine $\pol_d(\A;K)\cap\F^\df(\A;K)$ de $\F(\A;K)$ (qui est épaisse), et $\overline{\pol}^\df_d(\A;K)$ la sous-catégorie prélocalisante de $\F(A,K)$ engendrée par $\pol^\df_d(\A;K)$.

Lorsque $\A=\mathbf{P}(A)$, ces sous-catégories de $\F(A,K)$ seront notées simplement $\pol^\df_d(A,K)$ et $\overline{\pol}^\df_d(A,K)$ respectivement.
\end{nota}

Autrement dit, un foncteur $F$ de $\F(\A;K)$ appartient à $\overline{\pol}^\df_d(\A;K)$ si et seulement s'il appartient à $\pol_d(\A;K)$ et est localement à valeurs de dimensions finies, i.e. colimite de sous-foncteurs à valeurs de dimensions finies. 

L'observation suivante  découle par exemple de \cite[lemme~11.10]{DTV}.

\begin{lemm}\label{lm-lfevprel} La catégorie $\overline{\pol}^\df_d(A,K)$ est une catégorie de Grothendieck localement finie.
\end{lemm}

Nous déduirons le théorème de structure principal de cette section de ceux de la section~\ref{ss-tshtr} et du lemme formel simple suivant sur l'extension des scalaires dans les catégories abéliennes.

\begin{lemm}\label{lm-cbcatablpet} Soient $\E_1$, $\E_2$ des catégories de Grothendieck $K$-linéaires, avec $\E_1$ localement de type fini, $\Phi : \E_1\to\E_2$ un foncteur exact, cocontinu et $K$-linéaire, et $K\subset L$ une extension de corps. On suppose que le foncteur $\Phi\underset{K}{\boxtimes}L : \E_1\underset{K}{\boxtimes}L\to\E_2\underset{K}{\boxtimes}L$ est pleinement fidèle et que son image essentielle est une sous-catégorie semi-épaisse de $\E_2\underset{K}{\boxtimes}L$.

Alors $\Phi$ est pleinement fidèle et son image essentielle $\E'$ est une sous-catégorie semi-épaisse de $\E_2$. De plus, tout objet $T$ de $\E_2$ tel que $T\underset{K}{\otimes}L$ appartienne à l'image essentielle de $\Phi\underset{K}{\boxtimes}L$ appartient à $\E'$.
\end{lemm}

\begin{proof}Soient $X$ et $Y$ des objets de $\E_1$.
L'application $L$-linéaire canonique $\E_1(X,Y)\underset{K}{\otimes} L\to (\E_1\underset{K}{\boxtimes}L)(X\underset{K}{\otimes}L,Y\underset{K}{\otimes}L)$ est injective, et bijective si $X$ est de type fini, en raison de l'isomorphisme d'adjonction $(\E_1\underset{K}{\boxtimes}L)(X\underset{K}{\otimes}L,Y\underset{K}{\otimes}L)\simeq\E_1(X,\mathrm{Ou}(Y\underset{K}{\otimes}L))=\E_1(X,L\underset{K}{\otimes}Y)$, où $\mathrm{Ou} : \E_1\underset{K}{\boxtimes}L\to\E_1$ désigne le foncteur d'oubli.

Comme $\Phi$ est $K$-linéaire et cocontinu, on dispose par ailleurs d'un isomorphisme canonique $(\Phi\underset{K}{\boxtimes}L)(X\underset{K}{\otimes}L)\simeq\Phi X\underset{K}{\otimes}L$. On en déduit un diagramme commutatif
$$\xymatrix{\E_1(X,Y)\underset{K}{\otimes} L\ar[r]\ar[d] & \E_2(\Phi X,\Phi Y)\underset{K}{\otimes} L\ar[d]\\
(\E_1\underset{K}{\boxtimes}L)(X\underset{K}{\otimes}L,Y\underset{K}{\otimes}L)\ar[r] & (\E_2\underset{K}{\boxtimes}L)(\Phi X\underset{K}{\otimes}L,\Phi Y\underset{K}{\otimes}L)
}$$
dont les lignes sont induites par $\Phi$ et $\Phi\underset{K}{\boxtimes}L$ respectivement et dont les colonnes sont bijectives si $X$ est de type fini. Comme $\Phi\underset{K}{\boxtimes}L$ est pleinement fidèle, la ligne inférieure est également bijective. Il s'ensuit que la ligne supérieure est bijective, donc que le morphisme $\E_1(X,Y)\to\E_2(\Phi X,\Phi Y)$ induit par $\Phi$ est bijectif, dès lors que $X$ est de type fini. Comme $\Phi$ est cocontinu, l'application $\E_1(X,Y)\to\E_2(\Phi X,\Phi Y)$ est aussi bijective si $X$ est localement de type fini. Comme $\E_1$ est localement de type fini, cela montre que $\Phi$ est pleinement fidèle.

Comme $\Phi$ est exact et pleinement fidèle, pour montrer que son image essentielle $\E'$ est semi-épaisse dans $\E_2$, il suffit d'établir qu'elle est stable par sous-objet. Soient $X$ un objet de $\E_1$ et $T$ un sous-objet de $\Phi(X)$. Comme l'image essentielle de $\Phi\underset{K}{\boxtimes}L$ est semi-épaisse, il existe un objet $Y$ de $\E_1\underset{K}{\boxtimes}L$ tel que $(\Phi\underset{K}{\boxtimes}L)(Y)\simeq T\underset{K}{\otimes}L\subset\Phi(X)\underset{K}{\otimes}L\simeq (\Phi\underset{K}{\boxtimes}L)(X\underset{K}{\otimes}L)$. En appliquant le foncteur d'oubli $\E_2\underset{K}{\boxtimes}L\to\E_2$, on en déduit un isomorphisme $\Phi(\mathrm{Ou}(Y))\simeq L\underset{K}{\otimes}T$ dans $\E_2$. En particulier, $T$ est facteur direct de $\Phi(\mathrm{Ou}(Y))$ : il existe un idempotent de $\mathrm{End}_{\E_2}(\Phi(\mathrm{Ou}(Y)))$ dont l'image est $T$. Comme $\Phi$ est pleinement fidèle, cet idempotent provient d'un idempotent de $\mathrm{End}_{\E_1}(\mathrm{Ou}(Y))$. L'image de cet idempotent fournit un objet $Z$ de $\E_1$ tel que $\Phi(Z)\simeq T$.

Considérons maintenant un objet $T$ de $\E_2$ tel que $T\underset{K}{\otimes}L\simeq (\Phi\underset{K}{\boxtimes}L)(X)$ pour un objet $X$ de $\E_1$ : appliquant de nouveau le foncteur d'oubli, on obtient $L\underset{K}{\otimes}T\simeq\Phi(\mathrm{Ou}(X))$ dans $\E_2$ ; ainsi, $T$, qui est isomorphe à un sous-objet de $L\underset{K}{\otimes}T$, appartient à $\E'$ d'après ce qui précède, ce qui termine la démonstration.
\end{proof}

\subsection{Résultat fondamental}

\begin{prdef}\label{prdf-auxID} Soient $n\in\mathbb{N}$, $I_1,\dots,I_n$ des idéaux de $A$ et $d\in\mathbb{N}$. On note $\mathfrak{p}_i:=\sqrt{I_i}$, pour $i\in\llbracket 1,n\rrbracket$, la racine de l'idéal $I_i$ ; on suppose que $\mathfrak{p}_1,\dots,\mathfrak{p}_n$ sont des idéaux premiers deux à deux distincts. Soient
$$E:=\{i\in\llbracket 1,n\rrbracket\,|\,\cd(A/\mathfrak{p}_i)=\infty\}$$
et $\C_d[I_1,\dots,I_n]$ la sous-catégorie pleine de $\F^\df\big(\prod_{i=1}^n(\A\otimes_A (A_{\mathfrak{p}_i}/I_i.A_{\mathfrak{p}_i}));K\big)$ constituée des multifoncteurs qui sont polynomiaux de degré au plus $d$ par rapport à la $i$-ème variable lorsque $i\in E$.

Désignons par $\Phi_d^{I_1,\dots,I_n} : \C_d[I_1,\dots,I_n]\to\F(\A;K)$ la restriction à $\C_d[I_1,\dots,I_n]$ de la précomposition par le foncteur canonique $\A\to\prod_{i=1}^n(\A\otimes_A (A_{\mathfrak{p}_i}/I_i.A_{\mathfrak{p}_i}))$. Alors $\Phi_d^{I_1,\dots,I_n}$ est pleinement fidèle et son image essentielle est une sous-catégorie semi-épaisse de $\F(\A;K)$.
\end{prdef}

Lorsqu'une ambiguïté est possible sur le corps $K$ au but, nous noterons $\C_d^K[I_1,\dots,I_n]$ (resp. $\Phi_d^{I_1,\dots,I_n;K}$) pour $\C_d[I_1,\dots,I_n]$ (resp. $\Phi_d^{I_1,\dots,I_n}$). La proposition précédente sera démontrée en même temps que le théorème fondamental suivant.

\begin{theo}\label{th-tfdf_gal} Soit $F$ un foncteur de type fini de $\F^\df(\A;K)$. Il existe des entiers naturels $n, d,r$, des idéaux $I_1,\dots,I_n$ de $A$ et $\mathfrak{p}_1,\dots,\mathfrak{p}_n\in\spec(A)$ tels que :
\begin{enumerate}
    \item\label{itthps1} $\forall i\in\llbracket 1,n\rrbracket\qquad \mathfrak{p}_i^r\subset I_i\subset\mathfrak{p}_i$ ;
    \item\label{itthps2} si $\mathrm{car}(A/\mathfrak{p}_i)\ne p$, alors $A/I_i$ est fini ;
    \item $F$ appartient à l'image essentielle de $\Phi_d^{I_1,\dots,I_n}$.
\end{enumerate}
\end{theo}

\begin{proof}[Démonstration de la proposition~\ref{prdf-auxID} et du théorème~\ref{th-tfdf_gal}.] La sous-catégorie $\C_d^K[I_1,\dots,I_n]$ de $\F\big(\prod_{i=1}^n(\A\otimes_A (A_{\mathfrak{p}_i}/I_i.A_{\mathfrak{p}_i}));K\big)$ est semi-épaisse (et même épaisse) ; notons $\bar{\C}_d^K[I_1,\dots,I_n]$ la sous-catégorie prélocalisante qu'elle engendre, qui est donc une catégorie de Grothendieck localement de type fini.

Supposons d'abord le corps $K$ algébriquement clos. Tout foncteur de $\C_d^K[I_1,\dots,I_n]$, et donc de $\bar{\C}_d^K[I_1,\dots,I_n]$, est alors \htr\ (proposition~\ref{pr-df_htr}). Le corollaire~\ref{cor-septorpol}, la proposition~\ref{pr-indep_alpha} et le théorème~\ref{th-decPPart} montrent que la restriction $\bar{\Phi}_d^{I_1,\dots,I_n;K}$ à $\bar{\C}_d^K[I_1,\dots,I_n]$ de la précomposition par le foncteur canonique $\A\to\prod_{i=1}^n(\A\otimes_A (A_{\mathfrak{p}_i}/I_i.A_{\mathfrak{p}_i}))$ est pleinement fidèle et a pour image essentielle une sous-catégorie \emph{semi-épaisse} de $\F(\A;K)$ (constituée de foncteurs \hts\ dont les poids partiels appartiennent au sous-monoïde $\mon(A_\mu,K_\mu)$ engendré par les images canoniques de $\mon((A/\mathfrak{p}_i)_\mu,K_\mu)$ pour $i\in\llbracket 1,n\rrbracket$). La proposition~\ref{pr-df_htr} et le théorème~\ref{th-struct_htr} (ainsi que les propositions \ref{pr-tordI} et \ref{pr-poli_alpha}) établissent pour leur part le théorème~\ref{th-tfdf_gal} pour $K$ algébriquement clos.

On note par ailleurs que $\Phi_d^{I_1,\dots,I_n;K}$ est un foncteur exact, cocontinu et $K$-linéaire. De plus, pour toute extension de corps $K\subset L$, on dispose d'un diagramme commutatif
$$\xymatrix{\bar{\C}_d^K[I_1,\dots,I_n]\underset{K}{\boxtimes}L\ar[r]^-\simeq\ar[d]_{\bar{\Phi}_d^{I_1,\dots,I_n;K}\underset{K}{\boxtimes}L} & \bar{\C}_d^L[I_1,\dots,I_n]\ar[d]^{\bar{\Phi}_d^{I_1,\dots,I_n;L}} \\
\F(\A;K)\underset{K}{\boxtimes}L\ar[r]^-\simeq & \F(\A;L)
\;.}$$

En prenant pour $L$ une clôture algébrique de $K$, on voit que la proposition~\ref{prdf-auxID} et le théorème~\ref{th-tfdf_gal} découlent du cas particulier établi ci-dessus et du lemme~\ref{lm-cbcatablpet}.
\end{proof}

Lorsque la condition (IMQ) est vérifiée, avec les notations du théorème~\ref{th-tfdf_gal} et de la proposition-définition~\ref{prdf-auxID}, tous les idéaux $I_i$ tels que $i\notin E$ sont cofinis, d'où le corollaire suivant.

\begin{coro}\label{cor-imq} Supposons que la condition \textnormal{(IMQ)} est vérifiée. Soit $F$ un foncteur de type fini de $\F^\df(\A;K)$. Alors il existe un idéal cofini $I$ de $A$, un entier naturel $d$ et une sous-catégorie épaisse $\C$ de $\fct(\A/I,\pol^\df_d(\A;K))$ tels que la restriction à $\C\subset\fct(\A/I,\F(\A;K))\simeq\F(\A/I\times\A;K)$
de la précomposition par le foncteur canonique $\A\to\A/I\times\A$ induise une équivalence entre $\C$ et une sous-catégorie semi-épaisse de $\F(\A;K)$ contenant $F$. 
\end{coro}

Le corollaire précédent dit en particulier que, sous l'hypothèse (IMQ), tout foncteur de type fini de $\F(A,K)$ possède une décomposition de type $FP$. Nous verrons à l'exemple~\ref{cex-tFP} que ce n'est pas toujours vrai sans la condition (IMQ).

L'énoncé suivant, qui est à rapprocher du théorème~\ref{th-htr_ln}, sera bientôt généralisé, au théorème~\ref{th-noeth_cor_gal}.

\begin{coro}\label{cor-noecorpg} Si la condition \textnormal{(IMQ)} est vérifiée, alors tout foncteur de type fini de $\F^\df(A,K)$ est noethérien.
\end{coro}

\begin{proof}
La catégorie de Grothendieck $\overline{\pol}^\df_d(A,K)$ étant localement noethérienne (lemme~\ref{lm-lfevprel}), la conclusion résulte du corollaire~\ref{cor-imq} et du théorème~\ref{th-PSS}.
\end{proof}

\subsection{Application à la propriété $pf_\infty$}

Le théorème suivant fournit une généralisation substantielle de \cite[cor.~8.3]{DT-schw}.

\begin{theo}\label{th-pfi_corps} Supposons que les anneaux $A$ et $A\otimes_\mathbb{Z}K$ sont noethériens (par exemple, que l'anneau $A$ est essentiellement de type fini). Alors tout foncteur de type fini de $\F^\df(A,K)$ appartient à $\pf_\infty(\F(A,K))$.\index{nota}{pf@$\pf_n$}
\end{theo}

(On rappelle que $\pf_\infty(\F(A,K))$ désigne la classe des foncteurs de $\F(A,K)$ qui possèdent une résolution projective dont tous les termes sont de type fini --- cf. page~\pageref{pagepfn}, ainsi que \cite[section~1]{DT-schw}.)

\begin{proof} Pour tout idéal $I$ de $A$, le foncteur $\F(A/I\times A,K)\to\F(A,K)$ de précomposition par le foncteur canonique $\mathbf{P}(A)\to\mathbf{P}(A/I)\times\mathbf{P}(A)$ préserve la propriété $pf_\infty$ grâce à \cite[cor.~7.3 et prop.~1.3]{DT-schw} et à la noethérianité de l'anneau $A$. Le corollaire~\ref{cor-imq} montre alors qu'il suffit d'établir que, pour tout $d\in\mathbb{N}$, tout idéal cofini $I$ de $A$ et tout foncteur de type fini $F$ de $\fct(\mathbf{P}(A/I),\F(A,K))$ à valeurs dans $\pol_d^\df(A,K)$ possède la propriété $pf_\infty$.

Comme $F$ est un foncteur de type fini de la catégorie de Grothendieck localement noethérienne $\fct(\mathbf{P}(A/I),\overline{\pol}^\df_d(A,K))$ (cf. démonstration du corollaire~\ref{cor-noecorpg}), il vérifie la propriété $pf_\infty$, et donc a fortiori $psf_\infty$ (notion plus générale que $pf_\infty$ discutée dans \cite[§\,2]{DT-schw}), dans $\fct(\mathbf{P}(A/I),\overline{\pol}^\df_d(A,K))$ (cf. \cite[prop.~1.3]{DT-schw}). Or la propriété $psf_\infty$ est préservée par post-composition par tout foncteur exact cocontinu (grâce à \cite[prop.~2.8 et 2.6]{DT-schw} par exemple) : appliquant ceci à l'inclusion $\overline{\pol}^\df_d(A,K)\hookrightarrow\F(A,K)$, on voit que le foncteur $F$ vérifie aussi $psf_\infty$ dans $\fct(\mathbf{P}(A/I),\F(A,K))$. Ses valeurs sont des foncteurs polynomiaux de type fini (car le foncteur $F$ est de type fini et l'anneau $A/I$ fini) de $\F^\df(A,K)$, elles appartiennent donc à $\pf_\infty(\F(A,K))$ par \cite[th.~6.14]{DT-schw} (dont les hypothèses sont vérifiées parce que $A$ et $A\otimes_\mathbb{Z}K$ sont des anneaux noethériens). On conclut en appliquant \cite[cor.~2.15]{DT-schw}.
\end{proof}

\begin{rema}
On peut obtenir la même conclusion (ou la propriété $pf_n$ pour les foncteurs de type fini de $\F^\df(A,K)$, pour un $n\in\mathbb{N}$) sous des hypothèses un peu plus faibles sur $A$, mais techniques, à savoir celles de \cite[cor.~8.2]{DT-schw}. Toutefois, les hypothèses de ce dernier énoncé semblent le plus souvent difficiles à vérifier hors du cas particulier fondamental où les anneaux $A$ et $A\otimes_\mathbb{Z}K$ sont noethériens.

On peut aussi affaiblir légèrement l'hypothèse de commutativité de $A$ : la conclusion du théorème~\ref{th-pfi_corps} persiste si l'on remplace $A$ un anneau non commutatif qui est de type fini comme module à droite sur un sous-anneau commutatif vérifiant les hypothèses du théorème (la démonstration, à partir du théorème~\ref{th-pfi_corps}, est essentiellement la même que celle du corollaire~\ref{cor-noeth_corg} ci-après, en utilisant \cite[cor.~7.3]{DT-schw}).
\end{rema}

\section[La propriété localement noethérienne]{La propriété localement noethérienne pour $\F^\df(A,K)$ dans le cas général}\label{sAgal}

L'énoncé suivant ne nécessite \emph{aucune} hypothèse sur l'anneau $A$ (sinon la commutativité).

\begin{theo}\label{th-noeth_cor_gal}
Tout foncteur de type fini $F$ de $\F^\df(A,K)$ est noethérien.
\end{theo}

\begin{proof}
Le corollaire~\ref{cor-acdf} fournit un sous-anneau $A'$ de type fini, donc noethérien et a fortiori vérifiant (IMQ), de $A$ tel que le foncteur $\Phi : \F(A,K)\to\F(A',K)$ de précomposition par l'extension des scalaires $\mathbf{P}(A')\to\mathbf{P}(A)$ envoie $F$ sur un foncteur \emph{de type fini} de $\F^\df(A',K)$. Le corollaire~\ref{cor-noecorpg} montre que $\Phi(F)$ est noethérien. Comme $\Phi$ est exact et fidèle, cela entraîne que $F$ est lui-même noethérien.
\end{proof}

\begin{rema}\label{rq-cas_embetant}
L'argument de changement de base à la source qu'on vient d'utiliser est efficace pour établir la propriété noethérienne mais échoue à dire beaucoup plus sur la structure des foncteurs de type fini de $\F^\df(A,K)$. Supposons ainsi que $A$ est un anneau local, semi-primaire, infini, dont le corps résiduel $\kappa$ est fini de caractéristique $p$. Alors les décompositions en poids, même des effets croisés, d'un foncteur de $\F^\df(A,K)$ ne disent essentiellement rien de lui en général, car la projection $A\twoheadrightarrow\kappa$ induit un isomorphisme $\mon(\kappa_\mu,K_\mu)\xrightarrow{\simeq}\mon(A_\mu,K_\mu)$, qui est donc \emph{fini}.
Dans cette situation, nous ignorons par exemple si tout foncteur de type fini de $\F^\df(A,K)$ possède une filtration finie dont les sous-quotients possèdent une décomposition de type $FP$. (Cf. exemple~\ref{cex-tFP}.)
\end{rema}

Le théorème~\ref{th-noeth_cor_gal} répond positivement à la conjecture~11.11 (1) de \cite{DTV} dans le cas d'un anneau commutatif à la source. Les méthodes du présent article semblent toutefois impuissantes à s'affranchir de l'hypothèse de commutativité, au-delà de la généralisation légère et immédiate suivante.

\begin{coro}\label{cor-noeth_corg}
Soit $R$ un anneau non commutatif. On suppose que $R$ est de type fini comme module à droite sur un sous-anneau commutatif $A$. Alors tout foncteur de type fini $F$ de $\F^\df(R,K)$ est noethérien.
\end{coro}

\begin{proof}
Comme $R$ est un $A$-module à droite de type fini, le foncteur de précomposition par l'extension des scalaires $\F(R,K)\to\F(A,K)$ préserve les foncteurs de type fini. La conclusion se déduit donc de l'exactitude et de la fidélité de ce foncteur ainsi que du théorème~\ref{th-noeth_cor_gal}.
\end{proof}

Combiné au corollaire~\ref{cor-htn_df}, le théorème~\ref{th-noeth_cor_gal} fournit l'équivalence suivante :
\begin{coro}\label{cor-car_df} Soit $F$ un foncteur de type fini de $\F(A,K)$. Les assertions suivantes sont équivalentes :
\begin{enumerate}
\item $F$ appartient à $\F^\df(A,K)$ ;
\item $F$ est noethérien et il existe une extension de corps finie $K\subset L$ telle que le foncteur $F\otimes L$ de $\F(A,L)$ soit \htr.
\end{enumerate}
\end{coro}

\begin{proof}
La première condition implique la deuxième en vertu du théorème~\ref{th-noeth_cor_gal} et de la proposition~\ref{pr-tfdf_htr-ext}.

Réciproquement, si $F$ est noethérien et que $K\subset L$ est une extension de corps \emph{finie}, alors $F\otimes L$ est noethérien dans $\F(A,K)$, donc a fortiori dans $\F(A,L)$, on voit donc par le corollaire~\ref{cor-htn_df} que $F\otimes L$ appartient à $\F^\df(A,L)$. Aussi $F$ appartient-il à $\F^\df(A,K)$.
\end{proof}

\begin{rema}\label{rq-pfiag} Contrairement à la propriété noethérienne ici considérée, la propriété $pf_\infty$ pour les foncteurs de type fini et à valeurs de dimensions finies de $\F(A,K)$, montrée au théorème~\ref{th-pfi_corps} lorsque $A$ vérifie des hypothèses de finitude satisfaites en particulier lorsque $A$ est un anneau de type fini, ne s'étend pas au cas d'un anneau (même commutatif) arbitraire $A$. De fait, \cite[th.~13.17]{DTV} donne une condition nécessaire et suffisante sur $A$ (et $K$) pour que tout foncteur \emph{simple} à valeurs de dimensions finies de $\F^\df(A,K)$ soit de présentation finie. Il existe des anneaux, et même des corps (commutatifs), qui ne la vérifient pas : si $K$ n'est pas un corps de type fini, le foncteur \guillemotleft~identité~\guillemotright\ de $\F(K,K)$, qui est additif, simple et à valeurs de dimensions finies, n'est pas de présentation finie.
\end{rema}

\chapter{Foncteurs de type fini à valeurs de type fini}\label{stftf}

\begin{cvi}
Dans ce chapitre, $\A$ désigne toujours une petite catégorie additive $A$-linéaire essentiellement petite. La lettre $k$ désigne un anneau \emph{noethérien}, et $\E$ une catégorie de Grothendieck.
\end{cvi}

Nous généralisons ici plusieurs des résultats fondamentaux du chapitre précédent aux foncteurs de type fini de $\mathbf{P}(A)$ vers les $k$-modules de type fini. En revanche, le théorème de décomposition primaire~\ref{th-tfdf_gal} ne s'étend pas à cette situation (voir l'exemple~\ref{cex_dec-prim}). Les démonstrations consistent toujours à se ramener au cas où $k$ est un corps par récurrence noethérienne sur cet anneau, à l'aide d'arguments en partie formels de localisation et de quotient au but. Nous commençons toutefois par constater que la propagation de la propriété localement noethérienne dans ce type de situation (passage d'espaces vectoriels de dimension finie à des modules de type fini sur un anneau noethérien au but) n'est pas automatique, même pour des foncteurs sur une catégorie additive, avant de donner des critères élémentaires pour pouvoir le faire (section~\ref{scbpf}). Nous suivons ensuite un plan analogue à celui du chapitre~\ref{shts}, avec parfois des détours rendus nécessaires par l'absence de décomposition primaire.

\section{Changement de base et propriétés de finitude}\label{scbpf}

Les résultats de ce mémoire concernant presque toujours les catégories de foncteurs à valeurs dans des espaces vectoriels sur un corps $K$, pour en déduire des résultats, comme la noethérianité, sur des foncteurs de $\F(\A;k)$, nous sommes conduits à aborder la question plus générale suivante : soient $\E$ une catégorie de Grothendieck $k$-linéaire et $X$ un objet de $\E$. Supposons que, pour tout morphisme d'anneaux $k\to K$ où $K$ est un corps, l'objet $X\underset{k}{\otimes}K$ de $\E\underset{k}{\boxtimes}K$ est noethérien. Sous quelles hypothèses supplémentaires peut-on conclure que $X$ est noethérien ?

(Sans hypothèse supplémentaire, il est trivial que $X$ peut ne pas être noethérien, même si $\E$ est localement noethérienne, comme l'illustre l'exemple du groupe abélien $\mathbb{Q}/\mathbb{Z}$, qui devient nul après tensorisation par n'importe quel corps.)

Même lorsque $\E$ est une catégorie de foncteurs du type $\F(\C;k)$ et que $X$ est un foncteur de type fini et à valeurs de type fini, nous allons voir que la réponse est généralement négative. Nous donnerons toutefois ensuite des critères qui permettront de parvenir à nos fins lorsque la catégorie source est $\mathbf{P}(A)$, à l'aide des résultats de structure (beaucoup plus forts que les propriétés de finitude qu'ils impliquent) établis antérieurement dans ce travail.

\begin{rema}
Ces questions de comportement par changement de base des propriétés de finitude s'étendent de façon naturelle (et manifestement délicate) au produit tensoriel de catégories de Grothendieck introduit dans \cite{LRGS}. 
\end{rema}

\subsection{Un exemple de résultat négatif}

Dans cette section, nous allons construire une petite catégorie additive $\A$ telle que $\F^\df(\A;K)$ est localement noethérienne pour tout corps $K$ mais qu'il existe des foncteurs de $\F(\A;\mathbb{Z})$ de type fini, à valeurs de type fini et non noethériens.

Commençons par une construction générale. Si $\Phi : \B\to\C$ est un foncteur additif entre catégories additives, on définit une catégorie additive $\B\ltimes_\Phi\C$ par ${\rm Ob}\,(\B\ltimes_\Phi\C)={\rm Ob}\,\B\times{\rm Ob}\,\C$ et $(\B\ltimes_\Phi\C)((x,y),(x',y'))=\B(x,x')\oplus\C(y,y')\oplus\C(\Phi(x),y')$, la composition étant la composition évidente, de sorte que l'on dispose d'un diagramme de foncteurs additifs $\B\times\C\to\B\ltimes_\Phi\C\to\B\times\C$ dont la composée est l'identité.

Considérons la situation suivante : $\B=\mathbf{P}(\mathbb{Z})$, $\C=\prod_p\mathbf{P}(\mathbb{Z}/p)$, où le produit est pris sur les nombres premiers $p$, et $\Phi : \B\to\C$ est le foncteur additif dont chaque composante $\mathbf{P}(\mathbb{Z})\to\mathbf{P}(\mathbb{Z}/p)$ est la réduction modulo $p$. On considère la sous-catégorie pleine $\A$ de $\B\ltimes_\Phi\C$ des objets dont les composantes dans $\mathbf{P}(\mathbb{Z}/p)$ sont nulles sauf pour un nombre fini de nombres premiers $p$. (Noter que $\A$ est équivalente à la sous-catégorie pleine de $\mathbf{Ab}$ des groupes abéliens de type fini dont la torsion $l$-primaire est annulée par $l$ pour tout nombre premier $l$.)

Si $\I$ est un idéal $K$-cotrivial de $\A$, alors la trace de $\I$ sur $\B=\mathbf{P}(\mathbb{Z})$ est aussi un idéal cotrivial, donc de la forme $N.\B$ pour un entier $N>0$. Il s'ensuit que $\A/\I$ est un quotient de $\A/N.\B$ ; or on vérifie aisément que cette catégorie additive s'identifie à $(\mathbf{P}(\mathbb{Z}/N)\ltimes\bigoplus_{p|N}\mathbf{P}(\mathbb{Z}/p))\oplus\bigoplus_{p\nmid N}\mathbf{P}(\mathbb{Z}/p)$. Maintenant, si $\E$ est une catégorie localement noethérienne, alors $\fct(\A/N.\B,\E)$ est localement noethérienne en raison du théorème~\ref{th-PSS} et des observations suivantes :
\begin{itemize}
\item de façon générale, si $(\A_i)_{i\in E}$ est une famille de petites catégories additives, tout foncteur de type fini de $\fct(\bigoplus_{i\in E}\A_i,\E)$ se factorise par la projection $\bigoplus_{i\in E}\A_i\twoheadrightarrow\bigoplus_{i\in I}\A_i$ pour un sous-ensemble {\em fini} $I$ de $E$ ; 
\item si $I$ est un ensemble fini de nombres premiers, alors $(\mathbf{P}(\mathbb{Z}/N)\ltimes\bigoplus_{p|N}\mathbf{P}(\mathbb{Z}/p))\oplus\bigoplus_{p\in I}\mathbf{P}(\mathbb{Z}/p)$ est équivalente à $\mathbf{P}(A_I)$ pour un certain anneau non commutatif fini $A_I$.
\end{itemize}

Pour conclure que $\F^\df(\A;K)$ est localement noethérienne, il suffit donc, par le théorème~\ref{th-DTVglob}, de vérifier que les foncteurs {\em \phs} de $\F(\A;K)$ sont localement noethériens. Ce résultat découle de ce que tout foncteur \ph\ de $\F(\A;K)$ se factorise par :
\begin{itemize}
    \item la projection $\A\twoheadrightarrow\mathbf{P}(\mathbb{Z})$ si $K$ est de caractéristique nulle ;
    \item la projection $\A\twoheadrightarrow\mathbf{P}(\mathbb{Z})\ltimes\mathbf{P}(\mathbb{Z}/p)$ si $K$ est de caractéristique $p>0$.
\end{itemize}

Mais il existe des foncteurs, même additifs, à valeurs de type fini de $\F(\A;\mathbb{Z})$ qui ne sont pas localement noethériens. Considérons le foncteur additif représenté par l'objet $\mathbb{Z}$ de $\A$ : il est donné sur les objets par $(U,(V_p))\mapsto U\oplus\bigoplus V_p$ (où $U$ est un objet de $\mathbf{P}(\mathbb{Z})$ et les $V_p$, pour $p$ premier, des objets de $\mathbf{P}(\mathbb{Z}/p)$) et contient comme sous-foncteur la somme directe infinie non triviale des $(U,(V_p))\mapsto V_l$ (où $l$ est un nombre premier fixé). Il est donc de type fini mais non noethérien.

\subsection{Un résultat positif pour des foncteurs à valeurs noethériennes}\label{par-cbn}

\begin{nota}
Si $x$ est un élément de $k$, $\E$ une catégorie de Grothendieck $k$-linéaire et $T$ un objet de $\E$, on note $_x T$ le noyau de l'endomorphisme de $T$ de multiplication par $x$, et $T/x$ son conoyau. On note également $_{(x)}T$ la réunion filtrante croissante sur $n\in\mathbb{N}$ des $_{x^n}T$.
\end{nota}

\begin{hyp}
Dans toute la section~\ref{par-cbn}, on suppose que $\E$ est une catégorie de Grothendieck $k$-linéaire ; $x$ et $y$ désignent des éléments de $k$.
\end{hyp}

\begin{lemm}\label{lm-noethqxy} Soit $T$ un objet de $\E$. Si les objets $T/x$ et $T/y$ de $\E$ sont noethériens, alors $T/xy$ est noethérien.
\end{lemm}

\begin{proof}
Cela découle de la suite exacte canonique $T/x\to T/xy\to T/y\to 0$.
\end{proof}

\begin{lemm}\label{lm-noestpn} Soient $T$ un objet noethérien de $\E$ et $n$ un entier naturel. Si les sous-objets $_{x^n}T$ et $_{(x)}T$ de $T$ ont la même image dans $T/x$, alors le morphisme canonique $_{(x)}T\hookrightarrow T\twoheadrightarrow T/x^n$ est un monomorphisme.
\end{lemm}

\begin{proof}
L'hypothèse s'écrit $_{(x)}T\subset\,_{x^n}T+x.T$, d'où $_{(x)}T\subset\,_{x^n}T+x._{(x)}T$.

Comment $T$ est noethérien, il existe $N\in\mathbb{N}$ tel que $_{x^N}T=_{(x)}T$. Montrons que le plus petit entier $N$ vérifiant cette propriété est inférieur à $n$. En effet, si $N>0$, alors $x._{x^N}T\subset\,_{x^{N-1}}T$, donc l'inclusion précédente entraîne $_{(x)}T\subset\,_{x^{N-1}}T$ si $N>n$, ce qui contredit la minimalité de $N$.

Ainsi, $_{x^n}T=\,_{(x)}T$, donc le noyau $x^n.T\cap\,_{(x)}T=x^n._{(x)}T (\subset T)$ du morphisme $_{(x)}T\hookrightarrow T\twoheadrightarrow T/x^n$ est nul.
\end{proof}

\begin{prop}\label{pr-noethvn} Soient $\C$ une petite catégorie et $F$ un foncteur de $\fct(\C,\E)$. On suppose que les valeurs de $F$ sont des objets noethériens de $\E$ et que $F/x$ est un foncteur noethérien. Alors $_{(x)}F$ est noethérien.
\end{prop}

\begin{proof}
Comme $F/x$ est noethérien, il existe $n\in\mathbb{N}$ tel que les sous-foncteurs $_{x^n}F$ et $_{(x)}F$ aient la même image dans $F/x$. Comme les valeurs de $F$ sont noethériennes, le lemme~\ref{lm-noestpn} montre que $_{(x)}F$ est un sous-foncteur de $F/x^n$, qui est noethérien par le lemme~\ref{lm-noethqxy}, d'où la conclusion.
\end{proof}

\begin{prop}\label{pr-noethfploc} Soient $\C$ une petite catégorie et $F$ un foncteur de $\fct(\C,\E)$ à valeurs dans les objets noethériens de $\E$. On suppose que $F/x$ est un objet noethérien de $\fct(\C,\E)$, et que $F[x^{-1}]$ est un objet noethérien de $\fct(\C,\E[x^{-1}])$. Alors $F$ est un objet noethérien de $\fct(\C,\E)$.
\end{prop}

\begin{proof} Soit $(G_n)_{n\in\mathbb{N}}$ une suite croissante de sous-foncteurs de $F$. Comme la suite $(G_n[x^{-1}])$ de sous-foncteurs de $F[x^{-1}]$ stationne, il existe $N\in\mathbb{N}$ tel que $G_i/G_N\subset \,_{(x)}F'$ pour tout $i\ge N$, où $F':=F/G_N$.

Comme $F'$ est un quotient de $F$, c'est un foncteur à valeurs noethériennes, et $F'/x$ est noethérien. Par conséquent, la proposition~\ref{pr-noethvn} montre que $_{(x)}F'$ est noethérien, donc $(G_i/G_N)_{i\ge N}$, et par conséquent $(G_i)$, stationne, d'où la proposition.
\end{proof}

\subsection{Cas de la propriété $pf_\infty$}

Le résultat suivant ne prétend pas à l'optimalité, mais il permet de montrer simplement la propriété $pf_\infty$ à partir d'hypothèses que nous vérifierons, pour des foncteurs appropriés, à la section~\ref{ssct-Aetf}. 

\begin{prop}\label{pr-pfi_chgtb} Supposons que l'anneau $k$ est intègre, que la catégorie de Grothendieck $\E$ est $k$-linéaire, et que $t$ est un élément non nul de $k$. Soient $X$ et $Y$ des objets de $\E$. On suppose que $X$ est de type fini, que $Y$ vérifie $pf_\infty$, et que $X$ et $Y$ sont sans $t$-torsion. On suppose de plus que $X[t^{-1}]$ est facteur direct de $Y[t^{-1}]$ dans $\E[t^{-1}]$, et que tout quotient de $X/t$ appartient à $\pf_\infty(\E)$.

Alors $X$ vérifie la propriété $pf_\infty$.
\end{prop}

\begin{proof}
Comme $X$ est de type fini et $Y$ sans $t$-torsion, le morphisme canonique $\E(X,Y)[t^{-1}]\to\E[t^{-1}](X[t^{-1}],Y[t^{-1}])$ est un isomorphisme ; de même $\E(Y,X)[t^{-1}]\xrightarrow{\simeq}\E[t^{-1}](Y[t^{-1}],X[t^{-1}])$. Par conséquent, le fait que $X[t^{-1}]$ est facteur direct de $Y[t^{-1}]$ dans $\E[t^{-1}]$ se traduit par l'existence de morphismes $u : X\to Y$ et $v : Y\to X$ de $\E$ tels quel $vu=t^i\mathrm{Id}_X$ pour un $i\in\mathbb{N}$.

Comme $X$ est sans $t$-torsion, on dispose de suites exactes
\begin{equation}\label{eq-sevrmc}
0\to X\xrightarrow{u} Y\to\mathrm{Coker}(u)\to 0\,,
\end{equation}
\begin{equation}\label{eq-seprc}
0\to\mathrm{Ker}(v)\to Y\xrightarrow{v} X\to\mathrm{Coker}(v)\to 0
\end{equation}
\begin{equation}\label{eqsepc}\text{et}\qquad 0\to\mathrm{Ker}(v)\to\mathrm{Coker}(u)\to X/t^i\to\mathrm{Coker}(v)\to 0\,.
\end{equation}
La suite exacte \eqref{eq-sevrmc}, combinée à \cite[prop.~1.5]{DT-schw}, montre l'équivalence
\begin{equation}\label{eqeqv1}
X\in\pf_m(\E)\Leftrightarrow\mathrm{Coker}(u)\in\pf_{m+1}(\E),
\end{equation}
pour tout $m\in\mathbb{N}$, puisque $Y\in\pf_\infty(\E)$.

Comme tout quotient de $X/t$ vérifie $pf_\infty$, on vérifie par récurrence sur $j\in\mathbb{N}^*$, toujours à partir de \cite[prop.~1.5]{DT-schw} (qu'on utilisera encore, sans le mentionner, dans la suite de cette démonstration), que tout sous-quotient de $X/t^j$ vérifie aussi $pf_\infty$, à partir des suites exactes $X/t^{j-1}\to X/t^j\to X/t\to 0$. Par conséquent, \eqref{eqsepc} montre la relation $\mathrm{Coker}(v)\in\pf_\infty(\E)$ et l'équivalence
\begin{equation}\label{eqeqv2}
\mathrm{Ker}(v)\in\pf_m(\E)\Leftrightarrow\mathrm{Coker}(u)\in\pf_m(\E)
\end{equation}
pour tout $m\in\mathbb{N}$.

Comme $Y$ et $\mathrm{Coker}(v)$ vérifient $pf_\infty$, \eqref{eq-seprc} entraîne l'équivalence
\begin{equation}\label{eqeqv3}
\mathrm{Ker}(v)\in\pf_m(\E)\Leftrightarrow X\in\pf_{m+1}(\E)
\end{equation}
pour tout $m\in\mathbb{N}$. Combinant \eqref{eqeqv1}, \eqref{eqeqv2} et \eqref{eqeqv3}, on obtient l'équivalence $X\in\pf_m(\E)\Leftrightarrow X\in\pf_{m+2}(\E)$ pour tout $m\in\mathbb{N}$. Comme $X$ est de type fini (i.e. $pf_0$), on en déduit qu'il vérifie $pf_\infty$ comme souhaité.
\end{proof}

\section{Cas d'un anneau $A$ sans quotient fini}

Le théorème~\ref{th-asqfc} se généralise directement à la catégorie $\F(\A;k)$ :

\begin{theo}\label{th-asqf}
Si l'anneau $A$ est sans quotient fini, alors tout foncteur $F$ de type fini et à valeurs de type fini de $\F(\A;k)$ est polynomial.
\end{theo}

\begin{proof} Par récurrence noethérienne, on peut supposer que, pour tout idéal non nul $I$ de $k$, tout foncteur de type fini et à valeurs de type fini de $\F(\A;k/I)$ est polynomial.

Si $k$ n'est pas intègre, soient $x, y\in k\setminus\{0\}$ tels que $xy=0$ : les foncteurs $F/x$ et $F/y$ sont polynomiaux par l'hypothèse de récurrence, de sorte que la suite exacte $F/x\to F/xy\to F/y\to 0$ donne la conclusion.

Si $k$ est intègre, le foncteur $F_K:=F\otimes_k K$, où $K$ est le corps des fractions de $k$, de $\F(\A;K)$ est polynomial d'après le théorème~\ref{th-asqfc}. Il existe donc un entier $d\in\mathbb{N}$ tel que le morphisme canonique $\delta_{d+1}^* cr_{d+1}(F)\to F$, dont le conoyau est $\qpol_d(F)$ (cf. \eqref{eq-qpol}, page~\pageref{eq-qpol}), devienne nul après tensorisation par $K$. Autrement dit, son image $G$, qui est de type fini comme $\delta_{d+1}^* cr_{d+1}(F)$, est à valeurs de torsion sur $k$. Il s'ensuit qu'il existe $x\in k\setminus\{0\}$ tel que $G$ soit annulé par $x$ ; ainsi, $G$ définit un foncteur de type fini, et à valeurs de type fini comme $F$, de $\F(\A;k/(x))$. L'hypothèse de récurrence montre que $G$ est polynomial, de sorte que la suite exacte  $0\to G\to F\to\qpol_d(F)\to 0$ permet de conclure.
\end{proof}

Le théorème~\ref{th-PAP_df} (hors de la propriété localement finie, sauf si elle est satisfaite dans $k\Md$, i.e. si $k$ est artinien) s'étend également à la catégorie $\F(A,k)$ ; le montrer nous demandera toutefois un petit peu de préparation.

Comme dans le cas où $k$ est un corps, on dit que l'anneau $A$ est \emph{$k$-trivial} s'il est fini et que $A\otimes_\mathbb{Z}k=0$, c'est-à-dire que le cardinal de $A$ est inversible dans $k$. Il revient au même de demander que les seuls foncteurs polynomiaux de $\F(A,k)$ soient les foncteurs constants \cite[prop.~2.13]{DTV}. Un idéal $I$ de $A$ est dit \emph{$k$-cotrivial}\index{termin}{cotrivial \emph{(idéal)}} si l'anneau quotient $A/I$ est $k$-trivial. Les foncteurs \emph{antipolynomiaux}\index{termin}{antipolynomial \emph{(foncteur)}} de $\F(A,k)$ sont, comme lorsque $k$ est un corps, les foncteurs se factorisant à travers la réduction modulo un idéal $k$-cotrivial.

 Dans ce qui suit, on identifie, comme on l'a déjà fait de nombreuses fois, la catégorie $\F(A/I\times A,k)$ à $\fct(\mathbf{P}(A/I)\times\mathbf{P}(A);k)$.

\begin{prdef}\label{pr-ktriv_ann} Soient $I$ un idéal $k$-cotrivial de $A$ et $d\in\mathbb{N}$. La restriction du foncteur $\F(A/I\times A,k)\to\F(A,k)$ de précomposition par le changement de base le long du morphisme d'anneaux canonique $A\to A/I\times A$ à la sous-catégorie des foncteurs polynomiaux de degré au plus $d$ par rapport à la première variable est pleinement fidèle. De plus, son image essentielle est une sous-catégorie bilocalisante de $\F(A,k)$.

On notera $\F_{I,d}(A,k)$ cette sous-catégorie et $\qpol_{I,d}^k : \F(A,k)\to\F_{I,d}(A,k)$ l'adjoint à gauche de l'inclusion.
\end{prdef}

\begin{proof}
Le fait que ce foncteur soit pleinement fidèle et ait une image stable par sous-quotient est démontré, dans le cas où $k$ est un corps, dans \cite[prop.~4.9]{DTV}. La démonstration (\cite[§\,4.3]{DTV}) s'étend sans changement dans le cas où $k$ est un anneau quelconque.

La stabilité par extensions de $\F_{I,d}(A,k)$ découle de \cite[th.~4.4 et~5.4]{DT-ext}. Sa stabilité par limites et colimites résulte de la même propriété pour $\pol_d(A,k)$ et de ce qu'un foncteur de précomposition commmute aux limites et colimites.
\end{proof}

 Comme dans le cas où $k$ est un corps, on dira qu'un foncteur $F$ de $\F(A,k)$ possède une décomposition à la Steinberg de type $PAP$ s'il existe un bifoncteur $B$ de $\F(A\times A,k)$ polynomial par rapport à la première variable et antipolynomial par rapport à la deuxième tel que $F\simeq\delta^*B$. Si $F$ est de type fini, il revient au même de demander qu'il existe un idéal $k$-cotrivial $I$ et un entier $d\ge 0$ tels que $F$ appartienne à la sous-catégorie $\F_{I,d}(A,k)$ (cf. \cite[lm~4.6]{DT-ext}).
 
 \begin{lemm}\label{lm-adjqId} Soient $I$ un idéal $k$-cotrivial de $A$, $d\in\mathbb{N}$ et $k'$ une $k$-algèbre. Pour tout foncteur $F$ de $\F(A,k)$, on dispose dans $\F(A,k')$ d'un isomorphisme naturel $\qpol_{I,d}^{k'}(F\otimes_k k')\simeq\qpol_{I,d}^k(F)\otimes_k k'$.
 \end{lemm}

\begin{proof} La postcomposition $Res : \F(A,k')\to\F(A,k)$ par la foncteur d'oubli et son adjoint à gauche $-\otimes_k k'$ préservent les sous-catégories $\F_{I,d}$ : on dispose d'un diagramme commutatif (à isomorphisme près)
$$\xymatrix{\F_{I,d}(A,k')\ar@{^{(}->}[r]\ar[d] & \F(A,k')\ar[d]\\
\F_{I,d}(A,k)\ar@{^{(}->}[r] & \F(A,k)
}$$
où les flèches verticales sont données par $Res$, d'où en prenant les adjoints à gauche un diagramme commutatif (à isomorphisme près)
$$\xymatrix{\F_{I,d}(A,k') & \F(A,k')\ar[l]_-{\qpol_{I,d}^{k'}}\\
\F_{I,d}(A,k)\ar[u] & \F(A,k)\ar[u]\ar[l]_-{\qpol_{I,d}^k} 
}$$
où les flèches verticales sont données par $-\otimes_k k'$.
\end{proof}

\begin{rema}\label{rq-qpId_form}
On vérifie aussitôt que l'on dispose de la \guillemotleft~formule explicite~\guillemotright\ suivante pour $\qpol_{I,d}^k$ : si $E$ est une partie génératrice de l'idéal $I$, $\qpol_{I,d}^k(F)$ est la précomposition par la diagonale du bifoncteur donné sur un objet $(U,V)$ par le conoyau du morphisme
$$\delta_{d+1}^*cr_{d+1}(F(U\oplus -))(V)\oplus\bigoplus_{a\in E}F(U\oplus U\oplus V)\to F(U\oplus V)$$
dont la première composante est donnée par la transformation naturelle canonique $\delta_{d+1}^*cr_{d+1}\to\mathrm{Id}$ (cf. \eqref{eq-qpol}) et la composante $F(U\oplus U\oplus V)\to F(U\oplus V)$ d'indice $a\in E$ est $F(p_a)-F(p_0)$, où $p_a : U\oplus U\oplus V\to U\oplus V$ envoie $(r,s,t)$ sur $(r+as,t)$.
\end{rema}

Nous pouvons maintenant donner une première version, faible (l'hypothèse de type fini est renforcée en une hypothèse noethérienne, pour une raison mentionnée à la remarque ci-dessous~\ref{rq-cotrivTF}), de la généralisation annoncée du théorème~\ref{th-PAP_df}. On renvoie au corollaire~\ref{cor-pqfe_ann} ci-après pour un énoncé plus fort.

\begin{theo}\label{th-pqfe_ann}
Supposons que tout idéal cofini de $A$ est $k$-cotrivial. Alors tout foncteur noethérien et à valeurs de type fini $F$ de $\F(A,k)$ possède une décomposition à la Steinberg de type $PAP$.
\end{theo}

\begin{proof} On suit le même schéma de démonstration que pour le théorème~\ref{th-asqf} : par récurrence noethérienne, on peut supposer que, pour tout idéal non nul $I$ de $k$, tout foncteur noethérien et à valeurs de type fini de $\F(\A;k/I)$ admet une décomposition à la Steinberg de type $PAP$. Si $k$ n'est pas intègre, soient $x, y\in k\setminus\{0\}$ tels que $xy=0$ : les foncteurs $F/x$ et $F/y$ possèdent une décomposition à la Steinberg de type $PAP$ par l'hypothèse de récurrence, de sorte que la suite exacte canonique $F/x\to F\to F/y\to 0$ et la proposition~\ref{pr-ktriv_ann} montrent la même propriété pour $F$.

Si $k$ est intègre, soit $K$ son corps des fractions : le théorème~\ref{th-PAP_df} montre que le foncteur de type fini $F\otimes_k K$ de $\F^\df(A,K)$ possède une décomposition de type $PAP$, donc qu'il existe un idéal $K$-cotrivial $I$ de $A$ et $d\in\mathbb{N}$ tels que $F\otimes_k K\in\F_{I,d}(A,K)$. Autrement dit, l'épimorphisme canonique $F\otimes_k K\twoheadrightarrow\qpol_{I,d}^K(F\otimes_k K)$ est un isomorphisme. Comme $I$ est un idéal cofini de $A$, c'est aussi un idéal $k$-cotrivial, de sorte que le lemme~\ref{lm-adjqId} montre que le noyau $G$ de l'épimorphisme canonique $F\twoheadrightarrow\qpol_{I,d}^k(F)$ vérifie $G\otimes_k K=0$. Autrement dit, $G$ est à valeurs de torsion sur $k$. Or $G$ est noethérien, comme $F$, donc il existe $x\ne 0$ dans $k$ tel que  les valeurs de $G$ soient annulées par $x$. On peut ainsi voir $G$ comme un foncteur noethérien à valeurs de type fini de $\F(A,k/(x))$, et l'hypothèse de récurrence montre que $G$ possède une décomposition à la Steinberg de type $PAP$ (dans $\F(A,k/(x))$ ou $\F(A,k)$, cela ne change rien). La proposition~\ref{pr-ktriv_ann} permet alors de conclure.
\end{proof}

\begin{rema}\label{rq-cotrivTF}
Si l'on fait l'hypothèse supplémentaire que tout idéal cofini de $A$ est de type fini (par exemple, que $A$ est noethérien), alors la description explicite des $\qpol_{I,d}^k$ donnée à la remarque~\ref{rq-qpId_form} montre que le noyau de l'épimorphisme canonique $F\twoheadrightarrow\qpol_{I,d}^k(F)$ est de type fini si $F$ est de type fini, de sorte que, comme dans la démonstration du théorème~\ref{th-asqf}, on peut obtenir la conclusion du théorème~\ref{th-pqfe_ann} en supposant seulement $F$ de type fini. Sans hypothèse supplémentaire sur $A$, il ne semble toutefois pas possible de procéder directement de la même manière ; nous résoudrons le problème au corollaire~\ref{cor-pqfe_ann} en utilisant un résultat de noethérianité général pour les foncteurs de type fini à valeurs de type fini de $\F(A,k)$.
\end{rema}

\section{Cas d'un anneau $A$ de type fini}\label{ssct-Aetf}

\subsection{Contre-exemple à la décomposition de type $FP$}

Le corollaire~\ref{cor-imq}, et donc a fortiori le théorème~\ref{th-tfdf_gal}, ne subsistent généralement pas lorsque le corps $K$ au but est remplacé par un anneau noethérien, même celui des entiers, comme nous allons dans le voir dans l'exemple~\ref{cex_dec-prim} ci-dessous.

 Nous commençons par une observation générale sur les décompositions de type $FP$ (la définition~\ref{def-decFP} s'étend de façon évidente à une catégorie but quelconque au lieu des $K$-espaces vectoriels).

\begin{rema}\label{rq-psi_dI} Soient $I$ un idéal de $A$ et $d\in\mathbb{N}$. Notons
$$\varphi_{I,d} : \fct(\mathbf{P}(A/I),\pol_d(A,k))\to\F(A,k)$$
le foncteur composé de l'inclusion
$$\fct(\mathbf{P}(A/I),\pol_d(A,k))\hookrightarrow\fct(\mathbf{P}(A/I),\F(A,k))\simeq\F(A/I\times A,k)$$
et de la précomposition par le foncteur canonique $\mathbf{P}(A)\to\mathbf{P}(A/I)\times\mathbf{P}(A)$.

Par adjonction somme/diagonale, le foncteur $\varphi_{I,d}$ admet un adjoint à droite $\psi_{I,d}$ donné comme suit : pour un foncteur $F$ de $\F(A,k)$, l'image de $\psi_{I,d}(F)$ par l'inclusion canonique $\fct(\mathbf{P}(A/I),\pol_d(A,k))\hookrightarrow\fct(\mathbf{P}(A/I),\F(A,k))\hookrightarrow\fct(\mathbf{P}(A),\F(A,k))\simeq\F(A\times A,k)$ est le plus grand sous-foncteur de $F\circ\oplus$ qui se factorise à travers la réduction modulo $I$ par rapport à la première variable et est polynomial de degré au plus $d$ par rapport à la seconde (sous-foncteur qui peut s'exprimer comme le noyau d'un morphisme explicite, à partir d'une partie génératrice de l'idéal $I$,  de façon duale de la formule pour $\qpol_{I,d}^k$ donnée à la remarque~\ref{rq-qpId_form}).
\end{rema}

\begin{exem}\label{cex_dec-prim} Dans la catégorie $\F(\mathbb{Z},\mathbb{Z})$, considérons le foncteur $T : V\mapsto\mathbb{Z}[V/2]\otimes V$ (où le produit tensoriel, comme tous ceux considérés dans cet exemple, est pris sur $\mathbb{Z}$) et le foncteur $F$ image du morphisme $P^{\mathbb{Z}}\to T$ donné sur $V$ par la composée de la diagonale $\mathbb{Z}[V]\to\mathbb{Z}[V]\otimes\mathbb{Z}[V]$ et du produit tensoriel des épimorphismes canoniques $\mathbb{Z}[V]\twoheadrightarrow\mathbb{Z}[V/2]$ et $\mathbb{Z}[V]\twoheadrightarrow V$. Alors $T$ et $F$ sont des foncteurs de type fini de $\F(\mathbb{Z},\mathbb{Z})$, à valeurs dans les groupes abéliens libres de rang fini. Le foncteur $T$ possède une décomposition de type $FP$, mais nous allons voir que son sous-foncteur $F$ n'en possède pas.

Tout d'abord, on note que la restriction à $F$ de l'épimorphisme $\alpha : T\twoheadrightarrow\pi^*\Lambda^2_{\FF_2}$ (où $\pi : \mathbf{P}(\mathbb{Z})\to\mathbf{P}(\FF_2)$ désigne la réduction modulo $2$ et $\Lambda^2$ la deuxième puissance extérieure) donné sur $V$ par la composée des deux projections canoniques
$$\mathbb{Z}[V/2]\otimes V\twoheadrightarrow (V/2)\otimes (V/2)\twoheadrightarrow\Lambda^2_{\FF_2}(V/2)$$
est nulle.

Soient $N>0$ et $d\ge 2$ des entiers. On vérifie aussitôt, en utilisant la remarque précédente, que $\psi_{(N),d}(\pi^*\Lambda^2_{\FF_2})(U)(V)$ est donné par
\begin{equation}\label{eqpsi_ev}\psi_{(N),d}(\pi^*\Lambda^2_{\FF_2})(U)(V)\simeq\Lambda^2_{\FF_2}(U/2)\oplus\Lambda^2_{\FF_2}(V/2)\oplus (U\otimes V)/2
\end{equation}
pour $U\in\mathbf{P}(\mathbb{Z}/N)$ et $V\in\mathbf{P}(\mathbb{Z})$ si $N$ est pair, ce qu'on suppose désormais ; $\psi_{(N),d}(T)$ est quant à lui donné par
\begin{equation}\label{eqpsi_T}
\psi_{(N),d}(T)(U)(V)\simeq\mathbb{Z}[U/2]\otimes V\,.
\end{equation}
Pour le voir, on écrit
$$T(V\oplus W)\simeq\big(\mathbb{Z}[W/2]\otimes (V\oplus W)\big)\oplus \big(\bar{\mathbb{Z}}[V/2]\otimes\mathbb{Z}[W/2]\otimes (V\oplus W)\big)$$
et l'on note tout d'abord que, pour tout $W$, le foncteur $V\mapsto\bar{\mathbb{Z}}[V/2]\otimes\mathbb{Z}[W/2]\otimes (V\oplus W)$ ne possède aucun sous-foncteur polynomial non nul (comme, plus généralement, tout foncteur de la forme $\bar{\mathbb{Z}}[\Phi]\otimes X$, où $\Phi$ est un foncteur additif de $\F(\mathbb{Z},\mathbb{Z})$ et $X$ un foncteur quelconque de cette catégorie --- il s'agit d'une conséquence classique de la suite exacte \eqref{eq-ppol}, page~\pageref{eq-ppol}). Ensuite, on on constate que, pour tout $V$, le foncteur $G : W\mapsto\mathbb{Z}[W/2]\otimes W$ ne contient aucun sous-foncteur non nul se factorisant à travers la réduction modulo $N$ : grâce à la proposition~\ref{pr-radevtau}, il suffit de vérifier que 
$$G(W)\xrightarrow{G\left(\begin{array}{c} 1 \\
                                                N 
                                               \end{array}\right)-G\left(\begin{array}{c} 1 \\
                                                0 
                                               \end{array}\right)} G(W\oplus W)$$
                                               est injectif, ce qui est immédiat ($W\in\mathbf{P}(\mathbb{Z})$ est un groupe abélien sans torsion).
                                               
Les différentes observations qui précèdent établissent \eqref{eqpsi_T} (en utilisant la remarque~\ref{rq-psi_dI}). De plus, on vérifie aussitôt que la coünité $\varphi_{(N),d}\psi_{(N),d}(T)\to T$ est un isomorphisme et que l'épimorphisme $\alpha : T\twoheadrightarrow\pi^*\Lambda^2_{\FF_2}$ induit, via les isomorphismes \eqref{eqpsi_ev} et \eqref{eqpsi_T},
$$\psi_{(N),d}(\alpha)(U)(V) : \mathbb{Z}[U/2]\otimes V\twoheadrightarrow (U/2)\otimes V\hookrightarrow\Lambda^2_{\FF_2}(U/2)\oplus\Lambda^2_{\FF_2}(V/2)\oplus (U\otimes V)/2\,.$$

Comme la composée $F\hookrightarrow T\xrightarrow{\alpha}\pi^*\Lambda^2_{\FF_2}$ est nulle, il s'ensuit que l'image de $\psi_{(N),d}(F\hookrightarrow T)$ est incluse dans $X\boxtimes\mathrm{Id}$, où $X(U):=\mathrm{Ker}\,(\mathbb{Z}[U/2]\twoheadrightarrow U/2)$ et $\mathrm{Id}$ désigne ici, par abus, le foncteur d'inclusion $\mathbf{P}(\mathbb{Z})\hookrightarrow\mathbf{Ab}$. Par conséquent, l'image de la composée de la coünité $\varphi_{(N),d}\psi_{(N),d}(F)\to F$ et de l'inclusion $F\hookrightarrow T$ est incluse dans $\pi^*X\otimes\mathrm{Id}$, qui est manifestement \emph{strictement} inclus dans $F$ en tant que sous-foncteur de $T$ (l'élément $[1]$ de $T(\mathbb{Z})\simeq\mathbb{Z}[\mathbb{Z}/2]$ appartient à $F(\mathbb{Z})$ mais pas à $(\pi^*X\otimes\mathrm{Id})(\mathbb{Z})$). Ainsi, la coünité $\varphi_{(N),d}\psi_{(N),d}(F)\to F$ n'est pas un épimorphisme, ce qui entraîne que $F$ n'appartient pas à l'image essentielle de $\varphi_{(N),d}$. Comme ceci est vrai pour tout entier $d\ge 2$ et tout entier pair $N>0$, $F$ n'admet aucune décomposition de type $FP$.
\end{exem}

\subsection{Retour sur le spectre, l'algèbre caractéristique et les poids}

Nous avons défini les décompositions en poids, l'algèbre caractéristique et le spectre d'un foncteur de $\fct(\A,\E)$ pour une petite catégorie additive $A$-linéaire $\A$ et une catégorie de Grothendieck $K$-linéaire $\E$. Ces définitions (et d'autres reliées, comme celle de $\ti$ et de $\hi$) s'étendent sans changement lorsque $\E$ est une catégorie de Grothendieck $k$-linéaire, où $k$ est un anneau quelconque. (Toutefois, les propriétés fondamentales des décompositions en poids \emph{ne} s'étendent \emph{pas} sans lourde modification à cette situation, c'est pourquoi nous ne les avons étudiées que sur un corps de base $K$ jusqu'alors.)

Indiquons sommairement les quelques notions dont nous aurons besoin dans cette section~\ref{ssct-Aetf}. Si $F$ est un foncteur de $\fct(\A,\E)$, où $\E$ est $k$-linéaire, son spectre\index{termin}{spectre \emph{(d'un foncteur)}} est un quotient du foncteur $P^{A;k}(=k[-])$ de $\F(A,k)$ ; on le notera $\Sp(F)$,\index{nota}{sp@$\Sp$, $\Sp_A$ \emph{(spectre d'un foncteur)}} ou $\Sp^k(F)$ s'il y a besoin de spécifier l'anneau $k$ (nous n'aurons pas besoin de préciser l'anneau $A$). On dispose d'isomorphismes canoniques $F\simeq\Sp(F)\ti F$ et $\hi(\Sp(F),F)\simeq F$ (induits par la projection $P^{A;k}\twoheadrightarrow\Sp(F)$), et $\Sp(F)$ est le plus petit quotient de $P^{A;k}$ vérifiant l'une ou l'autre de ces propriétés. L'algèbre caractéristique\index{termin}{algebre caracteristique@algèbre caractéristique} de $F$ (relativement à $A$ et $k$) sera notée $\ac_A^k(F)$\index{nota}{A@$\ac$, $\ac_A$ \emph{(algèbre caractéristique d'un foncteur)}} (l'indice $A$ ou l'exposant $k$ seront omis si aucune confusion ne peut en résulter) : c'est le quotient de la $k$-algèbre $k[A_\mu]$ par l'idéal des $\sum_{i=1}^n\lambda_i[a_i]$ (où $\lambda_i\in k$ et $a_i\in A$) tels que $\sum_{i=1}^n\lambda_i F(1\oplus a_i)=0$.

En matière de poids, nous aurons seulement besoin de généraliser une partie de la définition~\ref{defipoi} et la notation~\ref{not-Tspc} :
\begin{nota}\label{not-Tspa}\begin{enumerate}
\item Si $F$ est un foncteur de $\F(A,k)$, $r\ge 0$ un entier et $w : A\to k$ un poids (c'est-à-dire un morphisme de monoïdes multiplicatifs), on note $F^{[r]}_w$ le plus grand quotient de $F$ appartenant à $U(\F(A,k);\mathfrak{m}_w^r)$, où $\mathfrak{m}_w:=\mathrm{Ker}\,(k[A_\mu]\to k\qquad [a]\mapsto w(a))$.
\item Pour $d,r,n\in\mathbb{N}$, $w_1,\dots,w_n : A\to k$ des poids et $I$ un idéal de $A$, on définit dans $\F(A,k)$
$$\mathrm{T}_{A,I,d,r,(w_1,\dots,w_n)}^k:=\pi_I^*P^{A/I;k}\otimes_k\Big(\bigoplus_{i=1}^n\qpol_d(P^{A;k})_{w_i}^{[r]}\Big)\,,$$
où $\pi_I$ désigne la réduction modulo $I$.
\end{enumerate}
\end{nota}

\begin{lemm}\label{lm-techn_TJM} Supposons que l'anneau $A$ est de type fini. Dans la situation de la notation précédente, le foncteur $\mathrm{T}_{A,I,d,r,(w_1,\dots,w_n)}^k$ est de présentation finie dans $\F(A,k)$.
\end{lemm}

\begin{proof}
Comme l'anneau $A$ est noethérien, \cite[prop.~7.1]{DT-schw} montre que $\pi_I^*P^{A/I;k}$ est de présentation finie (et vérifie même $pf_\infty$) dans $\F(A,k)$. Le foncteur $\bigoplus_{i=1}^n\qpol_d(P^{A;k})_{w_i}^{[r]}$ est polynomial de type fini dans $\F(A,k)$ ; comme $A$ est de type fini et $k$ noethérien, \cite[th.~2]{DT-TJM} montre que ce foncteur est de présentation finie (et vérifie même $pf_\infty$). Le produit tensoriel sur $k$ préserve les générateurs projectifs standards de $\F(A,k)$, donc les foncteurs de présentation finie, de $\F(A,k)$, d'où le lemme.
\end{proof}

\paragraph*{Extension des scalaires au but} Soit $k'$ une $k$-algèbre. On dispose de certaines compatibilités évidentes des constructions considérées précédemment à l'extension des scalaires au but le long de $k\to k'$ : on a par exemple des isomorphismes canoniques $P^{A;k'}\simeq P^{A;k}\otimes_k k'$ et, dans la situation de la notation~\ref{not-Tspa}, $\mathrm{T}_{A,I,d,r,(\tilde{w}_1,\dots,\tilde{w}_n)}^{k'}\simeq \mathrm{T}_{A,I,d,r,(w_1,\dots,w_n)}^k\otimes_k k'$, où $\tilde{w_i}$ désigne la composée de $w_i$ et de $k\to k'$.

Le foncteur $\ti : \F(A,k)\times\fct(\A,\E)\to\fct(\A,\E)$ commute également à l'extension des scalaires au but : si $T$ est un foncteur de $\F(A,k)$ et $F$ un foncteur de $\fct(\A,\E)$ (où $\E$ est toujours une catégorie de Grothendieck $k$-linéaire), on a dans $\fct(\A,\E\underset{k}{\boxtimes}k')$ un isomorphisme canonique
\begin{equation}\label{eq-ti_cbb}
(T\otimes_k k')\ti (F\underset{k}{\boxtimes}k')\simeq (T\ti F)\underset{k}{\boxtimes}k'\,.
\end{equation}

En revanche, le foncteur $\hi$, comme le foncteur Hom usuel, qu'il généralise, ne commute à l'extension des scalaires que sous certaines hypothèses. Nous aurons seulement besoin à ce propos de l'énoncé suivant.

\begin{lemm}\label{lm-pfhi} Supposons que la $k$-algèbre $k'$ est \textbf{plate}.
Soient $T$ un foncteur de $\F(A,k)$, $F$ et $G$ des foncteurs de $\fct(\A,\E)$.

On dispose dans $\fct(\A,\E\underset{k}{\boxtimes}k')$ d'un morphisme naturel
$$\hi(T,F)\underset{k}{\boxtimes}k'\to\hi(T\otimes_k k',F\underset{k}{\boxtimes}k')$$
qui est un isomorphisme si $T$ est de présentation finie (dans $\F(A,k)$).
\end{lemm}

\begin{proof}
Le résultat se déduit de la définition de $\hi$, de l'adjonction entre restriction et extension des scalaires, du fait qu'un $k$-module plat est colimite filtrante de $k$-modules projectifs de type fini et de ce que $\operatorname{Hom}(X,-)$ commute aux colimites filtrantes
lorsque $X$ est de présentation finie.
\end{proof}

\subsection{Lemmes de localisation}\label{par-an_integre}

\begin{conv} Dans tout le §\,\ref{par-an_integre}, on suppose l'anneau $k$ \textbf{intègre}. On désigne par $K$ son corps des fractions.
\end{conv}

Nous commençons par deux lemmes simples d'algèbre commutative. Le premier, particulièrement classique, est laissé en exercice.

%trouver une réf., voire se passer de ce petit lemme ?
\begin{lemm}\label{lm-modfloc} Soit $M$ un module de type fini sur un anneau intègre. Il existe un élément non nul $x$ de ce dernier tel que $M[x^{-1}]$ soit un module libre de rang fini sur l'anneau localisé.
\end{lemm}

\begin{lemm}\label{lm-loc_atf} Supposons que l'anneau $A$ est de type fini. Soient $L$ une extension algébrique du corps $K$ et $\alpha_1,\dots,\alpha_r : A\to L$ une collection finie de morphismes d'anneaux.

Alors il existe $x\in k\setminus\{0\}$, $n\in\mathbb{N}^*$ et un sous-anneau $k'$ de $L$ contenant $k[x^{-1}]$ tels que les $\alpha_i$ soient à valeurs dans $k'$ et que $k'$ soit un $k[x^{-1}]$-module libre de rang $n$.
\end{lemm}

\begin{proof}
Considérant l'image par les $\alpha_i$ d'un ensemble fini de générateurs de l'anneau $A$, on obtient un sous-ensemble fini $E$ de $L$ tel que la condition que les $\alpha_i$ soient à valeurs dans $k'$ équivaille à ce que $k'$ contienne $E$. Comme $L$ est algébrique sur $K$ et $E$ fini, il existe $t\in k\setminus\{0\}$ tel que chaque élément élément de $E$ soit annulé par un polynôme unitaire à coefficients dans $k[t^{-1}]$. Il s'ensuit que la sous-$k[t^{-1}]$-algèbre $R$ de $L$ engendrée par $E$ est un module de type fini sur $k[t^{-1}]$. La conclusion résulte donc du lemme~\ref{lm-modfloc}.
\end{proof}

\begin{lemm}\label{lm-epi_cbint} Si $u : F\to G$ est un morphisme de $\F(\A;k)$, avec $G$ de type fini, et $L$ un surcorps de $K$ tels que $u\otimes_k L$ soit un épimorphisme de $\F(\A;L)$, alors il existe $x\in k\setminus\{0\}$ tel que $u\otimes_k k[x^{-1}]$ soit un épimorphisme de $\F(\A;k[x^{-1}])$.
\end{lemm}

\begin{proof} Comme $u\otimes_k L$ est un épimorphisme, il en est de même pour $u\otimes_k K$ dans $\F(\A;K)$ (puisque $K\subset L$ est fidèlement plat), de sorte que $\mathrm{Coker}(u)$ est à valeurs de torsion sur $k$. Or ce foncteur est de type fini, car c'est un quotient de $G$, il existe donc $x\in k\setminus\{0\}$ tel que les valeurs de $\mathrm{Coker}(u)$ soient annulées par $x$. Ainsi $u\otimes_k k[x^{-1}]$ est-il un épimorphisme.
\end{proof}

Nous aurons besoin de la généralisation suivante de la proposition~\ref{pr-spht}.

\begin{prop}\label{pr-spht2} Soient $d,r,n\in\mathbb{N}$, $w_1,\dots,w_n$ des éléments deux à deux distincts de $\mon(A_\mu,k_\mu)$ s'exprimant chacun comme produit de morphismes d'anneaux d'images infinies et $I$ un idéal cofini de $A$. Soit $u : P^{A;k}\to\mathrm{T}^k_{A,I,d,r,(w_1,\dots,w_n)}$
le morphisme de $\F(A,k)$ composé de la diagonale $P^{A;k}\to P^{A;k}\otimes_k P^{A;k}$ et du morphisme canonique $P^{A;k}\to\pi_I^*P^{A/I;k}$ tensorisé par le morphisme $P^{A;k}\to\bigoplus_{i=1}^n\qpol_d(P^{A;k})_{w_i}^{[r]}$ dont les composantes sont les projections.

Alors il existe $x\in k\setminus\{0\}$ tel que $u\otimes_k k[x^{-1}]$ soit un épimorphisme de $\F(A,k[x^{-1}])$.
\end{prop}

\begin{proof}
Comme le foncteur $\mathrm{T}^k_{A,I,d,r,(w_1,\dots,w_n)}$ de $\F(A,k)$ est de type fini (en tant que quotient d'une somme directe finie de copies de $P^{A;k}\otimes_k P^{A;k}\simeq P^{A^2;k}$), le lemme~\ref{lm-epi_cbint} montre qu'il suffit de voir que $u\otimes_k K$ est un épimorphisme. C'est précisément ce que nous apprend la proposition~\ref{pr-spht}.
\end{proof}

Nous pouvons désormais donner le résultat principal de ce paragraphe, qui constitue l'étape clef de la récurrence noethérienne qui démontrera le théorème fondamental~\ref{thf-Aesstf} ci-après.

\begin{prop}\label{pr-factaidempres} Supposons que l'anneau $A$ est de type fini. Soit $F$ un foncteur de type fini de $\F(A,k)$ dont les valeurs sont des $k$-modules de type fini. 

Alors il existe un élément non nul $x$ de $k$, des entiers naturels $r$, $n$, $d$, une $k[x^{-1}]$-algèbre $k'$, libre de rang fini non nul sur $k[x^{-1}]$, des morphismes de monoïdes $w_1,\dots,w_n : A_\mu\to k'_\mu$, un idéal cofini $I$ de $A$ et foncteur $B$ de $\F(A/I\times A,k)$ vérifiant les propriétés suivantes :
\begin{enumerate}
\item le morphisme $P^{A;k'}\to\mathrm{T}^{k'}_{A,I,d,r,(w_1,\dots,w_n)}$ de la proposition~\ref{pr-spht2} est un épimorphisme, et $\Sp^{k'}(F\otimes_k k')\le\mathrm{T}^{k'}_{A,I,d,r,(w_1,\dots,w_n)}$ ;
\item $B$ est de type fini et à valeurs de type fini, polynomial de degré au plus $d$ par rapport à la deuxième variable ;
\item $F[x^{-1}]$ est facteur direct dans $\F(A,k[x^{-1}])$ de la composée de $B[x^{-1}]$ et du changement de base à la source le long du morphisme d'anneaux canonique $A\to A/I\times A$.
\end{enumerate}
\end{prop}

\begin{proof}
Il existe un extension finie de corps $K\subset L$ tel que le foncteur $F\otimes_k L$ appartienne à $\Tt^\cf(A,L)$ (proposition~\ref{pr-tfdf_htr-ext}). Comme $A$ est noethérien, et vérifie donc en particulier la condition (IMF), le corollaire~\ref{cor-hpapfor} fournit un idéal cofini $I$ de $A$, des entiers naturels $n$, $r$ et $d$ et des poids deux à deux distincts $w_1,\dots,w_n : A\to L$ de $\mathfrak{M}_{\mathrm{poli}}$ tels que $\Sp^L(F\otimes_k L)\le\mathrm{T}^L_{A,I,d,r,(w_1,\dots,w_n)}$.

Le lemme~\ref{lm-loc_atf} fournit $x\in k\setminus\{0\}$ et une sous-$k[x^{-1}]$-algèbre $k'$ de $L$, libre de rang fini non nul sur $k[x^{-1}]$, tels que les $w_i$ prennent leurs valeurs dans $k'$ (on notera, par abus, de la même façon les $w_i$ et les poids $A\to k'$ correspondants). Quitte à remplacer $x$ par $xt$ pour un élément $t$ approprié de $k\setminus\{0\}$, et $k'$ par $k'[t^{-1}]$, on peut supposer que le morphisme canonique $P^{A;k'}\to\mathrm{T}^{k'}_{A,I,d,r,(w_1,\dots,w_n)}$ est un épimorphisme, par la proposition~\ref{pr-spht2}.

Notons $F'$ le foncteur $F\otimes_k k'$ de $\F(A,k')$. On dispose dans $\F(A,L)$ d'isomorphismes canoniques
$$\hi(\mathrm{T}^{k'}_{A,I,d,r,(w_1,\dots,w_n)},F')\otimes_{k'} L\simeq\hi(\mathrm{T}^{k'}_{A,I,d,r,(w_1,\dots,w_n)}\otimes_{k'} L,F'\otimes_{k'} L)\simeq\cdots$$
$$\hi(\mathrm{T}^L_{A,I,d,r,(w_1,\dots,w_n)},F\otimes_k L)\simeq F\otimes_k L\simeq F'\otimes_{k'} L$$
dont le premier découle des lemmes~\ref{lm-pfhi} et~\ref{lm-techn_TJM} (le morphisme d'anneaux $k'\to L$ est plat car il est injectif et $L$ est un corps) et l'avant-dernier de l'inégalité $\Sp^L(F\otimes_k L)\le\mathrm{T}^L_{A,I,d,r,(w_1,\dots,w_n)}$. Comme cet isomorphisme est induit par le monomorphisme canonique $\hi(\mathrm{T}^{k'}_{A,I,d,r,(w_1,\dots,w_n)},F')\hookrightarrow F'$ et que $F'$ est de type fini dans $\F(A,k')$ (puisque $F$ est de type fini dans $\F(A,k)$), le lemme~\ref{lm-epi_cbint} fournit un élément non nul $t$ de $k$ tel que le monomorphisme canonique $\hi(\mathrm{T}^{k'}_{A,I,d,r,(w_1,\dots,w_n)},F')[t^{-1}]\hookrightarrow F'[t^{-1}]$ de $\F(A,k'[t^{-1}])$ soit un isomorphisme. En utilisant de nouveau le lemme~\ref{lm-pfhi}, on en déduit l'inégalité $\Sp^{k'[t^{-1}]}(F'[t^{-1}])\le\mathrm{T}^{k'[t^{-1}]}_{A,I,d,r,(w_1,\dots,w_n)}$. Ainsi, quitte à remplacer $x$ par $xt$ et $k'$ par $k'[t^{-1}]$, on peut supposer que $\Sp^{k'}(F')\le\mathrm{T}^{k'}_{A,I,d,r,(w_1,\dots,w_n)}$. Cela établit la première assertion.

Notons à présent $G$ le foncteur $\big(\pi_I^*P^{A/I;k}\otimes_k\qpol_d(P^{A;k})\big)\ti F$ de $\F(A,k)$. Si $B$ désigne le foncteur de $\F(A/I\times A,k)$ donné par
$$B(U,V)=\operatorname{Coend}\Big(\big(\pi_I^*P^{U;k}_{\mathbf{P}(A/I)^\op}\boxtimes_k\qpol_d\big(P^{V;k}_{\mathbf{P}(A)^\op}\big)\big)\boxtimes_k (F\circ\oplus) :$$
$$(\mathbf{P}(A)\times\mathbf{P}(A))^\op\times (\mathbf{P}(A)\times\mathbf{P}(A))\to k\Md\Big)\,,$$
on vérifie facilement que $B$ est polynomial de degré au plus $d$ par rapport à la deuxième variable et, en utilisant l'adjonction somme/diagonale, que $G$ est isomorphe à la composée de $B$ et du morphisme canonique $\mathbf{P}(A)\to\mathbf{P}(A/I)\times\mathbf{P}(A)$. En outre, l'épimorphisme canonique
$$P^{U;k}_{\mathbf{P}(A)^\op}\boxtimes_k P^{V;k}_{\mathbf{P}(A)^\op}\twoheadrightarrow \pi_I^*P^{\pi_I(U);k}_{\mathbf{P}(A/I)^\op}\boxtimes_k\qpol_d\big(P^{V;k}_{\mathbf{P}(A)^\op}\big)$$
naturel en les objets et $V$ de $\mathbf{P}(A)$ induit un épimorphisme naturel
$$F(U\oplus V)\twoheadrightarrow B(\pi_I(U),V)$$
qui montre que $B\circ (\pi_I\times\operatorname{Id}_{\mathbf{P}(A)})$ est comme $F$, de type fini et à valeurs de type fini sur $k$. Comme le foncteur $\pi_I$ est essentiellement surjectif, cela entraîne les mêmes propriétés pour $B$, d'où la deuxième assertion.

Quitte à remplacer $B$ par $B^{\oplus m}$, il suffit maintenant pour conclure de démontrer que $F[x^{-1}]$ est facteur direct de $G[x^{-1}]^{\oplus m}$, où $m$ est le rang de $k'$ sur $k[x^{-1}]$.

Comme la projection $P^{A;k'}\twoheadrightarrow\mathrm{T}^{k'}_{A,I,d,r,(w_1,\dots,w_n)}$ se factorise par le morphisme $P^{A;k'}\to\pi_I^*P^{A/I;k'}\otimes_k\qpol_d(P^{A;k'})$ composé de la diagonale $P^{A;k'}\to P^{A;k'}\otimes_{k'} P^{A;k'}$ et du produit tensoriel des projections canoniques, appliquant le foncteur $-\ti F'$, on obtient une composée
$$F'\simeq P^{A;k'}\ti F'\to\big(\pi_I^*P^{A/I;k'}\otimes_k\qpol_d(P^{A;k'})\big)\ti F'\to\mathrm{T}^{k'}_{A,I,d,r,(w_1,\dots,w_n)}\ti F'\simeq F'$$
égale à l'identité, en vertu de la première assertion. Or
$$\big(\pi_I^*P^{A/I;k'}\otimes_k\qpol_d(P^{A;k'})\big)\ti F'\simeq\big(\big(\pi_I^*P^{A/I;k}\otimes_k\qpol_d(P^{A;k})\big)\otimes_k k'\big)\ti\big(F\otimes_k k'\big)$$
est isomorphe à $G\otimes_k k'$ par \eqref{eq-ti_cbb} (page~\pageref{eq-ti_cbb}). Ainsi, $F'$ est facteur direct de $G\otimes_k k'$. Comme $k'\simeq k[x^{-1}]^m$ comme $k[x^{-1}]$-module, on conclut en postcomposant par le foncteur de restriction des scalaires $k'\Md\to k[x^{-1}]\Md$.
\end{proof}

\begin{coro}\label{cor-pfinoe_aux}  Supposons que l'anneau $A$ est de type fini. Soit $F$ un foncteur de type fini de $\F(A,k)$ dont les valeurs sont des $k$-modules de type fini. 

Alors il existe un élément non nul $x$ de $k$ tel que le foncteur $F[x^{-1}]$ soit noethérien dans $\F(A,k[x^{-1}])$ et un foncteur $G$ de $\pf_\infty(\F(A,k))$, à valeurs sans $x$-torsion, tel que $F[x^{-1}]$ soit facteur direct de $G[x^{-1}]$.
\end{coro}

\begin{proof} On choisit pour $x$ l'élément non nul de $k$ fournit par la proposition~\ref{pr-factaidempres} (dont on conserve aussi les notations de la conclusion pour $k'$, $I$, $B$, etc.).

Comme $k'$ est un $k[x^{-1}]$-module libre de rang fini non nul, $F[x^{-1}]$ est noethérien dans $\F(A,k[x^{-1}])$ si et seulement si $F'=F\otimes_k k'$ est noethérien dans $\F(A,k')$.

La sous-catégorie prébilocalisante $\F_{\mathrm{T}^{k'}_{A,I,d,r,(w_1,\dots,w_n)}}(A,k')$ de $\F(A,k')$ constituée des foncteurs $X$ tels que $\Sp^{k'}(X)\le\mathrm{T}^{k'}_{A,I,d,r,(w_1,\dots,w_n)}$ est équivalente à la sous-catégorie prébilocalisante $\fct(\mathbf{P}(A/I),\F_{\bigoplus_{i=1}^n\qpol_d(P^{A;k})_{w_i}^{[r]}}(A,k'))$ de $\fct(\mathbf{P}(A/I),\F(A,k'))\simeq\F(A/I\times A,k')$ (en effet, le corollaire~\ref{cor-sp_pt} s'étend, avec la même démonstration, au cas où le corps $K$ au but est remplacé par un anneau).
Comme $\F_{\bigoplus_{i=1}^n\qpol_d(P^{A;k})_{w_i}^{[r]}}(A,k')$ est une sous-catégorie prébilocalisante de $\pol_d(A,k)$, qui est localement noethérienne \cite[th.~1]{DT-TJM}, et que $A/I$ est fini, le théorème~\ref{th-PSS} permet de conclure pour la noethérianité.

Soit $G$ la composée du foncteur $B':=B/_{(x)}B$ de $\F(A/I\times A,k)$ et du foncteur canonique $\mathbf{P}(A)\to\mathbf{P}(A/I)\times\mathbf{P}(A)$ : $G$ est à valeurs sans $x$-torsion, et $F[x^{-1}]$ est facteur direct de $B[x^{-1}]\simeq B'[x^{-1}]$. De surcroît, $B'$ est un foncteur de type fini à valeurs dans la sous-catégorie $\fct(\mathbf{P}(A/I),\pol_d(A,k))$ de $\fct(\mathbf{P}(A/I),\F(A,k))\simeq\F(A/I\times A,k)$. Comme la catégorie $\fct(\mathbf{P}(A/I),\pol_d(A,k))$ est localement noethérienne (cf. supra), $B'$ appartient à $\pf_\infty(\fct(\mathbf{P}(A/I),\pol_d(A,k)))$. Or \cite[th.~2]{DT-TJM} montre que le foncteur d'inclusion $\pol_d(A,k)\hookrightarrow\F(A,k)$ préserve la propriété $pf_\infty$. Le lemme~\ref{lm-pfn_but} ci-dessous permet de conclure que $B'$ appartient à $\pf_\infty(\F(A/I\times A,k))$. Enfin, comme $A$ est un anneau noethérien, le foncteur $\mathbf{Add}(\mathbf{P}(A),\mathbf{Ab})\simeq A\Md\to(A/I)\Md\times A\Md\simeq\mathbf{Add}(\mathbf{P}(A/I)\times\mathbf{P}(A),\mathbf{Ab})$ induit par le foncteur canonique $\mathbf{P}(A)\to\mathbf{P}(A/I)\times\mathbf{P}(A)$ préserve la propriété $pf_\infty$ (qui équivaut à la type-finitude à la source comme au but, qui sont localement noethériens), de sorte que \cite[cor.~7.3]{DT-schw} permet d'en déduire que $G$ appartient à $\pf_\infty(\F(A,k))$, d'où le résultat.
\end{proof}

Le lemme formel suivant complète les résultats de \cite{DT-schw} utiles dans cette section.
% placer ce lemme en tout début d'article ?
\begin{lemm}\label{lm-pfn_but} Soient $\C$ une petite catégorie, $\Phi : \E\to\E'$ un foncteur exact entre catégories de Grothendieck et $n\in\mathbb{N}\cup\{\infty\}$. On suppose que :
\begin{enumerate}
\item la catégorie $\E$ est engendrée par des objets projectifs de type fini ;
\item le foncteur $\Phi$ préserve la propriété $pf_n$.\index{nota}{pf@$\pf_n$}
\end{enumerate}

Alors le foncteur de postcomposition $\Phi_* : \fct(\C,\E)\to\fct(\C,\E')$ préserve la propriété $pf_n$.
\end{lemm}

\begin{proof}
Si $P$ est un objet projectif de type fini de $\E$ et $c$ un objet de $\C$, alors le foncteur $P[\C(c,-)]$ de $\fct(\C,\E)$ est projectif de type fini, par le lemme de Yoneda. De plus, les foncteurs de ce type engendrent la catégorie $\fct(\C,\E)$. Il suffit donc \cite[prop.~1.9\,(1)]{DT-schw} de montrer que $\Phi_*(P[\C(c,-)])=\Phi(P)[\C(c,-)]$ vérifie la propriété $pf_n$ dans $\fct(\C,\E')$. Comme $\Phi(P)$ vérifie par hypothèse $pf_n$ dans $\E'$, cela découle encore du lemme de Yoneda, qui fournit des isomorphismes naturels $\mathrm{Ext}^*_{\fct(\C,\E')}(\Phi(P)[\C(c,-)],T)\simeq\mathrm{Ext}^*_{\E'}(\Phi(P),T(c))$.
\end{proof}

\subsection{Résultat fondamental}

Nous sommes maintenant en mesure d'établir le résultat principal de cette section~\ref{ssct-Aetf} :

\begin{theo}\label{thf-Aesstf}
Si l'anneau $A$ est de type fini, alors tout foncteur $F$ de $\F(A,k)$, de type fini et à valeurs de type fini, est noethérien et vérifie la propriété $pf_\infty$.\index{nota}{pf@$\pf_n$}
\end{theo}

\begin{proof} Par récurrence noethérienne, on peut supposer que, pour tout idéal non nul $I$ de $k$, tout foncteur de type fini et à valeurs de type fini de $\F(A,k/I)$ est noethérien et possède la propriété $pf_\infty$.

Si $k$ n'est pas intègre, soient $x, y\in k\setminus\{0\}$ tels que $xy=0$ : on dispose d'une suite exacte courte $0\to G\to F\to F\otimes k/(y)\to 0$ où $G$ (resp. $F\otimes k/(y)$) est un foncteur appartenant à l'image dans $\F(A,k)$ d'un foncteur de $\F(A,k/(x))$ (resp. $\F(A,k/(y))$, de type fini et à valeurs de type fini. Par l'hypothèse de récurrence, $G$ (resp. $F\otimes k/(y)$) est noethérien et vérifie $pf_\infty$ dans $\F(A,k/(x))$ (resp. $\F(A,k/(y))$, ce qui implique qu'ils vérifient les mêmes propriétés dans $\F(A,k)$ : pour la propriété noethérienne, cela résulte de ce que $\F(A,k/(x))$ est une sous-catégorie abélienne pleine stable par sous-objet de $\F(\A,k)$, et pour la propriété $pf_\infty$, cela provient du lemme~\ref{lm-pfn_but}, puisque le foncteur $k/(x)\Md\to k\Md$ induit par la projection $k\twoheadrightarrow k/(x)$ la préserve (car $k$ est un anneau noethérien). Il s'ensuit que $F$ est noethérien et appartient à $\pf_\infty(\F(A,k))$.

Supposons maintenant $k$ intègre : le corollaire~\ref{cor-pfinoe_aux} fournit $x\in k\setminus\{0\}$ tel que le foncteur $F[x^{-1}]$ soit noethérien dans $\F(A,k[x^{-1}])$ et un foncteur $G$ de $\pf_\infty(\F(A,k))$, à valeurs sans $x$-torsion, tel que $F[x^{-1}]$ soit facteur direct de $G[x^{-1}]$.

L'hypothèse de récurrence montre que $F/x$ est noethérien. La proposition~\ref{pr-noethfploc} permet donc de conclure que $F$ est noethérien.

Montrons que $F$ vérifie $pf_\infty$ d'abord lorsque $F$ est sans $x$-torsion. Comme l'hypothèse de récurrence entraîne que tout quotient de $F/x$ vérifie $pf_\infty$ dans $\F(A,k/(x))$, donc dans $\F(A,k)$ (cf. supra), cela découle de la proposition~\ref{pr-pfi_chgtb}.

Dans le cas général, comme $F$ est noethérien d'après ce qui précède, son sous-foncteur de $x$-torsion $_{(x)}F$ est de type fini, donc annulé par $x^n$ pour $n\in\mathbb{N}$ assez grand. L'hypothèse de récurrence montre ainsi que $_{(x)}F$ vérifie $pf_\infty$ dans $\F(A,k/(x^n))$, donc dans $\F(A,k)$. Comme $F/_{(x)}F$ appartient à $\pf_\infty(\F(A,k))$ d'après le cas précédent, il en est de même pour $F$, d'où le théorème.
\end{proof}

Nous verrons bientôt, au théorème~\ref{th-noethgal}, que l'on peut obtenir la propriété noethérienne pour les foncteurs de type fini et à valeurs de type fini de $\F(A,k)$ sans aucune hypothèse sur $A$. Ce n'est pas le cas pour la propriété $pf_\infty$ (même si $A$ et $k$ sont des corps --- cf. remarque~\ref{rq-pfiag}) ; on peut toutefois affaiblir légèrement l'hypothèse que $A$ est un anneau de type fini :

\begin{coro}\label{cor-pfi_an_gal}
Si l'anneau $A$ est essentiellement de type fini, alors tout foncteur de type fini et à valeurs de type fini $F$ de $\F(A,k)$ vérifie la propriété $pf_\infty$.\index{nota}{pf@$\pf_n$}
\end{coro}

\begin{proof}
Soit $A'$ un anneau de type fini tel que $A\simeq A'[S^{-1}]$ pour une certaine partie multiplicative $S$ de $A'$. La proposition~\ref{pr-loc_fcttftf} montre que l'image de $F$ dans $\F(A',k)$ par la précomposition $\iota$ par le changement de base le long du morphisme canonique $A'\to A$ est de type fini. Il s'ensuit $\iota(F)\in\pf_\infty(\F(A',k))$, par le théorème~\ref{thf-Aesstf}. L'adjoint à gauche $\Phi$ à $\iota$ est \emph{exact} (il est en effet donné par une colimite \emph{filtrante} associée à la partie multiplicative $S$) ; comme $\iota$ est exact et cocontinu, on en déduit par \cite[prop.~1.10]{DT-schw} que $\Phi$ préserve la propriété $pf_\infty$. Ainsi $F\simeq\Phi(\iota(F))\in\pf_\infty(\F(A,k))$.
\end{proof}

%On peut même gagner un chouilla en plus, en remplaçant $A$ par une algèbre non commutative finie sur un anneau commutatif essentiellement de type fini (ou plus généralement une algèbre non commutative de type fini comme module du bon côté sur un sous-anneau commutatif essentiellement de type fini) en appliquant \cite[prop.~1.9\,(2)]{DT-schw}, mais il n'est sans doute pas utile de le dire.

\begin{rema} Nous ne savons pas s'il est possible d'affaiblir l'hypothèse que $A$ est un anneau essentiellement de type fini en celle que les anneaux $A$ et $A\otimes_\mathbb{Z}k$ sont noethériens, comme c'est le cas lorsque $k$ est un corps (théorème~\ref{th-pfi_corps}).
\end{rema}

\section{Cas général}\label{ssct-noeth}

Le résultat suivant constitue une généralisation du corollaire~\ref{cor-acdf} qui se démontre exactement comme ce dernier.

\begin{lemm}\label{lm-acdf2} Soit $F$ un foncteur de type fini de $\F(A,k)$ tel que $\ac_A^k(F)$ soit une $k$-algèbre de type fini. Il existe un sous-anneau de type fini $A'$ de $A$ tel que l'image de $F$ dans $\F(A',k)$ (par la précomposition par l'extension des scalaires $\mathbf{P}(A')\to\mathbf{P}(A)$) est de type fini.
\end{lemm}

Le lemme précédent permet aussitôt d'étendre le théorème~\ref{thf-Aesstf} à un anneau quelconque à la source :

\begin{theo}\label{th-noethgal}
Tout foncteur de type fini et à valeurs de type fini de $\F(A,k)$ est noethérien.
\end{theo}

\begin{proof}
On se ramène au cas où $A$ est un anneau de type fini, traité par le théorème~\ref{thf-Aesstf}, à partir du lemme~\ref{lm-acdf2}, de la même façon que dans la démonstration du théorème~\ref{th-noeth_cor_gal}.
\end{proof}

On en déduit, de la même manière qu'au corollaire~\ref{cor-noeth_corg}, la légère généralisation suivante :

\begin{coro}\label{cor-noeth_anng}
Soit $R$ un anneau non commutatif. On suppose que $R$ est de type fini comme module à droite sur un sous-anneau commutatif $A$. Alors tout foncteur de type fini et à valeurs de type fini de $\F(R,k)$ est noethérien.
\end{coro}

En combinant les théorèmes \ref{th-noethgal} et \ref{th-pqfe_ann}, on obtient aussitôt :

\begin{coro}\label{cor-pqfe_ann}
Supposons que tout idéal cofini de $A$ est $k$-cotrivial. Alors tout foncteur de type fini et à valeurs de type fini de $\F(A,k)$ possède une décomposition à la Steinberg de type $PAP$.
\end{coro}

\appendix

\chapter{Action d'un groupe fini sur une catégorie abélienne}\label{app_act-grp}

La construction classique du produit semi-direct d'une catégorie abélienne $\E$ (de Grothendieck ici) par un groupe $G$ (fini ici) agissant sur celle-ci intervient dans la théorie des foncteurs polynomiaux par l'intermédiaire du théorème~\ref{th-recPiranv} de Pirashvili. Après des rappels généraux sur cette construction, cet appendice présente deux résultats (qui ne prétendent pas à l'originalité mais dont nous ne connaissons pas de référence dans la littérature sous la forme exacte dont nous avons besoin) qui interviendront à deux endroits de ce mémoire :
\begin{itemize}
    \item la préservation de la noethérianité par le foncteur d'oubli $\E\rtimes G\to\E$ (proposition~\ref{pr-onpt}), indispensable pour démontrer les théorèmes principaux du chapitre~\ref{sec-finpolp} ;
    \item la finitude de la dimension des corps (non commutatifs) d'endomorphismes d'objets simples de $\E\rtimes G$ sur leur centre, sous l'hypothèse que la même condition vaut dans $\E$, indispensable pour démontrer le théorème~\ref{th-chipol} et son important corollaire~\ref{cor-EMLpol-pol} (la proposition~\ref{pr-endofcentr} permet d'éviter toute hypothèse de valeurs de dimensions finies).
\end{itemize}

\section{Généralités}

\begin{hyp} Dans tout cet appendice, $\E$ désigne une catégorie de Grothendieck munie d'une action d'un groupe \textbf{fini} $G$, c'est-à-dire d'une famille $(T_g)_{g\in G}$ d'endofoncteurs de $\E$ telle que $T_1=\mathrm{Id}_\E$ et $T_{gh}=T_g.T_h$ pour tous $g, h\in G$.
\end{hyp}

Ainsi, les $T_g$ sont des automorphismes de $\E$, donc des foncteurs continus et cocontinus, et en particulier exacts. 

\begin{defi}\label{def-psdc}
On appelle produit semi-direct de $\E$ par $G$, et l'on note $\E\rtimes G$ la catégorie dont les objets sont les objets $X$ de $\E$ munis d'isomorphismes $\alpha_g : X\xrightarrow{\simeq}T_g(X)$ pour tout $g\in G$ de sorte que, pour tous $g, h\in G$, le diagramme suivant commute :
$$\xymatrix{X\ar[r]^-{\alpha_g}\ar[d]_-{\alpha_{gh}} & T_g(X)\ar[d]^-{T_g(\alpha_h)}\\
T_{gh}(X)\ar@{=}[r] & T_g(T_h(X))
}.$$
Un morphisme $(X,(\alpha_g)_{g\in G})\to (Y,(\beta_g)_{g\in G})$ de $\E\rtimes G$ est un morphisme $f : X\to Y$ tel que, pour tout $g\in G$, le diagramme
$$\xymatrix{X\ar[r]^-{\alpha_g}\ar[d]_f  & T_g(X)\ar[d]^-{T_g(f)}\\
Y\ar[r]^-{\beta_g} & T_g(Y)
}$$
commute. La composition des morphismes est induite par celle de $\E$.
\end{defi}

L'exemple typique suivant (les résultats de cet appendice constituent des généralisations directes de propriétés des anneaux de groupe tordus) justifie la terminologie et la notation.

\begin{exem}\label{ex-agtC}
Si $R$ est un anneau (non nécessairement commutatif) muni d'une action d'un groupe $G$ (par automorphismes d'anneaux), alors $G$ opère sur la catégorie $R\Md$, via les foncteurs $T_g$ de restriction des scalaires le long de l'automorphisme de $R$ correspondant à l'action de $g^{-1}$. On vérifie aussitôt qu'on dispose d'une équivalence de catégories canonique $(R\Md)\rtimes G\simeq R\langle G\rangle\Md$, où $R\langle G\rangle$ désigne l'{\em anneau de groupe tordu} de $G$ sur $R$, dont le groupe additif sous-jacent est $R[G]$ et la multiplication est donnée par $(r[g]).(s[h]):=(r.g_*(s))[gh]$ (où $g_*(s)$ désigne l'effet de l'action de $g$ sur $s$). Si $R=A[H]$ ($A$-algèbre de groupe usuelle), où $H$ est un groupe muni d'une action de $G$, alors $R\langle G\rangle=A[H\rtimes G]$.
\end{exem}

\begin{rema}\label{rq-wrca}
Plus généralement, si $k$ est un anneau et $R$ une $k$-algèbre (non nécessairement commutative), toute action d'un groupe $G$ sur $R$ (par automorphismes de $k$-algèbres) sur $R$ induit une action de $G$ sur $\E\underset{k}{\boxtimes}R$, si la catégorie de Grothendieck $\E$ est $k$-linéaire, et l'on dispose d'une équivalence de catégories canonique $$\big(\E\underset{k}{\boxtimes}R\big)\rtimes G\simeq\E\underset{k}{\boxtimes}R\langle G\rangle\,.$$
\end{rema}

\begin{nota}
On note $\mathcal{O} : \E\rtimes G\to\E\quad (X,(\alpha_g)_{g\in G})\mapsto X$ le foncteur d'oubli.
\end{nota}

\begin{rema}
\begin{itemize}
    \item Par définition de $\E\rtimes G$, on dispose d'un isomorphisme canonique $T_g\circ\mathcal{O}\simeq\mathcal{O}$ pour tout $g\in G$.
    \item Dans la situation de l'exemple~\ref{ex-agtC}, $\mathcal{O}$ s'identifie à la restriction des scalaires le long du morphisme d'anneaux canonique $R\to R\langle G\rangle$.
    \item Si $X$ et $Y$ sont des objets de $\E\rtimes G$, le groupe abélien $\E(\mathcal{O}(X),\mathcal{O}(Y))$ est muni d'une action naturelle de $G$, et $(\E\rtimes G)(X,Y)$ s'identifie au sous-groupe $\E(\mathcal{O}(X),\mathcal{O}(Y))^G$ des invariants.
\end{itemize}
\end{rema}

\begin{prop}
La catégorie $\E\rtimes G$ est une catégorie de Grothendieck ; le foncteur $\mathcal{O}$ est exact et fidèle.
\end{prop}

La proposition précédente, facile et classique, est laissée en exercice. Il en est de même pour l'énoncé suivant (cf. \cite[lemme~2.18 et le paragraphe le précédant]{DTV}).

\begin{prdef}\label{prdf-adjOL}
On définit un foncteur $\mathcal{L} : \E\to\E\rtimes G$ par $\mathcal{L}(X)=\Big(\bigoplus_{a\in G}T_a(X),(\alpha_g)_{g\in G}\Big)$, où $\alpha_g$ est l'isomorphisme canonique
$$\alpha_g : \bigoplus_{a\in G}T_a(X)\to T_g\Big(\bigoplus_{a\in G}T_a(X)\Big)\simeq\bigoplus_{a\in G}T_{ga}(X)$$
induit par la permutation $a\mapsto ga$ de $G$ ; l'effet sur les morphismes de $\mathcal{L}$ est donné par la fonctorialité des $T_a$.

Le foncteur $\mathcal{L}$ est adjoint des deux côtés à $\mathcal{O}$. En particulier, $\mathcal{L}$ et $\mathcal{O}$ commutent aux limites et aux colimites.
\end{prdef}

\section{Préservation de propriétés de finitude}

Le but principal de cette section consiste à établir le résultat suivant, qui constitue une généralisation directe de \cite[thm~4]{Grz85} (qui traite le cas d'une catégorie de modules sur un anneau de groupe tordu --- cf. exemple~\ref{ex-agtC}). Toutefois, comme nous ne connaissons pas de référence dans le cas général, nous en donnons la démonstration complète.

\begin{prop}\label{pr-onpt} Le foncteur d'oubli $\mathcal{O} : \E\rtimes G\to\E$ préserve les objets noethériens.
\end{prop}

Pour ce faire, on commence par montrer, comme \cite{Grz85}, le résultat suivant :

\begin{lemm}\label{lm-essentiel} Soient $X$ un objet non nul de $\E\rtimes G$ et $U$ un sous-objet \textbf{essentiel} de $\mathcal{O}(X)$. Alors il existe un sous-objet non nul $Y$ de $X$ tel que $\mathcal{O}(Y)\subset U$.
\end{lemm}

\begin{proof} Pour tout sous-objet $V$ de $\mathcal{O}(X)$ et tout $g\in G$, notons $_gV$ l'image de
$$T_g(V)\xrightarrow{T_g(V\hookrightarrow\mathcal{O}(X))}T_g(\mathcal{O}(X))\simeq\mathcal{O}(X).$$
On note que $_g(_hV)=_{gh}V$ pour tout $h\in G$, et que $V$ définit un sous-objet de $X$ dans $\E\rtimes G$ si et seulement si $_gV=V$ pour tout $g\in G$.

Par ailleurs, pour tout $g\in G$, du fait que $T_g$ est un automorphisme de $\E$, $U_g$ est comme $U$ un sous-objet essentiel de $\mathcal{O}(X)$. Il s'ensuit que $Y:=\bigcap_{g\in G}U_g$ est aussi un sous-objet essentiel, et en particulier non nul, de $\mathcal{O}(X)$, puisque $G$ est fini. Il résulte maintenant des observations précédentes que $Y$ définit un sous-objet de $X$ dans $\E\rtimes G$.
\end{proof}

\begin{lemm}\label{lm-otfpt} Le foncteur d'oubli $\mathcal{O} : \E\rtimes G\to\E$ préserve les objets de type fini.
\end{lemm}

\begin{proof}
En effet, l'adjoint à droite $\mathcal{L}$ de $\mathcal{O}$ commute aux colimites (cf. proposition-définition~\ref{prdf-adjOL}).
\end{proof}

\begin{proof}[Démonstration de la proposition~\ref{pr-onpt}] Soit $X$ un objet noethérien de $\E\rtimes G$. Pour montrer que $\mathcal{O}(X)$ est noethérien, on peut supposer que $X$ est non nul et, par récurrence noethérienne, que $\mathcal{O}(X/Y)$ est noethérien pour tout sous-objet non nul $Y$ de $X$.

Soit $(S_n)_{n\in\mathbb{N}}$ une suite croissante de sous-objets de $\mathcal{O}(X)$ et $S_\infty$ sa réunion. Comme $\E$ est une catégorie de Grothendieck, le lemme de Zorn montre qu'il existe un sous-objet \emph{maximal} $R$ de $\mathcal{O}(X)$ tel que $S_\infty\cap R=0$ ; $S_\infty\oplus R$ est alors un sous-objet \emph{essentiel} de $\mathcal{O}(X)$. Le lemme~\ref{lm-essentiel} fournit un sous-objet non nul $Y$ de $X$ tel que $\mathcal{O}(Y)\subset S_\infty\oplus R$. 

Comme $X$ est noethérien, $Y$ est de type fini, donc $\mathcal{O}(Y)$ est de type fini (lemme~\ref{lm-otfpt}). Il s'ensuit qu'il existe $N\in\mathbb{N}$ tel que $\mathcal{O}(Y)\subset S_N\oplus R$. L'image dans $\mathcal{O}(X)/\mathcal{O}(Y)\simeq\mathcal{O}(X/Y)$ de la suite croissante $(S_n\oplus R)_{n\ge N}$ de sous-objets de $\mathcal{O}(X)$ stationne, puisque $\mathcal{O}(X/Y)$ est noethérien. Par conséquent, $(S_n\oplus R)_{n\ge N}$, et donc $(S_n)_{n\in\mathbb{N}}$, stationnent.
\end{proof}

Nous montrons maintenant que le foncteur $\mathcal{O}$ préserve également les objets finis et les objets semi-simples (cf. \cite{Grz85}, ou \cite[§\,2.6]{DTV}). Ce résultat, plus facile que la préservation de la noethérianité, nous servira à la fin de cet appendice, pour établir la proposition~\ref{pr-endofcentr}, ainsi qu'au chapitre~\ref{sec-finpolp}.

\begin{lemm}\label{lm-Oss} Si $X$ est un objet simple de $\E\rtimes G$, alors $\mathcal{O}(X)$ est un objet semi-simple fini de $\E$.

Plus précisément, si $S$ désigne un facteur de composition de $\mathcal{O}(X)$ et que $\Gamma$ désigne le sous-groupe de $G$ des éléments $g$ tels que $T_g(S)\simeq S$, alors $T_g(S)$ ne dépend, à isomorphisme près, que de la classe $\bar{g}$ de $g$ dans $G/\Gamma$, et il existe un entier $n\ge 1$ tel que
$$\mathcal{O}(X)\simeq\bigoplus_{\bar{g}\in G/\Gamma}T_g(S)^{\oplus n}\;.$$
\end{lemm}

\begin{proof}
Comme $X$ est de type fini et non nul, il en est de même pour $\mathcal{O}(X)$, par le lemme~\ref{lm-otfpt}. Par conséquent, $\mathcal{O}(X)$ possède un quotient simple $S$. Comme $S$ est simple, l'adjoint $X\to\mathcal{L}(S)$ (cf. proposition-définition~\ref{prdfla}) de la projection $\mathcal{O}(X)\twoheadrightarrow S$ est un monomorphisme. Par conséquent, $\mathcal{O}(X)$ est un sous-objet de $\mathcal{O}\mathcal{L}(S)\simeq\bigoplus_{g\in G}T_g(S)$. Comme les foncteurs $T_g$ sont des automorphismes de la catégorie $\E$, les $T_g(S)$ sont des objets simples de $\E$. Il s'ensuit que $\mathcal{O}\mathcal{L}(S)$, et donc $\mathcal{O}(X)$, sont semi-simples finis. On obtient la deuxième partie de la conclusion en regroupant les facteurs simples isomorphes de $\mathcal{O}(X)$ et en utilisant l'isomorphisme $T_g(\mathcal{O}(X))\simeq\mathcal{O}(X)$ pour tout $g\in G$.
\end{proof}

\begin{prop}\label{pr-psspsd} Soit $X$ un objet de $\E\rtimes G$.
\begin{enumerate}
\item L'objet $\mathcal{O}(X)$ de $\E$ est semi-simple si $X$ est un objet semi-simple de $\E\rtimes G$.
\item Si $\E$ est $\mathbb{Q}$-linéaire et que $\mathcal{O}(X)$ de $\E$ est semi-simple fini, alors $X$ est semi-simple fini.
\item L'objet $\mathcal{O}(X)$ de $\E$ est fini si et seulement si $X$ est un objet fini de $\E\rtimes G$.
\end{enumerate}
\end{prop}

\begin{proof}
Comme $\mathcal{O}$ est un foncteur exact, le lemme~\ref{lm-Oss} entraîne que $\mathcal{O}$ préserve les objets finis ainsi que les objets semi-simples finis. Comme $\mathcal{O}$ commute également aux colimites arbitraires, il préserve les objets semi-simples.

Comme $\mathcal{O}$ est exact et fidèle, si $\mathcal{O}(X)$ est fini, $X$ est fini et de longueur inférieure ou égale à celle de $\mathcal{O}(X)$.

 Supposons maintenant que $\E$ est $\mathbb{Q}$-linéaire et $\mathcal{O}(X)$ est semi-simple fini : comme $X$ est fini, il suffit de montrer que tout sous-objet $Y$ de $X$ est facteur direct. C'est le théorème de Maschke : $\mathcal{O}(Y)$ est facteur direct de $\mathcal{O}(X)$ ; si $p\in\mathrm{Hom}_\E(\mathcal{O}(X),\mathcal{O}(Y))$ est une rétraction de l'inclusion $\mathcal{O}(Y)\hookrightarrow\mathcal{O}(X)$, alors $\frac{1}{\cd(G)}\sum_{g\in G}g.p$ en est une rétraction $G$-invariante, autrement dit elle définit une rétraction $X\twoheadrightarrow Y$ de l'inclusion, d'où la proposition.
\end{proof}

\section{Endomorphismes des objets simples}

\begin{nota}\label{nota-wrt} Soit $E$ un ensemble.
\begin{enumerate}
    \item On note $\Si(E)$ le groupe des permutations de $E$.
    \item Si $L$ est un anneau, on note $L^{\times E}$ le produit sur $E$ de copies de $L$, c'est-à-dire l'anneau des fonctions (ensemblistes) de $E$ dans $L$.
    \item Si $T$ est un groupe, on note $T\wr\Si(E):=T^E\rtimes\Si(E)$\index{nota}{$\wr$ \emph{(produit en couronne)}|textbf} le produit en couronne de $T$ par $\Si(E)$.
\end{enumerate}
\end{nota}

(La mention au symbole de produit dans le terme $L^{\times E}$ est destinée à le distinguer du sous-anneau $L^G$ des invariants de $L$ par l'action d'un groupe $G$.)

Le résultat suivant est un exercice classique qui découle de l'étude des idempotents de $L^{\times E}$, lorsque $L$ est corps. 

\begin{lemm}\label{lm-autprodc} Soient $L$ un corps et $E$ un ensemble fini. Le morphisme de groupes canoniques $\mathrm{Aut}_\mathbf{Ann}(L)\wr\Si(E)\to\mathrm{Aut}_\mathbf{Ann}(L^{\times E})$ est un isomorphisme.
\end{lemm}

\begin{lemm}\label{lm-invautprc} Soient $L$ un corps, $E$ un ensemble fini et $G$ un groupe fini opérant sur l'anneau $L^{\times E}$. Alors $L^{\times E}$ est un module de type fini sur son sous-anneau $(L^{\times E})^G$ d'invariants sous l'action de $G$.
\end{lemm}

\begin{proof} On ne perd pas en généralité à supposer que $G$ est un sous-groupe de $\mathrm{Aut}_\mathbf{Ann}(L^{\times E})$.
Utilisant l'isomorphisme du lemme~\ref{lm-autprodc}, notons $T$ l'image du morphisme de groupes $G\hookrightarrow\mathrm{Aut}_\mathbf{Ann}(L^{\times E})\simeq\mathrm{Aut}_\mathbf{Ann}(L)\wr\Si(E)\twoheadrightarrow\Si(E)$, et $H$ son noyau, c'est-à-dire l'intersection des sous-groupes $G$ et $\mathrm{Aut}_\mathbf{Ann}(L)^{\times E}$ de $\mathrm{Aut}_\mathbf{Ann}(L^{\times E})\simeq\mathrm{Aut}_\mathbf{Ann}(L)\wr\Si(E)$.

Pour $i\in E$, notons $H_i$ le sous-groupe de $\mathrm{Aut}_\mathbf{Ann}(L)$ image de la composée de l'inclusion $H\hookrightarrow\mathrm{Aut}_\mathbf{Ann}(L)^{\times E}$ et de la projection $\mathrm{Aut}_\mathbf{Ann}(L)^{\times E}\twoheadrightarrow\mathrm{Aut}_\mathbf{Ann}(L)$ donnée par l'évaluation en $i$ : $G$, donc $H$ et les $H_i$ sont des groupes finis, et
$$H\subset\prod_{i\in E}H_i\subset\mathrm{Aut}_\mathbf{Ann}(L)^{\times E}\;.$$
Soit $k_i$ le corps des invariants de $L$ sous l'action du groupe fini $H_i$ : $L$ est de dimension finie sur chaque $k_i$, et l'inclusion de groupes précédente fournit l'inclusion de sous-anneaux de $L^{\times E}$
$$\prod_{i\in E}k_i=\big(L^{\times E}\big)^{\prod_{i\in E}H_i}\subset (L^{\times E})^H\;.$$
L'action de $T\simeq G/H$ sur $(L^{\times E})^H$ induite par l'action de $G$ sur $L^{\times E}$ est la restriction à $T\subset\Si(E)$ de l'action de $\Si(E)$ sur $L^{\times E}$ par permutation des facteurs du produit direct.

Soit $\gamma : E/T\to E$ une section (ensembliste) de la surjection canonique, notée $i\mapsto\bar{i}$. La stabilité de $(L^{\times E})^H$ par l'action de $T$ et l'inclusion précédente fournissent
$$\prod_{i\in E}k_{\gamma(\bar{i})}\subset (L^{\times E})^H$$
puis
$$\prod_{x\in E/T}k_{\gamma(x)}\subset\big((L^{\times E})^H)^T=(L^{\times E})^G\;,$$
où $\prod_{x\in E/T}k_{\gamma(x)}$ est plongé dans $\prod_{i\in E}k_{\gamma(\bar{i})}$ via la précomposition par $E\twoheadrightarrow E/T$.

Comme $L$ est un $k_{\gamma(x)}$-espace vectoriel de dimension finie pour tout $x\in E/T$, et que $E/T$ est fini, $L^{\times E}$ est un module de type fini sur $\prod_{x\in E/T}k_{\gamma(x)}$, et donc a fortiori sur $(L^{\times E})^G$ d'après ce qui précède.
\end{proof}

\begin{prop}\label{pr-endofcentr} Soient $\E$ une catégorie de Grothendieck et $G$ un groupe fini opérant sur $\E$. Supposons que, pour tout objet $S$ de $\E$, le corps (non commutatif) $\mathrm{End}_\E(S)$ soit de dimension finie sur son centre.

Alors pour tout objet simple $X$ de $\E\rtimes G$, le corps (non commutatif) $\mathrm{End}_{\E\rtimes G}(X)$ est de dimension finie sur son centre.
\end{prop}

\begin{proof} Par le lemme~\ref{lm-Oss}, $\mathcal{O}(X)$ est semi-simple fini ; on se donne $S$, $n$ et $\Gamma$ comme dans la conclusion de ce lemme.

Puisque les $T_g(S)$ sont des objets simples deux à deux non isomorphes de $\E$, mais de corps (non commutatifs) d'endomorphismes isomorphes, lorsque $g$ parcourt $G/\Gamma$, on a un isomorphisme d'anneaux (non commutatifs)
$$\mathrm{End}_{\mathcal{O}(X)}\simeq\M_n(\mathrm{End}\,S)^{\times G/\Gamma}\;.$$
Le groupe $G$ opère (via la structure d'objet de $\E\rtimes G$ de $X$) sur cet anneau, donc sur son centre $Z(\mathrm{End}\,S)^{\times G/\Gamma}$, et
$$\mathrm{End}(X)\simeq (\M_n(\mathrm{End}\,S)^{\times G/\Gamma})^G\supset Z\big((\M_n(\mathrm{End}\,S)^{\times G/\Gamma})^G\big)\supset\big(Z(\mathrm{End}\,S)^{\times G/\Gamma}\big)^G\;.$$

Considérons le diagramme commutatif d'inclusions d'anneaux (non commutatifs)
$$\xymatrix{\big(Z(\mathrm{End}\,S)^{\times G/\Gamma}\big)^G\ar[r]\ar[d] & Z\big((\M_n(\mathrm{End}\,S)^{\times G/\Gamma})^G\big)\ar[d] \\
Z(\mathrm{End}\,S)^{\times G/\Gamma}\ar[r] & \M_n(\mathrm{End}\,S)^{\times G/\Gamma}
}\;:$$
$Z(\mathrm{End}\,S)^{\times G/\Gamma}$ est un module de type fini sur $\big(Z(\mathrm{End}\,S)^{\times G/\Gamma}\big)^G$ d'après le lemme~\ref{lm-invautprc}, et $\M_n(\mathrm{End}\,S)^{\times G/\Gamma}$ est un $Z(\mathrm{End}\,S)^{\times G/\Gamma}$-module de type fini parce que $\mathrm{End}\,S$ est par hypothèse un espace vectoriel de dimension finie sur $Z(\mathrm{End}\,S)$. Il s'ensuit que $\M_n(\mathrm{End}\,S)^{\times G/\Gamma}$ est un $\big(Z(\mathrm{End}\,S)^{\times G/\Gamma}\big)^G$-module de type fini, et a fortiori un $Z\big((\M_n(\mathrm{End}\,S)^{\times G/\Gamma})^G\big)\simeq Z(\mathrm{End}(X))$-espace vectoriel de dimension finie. Son sous-anneau $(\M_n(\mathrm{End}\,S)^{\times G/\Gamma})^G\simeq\mathrm{End}(X)$ est donc aussi un $Z(\mathrm{End}(X))$-espace vectoriel de dimension finie, comme souhaité.
\end{proof}

\chapter{Décompositions en poids et groupes d'extensions}\label{ch-apExt}

Nous examinons ici sommairement le comportement cohomologique de plusieurs des décompositions qui jouent un rôle important dans ce mémoire, à savoir les décompositions en poids et celles qui apparaissent dans nos théorèmes de structure. Nous ne faisons qu'effleurer le sujet, qui pose rapidement des problèmes difficiles dès lors qu'on manque de propriétés de finitude sur les foncteurs considérés.

Les méthodes mises en jeu dans cet appendice recoupent, ou reposent souvent sur \cite{DT-ext,Dja-JLMS}, dans un esprit qui rappelle la cohomologie locale \cite{Sch-Si} en algèbre commutative (notamment pour le théorème~\ref{th-comp_Ext-tordI}).

\section{Décomposition selon les poids partiels et groupes d'extensions}

Plusieurs des décompositions fondamentales qui apparaissent dans ce mémoire possèdent un bon comportement relativement aux groupes d'extensions. C'est déjà le cas pour les décompositions de type $PAP$, grâce à \cite[th.~4.4 et cor.~5.7]{DT-ext}. Nous ignorons toutefois si \cite[th.~4.4]{DT-ext} se généralise aux foncteurs ayant une décomposition de type $HPAP$.

L'assertion~\ref{ithp3} du théorème~\ref{th-decPPart} fournit également un résultat puissant de décomposition des groupes d'extensions, qui recoupe les résultats de \cite{DT-ext} qu'on vient de mentionner (et les généralisent aux décompositions de type $HPAP$ \textit{pour les foncteurs \hts}). On prendra toutefois garde qu'on ne peut pas en général s'affranchir de l'hypothèse de finitude sur le foncteur apparaissant au but de l'Ext dans le théorème~\ref{th-decPPart}, comme l'illustrera l'exemple~\ref{exextsdi}.

À l'aune de ces résultats, et de ceux de la section~\ref{ss-tshtr}, il apparaît très naturel d'étudier le comportement cohomologique des foncteurs d'inclusion $\Tt_{[\mathfrak{M}_{\mathrm{tord},I}]}(\A,\E)\to\fct(\A,\E)$, où $I$ est un idéal de $A$ tel que $A/I$ soit un corps fini de caractéristique $p$, et $\Tt_{[\mathfrak{M}_{\mathrm{poli},\alpha}]}(\A,\E)\to\fct(\A,\E)$, où $\alpha : A\to K$ est un morphisme d'anneaux d'image infinie. Le problème semble ardu en général ; nous nous contenterons ici de donner quelques conditions nécessaires et/ou suffisantes pour que ce soient des \textit{plongements homologiques}, i.e. qu'ils induisent des isomorphismes entre groupes d'extensions.

\begin{theo}\label{th-isoExt_Fpbar} Soit $\alpha : A\to K$ est un morphisme d'anneaux d'image infinie. Supposons que $\mathrm{Ext}^i_{A_K}(K_\alpha,K_\alpha)=0$ pour tout $i\in\mathbb{N}^*$, où $K_\alpha$ désigne le $(A,K)$-bimodule $K$ muni de l'action de $A$ déduite de $\alpha$. Si $p>0$, on suppose également que $A$ est une $\mathbb{F}_p$-algèbre et que le corps $K$ est parfait.

Alors l'inclusion $\Tt_{[\mathfrak{M}_{\mathrm{poli},\alpha}]}(A,K)\to\F(A,K)$ induit un isomorphisme
$$\mathrm{Ext}^*_{\Tt_{[\mathfrak{M}_{\mathrm{poli},\alpha}]}(A,K)}(F,G)\xrightarrow{\simeq}\mathrm{Ext}^*_{\F(A,K)}(F,G)$$
lorsque $F$ est un foncteur de $\Tt_{[\mathfrak{M}_{\mathrm{poli},\alpha}]}(A,K)$ et $G$ un foncteur de $\Tt^\cf_{[\mathfrak{M}_{\mathrm{poli},\alpha}]}(A,K)$.
\end{theo}

\begin{proof} Un argument formel de colimite montre qu'il suffit d'établir le résultat lorsque $F$ est \textit{de type fini}, ce qu'on suppose désormais.

L'annulation de $\mathrm{Ext}^1_{A_K}(K_\alpha,K_\alpha)\simeq\mathrm{Der}_\alpha(A,K)$ entraîne que tout foncteur de $\Tt_{[\mathfrak{M}_{\mathrm{poli},\alpha}]}(A,K)$ est \hd. En effet, la proposition~\ref{pr-poli_alpha} (et la remarque qui la suit) permet de se ramener au cas où $A$ est local, de sorte que la condition (CAD) est alors vérifiée, ce qui permet de conclure par la proposition~\ref{pr-hdpolpoi} lorsque $p=0$ et par la proposition~\ref{pr-htshd_polp} (dont on peut également reprendre directement la démonstration, sans supposer $A$ local) sinon. En particulier, la catégorie $\Tt_{[\mathfrak{M}_{\mathrm{poli},\alpha}]}(A,K)$ est localement noethérienne, et semi-simple si $p=0$ (théorème~\ref{th-finitude-diago}).

En caractéristique $p=0$, il s'agit de démontrer l'annulation de $\mathrm{Ext}^n_{\F(A,K)}(F,G)$ pour $n>0$ et $F$ et $G$ dans $\Tt^\cf_{[\mathfrak{M}_{\mathrm{poli},\alpha}]}(A,K)$. Lorsque $F$ et $G$ sont additifs, cela vient du fait classique de le morphisme canonique $\mathrm{Ext}^*_{\mathbf{Add}(A,K)}(F,G)\to\mathrm{Ext}^*_{\F(A,K)}(F,G)$ est un isomorphisme (cf. par exemple \cite[th.~1.2]{Dja-JLMS}) et de l'hypothèse d'annulation de $\mathrm{Ext}^i_{A_K}(K_\alpha,K_\alpha)$ pour $i>0$, via l'équivalence de catégories \eqref{eq-pteAdda} (page~\pageref{eq-pteAdda}). Le cas général (dans lequel $F$ et $G$ sont nécessairement polynomiaux) se ramène au cas additif par des arguments classiques de dévissage (cf. par exemple \cite[section~8]{DT-add}).

Supposons maintenant $p>0$. Pour $n\in\mathbb{N}$, notons $\C_n$ la sous-catégorie pleine de $\F(A,K)$ image essentielle du foncteur exact $\Pp(K,K)\to\F(A,K)$ de précomposition par l'extension des scalaires le long du morphisme $A\to K$ composé de $\alpha$ et de $K\to K\quad x\mapsto x^{1/p^n}$. Alors $(\C_n)$ est une suite croissance de sous-catégories prélocalisantes de $\Tt_{[\mathfrak{M}_{\mathrm{poli},\alpha}]}(A,K)$, et tout foncteur de $\Tt^\cf_{[\mathfrak{M}_{\mathrm{poli},\alpha}]}(A,K)$ appartient à $\C_n$ pour $n\in\mathbb{N}$ assez grand, en vertu du théorème~\ref{th-struct_ord}. La conclusion découle alors de \cite[th.~7.6]{DT-str_pol}, du lemme~\ref{lm-class_ext_ln} ci-dessous et de ce que l'image de $\alpha$ est infinie (ce qui implique que les torsions de Frobenius itérées du foncteur additif associé à $K_\alpha$ via \eqref{eq-pteAdda} sont des foncteurs absolument simples deux à deux non isomorphes et sans extensions entre eux).
\end{proof}

Le lemme suivant est classique et essentiellement contenu dans \cite[prop.~2.1 et rem.~2.2]{Psch}, par exemple ; nous en rappellerons toutefois pour mémoire la démonstration.

\begin{lemm}\label{lm-class_ext_ln} Soient $\E$ une catégorie de Grothendieck localement noethérienne, $I$ un ensemble ordonné filtrant et $(\C_i)_{i\in I}$ une famille croissante de sous-catégories prélocalisantes\index{termin}{prelocalisante@prélocalisante \emph{(sous-catégorie)}} de $\E$. On suppose que tout objet noethérien de $\E$ appartient à l'une des catégories $\C_i$.

Si $i$ est un élément de $I$ et que $X$ et $Y$ sont des objets de $\C_i$, avec $X$ de type fini, alors le morphisme canonique
$$\underset{j\ge i}{\col}\mathrm{Ext}^*_{\C_j}(X,Y)\to\mathrm{Ext}^*_\E(X,Y)$$
est bijectif.
\end{lemm}

\begin{proof} Soit $0\to Y\to T_n\to\dots\to T_1\xrightarrow{\varphi} X\to 0$ une suite exacte de $\E$. Comme $X$ est de type fini et que $\E$ est localement noethérienne, il existe un sous-objet noethérien $U_1$ de $T_1$ tel que la restriction à $U_1$ de $\varphi$ reste un épimorphisme. On peut de même définir par récurrence sur $m\in\llbracket 2,n\rrbracket$ un sous-objet noethérien $U_m$ de $T_m$ tel que $\mathrm{Im}(U_m\hookrightarrow T_m\to T_{m-1})=\mathrm{Im}(T_m\to T_{m-1})\cap U_{m-1}$.

L'hypothèse montre qu'il existe $j\in E$ tel que $j\ge i$ et $U_m\in\C_j$ pour tout $m\in\llbracket 1,n\rrbracket$.

Soit $Y'$ le produit fibré des monomorphismes $Y\to T_n$ et $U_n\hookrightarrow T_n$ : on dispose dans $\E$ d'un diagramme commutatif aux lignes exactes
$$\xymatrix{0\ar[r] & Y'\ar[r]\ar[d] & U_n\ar[r]\ar[d] & \dots\ar[r] & U_1\ar[r]\ar[d] & X\ar@{=}[d]\ar[r] & 0\\
0\ar[r] & Y\ar[r] & T_n\ar[r] & \dots\ar[r] & T_1\ar[r] & X\ar[r] & 0
}$$
dont les flèches verticales sont les monomorphismes canoniques. On en déduit un diagramme commutatif aux lignes exactes
$$\xymatrix{0\ar[r] & Y\ar[r]\ar@{=}[d] & Y\underset{Y'}{+}U_n\ar[r]\ar[d] & U_n\ar[r]\ar[d] & \dots\ar[r] & U_1\ar[r]\ar[d] & X\ar@{=}[d]\ar[r] & 0\\
0\ar[r] & Y\ar[r] & T_n\ar[r] & T_{n-1}\ar[r] & \dots\ar[r] & T_1\ar[r] & X\ar[r] & 0
}$$
dans lequel tous les termes de la ligne supérieure appartiennent à $\C_j$, ce qui permet facilement de conclure.
\end{proof}

\begin{rema} On ne peut en général pas s'affranchir de toute hypothèse de finitude sur le foncteur $G$ (voir l'exemple~\ref{exextsdi} ci-après) dans l'énoncé du théorème~\ref{th-isoExt_Fpbar}. C'est néanmoins le cas en caractéristique $p=0$ grâce à l'exemple~\ref{ex-ann_poiinf-disc} ci-après. C'est aussi le cas pour $p>0$ si $A$ est un anneau de type fini grâce à l'exemple~\ref{ex-ann_poiinf2} ci-après.

On peut voir également que, sans hypothèse sur $A$ ni $K$, le morphisme canonique du théorème~\ref{th-isoExt_Fpbar} est toujours un isomorphisme si le foncteur \htr\ (pas nécessairement à caractère fini) $G$ est à valeurs de dimensions finies. Cela se voit par un argument de dualité, car l'analogue du théorème~\ref{th-isoExt_Fpbar} pour les groupes de torsion, qui commutent aux colimites filtrantes, vaut sans aucune hypothèse de finitude sur les foncteurs \hts\ qui apparaissent. 
\end{rema}

L'énoncé et la démonstration du résultat suivant présentent de grandes similitudes avec certaines parties de \cite{Dja-JLMS} et \cite{DT-ext}. Comme nous n'en avons de plus besoin nulle part ailleurs dans ce mémoire, nous ne donnerons qu'une esquisse de la démonstration.

\begin{theo}\label{th-comp_Ext-tordI} Soit $I$ un idéal de $A$ tel que $A/I$ soit isomorphe à un sous-corps fini de $K$ et que $I/I^2$ soit fini. Les assertions suivantes sont équivalentes :
\begin{enumerate}
\item\label{it1-compExt} l'inclusion $\Tt_{[\mathfrak{M}_{\mathrm{tord},I}]}(A,K)\to\F(A,K)$ induit un isomorphisme
$$\mathrm{Ext}^*_{\Tt_{[\mathfrak{M}_{\mathrm{tord},I}]}(A,K)}(F,G)\xrightarrow{\simeq}\mathrm{Ext}^*_{\F(A,K)}(F,G)$$
pour tous foncteurs $F$ de $\Tt_{[\mathfrak{M}_{\mathrm{tord},I}]}(A,K)$ et $G$ de $\Tt^\cf_{[\mathfrak{M}_{\mathrm{tord},I}]}(A,K)$ ;
\item\label{it2-compExt} la cohomologie locale $\underset{n\in\mathbb{N}}{\col}\mathrm{Ext}^i_A(A/I^n,A/I)$ est nulle pour tout $i\in\mathbb{N}^*$.
\end{enumerate}
\end{theo}

\begin{proof}[Esquisse de démonstration]
Les hypothèses faites sur $I$ garantissent que $\Tt^\cf_{[\mathfrak{M}_{\mathrm{tord},I}]}(A,K)$ est la réunion de la suite croissante de sous-catégories prébilocalisantes $(\F(A/I^n,K))_{n\in\mathbb{N}}$ de $\F(A,K)$, et que la catégorie de Grothendieck $\Tt_{[\mathfrak{M}_{\mathrm{tord},I}]}(A,K)$ est localement noethérienne, grâce au théorème~\ref{th-PSS}. En particulier, le lemme~\ref{lm-class_ext_ln} montre que le morphisme canonique $\underset{n}{\col}\mathrm{Ext}^*_{\F(A/I^n,K)}(F,G)\to\mathrm{Ext}^*_{\Tt_{[\mathfrak{M}_{\mathrm{tord},I}]}(A,K)}(F,G)$ est bijectif lorsque $F$ est un foncteur de type fini de $\Tt_{[\mathfrak{M}_{\mathrm{tord},I}]}(A,K)$ et $G$ un foncteur de $\Tt^\cf_{[\mathfrak{M}_{\mathrm{tord},I}]}(A,K)$. L'assertion \ref{it1-compExt} (dans laquelle on peut se contenter de traiter des foncteurs $F$ de type fini par un argument formel de colimite) équivaut donc à dire que le morphisme canonique $\underset{n}{\col}\mathrm{Ext}^*_{\F(A/I^n,K)}(F,G)\to\mathrm{Ext}^*_{\F(A,K)}(F,G)$ est bijectif lorsque $F$ est un foncteur de type fini de $\Tt_{[\mathfrak{M}_{\mathrm{tord},I}]}(A,K)$ et $G$ un foncteur de $\Tt^\cf_{[\mathfrak{M}_{\mathrm{tord},I}]}(A,K)$. Des arguments formels de dualité entre Ext et Tor, de comparaison de foncteurs homologiques, et le lemme~\ref{lm-alg_lin} ci-dessous, montrent ensuite que la condition~\ref{it1-compExt} équivaut à demander que pour tous $n, i\in\mathbb{N}$, il existe $N\ge n$ tel que pour tout entier $m\ge N$ et tous foncteurs $F$ de $\F(A/I^n,K)$ et $X$ de $\F(\mathbf{P}(A/I^n)^\op;K)$, le morphisme canonique
$$\mathrm{Tor}^{\mathbf{P}(A)}_*(\pi^*X,\pi^*F)\to\mathrm{Tor}^{\mathbf{P}(A/I^m)}_*(\alpha_m^*X,\alpha_m^*F)$$
soit un isomorphisme en degré $*\le i$, où $\pi : \mathbf{P}(A)\to\mathbf{P}(A/I^n)$ et $\alpha_m : \mathbf{P}(A/I^m)\to\mathbf{P}(A/I^n)$ désignent les foncteurs de réduction modulo $I^n$, ou encore à demander la même condition lorsque $X$ et $F$ sont la linéarisation de foncteurs additifs de $\mathbf{P}(A/I^n)$ (contravariant pour $X$) vers $\mathbf{Ab}$. Une variante des résultats d'excision de \cite[§\,5]{DT-ext} (reposant sur la correspondance de Dold-Kan) permet de voir que cette condition équivaut à la suivante : pour tous $n, i\in\mathbb{N}$, il existe $N\ge n$ tel que pour tout entier $m\ge N$ et tous $A/I^n$-modules $U$ et $V$, le morphisme canonique
$$\mathrm{Tor}^A_*(U,V)\to\mathrm{Tor}^{A/I^m}_*(U,V)$$
est un isomorphisme en degré $*\le i$. Des arguments formels de comparaison de foncteurs homologiques et de filtrations montrent que cette condition revient à demander que pour tout $i\in\mathbb{N}$, il existe $N\in\mathbb{N}$ tel que $\mathrm{Tor}^A_m(A/I^n,A/I)=0$ pour $0<m\le i$. L'équivalence avec la condition~\ref{it2-compExt} s'obtient de nouveau par dualité (d'espaces vectoriels sur $A/I$, et non plus sur $K$ comme précédemment) entre Ext et Tor à l'aide du lemme~\ref{lm-alg_lin}.
\end{proof}

L'énoncé simple suivant est laissé en exercice.

\begin{lemm}\label{lm-alg_lin} Soient $V$ un espace vectoriel, $\cdots\to U_{n+1}\xrightarrow{u_n} U_n\to\dots\to U_0$ une tour d'espaces vectoriels et $f_n : V\to U_n$ des applications linéaires telles que $f_n=u_n f_{n+1}$. Les assertions suivantes sont équivalentes :
\begin{enumerate}
%\item pour tout ensenble $E$, l'application linéaire canonique $\underset{n\in\mathbb{N}}{\col}(U_n^{\oplus E})^\vee\to (V^{\oplus E})^\vee$ induite par les précédentes est un isomorphisme, où $(-)^\vee$ désigne la dualité ;
\item l'application linéaire $\underset{n\in\mathbb{N}}{\col}U_n^\vee\to V^\vee$ induite par les précédentes est un isomorphisme, où $(-)^\vee$ désigne la dualité ;
\item pour $n$ assez grand grand, $f_n$ est injective, et il existe $m\ge n$ tel que le $U_m\xrightarrow{u_n\dots u_{m-1}}U_n$ se factorise par $f_n$.
\end{enumerate}
\end{lemm}

Comme la cohomologie locale qui apparaît dans la condition~\ref{it2-compExt} du théorème~\ref{th-comp_Ext-tordI} est nulle lorsque l'anneau $A$ est noethérien (cf. par exemple \cite[lemme~2.4.1]{Sch-Si} ; cela découle aussi directement de la proposition~\ref{preserv-inj}), on en déduit l'important cas particulier suivant :
\begin{coro} Supposons l'anneau $A$ noethérien. Soit $I$ un idéal de $A$ tel que $A/I$ soit isomorphe à un sous-corps fini de $K$. Alors l'inclusion $\Tt_{[\mathfrak{M}_{\mathrm{tord},I}]}(A,K)\to\F(A,K)$ induit un isomorphisme
$$\mathrm{Ext}^*_{\Tt_{[\mathfrak{M}_{\mathrm{tord},I}]}(A,K)}(F,G)\xrightarrow{\simeq}\mathrm{Ext}^*_{\F(A,K)}(F,G)$$
pour tous foncteurs $F$ de $\Tt_{[\mathfrak{M}_{\mathrm{tord},I}]}(A,K)$ et $G$ de $\Tt^\cf_{[\mathfrak{M}_{\mathrm{tord},I}]}(A,K)$.
\end{coro}

\section{Décompositions en poids fortes infinies}\label{sscexext}

\begin{conv} Dans toute la section~\ref{sscexext}, $\C$ désigne une catégorie essentiellement petite munie d'une action d'un monoïde $M$ et $E$ un sous-ensemble de $\mon(M,K_\mu)$. On munit $\mon(M,K_\mu)\subset K^M$ de la topologie induite par la topologie produit, où chaque facteur $K$ est muni de la topologie discrète.
\end{conv}

On fait également, dans la section~\ref{sscexext}, l'hypothèse suivante :
\begin{hyp}\label{hyp-groprod} Il existe un objet $x$ de $\C$ et des objets $T_\pi$ de $\F(\C;K)_\pi$, pour tout $\pi\in E$, tels que $T_\pi(x)\ne 0$.
\end{hyp}

Cette hypothèse est notamment vérifiée pour $\A=\mathbf{P}(A)$ muni de l'action canonique de $A_\mu$, car pour tout $w\in\mon(A_\mu,K_\mu)$, $\dim_K P^A_w(A)=1$.

\begin{lemm}\label{lm-topdisc} Soit $(\pi_n)_{n\in\mathbb{N}}$ une suite de $\mon(M,K_\mu)$ convergeant vers $w$. Pour toute suite $(T_n)_{n\in\mathbb{N}}$ de foncteurs de $\F(\C;K)$, si $T_n$ est homogène de poids fort $\pi_n$, alors $\prod_{n\in\mathbb{N}}T_n/\bigoplus_{n\in\mathbb{N}}T_n$ est homogène de poids fort $w$.
\end{lemm}

\begin{proof}
Soit $a\in M$. Comme la suite $(\pi_n)$ converge vers $w$, il existe $N\in\mathbb{N}$ tel que $\pi_n(a)=w(a)$ pour $n\ge N$. Puisque chaque $T_n$ est homogène de poids fort $\pi_n$, il s'ensuit que l'action de $a$ sur $\prod_{n\in\mathbb{N}}T_n/\bigoplus_{n\in\mathbb{N}}T_n=\prod_{n\ge N}T_n/\bigoplus_{n\ge N}T_n$ coïncide avec la multiplication par $w(a)$.
\end{proof}

\begin{lemm}\label{lm-prodsom} Soient $w\in\mon(M,K_\mu)\setminus E$, $T_\pi$ (pour $\pi\in E$) et $X$ des foncteurs de $\F(\C;K)$. On suppose que $X$ (resp. chaque $T_\pi$) est homogène de poids fort $w$ (resp. $\pi$). Alors
$$\mathrm{Ext}^1_{\F(\C;K)}\Big(X,\bigoplus_{\pi\in E}T_\pi\Big)\simeq\mathrm{Hom}_{\F(\C;K)}\Big(X,\prod_{\pi\in E}T_\pi\Big/\bigoplus_{\pi\in E}T_\pi\Big).$$
\end{lemm}

\begin{proof} On a
$$\mathrm{Ext}^*_{\F(\C;K)}\Big(X,\prod_{\pi\in E}T_\pi\Big)\simeq\prod_{\pi\in E}\mathrm{Ext}^*_{\F(\C;K)}(X,T_\pi)=0$$
d'après la proposition~\ref{pr-etr}. La suite exacte longue associée au foncteur cohomologique $\mathrm{Ext}^*(X,-)$ donne alors la conclusion.
\end{proof}

\begin{prop}\label{pr-Ext1-SDI} Soit $w : M\to K_\mu$ un poids n'appartenant pas à $E$. Les assertions suivantes sont équivalentes :
\begin{enumerate}
\item\label{itpApp1} pour tout objet $X$ de $\F(\C;K)_w$ et tous objets $T_\pi$ de $\F(\C;K)_\pi$ pour $\pi\in E$, $\mathrm{Ext}^1_{\F(\C;K)}(X,\bigoplus_{\pi\in E}T_\pi)=0$ ;
\item\label{itpApp2} pour tous objets $T_\pi$ de $\F(\C;K)_\pi$, $(\prod_{\pi\in E}T_\pi/\bigoplus_{\pi\in E}T_\pi)^w$ est nul ;
\item\label{itpApp3} $w$ n'appartient pas à l'adhérence  séquentielle de $E$ dans $\mon(M,K_\mu)$.
\end{enumerate}
\end{prop}

\begin{proof}
L'équivalence entre \ref{itpApp1} et \ref{itpApp2} découle du lemme~\ref{lm-prodsom}.

Le lemme~\ref{lm-topdisc} et l'hypothèse~\ref{hyp-groprod} montrent que \ref{itpApp2} implique \ref{itpApp3}.

Pour montrer la réciproque, supposons qu'il existe des objets $T_\pi$  de $\F(\C;K)_\pi$ tels que $(\prod_{\pi\in E}T_\pi/\bigoplus_{\pi\in E}T_\pi)_w$ soit non nul : soient $c$ un objet de $\C$ sur lequel ce foncteur est non nul, et $(u_\pi)_{\pi\in E}$ un élément de $\prod_{\pi\in E}T_\pi(c)$ dont l'image dans $(\prod_{\pi\in E}T_\pi/\bigoplus_{\pi\in E}T_\pi)(c)$ appartienne à $(\prod_{\pi\in E}T_\pi/\bigoplus_{\pi\in E}T_\pi)^w(c)$ et soit non nulle. L'ensemble $E'$ des $\pi$ tels que $u_\pi\ne 0$ est donc infini ; soit $S$ une partie infinie \emph{dénombrable} de $E'$.

 Pour tout $a\in M$, le morphisme induit par l'endomorphisme $a$ de $c$ envoie $(u_\pi)$ sur $(\pi(a)u_\pi)$, puisque $T_\pi$ appartient à $\F(\C;K)_\pi$ pour tout $\pi\in E$. Mais dans $(\prod_{\pi\in E}T_\pi/\bigoplus_{\pi\in E}T_\pi)(c)$, cette image doit coïncider avec $(w(a)u_\pi)$ ; autrement dit, la famille $((w(a)-\pi(a)).u_\pi)_{\pi\in E}$ est presque nulle. En particulier, l'ensemble $\{\pi\in S\,|\,w(a)\ne\pi(a)\}$ est fini pour tout $a\in M$. Cela signifie que $w$ appartient à l'adhérence de $S$ dans $\mon(M,K_\mu)$. Comme $S$ est une partie infinie dénombrable de $E$, cela montre que $w$ appartient à l'adhérence séquentielle de $E$, d'où la proposition.
\end{proof}

\begin{exem}\label{exextsdi} Supposons que le corps $K$ est infini et localement fini ; notons $\mathrm{I}$ le foncteur de $\F(K,K)$ d'inclusion de $\mathbf{P}(K)$ dans $K\Md$ et, pour $r\in\mathbb{N}$, $\mathrm{I}^{(r)}$ sa $r$-ième torsion de Frobenius (i.e. l'extension des scalaires le long du morphisme de Frobenius itéré $r$ fois). Notons $\varphi\in\mon(K_\mu,K_\mu)$ le poids $x\mapsto x^p$. Comme $K$ est infini, la suite $(\varphi^r)_{r\in\mathbb{N}}$ est injective. Du fait que $K$ est localement fini, la suite $(\varphi^{n!})_{n\in\mathbb{N}}$ converge vers $1$.

En utilisant les lemmes~\ref{lm-topdisc} et~\ref{lm-prodsom}, on en déduit que la dimension de l'espace vectoriel $\mathrm{Ext}^1_{\F(K,K)}(\mathrm{I},\bigoplus_{n\ge 1}\mathrm{I}^{(n!)})$ égale le continu. Tout élément non trivial de ce groupe d'extensions produit un foncteur additif de $\F(K,K)$ extension d'un foncteur simple \hd\ par un foncteur \hd, mais qui ne possède pas de décomposition en poids faible. En particulier, on voit que la sous-catégorie prélocalisante $\Tt(K,K)=\Oo(K,K)$ (proposition~\ref{prto}) de $\F(K,K)$ n'est pas stable par extensions.
\end{exem}

On généraliser la proposition~\ref{pr-Ext1-SDI} au-delà des groupes $\mathrm{Ext}^1$, au moins lorsque $M$ est au plus dénombrable --- sinon, adhérence et adhérence séquentielle diffèrent a priori dans $\mon(M,K_\mu)$, et nous ne savons pas si la condition suffisante à l'annulation cohomologique dans l'énoncé ci-dessous est nécessaire.

\begin{prop}\label{pr-ann_Ext-discr} Soit $w$ un élément de $\mon(M,K_\mu)$ n'appartenant pas à l'adhérence de $E$. Si $X$ est un foncteur de $\F(\C;K)_{<w>}$ et $T$ un foncteur de $\F(\C;K)$ possédant une décomposition en poids faible tel que $\Pi(T)\subset E$, alors $\mathrm{Ext}^*_{\F(\C;K)}(X,T)$ est nul.
\end{prop}

\begin{proof} L'hypothèse $w\notin\bar{E}$ signifie qu'il existe une partie finie de $M$ sur laquelle aucun élément de $E$ ne coïncide avec $w$. Il suffit évidemment de montrer le résultat lorsque cette partie est réduite à un élément $a$ de $M$. La démonstration est alors essentiellement la même que celle de la proposition~\ref{pr-etr} : on commence par supposer que $X$ est homogène de poids \textit{fort} $w$. Dans ce cas, l'élément $[a]-w(a)$ de $K[M]$ opère par $0$ sur $F$, donc sur $\mathrm{Ext}^*_{\F(\C;K)}(X,T)$ (la catégorie $\F(\C;K)$ est $K[M]$-linéaire), mais par un automorphisme sur $T$, donc aussi sur $\mathrm{Ext}^*_{\F(\C;K)}(X,T)$. Ainsi, $\mathrm{Ext}^*_{\F(\C;K)}(X,T)=0$ lorsque $X$ appartient à $\F(\C;K)_w$. Le cas où $X$ appartient à $\F(\C;K)_{<w>}$ s'en déduit par un argument formel facile, comme pour la proposition~\ref{pr-etr}.
\end{proof}

\begin{exem}\label{ex-ann_poiinf-disc}
En caractéristique $p=0$, pour $M=\mathbb{Z}_\mu$, l'évaluation en $2$ est un morphisme continu et injectif $\mon_\mathrm{pol}(\mathbb{Z}_\mu,K_\mu)\simeq\mathbb{N}\to K_\mu$, de sorte que $\mon(M,K)$ est \textit{discret}. Il s'ensuit que la décomposition en poids des foncteurs analytiques de $\F(\A;K)$ (où $\A$ est une catégorie additive essentiellement petite) --- cf. proposition~\ref{pr-decpolcar0} --- induit une décomposition entre groupes d'extensions.
\end{exem}

\begin{exem}\label{ex-ann_poiinf2} Supposons que $p$ est premier, que $A$ est un anneau de type fini, que le corps $K$ est parfait et que $\alpha : A\to K$ est un morphisme d'anneaux d'image infinie. Soit $E$ le sous-monoïde de $\mon(A_\mu,K_\mu)$ constitué des $\alpha^{n/p^r}$ pour $n, r\in\mathbb{N}$. Alors $E$ est discret dans $K^M$, car il existe $a\in A$ tel que $\alpha(a)$ soit transcendant sur $\FF_p$ (sinon, l'image de $\alpha$ serait un sous-corps de $K$ algébrique sur $\FF_p$, donc fini, puisque $A$ est un anneau de type fini), de sorte que l'évaluation en $a$ définit un morphisme continu injectif $E\to K_\mu$, ce qui entraîne que $E$ est discret. Ainsi, la proposition~\ref{pr-ann_Ext-discr} fournit une décomposition des groupes d'extensions entre foncteurs admettant des décompositions en poids avec des poids appartenant à $E$.
\end{exem}

%Reprendre aussi les arguments de cohérence de $\F(k,k)$ apparaissant dans l'appendice des vieilles versions de \cite{DT-schw}, y compris dans un cadre plus général ???
%
%AUSSI : exemple d'idéal $I$ de $A$ avec $I=I^2$ mais $\mathrm{Ext}^{>0}_A(A/I,A/I)\ne 0$ pour d'autres contre-exemples ??? -> $A=\FF_p[\mathbb{Z}/p^\infty]$, $I$ idéal d'augmentation donne un tel cas, le $\mathrm{Ext}^2$ est non nul.
%
%\section{Décompositions en poids faibles infinies}
%
%Mettre a priori des choses plus courtes, peut-être en se contentant de faire référence à l'exemple~\ref{ex-Vpasepais}.

%AUTRES APPENDICES ou compléments ENVISAGEABLES : sur le centre d'une catégorie (y compris en recyclant une partie du truc sur les blocs)...

\chapter{Le cas où $A$ est fini modulo nilpotents}\label{ap_afmn}

Lorsque l'anneau quotient $A/\nil(A)$ est fini, les résultats connus sur $\F(A,K)$ se ramènent essentiellement aux deux situations suivantes (le cas général s'y ramenant en décomposant $A$ en produit direct de deux anneaux entrant dans chacune d'entre elle). Lorsque la caractéristique de $A$ est inversible dans $K$ (cas d'inégale caractéristique), tous les foncteurs de type fini de $\F^\df(A,K)$ sont antipolynomiaux (cela découle de \cite[th.~4.10]{DTV}, cité au théorème~\ref{th-DTVglob} du présent mémoire). Lorsque la caractéristique de $A$ est une puissance de $p$, tous les foncteurs de type fini de $\F^\df(A,K)$ sont \phs, grâce à la proposition~\ref{pr-htr_impl_ph}, mais cela dit assez peu sur la structure de ces foncteurs. Nous verrons à la section~\ref{scxFP} qu'ils ne possèdent pas nécessairement de décomposition de type $FP$.

Nous complétons l'étude du cas d'inégale caractéristique en démontrant le théorème~\ref{th-ttantipol}, qui caractérise l'antipolynomialité de tous les poids de $\mon(A_\mu,K_\mu)$ par l'existence d'un idéal cofini de carré nul, lorsque $K$ possède suffisamment de racines de l'unité. Il nous a semblé instructif de montrer, pour compléter la section~\ref{scarren}, qu'on peut déduire facilement de résultats classiques de Hartley \cite{Hartley77,Hartley83} que cette condition d'existence d'un idéal cofini de carré nul (toujours en inégale caractéristique) équivaut au fait que toutes les représentations $K$-linéaires irréductibles de $\GL_n(A)$ ou de $\SL_n(A)$ sont de dimension finie (et donc antipolynomiales en un sens directement adapté de celui des foncteurs), où $n\ge 3$ est un entier fixé. Ainsi une condition de finitude ou d'antipolynomialité pour les représentations $K$-linéaires génériques des groupes linéaires sur $A$ équivaut-elle à une condition partiellement analogue pour les représentations usuelles des groupes linéaires, même si les méthodes pour le démontrer sont totalement indépendantes.

Dans la section~\ref{spfsimp}, nous étudions le pendant polynomial de ces considérations (pour lequel on peut aborder de façon simultanée la situation fonctorielle et les représentations usuelles des groupes linéaires). Nous montrons précisément que $A/\nil(A)$ est fini de caractéristique $p$
 si et seulement si tout foncteur simple de $\F(A,K)$ (sans restriction a priori sur les dimensions de ses valeurs), ou encore si et seulement si toutes les représentations $K$-linéaires irréductibles de $\GL_n(A)$ (sans hypothèse de dimension finie a priori), pour un entier $n\ge 3$ fixé, vérifient une propriété polynomiale.

\section{La condition d'existence d'un idéal cofini de carré nul}\label{scarren}

\subsection{Antipolynomialité des poids : démonstration du théorème~\ref{th-ttantipol}}\label{par-thIc2}

Nous commençons par établir la partie facile du théorème.

\begin{proof}[Démonstration de l'assertion~\ref{it-thPfac} du théorème~\ref{th-ttantipol}] Pour $\p\in\spec(A)$, la considération du poids $\xi_{\p} : A\to K$ (où la notation $\xi$ est celle de la proposition-définition~\ref{prdfla}) montre, en appliquant la proposition~\ref{pr-zerantipol}, que $\p$ est maximal et $K$-cotrivial.

Considérant maintenant le poids $\xi_{A\setminus A^\times}$, la proposition~\ref{pr-zerantipol} montre que la réunion des idéaux maximaux de $A$ est égale à une réunion finie d'idéaux maximaux, ce qui implique (cf. par exemple \cite[ch.~II, §\,1, prop.~2]{Bki-AC1}) que $A$ n'a qu'un nombre fini d'idéaux maximaux. Il s'ensuit, par la proposition~\ref{pr-qpsi} (implication \ref{itqpf3}$\Rightarrow$\ref{itqpf1}), que $A$ est quasi-parfait et que $A/\rj(A)$ est $K$-trivial.
\end{proof}

Pour la suite de la démonstration, nous aurons besoin du lemme qui suit sur le foncteur de dualité des groupes abéliens $(-)^\vee:=\mathrm{Hom}(-,\mathbb{Q}/\mathbb{Z}) : \mathbf{Ab}^\op\to\mathbf{Ab}$, qui est exact, fidèle et continu.

%il serait mieux de trouver une réf !
\begin{lemm}\label{lm-dual-grab} Soient $l$ un nombre premier et $V$ un $l$-groupe abélien. On suppose que $V^\vee$ est un groupe abélien de torsion.
\begin{enumerate}
\item Il existe $N\in\mathbb{N}$ tel que $V$ soit somme directe de groupes cycliques $\mathbb{Z}/l^i$ avec $i\le N$.
\item Si $(U_n)_{n\in\mathbb{N}^*}$ une suite décroissante de sous-groupes de $V$ telle que le monomorphisme canonique
$\underset{n\in\mathbb{N}^*}{\col}(V/U_n)^\vee\to V^\vee$
est un isomorphisme, alors il existe $n\in\mathbb{N}^*$ tel que $U_n=0$.
\end{enumerate}
\end{lemm}

\begin{proof} Comme $V$ est un $l$-groupe abélien, il existe une suite exacte $0\to V'\to V\to V''\to 0$ où $V'$ est une somme directe de $l$-groupes cycliques et $V''$ est un $l$-groupe abélien divisible \cite[{\em Theorem}~32.3]{Fuchs}. On en déduit aussitôt le premier point, puisque le dual d'un groupe abélien divisible non nul n'est pas de torsion et qu'un produit de $l$-groupes cycliques n'est de torsion que si l'ensemble de leurs cardinaux est borné.

On établit le deuxième point par récurrence sur l'entier $N$ donné par le premier. Si $N=1$, $V$ est un $\FF_l$-espace vectoriel. Si $(U_n)$ est stationnaire, elle stationne forcément en $0$, puisque $(-)^\vee$ est exact et fidèle. Sinon, quitte à extraire une sous-suite, on peut supposer $(U_n)$ {\em strictement} décroissante. Choissons $e_n\in U_n\setminus U_{n+1}$ pour tout $n\in\mathbb{N}^*$ : alors $(e_n)$ est une famille libre du $\FF_l$-espace vectoriel $V$, donc il existe $\varphi\in V^\vee$ tel que $\varphi(e_n)\ne 0$ pour tout $n$ ; ainsi $\varphi$ n'est nul sur aucun $U_n$, contradiction qui termine la démonstration dans le cas $N=1$.

Supposons maintenant $N>1$ et l'assertion établie pour les groupes annulés par $l^{N-1}$. Il existe un sous-groupe $W$ de $V$ tel que $W$ soit annulé par $l$ et $V/W$ par $l^{N-1}$. L'hypothèse de récurrence appliquée à $V/W$ montre que $U_n\subset W$ pour $n$ assez grand. Comme l'inclusion induit un morphisme surjectif $V^\vee\twoheadrightarrow W^\vee$, le cas $N=1$ déjà traité montre maintenant que $U_n=0$ pour $n$ assez grand, d'où le lemme.
\end{proof}

\begin{proof}[Démonstration de l'assertion~\ref{it-thPsp} du théorème~\ref{th-ttantipol}] 
Comme l'anneau $A$ est quasi-parfait d'après l'assertion~\ref{it-thPfac} qu'on a déjà établie, c'est un produit fini d'anneaux locaux, de sorte qu'on peut supposer $A$ local. Le groupe des unités de $A$ est la somme directe du groupe (cyclique fini) des unités de son corps résiduel, de caractéristique $l$, et du groupe multiplicatif $1+\rj(A)$, qui est un $l$-groupe abélien que nous noterons simplement $V$.

On a $\mon(A_\mu,K_\mu)\simeq\mathbf{Ab}(A^\times,K^\times)_+$ (cf. exemple~\ref{exev-mormon}), et $\mathbf{Ab}(A^\times,K^\times)\simeq\mathbf{Ab}(A^\times,K^\times_{\mathrm{tor}})$ contient un facteur direct isomorphe à $V^\vee$, puisque $K$ est de caractéristique différente de $l$ et contient assez de racines de l'unité. 

Par conséquent, tout morphisme $V\to\mathbb{Q}/\mathbb{Z}$ se factorise par la restriction à $V$ d'un morphisme d'anneaux de $A$ sur un anneau fini. Comme le radical d'un anneau fini est nilpotent, il s'ensuit que, en notant $U_n$ le groupe multiplicatif $1+\rj(A)^n$ (pour $n\in\mathbb{N}^*$), tout morphisme $V\to\mathbb{Q}/\mathbb{Z}$ est nul sur $U_n$ pour $n$ assez grand ; autrement dit, le morphisme canonique $\underset{n\in\mathbb{N}^*}{\col}(V/U_n)^\vee\to V^\vee$ est un isomorphisme.

Comme chaque groupe $V/U_n$ est un $l$-groupe abélien de torsion bornée (car $U_i/U_{i+1}$ est isomorphe au groupe additif $\rj(A)^i/\rj(A)^{i+1}$, qui est annulé par $l$), les $(V/U_n)^\vee$ sont de torsion, et $V^\vee$ également. Le lemme~\ref{lm-dual-grab} montre alors que $U_n$ est trivial pour $n$ assez grand, c'est-à-dire que $\rj(A)$ est nilpotent, ainsi $A$ est semi-primaire.
\end{proof}

Nous passons maintenant à la démonstration de la dernière assertion du théorème, qui est également la plus longue.

\begin{proof}[Démonstration de l'assertion~\ref{it-thPdur} du théorème~\ref{th-ttantipol}] Au vu de ce qu'on a déjà montré, on peut supposer que $A$ est local ; il s'agit de montrer que $A$ possède un idéal cofini de carré nul (comme $A/\rj(A)$ est $K$-trivial, tout idéal cofini est nécessairement $K$-cotrivial). Définissons
$$\mathfrak{t}:=\{x\in\rj(A)\,|\,\cd(Ax)<\infty\}\;;$$
c'est un idéal de $A$ ; comme $A/\rj(A)$ est fini, on a aussi $\mathfrak{t}=\{x\in A\,|\,\cd(\rj(A)x)<\infty\}$. On va procéder en trois étapes, mais on fait quelques rappels et observations préliminaires auparavant.

Tout d'abord, comme $K$ contient assez de racines de l'unité et que le sous-groupe $(1+\rj(A))_\mu$ de $A^\times$ est de torsion, et sans $p$-torsion, $K^\times$ est un groupe abélien injectif et $\mathbf{Ab}((1+\rj(A))_\mu,K^\times)\simeq (1+\rj(A))_\mu^\vee$.

Si $T$ est une partie de $\mathbb{N}$, disons qu'une suite $(e_i)_{i\in T}$ d'éléments d'un groupe abélien $V$ est \emph{pseudo-libre} si, pour tout $n\in T$, $e_n$ n'appartient pas au sous-groupe de $V$ engendré par les $e_i$ pour $i\in T$ tels que $i<n$. Si $(e_i)$ est une suite pseudo-libre de $V$, alors il existe $l\in V^\vee$ tel que $l(e_i)\ne 0$ pour tout $i$.
\begin{enumerate}
\item[(A)] Si $X$ est une partie finie de $A\setminus\{0\}$, alors il existe un idéal cofini $I$ de $A$ tel que $I\cap X=\varnothing$.

(Comme une intersection finie d'idéaux cofinis est cofinie, on peut se borner à traiter le cas où $X$ est réduit à un élément $x$. Si $x\notin\rj(A)$, l'assertion est triviale : $I=\rj(A)$ convient. Sinon, $1+x\in A^\times\setminus\{1\}$, donc il existe $l\in\mathbf{Ab}(A^\times,K^\times)\simeq (A^\times)^\vee$ tel que $l(1+x)\ne 1$ ; par hypothèse, il existe un idéal cofini $I$ de $A$ tel que $l$ se factorise par la projection canonique $A^\times\twoheadrightarrow (A/I)^\times$, d'où $x\notin I$.)
\item[(B)] Si $(e_i)_{i\in\mathbb{N}}$ est une suite d'éléments de $A\setminus\mathfrak{t}$, il existe une suite $(u_i)_{i\in\mathbb{N}}$ d'éléments de $\rj(A)$ telle que $(1+e_i u_i)_{i\in\mathbb{N}}$ soit une suite pseudo-libre du groupe muliplicatif $1+\rj(A)$.

(Par récurrence, on peut supposer les $u_j$ construits, pour $j<i$, de sorte que $(1+e_j u_j)_{j<i}$ soit pesudo-libre dans $(1+\rj(A))_\mu$. Le sous-groupe de $(1+\rj(A))_\mu$ engendré par les $1+e_j u_j$ pour $j<i$ est fini, et $\rj(A)e_i$ est infini puisque $e_i\notin\mathfrak{t}$. On peut donc trouver $u_i\in\rj(A)$ tel que $1+e_i u_i$ n'appartienne pas à ce sous-groupe, de sorte que $(1+e_j u_j)_{j\le i}$ est pseudo-libre dans $(1+\rj(A))_\mu$.)
\end{enumerate}

\noindent
\textbullet\
\emph{Montrons que l'idéal $\mathfrak{t}$ de $A$ est cofini.}

Supposons le contraire : il existe alors un sous-groupe infini dénombrable $V$ de $A_\mathrm{add}$ tel que $V\cap\mathfrak{t}=\{0\}$. L'observation (B) fournit une fonction (ensembliste) $f : V\setminus\{0\}\to\rj(A)$ telle que la suite $(1+v_n f(v_n))_{n\in\mathbb{N}}$, où $\mathbb{N}\to V\setminus\{0\}\quad n\mapsto v_n$ est une bijection, soit pseudo-libre dans le groupe multiplicatif $1+\rj(A)$. Il existe donc $l\in\mathbf{Ab}(A^\times,K^\times)$ tel que $l(1+vf(v))\ne 1$ pour tout $v\in V\setminus\{0\}$, puis un idéal cofini $I$ de $A$ tel que $vf(v)\notin I$, donc $v\notin I$,  pour tout $v\in V\setminus\{0\}$. Autrement dit, $V\cap I=\{0\}$. Ainsi, $V$ s'injecte dans $A/I$, qui est fini, ce qui est absurde.

\noindent
\textbullet\
\emph{Montrons que l'idéal $\mathfrak{t}^2$ de $A$ est fini.}

Supposons le contraire. On note au préalable que si $\mathfrak{u}\subset\mathfrak{t}$ est un idéal cofini de $A$, alors $\mathfrak{u}^2$ de $A$ est également infini. En effet, comme $\mathfrak{u}$ est cofini, il existe un sous-groupe fini $V$ de $A_\mathrm{add}$ tel que $\mathfrak{t}=\mathfrak{u}+V$ ; on a alors $\mathfrak{t}^2=\mathfrak{u}^2+\mathfrak{t}.V\subset\mathfrak{u}^2+A.V$. Comme $Av$ est fini pour tout $v\in \mathfrak{t}$, la finitude de $V\subset\mathfrak{t}$ entraîne celle de $A.V$, de sorte que $\mathfrak{t}^2$ ne peut être infini que si $\mathfrak{u}^2$ l'est.

Maintenant, on construit par récurrence sur $n\in\mathbb{N}$ des éléments $r_n$ et $s_n$ de $\mathfrak{t}$ tels que $(1+r_n s_n)$ soit une suite pseudo-libre du groupe multiplicatif $1+\rj(A)$ et que $r_n s_m=0$ pour $n\ne m$. Soit $i\in\mathbb{N}$ ; supposons donnés $r_n$ et $s_n$, pour $n<i$, avec ces propriétés. Alors $\mathfrak{u}:=\mathfrak{t}\cap\bigcap_{n<i}(\mathrm{ann}(r_n)\cap\mathrm{ann}(s_n))$ (où ann désigne l'annulateur) est un idéal cofini de $A$, car $\mathrm{ann}(x)$ cofini pour tout $x\in\mathfrak{t}$ par définition ; ainsi, $\mathfrak{u}^2$ est infini d'après ce qui précède. On peut donc trouver $r_i, s_i\in\mathfrak{u}$ tels que $1+r_i s_i$ n'appartienne pas au sous-groupe (qui est fini) du groupe multiplicatif $1+\rj(A)$ engendré par les $1+r_n s_n$ pour $n<i$. On a par ailleurs $r_i s_n=0$ et $r_n s_i=0$ pour tout $n<i$, ce qui achève la construction des suites $(r_n)$ et $(s_n)$.

Comme $(1+r_n s_n)$ est une suite pseudo-libre de $(1+\rj(A))_\mu$, il existe $l\in\mathbf{Ab}(A^\times,K^\times)$ tel que $l(1+r_n s_n)\ne 1$ pour tout $n\in\mathbb{N}$, puis un idéal cofini $I$ de $A$ ne contenant aucun $r_n s_n$. Comme $A/I$ est fini, il existe $n\in\mathbb{N}$ tel que $\bar{r}_n\in\sum_{i<n}<\bar{r}_i>$ dans le groupe additif $A/I$, où $\bar{r}_i$ désigne la classe de $r_i$ modulo $I$ et $<\bar{r}_i>$ le sous-groupe qu'elle engendre. Autrement dit, $r_n\in I+\sum_{i<n}<r_i>$ dans le groupe $A_\mathrm{add}$. En multipliant par $s_n$ et en utilisant la nullité de $r_i s_n$ pour $i\ne n$, il vient $r_n s_n\in I$, d'où la contradiction recherchée.

\noindent
\textbullet\
\emph{Fin de la démonstration.}

Comme $\mathfrak{t}^2$ est fini, il existe un idéal cofini $I$ de $A$ tel que $I\cap\mathfrak{t}^2=\{0\}$ (observation (A) supra). Du fait que $\mathfrak{t}$ est un idéal cofini, l'idéal $I\cap\mathfrak{t}$ est encore cofini. Comme $(I\cap\mathfrak{t})^2\subset I\cap\mathfrak{t}^2=\{0\}$, cela termine la démonstration.
\end{proof}

\subsection{Endofinitude des représentations linéaires irréductibles des groupes linéaires}

Nous commençons par quelques rappels de théorie des groupes. Si $G$ est un groupe, alors pour tout nombre premier $l$, $G$ possède un plus grand $l$-sous-groupe (i.e. sous-groupe dont chaque élément a pour ordre une puissance de $l$) distingué, noté $O_l(G)$. Par convention, $O_0(G)$ désigne le sous-groupe trivial de $G$.

On rappelle qu'un groupe est dit \emph{virtuellement abélien} s'il possède un sous-groupe abélien (qu'on peut supposer distingué) d'indice fini.

Enfin, si $R$ est un anneau (non nécessairement commutatif) et $M$ un $R$-module à gauche, on rappelle que $M$ est dit \emph{endofini}\index{termin}{endofini|textbf} si $M$ est un module à gauche fini sur l'anneau (non nécessairement commutatif) $\mathrm{End}_R(M)$.

\begin{theo}[Hartley]\label{th-Hartley} Soit $G$ un groupe \emph{localement fini} (i.e. dont tout sous-groupe de type fini est fini). Les assertions suivantes sont équivalentes :
\begin{enumerate}
    \item toute représentation $K$-linéaire irréductible de $G$ est endofinie ;
    \item le groupe $G/O_p(G)$ est virtuellement abélien.
\end{enumerate}
\end{theo}

Ce résultat est donné par \cite[th.~B]{Hartley77} en caractéristique $p=0$, et est le théorème principal de \cite{Hartley83} pour $p>0$.

Pour d'appliquer le théorème~\ref{th-Hartley} aux groupes linéaires, nous aurons besoin du lemme élémentaire suivant.

\begin{lemm}\label{lm-Op}
\begin{enumerate}
    \item\label{it-repev1} Soient $G_1$ et $G_2$ des groupes.
    
    On a $O_p(G_1\times G_2)=O_p(G_1)\times O_p(G_2)$.
    \item\label{it-repev2} Si $H$ est un sous-groupe d'un groupe $G$, alors $O_p(G)\cap H\subset O_p(H)$. Si de plus $H$ est distingué dans $G$, alors la projection $G\twoheadrightarrow G/H$ induit un morphisme de groupes surjectif $G/O_p(G)\twoheadrightarrow (G/H)/O_p(G/H)$.
    \item\label{it-repev3} Si $k$ est un corps et $n\ge 2$ un entier, alors $O_p(\PSL_n(k))$ est trivial sauf dans les deux cas suivants : $n=2$, $p=3$, $k\simeq\FF_2$, et $n=2$, $p=2$ et $k\simeq\FF_3$.
    \item\label{it-annlocf} Supposons l'anneau $A$ localement fini. Soit $n\ge 3$ un entier.
    \begin{enumerate}
        \item\label{itegcar} Si $p$ est nilpotent dans $A$, alors $\PSL_n(A)/O_p(\PSL_n(A))=\PSL_n(A/\nil(A))$.
        \item\label{itcarcr} Si $p$ est inversible dans $A$, alors $O_p(\PSL_n(A))$ est trivial.
    \end{enumerate}
\end{enumerate}
\end{lemm}

\begin{proof}
Les assertions \ref{it-repev1} et \ref{it-repev2} sont immédiates.

Si $k$ est un corps et $n\ge 2$ un entier, alors $\PSL_n(k)$ est un groupe simple non abélien, d'où $O_p(\PSL_n(k))=\{1\}$, sauf si $n=2$ et que $k$ est isomorphe à $\FF_2$ ou $\FF_3$. Comme $\PSL_2(\FF_2)$ (resp. $\PSL_2(\FF_3)$) est isomorphe au groupe symétrique (resp. alterné) sur $3$ (resp. $4$) lettres, on voit par inspection que $O_3(\PSL_2(\FF_2))\simeq\mathbb{Z}/3$ et $O_p(\PSL_2(\FF_2))=\{1\}$ pour $p\ne 3$, et que $O_2(\PSL_2(\FF_3))\simeq (\mathbb{Z}/2)^2$ et $O_p(\PSL_2(\FF_3))=\{1\}$ pour $p\ne 2$. Cela montre le point~\ref{it-repev3}.

Soit $A$ un anneau localement fini. Pour tout idéal maximal $\mathfrak{m}$ de $A$, l'image de $O_p(\PSL_n(A))$ par le morphisme canonique $\PSL_n(A)\to\PSL_n(A/\mathfrak{m})$, qui est surjectif ($A/\mathfrak{m}$ est un corps, donc $\SL_n(A/\mathfrak{m})$ est engendré par les matrices élémentaires), est incluse dans $O_p(\PSL_n(A/\mathfrak{m}))$ (par le point \ref{it-repev2}), qui est trivial pour $n\ge 3$, grâce au point~\ref{it-repev3}. Comme l'intersection des idéaux maximaux de $A$ égale $\nil(A)$ (dans un anneau localement fini, tout idéal premier est maximal), on en déduit que $O_p(\PSL_n(A))$ est inclus dans le noyau $N$ du morphisme canonique $\PSL_n(A)\twoheadrightarrow\PSL_n(A/\nil(A))$.

Si $p$ est nilpotent dans $A$, alors le sous-groupe distingué $N$ de $\PSL_n(A)$ est un $p$-groupe, il est donc inclus dans $O_p(\PSL_n(A))$. Cela achève de de montrer \ref{itegcar}.

Si $p$ est inversible dans $A$, alors $N$ ne contient aucun élément d'ordre $p$, donc l'inclusion $O_p(\PSL_n(A))\subset N$ entraîne la trivialité de $O_p(\PSL_n(A))$, d'où le lemme.
\end{proof}

Le lemme suivant est classique (et élémentaire).
% Idéalement, trouver une réf !
\begin{lemm}\label{lm-grp_lf} Soient $G$ un groupe localement fini tel que $G/O_p(G)$ soit virtuellement abélien et $V$ une représentation $K$-linéaire irréductible de $G$. Alors il existe un sous-groupe distingué d'indice fini de $G$ contenant $O_p(G)$ opérant trivialement sur $V$.
\end{lemm}

\begin{prop}\label{pr-repGLendofini} Soit $G$ un groupe de la forme $\GL_n(A)$, $\SL_n(A)$, $\operatorname{PGL}_n(A)$ ou $\PSL_n(A)$ avec $n\ge 3$, ou, lorsque $p=0$, $\GL_2(A)$ ou $\SL_2(A)$.

Les assertions suivantes sont équivalentes :
\begin{enumerate}
    \item\label{it-fqf} Toute représentation $K$-linéaire irréductible de $G$ se factorise à travers la réduction modulo un idéal cofini de $A$.
    \item\label{it-dfK} toute représentation $K$-linéaire irréductible de $G$ est de dimension finie sur $K$ ;
    \item\label{it-dfC} toute représentation $K$-linéaire irréductible de $G$ est de dimension finie sur le centre de son corps d'endomorphismes ;
    \item\label{it-dfE} toute représentation $K$-linéaire irréductible de $G$ est endofinie\index{termin}{endofini} ;
    \item\label{it-qpftP} $A$ est isomorphe au produit direct d'un anneau $A_1$ possédant un idéal $K$-cotrivial de carré nul et d'un anneau $A_2$ tel que $A_2/\nil(A_2)$ soit fini de caractéristique $p$.
\end{enumerate}
\end{prop}

\begin{proof} Les implications \ref{it-fqf}$\Rightarrow$\ref{it-dfK}$\Rightarrow$\ref{it-dfC}$\Rightarrow$\ref{it-dfE} sont immédiates.

Nommons, comme Hartley \cite{Hartley77,Hartley83}, \emph{restreint} un groupe dont toutes les représentations $K$-linéaires irréductibles sont endofinies. Alors tout sous-groupe et tout quotient d'un groupe restreint est restreint \cite[lemme~2.2]{Hartley77}, ce que nous utiliserons abondamment dans ce qui suit. En particulier, pour montrer l'implication \ref{it-dfE}$\Rightarrow$\ref{it-qpftP}, on peut se contenter des cas de $G=\PSL_3(A)$ et, lorsque $p=0$, $G=\SL_2(A)$.

Supposons que $G$ est un groupe restreint. Commençons par montrer que l'anneau $A$ est localement fini. Comme le sous-quotient $\PSL_2(A)$ de $G$ est restreint, $A$ ne peut pas être de caractéristique nulle, car sinon $\PSL_2(A)$ contiendrait $\PSL_2(\mathbb{Z})$, qui lui-même contient un groupe libre non abélien, or un tel groupe n'est pas restreint\,\footnote{On peut le déduire directement du théorème~\ref{th-Hartley} de la façon suivante : comme un groupe libre non abélien $G$ contient un groupe libre de rang infini dénombrable, si $G$ était restreint, tout groupe dénombrable serait restreint. Or ce n'est pas le cas de $\PSL_2(k)$ si $k$ est la clôture algébrique d'un corps fini, puisque ce groupe dénombrable est localement fini mais simple et non abélien.}. Par décomposition primaire, on peut donc, pour montrer \ref{it-dfE}$\Rightarrow$\ref{it-qpftP}, supposer que la caractéristique de $A$ est une puissance d'un nombre premier $l$.

Montrons maintenant que $A$ est localement fini. On peut pour cela supposer que $A$ est de caractéristique $l$, quitte à remplacer $A$ par $A/l$, car la finitude locale de $A/l$ entraîne celle de $A$, puisque $l$ est nilpotent dans $A$. Il s'agit alors de montrer que $A$ est une algèbre entière sur le corps $\FF_l$. Si ce n'était pas le cas, $A$ contiendrait un sous-anneau isomorphe à la $\FF_l$-algèbre de polynômes $\mathbb{F}_l[X]$, or $\PSL_2(\mathbb{F}_l[X])$ contient un groupe libre non abélien, ce qui permet de conclure comme précédemment.

Ainsi l'anneau $A$ est localement fini, donc $G$ est un groupe localement fini, et le théorème de Hartley (\ref{th-Hartley}) montre que $G/O_p(G)$ est virtuellement abélien. Ainsi, par le lemme~\ref{lm-Op}.\,\ref{it-annlocf}, $\PSL_n(A)$ (resp. $\PSL_n(A/\nil(A))$) est virtuellement abélien si $l\ne p$ (resp. $l=p$). Quitte à remplacer $A$ par $A/\nil(A)$ si $l=p$, il suffit donc de montrer que $A$ possède un idéal cofini de carré nul si $\SL_2(A)$ ou $\PSL_3(A)$ est virtuellement abélien.

Traitons d'abord le cas où $\SL_2(A)$ possède un sous-groupe abélien distingué d'indice fini  $\Gamma$ ; définissons
$V:=\{x\in A\,|\,\mathrm{e}_{1,2}(x)\in\Gamma\}$, où, de manière générale, $\mathrm{e}_{i,j}(x)$ désigne, pour $1\le i\ne j\le n$ et $x\in A$, la matrice élémentaire de $\SL_n(A)$ avec $x$ en position $(i,j)$.
Alors $V$ est clairement un sous-groupe de $A_\mathrm{add}$ ; comme $\Gamma$ est d'indice fini dans $\SL_2(A)$, $V$ est d'indice fini dans $A_\mathrm{add}$. De plus, comme le conjugué de $\mathrm{e}_{1,2}(x)$ par $\left(\begin{array}{cc}0 & -1\\
                                                1 & 0
                                               \end{array}\right)\in\SL_2(A)$ est $\mathrm{e}_{2,1}(-x)$,
le fait que $G$ soit distingué et abélien montre que pour tous $x, y\in V$, les matrices $\mathrm{e}_{1,2}(x)$ et $\mathrm{e}_{2,1}(y)$ commutent, ce qui entraîne $xy=0$.
Soient $a_1,\dots,a_r$ des relevés dans $A$ des éléments du groupe (additif) fini $A/V$ : le sous-groupe $I:=\bigcap_{i=1}^r a_i^{-1}V$ de $A_\mathrm{add}$ est cofini, comme $V$. De plus, pour tout $t\in A$, il existe $i\in\llbracket 1,r\rrbracket$ tel que $a_i-t\in V$, d'où, pour $x\in I\subset V$, $tx=a_i x\in V$ d'après ce qui précède, puis $tx\in I$. Il s'ensuit que $I$ est un idéal de $A$, nécessairement de carré nul.

Traitons maintenant le cas où $\PSL_3(A)$ possède un sous-groupe abélien distingué d'indice fini  $\Gamma$ ; définissons $I:=\{x\in A\,|\,\pi(\mathrm{e}_{1,2}(x))\in\Gamma\}$, où $\pi : \SL_3(A)\twoheadrightarrow\PSL_3(A)$ est la projection. Comme $\Gamma$ est distingué, on a $\pi(\mathrm{e}_{i,j}(x))\in\Gamma$ pour tout $x\in I$ et tous $1\le i\ne j\le 3$. Il est clair que $I$ est un sous-groupe d'indice fini de $A_\mathrm{add}$. De plus, pour tout $t\in A$ et tout $x\in I$, l'image par $\pi$ du commutateur $\mathrm{e}_{1,2}(xt)=[\mathrm{e}_{1,3}(x),\mathrm{e}_{3,2}(t)]$ appartient à $I$, puisque c'est le cas pour $\mathrm{e}_{1,3}(x)$ et que $\Gamma$ est distingué. Ainsi, $I$ est un idéal de $A$. Comme $\Gamma$ est abélien, pour tous $x, y\in I$, l'image par $\pi$ du commutateur $\mathrm{e}_{1,2}(xy)=[\mathrm{e}_{1,3}(x),\mathrm{e}_{3,2}(y)]$ est l'identité ; autrement dit, $\mathrm{e}_{1,2}(xy)$ est une homothétie, ce qui entraîne $xy=0$. Ainsi $I$ est un idéal cofini de carré nul de $A$, ce qui achève de montrer l'implication \ref{it-dfE}$\Rightarrow$\ref{it-qpftP}.

Montrons l'implication \ref{it-qpftP}$\Rightarrow$\ref{it-fqf} par exemple dans le cas où $G=\SL_n(A)$ (celui de $\PSL_n(A)$ s'en déduit directement, et le cas de $\GL_n(A)$, qui implique celui de $\operatorname{PGL}_n(A)$, est entièrement similaire). Les lemmes \ref{lm-Op} et \ref{lm-grp_lf} montrent qu'il suffit d'établir (quitte à remplacer $A$ par $A/\nil(A_2)$) que, si $A$ possède un idéal cofini $I$ de carré nul (ce qui implique que $G$ est virtuellement abélien : le noyau du morphisme canonique $\SL_n(A)\to\SL_n(A/I)$ est abélien d'indice fini), alors tout sous-groupe d'indice fini $\Gamma$ de $G$ contient le noyau du morphisme canonique $\SL_n(A)\to\SL_n(A/T)$ pour un certain idéal cofini $T$. Soit $S:=\{x\in I\,|\,\mathrm{e}_{1,2}(x)\in\Gamma\}\cap\{x\in I\,|\,\mathrm{diag}(1+x,1-x,1,\dots,1)\in\Gamma\}$ : $S$ est un sous-groupe cofini de $A_\mathrm{add}$, de sorte que $T:=\bigcap_{i=1}^r a_i^{-1}S$, où $a_1,\dots,a_r$ sont des relevés dans $A$ des éléments de $A/S$, est un idéal cofini. Par conjugaison appropriée, on voit que $\mathrm{e}_{i,j}(x)\in\Gamma$ pour tous $1\le i\ne j\le n$, et que $\mathrm{diag}(1+x_1,\dots,1+x_n)\in\Gamma$ lorsque $x_1,\dots,x_n$ dont des éléments de $T$ de somme nulle. Comme l'idéal $T$ est inclus dans $\rj(A)$, on voit par opérations élémentaires sur les lignes et les colonnes que toute matrice de $\SL_n(A)$ dont l'image dans $\SL_n(A/T)$ est l'identité est un produit de matrices de la forme $\mathrm{e}_{i,j}(x)$ avec $x\in T$ et $\mathrm{diag}(1+x_1,\dots,1+x_n)$, où les $x_i$ sont des éléments de $T$ de somme nulle. Cela termine la démonstration.
\end{proof}

\begin{rema} La proposition~\ref{pr-repGLendofini} ne s'étend à aucun des groupes $\GL_2(A)$, $\SL_2(A)$, $\PSL_2(A)$, ni $\operatorname{PGL}_2(A)$, lorsque $p$ n'est pas nul. Cela provient des cas exceptionnels du lemme~\ref{lm-Op}\,.\ref{it-repev3}. Supposons ainsi que $A$ est un produit infini de copies de $\FF_2$ : alors $\GL_2(A)=\SL_2(A)=\PSL_2(A)=\operatorname{PGL}_2(A)$ est un produit infini de copies de $\GL_2(\FF_2)\simeq\Si_3$. Comme $\Si_3/O_3(\Si_3)\simeq\mathbb{Z}/2$ est abélien, on en déduit que $\GL_2(A)/O_3(\GL_2(A))$ est abélien. Il s'ensuit que, si $K$ est de caractéristique $3$, alors toute représentation irréductible sur $K$ de $\GL_2(A)$ est de dimension finie, bien que $A$ n'admette aucun idéal cofini de carré nul.

Même en caractéristique $p=0$, la proposition~\ref{pr-repGLendofini} ne s'étend pas à $\PSL_2(A)$. Considérons ainsi l'anneau $A=\FF_2[x_1,\dots,x_n,\dots]/(x_1^2,\dots,x_n^2,\dots)$ : il n'est pas semi-primaire, donc a fortiori ne possède aucun idéal cofini de carré nul. Pour autant, on vérifie sans peine que le sous-groupe d'indice fini $\mathrm{Ker}(\PSL_2(A)\twoheadrightarrow\PSL_2(\FF_2))$ (morphisme induit par la réduction modulo l'idéal d'augmentation) est abélien, de sorte que toute représentation irréductible de $\PSL_2(A)$ est de dimension finie. Toutefois, si $2\in A^\times$ et que $\PSL_2(A)$ est virtuellement abélien, alors $A$ possède un idéal cofini de carré nul, comme on le voit à l'aide du même argument que celui utilisé pour $\SL_2(A)$ dans la démonstration de la proposition~\ref{pr-repGLendofini} --- en effet, si le commutateur des matrices élémentaires $\mathrm{e}_{1,2}(x)$ et $\mathrm{e}_{2,1}(y)$ est une homothétie, alors $2xy=0$.
\end{rema}

Nous ignorons si, lorsque les conditions de la proposition~\ref{pr-repGLendofini}, toute représentation $K$-linéaire simple du \emph{monoïde} (non commutatif) $\M_n(A)_\mu$ est de dimension finie. Noter qu'il est équivalent, pour des raisons formelles (cf. par exemple \cite[§\,1.3]{DTV}), que les représentations $K$-linéaires simples des $\M_n(A)_\mu$ soient de dimension finie pour tout $n\in\mathbb{N}$ et que tous les foncteurs simples de $\F(A,K)$ soient à valeurs de dimensions finies.
Il résultera toutefois des résultats de la section suivante que c'est le cas si l'anneau $A_1$ qui apparaît dans la condition~\ref{it-qpftP} est fini.

\section{Polynomialité des foncteurs simples}\label{spfsimp}

Cette section alternera les raisonnements sur les représentations usuelles des groupes linéaires et leurs représentations génériques, en jouant sur les liens entre les deux.

\begin{lemm}\label{lm-rad_Mn-finitude} Soit $n\in\mathbb{N}$. Supposons que, pour tout sous-anneau de type fini $A'$ de $A$, toute représentation $K$-linéaire simple du monoïde multiplicatif (non commutatif) $\M_n(A')$ se factorise par la réduction modulo $\nil(A')$. Alors toute représentation $K$-linéaire simple du monoïde multiplicatif (non commutatif) $\M_n(A)$ se factorise par la réduction modulo $\nil(A)$.
\end{lemm}

\begin{proof} Soient $\xi$ un élément de $K[\M_n(A)]$ dont l'image dans $K[\M_n(A)/\nil(A)]$ par l'application induite par la réduction modulo $\nil(A)$ est nulle et $t$ un élément arbitraire de $K[\M_n(A)]$. Il existe un sous-anneau de type fini $A'$ de $A$ tel que $\xi$ et $t$ appartiennent à $K[\M_n(A')]$. On a clairement $\xi\in\mathrm{Ker}(K[\M_n(A')]\twoheadrightarrow K[\M_n(A')/\nil(A')])$. Ainsi, l'hypothèse implique que $\xi$ appartient au radical de Jacobson de l'algèbre (non commutative) $K[\M_n(A')_\mu]$. Par conséquent, l'élément $1-t\xi$ de $K[\M_n(A')]$ est inversible ; il est donc aussi inversible dans $K[\M_n(A)]$. Cela montre que $\xi$ appartient au radical de Jacobson de l'algèbre (non commutative) $K[\M_n(A)_\mu]$, d'où le lemme.
\end{proof}

\begin{lemm}\label{lm-rad_pol-simple} Soient $\A$ une catégorie additive $A$-linéaire essentiellement petite, $\E$ une catégorie de Grothendieck et $S$ un foncteur polynomial semi-simple de $\fct(\A,\E)$. Alors $\nil(A)\subset\rf(S)$. 
\end{lemm}

\begin{proof} Il suffit de montrer le résultat lorsque $S$ est simple de degré polynomial $d>0$.
Commençons par le cas où $d=1$ : si $a\in\nil(A)$, $S(a)$ est un endomorphisme nilpotent, donc nul, de $S$. Comme $S$ est additif, cela entraîne $a\in\rf(S)$.

Dans le cas général, le multifoncteur $cr_d(S) : \A^d\to\E$ est additif par rapport à chaque variable, et $cr_d(S)$ est semi-simple par la proposition~\ref{pr-psspsd} et le théorème~\ref{th-recPiranv}. Le cas précédent montre que ainsi que le radical de $cr_d(S)$ par rapport à chaque variable contient $\nil(A)$. Il s'ensuit que $\nil(A)\subset\rf(\delta_d^*cr_d(S))$. Or l'unité $S\to\delta_d^*cr_d(S)$ de l'adjonction est non nulle (puisque $cr_d(S)\ne 0$), c'est donc un monomorphisme puisque $S$ est simple. Par conséquent, on a $\rf(S)\supset\rf(\delta_d^*cr_d(S))\supset\nil(A)$.
\end{proof}

Les énoncés suivants font usage de la notion de représentation \emph{polynomiale à la Eilenberg-MacLane} (EML-polynomiale en abrégé) introduite dans \cite[§\,9.1]{DTV}.

\begin{lemm}\label{lm-EMLpol_endofini} Soit $n\in\mathbb{N}^*$.
\begin{enumerate}
    \item Toute représentation $K$-linéaire simple EML-polynomiale du monoïde (non commutatif) $\M_n(A)_\mu$ est de dimension finie sur le centre de son corps d'endomorphismes.
    \item Toute représentation $K$-linéaire irréductible et EML-polynomiale de $\GL_n(A)$ ou de $\SL_n(A)$ est endofinie.\index{termin}{endofini}
\end{enumerate}
\end{lemm}

\begin{proof}
Soient $M$ un $K[\M_n(A)_\mu]$-module à gauche simple et $F$ le \emph{prolongement intermédiaire} associé : alors $F$ est un foncteur simple de $\F(A,K)$ ayant même corps (non commutatif) d'endomorphismes que $M$ et tel que $F(A^n)\simeq M$ dans $K[\M_n(A)_\mu]\Md$ (cf. par exemple \cite[prop.~1.8 et 1.9]{DTV}). De plus, si $M$ est une représentation EML-polynomiale de $\M_n(A)_\mu$, alors $F$ est un foncteur polynomial (cf. par exemple \cite[cor.~1.10 et lemme~2.12]{DTV}). Dans ce cas, $F$ est à valeurs de dimensions finies sur le centre de son corps d'endomorphismes par la proposition~\ref{pr-simpoldfext}. Cela établit le premier point.

Supposons maintenant que $V$ est une représentation EML-polynomiale de $G=\GL_n(A)$ ou $\SL_n(A)$. Par définition, $V$ est la restriction à $G$ d'une représentation EML-polynomiale $M$ du monoïde (non commutatif) $\M_n(A)_\mu$. Si $V$ est irréductible, il en est de même pour $M$, et l'on a des inclusions de corps (non commutatifs)
$$Z(\mathrm{End}_{\M_n(A)}(M))\subset\mathrm{End}_{\M_n(A)}(M)\subset\mathrm{End}_G(V)\,.$$
Comme $M$ est de dimension finie sur $Z(\mathrm{End}_{\M_n(A)}(M))$ d'après le premier point, il s'ensuit que $V$ est de dimension finie sur $\mathrm{End}_G(V)$, ce qui établit le lemme.
\end{proof}

Le résultat principal de cette section est le suivant :

\begin{theo}\label{th-ts_simples_pol} Soient $n\in\mathbb{N}^*$ et $G$ l'un des groupes $\GL_n(A)$ ou $\SL_n(A)$. Considérons les conditions suivantes :
\begin{enumerate}
    \item\label{it-thapC1} tous les foncteurs simples de $\F(A,K)$ sont polynomiaux ;
    \item\label{it-thapC2} toutes les représentations $K$-linéaires simples du monoïde multiplicatif (non commutatif) $\M_n(A)$ sont EML-polynomiales ;
    \item\label{it-thapC3} tout poids $A_\mu\to K_\mu$ est polynomial ;
    \item\label{it-thapC4} l'anneau $A/\nil(A)$ est isomorphe à un produit fini de corps finis de caractéristique $p$.
    \item\label{it-thapC5} toutes les représentations $K$-linéaires irréductibles de $G$ sont EML-polynomiales.
\end{enumerate}

 Alors :
\begin{itemize}
    \item les conditions \ref{it-thapC1}, \ref{it-thapC2}, \ref{it-thapC3} et \ref{it-thapC4} sont équivalentes ;
    \item elles impliquent la condition \ref{it-thapC5} ;
    %si $G=\GL_n(A)$ ou que $K$ contient toutes les racines de l'unité ;
    \item si $n\ge 3$, ou que $n=2$ et $p=0$, alors \ref{it-thapC5} entraîne les autres conditions ;
    \item lorsque ces conditions sont remplies, tout foncteur (resp. toute représentation $K$-linéaire) simple de $\F(A,K)$ (resp. de $\M_n(A)_\mu$, $\GL_n(A)$ ou $\SL_n(A)$) est à valeurs de dimensions finies (resp. de dimension finie) et se factorise par la réduction modulo $\nil(A)$.
\end{itemize}
\end{theo}

\begin{proof}
L'implication \ref{it-thapC1}$\Rightarrow$\ref{it-thapC2} résulte de ce que toute représentation $K$-linéaire simple du monoïde (non commutatif) $\M_n(A)_\mu$ est l'évaluation en $A^n$ d'un foncteur simple de $\F(A,K)$ (cf. par exemple \cite[prop.~1.8 et 1.9]{DTV}) et de ce que l'évaluation en $A^n$ d'un foncteur polynomial définit une représentation EML-polynomiale de $\M_n(A)_\mu$ \cite[prop.~9.3]{DTV}.

Pour tout $w\in\mon(A_\mu,K_\mu)$, la $K$-algèbre (non commutative) $K[\M_n(A)_\mu]/<[\lambda.g]-w(\lambda)[g]\,|\,\lambda\in A, g\in\M_n(A)>$ est non nulle et possède donc au moins un module à gauche simple. Cela fournit une représentation $K$-linéaire simple de $\M_n(A)_\mu$ homogène de poids $w$ (i.e. sur laquelle tout $\lambda\in A\subset\M_n(A)$ opère par la multiplication par $w(\lambda)$). Si une telle représentation est EML-polynomiale, le poids $w$ est nécessairement polynomial, d'où l'implication \ref{it-thapC2}$\Rightarrow$\ref{it-thapC3}.

L'équivalence entre \ref{it-thapC3} et \ref{it-thapC4} a été montrée à la proposition~\ref{pr-poitjspol}.

Le lemme~\ref{lm-rad_pol-simple} montre que tout foncteur simple polynomial de $\F(A,K)$ se factorise par la réduction modulo $\nil(A)$ ; cela vaut également pour les représentations linéaires simples EML-polynomiales de $\M_n(A)_\mu$ grâce à \cite[prop.~9.3]{DTV}.

Montrons tout d'abord l'implication \ref{it-thapC4}$\Rightarrow$\ref{it-thapC1} lorsque l'anneau $A$ est de type fini. Comme $A$ est noethérien, le fait que $A/\nil(A)$ soit fini implique que $A$ est lui-même fini. Sa caractéristique est une puissance de $p$, donc pour tous objets $U$ et $V$ de $\mathbf{P}(A)$, $\mathrm{Hom}(U,V)$ est un $p$-groupe abélien fini, ce qui entraîne que toute fonction de $\mathrm{Hom}(U,V)$ vers un $K$-espace vectoriel est polynomiale \cite[chap.~VI, th.~1.2]{Passi}. Ainsi tout foncteur de $\F(A,K)$ est \ph, ce qui implique que tout foncteur simple de $\F(A,K)$ est polynomial \cite[prop.~2.11]{DTV}.

Montrons maintenant \ref{it-thapC4}$\Rightarrow$\ref{it-thapC1} dans le cas général. Si $A'$ est un sous-anneau de type fini de $A$, alors $A'/\nil(A')$ est également un produit fini de corps finis de caractéristique $p$, de sorte que ce qui précède montre que toutes les représentations $K$-linéaires simples de $\M_m(A')$ sont EML-polynomiales pour tout $m\in\mathbb{N}$, et donc se factorisent par la réduction modulo les éléments nilpotents. Le lemme~\ref{lm-rad_Mn-finitude} permet d'en déduire le même résultat pour les représentations simples des $\M_m(A)$. Comme tout foncteur simple $S$ de $\F(A,K)$ est le prolongement intermédiaire d'une représentation $K$-linéaire simple de $\M_m(A)_\mu$ pour un certain $m\in\mathbb{N}$ (tout $m$ tel que $S(A^m)\ne 0$ convient --- cf. par exemple \cite[cor.~1.10]{DTV}) et que le prolongement intermédiaire d'une représentation EML-polynomiale de $\M_m(A)_\mu$ est un foncteur polynomial (par \cite[cor.~1.10 et lemme~2.12]{DTV} par exemple), cela montre \ref{it-thapC1}.

Nous avons montré l'équivalence entre les quatre premières conditions. Supposons maintenant la condition \ref{it-thapC4} vérifiée. Alors les représentations irréductibles sur $K$ du groupe localement fini $\GL_n(A)$ (resp. $\SL_n(A)$) sont les mêmes que celles de $\GL_n(A/\nil(A))$ (resp. $\SL_n(A/\nil(A))$), grâce au lemme~\ref{lm-grp_lf} et à l'observation que le noyau du morphisme de groupes canonique $\GL_n(A)\twoheadrightarrow\GL_n(A/\nil(A))$ (resp. $\SL_n(A)\twoheadrightarrow\SL_n(A/\nil(A))$) est un $p$-groupe. Autrement dit, pour montrer la deuxième assertion de l'énoncé, il suffit d'établir que toutes les représentations irréductibles sur $K$ de $\GL_n(A)$ (resp. $\SL_n(A)$) sont EML-polynomiales lorsque $A$ est semi-simple fini de caractéristique $p$. Sous cette hypothèse sur $A$, toutes les représentations du monoïde $\M_n(A)$ sont polynomiales (puisque toute fonction de $\M_n(A)$ vers un $\mathbb{Z}/p$-espace vectoriel est polynomiale), donc aussi toutes les représentations de $\GL_n(A)$, puisqu'on peut les prolonger par $0$ en des représentations de $\M_n(A)$. Les représentations $K$-linéaires irréductibles de $\SL_n(A)$ proviennent de représentations irréductibles du groupe algébrique $\SL_n$, dont il est classique qu'elles sont la restriction de représentations de $\GL_n$, ce qui achève de montrer l'implication \ref{it-thapC4}$\Rightarrow$\ref{it-thapC5}.
% Le cas où $K$ n'est pas assez gros se ramène à celui où $K$ est assez gros par arguments standards d'extensions de corps finis.

Supposons maintenant la condition \ref{it-thapC5} vérifiée, et $n\ge 3$, ou $n\ge 2$ et $p=0$. Le lemme~\ref{lm-EMLpol_endofini} et l'implication \ref{it-dfE}$\Rightarrow$\ref{it-qpftP} de la proposition~\ref{pr-repGLendofini} montrent alors que $A$ est le produit d'un anneau $A_1$ contenant un idéal $K$-cotrivial de carré nul et d'un anneau $A_2$ tel que $A_2/\nil(A_2)$ soit un produit fini de corps finis de caractéristique $p$. L'anneau $A_1$ est nécessairement nul, car sinon il existerait un morphisme surjectif d'anneau $A\twoheadrightarrow k$ où $k$ est un corps fini de caractéristique distincte de $p$, auquel cas $\GL_n(k)$ (resp. $\SL_n(k)$), donc $G$ possède manifestement une représentation irréductible non EML-polynomiale (seule la représentation triviale est EML-polynomiale dans ce contexte, et $\SL_n(k)$ n'est pas un $p$-groupe vu que $n\ge 2$). 

La dernière assertion du théorème résulte pour partie de l'observation faite plus haut que les foncteurs simples polynomiaux de $\F(A,K)$, de même que les représentations irréductibles EML-polynomiales de $\M_n(A)$, se factorisent par la réduction modulo $\nil(A)$ ; dans le cas de $\GL_n(A)$ et $\SL_n(A)$, elle se déduit du lemme~\ref{lm-grp_lf}.
Cela termine la démonstration.
\end{proof}

Comme un foncteur polynomial à valeurs de dimensions finies de $\F(A,K)$ est toujours fini, on en déduit :

\begin{coro} Si l'anneau $A/\nil(A)$ est isomorphe à un produit fini de corps finis de caractéristique $p$, alors un foncteur de $\F(A,K)$ est fini si et seulement s'il est polynomial et à valeurs de dimensions finies.
\end{coro}

\begin{rema}
Sous la même hypothèse sur $A$, on ne peut rien conclure quant à $\rf(F)$ (sinon que cet idéal contient une \emph{puissance} de $\nil(A)$ --- cf. corollaire~\ref{cor-kerpp}) si $F$ est un foncteur fini de $\F(A,K)$, comme l'illustre l'exemple~\ref{ex-poiptpAFfini}.
\end{rema}

\begin{rema} Un résultat partiellement analogue au théorème~\ref{th-ts_simples_pol} est le suivant : si l'idéal $\nil(A)$ est égal à son carré (par exemple, $A=k[\mathbb{Z}[1/l]/\mathbb{Z}]$, où $l$ est un nombre premier et $k$ un corps de caractéristique $l$), alors tout foncteur de $\F^\df(A,K)$ se factorise par la réduction modulo $\nil(A)$, sans qu'on ait besoin que $A/\nil(A)$ soit fini, ni de caractéristique $p$ (voir la remarque~\ref{rq-nilcar}).
\end{rema}

\section[Contre-exemple à la décomposition de type $FP$]{Contre-exemple à la décomposition de type $FP$ pour un foncteur de type fini de $\F^\df(A,K)$}\label{scxFP}

Le contre-exemple que nous allons donner est très analogue à celui présenté, dans $\F(\mathbb{Z},\mathbb{Z})$, dans l'exemple~\ref{cex_dec-prim}.

Commençons par une situation générale : supposons que $\A$ est une catégorie additive $A$-linéaire essentiellement petite et que $\alpha$, $\rho$ et $\pi$ sont des foncteurs de $\mathbf{Add}(\A;K)$ s'insérant dans une suite exacte $0\to\alpha\to\rho\to\pi\to 0$. Supposons également que $I$ est un idéal de $A$ inclus dans $\rf(\pi)$ et $\rf(\alpha)$, et que $\alpha$ est le plus grand sous-foncteur de $\rho$ dont le radical contienne $I$.

On considère les foncteurs de $\F(\A;K)$ définis par $Y:=K[\pi]\otimes\rho$, $Z:=\Lambda^2_K\circ\pi$, et $X:=\mathrm{Ker}(Y\xrightarrow{u} Z)$, où $u$ désigne l'épimorphisme
$$K[\pi]\otimes\rho\twoheadrightarrow K[\pi]\otimes\pi\twoheadrightarrow \pi\otimes\pi\twoheadrightarrow\Lambda^2(\pi)\;.$$

Soit $d>1$ un entier. Le foncteur $\varphi_{I,d} : \fct(\A/I,\pol_d(\A;K))\to\F(\A;K)$ composé de l'inclusion $\iota_{I,d} : \fct(\A/I,\pol_d(\A;K))\hookrightarrow\fct(\A/I,\F(\A;K))\simeq\F(\A/I\times\A;K)$ et de la précomposition par le foncteur canonique $\A\to\A/I\times\A$ possède un adjoint à droite $\psi_{I,d}$ ; explicitement, pour $F$ dans $\F(\A;K)$, l'image de $\psi_{I,d}(F)$ par l'inclusion canonique $\fct(\A/I,\pol_d(\A;K))\hookrightarrow\fct(\A/I,\F(\A;K))\simeq\F(\A/I\times\A;K)\hookrightarrow\F(\A\times\A;K)$ (la dernière inclusion étant la précomposition par le foncteur canonique $\A\times\A\twoheadrightarrow\A/I\times\A$) est le plus grand sous-foncteur de $F\circ\oplus$ qui par rapport à la première variable se factorise à travers la réduction modulo $I$ et est polynomial de degré au plus $d$ par rapport à la deuxième variable (cf. remarque~\ref{rq-psi_dI}).

Pour déterminer $\psi_{I,d}(Y)$, on écrit
$$Y(U\oplus V)\simeq K[\pi(U)]\otimes K[\pi(V)]\oplus (\rho(U)\oplus\rho(V))\;;$$
comme, pour tout $U$, le foncteur $V\mapsto K[\pi(U)]\otimes\bar{K}[\pi(V)]\oplus (\rho(U)\oplus\rho(V))$ ne possède aucun sous-foncteur polynomial non nul (cf. exemple~\ref{cex_dec-prim}), le plus grand sous-bifoncteur de $Y\oplus\circ$ polynomial par rapport à la deuxième variable est donné par $(U,V)\mapsto K[\pi(U)]\oplus (\rho(U)\oplus\rho(V))$. Comme $I\subset\rf(\pi)\cap\rf(\alpha)$, pour tout $V$, le foncteur $U\mapsto K[\pi(U)]\oplus (\alpha(U)\oplus\rho(V))$ a également un radical qui contient $I$. Pour voir que c'est le plus grand sous-foncteur de $U\mapsto K[\pi(U)]\oplus (\rho(U)\oplus\rho(V))$ vérifiant cette propriété, on raisonne comme dans l'exemple~\ref{cex_dec-prim}, en invoquant la proposition~\ref{pr-radevtau} et en utilisant que $\alpha$ est le plus grand sous-foncteur de $\rho$ dont le radical contient $I$.

On en déduit aisément un diagramme commutatif
\begin{equation}\label{eq-cou30}
\xymatrix{\varphi_{I,d}\psi_{I,d}(Y)\ar[r]^-\simeq\ar[d] &  K[\pi]\otimes (\alpha\oplus\rho)\ar[d] \\
Y\ar@{=}[r] & K[\pi]\otimes\rho
}
\end{equation}
dont la flèche verticale de gauche est la coünité de l'adjonction et la flèche de droite le produit tensoriel de l'identité de $K[\pi]$ et du morphisme $\alpha\oplus\rho\to\rho$ dont la première composante est l'inclusion et la deuxième l'identité.

On obtient également facilement (de façon analogue à l'exemple~\ref{cex_dec-prim}) un isomorphisme naturel $\iota_{I,d}\psi_{I,d}(Z)(U,V)\simeq\Lambda^2(\pi(U))\oplus\Lambda^2(\pi(V))\oplus (\pi(U)\otimes\pi(V))$ ; et le morphisme $\iota_{I,d}\psi_{I,d}(Y)\to\iota_{I,d}\psi_{I,d}(Z)$ qu'induit $u$ s'identifie, évalué sur $(U,V)$, à
$$K[\pi(U)]\otimes (\alpha(U)\oplus\rho(V))\twoheadrightarrow K[\pi(U)]\otimes\rho(V)\twoheadrightarrow\pi(U)\otimes\pi(V)\hookrightarrow\Lambda^2(\pi(U))\oplus\Lambda^2(\pi(V))\oplus (\pi(U)\otimes\pi(V)).$$
Par conséquent, $\iota_{I,d}\psi_{I,d}(X)$, qui est le noyau du morphisme précédent, s'identifie à la somme directe de $K[\pi]\otimes\alpha$ et du sous-foncteur de $K[\pi]\boxtimes\rho$ somme (non directe) de $(\mathrm{Ker}(K[\pi]\twoheadrightarrow\pi)\boxtimes\rho$ et de $K[\pi]\boxtimes\alpha$.

En utilisant \eqref{eq-cou30}, on en déduit que l'image $F$ de la coünité $\varphi_{I,d}\psi_{I,d}(X)\to X$ est la somme (non directe) des sous-foncteurs $K[\pi]\otimes\alpha$ et $\mathrm{Ker}(K[\pi]\twoheadrightarrow\pi)\otimes\rho$.

On en déduit une suite exacte courte
$$0\to F\to X\to S^2(\pi)\to 0\;;$$
le fait que $F$ soit un sous-objet \emph{strict} de $X$ montre en particulier que $X$ n'appartient pas à l'image essentielle de $\varphi_{I,d}$.

Néanmoins, on a $I\subset\rf(S^2(\pi))$ (car $I\subset\rf(\pi)$ par hypothèse) ; en particulier, $S^2(\pi)$ appartient à l'image essentielle de $\varphi_{I,d}$. Cela montre que $F$ est extension d'un objet de cette image essentielle par un autre.

\begin{exem}\label{cex-tFP} Reprenons l'exemple~\ref{ex-poiptpAFfini} : il fournit, dès lors que $p$ est premier et $K$ infini, un anneau $A$ local, de radical $\mathfrak{m}$ nilpotent et de corps résiduel fini, ainsi qu'une suite exacte courte $0\to\pi\to\rho\to\pi\to 0$, où $\pi$ est un (et en fait le seul) foncteur additif simple de $\F(A,K)$ et $\rf(\pi)=\mathfrak{m}$ mais $\rf(\rho)=0$ ($\rho$ est le foncteur noté $F$ dans l'exemple~\ref{ex-poiptpAFfini}). En particulier, dès que $I$ est un idéal \emph{non nul} de $A$, $\pi$ est le plus grand sous-foncteur de $\rho$ dont le radical contienne $I$. Ainsi, les hypothèses posées en début de section sont vérifiées. Le foncteur $F$ construit ci-avant est à valeurs de dimensions finies, car $\pi$ est à valeurs finies, et de type fini car c'est un sous-foncteur de $K[\pi]\otimes\rho$, qui est de type fini et à valeurs de dimensions finies, donc noethérien (théorème~\ref{th-noeth_cor_gal}).

Ce qui précède montre que $F$ ne possède pas de décomposition de type $FP$ --- mais que $F$ est extension d'un foncteur possédant une décomposition de type $FP$ par un autre foncteur possédant une telle décomposition.
\end{exem}

\backmatter

\printindex{nota}{Index des notations}
\printindex{termin}{Index terminologique}
 
\bibliographystyle{smfalpha}
\bibliography{bib-poids.bib}

\providecommand{\bysame}{\leavevmode ---\ }
\providecommand{\og}{``}
\providecommand{\fg}{''}
\providecommand{\smfandname}{\&}
\providecommand{\smfedsname}{\'eds.}
\providecommand{\smfedname}{\'ed.}
\providecommand{\smfmastersthesisname}{M\'emoire}
\providecommand{\smfphdthesisname}{Th\`ese}
\begin{thebibliography}{LRGS18}

\bibitem[Alb15]{Alb}
{\scshape T.~Albu} -- \emph{Chain conditions in modular lattices with
  applications to {G}rothendieck categories and torsion theories}, Monograph
  Series of the Parana's Mathematical Society, vol.~1, Sociedade Paranaense de
  Matem\'{a}tica, Maring\'{a}, 2015.

\bibitem[Bou81]{Bki2}
{\scshape N.~Bourbaki} -- \emph{\'{E}l\'{e}ments de math\'{e}matique}, Masson,
  Paris, 1981, Alg\`ebre. Chapitres 4 \`a 7.

\bibitem[Bou85]{Bki-AC1}
\bysame , \emph{\'{E}l\'{e}ments de math\'{e}matique}, Masson, Paris, 1985,
  Alg\`ebre commutative. Chapitres 1 \`a 4, Reprint.

\bibitem[Bou12]{Bki}
\bysame , \emph{\'el\'ements de math\'ematique. {A}lg\`ebre. {C}hapitre 8.
  {M}odules et anneaux semi-simples}, Springer, Berlin, 2012, Second revised
  edition of the 1958 edition.

\bibitem[BT73]{Borel-Tits}
{\scshape A.~Borel {\normalfont \smfandname} J.~Tits} -- {\og Homomorphismes
  ``abstraits'' de groupes alg\'{e}briques simples\fg}, \emph{Ann. of Math.
  (2)} \textbf{97} (1973), p.~499--571.

\bibitem[DG24]{DG}
{\scshape A.~Djament {\normalfont \smfandname} T.~Gaujal} -- {\og
  Repr\'esentations g\'en\'eriques des groupes lin\'eaires finis en in\'egale
  caract\'eristique\fg}, \emph{Bull. Soc. Math. France} \textbf{152} (2024),
  no.~2, p.~295--354.

\bibitem[Dja07]{Dja-gr}
{\scshape A.~Djament} -- {\og Foncteurs en grassmanniennes, filtration de
  {K}rull et cohomologie des foncteurs\fg}, \emph{M\'{e}m. Soc. Math. Fr.
  (N.S.)} (2007), no.~111, p.~xxii+213 pp. (2008).

\bibitem[Dja15]{Dja-JLMS}
\bysame , {\og Groupes d'extensions et foncteurs polynomiaux\fg}, \emph{J.
  Lond. Math. Soc. (2)} \textbf{92} (2015), no.~1, p.~63--88.

\bibitem[Dja16]{Dja-bki}
\bysame , {\og La propri\'{e}t\'{e} noeth\'{e}rienne pour les foncteurs entre
  espaces vectoriels [d'apr\`es {A}. {P}utman, {S}. {S}am et {A}.
  {S}nowden]\fg}, \emph{Ast\'{e}risque} (2016), no.~380, S\'{e}minaire
  Bourbaki. Vol. 2014/2015, p.~Exp. No. 1090, 35--60.

\bibitem[DTa]{DT-ext}
{\scshape A.~Djament {\normalfont \smfandname} A.~Touz\'{e}} -- {\og Functor
  homology over an additive category\fg}, Pr\'epublication disponible sur
  https://hal.archives-ouvertes.fr/hal-03432824.

\bibitem[DTb]{DT-add}
\bysame , {\og The homology additive functors in prime characteristic\fg},
  Pr\'epublication disponible sur
  https://hal.archives-ouvertes.fr/hal-04639695.

\bibitem[DTc]{DT-str_pol}
\bysame , {\og Homology of strict polynomial functors over
  $\mathbb{F}_p$-linear additive categories\fg}, Pr\'epublication disponible
  sur https://hal.archives-ouvertes.fr/hal-04646079.

\bibitem[DT23]{DT-schw}
\bysame , {\og Finitude homologique des foncteurs sur une cat\'{e}gorie
  additive et applications\fg}, \emph{Trans. Amer. Math. Soc.} \textbf{376}
  (2023), no.~2, p.~1113--1154.

\bibitem[DT24]{DT-TJM}
\bysame , {\og Sur la noeth\'{e}rianit\'{e} locale des foncteurs
  polynomiaux\fg}, \emph{Tunis. J. Math.} \textbf{6} (2024), no.~1, p.~97--113.

\bibitem[DTV23]{DTV}
{\scshape A.~Djament, A.~Touz\'{e} {\normalfont \smfandname} C.~Vespa} -- {\og
  D\'{e}compositions à la {S}teinberg sur une catégorie additive\fg},
  \emph{Ann. Sci. {\'E}c. Norm. Sup{\'e}r. (4)} \textbf{56} (2023), no.~2,
  p.~427--516.

\bibitem[EK69]{EK69}
{\scshape D.~B.~A. Epstein {\normalfont \smfandname} M.~Kneser} -- {\og
  Functors between categories of vector spaces\fg}, in \emph{Category {T}heory,
  {H}omology {T}heory and their {A}pplications, {III} ({B}attelle {I}nstitute
  {C}onference, {S}eattle, {W}ash., 1968, {V}ol. {T}hree)}, Springer, Berlin,
  1969, p.~154--170.

\bibitem[EML54]{EML}
{\scshape S.~Eilenberg {\normalfont \smfandname} S.~Mac~Lane} -- {\og On the
  groups {$H(\Pi,n)$}. {II}. {M}ethods of computation\fg}, \emph{Ann. of Math.
  (2)} \textbf{60} (1954), p.~49--139.

\bibitem[Eps69]{Ep69}
{\scshape D.~B.~A. Epstein} -- {\og Group representations and functors\fg},
  \emph{Amer. J. Math.} \textbf{91} (1969), p.~395--414.

\bibitem[FS97]{FS}
{\scshape E.~M. Friedlander {\normalfont \smfandname} A.~Suslin} -- {\og
  Cohomology of finite group schemes over a field\fg}, \emph{Invent. Math.}
  \textbf{127} (1997), no.~2, p.~209--270.

\bibitem[Fuc70]{Fuchs}
{\scshape L.~Fuchs} -- \emph{Infinite abelian groups. {V}ol. {I}}, Pure and
  Applied Mathematics, Vol. 36, Academic Press, New York-London, 1970.

\bibitem[Gab62]{Gabriel}
{\scshape P.~Gabriel} -- {\og Des cat\'{e}gories ab\'{e}liennes\fg},
  \emph{Bull. Soc. Math. France} \textbf{90} (1962), p.~323--448.

\bibitem[Gil84]{Gilm}
{\scshape R.~Gilmer} -- \emph{Commutative semigroup rings}, Chicago Lectures in
  Mathematics, University of Chicago Press, Chicago, IL, 1984.

\bibitem[Gre80]{Green}
{\scshape J.~A. Green} -- \emph{Polynomial representations of {${\rm
  GL}_{n}$}}, Lecture Notes in Mathematics, vol. 830, Springer-Verlag,
  Berlin-New York, 1980.

\bibitem[Grz85]{Grz85}
{\scshape P.~Grzeszczuk} -- {\og On {$G$}-systems and {$G$}-graded rings\fg},
  \emph{Proc. Amer. Math. Soc.} \textbf{95} (1985), no.~3, p.~348--352.

\bibitem[Har77]{Hartley77}
{\scshape B.~Hartley} -- {\og Injective modules over group rings\fg}, \emph{Q.
  J. Math., Oxf. II. Ser.} \textbf{28} (1977), p.~1--29.

\bibitem[Har83]{Hartley83}
\bysame , {\og Locally finite groups whose irreducible modules are finite
  dimensional\fg}, \emph{Rocky Mountain J. Math.} \textbf{13} (1983), no.~2,
  p.~255--263.

\bibitem[IMR16]{IMR16}
{\scshape M.~C. Iovanov, Z.~Mesyan {\normalfont \smfandname} M.~L. Reyes} --
  {\og Infinite-dimensional diagonalization and semisimplicity\fg},
  \emph{Israel J. Math.} \textbf{215} (2016), no.~2, p.~801--855.

\bibitem[Isb59]{Isb59}
{\scshape J.~R. Isbell} -- {\og On the multiplicative semigroup of commutative
  ring\fg}, \emph{Proc. Amer. Math. Soc.} \textbf{10} (1959), p.~908--909.

\bibitem[Kuh94a]{Ku1}
{\scshape N.~J. Kuhn} -- {\og Generic representations of the finite general
  linear groups and the {S}teenrod algebra. {I}\fg}, \emph{Amer. J. Math.}
  \textbf{116} (1994), no.~2, p.~327--360.

\bibitem[Kuh94b]{Ku2}
\bysame , {\og Generic representations of the finite general linear groups and
  the {S}teenrod algebra. {II}\fg}, \emph{$K$-Theory} \textbf{8} (1994), no.~4,
  p.~395--428.

\bibitem[Kuh95]{Ku3}
\bysame , {\og Generic representations of the finite general linear groups and
  the {S}teenrod algebra. {III}\fg}, \emph{$K$-Theory} \textbf{9} (1995),
  no.~3, p.~273--303.

\bibitem[Kuh97]{K97}
\bysame , {\og Invariant subspaces of the ring of functions on a vector space
  over a finite field\fg}, \emph{J. Algebra} \textbf{191} (1997), no.~1,
  p.~212--227.

\bibitem[Kuh15]{Ku-adv}
\bysame , {\og Generic representation theory of finite fields in nondescribing
  characteristic\fg}, \emph{Adv. Math.} \textbf{272} (2015), p.~598--610.

\bibitem[LRGS18]{LRGS}
{\scshape W.~Lowen, J.~Ramos~Gonz\'{a}lez {\normalfont \smfandname}
  B.~Shoikhet} -- {\og On the tensor product of linear sites and {G}rothendieck
  categories\fg}, \emph{Int. Math. Res. Not. IMRN} (2018), no.~21,
  p.~6698--6736.

\bibitem[Mes18]{Mes18}
{\scshape Z.~Mesyan} -- {\og Infinite-dimensional triangularization\fg},
  \emph{J. Pure Appl. Algebra} \textbf{222} (2018), no.~7, p.~1529--1547.

\bibitem[ML98]{ML-cat}
{\scshape S.~Mac~Lane} -- \emph{Categories for the working mathematician},
  second \smfedname, Graduate Texts in Mathematics, vol.~5, Springer-Verlag,
  New York, 1998.

\bibitem[Pas77]{Passman}
{\scshape D.~S. Passman} -- \emph{The algebraic structure of group rings},
  {John} {Wiley} \& {Sons}, New {York} etc., 1977.

\bibitem[Pas79]{Passi}
{\scshape I.~B.~S. Passi} -- \emph{Group rings and their augmentation ideals},
  Lecture Notes in Mathematics, vol. 715, Springer, Berlin, 1979.

\bibitem[Pir88]{Pira88}
{\scshape T.~I. Pirashvili} -- {\og Polynomial functors\fg}, \emph{Trudy
  Tbiliss. Mat. Inst. Razmadze Akad. Nauk Gruzin. SSR} \textbf{91} (1988),
  p.~55--66.

\bibitem[Pop73]{Pop}
{\scshape N.~Popescu} -- \emph{Abelian categories with applications to rings
  and modules}, Academic Press, London-New York, 1973, London Mathematical
  Society Monographs, No. 3.

\bibitem[Pow98]{GP-cow}
{\scshape G.~M.~L. Powell} -- {\og The structure of indecomposable injectives
  in generic representation theory\fg}, \emph{Trans. Amer. Math. Soc.}
  \textbf{350} (1998), no.~10, p.~4167--4193.

\bibitem[PS17]{PSam}
{\scshape A.~Putman {\normalfont \smfandname} S.~V. Sam} -- {\og Representation
  stability and finite linear groups\fg}, \emph{Duke Math. J.} \textbf{166}
  (2017), no.~13, p.~2521--2598.

\bibitem[PS21]{Psch}
{\scshape L.~Positselski {\normalfont \smfandname} O.~M. Schn{\"u}rer} -- {\og
  Unbounded derived categories of small and big modules: is the natural functor
  fully faithful?\fg}, \emph{J. Pure Appl. Algebra} \textbf{225} (2021),
  no.~11, p.~23, Id/No 106722.

\bibitem[Sco00]{Sco}
{\scshape A.~Scorichenko} -- {\og Stable {K}-theory and functor homology over a
  ring\fg}, \smfphdthesisname, Evanston, 2000.

\bibitem[SGA77]{SGA5}
\emph{Cohomologie {$l$}-adique et fonctions {$L$}} -- Lecture Notes in
  Mathematics, Vol. 589, Springer-Verlag, Berlin-New York, 1977, S\'{e}minaire
  de G\'{e}ometrie Alg\'{e}brique du Bois-Marie 1965--1966 (SGA 5), Edit\'{e}
  par Luc Illusie.

\bibitem[Sha77]{Sharp}
{\scshape R.~Y. Sharp} -- {\og The dimension of the tensor product of two field
  extensions\fg}, \emph{Bull. Lond. Math. Soc.} \textbf{9} (1977), p.~42--48.

\bibitem[SS17]{SamSn}
{\scshape S.~V. Sam {\normalfont \smfandname} A.~Snowden} -- {\og Gr\"obner
  methods for representations of combinatorial categories\fg}, \emph{J. Amer.
  Math. Soc.} \textbf{30} (2017), no.~1, p.~159--203.

\bibitem[SS18]{Sch-Si}
{\scshape P.~Schenzel {\normalfont \smfandname} A.-M. Simon} --
  \emph{Completion, \v{C}ech and local homology and cohomology}, Springer
  Monographs in Mathematics, Springer, Cham, 2018, Interactions between them.

\bibitem[Tou13]{TouzeRingel}
{\scshape A.~Touz\'{e}} -- {\og Ringel duality and derivatives of non-additive
  functors\fg}, \emph{J. Pure Appl. Algebra} \textbf{217} (2013), no.~9,
  p.~1642--1673.

\bibitem[Tou14]{T-AIF}
\bysame , {\og Bar complexes and extensions of classical exponential
  functors\fg}, \emph{Ann. Inst. Fourier (Grenoble)} \textbf{64} (2014), no.~6,
  p.~2563--2637.

\bibitem[Tou17]{T-torsion}
\bysame , {\og A functorial control of integral torsion in homology\fg},
  \emph{Fund. Math.} \textbf{237} (2017), no.~2, p.~135--163.

\bibitem[Tou21]{T-expo}
\bysame , {\og On the structure of graded commutative exponential functors\fg},
  \emph{Int. Math. Res. Not. IMRN} (2021), no.~17, p.~13305--13415.

\bibitem[Vam78]{Vamos}
{\scshape P.~Vamos} -- {\og On the minimal prime ideals of a tensor product of
  two fields\fg}, \emph{Math. Proc. Camb. Philos. Soc.} \textbf{84} (1978),
  p.~25--35.

\bibitem[WJ87]{WJ87}
{\scshape P.~Wauters {\normalfont \smfandname} E.~Jespers} -- {\og When is a
  semigroup ring of a commutative semigroup local or semilocal?\fg}, \emph{J.
  Algebra} \textbf{108} (1987), no.~1, p.~188--194.

\bibitem[Woo71]{Woo}
{\scshape S.~M. Woods} -- {\og On perfect group rings\fg}, \emph{Proc. Amer.
  Math. Soc.} \textbf{27} (1971), p.~49--52.

\end{thebibliography}

 \end{document}